# HF=HM IV: The Seiberg-Witten Floer homology and ech correspondence

Cagatay Kutluhan, Yi-Jen Lee and Clifford Henry Taubes

ABSTRACT: This is the fourth of five papers that construct an isomorphism between the Seiberg-Witten Floer homology and the Heegaard Floer homology of a given compact, oriented 3-manifold. The isomorphism is given as a composition of three isomorphisms; the first of these relates a version of embedded contact homology on an an auxillary manifold to the Heegaard Floer homology on the original. The second isomorphism relates the relevant version of the embedded contact homology on the auxilliary manifold with a version of the Seiberg-Witten Floer homology on this same manifold. The third isomorphism relates the Seiberg-Witten Floer homology on the auxilliary manifold with the appropriate version of Seiberg-Witten Floer homology on the original manifold. The paper describes the second of these isomorphisms.

This is the fourth of five papers that together supply an isomorphism between Ozsvath and Szabo's [OS1], [OS2] Heegaard Floer homology of any given compact, oriented 3-manifold and a version of the Seiberg-Witten Floer homology of the same 3-manifold. The existence of an isomorphism between these respective Floer homologies is stated as the Main Theorem in [KLTI]. A particular isomorphism is described in [KLTI] that can be written as the concatenation of three separate isomorphisms which involve an auxilliary manifold that is obtained from the original by connect summing a certain number of copies of $S^1 \times S^2$. As noted in Section 3 of [KLTI], the middle isomorphism in this concatenation identifies a version of the Seiberg-Witten Floer homology on the connect sum manifold with a version of Michael Hutchings' embedded contact homology [Hu1] on this same connect sum The relevant version of this embedded contact homology is described in Section 2.2 of [KLTI] and in the Appendix of [KLTII]. This homology is also the focus in [KLTIII]. The isomorphism between the respective Seiberg-Witten Floer and embedded contact homologies on the connect sum manifold is asserted by Theorems 3.3 and 3.4 in [KLTI]. The latter reference deduces its Theorem 3.4 from Theorem 3.3; this paper proves Theorem 3.3 from [KLTI]. Section 1d explains why the latter theorem is consequence of Theorem 1.5 in Section 1d. Respective proofs are also given below for Theorems 3.1 and 3.2 in [KTLI]. These are seen in Section 1c to be consequences of Propositions 1.1-1.4.

What follows is a table of contents for this article.

Most of what is done in [KLTII] and [KLTIII] is not relevant for what follows. Even so, certain results and constructions from these papers are needed. In particular, the geometry needed to define the appropriate versions of the Seiberg-Witten Floer homology and the embedded contact homology is described in [KLTII] and [KLTIII]. Section 1a provides a summary of this geometry.

The following conventions are used throughout the remainder of this paper: Section numbers, equation numbers, and other references from [KLTI], [KLTII] and [KLTIII] are distinguished from those in this paper by the use of the Roman numerals I, II and III as a prefix. For example, 'Section II.1' refers to Section 1 in [KLTII]. Note also that the convention here as in [KLTII] and [KLTIII] is to use $c_0$ to denote a constant in $(1, \infty)$ whose value is independent of all relevant parameters. The value of $c_0$ can increase between subsequent appearances. A second convention used here and in [KLTII] and [KLTIII] concerns a function that is denoted by $\chi$. The latter is a fixed, non-increasing function on $\mathbb{R}$ that equals 1 on $(-\infty, 0]$ and equals 0 on $[1, \infty)$.

**Acknowledgements**

The first author was supported in part by a National Science Foundation Postdoctoral Fellowship under Award 1103795. The second and third authors are also supported in part by grants from the National Science Foundation.

## 1. Embedded contact homology and Seiberg-Witten Floer homology

Let M denote the given, compact oriented 3-manifold. A self-indexing Morse function for M and certain auxilliary data is used in [KLTII] to construct a second manifold, Y, as a connect sum of M with copies of $S^1 \times S^2$. The manifold Y has two orientations, one coming from the part it shares with M and then the opposite orientation. The theorems in [KLTI] compare certain Seiberg-Witten Floer homologies on the M-oriented version with a certain sort of embedded contact homology on the oppositely oriented version. These respective homologies need certain geometric data for their definition. This section starts with a description of the necessary data. It then briefly describes the relevant version of embedded contact homology and the relevant version of Seiberg-Witten Floer homology. It ends by restating Theorems 3.1-3.4 from [KLTI].

### a) The geometry of Y

The construction of Y from M is described in Section II.1. This subsection summarizes the salient feature of Y.

*Part 1*: The construction of Y starts with the choice of a self-indexing Morse function, $f$: $M \to [0,3]$ with one index 0 critical points, one index 3 critical point and some index 1 and index 2 critical points. The number of index 1 (and thus index 2) critical points is denoted by $G$. The manifold Y is diffeomorphic to the connect sum of M with $G+1$ copies of $S^1 \times S^2$. The manifold Y is oriented so that the part from M has the orientation opposite from M's orientation. Note that [KLTI] uses $\overline{Y}$ to denote this orientation of the connect sum and uses Y to denote the connect sum with the orientation induced from M.

Also needed from M is the choice of a class in $H^2(M;\mathbb{Z})$ which defines a homomorphism from $H_2(M;\mathbb{Z})$ to $2\mathbb{Z}$. This class is denoted in what follows by $c_{1M}$. A $\text{Spin}^{\mathbb{C}}$ structure will be chosen momentarily, and its first Chern class will play the role of $c_{1M}$. Needed also is a chosen pairing between the set of index 1 critical points of $f$ and the set of index 2 critical points of $f$. The resulting set of G pairs is denoted by $\Lambda$. An element $\mathfrak{p} \in \Lambda$ is written as an ordered pair of points with it understood that the first entry is the index 1 critical point of $f$ and the second entry is the index 2 critical point of $f$. Various other Morse-theoretic items from M are needed to constuct Y and its geometry, but these others play minor roles in this paper.

The definition of Y required the choice of positive numbers which are denoted by $\delta_*$ and R. This $\delta_*$ is from (0, 1); it is determined by the chosen function $f$. Meanwhile, R has the lower bound $-100\ln\delta_*$. This constant R has no apriori upper bound, and the freedom to take R as large as needed is exploited in [KLTII], [KLTIII] and in the constructions to come in this article.

The construction of the geometry needed for the embedded contact geometry chain complex required the choice of two additional positive numbers which are denoted by $\delta$ and $x_0$. The trio $(\delta, x_0, R)$ are constrained by the requirements that $\delta < \delta_*$, $x_0 < \delta^3$ and $R \geq -c_0 \ln x_0$. The choice of $\delta$ determines an upper bound for $x_0$, and that the choice of $x_0$ subject to this upper bound then determines a lower bound for R. Constants $\delta$, $x_0$, and R that satisfy these bounds are said to be *appropriate*. The freedom to take $\delta$ as small as desired is also exploited in [KLTII], [KLTIII] and in what follows.

*Part 2*: The manifold Y is constructed by attaching G+1 handles to M. In particular, Y is written as the union of sets $M_\delta \cup \mathcal{H}_0 \cup (\cup_{\mathfrak{p}\in\Lambda} \mathcal{H}_\mathfrak{p})$ where the notation is as follows: First, $M_\delta$ is the complement in M of $2(G+1)$ disjoint balls about the critical points of $f$. What is written as $\mathcal{H}_0$ is a 1-handle and so diffeomorphic to $[-1, 1] \times S^2$. It intersects $M_\delta$ near $\{-1\} \times S^2$ as an annulus in a ball centered on the index 3 critical point of $f$; and it intersects $M_\delta$ near $\{1\} \times S^2$ as an annulus in a ball centered on the index 0 critical point of $f$. Meanwhile, the various $\mathfrak{p} \in \Lambda$ version of $\mathcal{H}_\mathfrak{p}$ are 1-handles and so each is diffeomorphic to $[-1, 1] \times S^2$ also. These are pairwise disjoint and disjoint from $\mathcal{H}_0$. Any given $\mathfrak{p} \in \Lambda$ version of $\mathcal{H}_\mathfrak{p}$ intersects $M_\delta$ near $\{-1\} \times S^2$ as an annulus in a ball centered around $\mathfrak{p}$'s index 2 critical point and it intersects $M_\delta$ near $\{1\} \times S^2$ as an annulus in a ball cented around $\mathfrak{p}$'s index 1 critical point.

The handle $\mathcal{H}_0$ and those from the set $\{\mathcal{H}_\mathfrak{p}\}_{\mathfrak{p}\in\Lambda}$ have preferred coordinates, these denoted by $(u, (\theta, \phi))$ where $(\theta, \phi)$ are spherical coordinates for the $S^2$ factor and where u is the Euclidean coordinate for the closed interval $[-R - \ln(7\delta_*), R + \ln(7\delta_*)]$. The function $f$ appears in these coordinates near the $u < 0$ end of $\mathcal{H}_0$ as $f = 3 - e^{2(u+R)}$ and near the $u > 0$ end of $\mathcal{H}_0$ as $f = e^{2(u-R)}$. Meanwhile, the function $f$ appears near the respective negative and positive u ends of any given $\mathfrak{p} \in \Lambda$ version of $\mathcal{H}_\mathfrak{p}$ as

$$f = 2 - e^{2(u+R)}(1 - 3\cos^2\theta) \quad \textit{and} \quad f = 1 + e^{2(u-R)}(1 - 3\cos^2\theta) \tag{1.1}$$

The 3-form $du \sin\theta \, d\theta \, d\phi$ gives the Y-orientation to each handle. Orient the cross-sectional spheres in each 1-handle using the 2-form $\sin\theta \, d\theta \, d\phi$.

*Part 3*: The definition of the relevant version of Hutchings' embedded contact homology uses a pair $(w, a)$ of 2-form and 1-form on Y. The latter define a *stable Hamiltonian* structure, which is to say that $w$ is closed, $a \wedge w$ is nowhere zero and defines Y's orientation, and $da \subset \mathrm{span}(w)$. The vector field that generates the kernel of $w$ and has pairing 1 with $a$ is denoted by $v$. The salient features of $w$, $a$ and $v$ are listed in the upcoming (1.3). This equation refers to auxilliary functions $x_0, \chi_+, \chi_-$, f and g. These are functions $[-R - \ln(7\delta_*), R + \ln(7\delta_*)]$ that are defined using the chosen function $\chi$. By way of a reminder, $\chi$ is a smooth, non-increasing function on $\mathbb{R}$ that is 1 on $(-\infty, 0]$ and equals 0 on $[1, \infty)$. The aforementioned five functions are

- $x = x_0 \chi(|u| - R - \ln\delta + 12)$
- $\chi_+ = \chi(-u - \frac{1}{4}R)$ *and* $\chi_- = \chi(u - \frac{1}{4}R)$.
- $f = x + 2(\chi_+ e^{2(u-R)} + \chi_- e^{-2(u+R)})$ *and* $g = (\chi_+ e^{2(u-R)} - \chi_- e^{-2(u+R)})$.

(1.2)

What follows is the promised list.

- ON $M_\delta$: *The 2-form* $w$ *on* $M_\delta$ *is nowhere zero on the kernel of the 1-form* $df$ *and* $v$ *here is a certain pseudogradient vector field for* $f$.

- IN THE HANDLE $\mathcal{H}_0$: *The 2-form* $w$ *and the vector field* $v$ *on* $\mathcal{H}_0$ *are*

$$w = \sin\theta \, d\theta \wedge d\phi \quad \textit{and} \quad v = \frac{1}{2(\chi_+ e^{2(u-R)} + \chi_- e^{-2(u+R)})} \frac{\partial}{\partial u} .$$

- IN THE HANDLES $\{\mathcal{H}_{\mathfrak{p}}\}_{\mathfrak{p} \in \Lambda}$: *Fix* $\mathfrak{p} \in \Lambda$. *The trio* $a$, $w$ *and* $v$ *on* $\mathcal{H}_{\mathfrak{p}}$ *are*

$$\begin{aligned} a &= (x + g')(1 - 3\cos^2\theta)\, du - \sqrt{6}\, f \cos\theta \sin^2\theta \, d\phi + 6 g \cos\theta \sin\theta \, d\theta , \\ w &= 6 x \cos\theta \sin\theta \, d\theta \wedge du - \sqrt{6}\, d\{f \cos\theta \sin^2\theta \, d\phi\} , \\ v &= \hat{c}_v^{-1}\{f(1 - 3\cos^2\theta)\, \partial_u - \sqrt{6}\, x \cos\theta \, \partial_\phi + f' \cos\theta \sin\theta \, \partial_\theta\} . \end{aligned}$$

(1.3)

Here, $\hat{c}_v = (x + g')f(1 - 3\cos^2\theta)^2 + 6(xf + gf')\cos^2\theta \sin^2\theta$ is a positive function of $(u, \theta)$.

An additional property of $w$ plays a central role in the story to come. To say more about this, introduce the direct sum decomposition

$$H_2(Y; \mathbb{Z}) = H_2(M; \mathbb{Z}) \oplus H_2(\mathcal{H}_0; \mathbb{Z}) \oplus (\oplus_{\mathfrak{p} \in \Lambda} H_2(\mathcal{H}_{\mathfrak{p}}; \mathbb{Z}))$$

(1.4)

that comes via Mayer-Vietoris by writing $Y = \mathcal{M}_\delta \cup \mathcal{H}_0 \cup (\cup_{\mathfrak{p}\in\Lambda}\mathcal{H}_\mathfrak{p})$. The summands in (1.4) that correspond to the various 1-handles are isomorphic to $\mathbb{Z}$; and any oriented, cross-sectional sphere is a generator.

The additional property concerns the cohomology class defined by *w*. This class is determined by what follows: Integration of *w* over closed 2-cycles defines the linear map from $H_2(Y;\mathbb{Z})$ to $\mathbb{Z}$ that has value 2 on the generator of $H_2(\mathcal{H}_0;\mathbb{Z})$; it has value zero on each $\mathfrak{p} \in \Lambda$ version of $H_2(\mathcal{H}_\mathfrak{p};\mathbb{Z})$; and it acts on the $H_2(M;\mathbb{Z})$ summand in (1.7) as the pairing with the chosen class $c_{1M}$.

*Part 4*: A particular closed integral curve of the vector field *v* plays a distinguished role in the embedded contact homology story. This curve is denoted by $\gamma^{(z_0)}$ here and in the other papers in this series. The curve $\gamma^{(z_0)}$ is disjoint from $\cup_{\mathfrak{p}\in\Lambda}\mathcal{H}_\mathfrak{p}$ and it crossed $\mathcal{H}_0$ so as to have intersection number 1 with each cross-sectional sphere. Note in this regard that the convention here and in what follows is to orient the integral curves of *v* using *v* for the oriented unit tangent vector. This curve intersects $\Sigma$ in precisely one point. The latter is denoted by $z_0$.

A pair of additional 1-forms enter the story. These are denoted by $\upsilon_\Diamond$ and $\hat{a}$.

- THE 1-FORM $\upsilon_\Diamond$: *The 1-form* $\upsilon_\Diamond$ *is closed and is such that* $\upsilon_\Diamond \wedge w \geq 0$. *Furthermore,* $\upsilon_\Diamond \wedge w = 0$ *only where both* $u = 0$ and $1 - 3\cos^2\theta = 0$ *on each* $\mathfrak{p} \in \Lambda$ *version of* $\mathcal{H}_\mathfrak{p}$. *This 1-form equals* $df$ *on* $M_\delta$, *it is given by* $\upsilon_\Diamond = 2(\chi_+ e^{2(|u|-R)} + \chi_- e^{-2(|u|+R)})\,du$ *on* $\mathcal{H}_0$, *and it is given by* $d((\chi_+ e^{2(u-R)} - \chi_- e^{-2(u+R)})(1 - 3\cos^2\theta))$ *on any given* $\mathfrak{p} \in \Lambda$ *version of* $\mathcal{H}_\mathfrak{p}$.

(1.5)

The definition $\hat{a}$ refers to the function $\chi_\delta$ that is defined on any given $\mathfrak{p} \in \Lambda$ version of $\mathcal{H}_\mathfrak{p}$ by the rule $\chi_\delta = \chi(|u| - R - \ln\delta + 10)$.

- THE 1-FORM $\hat{a}$: *The 1-form* $\hat{a}$ *has pairing* 1 *with* $v$ *and is such that* $\hat{a} \wedge w > 0$. *This 1-form is equal to* $\upsilon_\Diamond$ *on* $M_\delta \cup \mathcal{H}_0$ *and it is equal to* $\chi_\delta a + (1 - \chi_\delta)\upsilon_\Diamond$ *on any given* $\mathfrak{p} \in \Lambda$ *version of* $\mathcal{H}_\mathfrak{p}$.

(1.6)

The kernel of the 1-form $\hat{a}$ defines a 2-plane subbundle in TY on which *w* is non-degenerate. When oriented by *w*, the bundle $\ker(\hat{a})$ has an Euler class which evaluates as 2 on the generator of the $H_2(\mathcal{H}_0;\mathbb{Z})$ summand in (1.4) and evaluates as -2 on the generator of each $\mathfrak{p} \in \Lambda$ summand $H_2(\mathcal{H}_\mathfrak{p};\mathbb{Z})$. The vector field *v* has pairing 1 with $\hat{a}$ also.

Various other geometric properties of Y are introduced as needed in what follows.

*Part* 5: The almost complex geometry of $\mathbb{R} \times Y$ is define by an almost complex structure, this denoted by J. The latter is constrained in various ways; most of the constraints are given in Part 1 of Section II.3a and Section III.1c. The upcoming (1.7) reviews various features of J. This equation uses $s$ to denote the Euclidean coordinate on the $\mathbb{R}$ factor of $\mathbb{R} \times Y$.

- J *maps the Euclidean tangent vector* $\partial_s$ *to the* $\mathbb{R}$ *factor of* $\mathbb{R} \times Y$ *to* $v$.
- J *is not changed by constant translations of the coordinate* $s$ *on* $\mathbb{R} \times Y$.
- J *preserves the kernel of the 1-form* $\hat{a}$*; and its restriction to this 2-plane field defines the orientation given by* $w$.
- J *on* $\mathbb{R} \times \mathcal{H}_0$ *and on any given* $\mathfrak{p} \in \Lambda$ *version of* $\mathbb{R} \times \mathcal{H}_{\mathfrak{p}}$ *is invariant with respect to constant translations of the* $\mathbb{R}/(2\pi\mathbb{Z})$ *coordinate* $\phi$.

(1.7)

It is a consequence of (1.7) that the 2-form $\hat{\omega} = ds \wedge \hat{a} + w$ on $\mathbb{R} \times Y$ is compatible with J. This is to say that the bilinear form $\hat{\omega}(\cdot, J(\cdot))$ on $T(\mathbb{R} \times Y)$ defines a Riemannian metric. Note in particular that this metric has the form $ds^2 + g_Y$ with $g_Y$ being a metric on TY that makes $v$ a unit vector that is orthogonal to the kernel of $\hat{a}$. The corresponding metric on T*Y gives $\hat{a}$ norm one and is such that the Hodge star of $\hat{a}$ is $w$.

These respective metrics on $\mathbb{R} \times Y$ and Y are used implicitly in what follows.

### b) Embedded contact homology on Y

The appendix in [KLTII] describes the relevant version of embedded contact homology on Y. More is said about the chain complex and its homology in Section III.1b and Section III.9. This subsection provides a very brief summary of what is said in these section of [KLTII] and [KLTIII]. The summary here comprises Parts 2-4 of the five parts of this subsection. The first part of the subsection constitues a digression that concerns $\text{Spin}_{\mathbb{C}}$ structures on M and Y. The final part summarizes some observations from [KLTII] and [KLTIII] that are particularly relevant in the subsequent sections of this paper.

*Part 1*: A $\text{Spin}_{\mathbb{C}}$ structure on M is chosen whose associated first Chern class is the chosen class $c_{1M}$. The chosen $\text{Spin}_{\mathbb{C}}$ structure is fixed for the remainder of this article. The $\text{Spin}_{\mathbb{C}}$ structure on M determines in a canonical fashion a corresponding $\text{Spin}_{\mathbb{C}}$ structure on Y. This is done using a version of Mayer-Vietoris with the decomposition of Y as $M_\delta \cup \mathcal{H}_0 \cup (\cup_{p\in\Lambda}\mathcal{H}_p)$. The first Chern class of the resulting $Spin_{\mathbb{C}}$ structure on Y has pairing 2 with the generator of the $H_2(\mathcal{H}_0;\mathbb{Z})$ summand in (1.4) and it has pairing 0 with each of the $\mathfrak{p} \in \Lambda$ labled summands. The pairing with the $H_2(M;\mathbb{Z})$ summand is that of the first Chern class of the $Spin_{\mathbb{C}}$ structure on M, which is to say that of $c_{1M}$. It follows as a consequence that the image in $H^2(Y;\mathbb{R})$ of the first Chern class of the $Spin_{\mathbb{C}}$ structure on Y is the class define by the closed form *w*.

The image of $H_2(M;\mathbb{Z})$ in $\mathbb{Z}$ given by the pairing with $c_{1M}$ is a subgroup of $\mathbb{Z}$. If $c_{1M}$ is not torsion, use $p_M$ denote the largest integer that divides all of its elements. Note that $p_M$ is in all cases even.

*Part 2*: The $\mathbb{Z}$-module that serves as the embedded complex homology chain complex is defined using a certain principal $\mathbb{Z}$-bundle over a set that is denoted by $\mathcal{Z}_{ech,M}$. The set $\mathcal{Z}_{ech,M}$ is described in Proposition II.2.8. The principal $\mathbb{Z}$-bundle is denoted by $\hat{\mathcal{Z}}_{ech,M}$ and described in Section II.1f and in Part 4 of Section III.1b. The embedded contact homology chain complex is denoted by $\mathbb{Z}(\hat{\mathcal{Z}}_{ech,M})$.

By way of a reminder, an element in $\mathcal{Z}_{ech,M}$ is a set, $\Theta$, consisting of some number of closed integral curves of $v$ that lie entirely in the union of the $f \in (1, 2)$ part of $M_\delta$ and the various $\mathfrak{p} \in \Lambda$ versions of $\mathcal{H}_p$. In particular, the union of the curves that comprise such a set $\Theta$ intersect each $\mathfrak{p} \in \Lambda$ version of $\mathcal{H}_p$ in at most three components. There exists in all cases one component of this intersection that lies entirely in the $1 - 3\cos^2\theta > 0$ part of $\mathcal{H}_p$ as an arc that crosses $\mathcal{H}_p$ from the $u = -R - \ln(7\delta_*)$ end to the $u = R + \ln(7\delta_*)$ end. The locus in $\mathcal{H}_p$ where both $u = 0$ and $1 - 3\cos^2\theta = 0$ is a disjoint union of two closed integral curves of $v$, and one or both of these curves can also appear in $\Theta$. The curve with $u = 0$ and $\cos\theta = \frac{1}{\sqrt{3}}$ is denoted by $\hat{\gamma}_p^+$ and the curve where $u = 0$ and $\cos\theta = -\frac{1}{\sqrt{3}}$ by $\hat{\gamma}_p^-$.

If $\gamma$ is used to denote a closed integral curve of $v$, then $[\gamma]$ is used to denote both the oriented cycle defined by $\gamma$ and the corresponding element in $H_1(Y;\mathbb{Z})$ where it is understood that $\gamma$ is oriented by $v$. Meanwhile, $[\Theta] = \sum_{\gamma\in\Theta} [\gamma]$ is used to denote both a sum of oriented 1-cycles and the corresponding homology class. The latter is fixed by the chosen $Spin_{\mathbb{C}}$-structure; this class is the same for all elements in $\mathcal{Z}_{ech,M}$.

The principal $\mathbb{Z}$ bundle $\hat{\mathcal{Z}}_{ech,M} \to \mathcal{Z}_{ech,M}$ is defined after choosing a fiducial element $\Theta_0 \in \mathcal{Z}_{ech,M}$. The fiber of $\hat{\mathcal{Z}}_{ech,M}$ over a given element $\Theta \in \mathcal{Z}_{ech,M}$ is identified with the set of equivalence classes of pairs of the form $(\Theta, Z)$ where Z is a relative cycle in $H_2(Y; [\Theta] - [\Theta_0])$. The equivalence relation is defined using the pairing with the Poincare dual of the homology class of the closed integral curve $\gamma^{(z_0)}$. This pairing defines a homomorphism from the $\mathbb{Z}$-module of closed 2-cycles to $\mathbb{Z}$ that is denoted by $[\gamma^{(z_0)}]^{Pd}(\cdot)$. The equivalence relation that defines $\hat{\mathcal{Z}}_{ech,M}$ has $(\Theta, Z) \sim (\Theta', Z')$ if and only if $\Theta = \Theta'$ and also $[\gamma^{(z_0)}]^{Pd}(Z - Z') = 0$. The principal bundle projection map sends an equivalence class $(\Theta, Z)$ to $\Theta$. The element $1 \in \mathbb{Z}$ acts to send $(\Theta, Z)$ to $(\Theta, Z + [S_0])$ where $[S_0]$ is the u = 0 sphere in $\mathcal{H}_0$.

The module $\mathbb{Z}(\hat{\mathcal{Z}}_{ech,M})$ has a relative $\mathbb{Z}$ grading when $c_{1M}$ is torsion, and it has a relative $\mathbb{Z}/p_M\mathbb{Z}$ grading otherwise. The grading rule comes via a corresponding grading of the generating set $\hat{\mathcal{Z}}_{ech,M}$. The rule for assigning the relative grading of the generators is given by Hutchings in [Hu1] and [Hu2]. The rule is described briefly in Section III.9a.

*Part 3*: The appendix in [KLTII] and Sections III.1d and III.9b explain how the differential that defines the embedded contact homology is computed using J-holomorphic submanifolds in $\mathbb{R} \times Y$. Keep in mind that a J-holomorphic submanifold is properly embedded with J-invariant tangent space and such that the integral of *w* over the submanifold is finite. The particular J-holomorphic submanifolds that are used to define the differential comprise a topological space that is indexed by an ordered pair from $\hat{\mathcal{Z}}_{ech,M}$. Let $(\hat{\Theta}', \hat{\Theta})$ denote such a pair. The corresponding component of this topological space is denoted by $\mathcal{M}_1(\hat{\Theta}', \hat{\Theta})$. This space is a finite disjoint union of connected components, each being homeomorphic to $\mathbb{R}$. In fact, each component has a free $\mathbb{R}$ action that is induced by the constant translations along the $\mathbb{R}$ factor of $\mathbb{R} \times Y$.

Any given submanifold from $\mathcal{M}_1(\hat{\Theta}', \hat{\Theta})$ is characterized in part by the behavior of its $|s| >> 1$ part. To elaborate, suppose that C is a given submanifold from this space. There exists $s_* >> 1$ such that the $|s| \geq s_*$ portion is a disjoint union of embedded cylinders where d*s* is nowhere zero. Each such cylinder is said to be an *end* of the given submanifold. These ends have the following properties:

- *The* $s \geq s_*$ *ends are in 1-1 correspondence with the integral curves from* $\Theta$. *This correspondence is such that the set of constant s slices of any given end converge isotopically in* Y *as* $s \to \infty$ *to its partner in* $\Theta$.
- *The* $s \leq -s_*$ *ends are in 1-1 correspondence with the integral curves from* $\Theta'$. *This correspondence is such that the set of constant s slices of any given end converge isotopically in* Y *as* $s \to -\infty$ *to its partner in* $\Theta'$.

(1.8)

Section III.9b associates a sign, either 1 or -1, to each component of $\mathcal{M}_1(\hat{\Theta}', \hat{\Theta})$. This is done in accordance with the rules laid out by Hutchings in [Hu1] and [Hu2]. These signs determine the endomorphism of $\mathbb{Z}(\hat{\mathcal{Z}}_{ech,M})$ that supplies the embedded contact homology differential as follows: The relevant endomorphism of this $\mathbb{Z}$-module is given by its actions on the set of generators by the rule

$$\hat{\Theta} \to \sum_{\hat{\Theta}' \in \hat{\mathcal{Z}}_{ech,M}} N_{\hat{\Theta}',\hat{\Theta}} \hat{\Theta}' \,, \tag{1.9}$$

where any given $\hat{\Theta}' \in \hat{\mathcal{Z}}_{ech,M}$ version of $N_{\hat{\Theta}',\hat{\Theta}}$ is the sum of the +1's and -1's that are associated to the components of $\mathcal{M}_1(\hat{\Theta}', \hat{\Theta})$.

The differential on $\mathbb{Z}(\hat{\mathcal{Z}}_{ech,M})$ that defines the embedded contact homology decreases the relative grading by 1.

*Part 4*: The appendix in [KLTII] and Sections III.1d and III.9c describe a certain action of $\mathbb{Z}(\mathbb{U}) \otimes (\wedge^*(H_1(Y;\mathbb{Z})/\text{torsion}))$ on the embedded contact homology $\mathbb{Z}$-module. The endomorphism that generates the $\mathbb{Z}(\mathbb{U})$ factor is called the $\mathbb{U}$-map. The latter decreases the relative degree by 2 and it commutes with the generators of the $\wedge^*(H_1(Y;\mathbb{Z})/\text{torsion})$ factor. The generators of the latter decrease relative degree by 1. The $\mathbb{U}$-map generator and those of $\wedge^*(H_1(Y;\mathbb{Z})/\text{torsion})$ are given by corresponding endomorphisms of $\mathbb{Z}(\hat{\mathcal{Z}}_{ech,M})$. Each such endomorphism is defined by a version of (1.9) with the set $\{N_{\hat{\Theta}',\hat{\Theta}}\}_{\hat{\Theta}',\hat{\Theta} \in \hat{\mathcal{Z}}_{ech,M}}$ determined by certain sets of J-holomorphic submanifolds according to rules laid out by Hutchings in Section 12 of [HS]. See also Section 2.5 of [HT1] for a discussion of the $\mathbb{U}$-map generator.

An endomorphism of $\mathbb{Z}(\hat{\mathcal{Z}}_{ech,M})$ that defines the $\mathbb{U}$-map is defined in Section III.1d and Section III.9c with the help of chosen point in the handle $\mathcal{H}_0$. With the point chosen, the the set of coefficients in the corresponding version of (1.9) is denoted by $\{N^{\mathbb{U}}_{\hat{\Theta}',\hat{\Theta}}\}_{\hat{\Theta}',\hat{\Theta} \in \hat{\mathcal{Z}}_{ech,M}}$. These are such that any given $N^{\mathbb{U}}_{\hat{\Theta}',\hat{\Theta}}$ is non-zero only $\hat{\Theta}'$ and $\hat{\Theta}$ sit over the same element in $\mathcal{Z}_{ech,M}$. Moreover, there is precisely one non-zero $N^{\mathbb{U}}_{\hat{\Theta}',\hat{\Theta}}$ such that both $\hat{\Theta}'$ and $\hat{\Theta}$ sit over any given element in $\mathcal{Z}_{ech,M}$, and the corresponding $N^{\mathbb{U}}_{\hat{\Theta}',\hat{\Theta}}$ is equal to 1. A single J-holomorphic submanifold is used to compute this nonzero $N^{\mathbb{U}}_{\hat{\Theta}',\hat{\Theta}}$: If $\Theta \in \mathcal{Z}_{ech,M}$ is the given element, then the corresponding submanifold is the union of the cylinders from the set $\{\mathbb{R} \times \gamma\}_{\gamma \in \Theta}$ and $\{0\} \times S \subset \mathbb{R} \times Y$ with S being the u = constant sphere in $\mathcal{H}_0$ that contains the chosen point.

The endomorphisms of $\mathbb{Z}(\hat{\mathcal{Z}}_{ech,M})$ that define a set of generators for the $\wedge^*(H_1(Y;\mathbb{Z})/\text{torsion})$ action on the embedded contact homology are defined with the help of a chosen set of 1-cycles that supply a basis for $H_1(Y;\mathbb{Z})/\text{torsion}$. Section III.1d took this set to have the form that is described momentarily. To set the background, introduce $b_1(M)$ to denote the first Betti number of M. Section II.2a describes $1 + b_1(M)$ closed integral curves of $v$ in $M_\delta \cup \mathcal{H}_0$ that have intersection number 1 with each cross-sectional sphere in $\mathcal{H}_0$. One of these curves is the aforementioned $\gamma^{(z_0)}$. The curves in this set are labeled by the intersection point with the surface $f^{-1}(\frac{3}{2})$. This set of points is denoted by ¥ and the curve that contains a given z ∈ ¥ is denoted by $\gamma^{(z)}$. Pairing with the Poincaré duals of the homology classes of the cycles that comprise the set $\{[\gamma^{(z)}] - [\gamma^{(z_0)}]\}_{z \in ¥ - z_0}$ generates the dual in $\mathrm{Hom}(H_2(Y;\mathbb{Z});\mathbb{Z})$ of the $H_2(M;\mathbb{Z})$ summand in (1.7).

The basis used in Section III.1d contains the cycle $[\gamma^{(z_0)}]$, it contains a set of cycles that are labeled $\{\hat{\imath}^{(z)}\}_{z \in ¥ - z_0}$, and it is rounded out by a set of G cycles that are labled $\{\hat{\imath}_{\mathfrak{p}}\}_{\mathfrak{p} \in \Lambda}$. A given $z \in ¥ - z_0$ version of $\hat{\imath}^{(z)}$ lies entirely in the $M_{7\delta_*}$ part of Y. It is homologous to $[\gamma^{(z)}] - [\gamma^{(z_0)}]$ and it is obtained from the latter by first truncating the $\mathcal{H}_0$ portions of the curves $\gamma^{(z)}$ and $\gamma^{(z_0)}$ and then reconnecting the respective endpoints by arcs on the boundary of the radius $7\delta_*$ coordinate balls about the index 0 and index 3 critical points of $f$. A given $\mathfrak{p} \in \Lambda$ version of $\hat{\imath}_{\mathfrak{p}}$ is disjoint from the $f \in [1, 2]$ part of $M_{7\delta_*}$, and it intersects the rest of $M_{7\delta_*}$ and $\mathcal{H}_0$ as a smooth curve that is transverse to the level sets of $f$ in $M_\delta$ and the constant u spheres in $\mathcal{H}_0$; the orientation is such that it has intersection number 1 with the u = 0 sphere in $\mathcal{H}_0$. Meanwhile, $\hat{\imath}_{\mathfrak{p}}$ intersects $\cup_{\mathfrak{p}' \in \Lambda} \mathcal{H}_{\mathfrak{p}'}$ as the $\theta = 0$ arc in $\mathcal{H}_{\mathfrak{p}}$; its orientation gives it intersection number -1 with each u = 0 sphere in $\mathcal{H}_{\mathfrak{p}}$.

Suppose that $\hat{\imath}$ is a given cycle from the chosen basis of cycles. As noted above, the corresponding endomorphism of $\mathbb{Z}(\hat{\mathcal{Z}}_{ech,M})$ has the form given in (1.9). Denote the set of coefficients by $\{N^{\hat{\imath}}_{\hat{\Theta}',\hat{\Theta}}\}_{\hat{\Theta}',\hat{\Theta} \in \hat{\mathcal{Z}}_{ech,M}}$. Any given $N^{\hat{\imath}}_{\hat{\Theta}',\hat{\Theta}}$ is defined using the J-holomorphic submanifolds from $\mathcal{M}_1(\hat{\Theta}', \hat{\Theta})$. In particular, $N^{\hat{\imath}}_{\hat{\Theta}',\hat{\Theta}}$ is the value of a sum that is indexed by the components of $\mathcal{M}_1(\hat{\Theta}', \hat{\Theta})$ whereby the component of a given submanifold C contributes either +1 or -1 times the algebraic intersection number between C and $\mathbb{R} \times \hat{\imath}$. This intersection number is well defined because $\hat{\imath}$ is disjoint from the integral curves of $v$ that come from elements in $\mathcal{Z}_{ech,M}$. The +1 or -1 used here is the contribution of C's component to the version of $N_{\hat{\Theta}',\hat{\Theta}}$ that defines the embedded contact homology differential.

*Part 5*: This last part of the subsection introduces a certain filtration of the embedded contact homology chain complex that is preserved by the differential. The filtration is depicted by (I.2.3). What follows reviews what is involved. To start, invoke Proposition II.2.8 or Theorem I.2.1 to write the set $\mathcal{Z}_{ech,M}$ as $\mathcal{Z}_{HF} \times (\times_{\mathfrak{p} \in \Lambda}(\mathbb{Z} \times \mathrm{O}))$. By way of a reminder, $\mathcal{Z}_{HF}$ denotes a a certain set that is defined using data coming strictly from M and O is the 4-element set {0, 1, -1, {1,-1}}. The $\mathcal{Z}_{HF}$-label of any given element $\Theta \in \mathcal{Z}_{ech,M}$ characterizes the intersection of $\cup_{\gamma \in \Theta} \gamma$ with $M_\delta$. Meanwhile, each $\mathfrak{p} \in \Lambda$ factor of $\mathbb{Z} \times \mathrm{O}$ characterizes the intersection of $\cup_{\gamma \in \Theta} \gamma$ with $\mathcal{H}_{\mathfrak{p}}$. The integer component of this label characterizes the segment of $(\cup_{\gamma\in\Theta}\gamma)\cap\mathcal{H}_{\mathfrak{p}}$ that crosses $\mathcal{H}_{\mathfrak{p}}$ from its u = -R - ln(7$\delta_*$) end to its u = R + ln(7$\delta_*$) end. The label from the set {0, 1, -1, {1,-1}} signifies which, if any, integral curves from the set $\{\hat{\gamma}_{\mathfrak{p}}^{+}, \hat{\gamma}_{\mathfrak{p}}^{-}\}$ appear in $\Theta$. The +1 signifies $\hat{\gamma}_{\mathfrak{p}}^{+}$ and the -1 signifies $\hat{\gamma}_{\mathfrak{p}}^{-}$. Use the identification of $\mathcal{Z}_{\mathrm{ech,M}}$ with $\mathcal{Z}_{\mathrm{HF}}\times(\times_{\mathfrak{p}\in\Lambda}(\mathbb{Z}\times\mathrm{O}))$ to write a given element $\Theta$ as $(\hat{\upsilon}, (\mathfrak{k}_{\mathfrak{p}}, \mathrm{O}_{\mathfrak{p}})_{\mathfrak{p}\in\Lambda})$. For each $\mathfrak{p}\in\Lambda$, use $|\mathrm{O}_{\mathfrak{p}}|\in\{0, 1, 2\}$ to denote the sum of the absolute values of the elements in $\mathrm{O}_{\mathfrak{p}}$.

Associate to each non-negative integer L the subset $\mathcal{Z}_{\mathrm{ech,M}}{}^{L}\subset\mathcal{Z}_{\mathrm{ech,M}}$ whose elements are such that $\sum_{\mathfrak{p}\in\Lambda}(|\mathfrak{k}_{\mathfrak{p}}| + 2|\mathrm{O}_{\mathfrak{p}}|) < L$. These sets are such that $\mathcal{Z}_{\mathrm{ech,M}}{}^{L}\subset\mathcal{Z}_{\mathrm{ech,M}}{}^{L'}$ when L´ > L, and their union is the whole of $\mathcal{Z}_{\mathrm{ech,M}}$. Let $\hat{\mathcal{Z}}_{\mathrm{ech,M}}{}^{L}$ denote the inverse image of $\mathcal{Z}_{\mathrm{ech,M}}{}^{L}$ in $\hat{\mathcal{Z}}_{\mathrm{ech,M}}$. It follows from Theorem I.2.3 or Theorem III.1.1 that the embedded contact homology differential maps the submodule $\mathbb{Z}(\hat{\mathcal{Z}}_{\mathrm{ech,M}}{}^{L})\subset\mathbb{Z}(\hat{\mathcal{Z}}_{\mathrm{ech,M}})$ to itself and so the latter defines a subcomplex. The embedded contact homology is the direct limit of the homology for the filtered sequence of chain subcomplexes

$$\cdots\subset\mathbb{Z}(\hat{\mathcal{Z}}_{\mathrm{ech,M}}{}^{L})\subset\mathbb{Z}(\hat{\mathcal{Z}}_{\mathrm{ech,M}}{}^{L+1})\subset\cdots \tag{1.10}$$

of the chain complex $\mathbb{Z}(\hat{\mathcal{Z}}_{\mathrm{ech,M}}{}^{L})$.

**c) The Seiberg-Witten Floer homology on Y**

This subsection describes various versions of the Seiberg-Witten Floer homology on the manifold Y. The presentation that follows takes for granted the basic constructions and theorems about Seiberg-Witten Floer homology and focuses almost exclusively on those parts of the story that are specific to the geometry at hand. The book by Kronheimer and Mrowka [KM] is the recommended text book for those who are not familiar with the foundational background. There are nine parts to what follows.

Parts 1-5 introduce various geometric notions that are used in Part 6 to define the chain complex and differential whose homology groups constitute the desired versions Seiberg-Witten Floer homology. These groups are introduced in Part 8. The intervening Part 7 describes certain canonical endomorphisms of the chain complex that are used to generate an action of $\mathbb{Z}[\mathbb{U}]\otimes(\wedge^*(H^1(Y;\mathbb{Z})/\mathrm{torsion}))$ on the homology. Part 9 explains why Theorems I.3.1 and I.3.2 are direct consequences of what is said in Parts 6-8.

*Part 1*: Part 5 in Section 1a defined a Riemannian metric on Y; this being a metric with $*w = \hat{a}$ and $|\hat{a}| = 1$. Use this metric to define the bundle of oriented, orthonormal frames for Y. The given Spin$^{\mathbb{C}}$-structure on Y determines a lift of this bundle to a principal U(1) bundle. The defining action of U(2) on $\mathbb{C}^2$ supplies an

associated Hermitian $\mathbb{C}^2$-bundle. The latter is denoted by $\mathbb{S}$. Use $\det(\mathbb{S})$ to denote the complex line bundle $\wedge^2\mathbb{S}$.

There is a canonical homomorphism from $T^*Y$ into $\mathrm{End}(\mathbb{S})$, this being Clifford multiplication. The homomorphism is denoted by cl; it is characterized as follows: Let a and b denote a pair of covectors in a given fiber of $T^*Y$. Then

$$\mathrm{cl}(a)^\dagger = -\mathrm{cl}(a) \quad \textit{and} \quad \mathrm{cl}(a)\mathrm{cl}(b) = -\langle a, b\rangle - \mathrm{cl}(*(a \wedge b)) \tag{1.11}$$

where $\langle\ ,\ \rangle$ here denotes the metric inner product and $*$ denotes the metric's Hodge dual. This Clifford multiplication map induces two other useful endomorphisms. The first, denoted by $\hat{c}$, maps $\mathbb{S}\otimes T^*Y$ to $\mathbb{S}$. It is defined so as to send an reducible element $\psi \otimes a$ to $\mathrm{cl}(a)\psi$. The second is the $\mathbb{R}$-linear homomorphism from $\mathbb{S} \otimes \mathbb{S}$ to $T^*M \otimes_{\mathbb{R}} \mathbb{C}$ that is written as $\eta \otimes \psi \to \eta^\dagger\tau\psi$ and defined by the rule whereby $\langle a, \eta^\dagger\tau\psi\rangle = \eta^\dagger\mathrm{cl}(a)\psi$.

Clifford multiplication by $\hat{a}$ splits $\mathbb{S}$ as a direct sum of complex line bundles, this written as

$$\mathbb{S} = E \oplus EK^{-1}\ . \tag{1.12}$$

Here, $K^{-1}$ is isomorphic as a real 2-plane bundle to the oriented bundle $\ker(\hat{a})$, this being the kernel of $\hat{a}$ with the orientation defined by $w$. The convention in what follows takes the left most summand as the +i eigenspace of $\mathrm{cl}(\hat{a})$. The corresponding component decomposition of a given section of $\mathbb{S}$ is written as $(\alpha, \beta)$ with $\alpha$ being a section of E and $\beta$ being a section of $EK^{-1}$.

A Hermitian connection on $\det(\mathbb{S})$ with the Levi-Civita connnection on TY jointly define a Hermitian connection on $\mathbb{S}$ and thus a covariant deriviate, this being a map from $C^\infty(Y;\mathbb{S})$ to $C^\infty(Y;\mathbb{S}\otimes T^*Y)$. Meanwhile, the Clifford multiplication endomorphism defines the endomorphism, $\hat{\mathrm{cl}}: \mathbb{S} \otimes T^*Y \to \mathbb{S}$. The composition of the covariant derivative and then $\hat{\mathrm{cl}}$ defines a first order, elliptic operator from $C^\infty(Y;\mathbb{S})$ to itself, this being the *Dirac* operator.

The Dirac operator is used momentarily to define a canonical connection on $K^{-1}$. To do this, introduce $I_{\mathbb{C}}$ to denote a topologically trivial complex line bundle over Y, and let $\mathbb{S}_I$ denote the $\mathrm{Spin}_{\mathbb{C}}$ structure given by the $E = I_{\mathbb{C}}$ version of (1.12). Fix a unit norm section, $1_{\mathbb{C}}$ of $I_{\mathbb{C}}$ and view the pair $(1_{\mathbb{C}}, 0)$ as a section of $\mathbb{S}$ using the splitting in (1.12). Since $\det(\mathbb{S}_I) = K^{-1}$, a Hermitian connection on $K^{-1}$ defines an associated Dirac operator. The canonical connection on $K^{-1}$ is characterized by the fact that $(1_{\mathbb{C}}, 0)$ is annihilated by its associated Dirac operator. This connection on $K^{-1}$ is denoted by $A_K$.

Let $\mathbb{S}$ be the $\mathbb{C}^2$ bundle that comes from the given Spin$^{\mathbb{C}}$-structure. With (1.12) understood, $\det(\mathbb{S}) = E^2K^{-1}$ and thus any given Hermitian connection on $\det(\mathbb{S})$ can be written as $A_K + 2A$ where A is a connection on E. With A a given connection on E, the symbol $D_A$ is used in what follows to denote the Dirac operator on sections of $\mathbb{S}$ that is defined using $A_K + 2A$ for the connection on $\det(\mathbb{S})$. Use Conn(E) to denote the Frechêt space of smooth, Hermitian connections on E.

The Frechet Lie group $C^\infty(Y; S^1)$ acts smoothly on $\text{Conn}(E) \times C^\infty(Y; \mathbb{S})$ by the rule whereby any given map u sends any given pair $(A, \psi)$ to $(A - u^{-1}du, u\psi)$.

*Part 2*: The $\mathbb{Z}$-module that serves as the chain complex for the Seiberg-Witten Floer homology is constructed using solutions to certain versions of the Seiberg-Witten equations. These are equations for pairs $(A, \psi)$ with $\psi$ a section of $\mathbb{S}$ and $A \in \text{Conn}(E)$. The simplest of the relevant versions constitutes a family of equations whose members are labeled by a real number greater than $\pi$ and a smooth 1-form on Y. The version defined by a given $r \in [1, \infty)$ and 1-form $\mu$ asks that $(A, \psi)$ obey

- $B_A - r(\psi^\dagger\tau\psi - i\hat{a}) + \frac{1}{2} B_{A_K} - i{*}d\mu = 0.$
- $D_A\psi = 0.$

(1.13)

where the notation is such that $B_A = {*}F_A$ with $F_A$ being the curvature 2-form of the connection A. Likewise, $B_{A_K}$ denotes the Hodge dual of the curvature 2-form of $A_K$

Certain perturbed versions of (1.13) are needed to guarantee that the solutions to the associated equation and its instanton counter part are suitably generic. A given perturbed version of (1.13) is defined using a chosen element in a certain Banach space of $C^\infty(Y; S^1)$-invariant functions on $\text{Conn}(E) \times C^\infty(Y; \mathbb{S})$ (see Chapter 11 in [KM]). This Banach space is denoted by $\mathcal{P}$ and its norm is called the $\mathcal{P}$-norm. The Banach space is such that the differential of a given $\mathfrak{g} \in \mathcal{P}$ is a smooth map from $\text{Conn}(E) \times C^\infty(Y; \mathbb{S})$ to $C^\infty(Y; iT^*Y \oplus C^\infty(Y; \mathbb{S}))$. With $\mathfrak{g}$ chosen, write its differential at a given $(A, \psi)$ as $(\mathfrak{T}_{(A,\psi)}, \mathfrak{S}_{(A,\psi)})$, this being a pair consisting of an $i\mathbb{R}$-valued 1-form on Y and a section of $\mathbb{S}$. The 1-form $\mathfrak{T}_{(A,\psi)}$ is in the image of the operator ${*}d$. The $\mathfrak{g}$-perturbed version of the Seiberg-Witten equations are

- $B_A - r(\psi^\dagger\tau\psi - i\hat{a}) + \frac{1}{2} B_{A_K} - \mathfrak{T}_{(A,\psi)} = 0.$
- $D_A\psi - \mathfrak{S}_{(A,\psi)} = 0.$

(1.14)

The simplest but non-trivial perturbation has $\mathfrak{T} = i{*}d\mu$ and $\mathfrak{S} = 0$ with $\mu$ a smooth 1-form on Y taken from a certain Banach space of such forms; this perturbation gives the equation in (1.13). The Banach space is denoted by $\Omega$. The norm on this space is also called the $\mathcal{P}$-norm. The latter is such that the inclusion $\Omega \to C^\infty(M; iT^*Y)$ defines a bounded, compact mapping. The convention in this paper is to use 1-forms $\mu$ from $\Omega$ with $\mathcal{P}$-norm less than 1. All of the assertions hold (with the same proofs) given any a priori upper bound on the $\mathcal{P}$-norm.

When $\mu \in \Omega$, then the corresponding version of $\mathfrak{g}$ is denoted by $\mathfrak{e}_\mu$; it is the function that assigns to any given $(A, \psi)$ the integral over Y of the 3-form $-iF_A \wedge \mu$.

If $(A, \psi)$ is a solution to (1.14), then so is $(A - u^{-1}du, u\psi)$ for any $u \in C^\infty(Y; S^1)$. Use $\mathcal{Z}_{SW,r}$ in what follows to denote the $C^\infty(Y; S^1)$-quotient of the space of solutions to a given $\mathfrak{g} \in \mathcal{P}$ version of (1.14). Note in this regard that the group $C^\infty(Y; S^1)$ acts freely on the space of solutions to any given $r > \pi$ and $\mathfrak{g} \in \mathcal{P}$ version of (1.14). This is so because the group acts freely on $\mathrm{Conn}(E) \times (C^\infty(Y; S^1) – 0)$ and no $\psi = 0$ pair can solve (1.14) because both $\frac{i}{2\pi}(F_A + \frac{1}{2} B_{A_K})$ and $w$ represent the first Chern class of $\det(\mathbb{S})$. Note also that $\mathcal{Z}_{SW,r}$ is in all cases compact (see Chapter 29 in [KM].) By way of a warning, the notation does not indicate that $\mathcal{Z}_{SW,r}$ depends on the chosen perturbation $\mathfrak{g}$.

*Part 3*: The definition of the $\mathbb{Z}$-module for the Seiberg-Witten Floer homology requires the introduction of a subgroup of $C^\infty(Y; S^1)$ which is denoted by $\mathcal{G}_{M_\Lambda}$. A given map u sits in this subgroup when $-\frac{i}{2\pi}\int_{\gamma^{(z_0)}} u^{-1}du = 0$. Note in this regard that $-\frac{i}{2\pi}u^{-1}du$ has integer periods.

Use $\hat{\mathcal{Z}}_{SW,r}$ to denote the space of $\mathcal{G}_{M_\Lambda}$-orbits of solutions to a given $r \in [1, \infty)$ and $\mathfrak{g} \in \mathcal{P}$ version of (1.14). The space $\hat{\mathcal{Z}}_{SW,r}$ is a principal $\mathbb{Z} = H_2(\mathcal{H}_0; \mathbb{Z})$ bundle over $\mathcal{Z}_{SW,r}$. The action of $k \in \mathbb{Z}$ send the $\mathcal{G}_{M_\Lambda}$-equivalence class of $(A, \psi)$ to that of $(A - u^{-1}du, u\psi)$ with $u \in C^\infty(Y; S^1)$ any map with $-\frac{i}{2\pi}\int_{\gamma^{(z_0)}} u^{-1}du = k$.

The geometry of Y supplies a certain $\mathbb{Z}$-equivariant maps from $\hat{\mathcal{Z}}_{SW,r}$ to $\mathbb{R}$. The definition requires the choice of a fiducial connection on E, this denoted by $A_E$. This choice is constrained by the requirement that $A_E$ be flat on $\mathcal{H}_0$ and have holonomy 1 around the curve $\gamma^{(z_0)}$. The definition of the $\mathbb{Z}$-equivariant map from $\hat{\mathcal{Z}}_{SW,r}$ to $\mathbb{R}$ also requires the choice of a smooth function $\wp: [0,\infty) \to [0,\infty)$ which is non-decreasing, obeys $\wp(x) = 0$ for $x < \frac{7}{16}$ and $\wp(x) = 1$ for $x \geq \frac{9}{16}$. It proves convenient to choose $\wp$ so that its derivative, $\wp'$, is bounded by $2^{10}(1 - \wp)^{3/4}$. A function with these properties can be readily constructed from the function on $\mathbb{R}$ that is set equal to 0 for $t < 0$ and set equal to $e^{-1/t}$ for $t > 0$. Fix such a function $\wp$.

Given $\mathfrak{c} = (A, \psi)$ from $\mathrm{Conn}(E) \times C^\infty(Y; \mathbb{S})$, write $\psi = (\alpha, \beta)$ and define

$$\hat{A} = A - \tfrac{1}{2}\wp(|\alpha|^2)|\alpha|^{-2}(\bar{\alpha}\nabla_A\alpha - \alpha\nabla_A\bar{\alpha}), \tag{1.15}$$

this being a connection on E. The corresponding equivariant map to $\mathbb{R}$ is the restriction from $(\mathrm{Conn}(E) \times C^\infty(Y; \mathbb{S}))/\mathcal{G}_{M_\Lambda}$ of the $\mathcal{G}_{M_\Lambda}$-invariant map from $\mathrm{Conn}(E) \times C^\infty(Y;\mathbb{S})$ to $\mathbb{R}$ that sends any given pair $\mathfrak{c} = (A, \psi)$ to

$$x(\mathfrak{c}) = \tfrac{i}{2\pi}\int_{\gamma^{(z_0)}} (\hat{A} - A_E) \ . \tag{1.16}$$

A given element $\mathfrak{c} \in \mathrm{Conn}(E) \times C^\infty(Y;\mathbb{S})$ is deemed to be *holonomy nondegenerate* when $x(\mathfrak{c}) - \frac{1}{2}$ is not an integer. The locus where $x(\mathfrak{c}) - \frac{1}{2} \in \mathbb{Z}$ is a codimension 1 submanifold in $\mathrm{Conn}(E) \times C^\infty(Y;\mathbb{S})$. The element $\mathfrak{c}$ is holonomy nondegenerate if and only if all pairs in its $C^\infty(Y;S^1)$ orbit are holonomy nondegenerate.

*Part 4*: A certain symmetric, first-order elliptic operator on $C^\infty(Y; iT^*Y \oplus \mathbb{S} \oplus i\mathbb{R})$ is associated to each pair in $\mathrm{Conn}(E) \times C^\infty(Y;\mathbb{S})$. Fix $\mathfrak{c} = (A,\psi) \in \mathrm{Conn}(E) \times C^\infty(Y;\mathbb{S})$. The corresponding operator in the case when $\mathfrak{g} = \mathfrak{e}_\mu$ sends a section $\mathfrak{b} = (b, \eta, \phi)$ to the section with respective $iT^*Y$, $\mathbb{S}$ and $i\mathbb{R}$ components

- $*db - d\phi - 2^{-1/2}\, r^{1/2}\, (\psi^\dagger\tau\eta + \eta^\dagger\tau\psi)$
- $D_A\eta + 2^{1/2}\, r^{1/2}(cl(b)\psi + \phi\psi)$
- $*d{*}b - 2^{-1/2}\, r^{1/2}(\eta^\dagger\psi - \psi^\dagger\eta)$

(1.17)

In the generic $\mathfrak{g}$ case, the operator is obtained from (1.17) by adding $(\frac{2}{r})^{1/2}(\frac{d}{d\tau}\,\mathfrak{T}_{\mathfrak{c}+\tau\hat{\mathfrak{b}}})_{\tau=0}$ to the top bullet and $(\frac{d}{d\tau}\,\mathfrak{S}_{\mathfrak{c}+\tau\hat{\mathfrak{b}}})_{\tau=0}$ to the middle bullet with $\hat{\mathfrak{b}} = ((\frac{r}{2})^{1/2}b, \eta)$. The operator in all cases is denoted by $\mathfrak{L}_{\mathfrak{c},r}$. The pair $\mathfrak{c}$ is said to be *non-degenerate* when $\mathfrak{L}_{\mathfrak{c},r}$ has trivial kernel. A given pair $\mathfrak{c}$ is non-degenerate if and only all pairs in its $C^\infty(Y;S^1)$ orbit are likewise non-degenerate.

The following statements are analogs of what is asserted in Lemma 3.6 and Proposition 3.11 [T1] for the case when $\hat{a}$ is a contact 1-form and $w = d\hat{a}$. The proofs differ only slightly from those given in [T1].

- *Fix* $r \geq 1$ *and there is an open, dense set in* $\Omega$ *which is characterized as follows: If* $\mu$ *is from this set, then the corresponding* $\mathfrak{g} = \mathfrak{e}_\mu$ *version of* $\mathcal{Z}_{SW,r}$ *is a finite set of orbits of pairs in* $\mathrm{Conn}(E) \times C^\infty(Y;\mathbb{S})$ *and each such pair is non-degenerate and holonomy non-degenerate.*
- *There exists a residual set in* $\Omega$ *that is characterized as follows: Fix* $\mu$ *from this set. There is a countable, non-accumulating set in* $(\pi, \infty)$ *such that if* $r$ *is from its complement, then the corresponding* $(r, \mathfrak{g} = \mathfrak{e}_\mu)$ *version of* $\mathcal{Z}_{SW,r}$ *is a finite set of* $C^\infty(Y;S^1)$ *orbits in* $\mathrm{Conn}(E) \times C^\infty(Y;\mathbb{S})$ *and each such orbit contains only non-degenerate and holonomy non-degenerate pairs*.

(1.18)

Suppose now that $(r, \mathfrak{g})$ is such that the corresponding $\mathcal{Z}_{SW,r}$ is a finite set of orbits and that each orbit contains only non-degenerate and holonomy non-degenerate pairs. The principal bundle in this case has a canonical $\mathbb{Z}$-equivariant isomorphism

$$\hat{\mathcal{Z}}_{SW,r} = \mathcal{Z}_{SW,r} \times \mathbb{Z}$$

(1.19)

that is characterized as follows: The section $\mathcal{Z}_{SW,r} \times \{0\}$ of the product bundle corresponds to the set of the $\mathcal{G}_{M_A}$ orbits of solutions to (1.14) with $\mathrm{x}(\cdot) \in (-\frac{1}{2}, \frac{1}{2})$. Granted this identification, use $\hat{\mathcal{Z}}^{\geq}_{SW,r} \subset \hat{\mathcal{Z}}_{SW,r}$ to denote the subset that is identified via (1.19) to $\mathcal{Z}_{SW,r} \times \{0, 1, 2, \ldots\} \subset \mathcal{Z}_{SW,r} \times \mathbb{Z}$.

*Part 5*: Certain versions of the Seiberg-Witten equations on $\mathbb{R} \times Y$ play a central role in the story. As in the case of (1.14), the equations require the choice of $r \geq 1$ and a perturbation $\mathfrak{g}$ from $\mathcal{P}$. The corresponding equations is viewed here as system of differential equations for a map from $\mathbb{R}$ to $\mathrm{Conn}(E) \times C^\infty(Y;\mathbb{S})$. The equations demand that $s \to \mathfrak{d}(s) = (A, \psi)|_s$, obey

- $\frac{\partial}{\partial s}A + B_A - r(\psi^\dagger\tau\psi - i\hat{a}) + \frac{1}{2} B_{A_K} - \mathfrak{T}_{(A,\psi)} = 0$ .
- $\frac{\partial}{\partial s}\psi + D_A\psi - \mathfrak{S}_{(A,\psi)} = 0.$

(1.20)

A solution to (1.19) said to be an *instanton* if it has respective $s \to \pm\infty$ limits that obey (1.14). Any constant $\mathbb{R}$ translate of an instanton solution to (1.20) is also an instanton solution.

An instanton is said to be nondegenerate if a certain operator is Fredholm and has trivial cokernel. The relevant operator maps an $L^2_1$ Hilbert space completion of the space of compactly supported elements in $C^\infty(\mathbb{R} \times Y; iT^*Y \oplus \mathbb{S} \oplus i\mathbb{R})$ to an $L^2$ completion. This operator in the case $\mathfrak{g} = \mathfrak{e}_\mu$ sends a section $(b, \eta, \phi)$ to the section whose respective $iT^*Y$, $\mathbb{S}$ and $i\mathbb{R}$ components are

- $\frac{\partial}{\partial s}b + {*db} - d\phi - 2^{-1/2} r^{1/2} (\psi^\dagger\tau\eta + \eta^\dagger\tau\psi)$ ,
- $\frac{\partial}{\partial s}\eta + D_A\eta + 2^{1/2} r^{1/2}(cl(b)\psi + \phi\psi)$ ,
- $\frac{\partial}{\partial s}\phi + {*d{*b}} - 2^{-1/2} r^{1/2}(\eta^\dagger\psi - \psi^\dagger\eta)$ .

(1.21)

The operator in the generic $\mathfrak{g}$ case is obtained from (1.21) by adding $(\frac{2}{r})^{1/2}(\frac{d}{d\tau}\mathfrak{T}_{\mathfrak{d}+\tau\hat{\mathfrak{b}}})_{\tau=0}$ to the top bullet and $(\frac{d}{d\tau}\mathfrak{S}_{\mathfrak{d}+\tau\hat{\mathfrak{b}}})_{\tau=0}$ to the middle bullet. The operator in any case is denoted by $\mathfrak{D}_{\mathfrak{d}}$.

Let $\mathfrak{d}$: $\mathbb{R} \to \text{Conn}(E) \times C^\infty(Y;\mathbb{S})$ denote an instanton solution to some (r, $\mathfrak{g}$) version of (1.20). The corresponding version of (1.21) is a Fredholm operator if and only if the $s \to \pm\infty$ limits of $\mathfrak{d}$ are nondegenerate pairs in $\text{Conn}(E) \times C^\infty(Y; \mathbb{S})$. If this is the case, then $\mathfrak{D}_{\mathfrak{d}}$ has a corresponding Fredholm index, this denoted by $\iota_{\mathfrak{d}}$.

The assertion made by the upcoming (1.22) is used when considering perturbations. This fact has an almost verbatim analog stated in Part 5 from Section 3b of [T2] for the case when $\hat{a}$ is a contact 1-form on a given 3-manifold. The argument that prove (1.22) differs only cosmetically from those given for their [T2] analog.

To set the notation for (1.22), suppose that $\mathfrak{p} \in \mathcal{P}$ and that $\mathfrak{c}$ is a given pair from $\text{Conn}(E) \times C^\infty(Y;\mathbb{S})$. Then $\mathfrak{p}$ is said to *vanishes to second order* at $\mathfrak{c}$ if $\mathfrak{p}(\mathfrak{c}) = 0$, and if both the first and second derivatives at $\tau = 0$ of the function $\tau \to \mathfrak{p}(\mathfrak{c}+\tau\mathfrak{b})$ are zero for all pairs $\mathfrak{b} \in C^\infty(Y; iT^*Y \oplus \mathbb{S})$. With this term understood, introduce $\mathcal{P}_\mu \subset \mathcal{P}$ to denote the subset whose members have the following property: If $\mathfrak{p} \in \mathcal{P}_\mu$, then $\mathfrak{p} = 0$ to second order on any solution to (1.14) and $\mathfrak{p} = 0$ to second order at all points on any path $s \to \mathfrak{d}(s)$ with $\mathfrak{d}$ an $\iota_{\mathfrak{d}} \le 2$, non-degenerate instanton solution to the $(r, \mathfrak{g} = \mathfrak{e}_\mu)$ version of (1.20).

*Fix* $r \ge 1$ *and* $\mu$ *such that all solutions to the corresponding* $\mathfrak{g} = \mathfrak{e}_\mu$ *version of (1.14) are non-degenerate. There is a residual set in* $\mathcal{P}_\mu$ *characterized as follows: If* $\mathfrak{p}$ *is a member, then all instanton solutions to the* $(r, \mathfrak{g} = \mathfrak{e}_\mu + \mathfrak{p})$ *version of (1.20) are non-degenerate.*

(1.22)

Suppose that (r, $\mathfrak{g}$) is such that all solutions to (1.14) are non-degenerate and such that all $\iota_{\mathfrak{d}} \le 2$ instanton solutions to (1.20) are non-degenerate. If this is the case, then the set of $\iota_{\mathfrak{d}} = 0$ instanton solutions is the set of constant maps from $\mathbb{R}$ to the set of solutions to (1.14). To say something about the set of $\iota_{\mathfrak{d}} = 1$ instanton solutions to (1.20), suppose that $\mathfrak{c}_-$ and $\mathfrak{c}_+$ are given solutions to (1.14). Introduce $\mathcal{M}_1(\mathfrak{c}_-, \mathfrak{c}_+)$ to denote the set of instanton solutions to (1.20) with $s \to -\infty$ limit equal to $\mathfrak{c}_-$ and $s \to \infty$ limit a pair on the $\mathcal{G}_{M_\Lambda}$ orbit of $\mathfrak{c}_+$. This set $\mathcal{M}_1(\mathfrak{c}_-, \mathfrak{c}_+)$ has the structure of a smooth 1-dimensional manifold with a finite set of components, each being a copy of $\mathbb{R}$. Moreover, the group of constant translations of $\mathbb{R}$ induces a smooth, free $\mathbb{R}$-action on each component.

*Part 6*: This part defines chain complexes whose corresponding homology groups are of central concern here. To this end, fix $r \ge 1$ and fix $\mu \in \Omega$ with $\mathcal{P}$-norm bounded by 1 and such that all solutions to the (r, $\mathfrak{e}_\mu$) version of (1.14) are non-degenerate and also holonomy non-degenerate. A $\mathbb{Z}$-module that serves for the chain complex of interest is the free module generated by the elements in the $(r, \mathfrak{e}_\mu)$ version of $\hat{\mathcal{Z}}_{SW,r}$. This module is denoted in what follows by $\mathbb{Z}(\hat{\mathcal{Z}}_{SW,r})$. This module is written as $\mathbb{Z}(\hat{\mathcal{Z}}_{SW,r})$. The action of the generator of $H_2(\mathcal{H}_0;\mathbb{Z})$ gives this module a $\mathbb{Z}[t, t^{-1}]$ structure.

The $\mathbb{Z}$-module $\mathbb{Z}(\hat{\mathcal{Z}}_{SW,r})$ has a relative $\mathbb{Z}/(p_M\mathbb{Z})$ grading which is defined as follows: Let $\mathfrak{c}_0$ and $\mathfrak{c}_1$ denote two solutions to the relevant version of (1.14). Introduce $gr_{SW}(\mathfrak{c}_0)$ - $gr_{SW}(\mathfrak{c}_1)$ to denote the difference between the grading degrees of their respective $\mathcal{G}_{M_\Lambda}$-orbits in $\hat{\mathcal{Z}}_{SW,Y,r}$. This number is (-1) times the spectral flow for the [0,1]-parametrized family of operators $\tau \to \mathfrak{L}_{\mathfrak{c}(\tau)}$ with $\mathfrak{c}(0) = \mathfrak{c}_0$ and $\mathfrak{c}(1) = \mathfrak{c}_1$. Those unfamiliar with the notion of spectral flow can read about it in Chapter 14.2 of [KM] or in [T3]. The generator of the $\mathbb{Z}[t, t^{-1}]$ action on $\mathbb{Z}(\hat{\mathcal{Z}}_{SW,r})$ acts as a degree -2 endomorphism.

The relevant differential is a certain square zero endomorphism of $\mathbb{Z}(\hat{\mathcal{Z}}_{SW,r})$. This endomorphism is defined by its action on the generators. To say more, introduce by way of notation $[\mathfrak{c}]$ to denote the $\mathcal{G}_{M_\Lambda}$-equivalence class of a given pair $\mathfrak{c} = (A,\psi)$ from $\mathrm{Conn}(E) \times C^\infty(Y;\mathbb{S})$. Any given endomorphism of $\mathbb{Z}(\hat{\mathcal{Z}}_{SW,r})$ is defined by a rule

$$[\mathfrak{c}] \to \sum_{[\mathfrak{c}']\in\hat{\mathcal{Z}}_{SW,r}} w_{[\mathfrak{c}],[\mathfrak{c}']}\, [\mathfrak{c}'] \tag{1.23}$$

where each $[\mathfrak{c}']$ version of $w_{[\mathfrak{c}'],[\mathfrak{c}]}$ is an integer and only finitely many are non-zero.

In the case of the differential, the specification of the coefficient set $\{w_{[\mathfrak{c}'],[\mathfrak{c}]}\}_{[\mathfrak{c}],[\mathfrak{c}']\in\hat{\mathcal{Z}}_{SW,r}}$ requires first the choice of an element $\mathfrak{p} \in \mathcal{P}$ with norm much less than 1 such that the conclusions of (1.22) hold. Any given $[\mathfrak{c}], [\mathfrak{c}'] \in \hat{\mathcal{Z}}_{SW,r}$ version of $w_{[\mathfrak{c}'],[\mathfrak{c}]}$ is a sum that is indexed by the components of the $(r, \mathfrak{g} = \mathfrak{e}_\mu + \mathfrak{p})$ version of $\mathcal{M}_1(\mathfrak{c}', \mathfrak{c})$ with each component contributing either +1 or -1 to the sum. The sign is obtained by comparing two orientations for the component, one given by the generator of the $\mathbb{R}$ action and the other using Quillen's notion of a determinant line bundle for a family of Fredholm operators. This is done according to the rules given in Chapters 20-22 of [KM]; see also Section 3 of [T4].

**Proposition 1.1**: *There exists* $\kappa \geq 1$ *with the following significance: Fix* $r \geq \kappa$ *and an element* $\mu \in \Omega$ *with* $\mathcal{P}$*-norm less than 1 such that all solutions to the* $(r,\mu)$ *version of (1.13) are non-degenerate. Suppose that* $\mathfrak{p} \in \mathcal{P}_\mu$ *has small* $\mathcal{P}$*-norm and is described by (1.22).*

- *The rules given in [KM] for specifying the various* $[\mathfrak{c}], [\mathfrak{c}'] \in \hat{\mathcal{Z}}_{SW,r}$ *versions of* $w_{[\mathfrak{c},\mathfrak{c}']}$ *define a square zero endomorphism of* $\mathbb{Z}(\hat{\mathcal{Z}}_{SW,r})$.
- *Each solution to (1.13) is holonomy non-degenerate and so* $\mathbb{Z}(\hat{\mathcal{Z}}^{\geq}_{SW,r})$ *is well defined.*

- *The endomorphism given by the first bullet maps the submodule* $\mathbb{Z}(\hat{\mathcal{Z}}^{\geq}_{SW,r})$ *to itself.*

This proposition is proved in Section 7a. Take it on faith for now and use $\partial_{SW,Y}$ to denote the resulting endomorphism of $\mathbb{Z}(\hat{\mathcal{Z}}_{SW,r})$. It follows from the definition that $\partial_{SW,Y}$ decreases the $\mathbb{Z}/(p_M\mathbb{Z})$ grading by 1.

*Part 7*: This part of the subsection describes endomorphisms of $\mathbb{Z}(\hat{\mathcal{Z}}_{SW,r})$ that generate an action of $\mathbb{Z}[\mathbb{U}] \otimes (\wedge^*(H_1(Y;\mathbb{Z})/\text{torsion}))$ on the $\partial_{SW,Y}$-homology. Each such endomorphism is defined by the coefficients that appear in the relevant version of (1.23).

Consider first the endomorphism that generates the $\mathbb{Z}[\mathbb{U}]$ factor. Fix $[\mathfrak{c}]$ and $[\mathfrak{c}']$ so as to specify the corresponding version of $w_{[\mathfrak{c}],[\mathfrak{c}']}$. The specification of these coefficients requires the choice of a point $p \in \mathcal{H}_0$. Reintroduce $\mathfrak{p} \in \mathcal{P}_\mu$ from Part 6 and use $\mathcal{M}_2(\mathfrak{c},\mathfrak{c}')$ to denote the set of instanton solutions to the r and $\mathfrak{g} = \mathfrak{e}_\mu + \mathfrak{p}$ version of (1.20) with corresponding Fredholm index $\iota_{(\cdot)}$ is equal to 2, with $s \to -\infty$ limit equal to $\mathfrak{c}$ and with $s \to \infty$ limit in the $\mathcal{G}_{M_\Lambda}$ -orbit of $\mathfrak{c}'$. Use $\mathcal{M}_{2,p}(\mathfrak{c},\mathfrak{c}')$ to denote the subset of $\mathcal{M}_2(\mathfrak{c},\mathfrak{c}')$ that is characterized as follows: A given instanton $\mathfrak{d} = (A, \psi = (\alpha, \beta))$ is a member if and only if $\alpha|_{s=0}$ vanishes at p. The upcoming Proposition 1.2 asserts in part that $\mathcal{M}_{2,p}(\mathfrak{c},\mathfrak{c}')$ is a finite set if r is large. Granted that this is so, the coefficient $w_{[\mathfrak{c}],[\mathfrak{c}']}$ is given as a sum that is indexed by the instantons from $\mathcal{M}_{2,p}(\mathfrak{c},\mathfrak{c}')$. The contribution from each such instanton is specified using the rules in Chapter 23 of [KM]. Parts 3 and 4 of Section 1b in [T5] describe these same rules in the case when $\hat{a}$ is replaced by a contact 1-form and $w$ is replaced by the latter's exterior derivative.

The specification of the various endomorphisms that are meant to generate the $\wedge^*(H_1(Y;\mathbb{Z})/\text{torsion})$ action on the $\partial_{SW}$-homology requires the reintroduction of the set of 1-cycles $\{\{[\gamma^{(z)}]\}_{z\in¥}, \{\hat{\imath}_p\}_{p\in\Lambda}\}$ from Part 4 of Section 1b. Each cycle from this set labels a corresponding endomorphism. Let $\hat{\imath}$ denote such a cycle. Use $w^{\hat{\imath}}_{[\mathfrak{c}],[\mathfrak{c}']}$ to denote any given $[\mathfrak{c}]$, $[\mathfrak{c}']$ coefficient in $\hat{\imath}$'s version of (1.23). This coefficient is a weighted sum of intersection numbers that are defined using the elements in $\mathcal{M}_1(\mathfrak{c},\mathfrak{c}')$ whereby the contribution of a given instanton $(A, \psi = (\alpha, \beta))$ to the sum is either +1 or -1 times the algebraic intersection number between $\alpha^{-1}(0)$ and the locus $\mathbb{R} \times \hat{\imath}$ in $\mathbb{R} \times Y$. The rules for assigning a +1 or -1 weight to the intersection number are laid out in Chapter 23 of [KM]. Part 3 of Section 1b in [T5] describe these same rules in the case when $\hat{a}$ is replaced by a contact 1-form and $w$ is replaced by the latter's exterior derivative.

**Proposition 1.2**: *There exists* $\kappa \geq 1$ *with the following significance: Fix* $r \geq \kappa$ *and* $\mu \in \Omega$ *with* $\mathcal{P}$*-norm less than 1 such that all solutions to the* $(r,\mu)$ *version of (1.13) are non-* *degenerate and holonomy nondegenerate. If* $\mathfrak{p} \in \mathcal{P}$ *has small* $\mathcal{P}$*-norm and is described by (1.22), then the rules given in [KM] for specifying the coefficients for the just described endomorphisms of* $\mathbb{Z}(\hat{\mathcal{Z}}_{\mathrm{SW,r}})$ *define an action of* $\mathbb{Z}[\mathbb{U}] \otimes (\wedge^*(H_1(Y;\mathbb{Z})/\mathrm{torsion}))$ *on the* $\partial_{\mathrm{SW}}$ *homology. The generator of the action of the* $\mathbb{Z}[\mathbb{U}]$ *factor decreases the relative grading by* 2 *and those that generate the action of* $H_1(Y;\mathbb{Z})/\mathrm{torsion}$ *decrease the relative grading by 1. In addition, all of these endomorphisms map the submodule* $\mathbb{Z}(\hat{\mathcal{Z}}^{\geq}_{\mathrm{SW,r}})$ *to itself.*

This proposition is also proved in Section 7a.

*Part 8*: The formal adjoint of $\partial_{\mathrm{SW}}$ on the $\mathbb{Z}$-module $\mathrm{Hom}(\mathbb{Z}(\hat{\mathcal{Z}}_{\mathrm{SW,r}});\mathbb{Z})$ defines the differential for what is formally a version of Seiberg-Witten Floer cohomology. This formal adjoint of $\partial_{\mathrm{SW}}$ is denoted by $\partial^*_{\mathrm{SW}}$. The endomorphism $\partial^*_{\mathrm{SW}}$ sends a given basis element $[\mathfrak{c}]$ in $\hat{\mathcal{Z}}_{\mathrm{SW,r}}$ to

$$\partial^*_{\mathrm{SW}}[\mathfrak{c}] = \sum_{[\mathfrak{c}'] \in \hat{\mathcal{Z}}_{\mathrm{SW,r}}} w_{[\mathfrak{c}'],[\mathfrak{c}]}\,[\mathfrak{c}'] \tag{1.24}$$

This endomorphism increases the relative $\mathbb{Z}/p_M\mathbb{Z}$ grading by 1 and has square zero. The resulting $\partial^*_{\mathrm{SW}}$ homology groups enjoy an action of $\mathbb{Z}[\mathbb{U}] \otimes (\wedge^*(H_1(M;\mathbb{Z})/\mathrm{torsion})$ with the $\mathbb{Z}[\mathbb{U}]$ generator now increasing grading by 2 and the generators of $H_1(M;\mathbb{Z})/\mathrm{torsion}$ increasing grading by 1. The generators of this action come from the adjoints of the endomorphisms that are defined in Part 7.

**Proposition 1.3**: *There exists* $\kappa \geq 1$ *with the following significance: Fix* $r \geq \kappa$ *and* $\mu \in \Omega$ *with* $\mathcal{P}$*-norm less than 1 such that all solutions to the* $(r,\mu)$ *version of (1.13) are non-degenerate. Suppose that* $\mathfrak{p} \in \mathcal{P}_\mu$ *has small* $\mathcal{P}$*-norm and is described by (1.22). Then the expression on the right hand side in (1.24) defines* $\partial^*_{\mathrm{SW}}$ *as a square zero endomorphism of* $\mathbb{Z}(\hat{\mathcal{Z}}_{\mathrm{SW,r}})$*. The adjoints of the endomorphism of* $\mathbb{Z}(\hat{\mathcal{Z}}_{\mathrm{SW,r}})$ *that are defined in Part 7 likewise map* $\mathbb{Z}(\hat{\mathcal{Z}}_{\mathrm{SW,r}})$ *to itself and so define an action of* $\mathbb{Z}[\mathbb{U}] \otimes (\wedge^*(H_1(M;\mathbb{Z})/\mathrm{torsion})$ *on the homology groups of* $\partial^*_{\mathrm{SW}}$.

Proposition 1.3 is likewise proved in Section 7a.

Let $\hat{\mathcal{Z}}^{<}_{\mathrm{SW,r}} \subset \hat{\mathcal{Z}}_{\mathrm{SW,r}}$ denote the subset that corresponds via the identification in (1.19) to $\mathcal{Z}_{\mathrm{SW,r}} \times \{\ldots, -2, -1\}$. The endomorphism $\partial^*_{\mathrm{SW}}$ preserves the submodule of $\mathbb{Z}(\hat{\mathcal{Z}}^{<}_{\mathrm{SW,r}})$ as do those that give the generators of $\mathbb{Z}[\mathbb{U}] \otimes (\wedge^*(H_1(M;\mathbb{Z})/\mathrm{torsion}))$ action.

Granted what was just said, introduce $H^{\infty}_{\mathrm{SW,r}}$, $H^{-}_{\mathrm{SW,r}}$ and $H^{+}_{\mathrm{SW,r}}$ to denote the respective $\partial^*_{\mathrm{SW}}$-homology on the chain complexes $\mathbb{Z}(\hat{\mathcal{Z}}_{\mathrm{SW,r}})$, $\mathbb{Z}(\hat{\mathcal{Z}}^{<}_{\mathrm{SW,r}})$ and $\mathbb{Z}(\hat{\mathcal{Z}}_{SW,r})/\mathbb{Z}(\hat{\mathcal{Z}}^{<}_{SW,r})$. Each of these homology groups has a relative $\mathbb{Z}/p_M\mathbb{Z}$ grading, and each admits an action of $\mathbb{Z}[\mathbb{U}] \otimes (\wedge^*(H_1(Y;\mathbb{Z})/\text{torsion})$. Moreover, the latter are intertwined by the long exact sequence that is induced by the evident short exact sequence.

The next proposition speaks to the r-dependence of these $\partial^*_{SW}$-homology groups.

**Proposition 1.4**: *The versions of* $\kappa$ *that appear in Propositions 1.2 and 1.3 can be chosen so that the following is true: Suppose that* $r_1, r_2 \geq \kappa$, *and that* $(\mu_1, \mathfrak{p}_1)$ *and* $(\mu_2, \mathfrak{p}_2)$ *are pairs in* $\Omega \times \mathcal{P}$ *such that* $\mu_1$ *and* $\mu_2$ *have* $\mathcal{P}$*-norm less than 1, such that* $\mathfrak{p}_1$ *and* $\mathfrak{p}_2$ *have* $\mathcal{P}$*-norm much less than 1, and such that the conclusions of Propositions 1.1 and 1.2 hold for the data sets* $(r_1, \mu_1, \mathfrak{p}_1)$ *and* $(r_2, \mu_2, \mathfrak{p}_2)$. *Use these respective data sets to define the corresponding* $r = r_1$ *and* $r = r_2$ *versions of the groups* $H^{\infty}_{SW,r}$, $H^{-}_{SW,r}$ *and* $H^{+}_{SW,r}$.

- *There is a canonical isomorphism between the respective* $(r_1, \mu_1, \mathfrak{p}_1)$ *and* $(r_2, \mu_2, \mathfrak{p}_2)$ *versions of* $H^{\infty}_{SW,r}$ *that preserves the relative* $\mathbb{Z}/p_M\mathbb{Z}$ *gradings and intertwines the respective actions of* $\mathbb{Z}[\mathbb{U}] \otimes (\wedge^*(H_1(Y;\mathbb{Z})/\text{torsion}))$.
- *This canonical isomorphism maps the* $(r_1,\mu_1,\mathfrak{p}_1)$ *version of* $H^{-}_{SW,r}$ *isomorphically to the* $(r_2, \mu_2, \mathfrak{p}_2)$ *of* $H^{-}_{SW,r}$ *version, it induces an isomorphism between the two versions of* $H^{+}_{SW,r}$ *and it intertwines the respective long exact sequence homomorphisms.*
- *This canonical isomorphism is induced by a chain complex homomorphism from the* $(r_1,\mu_1,\mathfrak{p}_1)$ *version of* $\mathbb{Z}(\hat{\mathcal{Z}}_{SW,r})$ *to the* $(r_2,\mu_2,\mathfrak{p}_2)$ *version of* $\mathbb{Z}(\hat{\mathcal{Z}}_{SW,r})$ *that maps the* $(r_1,\mu_1,\mathfrak{p}_1)$ *version of* $\mathbb{Z}(\hat{\mathcal{Z}}^{<}_{SW,r})$ *to the* $(r_2,\mu_2,\mathfrak{p}_2)$ *version.*

This proposition is proved in Section 7c.

The canonical isomorphisms described by Proposition 1.4 are henceforth used to identify distinct $(r, \mu, \mathfrak{p})$ versions of $H^{\infty}_{SW,r}$, $H^{-}_{SW,r}$ and $H^{+}_{SW,r}$ and so write these groups respectively as $H^{\infty}_{SW}$, $H^{-}_{SW}$ and $H^{+}_{SW}$.

*Part 9*: This last part of the subsection brings the orientation reversed version of Y into the story so as to connect with what is said in [KLTI]. What is said here explains why Theorems I.3.1 and I.3.2 in [KLTI] follow directly from Propositions 1.1-1.4.

The orientation reversed twin of Y is denoted here by $\overline{Y}$. So as to be clear, the orientation on Y is defined so that the inclusion map $M_\delta \to M$ is orientation preserving and that of $M_\delta$ into Y is orientation reversing. The orientation for $\overline{Y}$ is such that both of the inclusion maps $M_\delta \to M$ and $M_\delta \to \overline{Y}$ are orientation *preserving*. As noted in the introduction, the convention used here for which orientation signifies Y and which signifies $\overline{Y}$ is opposite the convention used in [KLTI].

As explained below, the groups $H^{\infty}_{SW}$, $H^{-}_{SW}$ and $H^{+}_{SW}$ are canonically isomorphic to certain Seiberg-Witten Floer homology groups on Y, these being the respective groups $H^\infty_*$, $H^-_*$, and $H^+_*$ that are defined at the end of Section I.3.2. To see the connection, write the first line of (1.13) as

$$F_A - r(*\psi^\dagger\tau\psi - iw) + \tfrac{1}{2}F_{A_K} = 0 . \tag{1.25}$$

Now introduce $\bar{*}$ to denote the Hodge star as defined by the orientation for $\bar{Y}$. The latter is equal to -*. Likewise, intruduce $\overline{cl}$ to denote the $\bar{Y}$ version of the Clifford multiplication map. The latter is equal to -cl and as a consequence, the $\bar{Y}$ version of $\psi^\dagger\tau\psi$ is equal to (-1) times the Y version. Granted these last two observations, what is written in (1.13) is the equation that results when $\bar{*}$ is applied to both sides of the top line in (I.3.1). Meanwhile the lower line in (I.3.1) is -1 times the lower line in (1.13).

What was just said canonically identifies $\hat{\mathcal{Z}}_{SW,r}$ with a corresponding equivalence class of solutions to (I.3.1), this denoted in [KLTI] by $\hat{\mathcal{Z}}_{SW,Y,r}$. To see about the relation between $\partial^*_{SW}$ and the differential on $\hat{\mathcal{Z}}_{SW,Y,r}$, a look at (1.20) leads to the following observation: Let $\mathfrak{c}_-$ and $\mathfrak{c}_+$ denote solutions to (1.13) and suppose that $s \to \mathfrak{d}(s)$ is an instanton solution on Y with $s \to -\infty$ limit equal to $\mathfrak{c}_-$ and $s \to \infty$ limit equal to $\mathfrak{c}_+$. Then the map $s \to \mathfrak{d}(-s)$ is an instanton solution to (I.3.3) with $s \to -\infty$ limit $\mathfrak{c}_+$ and $s \to \infty$ limit $\mathfrak{c}_-$. This last observation implies that the identification just described between $\hat{\mathcal{Z}}_{SW,Y,r}$ and $\hat{\mathcal{Z}}_{SW,r}$ extends in a linear fashion to give an isomorphism between the chain complex on $\bar{Y}$ that is used to define the afore-mentioned groups $H^\infty_*$, $H^-_*$, and $H^+_*$ from Section I.3.2 and the chain complex $\mathbb{Z}(\hat{\mathcal{Z}}^<_{SW,r})$ with the differential $\partial^*_{SW}$.

The conclusions of the preceeding two paragraphs make Theorems I.3.1 and I.3.2 immediate consequences of Propositions 1.1−1.4.

**d) Seiberg-Witten Floer homology and embedded contact homology**

This subsection describes the relationship between the Seiberg-Witten Floer chain complex from Section 1c and its $\partial^*_{SW}$ homology and the embedded contact homology chain complex from Section 1b. This is the content of Theorem 1.5. This relationship is the analog of that described by Theorem 4.5 in [T2].

Theorem I.3.3 from [KLTI] follows directly from what is said in Theorem 1.5 and what is said in Part 9 of the previous subsection.

The upcoming Theorem 1.5 refers to the filtration of $\hat{\mathcal{Z}}_{ech,M}$ given in (1.10). The theorem also refers to a certain subsets in the various L ≥ 1 versions of what is denoted in Part 5 of Section 1b by $\hat{\mathcal{Z}}_{ech,M}{}^{L}$. The subset in question is denoted in the theorem and in what follows by $\hat{\mathcal{Z}}_{ech,M}{}^{L,<}$ and it is defined as follows: Part 4 in Section II.1b defines a principal $\mathbb{Z}$-bundle isomorphism $\hat{\mathcal{Z}}_{ech,M} = \mathcal{Z}_{ech,M} \times \mathbb{Z}$. This isomorphism sends the equivalence class $(\Theta, Z)$ to the pair $(\Theta, k)$ when $Z$ has intersection number $k$ with $\gamma^{(z_0)}$. The subset $\hat{\mathcal{Z}}_{ech,M}{}^{L,<}$ corresponds via this isomorphism to $\mathcal{Z}_{ech,M}{}^{L} \times \{-\infty, \ldots, -1\}$.

**Theorem 1.5**: *Let $\mathbb{H}^\infty$, $\mathbb{H}^-$ and $\mathbb{H}^+$ denote finitely generated subgroups of the respective groups $H^\infty{}_{SW}$, $H^-{}_{SW}$ and $H^+{}_{SW}$. Given these subgroups, there exists $L_\mathbb{H}$ and given $L \ge L_\mathbb{H}$ there exists $L' \ge L$ with the following significance: Fix $r$ sufficiently large, and then fix a pair $(\mu, \mathfrak{p}) \in \Omega \times \mathcal{P}$ such that $\mu$ has $\mathcal{P}$-norm less than 1, such that $\mathfrak{p}$ has sufficiently small $\mathcal{P}$-norm, and such that Propositions 1.1-1.3 can be invoked to define the chain complex $(\mathbb{Z}(\hat{\mathcal{Z}}_{SW,r}), \partial^*{}_{SW})$, the subcomplex $\mathbb{Z}(\hat{\mathcal{Z}}^{<}_{SW,r})$ and the $\mathbb{Z}[\mathbb{U}] \otimes (\wedge^*(H_1(Y;\mathbb{Z})/\text{torsion}))$ action on the homology. There exists an injective principal $\mathbb{Z}$ bundle map $\hat{\Phi}^r : \hat{\mathcal{Z}}_{ech,M}{}^{L'} \to \hat{\mathcal{Z}}_{SW,r}$ that defines a $\mathbb{Z}$-module homomorphism*

$$\mathbb{L}^r : \mathbb{Z}(\hat{\mathcal{Z}}_{ech,M}{}^{L'}) \to \mathbb{Z}(\hat{\mathcal{Z}}_{SW,r})$$

*with the properties listed below.*

- *$\mathbb{L}^r$ reverses the sign of relative grading degrees.*
- *$\mathbb{L}^r$ induces monomorphisms from $\mathbb{Z}(\hat{\mathcal{Z}}_{ech,M}{}^{L',<})$ to $\mathbb{Z}(\hat{\mathcal{Z}}^{<}_{SW,r})$ and from $\mathbb{Z}(\hat{\mathcal{Z}}_{ech,M}{}^{L'})/\mathbb{Z}(\hat{\mathcal{Z}}_{ech,M}{}^{L',<})$ to $\mathbb{Z}(\hat{\mathcal{Z}}_{SW,r})/\mathbb{Z}(\hat{\mathcal{Z}}^{<}_{SW,r})$.*
- *$\mathbb{L}^r$ intertwines $\partial_{ech}$ with $\partial^*{}_{SW}$ and it also intertwines the endomorphisms that define the generators of the respective $\mathbb{Z}[\mathbb{U}] \otimes (\wedge^*(H_1(Y;\mathbb{Z})/\text{torsion}))$ actions on the $\partial_{ech}$ homology and $\partial^*{}_{SW}$ homology.*
- *Let $\mathbb{Q}_{ech}{}^{L}$ denote either $\mathbb{Z}(\hat{\mathcal{Z}}_{ech,M}{}^{L})$, $\mathbb{Z}(\hat{\mathcal{Z}}_{ech,M}{}^{L,<})$ or $\mathbb{Z}(\hat{\mathcal{Z}}_{ech,M}{}^{L})/\mathbb{Z}(\hat{\mathcal{Z}}_{ech,M}{}^{L,<})$ and let $\mathbb{Q}_{ech}{}^{L'}$ denote the $L'$ version. Use $\mathbb{Q}_{SW}$ to denote the corresponding $\mathbb{Z}(\hat{\mathcal{Z}}_{SW,r})$, $\mathbb{Z}(\hat{\mathcal{Z}}^{<}_{SW,r})$ or $\mathbb{Z}(\hat{\mathcal{Z}}_{SW,r})/\mathbb{Z}(\hat{\mathcal{Z}}^{<}_{SW,r})$ as the case may be. If $\varsigma \in \mathbb{Q}_{ech}{}^{L}$ is such that $\mathbb{L}^r(\varsigma) = \partial^*{}_{SW}\mathfrak{z}$ for some $\mathfrak{z} \in \mathbb{Q}_{SW}$, then $\varsigma = \partial_{ech}\varsigma'$ for some $\varsigma' \in \mathbb{Q}_{ech}{}^{L'}$.*
- *The subgroups $\mathbb{H}^\infty$, $\mathbb{H}^-$, and $\mathbb{H}^+$ are represented by elements in the respective $\mathbb{L}^r$ images of $\mathbb{Z}(\hat{\mathcal{Z}}_{ech,M}{}^{L})$, $\mathbb{Z}(\hat{\mathcal{Z}}_{ech,M}{}^{L,<})$ and $\mathbb{Z}(\hat{\mathcal{Z}}_{ech,M}{}^{L})/\mathbb{Z}(\hat{\mathcal{Z}}_{ech,M}{}^{L,<})$.*

*Suppose that distinct choices for $(r, \mu, \mathfrak{p})$ are suitable for defining $\mathbb{Z}(\hat{\mathcal{Z}}_{SW,r})$, the differential $\partial^*{}_{SW}$, and the subcomplex $\mathbb{Z}(\hat{\mathcal{Z}}^{<}_{SW,r})$. If the respective values of $r$ are large enough to define $\mathbb{L}^r$ on $\mathbb{Z}(\hat{\mathcal{Z}}_{ech,M}{}^{L'})$, and in any event, sufficiently large, then homomorphism from the third bullet of Proposition 1.4 can be chosen that intertwine the resulting versions of $\mathbb{L}^r$.*

This theorem is proved in Sections 7d and 7f.

**e) Functions on $\mathrm{Conn}(E) \times C^\infty(Y; \mathbb{S})$**

This subsection introduces functions on $\mathrm{Conn}(E) \times C^\infty(Y; \mathbb{S})$ that play essential roles in the story. In what follows, $\mathfrak{c} = (A, \psi = (\alpha, \beta)) \in \mathrm{Conn}(E) \times C^\infty(Y; \mathbb{S})$ is a given element.

The first function is the *Chern-Simons function*. Reintroduce the chosen fiducial connection $A_E$ from Part 3 of Section 1c and write $A = A_E + \hat{a}_A$ with $\hat{a}_A$ an $i\mathbb{R}$-valued 1-form. The Chern-Simons function sends A to

$$\mathfrak{cs}(A) = -\int_Y \hat{a}_A \wedge d\hat{a}_A - 2\int_Y \hat{a}_A \wedge (F_{A_E} + \tfrac{1}{2} F_{A_K}) . \tag{1.26}$$

Note that $\mathfrak{cs}$ is invariant only under the action on Conn(E) of the subgroup in $C^\infty(Y; S^1)$ of maps u that define classes in $H^1(Y; S^1)$ that have cup product pairing zero with the first Chern class of $\det(\mathbb{S})$. This subgroup is denoted by $\mathcal{G}_{\mathbb{S}}$.

The second function is

$$\mathrm{W}(A) = i \int_Y \hat{a}_A \wedge w \tag{1.27}$$

This function is invariant only under the action on Conn(E) of $\mathcal{G}_{\mathbb{S}}$.

The next function is denoted by $\mathfrak{a}$. The critical points of $\mathfrak{a}$ are the solutions to (1.13). This function is given by

$$\mathfrak{a} = \tfrac{1}{2}\mathfrak{cs} - r\mathrm{W} + \mathfrak{e}_\mu + r\int_Y \psi^\dagger D_A \psi . \tag{1.28}$$

The *spectral flow* function is denoted by $\mathfrak{f}_\mathfrak{s}$. This function is constant on the components of the complement in $\mathrm{Conn}(E) \times C^\infty(Y; \mathbb{S})$ of the codimension 1 subvariety where the operator $\mathfrak{L}_{(\cdot)}$ in (1.17) has nontrivial kernel. It is discontinuous across this subvariety, but in any event it is locally bounded. A precise definition can be found in [T3]. What follows defines $\mathfrak{f}_\mathfrak{s}$ where $\mathfrak{L}_{(\cdot)}$ has trivial kernel. The definition of $\mathfrak{f}_\mathfrak{s}$ requires the choice of a section of $\mathbb{S}$, this denoted by $\psi_E$. This section must be chosen so that the $\mathfrak{c}_E = (A_E, \psi_E)$ and the r = 1 version of the operator $\mathfrak{L}_{(\cdot),r}$ has trivial kernel. The existence of such a section can be established using the Bochner-Weitzenboch formula in (A.12) of the appendix for the square of the operator $\mathfrak{L}_{(\cdot),r}$. Now suppose that $\mathfrak{c}$ is a given pair in $\mathrm{Conn}(E) \times C^\infty(Y; \mathbb{S})$ with the kernel of $\mathfrak{L}_{\mathfrak{c},r} = \{0\}$. Select a smooth map $\mathfrak{c}(\cdot)$ from [0,1] to $\mathrm{Conn}(E) \times C^\infty(Y; \mathbb{S})$ with $\mathfrak{c}(0) = \mathfrak{c}_E$ and $\mathfrak{c}(1) = \mathfrak{c}$ and a smooth map, $r(\cdot)$ from [0, 1] to [1, r] with $r(0) = 1$ and $r(1) = \mathrm{r}$. The function $\mathfrak{f}_\mathfrak{s}$ assigns to $\mathfrak{c}$ the spectral flow for the [0,1]-parametrized path of operators $\{\mathfrak{L}_{\mathfrak{c}(\tau),r(\tau)}\}_{\tau\in[0,1]}$. Note that $\mathfrak{f}_\mathfrak{s}(\mathfrak{c})$ is independent of the chosen maps $\mathfrak{c}(\cdot)$ and $r(\cdot)$. So defined, the function $\mathfrak{f}_\mathfrak{s}$ is also constant on the $\mathcal{G}_\mathbb{S}$ orbit of $\mathfrak{c}$.

Neither $\mathfrak{a}$, $\mathfrak{cs}$, $\mathrm{W}$, nor $\mathfrak{f}_\mathfrak{s}$ are invariant with respect to the action of $C^\infty(Y; S^1)$ on Conn(E) although all are invariant with respect to the action of the subgroup $\mathcal{G}_\mathbb{S}$. However, the following functions are invariant under the full action of $C^\infty(Y; S^1)$

$$\mathfrak{cs}^\mathfrak{f} = \mathfrak{cs} - 4\pi^2 \mathfrak{f}_\mathfrak{s}\,, \quad \mathrm{W}^\mathfrak{f} = \mathrm{W} - 2\pi\mathfrak{f}_\mathfrak{s} \quad \textit{and} \quad \mathfrak{a}^\mathfrak{f} = \mathfrak{a} + 2\pi(\mathrm{r} - \pi)\,\mathfrak{f}_\mathfrak{s}\,. \tag{1.29}$$

The last of the functions of interest is denoted by $\mathrm{M}$ and it is given by

$$\mathrm{M} = \mathrm{r}\int_Y (1 - |\alpha|^2)\ . \tag{1.30}$$

The question of bounding $\mathrm{M}$ on a solution to (1.14) or along the path of an instanton is a central concern in what follows.

## 2. Solutions to the Seiberg-Witten equations on Y

The solutions to the large r versions of (1.13) have certain properties that play central roles in many of the subsequent arguments that supply input to the proofs of Proposition 1.14 and Theorem 1.15. These properties are given by the various lemmas and propositions in the first two subsections that follow. The third subsection contains the proof of a proposition in the first subsection.

### a) A priori properties of solutions to (1.13)

The lemmas in the first parts of this subsection consider the pointwise behavior of solutions to (1.13). The second part of the subsection concerns the locus in Y where the curvature 2-form is large. This second part also talks about the function $\mathrm{M}$ in (1.30). The third part of what follows discusses the spectral flow function $\mathfrak{f}_\mathfrak{s}$ and the final part discusses the functions $\mathrm{W}$, $\mathfrak{cs}$ and $\mathfrak{a}$ from Section 1e.

*Part 1*: The upcoming Lemmas 2.1-2.3 have close analogs in Sections 2a of [T6], in Section 6 of [T1] and in Section 3 of [T7]. When the proof of a given lemma here differs only slightly from its partner in one of these references, then only the salient differences (if any) are noted.

The first lemma speaks to the size of the $C^\infty(Y; \mathbb{S})$ component of a solution.

**Lemma 2.1**: *There exists* $\kappa > 1$ *with the following significance: Fix* $\mu \in \Omega$ *with $\mathcal{P}$-norm less that* 1 *and* $r \geq \kappa$. *Let* $(A, \psi = (\alpha, \beta))$ *denote a solution to the corresponding version* $(r, \mu)$ *version of (1.13). Then*

- $|\alpha| \leq 1 + \kappa r^{-1}$.
- $|\beta|^2 \leq \kappa r^{-1}(1 - |\alpha|^2) + \kappa^3 r^{-2}$ .
- $|\nabla_A\alpha|^2 \leq \kappa r(1 - |\alpha|^2) + \kappa^3$.
- $|\nabla_A\beta|^2 \leq \kappa (1 - |\alpha|^2) + \kappa^3 r^{-1}$

*In addition, for each* $q \geq 1$, *there exist* $\kappa_q \in (1, \infty)$ *which is independent of* $(A, \psi)$, *r and* $\mu$, *and is such that*

- $|\nabla_A{}^q\alpha| + r^{1/2}|\nabla_A{}^q\beta| \leq \kappa_q r^{q/2}$.

***Proof of Lemma 2.1***: The lemma and its proof differ only in notation from Lemma 2.3 in [T6] and the latter's proof.

Given the equation in the top line of (1.13s), what follows is an immediate consequence of the first two bullets in Lemma 2.1:

$$|B_A| = i{*}(\hat{a} \wedge {*}B_A) + \mathfrak{e} \tag{2.1}$$

where $|\mathfrak{e}| \leq c_0$.

The second lemma addresses the size of the connection A.

**Lemma 2.2**: *There exists* $\kappa > 1$ *with the following significance: Fix* $\mu \in \Omega$ *with $\mathcal{P}$-norm less that* 1 *and* $r \geq \kappa$. *Let* $(A, \psi = (\alpha, \beta))$ *denote a solution to the corresponding version* $(r, \mu)$ *version of (1.13). There is a map* $u \in C^\infty(Y; \mathbb{S})$ *which is homotopic to the identity and such that* $A - u^{-1}du$ *can be written as* $A_E + \hat{a}^\perp + \mathfrak{p}_A$ *where* $\mathfrak{p}_A$ *is a harmonic 1-form and* $\hat{a}^\perp$ *is coclosed, $L^2$-orthogonal to the space of harmonic 1-forms, and is such that* $|\hat{a}^\perp| \leq \kappa(|M|^{1/3} r^{2/3} + 1)$.

***Proof of Lemma 2.2***: Given (2.1), the proof is identical to that of Lemma 2.4 in [T1].

The third lemma in this series extends what is said in Lemma 2.1 with some precise bounds for the size of $1 - |\alpha|^2$ and the covariant derivatives of $\alpha$ and $\beta$.

**Lemma 2.3**: *There exists* $\kappa > 1$ *with the following significance: Fix* $\mu \in \Omega$ *with $\mathcal{P}$-norm less that* 1 *and* $r \geq 1$. *Let* $(A, \psi = (\alpha, \beta))$ *denote a solution to the* $(r, \mu)$ *version of (1.13). Let* $Y_* \subset Y$ *denote the subset of points where* $1 - |\alpha|^2 \geq \kappa^{-1}$. *Then*

$$|1 - |\alpha|^2| \leq (e^{-\sqrt{r}\,\mathrm{dist}(\cdot, Y_*)/\kappa} + \kappa r^{-1}) .$$

***Proof of Lemma 2.3***: The manipulations done to prove Proposition 4.4 from the article *SW => Gr* in [T8] can be repeated here to obtain the desired inequality.

*Part 2*: The upcoming Proposition 2.4 describes the zero locus of the $\alpha$ part of a solution to a given large r version of (1.13) at the points in Y with distances greater than $c_0\, r^{-1/2}$ from the curves in the set $\cup_{p\in\Lambda}\{\hat\gamma_p^+, \hat\gamma_p^-\}$. The proposition refers to the connection $\hat{A}$ that is defined from any given pair of connection on E and section of E in (1.16).

**Proposition 2.4**: *There exists* $\kappa > 1$ *with the following significance: Fix* $r \ge \kappa$ *and* $\mu \in \Omega$ *with* $\mathcal{P}$*-norm bounded by* 1 *and suppose that* $(A, \psi = (\alpha, \beta))$ *is a solution to the corresponding* $(r, \mu)$ *version of (1.13). Let* $Y_r \subset Y$ *denote the set of points with distance greater than* $\kappa^2\, r^{-1/2}$ *from the curves in* $\cup_{p\in\Lambda}\{\hat\gamma_p^+, \hat\gamma_p^-\}$*. The zero locus of* $\alpha$ *in the closure of* $Y_r$ *is transversal and it consists of the disjoint union of at most* G *components with each a properly embedded arc or circle. The zero locus of* $\alpha$ *has the following additional properties*:

- *The tangent line to each component has distance at most* $\kappa r^{-1/2}$ *from* $\nu$.
- *Each component lies where* $1 - 3\cos^2\theta > 0$.
- *The intersection of the zero set with* $M_\delta$ *consists of* G *properly embedded segments that pair the index 1 and index 2 critical points of the incarnation of f as a function on* M *in the sense that distinct segments start on the boundary of the radius* $\delta$ *coordinate balls about distinct index 1 critical points of f and end on the boundary of the radius* $\delta$ *coordinate balls about distinct index 2 critical points.*
- *The absolute value of* $1 - |\alpha|^2$ *is less than* $\frac{1}{32}\kappa^{-1}$ *at all points with distance greater than* $\kappa r^{-1/2}$ *from the zero locus of* $\alpha$ *in* Y; *and less than* $\kappa r^{-1}$ *at all points with distance* $\kappa(\ln r)\, r^{-1/2}$ *or more from the zero locus of* $\alpha$ *in* Y.
- *The 2-form* $\frac{i}{2\pi} F_{\hat A}$ *has compact support and integral* 1 *on any disk in* $Y_r$ *that intersects* $\alpha^{-1}(0)$ *transversally at its center point, is otherwise disjoint from* $\alpha^{-1}(0)$, *and has closure with all boundary points at distance at least* $\kappa r^{-1/2}$ *from* $\alpha^{-1}(0)$.

The proof of Proposition 2.4 is given in Section 2c. The first assertion of the next lemma is little more than a corollary to Proposition 2.4. The second assertion refers to the 1-form $\upsilon_\Diamond$ from (1.5).

**Lemma 2.5**: *There exists* $\kappa > 1$ *with the following significance: Fix* $r \ge \kappa$ *and* $\mu \in \Omega$ *with* $\mathcal{P}$*-norm bounded by* 1 *and suppose that* $(A, \psi = (\alpha, \beta))$ *is a solution to the corresponding* $(r, \mu)$ *version of (1.13).*

- *Set* $\text{M} = r\int_Y (1 - |\alpha|^2)$. *Then* $-\kappa \le \text{M} \le \kappa \ln r$.
- $r\int_Y |\upsilon_\lozenge|^2 |1 - |\alpha|^2| \le \kappa$.

***Proof of Lemma 2.5***: The lower bound on M follows directly from Lemma 2.1. To obtain the asserted upper bound, use Proposition 2.4 to characterize the zero locus of $\alpha$ in $Y_r$. In particular, Lemma 2.3 with the third bullet of Proposition 2.4 bound $1 - |\alpha|^2$ at distance $\rho$ in $Y_r$ from $\alpha^{-1}(0) \cap Y_r$ by $c_0\, e^{-\sqrt{r}\rho/c_0}$. It follows from the first bullet of Proposition 2.4 and the formula for $v$ in (1.3) that the length $\alpha^{-1}(0) \cap Y_r$ is at most $c_0 \ln r$. These bounds together imply that the $Y_r$ contribution to the integral that defines M is at most $c_0 \ln r$. Meanwhile, the volume of $Y - Y_r$ is at most $c_0 r^{-1}$ and so the $Y - Y_r$ contribution to the integral that defines M is at most $c_0 \ln r$.

To prove the assertion of the second bullet, note that integral over Y of $\upsilon_\lozenge \wedge \frac{i}{2\pi} F_A$ is equal to the pairing between $c_1(E)$ and the class in $H_2(Y; \mathbb{R})$ that is Poincaré dual to the class that is defined by the De Rham cohomology by the closed form $\upsilon_\lozenge$. With this fact in mind, write $\upsilon_\lozenge$ as $q_\lozenge \hat{a} + \mathfrak{b}$ where $\mathfrak{b}$ annihilates the vector field $v$. Note in particular that what is said in Part 4 of Section 1a can be used to see that $|\upsilon_\lozenge|^2 \le c_0 q_\lozenge$. Granted this last point, use the top equation in (1.13) to see that

$$i{*}(\upsilon_\lozenge \wedge F_A) = r q_\lozenge (1 - |\alpha|^2) + \mathfrak{r} \tag{2.2}$$

where $|\mathfrak{r}| \le c_0 r |\alpha| |\beta|\, |\upsilon_\lozenge|$. Given that $|\upsilon_\lozenge| \le c_0 q_\lozenge^{1/2}$, the first and second bullets in Lemma 2.1 with (2.2) find

$$i{*}(\upsilon_\lozenge \wedge F_A) \ge \tfrac{1}{2} r q_\lozenge |1 - |\alpha|^2| - c_0 . \tag{2.3}$$

The lemma's assertion follows from (2.3) with second appeal to the bound $c_0^{-1}|\upsilon_\lozenge|^2 \le q_\lozenge$.

*Part 3*: The spectral flow function $\mathfrak{f}_\mathfrak{s}$ plays a central role in the proofs of the Proposition 1.4 and Theorem 1.5. The upcoming Proposition 2.6 supplies a crucial bound for its absolute value. To set the stage for this proposition, reintroduce from Part 4 of Section 1b the set $\{\gamma^{(z)}\}_{z \in \text{¥}}$ of closed integral curves of $v$. This set has $1 + b_1(M)$ elements. Each curve in this set lies in $M_\delta \cap \mathcal{H}_0$ and it has distance $c_0^{-1}$ or more from any segment of an integral curve of $v$ in the $f^{-1}(1,2)$ part of $M_\delta$ that starts on the boundary of the radius $\delta$ coordinate ball about an index 1 critical point of $f$ in M and ends on the boundary of the radius $\delta$ coordinate ball about an index 2 critical point of $f$. Associate to each $z \in ¥$ the map $\mathrm{x}^{(z)}$ from Conn(E) to $\mathbb{R}$ given by following rule: Let A denote any given connection on E, write A as $A = A_E + \hat{a}_A$ and set

$$\mathrm{x}^{(z)}(A) = \tfrac{i}{2\pi} \int_{\gamma^{(z)}} \hat{a}_A \, .$$

(2.4)

Use $[\gamma^{(z)}]$ in what follows to denote the class in $H_1(M;\mathbb{Z})$/torsion that is defined by a given $z \in ¥$ loop $\gamma^{(z)}$. The set of such cycles generates the image of the Poincaré dual of the classes in $H^2(Y; \mathbb{Z})$ that annihilate the $\oplus_{p\in\Lambda} H_2(\mathcal{H}_p;\mathbb{Z})$ summand in (1.4)'s depiction of $H_2(Y;\mathbb{Z})$. As the first Chern class of det($\mathbb{S}$) annihilates this summand, the image of its Poincaré dual in $H_1(M;\mathbb{Z})$/torsion can be written as $\sum_{z\in ¥} c_{\mathbb{S},z}[\gamma^{(z)}]$ with coefficients $\{c_{\mathbb{S},z}\}_{z\in ¥} \in \mathbb{Z}$. Use $\mathrm{x}_{\mathbb{S}}$ to denote the corresponding map $\sum_{z\in ¥} c_{\mathbb{S},z}\, \mathrm{x}^{(z)}$.

What follows is a consequence of the fact that the classes from the set $\{[\gamma^{(z)}]\}_{z\in ¥}$ are linearly independent in $H_1(Y;\mathbb{Z})$/torsion: Let A denote a connection on E. There is a smooth map $u: Y \to S^1$ such that $A - u^{-1}du$ obeys $0 \le \mathrm{x}^{(z)}(A - u^{-1}du) < 1$ for each $z \in ¥$. Note that u can be chosen so that $A - u^{-1}du - A_E = \hat{a}_A - u^{-1}du$ is a coclosed 1-form.

**Proposition 2.6**: *There exists* $\kappa > 1$ *with the following significance: Suppose that* $r > \kappa$ *and that* $\mu \in \Omega$ *has* $\mathcal{P}$*-norm less than 1. Let* $(A, \psi)$ *denote a non-degenerate solution to the* $(r, \mu)$ *version of (1.13). Then* $|\mathfrak{f}_{\mathfrak{s}}(A, \psi) - \mathrm{x}_{\mathbb{S}}(A, \psi)| \le \kappa$.

Section Bc of the appendix extends the function $|\mathfrak{f}_{\mathfrak{s}}|$ as a piecewise constant function on the whole of Conn(E) $\times\, C^\infty(Y;\mathbb{S})$. This understood, the assertion in Proposition 2.6 also holds in the case when the $(A, \psi)$ version of (1.17) has non-trivial kernel.

The proof of Proposition 2.6 is in the Appendix. The placement of the proof in an appendix is not a reflection of the import of Proposition 2.6; this proposition is absolutely crucial with regards to what is said subsequently about instanton solutions to (1.20). The proof is in the appendix as it is long and as the notions that enter are not used elsewhere

*Part 4*: The proposition that follows supplies a priori bounds for the functions $\mathfrak{cs}$ in (1.26), the function $\mathrm{W}$ in (1.27) and the function $\mathfrak{a}$ in (1.28).

**Proposition 2.7**: *There exists* $\kappa > 1$ *with the following significance: Fix* $r \ge \kappa$ *and* $\mu \in \Omega$ *with* $\mathcal{P}$*-norm bounded by* 1 *and suppose that* $(A, \psi)$ *is a solution to the* $(r, \mu)$ *version of (1.13). Then*

- $|\mathfrak{cs}^f| \le \kappa(r^{2/3} \mathrm{M}^{4/3} + \mathrm{M} + 1)$ ,
- $|\mathrm{W}^f - \mathrm{M}| < \kappa$ ,
- $-\kappa r(\mathrm{M}+1) \le \mathfrak{a}^f \le \kappa r$ .

***Proof of Proposition 2.7***: The assertion in the third bullet about $\mathfrak{a}^f$ follows from the assertions about $\mathfrak{cs}^f$ and $\mathrm{W}^f$. To prove the asserted bound for $\mathfrak{cs}^f$, use the Green's function for the operator $d+d^*$ to construct a smooth, coclosed 1-form on $Y-(\cup_{z\in¥}\gamma^{(z)})$ with the following properties: Let $B_{\mathbb{S}}$ denote this 1-form. Then $|B_{\mathbb{S}}| \le c_0 \sum_{z\in¥} \mathrm{dist}(\cdot,\gamma^{(z)})^{-1}$ and

$$\frac{i}{2\pi}\int_Y \hat{a}\wedge(F_{A_E} + \tfrac{1}{2}F_{A_K}) = \sum_{z\in¥} C_{\mathbb{S},z}\int_{\gamma^{(z)}} \hat{a} + \frac{i}{2\pi}\int_Y B_{\mathbb{S}}\wedge d\hat{a} \tag{2.5}$$

with â being any given 1-form on Y. Granted (2.5), write $\mathfrak{cs}(A)$ as

$$-\int_Y \hat{a}_A \wedge d\hat{a}_A - 2\int_Y B_{\mathbb{S}}\wedge d\hat{a}_A + 4\pi\, x_{\mathbb{S}}\,. \tag{2.6}$$

To bound the integral of $B_{\mathbb{S}} \wedge d\hat{a}_A$, use the top equation in (1.13) to see that the 2-form $d\hat{a}_A$ differs from $B_A$ by a smooth form. This understood, use this same equation with (2.1) and Lemma 2.1 to bound $d\hat{a}_A$ at distances $c_0^{-1}$ or less from any curve in the set $\{\gamma^{(z)}\}_{z\in¥}$ by $c_0 r(1-|\alpha|) + c_0$ and use Lemma 2.3 to bound the latter by $c_0$. The absolute value of the contribution to the integral of $B_{\mathbb{S}} \wedge d\hat{a}_A$ from the radius $c_0^{-1}$ tubular neighborhood of any curve from $\{\gamma^{(z)}\}_{z\in¥}$ is therefore bounded by $c_0$. Meanwhile, the absolute value of the contribution to the integral of $B_{\mathbb{S}} \wedge d\hat{a}_A$ from the complement in Y of the union of these neighborhoods is less than $c_0 M$.

To bound the left most integral in (2.6), remark first that both of the integrals over Y in (2.6) do not change when A is replaced by $A - u^{-1}du$ with u being any map from Y to $S^1$. This understood, choose a map u so that the $L^2$-orthogonal projection of $\hat{a}_A - u^{-1}du$ has $L^2$ norm bounded by $c_0$. Having done this, use Lemma 2.2 to bound the left most integral in (2.6) by $c_0 r^{2/3} M^{4/3}$. Granted these bounds, Proposition 2.7's bound of $\mathfrak{cs}^f$ follows from (2.6) and Proposition 2.6.

Consider next the assertion made by second bullet of the proposition. Look at (1.3) and (1.6) to see that $w$ on the $|u| \le R + c_0 \ln\delta$ part of any $\mathfrak{p} \in \Lambda$ version of $\mathcal{H}_{\mathfrak{p}}$ can be written as $d\hat{a}$. As a consequence, the function $\chi$ can be used with the Green's function for the operator $d+d^*$ to construct a smooth 1-form on $Y-(\cup_{z\in¥}\gamma^{(z)})$ with the properties listed next. Use $B_w$ to denote this 1-form. The form $B_w$ is zero on the $|u| \le R + c_0 \ln\delta$ part of each $\mathfrak{p} \in \Lambda$ version of $\mathcal{H}_{\mathfrak{p}}$. In addition, $|B_w| \le c_0 \sum_{z\in¥} \mathrm{dist}(\cdot,\gamma^{(z)})^{-1}$ and

$$i\int_Y \hat{a}\wedge w = i\sum_{z\in¥} C_{\mathbb{S},z}\int_{\gamma^{(z)}} \hat{a} + i\int_Y \hat{a}\wedge d\hat{a} + i\int_Y B_w \wedge d\hat{a}\ , \tag{2.7}$$

with $\hat{a}$ being any given 1-form on Y.

Take $\hat{a}$ in (2.7) to be the 1-form $\hat{a}_A$. The left hand side of the $\hat{a}_A$ version of (2.7) is $\mathrm{W}(A)$. The term on the right hand side with the sum indexed by ¥ is $2\pi\, \mathrm{x}_{\mathbb{S}}$. Use the top equation in (1.13) with (1.30) to see that the integral of $\hat{a} \wedge d\hat{a}_A$ can be written as $-i\mathrm{M} + \mathfrak{r}_A$ with $|\mathfrak{r}_A| \leq c_0$. Meanwhile, $\mathrm{B}_w \wedge d\hat{a}_A$ can be written as $\mathrm{B}_w \wedge F_A + \mathfrak{q}_A$ where $F_A$ denotes A's curvature 2-form and $\mathfrak{q}_A$ is a 2-form with $|\mathfrak{q}_A| \leq c_0$. Granted the preceding, the second bullet of the proposition follows with a bound by $c_0$ on the absolute value of the integral over Y of the form $\mathrm{B}_w \wedge F_A$ . Such a bound follows from the second bullet of Lemma 2.5 and Lemma 2.1.

Given what was said in the preceding two paragraphs, the bound for $|\mathrm{W}^f - \mathrm{M}|$ given in Proposition 2.7 follows from (2.7) and Proposition 2.6.

## b) The vortex equations I

The proof of Proposition 2.4 in Section 2c invokes various properties of the *vortex equations* on $\mathbb{C}$. Properties of these equations are also used to prove Proposition 2.6 and are used elsewhere as well. This section introduces these equations and supplies what is needed for the proof of Proposition 2.4. More is said about these equations in Sections 3 and 4.

The vortex equations ask that a pair $(A_0, \alpha_0)$ of connection on a complex line bundle over $\mathbb{C}$ and section of this bundle obey

- $* F_{A_0} = -i(1 - |\alpha_0|^2)$
- $\bar{\partial}_{A_0}\alpha_0 = 0$ ,
- $|\alpha_0| \leq 1$ .

(2.8)

The notation here is such that $*$ denotes the Euclidean Hodge dual on $\mathbb{C}$, while $F_{\hat{A}_0}$ and $\bar{\partial}_{A_0}$ denote the respective curvature 2-form of $A_0$ and the d-bar operator defined by $A_0$ on the space of sections of the given complex line bundle. The solutions with $1 - |\alpha_0|^2$ integrable are discussed at length in Sections 1 and 2 of [T9]. Solutions to (2.8) are also solutions to (4.1) in [T6] so what is said in Proposition 4.2 in [T6] applies as well.

Two properties of the solutions to (2.8) are needed for the proof of Proposition 2.4 that are not stated explicitly in Sections 1 and 2 of [T9] or by Proposition 4.1 in [T6]. These are given by

**Lemma 2.8**: *Let* $(A_0, \alpha_0)$ *denote a solution to the vortex equations. Then* $|\alpha_0|$ *can not have a local, non-zero minimum. Given* $\varepsilon > 0$, *there exists* $\kappa > 1$ *with the following* *significance*: *Suppose that* $(A_0, \alpha_0)$ *is a solution to the vortex equations and* $|\alpha_0| < 1 - \varepsilon$ *at the origin in* $\mathbb{C}$. *Then* $|\alpha_0| < \varepsilon$ *at a point with distance at most* $\kappa$ *from the origin.*

***Proof of Lemma 2.8***: The function $|\alpha_0|$ can be written as $e^u$ on a set where it is nowhere zero with $u < 0$ a smooth function. The top equation in (2.8) requires that $-\Delta u = (1 - e^{2u})$ where $\Delta$ here denotes the Laplacian on $\mathbb{R}^2$. This understood, the first assertion of the lemma follows from the maximum principal. To prove the second assertion, suppose to the contrary that it is false for some $\varepsilon > 0$. The equations in (2.8) are uniformly ellipitic, and thus by taking limits with counter examples for the successive cases $\kappa = 1, 2, \ldots$ finds a solution $(A_0, \alpha_0)$ with $|\alpha_0| < 1 - \varepsilon$ at the origin and with $|\alpha| > \varepsilon$ on $\mathbb{C}$. Introduce the function $\mathfrak{t}$ on $[0, \infty)$ whose value at any given $s \in [0, \infty)$ is the average value of $-u$ on the circle in $\mathbb{C}$ of radius s. The equation $-\Delta u = 1 - e^{2u}$ implies the equation $s\partial_s \mathfrak{t} = \mathfrak{h}$ where $\mathfrak{h}(s)$ is the integral of $1 - |\alpha|^2$ over the radius s disk in $\mathbb{C}$ centered at the origin. The fact that $1 - |\alpha|^2 < \varepsilon$ at the origin implies that $\mathfrak{h} \ge c_0^{-1}\varepsilon$ on $[1, \infty)$ and so $s\partial_s \mathfrak{t} \ge c_0^{-1}\varepsilon$ on $[1, \infty)$. This being the case, then $\mathfrak{t} \ge c_0\varepsilon(\ln s) - c_0^2$. On the other hand, $\mathfrak{t} \le |\ln \varepsilon|$ if $|\alpha_0| > \varepsilon$ and this bound is violated when $\ln s \ge c_0^{-1}\varepsilon^{-1}|\ln \varepsilon| + c_0^2$.

The vortex equations enter the Proposition 2.4 proof via the upcoming Lemma 2.9. The lemma refers to a *transverse disk* with a given radius through a given point in Y. Such a disk is the image via the metric's exponential map of the centered disk of the given radius in the 2-plane bundle kernel($\hat{a}$) at the given point. There exists $c_0 > 100$ such that any transverse disk with radius $c_0^{-1}$ is embedded with a priori bounds on the derivatives to any given order of its extrinsic curvature. In addition, the vector field $v$ along $D_0$ is everywhere $c_0^{-1}$ close to the normal vector to $D_0$. All transverse disks are assumed implicitly to have radius less than $c_0^{-1}$ so as to invoke these two properties.

Lemma 2.9 uses J to view kernel($\hat{a}$) as a complex line bundle and it uses the Riemannian metric to define a compatible Hermitian structure on kernel($\hat{a}$). Fix $p \in Y$ and an isometric isomorphism from kernel$(\hat{a})|_p$ to $\mathbb{C}$. Use $\varphi$ in what follows to denote the map from $\mathbb{C}$ to Y that is obtained by composing first the isomorphism with kernel$(\hat{a})|_p$ and then the metric's exponential map. With $r \ge 1$ given, Lemma 2.6 uses $\varphi_r$ to denote the composition of first multiplication by $r^{-1/2}$ on $\mathbb{C}$ and then $\varphi$.

To finish the notational preliminaries, suppose next that $(A, \psi = (\alpha, \beta))$ is a given solution to some $r \ge 1$ and $\mu \in \Omega$ version of (1.13). Use $(A_r, \psi_r)$ to denote $\varphi_r^*(A, \psi)$. Lemma 2.9 writes $\psi_r$ as $(\alpha_r, \beta_r)$.

**Lemma 2.9**: *Fix an integer* $k \ge 1$ *and there exists* $\kappa > 1$ *with the following properties: Fix* $r \ge \kappa$ *and* $\mu \in \Omega$ *with* $\mathcal{P}$*-norm bounded by* 1 *and suppose that* $(A, \psi)$ *is a solution to the corresponding* $(r, \mu)$ *version of (1.13). Fix a point in* Y *to obtain the pair* $(A_r, \alpha_r)$ *of* *connection on and section of a complex line bundle over $\mathbb{C}$. There exists a solution to the vortex equation on $\mathbb{C}$ whose restriction to the radius* k *disk about the origin in $\mathbb{C}$ has* $C^k$-*distance less than* $\frac{1}{k}$ *from* $(A_r, \alpha_r)$ *on this same disk.*

***Proof of Lemma 2.9***: The argument is essentially identical to that use to prove Lemma 6.1 in [T1].

## c) Proof of Proposition 2.4

The proof of the proposition has seven parts. By way of a look ahead, the arguments are much like those in Section 6d of [T1].

*Part 1*: Let $D_0$ denote an embedded disk in Y with the following properties: First, the disk has radius $\rho > c_0 r^{-1/2}$. Second, all points in the disk have distance at least $(c_0 + 10^8)\rho$ from $\cup_{p\in\Lambda}\{\hat{\gamma}_p^+, \hat{\gamma}_p^-\}$. Lie transport by $v$ moves $D_0$ to a new disk. For $t \in \mathbb{R}$, use $D_t$ to denote the new disk that is obtained by moving the points in $D_0$ a distance t along the integral curves of $v$. The formula for $v$ in (1.3) can be used to see that $t = t_1$ and $t = t_2$ versions of $D_t$ are disjoint unless $t_1 = t_2$.

Fix a compactly supported function on $D_0$ which is equal to 1 on the radius $\frac{1}{2}\rho$ concentric subdisk in $D_0$ and with derivative bounded by $c_0\rho^{-1}$. Use $\chi_0$ to denote this function and use $\chi_t$ to denote the Lie transport of $\chi_0$ by $v$.

*Part 2*: Fix $\kappa_0 \geq 1$. Fix $r \geq c_0$ and $\mu \in \Omega$ with $\mathcal{P}$-norm bounded by 1 and suppose that $(A, \psi = (\alpha, \beta))$ is a solution to the $(r, \mu)$ version of (1.13). Let $D_0$ denote a disk as described in Part 1 with $1 - |\alpha|^2 \geq \kappa_0^{-1}$ at the center point of $D_0$.

Use $\mathfrak{f}$ to denote the function on $[0, \infty)$ that is given by the rule

$$t \to \mathfrak{f}(t) = r\int_{D_t} \chi_t(1 - |\alpha|^2) \ . \tag{2.9}$$

Note that $\mathfrak{f}(0) \geq c_0^{-1}$. This lower bound follows from the upper bound on $|\nabla\alpha|$ given in Lemma 2.1 with the fact that $1 - |\alpha|^2 \geq \kappa_0^{-1}$ at the center point of $D_0$.

The derivative of $\mathfrak{f}$ is denoted in what follows by $\mathfrak{f}'$; it is given by

$$\mathfrak{f}' = r\int_{D_t} \chi_t(\bar{\alpha}(\nabla_A\alpha)_v + (\nabla_A\bar{\alpha})_v\alpha)\,, \tag{2.10}$$

where $(\nabla\alpha)_v$ is used here and subsequently to denote the section of E that is obtained by pairing $\nabla\alpha$ with the vector field $v$. As explained in the next paragraph, the norm of the derivative of $\mathfrak{f}$ is such that

$$|\mathfrak{f}'| \le c_0 r \int_{D_t} (|\partial_t\chi||\beta|^2 + (|v^\perp| + |d^\perp\chi_t|)|\alpha||\beta|) \tag{2.11}$$

where the notation uses $v^\perp$ and $d^\perp\chi_t$ to denote the orthogonal projections of $v$ and $d\chi_t$ to the respective tangent and contangent bundles of $D_t$. Granted (2.11), use Lemma 2.1 to see that

$$|\mathfrak{f}(t) - \mathfrak{f}(0)| \le c_0 t . \tag{2.12}$$

This last inequality implies that

$$\mathfrak{f}(t) \ge c_0^{-1}\kappa_0^{-1} \quad \textit{when} \quad |t| \le c_0^{-2}\kappa_0^{-1}. \tag{2.13}$$

To prove (2.11), note first that J defines an almost complex structure on the kernel of $\hat{a}$. Equation (1.13) identifies $(\nabla_A\alpha)_v$ with a constant multiple of the part of $\nabla_A\beta$ that comes from the (1,0) part of the dual to the kernel of $\hat{a}$. Meanwhile, it identifies $(\nabla_A\beta)_v$ with a constant multiple of the part of $\nabla_A\alpha$ that comes from the dual to the (0, 1) part of the kernel of $\hat{a}$. Equation (2.11) follows from these observations with an integration by parts. By way of a warning, these same identifications are used later in the proof without further comment.

*Part 3*: This part constitutes a digression of sorts to draw attention to some consequences of the bounds given by Lemma 2.1. The remarks that follow here concern the integral over disks in Y of the curvature of the connection $\hat{A}$ given in (1.16) and the curvatures of analogs of $\hat{A}$ that are defined using (1.16) with a different version of the function $\wp$. In particular, allow in (1.16) any function $\wp$ on $[0, \infty)$ that is non-decreasing and such that $\wp(x) \le c_0 x$ for x near 0.

Given $r \ge c_0$ and $\mu \in \Omega$ with $\mathcal{P}$-norm bounded by 1, let $(A, \psi = (\alpha, \beta))$ denote a solution to the $(r, \mu)$ version of (1.13). Use the pair $(A,\alpha)$ to define $\hat{A}$. The corresponding curvature 2-form is denoted by $F_{\hat{A}}$; it is given by the formula

$$F_{\hat{A}} = (1 - \wp)\, F_A - \wp'(\nabla_A\bar{\alpha} \wedge \nabla_A\alpha) , \tag{2.14}$$

where $F_A$ here denotes the curvature 2-form of the connection A. In the context at hand, $F_A = *B_A$. The following uses what is said in the last paragraph of Part 2 with the bounds provided by Lemma 2.1 to see that

$$F_{\hat{A}} = ((1-\wp)\, r\,(1-|\alpha|^2) + \wp'\, |\nabla_A\alpha|^2 + \mathfrak{e}_v)\, w + \hat{a} \wedge \mathfrak{e}^{\perp} \tag{2.15}$$

where $|\mathfrak{e}_v| \le c_0((1-\wp)+\wp')\,(|1-|\alpha|^2| + r^{-1})$ and $|\mathfrak{e}^{\perp}| \le c_0((1-\wp)+\wp')\,(r^{1/2}|1-|\alpha|^2|^{1/2}+1)$. This depiction of $F_{\hat{A}}$ is plays an important role in subsequent arguments.

An additional fact about Â is used extensively, this concerning the case when $\wp$ is chosen to equal 1 on a neighborhood of $[1, \infty)$ in $[0, \infty)$: If $\kappa_0 > 1$ is such that $\wp = 1$ on $[1-\kappa_0^{-1}, \infty)$, then Â is flat and $\alpha|\alpha|^{-1}$ is Â-covariantly constant where $1-|\alpha|^2 < \kappa_0^{-1}$.

*Part 4*: Fix $\kappa_0 > 4$ and a function $\wp$ on $[0, \infty)$ that is zero on $[0, 1-2\kappa_0^{-1}]$ and is equal to 1 on $[1-\kappa_0^{-1}, \infty)$. Use this version of $\wp$ to define the connection Â. The 2-form $\frac{i}{2\pi}F_{\hat{A}}$ represents the first Chern class of E in the De Rham cohomology of Y and so it has integral G on the $f \in [1+\delta, 2-\delta]$ level sets in $M_\delta$. It follows as a consequence that there are points on any such surface where $1-|\alpha|^2 > \kappa_0^{-1}$. Meanwhile, it has integral zero on any level set of $f$ with $f$ not in this range, and it has integral zero on the u = constant 2-spheres in $\mathcal{H}_0$. This last observation implies that $1-|\alpha|^2$ must be $\mathcal{O}(r^{-1})$ on much of Y. The next lemma describes this region.

To set the stage for the lemma, fix $q \ge 1$ and let $\mathcal{Y}_q$ denote the set of points in Y with the following property: A point is in $\mathcal{Y}_q$ if it lies on a segment of an integral curve of $v$ with length $q$ or less and with one end point in $\mathcal{H}_0$. Note that $\mathcal{Y}_q$ contains $\mathcal{H}_0$, and it contains both the $f \le 1$ and $f \ge 2$ parts of $M_\delta$ if $q > c_0$. If $q > c_0$, then it also contains a small radius tubular neighborhood of the integral curve segments of $v$ in $M_\delta$ that start on the boundary of the radius $\delta$ coordinate ball about an index 1 critical point of $f$ and end on the boundary of the radius $\delta$ coordinate ball about an index 2 critical point of $f$. It also contains much of the $1-3\cos^2\theta \le 0$ portion of any given $\mathfrak{p} \in \Lambda$ version of $\mathcal{H}_{\mathfrak{p}}$, the missing part is a small radius tubular neighborhood of $\hat{\gamma}_{\mathfrak{p}}^{+} \cup \hat{\gamma}_{\mathfrak{p}}^{-}$.

To say more about these last parts of $\mathcal{Y}_q$, fix $\varepsilon > 0$ and fix $\mathfrak{p} \in \Lambda$. Let $\mathcal{H}_{\mathfrak{p},\varepsilon} \subset \mathcal{H}_{\mathfrak{p}}$ denote the subset of points with distance greater than $\varepsilon$ from $\hat{\gamma}_{\mathfrak{p}}^{+} \cup \hat{\gamma}_{\mathfrak{p}}^{-}$ and with $(u, \theta)$ coordinates such that either $1-3\cos^2\theta \le 0$ or $f(u)\,|\cos\theta|\sin^2\theta > \frac{2}{3\sqrt{3}}(x_0+4e^{-2R})+\varepsilon$. By way of a reminder, the function f is given in (1.2). Lemma II.2.2 finds $q_\varepsilon > 1$ such that each point in $\mathcal{H}_{\mathfrak{p},\varepsilon}$ has distance $q_\varepsilon$ or less along an integral curve of $v$ from $\mathcal{H}_0$. For example, $\mathcal{H}_{\mathfrak{p}} \cap M_\delta$ is the part of $\mathcal{H}_{\mathfrak{p}}$ where $|u| > R + \ln\delta$ and so a given point in $\mathcal{H}_{\mathfrak{p}} \cap M_\delta$ is in $\mathcal{H}_{\mathfrak{p},\varepsilon}$ unless both $1-3\cos^2\theta > 0$ and $|\cos\theta| < c_0\delta^{-2}(x_0+\varepsilon)$. This has the following consequence when $\varepsilon \le x_0$. The complement of the radius $c_0 x_0\delta^{-2}$ tubular neighborhood of the $M_\delta$ part of the union of the ascending disks from the index 1 critical points of $f$ and the descending disks from the index 2 critical points of $f$ is in $\mathcal{Y}_q$ if $q > q_\varepsilon$.

**Lemma 2.10**: *Fix* $q \geq 1$ *and there exists* $\kappa > 1$ *with the following significance: Suppose that* $r \geq \kappa$ *and that* $\mu \in \Omega$ *has* $\mathcal{P}$*-norm less than one. Let* $(A, \psi = (\alpha, \beta))$ *denote a solution to the* $(r, \mu)$ *version of (1.13). Then* $1 - |\alpha|^2 \leq \kappa r^{-1}$ *at all points in* $\mathcal{Y}_q$.

The proof is given momentarily. The lemma that follows directly plays a central role in the proof of Lemma 2.10.

**Lemma 2.11**: *There exists* $\kappa > 1$ *with the following significance: Fix* $r \geq \kappa$ *and* $\mu \in \Omega$ *with* $\mathcal{P}$*-norm less that one. Let* $(A, \psi = (\alpha, \beta))$ *denote a solution to the* $(r, \mu)$ *version of (1.13). Then* $1 - |\alpha|^2 \leq \kappa r^{-1}$ *on* $\mathcal{H}_0$ *and on the part of* $M_\delta$ *where* $f \leq 1 - 2\delta^2 - \kappa(\ln r) r^{-1/2}$ *or where* $f \geq 2 + 2\delta^2 + \kappa(\ln r) r^{-1/2}$.

The proof Lemma 2.11 and subsequent parts of the proof of Proposition 2.4 use $\kappa_*$ to denote the constant that appears in Lemma 2.3.

***Proof of Lemma 2.11***: To prove the first assertion, let S denote either a constant u sphere in $\mathcal{H}_0$ or a compact, level set of $f$ in $M_\delta$ with $f$ either less than $1 - 2\delta^2$ or greater than $2 + 2\delta^2$. Suppose that $p \in S$ is a point where $1 - |\alpha|^2 > \frac{1}{4} \kappa_*^{-1}$. It follows from Lemma 2.1's bound on $|\nabla_A \alpha|$ that the integral of $\frac{i}{2\pi} F_A$ over the disk in S centered at this point with radius $r^{-1/2}$ is greater than $c_0^{-1}\kappa_*^{-1}$. This understood, use the formula for $B_A$ in (1.13) with the bounds on $|\beta|$ supplied by Lemma 2.1 to see that the integral of $\frac{i}{2\pi} F_A$ over S is larger than $c_0^{-1}\kappa_*^{-1}$. But this is impossible given that the 2-form $\frac{i}{2\pi} F_A$ represents the first Chern class of E. Granted this conclusion, use Lemma 2.3 to conclude that $1 - |\alpha|^2 \leq c_0 r^{-1}$ on $\mathcal{H}_0$ and on the parts of $M_\delta$ where $f \leq 1 - 2\delta^2 - c_0(\ln r) r^{-1/2}$ and where $f \geq 2 + 2\delta^2 + c_0(\ln r) r^{-1/2}$.

***Proof of Lemma 2.10***: Fix $z > 1$ to be specified shortly. It is enough to consider the cases when $q = n z^{-1}$ with $n \in \{0, 1, 2, \ldots\}$. Since Lemma 2.11 gives the case for $n = 0$, an induction argument will prove the lemma for the general case for a suitable $z = c_0$. This understood, suppose that the lemma holds for a given integer $n \geq 0$ and suppose for the sake of argument that there is a point in Y where $1 - |\alpha|^2 \geq \frac{1}{2} \kappa_*$ and with distance less than $(n+2) z$ along a segment of an integral curve of $v$ with an endpoint from the part of Y that is described in Lemma 2.10. Let $D_0$ denote the transverse disk centered at this point with radius $c_0 r^{-1/2}$. The function $\mathfrak{f}$ given in Part 2 is such that $\mathfrak{f}(0) \geq c_0^{-1}\kappa_*^{-1}$. It follows from (2.11) that the function $\mathfrak{f}(t)$ is greater than $\frac{1}{4} c_0^{-1}\kappa_*$ if $|t| \leq c_0^{-1}$. If $z < c_0^{-1}$, this last conclusion violates the induction hypothesis that $1 - |\alpha|^2 \le c_0 r^{-1}$ on integral curve segments of length nz or less with one end point in the set described by Lemma 2.11. Thus, $1 - |\alpha|^2$ must be less than $\frac{1}{2}\kappa_*$ at all points along an integral curve segment of length less than $(n+2)z^{-1}$ with one end point in this same set from Lemma 2.11. Granted this fact, then what is said in Lemma 2.3 completes the proof.

*Part 5*: Let $S \subset M_\delta$ denote a level set of $f$ with $f \in (1+\delta, 2-\delta)$. As noted in the first paragraph of Part 4, there are points on S where $1 - |\alpha|^2$ is greater than $\frac{1}{2}\kappa_*^{-1}$. Let p denote such a point. It follows from Lemmas 2.8 and 2.9 that there is a point in S with distance at most $c_0 r^{-1/2}$ from p where $|\alpha| < \frac{1}{100}$. In fact, there is a point in S with distance less than $c_0 r^{-1/2}$ from p where $\alpha = 0$. To prove the existence of such an $\alpha = 0$ point, define $\hat{A}$ as in the first paragraph of Part 4 using $\kappa_0 = 2\kappa_*$ and a suitable function $\wp$. It follows from Lemmas 2.8 and 2.9 that a disk of radius $c_0 r^{-1/2}$ centered on any point in S where $\hat{A}$ is not flat must contribute at least $c_0^{-1}$ to the integral of $\frac{i}{2\pi}F_{\hat{A}}$ over S. Use this last observation with (2.15) to see that there can be at most $c_0$ pairwise disjoint disks in S where $\hat{A}$ is not flat. It follows $1 - |\alpha|^2 < \frac{1}{2}\kappa_*^{-1}$ on the boundary of the disk centered at p with radius at most $c_0 r^{-1/2}$. In particular, the connection $\hat{A}$ is flat near the boundary of the latter disk and $\alpha/|\alpha|$ is $\hat{A}$-covariantly constant. As the integral of $\frac{i}{2\pi}F_{\hat{A}}$ over this disk is non-zero, it is a positive integer. These last two facts require a zero of $\alpha$ in this disk because the integral of $\frac{i}{2\pi}F_{\hat{A}}$ over the disk is the sum of the local Euler numbers at the zeros of $\alpha$ in the disk. What follows summarizes this. There exists $z \in [1, c_0]$ and at most G pairwise disjoint disks in S of radius $z r^{-1/2}$ with the following properties:

- $1 - |\alpha|^2$ *is less than* $\frac{1}{2}\kappa_*^{-1}$ *on the complement of the union of these disks.*
- *The integral of* $\frac{i}{2\pi}F_{\hat{A}}$ *over each disk is a positive integer; and this integer for any given disk is is the sum of the local Euler numbers of the zeros of* $\alpha$ *in the disk.*

(2.16)

Use the fact that there are at most G such disks to see there is a set of at most G disks of radius at most $(z + c_0)r^{-1/2}$ such that each disk in the set obeys (2.16) and such that the distance between pairwise distinct disks from the set is greater than $c_0^{-1}$. Let N denote the number of elements in this set. Enumerate this set of disks as $\{D^{(i)}_S\}_{1\le i\le N}$. For each index $i \in \{1, \ldots, N\}$ and for each $t \in \mathbb{R}$, use $D^{(i)}_{S,t}$ to denote the disk in Y that is obtained from $D^{(i)}_S$ by moving its points for time t along the integral curves of $v$. Let $t_S$ denote the value of $f$ on S. If $t_S + t \in (1+\delta, 2-\delta)$, then each $D^{(i)}_{S,t}$ is a disk in the $t_S + t$ level set of $f$. It follows from (2.15) and from the comment in the final paragraph of Part 2 that there exists $c_0 \in (1, c_0)$ with the following property: If $t_S + t \in (1+\delta, 2-\delta)$ and if $|t| < c_0^{-1}$, then $1 - |\alpha|^2 < \kappa_*^{-1}$ on the complement of $\cup_{1\le i\le N} D^{(i)}{}_{S,t}$ and the integral of $\frac{i}{2\pi} F_{\hat{A}}$ over any given $D^{(i)}{}_t$ is the same as its integral over $D^{(i)}{}_S$. Meanwhile, the diameter of $D^{(i)}{}_{S,t}$ is bounded by $c_0 r^{-1/2}$ and the pairwise separation between $D^{(i)}{}_{S,t}$ and $D^{(i')}{}_{S,t}$ when $i \ne i'$ is at least $c_0^{-1}$.

Granted the preceding observations, let t be such that $t_S + t \in [1+\delta, 2-\delta]$ and such that $|t| < c_0^{-1}$. Let S´ denote the $t_S + t$ level set. Define $N'$ and $\{D^{(i)}{}_{S'}\}_{1\le i\le N'}$ as done above for the case of S. It follows from what was said in the preceding paragraph that $N' = N$ and that the set $\{D^{(i)}{}_{S'}\}_{1\le i\le N}$ can be labeled so that any given $D^{(i)}{}_{S'}$ has non-empty intersection with $D^{(i)}{}_{S,t}$ and has distance at most $c_0^{-1}$ from any $i' \ne i$ version of $D^{(i)}{}_{S,t}$. If $r \ge c_0$, then these facts when applied sequentially to some $2n \le c_0$ level sets of $f$ starting with $\Sigma$, and then with $f$ value $\frac{3}{2} \pm c_0^{-1}$, then with $f$ value $\frac{3}{2} \pm 2c_0^{-1}$, and so on lead to the following:

- *The integer* N *is independent of the value of* $f \in [1+\delta, 2-\delta]$.
- *There is a set of* N *segments of integral curves of* $v$ *with the following properties:*
  a) *Each segment starts on the* $f = 1+\delta$ *and ends on the* $f = 2-\delta$ *level set.*
  b) *Distinct segments from the set have pairwise separation no less than* $c_0^{-1}$.
  c) *The intersection of each segment with a level set of f has distance at most* $c_0 r^{-1/2}$ *from a zero of* $\alpha$.
  d) $1 - |\alpha|^2 < \frac{1}{2}\kappa_*$ *at the points in the* $f \in (1+\delta, 2-\delta)$ *part of* $M_\delta \cup \mathcal{H}_0$ *with distance greater than* $c_0^2 r^{-1/2}$ *from the union of the segments in this set.*

(2.17)

Part 6 explains why N = G and it says more about the start and end points of (2.17)'s integral curve segments.

*Part 6*: Fix $\mathfrak{p} \in \Lambda$ and $\varepsilon > 0$ so as to reintroduce from Part 4 the subset $\mathcal{H}_{p,\varepsilon} \subset \mathcal{H}_p$. By way of a reminder, this is the subset of points with distance greater than $\varepsilon$ from $\hat{\gamma}_p^+ \cup \hat{\gamma}_p^-$ and such that either $1 - 3\cos^2\theta \le 0$ or $f(u)\,|\cos\theta|\sin^2\theta > \frac{2}{3\sqrt{3}}(x_0 + 4e^{-2R}) + \varepsilon$. Lemma 2.9 finds $c_\varepsilon > 1$ such that if $r > c_\varepsilon$, then $1 - |\alpha|^2 \le c_0 r^{-1/2}$ on $\cup_{p\in\Lambda}\mathcal{H}_{p,\varepsilon}$. Fix $\varepsilon \le x_0$ and assume henceforth that $r > c_\varepsilon$.

Given what was just said, the segments of integral curves of $v$ that arise in (2.17) intersect the $f \le 1 + \delta_*^2$ part of $M_\delta$ in the union of the radius $\delta_*$ coordinate balls about the index 1 critical points of $f$ and they likewise intersect the $f > 2 - \delta_*^2$ part of $M_\delta$ in the union of the radius $\delta_*^2$ coordinate balls about the index 2 critical points of $f$. Moreover, each such intersection lies where $1 - 3\cos^2\theta > 0$ and $|\cos\theta| \le c_0 x_0 \delta^{-2}$. This understood, Lemma 2.10 implies that $1 - |\alpha|^2 < c_0 r^{-1}$ on the $|\cos\theta| > c_0 x_0 \delta^{-2}$ portion of the $|u| \ge R + \ln\delta$ part of $\mathcal{H}_p$.

Define $\hat{A}$ as in Part 5. As the integral of $\frac{i}{2\pi} F_{\hat{A}}$ over any given constant u sphere in $\mathcal{H}_p$ is equal to 1, and as $\hat{A}$ is flat where $1 - |\alpha|^2 < \frac{1}{2}\kappa_*$, it follows from what was just said that the radius $\delta_*$ coordinate ball about any given index 1 critical point of $f$ must contain the starting point of one of (2.17)'s integral curve segments. By the same token, the radius $\delta_*$ coordinate ball about any give index 2 critical point of $f$ must contain the ending point of one of (2.17)'s integral curve segments. This can happen only if $\text{N} = \text{G}$.

Granted that $\text{N} = \text{G}$, then the following must also hold: Let D denote an embedded disk in $M_\delta$ that intersects just one of (2.17)'s integral curve segments and is such that its boundary has distance $c_0 r^{-1/2}$ or greater from all of (2.17)'s integral curve segments. Then the integral of $\frac{i}{2\pi} F_{\hat{A}}$ over D is equal to 1.

*Part 7*: What is said by Lemma 2.10, by (2.17) and in Part 6 verify all of Proposition 2.4 on $M_\delta \cup \mathcal{H}_0$ but for the assertion that $\alpha^{-1}(0)$ is transversal in $Y - T_\varepsilon$ with G components and with tangent line very close to $v$. As explained next, these as yet unproved assertions are all consequences of Lemmas 2.9 and 2.10.

To see about a proof, fix $R > c_0$ for the moment and let $D \subset M_\delta \cup \mathcal{H}_0$ denote a smoothly embedded, transverse disk of radius $2R r^{-1/2}$ and center on one of the (2.15)'s integral curve segments. Assume that $v$ is normal to D at its center point. Define $(A_r, \alpha_r)$ as in Lemma 2.10 and fix $\varepsilon' > 0$. According to Lemma 2.10, there is a constant $r_{\varepsilon'}$ that depends only on $\varepsilon'$ and is such that if $r \ge r_{\varepsilon'}$, then there is a solution to (2.8) with $C^2$ distance at most $\varepsilon'$ from $(A_r, \alpha_r)$ on the disk of radius R centered at the origin in $\mathbb{C}$. Let $(A_0, \alpha_0)$ denote such a solution.

To say more about $(A_0, \alpha_0)$, keep in mind that $\hat{A}$ is flat and $\alpha/|\alpha|$ is $\hat{A}$ covariantly constant along D at points with distance greater then some fixed multiple of $r^{-1/2}$ that is independent of R. Denote this multiple as $r$ and assume that $R >> r$. If $\varepsilon' < c_0^{-1}$, it then follows that $1 - |\alpha_0|^2 < \frac{3}{4} \kappa_*^{-1}$ on the radius $r$ disk about the origin in $\mathbb{C}$. It also follows that the $(A_0, \alpha_0)$ version of the connection $\hat{A}$ is flat on the annulus in $\mathbb{C}$ with inner radius $r$ and outer radius R, and is such that $\alpha_0/|\alpha_0|$ is $\hat{A}_0$-covariantly constant on this same annulus Moreover, the integral of the $(A_0, \alpha_0)$ version of $\frac{i}{2\pi} F_{\hat{A}_0}$ over any centered disk in $\mathbb{C}$ with radius between $r$ and R must equal 1. This understood, then $\alpha_0$ must have at least one zero in the centered, radius $r$ disk. Given that $\alpha_0$ is $\bar{\partial}_{A_0}$ -holomorphic, there can be at most one such zero and it is necessarily non-degenerate with local Euler number 1.

Use $z_0 \in \mathbb{C}$ to denote this zero. Given the aforementioned holomorphicity, $\alpha_0$ must appear near $z_0$ as $\alpha_0 = \varsigma(z - z_0) + \mathfrak{e}$ where $\varsigma \in \mathbb{C}$ and with $\mathfrak{e}$ such that $|\mathfrak{e}| \le \varsigma' |z - z_0|^2$. Note in this regard $|\varsigma| \ge c_0^{-1}$ and $|\varsigma'| < c_0$, this being a consequence of the fact that the equations in (2.8) are elliptic modulo the action $C^\infty(\mathbb{C}; S^1)$ and so any sequence of solutions has a subsequence that converges up to this group action in the $C^\infty$ Frechêt topology on compact subsets of $\mathbb{C}$. This same sequential compactness property implies that $|\alpha_0| \ge c_0^{-1} \frac{|z-z_0|}{1+|z-z_0|}$ in the radius R disk about $0 \in \mathbb{C}$.

These last facts about $\alpha_0$ have the implications that follow with regards to $\alpha$. First, if $\varepsilon' < c_0^{-1}$, then $\alpha$ has a single, transverse zero in D with distance at most $c_0\varepsilon' r^{-1/2}$ from D's center point. To give the second implication, use J as before to define the (1,0) and (0,1) parts of the complexification of the 2-plane bundle kernel($\hat{a}$). Introduce by way of notation $\partial_A\alpha$ to denote the (1,0) part of $\nabla_A\alpha$, this being the homomorphism from (1,0) part of this complexification to E that is defined by restricting the domain of $\nabla_A\alpha$. It must also be the case that $|\partial_A\alpha| \geq c_0 r^{1/2}$ at this zero of $\alpha$. Note in this regard that corresponding restriction $\bar{\partial}_A\alpha$ of $\nabla_A\alpha$ to the (0,1) part of kernel($\hat{a}$) $\otimes_{\mathbb{R}} \mathbb{C}$ is equal to the directional covariant derivative of $\beta$ along $v$ and so has norm bounded by $c_0$. By way of a reminder, the directional covariant derivative of $\alpha$ along $v$ was denoted by $(\nabla_A\alpha)_v$ and being a linear combination of covariant derivatives of $\beta$, it too has norm bounded by $c_0$.

What was just said as applied to transverse disks along the various components of (2.17)'s integral curve segments verifies the claim that $\alpha^{-1}(0)$ is transverse and it verifies the claim that each component of the radius $c_0 r^{-1/2}$ tubular neighborhood of (2.17)'s integral curve segments contains precisely one component. To see about the tangent line to a component, parametrize a neighborhood of a given point in a segment at unit speed by a map from a small interval about the origin in $\mathbb{R}$ to Y. Use x to denote this map. Let $\partial_t$ denote the Euclidean vector field on $\mathbb{R}$. Then $x_*\partial_t$ pairs with $\nabla_A\alpha$ to give zero. With this in mind, write $x_*\partial_t$ at the origin in $\mathbb{R}$ as $x_v v + x_{(1,0)} + x_{(0,1)}$ where $x_{(1,0)}$ is the projection of $x_*\partial_t$ to the (1,0) part of kernel($\hat{a}$) $\otimes_{\mathbb{R}} \mathbb{C}$, and $x_{(0,1)}$ is the complex conjugate of $x_{(1,0)}$. The fact that $\nabla_A\alpha$ annihilates $x_*\partial_t$ means that

$$x_v (\nabla_A\alpha)_v + x_{(1,0)}\partial_A\alpha + x_{(0,1)}\bar{\partial}_A\alpha = 0 \tag{2.18}$$

Given what is said in the preceding paragraph about the relative sizes of the various projections of the covariant derivative of $\alpha$, the equality in (2.18) can not hold unless

$$r^{1/2}|x_{(1,0)}| \leq c_0|x_v| \tag{2.19}$$

This last inequality implies the claim about the tangent vector to $\alpha^{-1}(0)$.

*Part 8*: What is said by Lemma 2.10 and by (2.17) with what is said in Parts 6 and 7 verify the assertions of Proposition 2.4 for the $M_\delta \cup \mathcal{H}_0$ part of Y. The upcoming Lemma 2.12 is used momentarily to extend the domain where these assertions hold into $\cup_{\mathfrak{p}\in\Lambda}\mathcal{H}_\mathfrak{p}$. To set the stage for this lemma, fix $\mathfrak{p} \in \Lambda$. Given $\mathrm{D} > 0$, use $H_{\mathfrak{p},\mathrm{D}}$ to denote the set of points in $\mathcal{H}_\mathfrak{p}$ with distance at least D from all points in $\hat{\gamma}_\mathfrak{p}^+ \cup \hat{\gamma}_\mathfrak{p}^-$.

**Lemma 2.12**: *There exists* $\kappa > 1$ *with the following significance: Fix* $r \ge \kappa$ *and fix* $\mu \in \Omega$ *with* $\mathcal{P}$*-norm less that one. Let* $(A, \psi = (\alpha, \beta))$ *denote a solution to the* $(r, \mu)$ *version of (1.13). If* n *is a given positive integer, set* $D(n) = (1 + \kappa^{-1})^n \kappa r^{-1/2}$. *Assume that* n *is such that the assertions of Proposition 2.4 hold in* $M_\delta \cup \mathcal{H}_0 \cup (\cup_{p\in\Lambda} \mathcal{H}_{p,D(n+1)})$. *Then the assertions of Proposition 2.4 also hold in* $M_\delta \cup \mathcal{H}_0 \cup (\cup_{p\in\Lambda} \mathcal{H}_{p,D(n)})$

This lemma is proved momentarily.

To finish the proof of Proposition 2.4, introduce $\kappa$ from Lemma 2.12 and let N denote the least integer such that $\mathcal{H}_{p,D(N+1)} \subset \mathcal{H}_p \cap M_\delta$ for all $p \in \Lambda$. The assertions of Proposition 2.4 have been verified on $M_\delta \cup \mathcal{H}_0 \cup (\cup_{p\in\Lambda} \mathcal{H}_{p,D(N+1)})$. This understood, invoke Lemma 2.12 a total of N times to prove sequentially that the assertions of Proposition 2.4 hold on $M_\delta \cup \mathcal{H}_0 \cup (\cup_{p\in\Lambda} \mathcal{H}_{p,D(N)})$, then on $M_\delta \cup \mathcal{H}_0 \cup (\cup_{p\in\Lambda} \mathcal{H}_{p,D(N-1)})$, etc.

***Proof of Lemma 2.12***: Fix $L > 1$ and set $D(n) = (1 + L^{-1})^n L r^{-1/2}$. Suppose that the assertions of Proposition 2.4 hold on $M_\delta \cup \mathcal{H}_0 \cup (\cup_{p\in\Lambda} \mathcal{H}_{p,D(n+1)})$. The proof that they also hold on $M_\delta \cup \mathcal{H}_0 \cup (\cup_{p\in\Lambda} \mathcal{H}_{p,D(n)})$ for a suitable, r-independent choice of L is given in seven steps.

<u>Step 1</u>: Only the third bullet of Proposition 2.4 needs comment when $1 - |\alpha|^2 < \kappa_*^{-1}$ on $\mathcal{H}_{p,D(n)} - \mathcal{H}_{p,D(n+1)}$. In any event, the third bullet restates part of Lemma 2.3 and so it holds whether or not $1 - |\alpha|^2 < \kappa_*^{-1}$ on the whole of $\mathcal{H}_{p,D(n)} - \mathcal{H}_{p,D(n+1)}$.

<u>Step 2</u>: Fix $p \in \Lambda$. Given $c \ge c_0$, suppose that $D_0$ denotes an embedded disk in $\mathcal{H}_p$ with radius $c r^{-1/2}$ whose points have distance at least $(c_0 + 10^8) c r^{-1/2}$ from both $\hat{\gamma}_p^+$ and $\hat{\gamma}_p^-$. Assume in addition that the vector field $v$ along $D_0$ is at all points $c_0^{-1}$ close to the normal vector. For example, a transverse disk has this last property.

Use $d_*(t)$ to denote the distance from $\hat{\gamma}_p^+ \cup \hat{\gamma}_p^-$ to the point at time t along the integral curve of $v$ that starts at the center point of $D_0$. It follows from (1.3)'s formula for $v$ (see Equation (II.2.3)) that there exists $\lambda \ge c_0^{-1}$ with the following property: Either one or both of the inequalities $d_*(t) \ge d_*(0) e^{\lambda t}$ and $d_*(-t) \ge d_*(0) e^{\lambda |t|}$ hold if t is such that the point relevant point at time t on the integral curve is in $\mathcal{H}_p$. The discussion that follows assumes the first of these conditions and t is assumed implicitly to be non-negative.

The rest of this step contains observations on the geometry of $D_0$ and the $t > 0$ versions of $D_t$. Assume for all of these that the center point of $D_t$ is in $\mathcal{H}_p$. Granted this assumption, then $d_*(t)$ can serve as a proxy of sorts for the distance between any given point in $D_t$ and $\hat{\gamma}_p^+ \cup \hat{\gamma}_p^-$. In particular, all points in $D_t$ have distance at least $(1 - c_0^{-1}) d_*(t)$ from $\hat{\gamma}_p^+ \cup \hat{\gamma}_p^-$ and distance at most $(1 + c_0^{-1}) d_*(t)$ from $\hat{\gamma}_p^+ \cup \hat{\gamma}_p^-$ because the points in $D_0$ have distance at least $10^8 c r^{-1/2}$ from $\hat{\gamma}_p^+ \cup \hat{\gamma}_p^-$ and at most $c r^{-1/2}$ from each other.

It is also the case that the image in $D_t$ of two points in $D_0$ with a given distance $\rho$ from each other have distance at most $\rho\, e^{\lambda t}$ apart in $D_t$. This is to say that the disk $D_t$ is not seriously distorted if $t < c_0^{-1}$. Moreover, if $t \le c_0^{-1}$, then the normal vector to $D_t$ at all points will be close to $v$.

There is one other point to keep in mind about $D_0$, this concerning the number of intersection points between $D_0$ and a given segment in $\mathcal{H}_p$ of an integral curve of $v$. In particular, $D_0$ has at most one intersection with any such section. This is proved using the formula for $v$ in (1.3) given the assumption that $d_*(0)$ is at least $c_0$ times $D_0$'s diameter.

Step 3: Assume in this step that the function $\mathfrak{f}$ from (2.7) is such that $\mathfrak{f}(0)$ is greater than $c^{1/3}$. Given that $d_*(t) \ge d_*(0)\, e^{\lambda t}$, so $d_*(t) \ge (1 + L^{-1})^3 d_*(0)$ if $t \ge 3\lambda^{-1} \ln(1+L^{-1})$. Set $t_* = 100\, \lambda^{-1} \ln(1+L^{-1})$ and use (2.6) to see that

$$\mathfrak{f}(t_*) \ge 10^{-2}\, c^{1/3} - c_0 \ln(1+L^{-1}). \tag{2.20}$$

Suppose that $D_0 \subset \mathcal{H}_{p,D(n-1)}$ and that $c > c_0 (1 + \ln(1+L^{-1}))^3$. If such is the case, then the inequality in (2.20) is not compatible with the assumption that the assertions of Proposition 2.4 hold on $M_\delta \cup \mathcal{H}_0 \cup (\cup_{p\in\Lambda} \mathcal{H}_{p,D(n+1)})$. It follows as a consequence that $\mathfrak{f}(0)$ can be no greater than $c^{1/3}$ if $c > c_0 (1 + \ln(1+L^{-1}))^3$. Assume this bound for $c$ in what follows and likewise assume that $D_0 \subset \mathcal{H}_{p,D(n-1)}$ so as to guarantee that $\mathfrak{f}(0) < c^{1/3}$.

Step 4: Fix $R \ge 2$ but less than $c_0^{-1} c^{2/3}$. Since $\mathfrak{f}(0) < c^{1/3}$, the bounds from Lemma 2.1 for $|\nabla\alpha|$ requires a point $x \in [1, Rc^{1/3}]$ such that $1 - |\alpha|^2 \le c_0 R^{-1}$ on the concentric annulus in $D_0$ with inner radius $x\, r^{-1/2}$ and outer radius $(x+1)\, r^{-1/2}$. With this understood, use Lemma 2.3 to deduce the following: If $R \ge c_0\kappa_*$, then $1 - |\alpha|^2 \le \frac{1}{16}\kappa_*^{-1}$ on such an annulus. Assume in what follows that $c > c_0$ is such that R can be chosen greater than $c_0\kappa_*$. Reintroduce the connection $\hat{A}$ from Part 4 as defined with $\kappa_0 = 2\kappa_*$. This connection is flat and $\alpha|\alpha|^{-1}$ is $\hat{A}$-covariantly constant on this annuli in $D_0$.

Step 5: Take $D_0$ to be a transverse disk. With $\kappa_*$ denoting as before the constant from Lemma 2.3, assume in addition that $1 - |\alpha|^2 \ge \kappa_*^{-1}$ at the center point of $D_0$. This assumption with Lemmas 2.8 and 2.9 guarantee a point with distance at most $c_0 r^{-1/2}$ from the origin in $D_0$ where $1 - |\alpha|^2 > \frac{9}{10}$. If $c > c_0$, then the existence of such a point implies that the integral of $\frac{i}{2\pi} F_{\hat{A}}$ over the subdisk in $D_0$ of radius $(x+1) r^{-1/2}$ is non-zero and positive. Moreover, this integral must be a positive integer because $\hat{A}$ is flat with a covariantly constant section near the boundary of this subdisk. Use $n_1$ in what follows to denote this integer. Keep in mind that $\alpha$ has a zero in this subdisk with positive local

Euler number because the sum of the local Euler numbers of the zeros of $\alpha$ in the subdisk is equal to this same $n_1$.

Use (1.13) with Lemma 2.1's bound on $|\nabla_A\beta|$ to draw the following conclusion: There exists $t_0 \geq c_0^{-1}$ so that $1 - |\alpha|^2 < \frac{1}{2}\kappa_*^{-1}$ at all times $t \leq t_0$ on the image in $D_t$ of the annulus with inner radius $x\,r^{-1/2}$ and outer radius $(x+1)r^{-1/2}$. This being the case, the integral of $\frac{i}{2\pi}F_{\hat{A}_0}$ over the image in $D_t$ of the radius $(x+1)r^{-1/2}$ subdisk in $D_0$ is still equal to $n_1$ and $\alpha$ still has at least one zero with positive local Euler number in the image in $D_t$ of the radius $(x+1)r^{-1/2}$ subdisk of $D_0$.

With the preceding understood, remark that if $L \geq c_0$, then $t_0 > t_*$ with $t_*$ as defined in Step 3. Assume that this is the case. Then $D_{t_0} \subset \mathcal{H}_{p,D(n+1)}$ and as a consequence, $n_1$ must equal 1 because Proposition 2.4's assertions hold on $M_\delta \cup \mathcal{H}_0 \cup (\cup_{p\in\Lambda}\mathcal{H}_{p,D(n+1)})$.

Step 6: Let $x \in [1, R\mathcal{c}^{1/3}]$ denote the constant from Step 4. It follows from what is said in Step 3 that there is a zero of $\alpha$ in $D_0$ with local Euler number 1 with distance at most $c_0\, r^{-1/2}$ from the center of $D_0$. Use $D_0$ now to denote the transverse disk through this point with radius $(1 - c_0^{-1})\mathcal{c}r^{-1/2}$. The conclusions of the preceding steps hold for this new version of $D_0$ also. In particular, there exists $x \in [1, R\mathcal{c}^{1/3}]$ such that the connection $\hat{A}$ is flat with covariantly constant section $\alpha|\alpha|^{-1}$ on the concentric annulus with inner and outer radii $x\,r^{-1/2}$ and $(x+1)r^{-1/2}$. Let $D_\Diamond$ denote the concentric subdisk in $D_0$ with radius $(x+1)r^{-1/2}$. As before, the integral of $\frac{i}{2\pi}F_{\hat{A}}$ over $D_\Diamond$ is equal to 1. This value of 1 for the integral of $\frac{i}{2\pi}F_{\hat{A}}$ demands the following:

*There exists* $z \geq 1$ *that is independent of* $(A, \alpha)$, $\mu$ *and* $r$ *and such that if* $r \geq z$, *then* $1 - |\alpha|^2 \leq \frac{1}{32}\kappa_*^{-1}$ *on the part of* $D_\Diamond$ *with distance greater than* $z r^{-1/2}$ *from the origin.*

(2.21)

What follows explains why (2.21) is true. To start, let $p \in D_\Diamond$ denote a given point where $1 - |\alpha|^2 \geq \frac{1}{32}\kappa_*^{-1}$. Use Lemmas 2.8 and 2.9 to find a point $p'$ in the transverse disk of radius $c_0 r^{-1/2}$ through p where $|\alpha|$ is less than $10^{-5}$. Since $v$ is almost normal to this transverse disk at $p'$ and also to $D_\Diamond$, so there is a point in $D_\Diamond$ on the integral curve of $v$ through $p'$ with distance at most $c_0 r^{-1/2}$ from $p'$. Let $p''$ denote the latter point. As noted previously, the Dirac equation in (1.13) identifies the covariant derivative of $\alpha$ in the direction of $v$ with a linear combination of covariant derivatives of $\beta$. This understood, then it follows from Lemma 2.1 that $|\alpha|$ at $p''$ is no greater than $10^{-4}$. Given this upper bound for $|\alpha|$, use Lemma 2.1's bound upper bound for $|\nabla_A\alpha|$ to see that the contribution to the integral over $D_\Diamond$ of $\frac{i}{2\pi}F_{\hat{A}}$ from the radius $c_0 r^{-1/2}$ disk in $D_\Diamond$ centered at $p''$ is greater than $c_0^{-1}$. This conclusion has the following consequence: There can be at most $c_0$ points

in any subset of $D_\Diamond$ with the following property: The distance between any two distinct points is greater than $c_0 r^{-1/2}$ and $1 - |\alpha|^2 > \frac{1}{32}\kappa_*^{-1}$ at each point.

Now, suppose that p is a point in $D_\Diamond$ with $1 - |\alpha|^2 > \frac{1}{32}\kappa_*^{-1}$. It follows from what was just said and from Lemmas 2.1 and 2.3 that there exists an $(A, \psi)$, $\mu$ and r independent constant $z_1 \geq 10^4$ and a subdisk $D_p \subset D_\Diamond$ centered at p with radius $z_1 r^{-1/2}$ with the following properties: First, $1 - |\alpha|^2 < \frac{1}{32}\kappa_*^{-1}$ on the annular neighborhood of the boundary of $D_p$ with inner and outer radii equal to $\frac{1}{2} z_1 r^{-1/2}$ and $z_1 r^{-1/2}$. Second, the integral of $\frac{i}{2\pi}F_{\hat{A}}$ over $D_p$ is at least $c_0^{-1}$. This last property implies that the integral of $\frac{i}{2\pi}F_{\hat{A}}$ over $D_p$ is at least 1 since the connection $\hat{A}$ is flat with $\alpha/|\alpha|$ covariantly constant on the annular boundary neighborhood.

Since the integral of $\frac{i}{2\pi}F_{\hat{A}}$ over the whole of $D_\Diamond$ is equal to 1, the conclusions of the preceding paragraph have the following consequence: Any two versions of $D_p$ are certain to overlap. It follows that $1 - |\alpha|^2$ is less than $\frac{1}{32}\kappa_*^{-1}$ at any point in $D_\Diamond$ with distance $c_0 r^{-1/2}$ or greater from the center point since $\alpha$ is zero at this point. This last observation verifies Proposition 2.4's fourth bullet for $M_\delta \cup \mathcal{H}_0 \cup (\cup_{p\in\Lambda}\mathcal{H}_{p,D(n)})$.

Step 7: Suppose that $t \in [0, t_0]$. Use $D_{\Diamond t} \subset D_t$ to denote the image of $D_\Diamond$. The definition of $t_0$ is such as to guarantee that $1 - |\alpha|^2$ is no greater than $\frac{1}{2}\kappa_*$ on the image in $D_{\Diamond t}$ of the set of points in $D_\Diamond$ with distance at most $c_0 r^{-1/2}$ from the origin. This understood, the integral of $\frac{i}{2\pi}F_{\hat{A}}$ over $D_{\Diamond t}$ is equal to 1. Granted these last conclusions, an essentially verbatim repetition of what is said in Part 7 proves that the zero locus of $\alpha$ in the cylinder $\cup_{1\leq t\leq t_0} D_{\Diamond t}$ is transverse and consists of a properly embedded arc whose tangent vector has angle at most $c_0 r^{-1/2}$ from $v$ and whose points have distance at most $c_0 r^{-1/2}$ from the integral curve of $v$ through the center point in $D_\Diamond$. Since $D_{\Diamond t_0} \subset \mathcal{H}_{p,D(n+1)}$, this arc smoothly extends the zero locus of $\alpha$ in $\mathcal{H}_{p,D(n+1)}$. The preceding observations verifies Proposition 2.4's first bullet for $M_\delta \cup \mathcal{H}_0 \cup (\cup_{p\in\Lambda}\mathcal{H}_{p,D(n)})$.

## 3. The map $\Phi^r$ from $\mathcal{Z}_{ech,M}^L$ to $\mathcal{Z}_{SW,r}$

Fix $\mu \in \Omega$ with $\mathcal{P}$-norm less then 1 and fix $L \geq 1$. The map $\hat{\Phi}^r : \hat{\mathcal{Z}}_{ech,M}^L \to \hat{\mathcal{Z}}_{SW,r}$ for Theorem 1.5 is a principal $\mathbb{Z}$ bundle covering map over a map from $\mathcal{Z}_{ech,M}^L$ into $\mathcal{Z}_{SW,r}$ that is denoted in what follows by $\Phi^r$.

The following proposition makes a formal assertion as to the existence of the desired map $\Phi^r$. It then says something about the form of the solutions to the relevant version of (1.13) that lie in the $C^\infty(Y; S^1)$ orbits in $\text{Conn}(E) \times C^\infty(Y; \mathbb{S})$ that comprise $\Phi^r$'s image. The proposition uses the isomorphism in (1.19) to identify $\hat{\mathcal{Z}}_{SW,r}$ with $\mathcal{Z}_{SW,r} \times \mathbb{Z}$ and it uses the isomorphism described in the paragraph preceding Theorem 1.5 to identify $\hat{\mathcal{Z}}_{ech,M}$ with $\mathcal{Z}_{ech,M} \times \mathbb{Z}$ and thus $\hat{\mathcal{Z}}_{ech,M}{}^{L}$ with $\mathcal{Z}_{ech,M}{}^{L} \times \mathbb{Z}$. The proposition also uses the following notation: When $\gamma$ denotes a closed integral curve of $v$, then $\ell_\gamma$ denotes its length.

**Proposition 3.1**: *There exists* $\kappa > \pi$, *and given* $E > 1$ *and* $L > \kappa E$, *there exists* $\kappa_L > \kappa$; *and* $\kappa$ *with* $\kappa_L$ *have the following significance: Fix* $r \geq \kappa_L$ *and an element* $\mu \in \Omega$ *with* $\mathcal{P}$*-norm less that 1. Use the solutions to the* $(r,\mu)$ *version of (1.13) to define* $\mathcal{Z}_{SW,r}$. *There exists a 1-1 map* $\Phi^r$: $\mathcal{Z}_{ech,M}{}^{L} \to \mathcal{Z}_{SW,r}$ *whose image contains the subset of* $\mathcal{Z}_{SW,r}$ *with* $M < E$. *Moreover, suppose that* $\Theta \in \mathcal{Z}_{ech,M}{}^{L}$ *and that* $\mathfrak{c} = (A, \psi = (\alpha, \beta))$ *is a solution to the* $(r, \mu)$ *version of (1.13) on the* $C^\infty(Y; S^1)$ *orbit defined by* $\Phi^r(\Theta)$. *Then*

- $\mathfrak{c}$ *is non-degenerate and holonomy non-degenerate.*
- $M(\mathfrak{c}) < 2\pi\sum_{\gamma\in\Theta} \ell_\gamma + \kappa^{-1}$.
- *The zero locus of* $\alpha$ *is a disjoint union of embedded circles whose components are in 1-1 correspondence with the integral curves of* $v$ *that comprise* $\Theta$. *This correspondence is such that*
  1) *Any given component of* $\alpha^{-1}(0)$ *lies in the radius* $\kappa r^{-1/2}$ *tubular neighborhood of its partner from* $\Theta$ *and it is isotopic in this neighborhood to this partner.*
  2) $|1 - |\alpha|| \leq \kappa(e^{-\sqrt{r}d/\kappa} + r^{-1})$ *at points with distance* $d$ *or more from* $\cup_{\gamma\in\Theta}\gamma$.
  3) *Let* $D \subset Y$ *denote an oriented, embedded disk with all boundary points at distance greater than* $\kappa r^{-1/2}$ *from* $\cup_{\gamma\in\Theta}\gamma$ *and with algebraic intersection number 1 with* $\cup_{\gamma\in\Theta}\gamma$. *Then* $\frac{i}{2\pi}\int_D F_{\hat{A}} = 1$.
- *View* $\Phi^r$ *as a map from* $\mathcal{Z}_{ech,M}{}^{L} \times \mathbb{Z}$ *to* $\mathcal{Z}_{SW,r} \times \mathbb{Z}$ *that acts as the identity on the* $\mathbb{Z}$ *factor. As such,* $\Phi^r$ *defines a* $\mathbb{Z}$*-equivariant covering map* $\hat{\Phi}^r$: $\hat{\mathcal{Z}}_{ech,M}{}^{L} \to \hat{\mathcal{Z}}_{SW,r}$ *via the isomorphisms described above that reverses the sign of the relative* $\mathbb{Z}$ *or* $\mathbb{Z}/p_M$ *degrees.*

The map $\Phi^r$ is constructed by copying in an almost verbatim fashion some of what is done in Section 3 of [T9] to construct an analogous map in the context where $\hat{a}$ is a contact 1-form and $w = \frac{1}{2}\, d\hat{a}$. The latter version of $\Phi^r$ is the map that is described in Theorems 4.2 of [T2] and Theorem 1.1 of [T9]. This contact 1-form incarnation of $\Phi^r$ is constructed in Section 3 of [T9] and its salient properties are stated as Theorem 1.1 in [T4] and Theorem 1.1 in [T6]. These theorems are proved respectively in Section 2 of [T4] and Section 2 of [T6]. As explained below only the simplest case of what is done in Section 3 of [T9], Section 2 of [T4] and Section 2 of [T6] are needed for what follows because of certain special features of the closed integral curves of $v$ that arise from elements in $\mathcal{Z}_{ech,M}$.

By way of a look ahead, Section 3a summarizes some basic facts about a particular subset of solutions to the vortex equations that are used to construct $\Phi^r$. The proof of Proposition 3.1 is given in Section 3b. Section 3c has the proof of Lemma 3.2 from Section 3a.

**a) The vortex equations II**

The proof Proposition 3.1 makes reference to (2.8)'s vortex equations. Of relevance here are the solutions which are such that $1 - |\alpha_0|^2$ is integrable on $\mathbb{C}$. The discussion of this subset of solutions to (2.8) has four parts. What is said in Parts 1-3 summarize various observations from Section 2 in [T9].

*Part 1*: As all complex line bundles over $\mathbb{C}$ are isomorphic to the product line bundle $\mathbb{C}\times\mathbb{C}$, no generality is lost by the focus in what follows on solutions $(A_0, \alpha_0)$ with $A_0$ a connection on this product bundle and $\alpha_0$ a complex function. Introduce $\theta_0$ to denote the product connection on the product line bundle $\mathbb{C} \times \mathbb{C}$. Write any given connection on $\mathbb{C} \times \mathbb{C}$ as $\theta_0 + \hat{a}$ with $\hat{a}$ being an $i\mathbb{R}$-valued 1-form on $\mathbb{C}$. Doing so identifies the set of solutions to (2.8) with a subset of the space $C^\infty(\mathbb{C}; iT^*\mathbb{C}) \times C^\infty(\mathbb{C}; \mathbb{C})$. This identification endows the set of solutions with a topology. Meanwhile, there is a free action of $C^\infty(\mathbb{C}; S^1)$ on the space of solutions to (2.8) whereby a given map u sends a given solution $(A_0, \alpha_0)$ to $(A_0 - u^{-1}du, u\alpha_0)$. This action is continuous, and so the set of $C^\infty(\mathbb{C}; S^1)$ equivalence classes of solutions has the induced quotient topology. The resulting subspace of solutions with $1 - |\alpha_0|^2$ being integrable is a disjoint union of components labeled by the non-negative integers. The integer m component consists of the set of equivalence classes of solutions that obey

$$\frac{1}{2\pi} \int_{\mathbb{C}} (1 - |\alpha_0|^2) = m . \tag{3.1}$$

The integer m component is denoted by $\mathfrak{C}_m$.

The space $\mathfrak{C}_m$ has the structure of a complex manifold that is holomorphically isomorphic to $\mathbb{C}^m$. The m complex functions $\{\sigma_q\}_{1\le q\le m}$ defined by

$$\sigma_q = \frac{1}{2\pi} \int_{\mathbb{C}} z^q (1 - |\alpha|^2) \tag{3.2}$$

define such an isomorphism. There are no convergence issues with regards to the integral in (3.2) by virtue of the fact that a solution to (2.8) and (3.1) obeys

- $|\alpha_0| < 1$ *with equality if and only if* $|\alpha_0| = 1$.

- $1 - |\alpha_0|^2 \le c_0\, e^{-\sqrt{2}\,\mathrm{dist}(\cdot,\alpha_0^{-1}(0))}$ .

(3.3)

Here, $c_0$ depends only on the integer m. As it turns out, the zeros of $\alpha$ are the roots of the polynomial $z \to \wp(z) = z^m + \sigma_1 z^{m-1} + \sigma_2 z^{m-2} + \cdots + \sigma_m$.

*Part 2*: Let $L \to \mathbb{C}$ denote a complex, Hermitian line bundle. Suppose for the moment that $(A_0,\alpha_0)$ defines a pair of Hermitian connection on L and section of L. Define the operator $\vartheta$ on $C^\infty(\mathbb{C},L)$ by the rule

$$(x, \iota) \to (\partial x + \tfrac{1}{\sqrt{2}}\,\bar{\alpha}_0\,\iota,\ \bar{\partial}_{A_0}\,\iota + \tfrac{1}{\sqrt{2}}\,\alpha_0\, x)$$

(3.4)

Here, $\partial$ is the holomorphic derivative on $C^\infty(\mathbb{C};\mathbb{C})$. The formal $L^2$ adjoint of $\vartheta$ is denoted by $\vartheta^\dagger$ and it is given by the rule

$$\vartheta^\dagger(c,\varsigma) = (-\bar{\partial}\, c + \tfrac{1}{\sqrt{2}}\,\bar{\alpha}_0\varsigma, -\partial_{A_0}\,\varsigma + \tfrac{1}{\sqrt{2}}\,\alpha_0\, c)\ .$$

(3.5)

The corresponding Laplacian $\vartheta\vartheta^\dagger$ can be written as

$$\vartheta\vartheta^\dagger(c,\varsigma) = ((-\partial\bar{\partial} + \tfrac{1}{2}\,|\alpha_0|^2)\, c, (-\bar{\partial}_{A_0}\partial_{A_0} + \tfrac{1}{2}\,|\alpha_0|^2)\,\varsigma) + \tfrac{1}{\sqrt{2}}\,(\partial_{A_0}\bar{\alpha}_0\,\varsigma,\ \bar{\partial}_{A_0}\alpha_0\, c)\ .$$

(3.6)

What follows is an important observation to keep in mind: The right most term in (3.6) is zero when $(A_0,\alpha_0)$ obeys the vortex equations.

*Part 3*: Suppose now that $(A_0,\alpha_0)$ is a solution to (2.8) and (3.1). The (1,0) tangent space to the orbit of $(A_0,\alpha_0)$ in $\mathfrak{C}_m$ is canonically isomorphic to the vector space of square integrable pairs $(x,\iota)$ of complex valued functions that are annihilated by $\vartheta$. This identification is used implicitly in what follows. The vector space of square integrable elements annihilated by $\vartheta$ is called the *kernel* of $\vartheta$.

Introduce the Hermitian inner product on the kernel of $\vartheta$ defined by the rule that sends an ordered pair $(\mathfrak{w} = (x, \iota), \mathfrak{w}' = (x',\iota'))$ in the kernel of $\vartheta$ to

$$\langle \mathfrak{w},\mathfrak{w}' \rangle = \tfrac{1}{\pi}\int_{\mathbb{C}} (\bar{x}x' + \bar{\iota}\iota')\ .$$

(3.7)

This Hermitian inner product is compatible with the complex structure and it defines a complete Kähler metric on $\mathfrak{C}_m$. In the case m = 1, this is the standard metric on $\mathfrak{C}_m = \mathbb{C}$, but such is not the case if m > 1.

Given a real number $\nu$ and a complex number $\mu$, define the function $\mathcal{h}$ on $\mathfrak{C}_m$ by the following rule: If $\mathfrak{c} = (A_0, \alpha_0)$, then

$$\mathcal{h}(\mathfrak{c}) = \tfrac{1}{4\pi} \int_{\mathbb{C}} (2\nu\,|z|^2 + \mu\bar{z}^2 + \bar{\mu}z^2)(1 - |\alpha_0|^2) \tag{3.8}$$

Suppose now that $\nu$ and $\mu$ depend on $t \in \mathbb{R}$ so that (3.6) defines a function on $\mathbb{R} \times \mathfrak{C}_m$. The Kähler metric on $\mathfrak{C}_m$ defines an associated symplectic form, and the latter with the $\mathbb{R}$ dependent function $\mathcal{h}$ define a corresponding 1-parameter family of Hamiltonian vector fields on $\mathfrak{C}_m$. An integral curve of this 1-parameter family of vector fields constitutes a map, $t \to \mathfrak{c}(t) \in \mathfrak{C}_m$, from $\mathbb{R}$ to $\mathfrak{C}_m$ that obeys the equation

$$\tfrac{i}{2}\,\mathfrak{c}' + \nabla^{(1,0)}\mathcal{h}|_{\mathfrak{c}} = 0 \ , \tag{3.9}$$

where $\mathfrak{c}'$ is shorthand for the $(1,0)$ part of $\mathfrak{c}_* \frac{d}{dt}$ , and where $\nabla^{(1,0)}\mathcal{h}$ denotes the $(1,0)$ part of the gradient of $\mathcal{h}$. Of interest in what follows are the solutions to (3.7) that obey the condition $\mathfrak{c}(t + T) = \mathfrak{c}(t)$ for some $T \geq 0$. Such a solution is said to be a *periodic* solution.

Let $\mathfrak{c}$: $S^1 \to \mathfrak{C}_m$ denote a given map. Associate to $\mathfrak{c}$ the bundle $\mathfrak{c}^*T_{1,0}\mathfrak{C}_m \to S^1$. The pull-back of the Riemannian connnection on $T\mathfrak{C}_m$ defines a Hermitian connection on $S^1$. The map $\mathfrak{c}$ is said to be *nondegenerate* when the operator

$$\zeta \to \tfrac{i}{2}\,\nabla_t\zeta + (\nabla_{\zeta_{\mathbb{R}}} \nabla^{1,0}\mathcal{h})|_{\mathfrak{c}} \tag{3.10}$$

on $C^\infty(S^1; \mathfrak{c}^*T_{1,0}\mathfrak{C}_m)$ has trivial kernel. The notation here is such that $\nabla_t$ denotes the covariant derivative of the afore-mentioned Hermitian connection. Also, $(\nabla_{\zeta_{\mathbb{R}}} \nabla^{1,0}\mathcal{h})|_{\mathfrak{c}}$ denotes the covariant derivative at $\mathfrak{c}$ along the vector defined by $\zeta$ in $T\mathfrak{C}_m|_{\mathfrak{c}}$ of the vector field $\nabla^{1,0}\mathcal{h} \in C^\infty(\mathfrak{C}_m; T_{1,0}\mathfrak{C}_m)$.

*Part 4*: Let $\gamma$ denote an integral curve of $v$. The proof of Proposition 3.1 refers to a certain pair of $\mathbb{R}$ valued and $\mathbb{C}$ valued functions $\gamma$ that are associated to a given Hermitian isomorphism between $K^{-1}|_\gamma$ and $\gamma \times \mathbb{C}$. To define these functions, fix a $\mathbb{C}$-linear, Hermitian isomorphism between $K^{-1}|_\gamma$ and $\gamma \times \mathbb{C}$. Let z denote the complex coordinate on $\mathbb{C}$ and let t denote an affine parameter for $\gamma$ such that $\gamma_* \frac{\partial}{\partial t} = \frac{\ell_\gamma}{2\pi} v$. Use the metric's exponential map along $\gamma$ to give a tubular neighborhood of $\gamma$ in Y the coordinates $(t, z)$ with it understood that these coordinates are only valid when z is a small radius disk about the origin in $\mathbb{C}$. Use these coordinates with the first order Taylor's expansion to write $w$ as $w = \frac{i}{2}\, dz\wedge d\bar{z} - 2(\nu z + \mu \bar{z})\, d\bar{z} \wedge dt - 2(\nu \bar{z} + \bar{\mu} z) dz\wedge dt + \cdots$ where $\nu$ and $\mu$ are respectively $\mathbb{R}$ and $\mathbb{C}$ valued functions on $S^1$, and where the unwritten terms are bounded by $c_0|z|^2$. Note that $\nu$ must be $\mathbb{R}$ valued because $w$ is closed.

The pair $(\nu,\mu)$ are the desired pair. Use this pair to define the function $\mathfrak{h}$. This done, fix a non-negative integer m and use $\mathfrak{C}_{(\gamma,m)}$ to denote the set of periodic solutions to (3.9) on $\mathfrak{C}_m$.

**Lemma 3.2**: *Suppose that* $\Theta \in \mathcal{Z}_{ech,M}$ *and that* $\gamma$ *is a closed integral curve of* $v$ *from* $\Theta$.

- *The space* $\mathfrak{C}_{(\gamma,1)}$ *consists of the constant map from* $S^1$ *to the* $\sigma_1 = 0$ *point in* $\mathfrak{C}_1$. *This is the equivalence of solutions to (2.8) and the* $m = 1$ *version of (3.1) with* $\alpha_0^{-1}(0) = 0$. *This solution is non-degenerate.*
- *Suppose that* $\mathfrak{p} \in \Lambda$ *and that* $\gamma \in \{\hat{\gamma}_{\mathfrak{p}}^{+}, \hat{\gamma}_{\mathfrak{p}}^{-}\}$. *Then* $\mathfrak{C}_{(\gamma,m)} = \emptyset$ *when* $m > 1$.

The proof of Lemma 3.2 constitutes Section 3c.

## b) Proof of Proposition 3.1

The proof differs only cosmetically from that of Theorem 4.2 in [T2]. As with the proof of the latter, there are three parts: Part 1 constructs the map $\Phi^r$ and verifies what is said in the second and third bullets. Part 2 proves what is said in the first and fourth bullets; and Part 3 verifies that image of $\Phi^r$ contains the $M < E$ subset of $\mathcal{Z}_{SW,r}$.

*Part 1*: The map $\Phi^r$ is constructed by copying what is done for its namesake in Theorem 1.1 from Section 1d in [T9]. The construction here constitutes what is perhaps the simplest of cases because only $\Theta \in \mathcal{Z}_{ech,M}{}^{L}$ and $\gamma \in \Theta$ versions of $\mathfrak{C}_{(\gamma,1)}$ are used. By way of a parenthetical remark, the first step in the construction of $\Phi^r$ uses the data from $\Theta$ to build a pair in $Conn(E) \times C^\infty(Y;\mathbb{S})$ that comes very close to solving (1.14). This construction is described in the first subsection of the appendix.

The first bullet of Lemma 3.2 guarantees that each map from $\{\{\mathfrak{C}_{(\gamma,1)}\}_{\gamma\in\Theta}\}_{\Theta\in\mathcal{Z}_{ech,M}}$ is non-degenerate; and so this set can be used as input for Theorem 1.1 in [T9]. The assertions made by the second and third bullets all follow directly from the construction and from Lemmas 2.1 and 2.3.

*Part 2*: To see about the first bullet of the proposition, suppose that $\Theta \in \mathcal{Z}_{ech,M}{}^{L}$ and that $\mathfrak{c} = (A,\alpha)$ is a solution to (1.13) that defines the equivalence class $\Phi^r(\Theta)$. The assertion that $\mathfrak{c}$ is non-degenerate can be had by copying almost verbatim the arguments in Section 2a of [T4] that prove the analogous assertion in Theorem 1.1 from [T4]. There are no substantive changes to these arguments from [T4]. The assertion that $\mathfrak{c}$ is holonomy non-degenerate follows from the third and fourth bullets of Proposition 2.4. To elaborate, these bullets imply that the connection $\hat{A}$ has a covariantly constant section near the integral curve $\gamma^{(z_0)}$. If $\hat{A}$ is flat near $\gamma^{(z_0)}$, then $\chi(\hat{A} - A_E)$ is an integer because $A_E$ was chosen to have holonomy 1 around $\gamma^{(z_0)}$.

To argue for the fourth bullet, fix $\Theta \in \mathcal{Z}_{ech,M}{}^L$ and suppose that $(A, \alpha)$ is a solution to (1.13) that defines the equivalence class $\Phi^r(\Theta)$. Fix a 2-cycle $Z \in H_2(M; [\Theta] - [\Theta_0])$ that has algebraic intersection number zero with $\gamma^{(z_0)}$. The $\mathbb{Z}$-equivariant covering map $\hat{\Phi}^r$ sends the equivalence class of a pair $(\Theta, Z)$ to the $\mathcal{G}_{M_\Lambda}$-orbit of a solution $(A, \psi)$ of (1.13) with the property that $\chi(\hat{A} - A_E) = 0$. With this point understood, the argument for the third bullet differ only cosmetically from those used [T4] to prove an analogous assertion from Theorem 4.2 in [T2]. The latter theorem follows directly from the relative degree assertion about the namesake $\Phi^r$ that appear in Theorem 1.1 in [T4]. The proof of this part of [T4]'s Theorem 1.1 constitute Sections 2b and 2c of [T4]. Note in this regard that the assumption that is made in Equation (2.56) in [T4] is not needed by virtue of the fact that the map $\Phi^r$ is constructed using only elements from the set $\{\{\mathfrak{C}_{(\gamma,1)}\}_{\gamma\in\Theta}\}_{\Theta\in\mathcal{Z}_{ech,M}}$.

*Part 3*: But for one additional substantive remark, the arguments for Theorem 1.1 in [T6] can copied with only notational changes to prove that if L is large, then $\Phi^r$ maps $\mathcal{Z}_{ech,M}{}^L$ onto the $M < E$ subset of $\mathcal{Z}_{SW,r}$ when r is large. The extra remark concerns the input to Theorem 1.1 of [T6] of a set denoted by $\mathfrak{C}\mathcal{Z}^L$ and a subset $\mathfrak{C}\mathcal{Z}^{L*} \subset \mathfrak{C}\mathcal{Z}^L$. Theorem 1.1 in [T6] requires $\mathfrak{C}\mathcal{Z}^{L*}$ to be the whole of $\mathfrak{C}\mathcal{Z}^L$. As explained below, $\mathfrak{C}\mathcal{Z}^L$ in this case is $\{\times_{\gamma\in\Theta} \mathfrak{C}_{(\gamma,1)}\}_{\Theta\in\mathcal{Z}_{ech,M}{}^L}$ and that $\mathfrak{C}\mathcal{Z}^{L*}$ in this case is indeed all of $\mathfrak{C}\mathcal{Z}^L$.

To define $\mathfrak{C}\mathcal{Z}^L$, introduce first $\mathcal{Z}$ to denote the set whose typical element consists of a finite collection of pairs whose first entry is a closed integral curve of $v$ and whose second entry is a positive integer. Let $\Theta$ denote such a collection. This set is constrained in two ways: Distinct pairs from $\Theta$ contain distinct closed orbits of $v$. The second constraint requires $[\Theta] = \sum_{(\gamma,m)\in\Theta} m[\gamma] \in H_1(Y; \mathbb{Z})$ to be the class that is defined by the elements in $\mathcal{Z}_{ech,M}$. The set $\mathcal{Z}_{ech,M}$ sits in $\mathcal{Z}$, but $\mathcal{Z}$ is strictly larger than $\mathcal{Z}_{ech,M}$; this can be seen using the parametrization given next.

Invoke Proposition II.2.8 or Theorem I.2.1 to write $\mathcal{Z}$ as $\mathcal{Z}_{HF} \times (\times_{\mathfrak{p}\in\Lambda}(\mathbb{Z} \times \hat{o}))$ where $\hat{o}$ is the set $\{0, 1, 2, \ldots\} \times \{0, 1, 2, \ldots\}$. This parametrization is such that the factor $\mathcal{Z}_{HF} \times (\times_{\mathfrak{p}\in\Lambda}\mathbb{Z})$ parametrizes pairs of the form $(\gamma, 1)$ with $\gamma \subset M_\delta \cup (\cup_{\mathfrak{p}\in\Lambda}\mathcal{H}_\mathfrak{p})$ crossing each $\mathfrak{p} \in \Lambda$ version of $\mathcal{H}_\mathfrak{p}$ once. To explain the factors of $\hat{o}$, fix $\mathfrak{p} \in \Lambda$. The element $\{0, 0\}$ in $\mathfrak{p}$'s version of $\hat{o}$ signifies that neither $\hat{\gamma}_\mathfrak{p}^+$ nor $\hat{\gamma}_\mathfrak{p}^-$ appears in $\Theta$. The element $(m_+, 0)$ from $\hat{o}$ with $m_+ > 0$ signifies that $\Theta$ contains $(\hat{\gamma}_\mathfrak{p}^+, m_+)$ and that $\Theta$ lacks a pair with $\hat{\gamma}_\mathfrak{p}^-$. By the same token, the element $(0, m_-)$ from $\hat{o}$ with $m_- > 0$ signifies that $\Theta$ contains the pair $(\hat{\gamma}_\mathfrak{p}^-, m_-)$ and that $\Theta$ lacks a pair with $\hat{\gamma}_\mathfrak{p}^+$. The element $(m_+, m_-)$ with both entries positive signifies the appearance in $\Theta$ of $(\hat\gamma_p^+, m_+)$ and $(\hat\gamma_p^-, m_-)$. Use $\mathcal{Z}^L$ to denote the subset of $\Theta = (\hat\upsilon, (\mathfrak{k}_p, (m_{p+}, m_{p-})_{p\in\Lambda}) \in \mathcal{Z}$ with $\sum_{p\in\Lambda}(\mathfrak{k}_p + 2m_{p+} + 2m_{p-}) < L$.

The set $\mathfrak{C}\mathcal{Z}^L$ maps to $\mathcal{Z}^L$ with fiber over any given $\Theta$ being $\times_{(\gamma,m)\in\Theta}\mathfrak{C}_{(\gamma,m)}$. The fiber over $\Theta$ of $\mathfrak{C}\mathcal{Z}^{L*}$ consists of the elements in $\times_{(\gamma,m)\in\Theta}\mathfrak{C}_{(\gamma,m)}$ whose entries are non-degenerate.

Granted these definitions, invoke the second bullet in Lemma 3.2 to see that $\mathfrak{C}\mathcal{Z}^L$ is indeed $\{\times_{\gamma\in\Theta}\mathfrak{C}_{(\gamma,1)}\}_{\Theta\in\mathcal{Z}_{ech,M}^L}$. This being the case, invoke the first bullet in Lemma 3.2 to see that $\mathfrak{C}\mathcal{Z}^{L*} = \mathfrak{C}\mathcal{Z}^L$.

**c) Proof of Lemma 3.2**

To prove the first bullet, view $\mathfrak{C}_1$ as $\mathbb{C}$ using (3.2)'s coordinate $\sigma_1$. Viewed in this way, then (3.9) is an equation for a function $t \to z(t)$ from $\mathbb{R}$ to $\mathbb{C}$, this being the equation

$$\tfrac{i}{2}\tfrac{d}{dt}z + \nu z + \mu\bar{z} = 0 . \tag{3.11}$$

Let $z(\cdot)$ denote a solution, but viewed as a map from $\mathbb{R}$ to $\mathbb{R}^2$. Then $z(2\pi)$ can be written as $U_\gamma z(0)$ where $U_\gamma \in SL(2;\mathbb{Z})$ is the linear return map that is described in Part 3 of Section II.1f. Proposition 2.7 of [KLTII] asserts that all of the relevant integral curves are hyperbolic, and by definition, this means that $U_\gamma$ has real eigenvalues with neither being 1 or -1. Thus (3.11) has a single solution, this being the constant map $t \to 0$. The operator (3.10) in this case is the operator that appears on the left side of (3.11), and so its kernel is trivial

The proof of the second bullet has ten parts. The Parts 1-3 say more about the solutions to the vortex equation. The remaining parts contain the proof proper. The arguments in Parts 4-10 focus on the case when $\gamma = \hat\gamma_p^+$. The arguments when $\gamma = \hat\gamma_p^-$ are essentially identical.

*Part 1*: The lemma that follows supplies three facts that play a central role in the proof of Lemma 3.2.

**Lemma 3.3**: *Fix* $m \geq 1$ *and suppose that* $(A_0, \alpha_0)$ *is a solution to (2.8) that defines a point in* $\mathfrak{C}_m$. *Then*

- $\frac{1}{\pi}\int_{\mathbb{C}}(\frac{1}{2}(1 - |\alpha_0|^2)^2 + |\partial_{A_0}\alpha_0|^2) = m.$
- $\frac{1}{2\pi}\int_{\mathbb{C}}(1 - |\alpha_0|^2)^2 \geq \frac{2}{5}m$ *and* $\frac{1}{\pi}\int_{\mathbb{C}}|\partial_A\alpha_0|^2 \leq \frac{3}{5}m$ .
- $|\partial_{A_0}\alpha_0| \leq \frac{\sqrt{3}}{2}(1 - |\alpha_0|^2)$.

***Proof of Lemma 3.3***: Use $\Delta$ in what follows to denote the Laplacian on $\mathbb{C}$. Meanwhile, let w denote the function $(1 - |\alpha_0|^2)$ and use g to denote $\partial_{A_0}\alpha_0$. It follows from (2.8) that

$$-\tfrac{1}{4}\Delta w + \tfrac{1}{2}|\alpha|^2 w = |g|^2 \quad \textit{and} \quad -\tfrac{1}{4}\nabla_A{}^{\dagger}\nabla_A g + \tfrac{1}{2}|\alpha|^2 g = \tfrac{3}{4} w g. \tag{3.12}$$

Write $|\alpha_0|^2 = -w + 1$ to write the left most equation in (3.12) as

$$-\tfrac{1}{4}\Delta w + \tfrac{1}{2} w = \tfrac{1}{2} w^2 + |g|^2 \ . \tag{3.13}$$

Integrate both sides of this equation and appeal to (3.1) to obtain the first bullet in the lemma.

The second bullet follows from the first and the third. To elaborate, use the third bullet of the lemma to see that the integral on the left hand side of the first bullet is less than $\frac{5}{4\pi}$ times the integral of $w^2$. As a consequence, the contribution of the term $\frac{1}{2} w^2$ to the integral on the left side of the first bullet is no less than $\frac{2}{5}$ m and so the contribution to this integral of $|g|^2$ is no greater than $\frac{3}{5}$ m.

To obtain the fourth bullet, use the right hand identity in (3.12) to see that

$$-\tfrac{1}{4}\Delta|g| + \tfrac{1}{2}|\alpha|^2|g| \le \tfrac{3}{4} w|g| \ . \tag{3.14}$$

To exploit (3.14), set $x = |g| - \frac{\sqrt{3}}{2} w$. The left most equation in (3.12) and (3.14) require

$$-\tfrac{1}{4}\Delta x + \tfrac{1}{2}|\alpha|^2 x \le |g|(\tfrac{3}{4} w - \tfrac{\sqrt{3}}{2}|g|) = -\tfrac{\sqrt{3}}{2}|g| x. \tag{3.15}$$

Granted (3.15), use the maximum principal to see that x can not have a positive local maximum. Given that x is square integrable, this implies that $x \le 0$ which is what is asserted by the second bullet of the lemma.

*Part 2*: Various additional facts about any given $m \in \{1, 2, \ldots.\}$ version of $\mathfrak{C}_m$ are needed for the proof. The first of these facts concerns an isometric and holomorphic action of the semi-direct product of $S^1$ and $\mathbb{C}$ on $\mathfrak{C}_m$. This action is induced by the group's action on $\mathbb{C}$ where $S^1$ acts as the group of rotations about the origin and $\mathbb{C}$ acts on itself by translation. The generator of the action of $\mathbb{C}$ on $\mathfrak{C}_m$ at the equivalence class of a solution $(A_0, \alpha_0)$ to (2.8) is the tangent vector that is defined by the element

$$\mathfrak{w}_1 = (\tfrac{1}{\sqrt{2}}(1 - |\alpha_0|^2), \partial_{A_0}\alpha_0) \ . \tag{3.16}$$

in the kernel of $\vartheta$. The action of $S^1$ on $\mathfrak{C}_m$ is such that $\eta \in S^1$ pulls back (3.2)'s functions $\{\sigma_q\}_{1\le q\le m}$ to $\{\eta^q\sigma_q\}_{1\le q\le m}$. The action has a unique fixed point in $\mathfrak{C}_m$, this given by the point where all $\sigma_q$ are zero. The latter point is the equivalence class of the solutions to (2.8) with $\alpha^{-1}(0) = 0$. The fixed point of the $S^1$ action is called the *symmetric vortex*.

*Part 3*: Define (3.4)'s operator $\vartheta$ using the solution $(A_0, \alpha_0)$. The absence of the right most term in (3.6) and the integrability of $1 - |\alpha_0|^2$ imply that $\vartheta\vartheta^\dagger$ has a bounded inverse that maps $L^2(\mathbb{C}; \mathbb{C} \times \mathbb{C})$ to the $L^2$-orthogonal complement of the kernel of $\vartheta$.

The other Laplacian, $\vartheta^\dagger\vartheta$, can be written as $\vartheta^\dagger\vartheta = \vartheta\vartheta^\dagger + \mathfrak{e}$ where $\mathfrak{e}$ is a zero'th order term that is bounded by $c_0(1 - |\alpha_0|^2)$. Given this last fact, the Bochner-Weitzenbock formula for $\bar\partial_A\partial_A$ can be used in conjunction with the left most equation in (3.12) and the maximum principal to see that any given element in the kernel of $\vartheta$ with $L^2$ norm 1 is bounded pointwise by $c_0(1 - |\alpha_0|^2)$. The argument also invokes the third bullet of Lemma 3.3 and (3.3). Granted the latter as input, the argument differs little from the argument in Part 1 proving the third bullet in Lemma 3.3. This being the case, the details are omitted.

*Part 4*: Use (3.2)'s coordinates $\{\sigma_q\}_{1\le q\le m}$ for $\mathfrak{C}_m$ so as to view the equation in (3.9) as an equation for a map, $t \to (\sigma_1(t), \sigma_2(t), \ldots, \sigma_m(t))$ from $\mathbb{R}$ to $\mathbb{C}^m$. In particular the map $t \to \sigma_2(t)$ must obey the equation

$$\frac{i}{2}\frac{d}{dt}\sigma_2 + 2\nu\sigma_2 + \mu\, g^{2\bar{2}} = 0 \tag{3.17}$$

where $g^{2\bar{2}}$ is the norm of $d\sigma_2$ as defined by the Kahler metric on $\mathfrak{C}_m$. Note that $g^{2\bar{2}}$ is a strictly positive function on $\mathfrak{C}_m$. Parts 6-10 explain why

$$g^{2\bar{2}} > 2|\sigma_2| \,. \tag{3.18}$$

Meanwhile, Part 5 explains why the functions $\nu$ and $\mu$ in (3.8) and (3.17) can be assumed constant, with $\mu$ real and such that $\mu > |\nu| + c_0^{-1}$. With the preceding understood, write $\sigma_2$ as $\sigma_x + i\sigma_y$ with $\sigma_x$ and $\sigma_y$ real valued functions. Then (3.17) and (3.18) require $-\frac{d}{dt}\sigma_y < 0$ and so there are no periodic solutions.

*Part 5*: The functions $\nu$ and $\mu$ that appear in (3.8) and (3.17) depend on a chosen unit length basis vector for the bundle $K^{-1}$ along the given closed integral curve, this being $\hat\gamma_p^+$. Even so, the question of existence or not of solutions to the corresponding version of (3.9) does not depend on the trivialization. This fact is exploited in what follows to choose a convenient trivialization.

The metric on $\mathcal{H}_p$ is invariant with respect to rotations of the coordinate $\phi$ and as $\nu$ is a constant multiple of $-\frac{\partial}{\partial\phi}$ along $\hat{\gamma}_p^+$ the basis vector for $K^{-1}$ along $\hat{\gamma}_p^+$ can be chosen so as to be covariantly constant along $\hat{\gamma}_p^+$. Choosing such a basis vector gives a pair $(\nu, \mu)$ with both being constants. As noted previously, $\nu$ must be real, and if $\mu$ is not real and non-negative to begin with, a suitable constant rotation of $\mathbb{C}$ changes the coordinates so that the resulting version of $\mu$ is real an non-negative.

The assertion that $\mu > |\nu|$ follows from the fact that $\hat{\gamma}_p^+$ is hyperbolic. By way of an explanation, the fact that $\mu$ and $\nu$ are constant can be used to solve (3.11) and thus write the matrix $U_\gamma$ and see directly its eigenvalues. These are real and neither 1 nor -1 if and only if $\mu > |\nu|$.

By way of a parenthetical remark, Part 2 of Section III.5a introduces the coordinates $(s_+, \phi_+, \theta_+, u_+)$ for the product of $\mathbb{R}$ with a tubular neighborhood in Y of $\hat{\gamma}_p^+$. These are such that the locus $\theta_+ = 0$, $u_+ = 0$ is the cylinder $\mathbb{R} \times \hat{\gamma}_p^+$ with $s_+$ being the Euclidean coordinate for the $\mathbb{R}$ factor and $\phi_+$ an $\mathbb{R}/(2\pi\mathbb{Z})$ valued coordinate for $\hat{\gamma}_p^+$. This understood, the differentials $d\theta_+$ and $du_+$ together define a trivialization of the normal bundle of $\hat{\gamma}_p^+$. Given this trivialization, the coefficients that appear in Equation (III.5.1) determine $\nu$ and $\mu$ as functions of the constants $x_0$ and R that are used in Section 1a to define the geometry of Y. A direct calculation using these coordinates will also verify the claim that $\nu$ and $\mu$ can be assumed constant, with $\mu$ real and greater than $|\nu|$.

*Part 6*: Suppose that $(A_0, \alpha_0)$ is a solution to (2.8) that defines a point in $\mathfrak{C}_m$. Let $\mathfrak{w} = (x, \iota)$ denote an element in the kernel of the operator $\vartheta$. The (1,0) part of $d\sigma_2$ pairs with the tangent vector defined by $\mathfrak{w}$ to give

$$-\frac{1}{2\pi}\int_{\mathbb{C}} z^2\bar{\alpha}\iota\ . \tag{3.19}$$

Since $\vartheta\mathfrak{w} = 0$, the first entry of (3.4) is zero and so the integrand in (3.19) can be replaced by $-\sqrt{2}\, z^2\partial x$. Having done so, integration by parts writes (3.19) as

$$-\frac{\sqrt{2}}{\pi}\int_{\mathbb{C}} z x\ . \tag{3.20}$$

Note in this regard that such an integration by parts is possible here (and in a subsequent integration by parts) by virtue of what is said in Part 3 to the effect that $|\mathfrak{w}|$ is bounded by a multiple of $(1 - |\alpha_0|^2)$, and thus is exponentially small where $|z|$ is large.

The integrand in (3.20) is the same as $z(1-|\alpha_0|^2)x + z\,\bar{\alpha}_0\alpha_0\,x$. As $\vartheta\mathfrak{w} = 0$, the left hand entry in (3.4) is zero, and so this is the same as $z(1-|\alpha_0|^2)x - \sqrt{2}\,z\,\bar{\alpha}_0\,\bar{\partial}_{A_0}\iota$. Use this fact with a second integration by parts to see that (3.19) is equal to

$$-\tfrac{2}{\pi}\int_{\mathbb{C}}(z\tfrac{1}{\sqrt{2}}(1-|\alpha|^2)x + z\bar{\partial}_A\bar{\alpha}\iota)\,. \tag{3.21}$$

Introduce $\Pi$ to denote the $L^2$ orthogonal projection from $L^2(\mathbb{C}; \mathbb{C}\oplus\mathbb{C})$ to the kernel of $\vartheta$. This last identity implies that $d\sigma_2$ acts on the kernel of $\vartheta$ as $\mathfrak{w} \to 2\langle\Pi(\bar{z}\,\mathfrak{w}_1),\mathfrak{w}\rangle$ with $\mathfrak{w}_1$ defined by (3.16). It follows as a consequence that

$$g^{2\bar{2}} = 4\langle\Pi(\bar{z}\,\mathfrak{w}_1),\Pi(\bar{z}\,\mathfrak{w}_1)\rangle\,. \tag{3.22}$$

Meanwhile, $\Pi(\bar{z}\,\mathfrak{w}_1)$ can be written as $\bar{z}\,\mathfrak{w}_1 + \Theta^\dagger\mathfrak{z}$ and so

$$\langle\Pi(\bar{z}\,\mathfrak{w}_1),\Pi(\bar{\alpha}\,\mathfrak{w}_1)\rangle = \tfrac{1}{\pi}\int_{\mathbb{C}}|z|^2\,(\tfrac{1}{2}(1-|\alpha_0|^2)^2 + |\partial_{A_0}\alpha_0|^2) - \langle\Theta^\dagger\mathfrak{z},\Theta^\dagger\mathfrak{z}\rangle\,. \tag{3.23}$$

With (3.23) in hand, write (3.22) as

$$g^{2\bar{2}} = \tfrac{2}{\pi}\int_{\mathbb{C}}|z|^2\,(\tfrac{1}{2}(1-|\alpha_0|^2)^2+|\partial_{A_0}\alpha_0|^2) + \tfrac{2}{\pi}\int_{\mathbb{C}}|z|^2\,(\tfrac{1}{2}(1-|\alpha_0|^2)^2+|\partial_{A_0}\alpha_0|^2) - 4\langle\vartheta^\dagger\mathfrak{z},\vartheta^\dagger\mathfrak{z}\rangle. \tag{3.24}$$

The comparison of $g^{2\bar{2}}$ with $2|\sigma_2|$ uses (3.24) with the rewriting of $\sigma_2$ as

$$\sigma_2 = \tfrac{1}{\pi}\int_{\mathbb{C}} z^2(\tfrac{1}{2}(1-|\alpha_0|^2)^2 + |\partial_{A_0}\alpha_0|^2)\,. \tag{3.25}$$

To obtain this last identity, multiply both sides of (3.13) by $\frac{1}{2\pi}z^2$ and integrate the resulting equation over $\mathbb{C}$. The integral of $\frac{1}{2\pi}z^2\Delta w$ is zero.

The inequality $g^{2\bar{2}} > 2|\sigma_2|$ follows directly from (3.24) and (3.25) if

$$\tfrac{1}{\pi}\int_{\mathbb{C}}|z|^2\,(\tfrac{1}{2}(1-|\alpha_0|^2)^2 + |\partial_{A_0}\alpha_0|^2) - 2\langle\vartheta^\dagger\mathfrak{z},\vartheta^\dagger\mathfrak{z}\rangle \ge 0. \tag{3.26}$$

The remaining Parts 7-10 supply a proof of this inequality.

*Part 7*: This step supplies an upper bound for $\langle\vartheta^\dagger\mathfrak{z}, \vartheta^\dagger\mathfrak{z}\rangle$. To this end, use (3.6) to see that $\mathfrak{z} = (0,\varsigma)$ with $\varsigma$ being the $L^2$ solution on $\mathbb{C}$ to the equation

$$-\bar{\partial}_{A_0}\partial_{A_0}\varsigma + \tfrac{1}{2}|\alpha_0|^2\varsigma = -\partial_{A_0}\alpha_0\,.$$

(3.27)

It follows as a consequence that

$$\langle\vartheta^\dagger\mathfrak{z}, \vartheta^\dagger\mathfrak{z}\rangle = \tfrac{1}{\pi}\int_{\mathbb{C}}(|\partial_{A_0}\varsigma|^2 + \tfrac{1}{2}|\alpha_0|^2|\varsigma|^2)\ .$$

(3.28)

Granted (3.28), it then follows from (3.27) that $\langle\vartheta^\dagger\mathfrak{z}, \vartheta^\dagger\mathfrak{z}\rangle \le \frac{1}{\pi}\|\varsigma\|_2\|\partial_{A_0}\alpha_0\|_2$ with $\|\cdot\|_2$ denoting here the $L^2$ norm. To see about the $L^2$ norm of $\varsigma$, commute derivatives using the top bullet in (2.2) to write the left hand side of (3.27) as $-\partial_{A_0}\bar{\partial}_{A_0}\varsigma + \frac{1}{2}\varsigma$. Take the $L^2$ inner product of both sides of the resulting equation with $\varsigma$. This leads to an equality between integrals. An application of Holder's inequality to the latter equality finds $\frac{1}{2}\|\varsigma\|_2^2 \le \|\varsigma\|_2\|\partial_{A_0}\alpha_0\|_2$ and so $\|\varsigma\|_2 \le 2\|\partial_{A_0}\alpha_0\|_2$. This being the case, then

$$\langle\vartheta^\dagger\mathfrak{z}, \vartheta^\dagger\mathfrak{z}\rangle \le \tfrac{2}{\pi}\int_{\mathbb{C}}|\partial_{A_0}\alpha_0|^2\ .$$

(3.29)

*Part 8*: This part of the proof exploits another identity coming from (3.13):

$$\tfrac{1}{\pi}\int_{\mathbb{C}}|z|^2\,(\tfrac{1}{2}(1-|\alpha_0|^2)^2 + |\partial_{A_0}\alpha_0|^2) = \tfrac{1}{2\pi}\int_{\mathbb{C}}|z|^2\,(1-|\alpha_0|^2) - 2m\,.$$

(3.30)

To derive (3.30), multiply both sides of (3.13) by $\frac{1}{\pi}|z|^2$ and then integrate the resulting equation over $\mathbb{C}$. An integration by parts identifies the integral over $\mathbb{C}$ of $-\frac{1}{4\pi}|z|^2\Delta w$ with that of $-\frac{1}{\pi}w$. According to (3.1), the latter integral is equal to -2m.

With (3.30) in mind, digress momentarily to consider a certain constrained minimization problem for a real valued, measurable function on $\mathbb{C}$. The problem asks for an infimum of the functional

$$u \to \mathfrak{s}(u) = \tfrac{1}{2\pi}\int_{\mathbb{C}}|z|^2\,u - 2m - 4(m - \tfrac{1}{2\pi}\int_{\mathbb{C}}u^2)$$

(3.31)

subject to the two constraints $0 \le u \le 1$ and $\frac{1}{2\pi}\int_{\mathbb{C}}u = m$. By way of an explanation, the function $u = (1-|\alpha_0|^2)$ obeys the constraints; and it follows from (3.13) with the first bullet of Lemma 3.3 that the value $\mathfrak{s}$ in this case is no greater than what is written on the left hand side of (3.26). As a consequence, (3.26) follows if the infimum of $\mathfrak{s}$ is positive.

As explained in momentarily, the functional $\mathfrak{s}$ takes on its minimum with the function, $u_*$, given as follows: Set $\lambda = 2m + 4$. Then

- $u_* = 1$ *where* $|z|^2 \le \lambda - 8$.
- $u_* = \frac{1}{8}(\lambda - |z|^2)$ *where* $\lambda - 8 \le |z|^2 \le \lambda$.
- $u_* = 0$ *where* $|z|^2 \ge \lambda$.

(3.32)

The value of $\mathfrak{s}$ on $u_*$ is $\frac{1}{4}\lambda^2 - \frac{16}{3} - 6m = m^2 - 2m - \frac{4}{3}$; and this is positive for all $m \ge 3$.

With regards to (3.32), note first that an averaging argument shows that any minimizer is a function of the radial coordinate on $\mathbb{C}$. Meanwhile, the variational equations for $\mathfrak{s}$ assert that a constrained minimizer, $u_*$, is such that $|z|^3 + 8|z|u_* = \lambda|z|$ where $0 \le u_* \le 1$. Here, $\lambda$ is the Lagrange multiplier for the constraint that the integral of $\frac{1}{2\pi}u$ is equal to m. Thus, the minimizer $u_*$ has the form that is depicted in (3.32) with $\lambda$ chosen so that this integral constraint is obeyed. A calculation finds that $\lambda = 2m+4$ and a second calculation finds that $\mathfrak{s}(u_*) = \frac{1}{4}\lambda^2 - \frac{16}{3} - 6m$.

*Part 9*: The verification of (3.26) when $m = 2$ requires more care with regards to the difference between $|\sigma_2|$ and the integral of $\frac{1}{\pi}|z|^2(\frac{1}{2}(1-|\alpha_0|^2)^2 + |\partial_{A_0}\alpha_0|^2)$. As explained in Part 10, this difference is no less than $2q$ where $q$ is the integral of this same function in the case when $(A_0, \alpha_0)$ is a symmetric solution to (2.8) from the space $\mathfrak{C}_1$. By way of a reminder, the space $\mathfrak{C}_1$ is diffeomorphic to $\mathbb{C}$ with the diffeomorphism given by the function $\sigma_1$. The symmetric solution is the $\sigma_1 = 0$ point, this the solution with $\alpha^{-1}(0) = 0$. This step proves that $q > \frac{2}{3}$. Granted the latter, then (3.24) and what is said in Part 8 find

$$g^{2\bar{2}} > 2|\sigma_2| + 4q - \tfrac{8}{3} > 2|\sigma_2| .$$

(3.33)

To derive the asserted lower bound for $q$, introduce w to denote $(1 - |\alpha_0|^2)$ and introduce g to denote $\partial_{A_0}\alpha_0$ for a $\mathfrak{C}_1$ version of $(A_0,\alpha_0)$ with $\alpha_0^{-1}(0) = 0$. Use $\rho = |z|$ to denote the radial coordinate on $\mathbb{C}$. Then $\partial_\rho w \le 0$ since w is rotationally symmetric and has no local maxima. Meanwhile, $|\partial_\rho w| = \sqrt{2}|\alpha_0||g| < \sqrt{2}|g|$. What with Lemma 3.3, this finds $|\partial_\rho w| < \frac{\sqrt{3}}{\sqrt{2}}|w|$. Keeping in mind that w is exponentially small at large $\rho$, integration by parts finds that

$$0 = \int_0^\infty \partial_\rho(\rho^2 w)d\rho = 2\int_0^\infty w\rho d\rho + \int_0^\infty (\partial_\rho w)\rho^2 d\rho .$$

(3.34)

Given what was said about $\partial_\rho w$, this last equation implies that

$$\int_0^\infty w\rho^2 d\rho > \tfrac{2\surd 2}{\surd 3}\int_0^\infty w\rho d\rho = \tfrac{2\surd 2}{\surd 3}\,.$$

(3.35)

By way of explanation, (3.1) asserts that the integral on the right side is equal to 1. To continue, use Hölder's inequality with (3.1) to see that the left hand side of (3.35) is no less than

$$\Big(\int_0^\infty w\rho^3 d\rho\Big)^{1/2}\Big(\int_0^\infty w\rho d\rho\Big)^{1/2} = \Big(\int_0^\infty w\rho^3 d\rho\Big)^{1/2}\,.$$

(3.36)

Taken together, (3.35) and (3.36) assert that

$$\int_0^\infty w\rho^3 d\rho > \tfrac{8}{3}\,.$$

(3.37)

This last equation with (3.30) say that $q > \frac{2}{3}$.

*Part 10*: Suppose now that $(A_0, \alpha_0)$ is a solution to (2.8) that defines a point in $\mathfrak{C}_2$. This step explains why

$$\tfrac{1}{\pi}\int_{\mathbb{C}} |z|^2 \big(\tfrac{1}{2}(1-|\alpha_0|^2)^2 + |\bar\partial_{A_0}\alpha_0|^2\big) > |\sigma_2| + 2q\,.$$

(3.38)

To this end, note (3.38) holds if the left hand side is greater than $2q$ plus the real part of $\bar{u}\sigma_2$ for any $u \in S^1$. Therefore, no generality is lost to prove that the left hand side of (3.38) is greater than $2q$ plus the real part of $\sigma_2$. Reintroduce $\sigma_x$ to denote this real part. Let $\sigma$: $\mathfrak{C}_2 \to (0, \infty)$ denote the function given by the left side integral in (3.38).

**Lemma 3.4**: *The function* $\sigma - \sigma_x$ *does not take on its infimum at any point in* $\mathfrak{C}_2$. *Furthermore, sequences in* $\mathfrak{C}_2$ *on which* $\sigma - \sigma_x$ *converges to its infimum have the following properties: Fix* $R \geq 1$ *and the all but a finite number of elements in the sequence are* $C^\infty(\mathbb{C}; S^1)$ *orbits of pairs* $(A_0, \alpha_0)$ *with* $\alpha_0$ *such that its two zeros have distance* R *or greater between them. Moreover, these zeros have distance* $\frac{1}{R}$ *or less from the real z axis in* $\mathbb{C}$.

Granted for the moment Lemma 3.4, write the coordinate z as x + iy with x and y being real, and then use (3.29) to write

$$\sigma - \sigma_x = \frac{1}{\pi}\int_{\mathbb{C}} y^2(1 - |\alpha|^2) - 2m.$$

(3.39)

It follows from (2.4) in [T4] that if R >> 1 and $(A_0,\alpha_0)$ is as described in Lemma 3.4, then what is written on the left hand side of (3.39) differs by at most $c_0R^{-1}$ from twice its value for the case where m = 1 and $(A_0,\alpha_0)$ is the $\sigma_1 = 0$ solution in $\mathfrak{C}_1$. Meanwhile, the $\sigma_1 = 0$ solution in $\mathfrak{C}_1$ is invariant with respect to the $S^1$ action on $\mathfrak{C}_1$ and so the m =1 and $\sigma_1 = 0$ version of (3.38) is equal to $q$.

***Proof of Lemma 3.4***: Fix an element in $\mathfrak{C}_2$ and write the zeros of any corresponding solution to (2.8) as an unordered pair $(z_1,z_2) \in \mathrm{Sym}^2(\mathbb{C})$. Fix pairs, $(A_1,\alpha_1)$ and $(A_2,\alpha_2)$, of m = 1 solutions to (2.8) with $\alpha_1^{-1}(0) = z_1$ and with $\alpha_2^{-1}(0) = z_2$. Part 4 in Section 2a of [T9] writes the given $\mathfrak{C}_2$ element as the $C^\infty(\mathbb{C};S^1)$ orbit of an m = 2 solution to (2.8) that can be written as $(A_0, \alpha_0)$ with

$$A_0 = A_1 + A_2 + (\bar\partial u\, d\bar z - \partial u\, dz) \quad \textit{and} \quad \alpha_0 = e^{-u}\alpha_1\alpha_2$$

(3.40)

such that u is a smooth, real valued function on $\mathbb{C}$ which obeys $|u| \le c_0 e^{-\mathrm{dist}(\cdot,\alpha^{-1}(0))/c_0}$. The top line in (2.8) requires u to obey

$$\Delta u = (1 - e^{-2u}|\alpha_1|^2|\alpha_2|^2) - (1 - |\alpha_1|^2) - (1 - |\alpha_2|^2)\ .$$

(3.41)

Were $u \le 0$, then the right hand side of (3.41) would be less than $-(1 - |\alpha_1|^2)(1 - |\alpha_2|^2)$ and thus not positive. This being the case, the maximum principal demands that $u > 0$. With this in mind, multiply both sides of (3.41) by $y^2$ and then integrate the result over $\mathbb{C}$. An integration by parts writes the integral of $y^2\Delta u$ as twice the integral of u. In particular, the integral of u is positive, and so

$$\int_{\mathbb{C}} y^2(1 - |\alpha|^2) \ge \int_{\mathbb{C}} y^2(1 - |\alpha_1|^2) + \int_{\mathbb{C}} y^2(1 - |\alpha_2|^2)\ .$$

(3.42)

Meanwhile, the bound on $u \le c_0 e^{-\mathrm{dist}(\cdot,\alpha^{-1}(0))/c_0}$ implies that

$$\int_{\mathbb{C}} y^2(1 - |\alpha|^2) \le \int_{\mathbb{C}} y^2(1 - |\alpha_1|^2) + \int_{\mathbb{C}} y^2(1 - |\alpha_2|^2) + c_0 e^{-|z_1 - z_2|/c_0}\ .$$

(3.43)

The assertions of Lemma 3.4 follow directly from (3.42) and (3.43).

## 4. Instantons

The purpose of this section is to provide various facts about the solutions to the r and $\mathfrak{g} = \mathfrak{e}_\mu$ versions of (1.20), this being the version reproduced below.

- $\frac{\partial}{\partial s} A + B_A - r(\psi^\dagger \tau \psi - i\hat{a}) + \frac{1}{2} B_{A_K} - i{*}d\mu = 0$ .
- $\frac{\partial}{\partial s} \psi + D_A \psi = 0.$

(4.1)

These facts assert a priori bounds on various integrals on pointwise norms.

### a) Apriori integral bounds

The analysis of (4.1) concerns the versions with $r > \pi$ and with $\mu \in \Omega$ a given element with $\mathcal{P}$-norm bounded by 1. Assume in what follows that $\mu$ is such that all solutions to (1.13) are non-degenerate.

To set some notation, suppose that $\mathfrak{d}$: $\mathbb{R} \to \text{Conn}(E) \times C^\infty(Y; \mathbb{S})$ is a given instanton solution to (4.1). The $s \to -\infty$ limit of $\mathfrak{d}$ is denoted by $\mathfrak{c}_-$ and the $s \to \infty$ limit by $\mathfrak{c}_+$. The latter are solutions to (1.13). The respective Conn(E) and $C^\infty(Y; \mathbb{S})$ components of $\mathfrak{d}$ are written as (A, $\psi$) and $\psi$ is often written in two component form as ($\alpha$, $\beta$). The lemmas that follow use $A_{\mathfrak{d}}$ to denote $\mathfrak{a}(\mathfrak{c}_-) - \mathfrak{a}(\mathfrak{c}_+)$.

Many of the lemmas here and rest of Section 4 have analogs in Section 3 of [T6]. Except for one item, the statement of a given lemma here is virtually identical to the statement of its partner in Section 3 of [T6]. Various lemmas in Section 3 of [T6] give the option of assuming the lower bound $\mathfrak{f}_\mathfrak{s}(\mathfrak{c}_+) - \mathfrak{f}_\mathfrak{s}(\mathfrak{c}_-) > -r^2$ in lieu of an upper bound on $A_{\mathfrak{d}}$. Their partners here do not give such an option. This difference is due solely to the term $2\pi r \mathfrak{f}_\mathfrak{s}$ in (1.29)'s formula for $\mathfrak{a}^f$. The version of $\mathfrak{a}^f$ used in Section 3 of [T6] has $\mathfrak{f}_\mathfrak{s}$ appearing only as $-2\pi^2 \mathfrak{f}_\mathfrak{s}$ while the version here has $2\pi(r - \pi)\mathfrak{f}_\mathfrak{s}$. Of relevance here is the sign difference when $r > \pi$.

Except for what was just said about $\mathfrak{f}_\mathfrak{s}(\mathfrak{c}_+) - \mathfrak{f}_\mathfrak{s}(\mathfrak{c}_-)$, the proof of most every lemma here is virtually identical to that of its partner in Section 3 of [T6]. When this is the case, the reader is referred to Section 3 of [T6] for the proofs. The correspondence between lemmas here and lemmas in Section 3 of [T6] are noted below. Be forwarned however that the lemmas in Section 3 of [T6] do not appear in the same order as those here.

The first lemma below supplies an inequality that relates $A_{\mathfrak{d}}$ to the change in $\mathfrak{f}_\mathfrak{s}$.

**Lemma 4.1**: *There exists a constant* $\kappa \geq 1$ *with the following significance: Suppose that* $r \geq \kappa$*, that* $\mu \in \Omega$ *with* $\mathcal{P}$*-norm less than 1. Suppose that* $\mathfrak{c}_+$ *and* $\mathfrak{c}_-$ *are solutions to the* $(r,\mu)$ *version of (1.13). Then* $\mathfrak{a}(\mathfrak{c}_-) - \mathfrak{a}(\mathfrak{c}_+) \leq 2\pi(r - \pi)(\mathfrak{f}_\mathfrak{s}(\mathfrak{c}_+) - \mathfrak{f}_\mathfrak{s}(\mathfrak{c}_-)) + \kappa r(M(\mathfrak{c}_+) + 1)$.

***Proof of Lemma 4.1***: Write $\mathfrak{a}(\mathfrak{c}_-) - \mathfrak{a}(\mathfrak{c}_+)$ as $\mathfrak{a}^f(\mathfrak{c}_-) - \mathfrak{a}^f(\mathfrak{c}_+) + 2\pi(r - \pi)(\mathfrak{f}_\mathfrak{s}(\mathfrak{c}_+) - \mathfrak{f}_\mathfrak{s}(\mathfrak{c}_-))$ and then appeal to the fourth bullet of Proposition 2.7.

The next lemma refers to a certain $i\mathbb{R}$-valued 1-form that can be associated to a given $(A,\psi) \in \text{Conn}(E) \times C^\infty(Y;\mathbb{S})$. This 1-form is denoted by $\mathfrak{B}_{(A,\psi)}$:

$$\mathfrak{B}_{(A,\psi)} = B_A - r\,(\psi^\dagger\tau^k\psi - i\,\hat{a}) - i*d\mu + \tfrac{1}{2} B_{A_K}\,. \tag{4.2}$$

The upcoming Lemma 4.2 gives an a priori bound for the $L^2$ norms of $\frac{\partial}{\partial s} A$, $\mathfrak{B}_{(A,\psi)}$, $\frac{\partial}{\partial s}\psi$ and $D_A\psi$. Lemma 4.2 is partnered with Lemma 3.4 in [T6] and its proof is identical to the latter's but for notation.

**Lemma 4.2**: *There exists a constant* $\kappa \geq 1$ *with the following significance: Suppose that* $r \geq \kappa$*, that* $\mu \in \Omega$ *has* $\mathcal{P}$*-norm less than 1 and that* $(A, \psi)$ *is an instanton solution to the* $(r, \mu)$ *version of (4.1). Let* $s' \geq s \in \mathbb{R}$*. Then*

$$\tfrac{1}{2}\int_{[s,s']\times Y}(|\tfrac{\partial}{\partial s}A|^2 + |\mathfrak{B}_{(A,\psi)}|^2 + 2r(|\tfrac{\partial}{\partial s}\psi|^2 + |D_A\psi|^2)) = \mathfrak{a}(\mathfrak{d}|_s) - \mathfrak{a}(\mathfrak{d}|_{s'})\,.$$

*Moreover,*

$$\tfrac{1}{2}\int_{\mathbb{R}\times Y}(|\tfrac{\partial}{\partial s}A|^2 + |\mathfrak{B}_{(A,\psi)}|^2 + 2r(|\tfrac{\partial}{\partial s}\psi|^2 + |D_A\psi|^2)) = \mathfrak{a}(\mathfrak{c}_-) - \mathfrak{a}(\mathfrak{c}_+).$$

Lemmas 4.1 and 4.2 with Lemma 2.5 have the following as a corollary: There is a constant $\kappa$ that is independent of $\mathfrak{d}$, r and $\mu$ and is such that $\mathfrak{f}_{\mathfrak{s}}(\mathfrak{c}_+) > \mathfrak{f}_{\mathfrak{s}}(\mathfrak{c}_-) - \kappa \ln r$.

The final lemma in this section speaks to the $L^2$ norms of $B_A$ and the covariant derivative of $\psi$ along the constant $s$ slices of $\mathbb{R} \times Y$. The latter is denoted by $\nabla^Y{}_A\psi$.

**Lemma 4.3**: *There exists* $\kappa \geq \pi$ *with the following significance: Suppose that* $r \geq \kappa$*, that* $\mu \in \Omega$ *has* $\mathcal{P}$*-norm less than 1. Let* $\mathfrak{d} = (A, \psi)$ *denote an instanton solution to the* $(r, \mu)$ *version of (4.1) with* $A_{\mathfrak{d}} < r^{2-z}$*. Fix a point* $s \in \mathbb{R}$*. Then*

$$\int_{[s,s+1]\times Y}(|\tfrac{\partial}{\partial s}A|^2 + |B_A|^2 + 2r|\tfrac{\partial}{\partial s}\psi|^2 + 2r|\nabla^Y_A\psi|^2) \leq \kappa r^2.$$

The proof of Lemma 4.3 is identical to its [T6] analog, Lemma 3.3 in [T6].

## b) A priori bounds on $\alpha$, $\beta$ and $B_A$ and $\frac{\partial}{\partial s}A$

The lemma that follows supplies the first of a series of a priori pointwise bounds on the size of the components of $\psi$, $B_A$ and $\frac{\partial}{\partial s}A$. The bounds in this first lemma are the fundamental ones from which all else follows. This upcoming Lemma 4.4 is the analog of Lemma 3.1 in [T6] and its proof essentially the same as that of the latter.

**Lemma 4.4**: *There exists* $\kappa \geq \pi$ *with the following significance: Fix* $r \geq \kappa$ *and fix* $\mu \in \Omega$ *with* $\mathcal{P}$*-norm bounded by* 1. *Suppose that* $\mathfrak{d} = (A, \psi)$ *is an instanton solution to the corresponding* $(r, \mu)$ *version of (4.1). Then*

- $|\alpha| \leq 1 + \kappa r^{-1}$
- $|\beta|^2 \leq \kappa r^{-1}(1 - |\alpha|^2) + \kappa^2 r^{-2}$.

The next set of bounds are for $|B_A|$ and $|\frac{\partial}{\partial s} A|$. Those stated by the next lemma are the analog of Lemma 3.2 in [T6]. The proof of the next lemma is virtually identical to the proof of the latter with Lemma 4.3 serving as the substitute for Lemma 3.3 in [T6].

**Lemma 4.5**: *There exists* $\kappa \geq \pi$ *with the following significance: Suppose that* $r \geq \kappa$ *and that* $\mu \in \Omega$ *has* $\mathcal{P}$*-norm less than 1. Suppose in addition that* $\mathfrak{d} = (A, \psi)$ *is an instanton solution to the* $(r, \mu)$ *version of (4.1) with* $A_{\mathfrak{d}} < r^2$. *Then* $|B_A| + |\frac{\partial}{\partial s} A| \leq \kappa r$.

The bound supplied by this lemma is used to prove the next one. This upcoming Lemma 4.6 is the analog of Lemma 3.6 in [T6] and its proof is identical with Lemma 4.5 serving as a substitute for Lemma 3.2 in [T6].

The notation used in Lemma 4.6 and subsquently has $\nabla_A$ denoting the covariant derivative on sections of the pull-back of E over $\mathbb{R} \times Y$ that is defined by viewing the connection A as an $\mathbb{R}$-dependent connection on this pull-back bundle. By way of an example: $\nabla_A \psi = \frac{\partial}{\partial s} \psi\, ds + \nabla^Y{}_A \psi$.

**Lemma 4.6**: *There exists* $\kappa \geq \pi$ *with the following significance: Suppose that* $r \geq \kappa$*, that* $\mu \in \Omega$ *has* $\mathcal{P}$*-norm less than 1. Suppose in addition that* $\mathfrak{d} = (A, \psi)$ *is an instanton solution to the* $(r, \mu)$ *version of (4.1) with* $A_{\mathfrak{d}} < r^2$. *Then*

- $|\nabla_A \alpha|^2 \leq \kappa r$.
- $|\nabla_A \beta|^2 \leq \kappa$.

*In addition, for each* $q \geq 1$*, there exists a constant* $\kappa_q$ *which is independent of* $\mathfrak{d}, \mu$ *and* $r$*, and is such that when* $r \geq \kappa$ *then*

- $|\nabla_A{}^q \alpha| + r^{1/2} |\nabla_A{}^q \beta| \leq \kappa_q r^{q/2}$.

The upcoming Lemma 4.7 is the analog of Lemma 3.7 in [T6]. This lemma and subsequent lemmas refer to the function $\underline{M}$ on $\mathbb{R}$ that is defined by the rule.

$$s \to \underline{M}(s) = r \int_{[s-1, s+1] \times Y} (1 - |\alpha|^2) \,. \tag{4.3}$$

The proof of the upcoming lemma differs little from that of Lemma 3.7 in [T6] with Lemma 4.5 serving in lieu of Lemma 3.2 in [T6].

**Lemma 4.7**: *There exists* $\kappa \geq \pi$ *with the following significance: Fix* $r \geq \kappa$ *and* $\mu \in \Omega$ *with* $\mathcal{P}$*-norm less than 1. Let* $\mathfrak{d} = (A, \psi)$ *denote an instanton solution to the* $(r, \mu)$ *version of (4.1) with* $A_{\mathfrak{d}} < r^2$. *Assume in addition that* $s_0 \in \mathbb{R}$ *and that* $\mathcal{K} \geq 1$ *are such that* $\sup_{s\in[s_0-2, s_0+2]} \underline{M}(s) \leq \mathcal{K}$. *Then*

$$|\tfrac{\partial}{\partial s} A - B_A| \leq r(1 + \kappa\mathcal{K}^{1/2} r^{-1/2})(1 - |\alpha|^2) + \kappa$$

*at all points in* $[s_0 - 1, s_0 + 1] \times Y$.

The upcoming Lemma 4.8 is a refinement of Lemma 3.8 in [T6] in that it makes no reference to $\underline{M}$. The proof given below works just as well in the context of Lemma 3.8 in [T6] and so the assertion of the latter lemma holds also with no reference to $\underline{M}$.

**Lemma 4.8**: *There exists* $\kappa \geq \pi$ *with the following significance: Fix* $r \geq \kappa$ *and* $\mu \in \Omega$ *with* $\mathcal{P}$*-norm less than 1. Let* $\mathfrak{d} = (A, \psi)$ *denote an instanton solution to the* $(r, \mu)$ *version of (4.1) with* $A_{\mathfrak{d}} < r^2$. *Fix* $s_0 \in \mathbb{R}$ *and let* $X_* \subset [s_0 - 2, s_0 + 2] \times Y$ *denote the set of points where* $1 - |\alpha| \leq \kappa^{-1}$. *The bounds stated below hold on* $X_*$.

- $|\nabla_A\alpha|^2 + r|\nabla_A\beta|^2 \leq \kappa r(1 - |\alpha|^2) + \kappa^2$
- $r(1 - |\alpha|^2) + |\nabla_A\alpha|^2 + r|\nabla_A\beta|^2 \leq \kappa(r^{-1} + r e^{-\sqrt{r}\,\mathrm{dist}(\cdot, X_*)/\kappa})$
- $|\beta|^2 \leq \kappa(r^{-2} + r^{-1} e^{-\sqrt{r}\,\mathrm{dist}(\cdot, X_*)/\kappa})$

***Proof of Lemma 4.8***: The proof of the top bullet is the same as the proof of the analog in Proposition 2.8 of the article *SW=> Gr* in [T8]. It uses only the bounds from Lemma 4.6 on $|B_A|$ and $|\frac{\partial}{\partial s} A|$. The bounds in the second and third bullets of Lemma 4.8 are derived in the three steps that follow.

<u>Step 1</u>: Mimic what is done in Step 2 of the proof of Proposition 4.4 in the article *SW=>Gr* from [T8] to find positive, $(A,\psi)$ and r independent constants $c > 1$ and $z_1$, $z_2$ such that function $y_1 = (|\nabla\alpha|^2 + z_1 r|\nabla\beta|^2 + z_2 r^2|\beta|^2)$ obeys a differential inequality of the form

$$d{*}d y_1 + c^{-2} r\, y_1 \leq c_0 r(1 - |\alpha|^2) y_1 + c_0. \tag{4.4}$$

Here d is used to denote the exterior derivative on $\mathbb{R} \times Y$ and $d^*$ is its formal $L^2$ adjoint. By definition, $\underline{M} \le c_0 r$ and (4.4) invokes the corresponding $\mathcal{K} = c_0 r$ version of Lemma 4.7 with the first bullet in (4.1) to bound $|B_A|$ and $|\frac{\partial}{\partial s} A|$.

Step 2: Fix $x_0 \in X_*$ and let $d_0$ denote the distance from $x_0$ to the boundary of $X_*$. There exists $c_0 > 1$ such that the function $x \to h_0(x) = e^{\sqrt{r}(\mathrm{dist}(x,x_0) - d_0)/2c}$ obeys the differential inequality $d^*dh_0 + \frac{1}{2} c^{-2} r h_0 \le 0$ when $\mathrm{dist}(x,x_0) \le c_0^{-1}$. Lemma 4.4 and the top bullet in Lemma 4.8 bound $y_1$ by $c_0 r$ in any event, and so $y_2 = y_1 - c_0(r h_0 + r^{-1})$ is non-positive on where $d(x, x_*) \ge d_0$. Meanwhile, (4.4) implies that $d^*dy_2 \le 0$ if $X_*$ is defined to be where $1 - |\alpha|^2 \le \frac{1}{2} c_0^{-1} c^{-2}$. Granted this definition, then the maximum principle asserts that $y_2 \le 0$ on $X_*$. In particular, this is true at $x_0$ and so $y_1 \le c_0(r\, e^{-\sqrt{r} d_0/2c} + r^{-1})$. This proves the third bullet in Lemma 4.8 and the assertions in the second bullet of Lemma 4.8 concerning $|\nabla\alpha|^2$ and $|\nabla\beta|^2$.

Step 3: Mimic the derivation of Equation (6.5) in Section 6a of [T1] to see that the function $w = 1 - |\alpha|^2$ obeys a differential inequality of the form

$$d^*dw + 2rw \le 2rw^2 + c_0(y_1 + 1) \tag{4.5}$$

Granted (4.5), then Step 2's maximum priniciple argument using $h_0$ can be repeated with only cosmetic changes to prove the assertion about $1 - |\alpha|^2$ in the lemma's second bullet.

Lemmas 4.9 and 4.10 are the respective analogs of Lemmas 3.9 and 3.10 of [T6]. To set the stage for these lemmas, suppose that $x \in \mathbb{R} \times Y$ and $\rho \in (r^{-1/2}, c_0^{-1})$ have been specified. The lemmas use $M_{(x,\rho)}$ to denote the integral of $r(1 - |\alpha|^2)$ over the radius $\rho$ ball in $\mathbb{R} \times Y$ centered at x.

**Lemma 4.9**: *There exists $\kappa \ge \pi$, and given data consisting of an open set $U \subset \mathbb{R} \times Y$, an open subset $V \subset U$ with compact closure and $\mathcal{K} \ge 1$, there exists $\kappa_{\mathcal{K},U,V} \ge 1$ with the following significance: Fix $r \ge \kappa$ and $\mu \in \Omega$ with $\mathcal{P}$-norm less than 1. Let $\mathfrak{d} = (A, \psi)$ denote an instanton solution to the $(r, \mu)$ version of (4.1) with $A_{\mathfrak{d}} < r^2$. Assume that $\mathfrak{d}$ is such that $\sup_{x\in U} \underline{M}_{x,1/\kappa} \le \mathcal{K}$. Then*

- $|\frac{\partial}{\partial s} A + B_A| \le r(1 - |\alpha|^2) + \kappa_{\mathcal{K},U,V}$,
- $|\frac{\partial}{\partial s} A - B_A| \le r(1 - |\alpha|^2) + \kappa_{\mathcal{K},U,V}$

*at all points in V.*

The proof of Lemma 4.9 is very similar to that of its analog in [T6], the latter being very similar to the proof of Proposition 3.4 in the article *SW=> Gr* from [T8].

The final lemma in this subsection is a monotonicity result and of a different flavor from the pointwise bounds given above. It plays a role in the proof of Lemma 4.9.

**Lemma 4.10**: *There exists* $\kappa \geq \pi$, *and given* $z \geq 1$, *there exists* $\kappa_z \geq 1$ *with the following significance: Fix* $r \geq \kappa$ *and* $\mu \in \Omega$ *has* $\mathcal{P}$*-norm less than 1. Let* $\mathfrak{d} = (A, \psi)$ *denote an instanton solution to the* $(r, \mu)$ *version of (4.1) with* $A_{\mathfrak{d}} < r^2$ *and* $\sup_{s\in\mathbb{R}} \underline{M}(s) \leq r^{1-1/z}$. *Given* $x \in \mathbb{R}\times Y$ *and* $\rho \in (r^{-1/2}, \kappa_z^{-1})$, *use* $M_{(x,\rho)}$ *to denote the integral of* $r(1 - |\alpha|^2)$ *over the radius* $\rho$ *ball in* $\mathbb{R}\times Y$ *centered at* $x$. *Then*

- *If* $\rho_1 > \rho_0$ *are in* $(r^{-1/2}, \kappa_z^{-1})$, *then* $M_{(x,\rho_1)} > \kappa_z^{-1} \rho_1^2/\rho_0^2\, M_{(x,\rho_0)}$.
- *Suppose that* $|\alpha| \leq \frac{3}{4}$ *at* $x$. *If* $\rho \in (r^{-1/2}, \kappa_z^{-1})$, *then* $M_{(x,\rho)} \geq \kappa^{-1} \rho^2$.
- *Suppose that* $\mathcal{K} \in (1, r^{1-1/z})$ *and suppose that* $d \in (r^{-1/2}, \kappa^{-1})$ *and* $x \in \mathbb{R}\times Y$ *are such that* $M_{(x,d)} \leq \mathcal{K} d^2$. *If* $\rho \in (r^{-1/2}, d)$, *then* $M_{(x,\rho)} \leq \kappa_z \mathcal{K} \rho^2$.

As with the proof of Lemma 3.10 in [T6], the proof of Lemma 4.10 differs little from the proof of Proposition 3.1 in the article *SW=> Gr* from [T8]. This lemma also plays a role in the subsequent sections.

## c) Instantons and holomorphic data on $\mathbb{C}$

The three parts of this section first introduce holomorphic notions on $\mathbb{C}^2$, and then explain how they model an instanton solution to (4.1) in a radius $\mathcal{O}(r^{-1/2})$ ball.

*Part 1*: This part introduces the relevant holomorphic data on $\mathbb{C}^2$. To this end, introduce complex coordinates $(x_0, x_1)$ for $\mathbb{C}^2 = \mathbb{R}^4$. Give $\mathbb{C}^2$ the standard metric with Kahler form $\omega_0 = \frac{i}{2}(dx_0 \wedge d\bar{x}_0 + dx_1 \wedge d\bar{x}_1)$. Use $P^+ \colon \wedge^2 T^*\mathbb{C}^2 \to \wedge^2 T^*\mathbb{C}^2$ to denote the projection to the self dual subspace and $P^-$ the projection to the anti-self dual subspace.

Of interest here are pairs $(A_0, \alpha_0)$ on $\mathbb{C}^2$ where $A_0$ is a unitary connection on the trivial bundle and $\alpha_0$ is a section of this bundle; and where these are such that

- $\bar{\partial}_{A_0} \alpha_0 = 0$.
- $P^+ F_{A_0} = -\frac{i}{2}(1 - |\alpha_0|^2)\,\omega_0$.
- $|\alpha_0| \leq 1$.
- $|P^- F_{A_0}| \leq |P^+ F_{A_0}| \leq 2^{-1/2}(1 - |\alpha_0|^2)$.

(4.6)

Proposition 4.1 in the article *SW=> Gr* from [T8] and Proposition 4.2 in [T6] describe the pairs $(A_0, \alpha_0)$ that satisfy these conditions. Except for the second bullet, the following proposition restates Proposition 4.2 in [T6]. The proof of the second bullet is the same as that of the second bullet of this same Proposition 4.2 in [T6].

**Proposition 4.11**: *Suppose that* $(A_0, \alpha_0)$ *obeys (4.6).*

- *If* $|\alpha_0| < 1$ *somewhere, then* $|\alpha_0|$ *is strictly less than 1, it has no positive local minimum and* $\inf_{\mathbb{C}^2} |\alpha_0| = 0$. *If* $\alpha_0^{-1}(0) \neq \emptyset$, *then* $\alpha_0^{-1}(0)$ *is either all of* $\mathbb{C}^2$ *or a complex analytic subvariety of complex dimension 1.*
- *There exists* $\kappa_0 > 1$ *that is independent of* $(A_0, \alpha_0)$ *and has the following signficance: Let* $X_* \subset \mathbb{C}^2$ *denote the set of points where* $1 - |\alpha_0| \geq \frac{3}{4}$. *Then*

$$1 - |\alpha_0| + |\nabla_{A_0}\alpha_0| \leq \kappa_0 e^{-\mathrm{dist}(\cdot, X_*)/\kappa_0} .$$

- *If* $|\alpha_0| < 1$ *somewhere, and if there exists* $m \geq 1$ *such that the integral of* $(1 - |\alpha_0|^2)$ *over the ball of any given radius* $R \geq 1$ *centered at the origin is less than* $mR^2$, *then*
  a) *The locus* $\alpha^{-1}(0)$ *is a non-empty, complex algebraic subvariety with complex dimension 1. As such, this locus near any given point is the zero locus of a holomorphic polynomial*
  b) *The order of the latter polynomial has a purely* $m$*-dependent upper bound.*

  *If, in addition, the integral over* $\mathbb{C}^2$ *of* $|P^+F_{A_0}|^2 - |P^-F_{A_0}|^2$ *is finite, then*
  c) *This integral is a non-negative integer multiple of* $4\pi^2$.
  d) *If the latter integral is zero, then* $(A_0, \alpha_0)$ *is the pull-back via a projection* $\mathbb{C}^2 \to \mathbb{C}$ *of a solution on* $\mathbb{C}$ *to the vortex equations in (1.4) and* $\alpha_0^{-1}(0)$ *is a union of planes.*
- *The set of gauge equivalence classes of pairs* $(A_0, \alpha_0)$ *that obey (4.1) is sequentially compact with respect to convergence on compact subsets of* $\mathbb{C}^2$ *in the* $C^\infty$ *topology.*

The solutions to (4.6) constitute the desired holomorphic data on $\mathbb{C}^2$.

*Part 2*: Fix a point $p \in \mathbb{R} \times Y$. A *complex Gaussian* coordinate system centered at p is a coordinate chart map from a ball about the origin in $\mathbb{C}^2$ to a neighborhood of p that takes the origin to x and defines a Gaussian coordinate chart when written in terms of the real coordinates for $\mathbb{C}^2$. In addition, the almost complex structure J at the point p must appear in these coordinates as the standard complex structure. The complex coordinates on $\mathbb{C}^2$ are written again as $(x_0, x_1)$. No generality is lost by assuming that any given such Gaussian coordinate chart is defined where $|x_0|^2 + |x_1|^2 \leq c_0^{-1}$ with $c_0$ being independent of p.

Introduce a new coordinate chart by composing the original one with the map from $\mathbb{C}^2$ to itself that sends $(x_0, x_1)$ to $(r^{-1/2}x_0, r^{-1/2}x_1)$. The new coordinate chart is defined on the ball of radius $c_0^{-1}r^{1/2}$ centered at the origin in $\mathbb{C}^2$. Use $\varphi_r$ in what follows to denote this coordinate chart map from the ball of radius $c_0^{-1}r^{1/2}$ in $\mathbb{C}^2$ to $\mathbb{R} \times Y$.

The $\varphi_r$ pull-back of the metric from $\mathbb{R} \times Y$ differs from the standard Euclidean metric by no more than $c_0 r^{-1}$ on the radius $2^4$ ball. The pull-back of the Riemannian curvature is also bounded in absolute value on this ball by $c_0 r^{-1}$, and the latter's derivatives to a given order $k \geq 1$ on this ball have norm bounded by $c_k r^{-1-k/2}$ with $c_k$ depending on k only. Meanwhile, the $\varphi_r$ pull-back of the almost complex structure on this ball differs from the standard one by at most $c_0 r^{-1/2}$ and its derivatives to order k have norm bounded by $c_k r^{-(1+k)/2}$.

*Part 3*: Let $\mathfrak{d} = (A, \psi = (\alpha, \beta))$ denote an instanton solution to (4.1) with $A_{\mathfrak{d}} \leq r^2$ and such that there exists $z > 1$ such that $\sup_{s\in\mathbb{R}} \underline{M}(s) \leq r^{1-1/z}$. Introduce $(A_r, \alpha_r)$ to denote the $\varphi_r$ pull-back of the $(A, \alpha)$. Use $F_{A_r}$ to denote the curvature 2-form of the connection $A_r$. Lemmas 4.4, 4.6 and 4.7 have implications with regards to $(A_r, \alpha_r)$ that are described in what follows. To say more, fix $R \geq 1$. Given Part 2's remarks about about the $\varphi_r$ pull-backs of the metric and almost complex structure, there exists $c_R > 1$ that is independent of p and such that if $r \geq c_R$, then the $\mathfrak{d}$ version of $(A_r, \alpha_r)$ is nearly a solution to (4.6) on the ball of radius R in $\mathbb{C}^2$ centered at the origin in the sense that

- $|\bar{\partial}_{A_r} \alpha_r| \leq c_R r^{-1/2}$.
- $|P^+ F_{A_r} + \frac{i}{2}(1 - |\alpha_r|^2)\omega_0| \leq c_R r^{-1}$.
- $|\alpha_r| \leq 1 + c_R r^{-1}$
- $|P^- F_{A_r}| \leq |P^+ F_{A_r}| + c_R r^{-1/2z} \leq 2^{-1/2}(1 - |\alpha_r|^2) + c_R r^{-1/2z}$.

(4.7)

Moreover, with the $\varphi_r$ pull-back of $\beta$, the pair $(A_r, \alpha_r)$ plus $\varphi_r{}^*\beta$ obey an equation on the radius R ball in $\mathbb{C}^2$ that gives bounds on the covariant derivatives of $\alpha_r$ and $F_{A_r}$ to any given order that are in dependent of p, $\mathfrak{d}$ and R. These bounds with (4.7) lead to the following lemma.

**Lemma 4.12**: *Given* $q \geq 1$, $R > 1$, $\varepsilon > 0$, $k \in \{1, 2, \ldots\}$ *and* $m > 1$, *there exists* $\kappa > 10R$ *with the following significance: Fix* $r \geq \kappa$ *and* $\mu \in \Omega$ *with* $\mathcal{P}$*-norm less than 1 and suppose that* $(A, \psi)$ *denotes an instanton solution to (4.1) with* $A_{\mathfrak{d}}$ *and* $\sup_{s\in\mathbb{R}} \underline{M}(s) \leq r^{1-1/q}$. *Given* $p \in \mathbb{R} \times Y$, *there exists a solution to* (4.6) *on* $\mathbb{C}^2$, *this denoted by* $(A_0, \alpha_0)$, *such that* $(A_r, \alpha_r) = (A_0 + \hat{a}, \alpha_0 + \eta)$ *with* $(\hat{a}, \eta)$ *having* $C^k$ *norm less than* $\varepsilon$ *on the ball of radius* R *in* $\mathbb{C}^2$ *centered at the origin. Moreover, suppose that the integral of* $r(1 - |\alpha|^2)$ *on each* *radius* $\rho \in (r^{-1/2}, \kappa r^{-1/2})$ *ball centered on* p *is less than* $m\rho^2$. *Then* $(A_0, \alpha_0)$ *can be chosen so as to obey Items a) and b) of the third bullet in Proposition 4.12.*

Lemma 4.12 is the analog here of Lemma 4.3 in [T6]. As with the latter, the proof differs little from that of Proposition 4.2 in the article *SW => Gr* from [T8].

## 5. A priori bounds for the function $\underline{M}$: The complement of $\cup_{p\in\Lambda}(\hat{\gamma}_p^+ \cup \hat{\gamma}_p^-)$

Write $\upsilon_\Diamond$ as $q_\Diamond \hat{a} + \mathcal{b}$ where the 1-form $\mathcal{b}$ annihilates $v$. By way of a reminder, the function $q_\Diamond$ differs from 1 only in $\cup_{p\in\Lambda}\mathcal{H}_p$, it vanishes only on $\cup_{p\in\Lambda}(\hat{\gamma}_p^+ \cup \hat{\gamma}_p^-)$, and it is such that $q_\Diamond \geq c_0^{-1}|\upsilon_\Diamond|^2$. Fix $r > c_0$ and $\mu \in \Omega$ with $\mathcal{P}$-norm less than 1 so as to define (4.1). Suppose that $\mathfrak{d} = (A, \psi = (\alpha, \beta))$ is an instanton solution to this $(r, \mu)$ version of (4.1). This section supplies a $\mathfrak{d}$ and r independent bound for the function on $\mathbb{R}$ given by the rule

$$s \to \underline{M}_\Diamond(s) = r \int_{[s-1,s+1]\times Y} q_\Diamond^6 (1 - |\alpha|^2) \ . \tag{5.1}$$

The proposition that follows makes a formal statement that such a bound exists.

**Proposition 5.1**: *There exists* $\kappa \geq \pi$ *and given* $c \geq 1$, *there exists* $\kappa_c > 1$ *with the following significance: Suppose that* $r \geq \kappa$ *and that* $\mu \in \Omega$ *has* $\mathcal{P}$*-norm less than 1. Suppose in addition that* $\mathfrak{d} = (A, \psi)$ *is an instanton solution to the* $(r, \mu)$ *version of (4.1) with* $A_\mathfrak{d} < c r$. *The corresponding function* $\underline{M}_\Diamond$ *obeys* $-\kappa < \underline{M}_\Diamond < \kappa_c$.

The lower bound follows directly from Lemma 4.4 so it holds without the bound for $A_\mathfrak{d}$. The proof of the upper bound occupies the rest of this section. By way of a parenthetical remark, the proof looks much like the proof of Lemma 5.8 in [LT].

### a) Preliminary bounds for $\underline{M}_\Diamond$ and $\underline{M}$

The lemma that follows supplies a preliminary and easy to come by bound for $\underline{M}$ that is used in the later subsections to invoke Lemma 4.7.

**Lemma 5.2**: *There exists* $\kappa \geq \pi$ *with the following significance: Fix* $r \geq \kappa$ *and* $\mu \in \Omega$ *with* $\mathcal{P}$*-norm less than 1. Let* $\mathfrak{d} = (A, \psi)$ *denote an instanton solution to the* $(r, \mu)$ *version of (4.1). The corresponding versions of* $\underline{M}_\Diamond$ *obeys* $-\kappa \leq \underline{M}_\Diamond \leq \kappa(A_\mathfrak{d} + 1)^{1/2}$ *and the corresponding version of* $\underline{M}$ *obeys* $-\kappa < \underline{M} < \kappa r^{2/3}(1 + A_\mathfrak{d})^{1/6}$.

***Proof of Lemma 5.2***: The lower bounds follow from Lemma 4.4. The first step of what follows establishes the upper bound for $\underline{M}_\Diamond$ and the second step establishes the upper bound for $\underline{M}$. The notation in these steps is that used earlier in the proof of the second bullet of Lemma 2.5.

Step 1: To prove the upper bound for $\underline{M}_\Diamond$, take the inner product on Y between $\upsilon_\Diamond$ and the 1-form on the right hand side of the top bullet in (4.1). Integrate the result over $[s-1,s+1] \times Y$. This integral is, of course, equal to zero. Thus,

$$\int_{[s-1,s+1]\times Y} (\upsilon_\Diamond \wedge * \tfrac{\partial}{\partial s} A) + \int_{[s-1,s+1]\times Y} \upsilon_\Diamond \wedge *B_A = r \int_{[s-1,s+1]\times Y} \upsilon_\Diamond \wedge *(\psi^\dagger \tau \psi - i\hat{a}) + \mathfrak{e} \tag{5.2}$$

where $|\mathfrak{e}| \leq c_0$. Write $\upsilon_\Diamond$ as $q_\Diamond \hat{a} + \mathfrak{b}$ with $\mathfrak{b}$ annihilating the vector field $v$ and use this rewriting for the integrand of the integral on the right hand side of (5.2). Then, use the bounds $|\mathfrak{b}| \leq |\upsilon_\Diamond|$ and $|\upsilon_\Diamond| \leq c_0 q_\Diamond^{1/2}$ with Lemma 4.4's bounds for $|\beta|$ to see that this integrand is greater than $\frac{1}{2}|\upsilon_\Diamond|^2(1-|\alpha|^2) - c_0 r^{-1}$. This the case, a bound on the integral on the left hand side of (5.2) supplies one for the integral on the right hand side of (5.1).

To obtain an upper bound for the left hand side of (5.2), use Lemma 4.2 to see that the integral of $*\frac{\partial}{\partial s}A$ that appears on the left hand side of (5.2) is no greater than $c_0(1+A_\eth)^{1/2}$. Meanwhile, the integral of $\upsilon_\Diamond \wedge *B_A$ is independent of A and r because it computes a pairing with the first Chern class of the bundle E. These last facts imply that the left hand side of (5.2) is no greater than $c_0(A_\eth + 1)^{1/2}$.

Step 2: Fix $\rho > 0$ and let $Y^\rho$ denote for the moment the set of points in Y with distance $\rho$ or more from the curves in the set $\cup_{p\in\Lambda}\{\hat{\gamma}_p^+ \cup \hat{\gamma}_p^-\}$. The integral in (5.1) is no less then the contribution from $Y^\rho$ and this is no less than $c_0^{-1}\rho^4(\underline{M} - c_0 r\rho^2) - c_0$. It follows as a consequence that $\underline{M} \leq c_0(\rho^{-4}(A_\eth+1)^{1/2} + r\rho^2)$. This understood, take $\rho^2 = r^{-1/3}(A_\eth+1)^{1/6}$ to obtain what is asserted by Lemma 5.2.

**b) A vortex-like inequality**

This subsection shows how Proposition 5.1 follows from Lemma 5.3. This lemma asserts an inequality that is reminiscent of the equality asserted by the first bullet of Lemma 3.3.

Lemma 5.3 refers to a certain function, $Q_\Diamond$, on Y which is specified in the next subsection. For the purposes of the lemma, it is enough to know that $Q_\Diamond \geq c_0^{-1} q_\Diamond^6$ and that $|dQ_\Diamond| \leq c_0$. Given $s \in \mathbb{R}$, this lemma uses $\chi_s$ to denote the function $\chi(2|s - (\cdot)| - 1)$ on $\mathbb{R}$. This function is 1 on $[s - \frac{1}{2}, s + \frac{1}{2}]$ and it is equal to 0 on the complement of $[s-1, s+1]$.

**Lemma 5.3**: *There exists* $\kappa > \pi$ *and given* $c \geq 1$, *there exists* $\kappa_c \geq 1$ *with the following significance: Fix* $r \geq \kappa$ *and* $\mu \in \Omega$ *with* $\mathcal{P}$*-norm less than 1. Let* $\mathfrak{d} = (A, \psi)$ *denote an instanton solution to the* $(r, \mu)$ *version of (4.1) with* $A_{\mathfrak{d}} < cr$. *Then*

$$\sup_{s\in\mathbb{R}} \underline{M}_{\Diamond} \leq \kappa_c \sup_{s\in\mathbb{R}} \int_{\mathbb{R}\times Y} \chi_s Q_{\Diamond}(2r\,|\,\alpha\,|^2\,(1 - |\,\alpha\,|^2) - |\nabla_A\alpha\,|^2) + \kappa_c^{\,2}\,.$$

This lemma is proved momentarily.

***Proof of Proposition 5.1***: The proof assumes that Lemma 5.3 is true so as to deduce a suitable upper bound for $\underline{M}_{\Diamond}$. To deduce such a bound from Lemma 5.3, introduce by way of notation $\mathcal{D}_A$ to denote $\frac{\partial}{\partial s} + D_A$, this being an operator on the space of sections of $\mathbb{S}$ over $\mathbb{R} \times Y$. Use $\mathcal{D}_A^{\dagger}$ to denote its formal $L^2$ adjoint. Given that $\mathcal{D}_A\psi = 0$, so $\mathcal{D}_A^{\dagger}\mathcal{D}_A\psi$ is also zero. Projecting the equation $\mathcal{D}_A^{\dagger}\mathcal{D}_A\psi = 0$ to the E summand of $\mathbb{S}$ gives an equation of the form

$$\tfrac{1}{4}\nabla_A^{\dagger}\nabla_A\alpha + \tfrac{1}{2}\,r\,\alpha(|\alpha|^2 - 1) + \mathfrak{r} = 0\,, \tag{5.3}$$

where $|\mathfrak{r}| \leq c_0(|\alpha| + |\beta| + |\nabla_A\beta|)$. Take the inner product of this equation with $\alpha$ to find an equation of the form

$$\tfrac{1}{4}d^{\dagger}d(1 - |\alpha|^2) + r|\alpha|^2(1 - |\alpha|^2) - \tfrac{1}{2}|\nabla_A\alpha|^2 + \mathfrak{e} = 0\,, \tag{5.4}$$

where $|\mathfrak{e}| \leq c_0(|\alpha|^2 + |\beta|^2 + |\nabla_A\beta|^2)$. Multiply both sides of this last equation by $\chi_s Q_{\Diamond}$ and integrate the result over $\mathbb{R} \times Y$. Integrate by parts and appeal to Lemma 4.4 and the bound on $|\nabla_A\beta|^2$ from Lemma 4.6 to see that the integral that appears on the right hand side of Lemma 5.3 has an $(A,\psi)$ and r independent upper bound.

**c) Proof of Lemma 5.3**

The four steps that follow derive Lemma 5.3 from the upcoming Lemma 5.4. The rest of the subsection supplies a proof of Lemma 5.4.

Step 1: This step specifies the function $Q_{\Diamond}$. To do this, first introduce the function $f_*$ that is defined on each $\mathfrak{p} \in \Lambda$ version of $\mathcal{H}_{\mathfrak{p}}$ by the rule $f_*(u,\theta) = g(u)(1 - 3\cos^2\theta)$ with g defined in the third bullet in (1.2). By way of a reminder from (1.5), this function is such that $\upsilon_{\Diamond} = df_*$ on $\mathcal{H}_{\mathfrak{p}}$. Choose a smooth, non-increasing function on $[0, \infty)$ with the properties listed next. This function is denoted by Q. It is such that $Q(t) = t^5$ for $t \in (0, \frac{1}{2}]$ and $Q(t) = 1$ for $t \geq 1$. With Q in hand, fix for the moment $\varepsilon \in (0, 1)$ and use $Q_\varepsilon$ to denote the function $Q(\varepsilon^{-2}q_\Diamond)$. Let $\nu(Q_\varepsilon)$ denote the pairing between $\nu$ and $dQ_\varepsilon$. The function $Q_\Diamond$ is the function $Q_\varepsilon q_\Diamond + f_*\nu(Q_\varepsilon)$ for a choice for $\varepsilon$ that guarantees it to be greater than $c_0^{-1}q_\Diamond^{\,6}$ and to have derivative norm bounded by $c_0$. This choice is such that $\varepsilon > c_0^{-1}$.

Step 2: The upcoming Equation (5.5) supplies an integral form of the Bochner-Weitzenboch identity for the operator $D_A(q_\Diamond D_A)$. The formula reintroduces from (1.11) the Clifford multiplication endomorphism $cl(\cdot)$. This formula is derived using integration by parts. Suppose for the moment that $(A,\psi)$ is any given pair in $Conn(E) \times C^\infty(Y;\mathbb{S})$. What follows is the promised identity.

$$\int_Y Q_\Diamond(|B_A|^2 + r^2|\psi^\dagger\tau\psi - i\hat{a}|^2 + r|\nabla_A^Y\psi|^2) = + 2ir\int_Y Q_\Diamond \hat{a} \wedge {*B_A} - r\int_Y (\psi^\dagger cl(dQ_\Diamond)D_A\psi - (D_A\psi)^\dagger cl(dQ_\Diamond)\psi) + \int_Y Q_\Diamond(|\mathfrak{B}_{(A,\psi)}|^2 + r|D_A\psi|^2) + \mathfrak{e}\,, \tag{5.5}$$

where $\mathfrak{e}$ obeys $|\mathfrak{e}| \leq c_0(1 + r)$. The proof of Lemma 5.3 uses the a priori bounds given by the next lemma on the first two integrals that appear on the right hand side of (5.5). Lemma 4.2 is used to bound the third, right most integral on the right hand side of (5.5).

**Lemma 5.4**: *There exists $\kappa > \pi$ and given $c \geq 1$, there exists $\kappa_c > 1$ with the following significance: Fix $r \geq \kappa$ and $\mu \in \Omega$ with $\mathcal{P}$-norm less than 1. Let $\mathfrak{d} = (A, \psi)$ denote an instanton solution to the $(r, \mu)$ version of (4.1) with $A_\mathfrak{d} < cr$ and $\sup_{s\in\mathbb{R}} \underline{M}_\Diamond \geq 1$.*

- $\sup_{s\in\mathbb{R}} 2r \int_{\mathbb{R}\times Y} \chi_s(i\int_Y Q_\Diamond \hat{a} \wedge {*B_A}) \leq \frac{1}{1000} r \sup_{s\in\mathbb{R}} \underline{M}_\Diamond + \kappa_c r$
- $\sup_{s\in\mathbb{R}} r|\int_{\mathbb{R}\times Y} \chi_s(\psi^\dagger cl(dQ_\Diamond)D_A\psi - (D_A\psi)^\dagger cl(dQ_\Diamond)\psi)| \leq \frac{1}{1000} r \sup_{s\in\mathbb{R}} \underline{M}_\Diamond + \kappa_c r$.

This lemma is proved momentarily. The remaining steps use Lemma 5.4 to complete the argument for Lemma 5.3.

Step 3: Take the $Conn(E) \times C^\infty(Y;\mathbb{S})$ pair $(A, \psi)$ in (5.5) to be the pair given in the statement of Lemma 5.4 at any given slice of $\mathbb{R} \times Y$ with constant $\mathbb{R}$ factor. Add the integral over this slice of $Q_\Diamond|\frac{\partial}{\partial s}\psi|^2$ to both sides of (5.5). View the result as an equality between functions on $\mathbb{R}$. Multiply this equality by $\chi_s$ and integrate over $\mathbb{R}$. Then use Lemmas 4.2 and 5.4 with (5.5) to see that

$$\int_{\mathbb{R}\times Y}\chi_s Q_\Diamond(|B_A|^2+r^2(1-|\alpha|^2)^2+r|\nabla_A\alpha|^2) \le \tfrac{1}{100}r\sup_{s\in\mathbb{R}} \underline{M}_\Diamond + c_{c*}r \ , \tag{5.6}$$

where $c_{c*}$ denotes the version of $\kappa_c$ that is given by Lemma 5.4. To make something of this, mimic what is done in Section 5.4 of [LT] by writing

$$\tfrac{\partial}{\partial s}A = -i(1-\sigma)(r(1-|\alpha|^2)+\mathfrak{z}_A)\hat{a} + \mathfrak{r} + \mathfrak{X} \quad \textit{and} \quad B_A = -i\sigma(r(1-|\alpha|^2)+\mathfrak{z}_B)\hat{a} + \mathfrak{r} - \mathfrak{X} \tag{5.7}$$

where the notation uses $\sigma$ to denote a function on $\mathbb{R}\times Y$. The notation has $\mathfrak{z}_A$ and $\mathfrak{z}_B$ denoting functions on $\mathbb{R}\times Y$ with norms bounded by 1, and it has both $\mathfrak{r}$ and $\mathfrak{X}$ annihilating $\nu$. Lemma 4.4 finds $|\mathfrak{r}| \le c_0(r^{1/2}|1-\alpha|^{1/2}+1)$. To say something more about $\mathfrak{X}$, use the top bullet in (4.1) and Lemma 4.7 with Lemma 5.2 to see that

$$4|\mathfrak{X}|^2 + (1-2\sigma)^2r^2(1-|\alpha|^2)^2 \le r^2(1+c_c r^{-1/12})^2(1-\alpha)^2 + c_c \ , \tag{5.8}$$

where $c_c$ here and in what follows denotes denotes a purely $c$ dependent constant with value greater than 1. The notation is such that $c_c$ increases between subsequent appearances. This last inequality implies that

$$|\mathfrak{X}|^2 \le r^2\sigma(1-\sigma)(1-|\alpha|^2)^2 + c_c r^{23/12}(1-|\alpha|^2)^2 + c_c. \tag{5.9}$$

Use (5.7) to write

- $\int_{\mathbb{R}\times Y}\chi_s Q_\Diamond|\tfrac{\partial}{\partial s}A|^2 = \int_{\mathbb{R}\times Y}\chi_s Q_\Diamond((1-\sigma)^2r^2(1-|\alpha|^2)^2+|\mathfrak{X}|^2) + \mathfrak{e}_A$ ,
- $\int_{\mathbb{R}\times Y}\chi_s Q_\Diamond|B_A|^2 = \int_{\mathbb{R}\times Y}\chi_s Q_\Diamond(\sigma^2r^2(1-|\alpha|^2)^2+|\mathfrak{X}|^2) + \mathfrak{e}_B$ ,

(5.10)

where $\mathfrak{e}_A$ and $\mathfrak{e}_B$ are such that $|\mathfrak{e}_A| + |\mathfrak{e}_B| \le \tfrac{1}{1000}r\sup_{s\in\mathbb{R}} \underline{M}_\Diamond + c_c r$.

Step 3: Let $s \in \mathbb{R}$ denote a point where the function of $\underline{M}_\Diamond$ is greater than $\tfrac{3}{4}$ times its supremum. Following along the lines of what is done in Section 5.4 of [LT], consider two cases, that when

$$\int_{\mathbb{R}\times Y}\chi_s Q_\Diamond|\tfrac{\partial}{\partial s}A|^2 \ge \tfrac{1}{100}\int_{\mathbb{R}\times Y}\chi_s Q_\Diamond r^2(1-|\alpha|^2)^2 \tag{5.11}$$

and otherwise. If (5.11) holds, add 100 times the right hand integral in (5.11) to both sides of (5.6) and invoke Lemma 4.2 to bound the resulting contribution to the right hand side. Doing so leads to the inequality

$$\int_{\mathbb{R}\times Y}\chi_s Q_\Diamond(2r^2(1-|\alpha|^2)^2+r|\nabla_A\alpha|^2) \le \tfrac{11}{100}\, r\sup_{s\in\mathbb{R}}\underline{M}_\Diamond + c_c r\ . \tag{5.12}$$

Write the left hand side of this inequality as

$$2r^2\int_{\mathbb{R}\times Y}\chi_s Q_\Diamond(1-|\alpha|^2) + r^2\int_{\mathbb{R}\times Y}\chi_s Q_\Diamond(-2|\alpha|^2(1-|\alpha|^2)+r|\nabla_A\alpha|^2)\ . \tag{5.13}$$

The left most integral in (5.13) with the factor of 2 is no less than $\frac{3}{2}\sup_{s\in\mathbb{R}}\underline{M}_\Diamond$ and so (5.11) and (5.13) imply what is asserted by Lemma 5.3.

Now suppose that (5.11) is not true. If this is so, then (5.8) and (5.10) imply that

$$\int_{\mathbb{R}\times Y}\chi_s Q_\Diamond|B_A|^2 \ge (1-\tfrac{1}{50})\int_{\mathbb{R}\times Y}\chi_s Q_\Diamond r^2(1-|\alpha|^2)^2\ . \tag{5.14}$$

Use this last inequality in (5.6) with the top bullet in Lemma 3.5 to see that (5.12) still holds. This being the case, then what is said at the end of the last paragraph can be repeated so as to complete the proof of Lemma 5.3.

***Proof of Lemma 5.4***: The first seven steps in the proof verify the top bullet and the eighth step verifies the lower bullet.

Step 1: As noted previously, $\upsilon_\Diamond = df_*$ on any given $\mathfrak{p}\in\Lambda$ version of $\mathcal{H}_\mathfrak{p}$. This understood, a homologous closed 1-form, this denoted by $\upsilon_\varepsilon$, is defined to be $\upsilon_\Diamond$ on $M_\delta\cup\mathcal{H}_0$ and to equal $d(Q_\varepsilon f_*)$ on each $\mathfrak{p}\in\Lambda$ version of $\mathcal{H}_\mathfrak{p}$. This form can be written as $\upsilon_\varepsilon = Q_\Diamond\hat{a} + \mathit{b}_\varepsilon$ with $\mathit{b}_\varepsilon$ annihilating the vector field $v$. Writing $\upsilon_\Diamond = q_\Diamond\hat{a} + \mathit{b}$ allows $\mathit{b}_\varepsilon$ to be written as $Q_\varepsilon\mathit{b} + f_* d^\perp Q_\varepsilon$ with $d^\perp$ denoting here the orthogonal projection of the exterior derivative to the annihilator of $v$ in $T^*Y$. Use the fact that $|\mathit{b}_\varepsilon|\le c_0|v_\Diamond|$ and $|v_\Diamond|\le c_0 q_\Diamond^{1/2}$ and that $|f_*|\le c_0 q_\Diamond$ to see that $|\mathit{b}_\varepsilon|\le c_0 q_\Diamond^{11/2}$.

Step 2: Write $Q_\Diamond\hat{a}$ as $\upsilon_\varepsilon - \mathit{b}_\varepsilon$. The integral of $\frac{i}{2\pi}\upsilon_\varepsilon\wedge *B_A$ computes the cup product pairing between the cohomology class defined by $\upsilon_\Diamond$ and the first Chern class of the line bundle E. This being the case, the lemma's top bullet follows given a suitable bound for the absolute value of the integral of $\mathit{b}_\varepsilon\wedge *B_A$. To obtain such a bound, use (5.9) to write this form as $\mathcal{b} \wedge *(\mathfrak{r} - \mathfrak{X})$. As Lemma 4.4 finds $|\mathfrak{r}| \le c_0$, an appropriate bound for the integral of $r|\mathcal{b} \wedge *\mathfrak{X}|$ will suffice.

To obtain such a bound, use (5.7) to see that

$$|\mathcal{b}_\varepsilon \wedge *\mathfrak{X}| \le c_0\, r\, q_\diamond^{11/2}(|1 - \sigma|^{1/2} + c_c r^{-1/12})\,|1 - |\alpha|^2| + c_c\,. \tag{5.15}$$

Fix $z \ge 1$ for the moment. Having done so, use the inequality

$$q_\diamond^{11/2} c_c r^{-1/12}|1 - |\alpha|^2| \le z^{-1} q_\diamond^{6}(1 - |\alpha|^2)^2 + z\, c_c r^{-1} \tag{5.16}$$

to see that the term with the factor $q_\diamond^{11/2} r^{-1/12}$ in (5.15) contributes at most $z^{-1} r\, \underline{M}_\diamond + z\, c_c r$ to the integral of $r\chi_s|\mathcal{b}_\varepsilon \wedge *\mathfrak{X}|$.

The inequality

$$q_\diamond^{11/2}|1 - \sigma|^{1/2}|1 - |\alpha|^2| \le z^{-1} q_\diamond^{22/3}|1 - |\alpha|^2|^{4/3} + z c_0 (1 - \sigma)^2 (1 - |\alpha|^2)^2\,, \tag{5.17}$$

implies that the term in (5.15) with the factor $q_\diamond^{11/2}|1 - \sigma|^{1/2}$ contributes at most

$$c_0 z^{-1} r^2 \int_{\mathbb{R}\times Y} \chi_s q_\diamond^{22/3}|1 - |\alpha|^2|^{2/3} + c_0 z\, r^2 \int_{\mathbb{R}\times Y} \chi_s (1 - \sigma)^2 (1 - |\alpha|^2)^2\,. \tag{5.18}$$

to the integral of $r\chi_s|\mathcal{b}_\varepsilon \wedge *\mathfrak{X}|$. Use (5.7) to see that the right most integral in (5.18) is no greater than $c_0 z \int_{\mathbb{R}\times Y} |\frac{\partial}{\partial s} A|^2$ and therefore no greater than $z\, c_c r$ on account of Lemma 4.2.

<u>Step 3</u>: The conclusions of Step 2 supply the bound

$$r\,\Big|\int_{\mathbb{R}\times Y} \chi_s \mathcal{b}_\varepsilon \wedge *B_A\Big| \le c_0 z^{-1} r^2 \int_{\mathbb{R}\times Y} \chi_s q_\diamond^{22/3}|1 - |\alpha|^2|^{2/3} + z^{-1} r \underline{M}_\diamond + z c_c r\,. \tag{5.19}$$

This step and the next supply an approriate upper bound for

$$r^2 \int_{\mathbb{R}\times Y} \chi_s q_\diamond^{22/3}|1 - |\alpha|^2|^{4/3}\,. \tag{5.20}$$

To start this task, introduce $\kappa_*$ to denote the version of $\kappa$ given by Lemma 4.8. Separate the domain of integration in (5.20) to the set of points where $1 - |\alpha|^2 < \kappa_*^{-1}$ and the part where this inequality is violated. Noting that $q_\diamond^{22/3} = q_\diamond^{6}\, q_\diamond^{4/3}$, the contribution to the integral in (5.20) from the $1 - |\alpha|^2 \ge \kappa_*^{-1}$ part of the domain is no greater than $c_0\, r\, \underline{M}_\diamond$.

This understood, it is enough to bound the integral of (5.20) with the domain restricted to the subset in $[s-2, s+2] \times Y$ where $1 - |\alpha|^2 < \kappa_*^{-1}$.

Step 4: Let $X_* \subset [s-2, s+2] \times Y$ denote the set described by Lemma 4.8. Use $X_{**}$ to denote the set of points in $[s-2, s+2] \times Y$ with distance at most $2\kappa_* r^{-1/2} \ln r$ from $X_*$. It follows from Lemma 4.8 that the contribution to (5.20) from the complement of $X_{**}$ is bounded by $c_0 r$.

It follows from Lemma 4.12 and the second bullet of Proposition 4.11 that there is a point where $|\alpha| < \frac{1}{2}$ with distance $r^{-1/2}$ or less from each point in $X_{**}$. Let p denote such a point. The function $q_\diamond$ in the ball of radius $8r^{-1/2}\ln r$ centered at p is no less than $\frac{1}{2} q_\diamond(p)$ and no greater than $2q_\diamond(p)$ unless p has distance less than $c_0 r^{-1/2} \ln r$ from the zero locus of $q_\diamond$, this being $\cup_{p\in\Lambda} (\hat{\gamma}_p^+ \cup \hat{\gamma}_p^-)$. The contribution to (5.20) and to $r\underline{M}(s)$ from the set of such points is no greater than $c_0 r^{-10}$.

Step 5: Fix $n \in \{1, 2, \ldots\}$ with $n < 8\ln r$. It follows from what was just said in Step 4 and from Lemmas 4.10 and 5.2 that the contribution to $r\underline{M}(s-1) + r\underline{M}(s+1)$ from the ball of radius $n r^{-1/2}$ centered at p is no less than $c_c n^2 q_\diamond(p)^6$.

Fix a maximal set $\mathfrak{U}_n \subset X_{**}$ with the following properties:

- *The function* $|\alpha|$ *is less than* $\frac{1}{2}$ *at all points in* $\mathfrak{U}_n$.
- *The union of the balls of radius* $2^8 n r^{-1/2}$ *centered at the points in* $\mathfrak{U}_n$ *cover the subset in* $X_{**}$ *with distance* $4nr^{-1/2}$ *or less from the subset where* $|\alpha| < \frac{1}{2}$ .
- *The respective balls of radius* $nr^{-1/2}$ *centered at distinct points in* $\mathfrak{U}_n$ *are disjoint.*

(5.21)

The conditions in the second and third bullets of (5.21) imply that any given point in $X_{**}$ with distance $2^8 n r^{-1/2}$ or less from where $|\alpha| < \frac{1}{2}$ is in at most $c_0$ balls of radius $2^8 n r^{-1/2}$ centered at the points in $\mathfrak{U}_n$.

It follows from what said in this step's opening paragraph and from the condition in the third bullet of (5.21) that $r\underline{M}(s-1) + r\underline{M}(s+1) \geq c_c^{-1} n^2 \sum_{p\in\mathfrak{U}_n} q_\diamond(p)^6$.

Step 6: Use $X_n$ for $n \in \{1, 2, \ldots\}$ to denote the subset of $X_{**} \cap ([s-\frac{3}{2}, s+\frac{3}{2}] \times Y)$ with distance between $nr^{-1/2}$ and $(n-1)r^{-1/2}$ from $X_*$ and with $1 - |\alpha|^2 < \kappa_*^{-1}$. Lemma 4.8 with Proposition 4.11 and Lemma 4.12 have the following corollary: Given that $r \geq c_0$, there is a point in $X_{**}$ where $|\alpha| < \frac{1}{2}$ and with distance less than $(n+c_0)r^{-1/2}$ from each point in $X_n$.

Let p denote a point in $X_{**}$ where $|\alpha| < \frac{1}{2}$ . Use Lemma 4.8 to see that the contribution to (5.20) from the ball of radius $(n+c_0)r^{-1/2}$ centered at p is no greater than

$c_0\, e^{-n/c_0}\, n^4 q_\Diamond{}^6$. This understood, use the second bullet in (5.21) with the conclusions of the preceding step to see that the contribution from $X_n$ to the integral in (5.20) is no greater than $c_c\, e^{-n/c_0}\, r \sup_{s\in\mathbb{R}} \underline{M}_\Diamond$.

Sum these various $n \in \{1, 2, \ldots\}$ contributions to (5.20) to see that the contribution to (5.20) from $X_{**}$ is at most $c_c r \sup_{s\in\mathbb{R}} \underline{M}_\Diamond$.

Step 7: Use the conclusion of the final paragraph of Step 6 with the conclusions in the final paragraph of Step 3 to see that (5.20) is no greater than $c_c r \sup_{s\in\mathbb{R}} \underline{M}_\Diamond$. This being the case, then the right hand side of (5.19) is no greater than $c_c(z^{-1} r \sup_{s\in\mathbb{R}} \underline{M}_\Diamond + z r)$. The $z \ge c_c{}^2$ version of this last bound with what is said at the outset of Step 2 give the top bullet of Lemma 5.4

Step 8: This step proves the second bullet of Lemma 5.4. To this end, use Lemmas 4.8 to see that $|D_A\psi| \le c_0 r^{1/2}((1-|\alpha|^2) + c_0 r^{-1})^{1/2}$. Meanwhile, $|d(Q_\varepsilon q_\Diamond)| \le c_0 q_\Diamond{}^{11/2}$ because $q_\Diamond{}^{1/2} \ge c_0{}^{-1}|dq_\Diamond|$. These observations have the following consequence: The supremum in the second bullet is no greater than $c_0(z^{-1} r \sup_{s\in\mathbb{R}} \underline{M}_Q + z r)$ for any $z \ge 1$. Any $z \ge 1000\, c_0$ version of this last fact gives the assertion of Lemma 5.4's second bullet.

## 6. A priori bounds for $\underline{M}$

Fix $r > c_0$ and $\mu \in \Omega$ with $\mathcal{P}$-norm less than 1. Suppose that $\mathfrak{d}$ is an instanton solution to the corresponding $(r, \mu)$ version of (4.1). This section uses Proposition 5.1's bound on $\underline{M}_\Diamond$ to derive a $\mathfrak{d}$ and r independent bound for the function $\underline{M}$. The proposition that follows makes the formal statement that such a bound exists.

**Proposition 6.1**: *There exists $\kappa \ge \pi$ and given $c \ge 1$, there exists $\kappa_c > 1$ with the following significance: Fix $r \ge \kappa$ and $\mu \in \Omega$ has $\mathcal{P}$-norm less than 1. Suppose that $\mathfrak{d} = (A, \psi)$ is an instanton solution to the $(r, \mu)$ version of (4.1) with $A_\mathfrak{d} < c r$ and $\lim_{s\to\infty} M(\mathfrak{d}|_s) < c$. Then the corresponding function $\underline{M}$ obeys $-\kappa < \underline{M} < \kappa_c$.*

But for one extra lemma, the proof of Proposition 6.1 is in Section 6b. The extra lemma is proved in Section 6d. Section 6a makes observations that are used in the proof of Proposition 6.1. Sections 6a and 6b borrow much from the proof of Lemma 5.2 in [T6]. Section 6c supply some facts that are used in Section 6d and again in Section 7's proof of Theorem 1.5. The assertions in Section 6c all have analogs in Section 4 of [T6].

**a) Functions E and $\underline{\text{E}}$**

Let $\mathfrak{d}$ = (A, ψ): $\mathbb{R} \to \text{Conn}(E) \times C^\infty(M; \mathbb{S})$ denote an instanton solution to the $(r,\mu)$ version of (4.1) with $A_{\mathfrak{d}} < cr$ and with $\lim_{s\to\infty} M(\mathfrak{d}|_s) < c$. Introduce the function E(·) on $\mathbb{R}$ given by

$$s \to E(s) = i \int_{\{s\}\times Y} \hat{a} \wedge * B_A \,;$$

(6.1)

and introduce the function $\underline{\text{E}}$ on $\mathbb{R}$ given by the rule

$$s \to \underline{E}(s) = \int_{[s-1,s+1]\times Y} E(\cdot)\,.$$

(6.2)

This section explains how bounds on $\underline{\text{E}}$ give bounds on $\underline{\text{M}}$.

To do this, differentiate E and use the top bullet of (4.1) and an integration by parts to see that

$$\tfrac{d}{ds} E = i \int_{\{s\}\times Y} d\hat{a} \wedge (- * B_A + (\psi^\dagger \tau \psi - i\hat{a})) + \mathfrak{e}$$

(6.3)

with $|\mathfrak{e}| \le c_0$. The 2-form $d\hat{a}$ is zero on $M_\delta \cup \mathcal{H}_0$ and it is equal to $w$ where $|u| < R + c_0 \ln \delta$ on each $\mathfrak{p} \in \Lambda$ version of $\mathcal{H}_{\mathfrak{p}}$. In particular, $d\hat{a} = w$ where $q_\diamond \le c_0^{-1}$. This understood, use Lemma 4.4 with Proposition 5.1 to write (6.3) as

$$\tfrac{d}{ds} E = - E + M + \mathfrak{r}_E$$

(6.4)

where $|\mathfrak{r}_E| \le c_c$ with $c_c$ denoting a purely $c$ dependent constant. By way of notation, $c_c$ will hence forth denote a purely $c$-dependent constant that is greater than 1. Its value can be assumed to increase between successive appearances.

Integrate (6.4) to see that

$$E(s) = e^{-s} \int_{-\infty}^{s} e^x (M(x) + \mathfrak{r}_E)\, dx$$

(6.5)

It follows from Lemma 5.2 that $M(\cdot) \ge -c_0$; and thus (6.4) leads to the bound

$$- c_0 \le E(s) \le e^t (E(s+t) + c_c) \ \textit{for any}\ t \ge 0.$$

(6.6)

It then follows from (6.5) and (6.6) that

$$\mathrm{E}(s) \le c_c + c_c^2\, \underline{\mathrm{E}}(s+2) \;\; \textit{and} \;\; \underline{\mathrm{M}}(s) \le c_c(\underline{\mathrm{E}}(s+4) + 1). \tag{6.7}$$

Thus, a bound on $\underline{\mathrm{E}}$ gives a bound on $\underline{\mathrm{M}}$. By way of a converse to (6.7), note that

$$\mathrm{E}(s) \le (1+c_c r^{-1/12})\mathrm{M}(s) + c_c \;\; \textit{and} \;\; \underline{\mathrm{E}}(s) \le (1+c_c r^{-1/12})\underline{\mathrm{M}}(s) + c_c\,; \tag{6.8}$$

this being a consequence of (5.7) and (5.9).

As explained next, the function E is closely related to the function $s \to \mathrm{W}(A|_s)$ with W as defined in (1.27). The discussion that follows uses $\mathrm{W}_A(s)$ to denote $\mathrm{W}(A|_s)$. To say more about E and $\mathrm{W}_A$, use (4.1) to see that

$$\tfrac{d}{ds}\mathrm{W}_A = -\mathrm{E} + \mathrm{M} + \mathfrak{r}_w\,, \tag{6.9}$$

where $|\mathfrak{r}_w| \le c_0$. In particular, a comparision with (6.4) finds that

$$|\tfrac{d}{ds}(\mathrm{E} - \mathrm{W}_A)| \le c_c\,. \tag{6.10}$$

Most of $|\mathrm{E} - \mathrm{W}_A|$ is accounted for by the restriction of A to $M_\delta \cup \mathcal{H}_0$ in the sense that the $\hat{a} = \hat{a}_A|_s$ version of (2.7) with (5.7), (5.9) and Proposition 5.1 can be used to write

$$\mathrm{W}_A = \mathrm{E} + i\sum_{z\in ¥} C_{\mathbb{S},z} \int_{\{s\}\times\gamma^{(z)}} \hat{a}_A + \mathfrak{r}_*\,, \tag{6.11}$$

where $\{C_{\mathbb{S},z}\}_{z\in ¥}$ are integers and where $\mathfrak{r}_*$ is a function on $\mathbb{R}$ with $|\mathfrak{r}_*| \le c_c$. A given $z \in ¥$ version of the integral that appears in (6.11) is the value at $A|_s$ of (2.4)'s function $X_z$.

The lemma that follows says more about the ¥-indexed sum in (6.11). This lemma writes the $s \to \infty$ limit of $\mathfrak{d}$ as $(A_+, \psi_+)$ and write $A_+$ as $A_E + \hat{a}_{A_+}$.

**Lemma 6.2**: *There exists* $\kappa \ge \pi$ *and given* $c \ge 1$, *there exists* $\kappa_c > 1$ *with the following significance: Fix* $r \ge \kappa$ *and* $\mu \in \Omega$ *with* $\mathcal{P}$*-norm less than 1. Let* $\mathfrak{d} = (A, \psi)$ *denote an instanton solution to the* $(r, \mu)$ *version of (4.1) with* $A_\mathfrak{d} < cr$ *and* $\lim_{s\to\infty} \mathrm{M}(\mathfrak{d}|_s) < c$. *Then*

$$\Big|\sum_{z\in ¥} C_{\mathbb{S},z}\Big(\int_{\{s\}\times\gamma^{(z)}} \hat{a}_A\Big) - \sum_{z\in ¥} C_{\mathbb{S},z}\Big(\int_{\gamma^{(z)}} \hat{a}_{A_+}\Big)\Big| < \kappa_c$$

*for all* $s \in \mathbb{R}$.

This lemma is proved in Section 6d. Accept it as true in the meantime. By way of a look a head, Proposition 5.1 plays a key role in the proof of this lemma; it plays no role otherwise.

**b) An algebraic inequality for $\underline{\mathrm{E}}$.**

The equation that follows asserts the desired algebraic inequality for $\underline{\mathrm{E}}$.

$$\underline{\mathrm{E}}(s) \leq c_c(1 + r^{-1/3}\sup_{x\geq s+1} |\underline{\mathrm{E}}(x)|^{4/3}) . \tag{6.12}$$

The derivation of this formula is given momentarily. What follows directly assumes (6.12) to prove Proposition 6.1.

***Proof of Proposition 6.1***: The $s \to \infty$ limit of $\underline{\mathrm{E}}$ is by assumption bounded by $c$. Fix for the moment $T > c$ and let $s_T$ denote the largest value of $s$ with $\underline{\mathrm{E}}(s) \geq T$. Let $c_*$ denote the version of $c_c$ that appears in (6.12). The $s = s_T$ version of (6.12) reads $T < c_*(1 + r^{-1/3}T^{4/3})$. This has no solutions for $T \in (2c_*, c_*^{-1}r)$ if $r > 2^8 c_*^8$. The bound $\underline{\mathrm{E}} < 2c_*$ leads via (6.7) to a purely $c$-dependent bound for $\underline{\mathrm{M}}$.

The four parts that follow derive the inequality in (6.12).

*Part 1*: The derivation starts with an inequality that involves the functions $\mathfrak{a}(\mathfrak{d}|_s)$, $\mathrm{E}$, $\mathrm{M}$ and a function, $\mathrm{O}$, on $\mathbb{R}$ that defined at any given $s$ by

$$\mathrm{O}(s) = \int_{\{s\}\times Y} (|\mathfrak{B}_{(A,\psi)}|^2 + r\,|D_A\psi|^2) \tag{6.13}$$

with $\mathfrak{B}_{(A,\psi)}$ defined in (4.2). The derivation of this first inequality mimics the derivation of an analogous inequality in Section 5b of [T6].

To start, use (1.28) to see that

$$\mathfrak{a}(\mathfrak{d}|_s) \leq \tfrac{1}{2}\,\mathfrak{cs}(A|_s) - r\,\mathrm{w}_A(s) + c_0 r^{1/2}\mathrm{O}(s)^{1/2} . \tag{6.14}$$

The next part of the derivation talks about the function $\mathfrak{cs}$.

*Part 2*: The formula for $\mathfrak{cs}$ is given in (1.26) as a sum of two integrals. To say more about the right most integral in (1.26), keep in mind that the $i\mathbb{R}$ valued 2-form $2F_{A_E} + F_{A_K}$ that appears this formula is cohomologous to $-2\pi i\, w$. This being the case, their difference is the exterior derivative of a fixed, smooth 1-form, this denoted by $\mathfrak{y}$. As a consequence, integration by parts equates the right most integral in (1.26) with

$$-2 \int_{\{s\}\times Y} \hat{a}_A \wedge (F_{A_E} + \tfrac{1}{2} F_{A_K}) = 2\pi \mathrm{w}_A + i \int_{\{s\}\times Y} *B_A \wedge \mathfrak{y} \ . \tag{6.15}$$

Use Lemma 5.2's preliminary bound for $\underline{\mathrm{M}}$ in Lemma 4.7 to bound the absolute value of the right most integral in (6.15) by $c_{c}(\mathrm{M}+1)$.

The remaining term in (1.26) is the integral of $\hat{a}_A \wedge d\hat{a}_A$. This term can be bounded by writing $\hat{a}_A = \hat{a}_A{}^{\perp} + \mathfrak{q}$ with $\hat{a}_A{}^{\perp}$ a coclosed 1-form that is orthogonal to the space of harmonic 1-forms on Y. Meanwhile, $\mathfrak{q}$ is a closed 1-form on Y. The integral of $\hat{a}_A \wedge d\hat{a}_A$ is the same as that of $\hat{a}_A{}^{\perp} \wedge d\hat{a}_A{}^{\perp}$. Meanwhile, the Green's function for the operator $*d$ can be used much as in the proof of Lemma 2.4 of [T1] to bound $|\hat{a}_A{}^{\perp}|$ by $c_0 r^{2/3} \mathrm{M}^{1/3}$. Given this last bound, use the fact that $d\hat{a}^{\perp} = B_A$ and use Lemma 5.2's preliminary bound for $\underline{\mathrm{M}}$ in Lemma 4.7 to see that

$$\Big| \int_{\{s\}\times Y} \hat{a}_A \wedge d\hat{a}_A \Big| \le c_{c}\, r^{2/3} \mathrm{M}^{4/3} . \tag{6.16}$$

Use the bound in (6.16) and the bound by $c_{c}(\mathrm{M} + 1)$ for the right most integral in (6.15) to see that

$$\mathfrak{cs}(A|_s) \le 2\pi \mathrm{w}_A + c_{c}(1 + \mathrm{M} + r^{2/3} \mathrm{M}^{4/3}) . \tag{6.17}$$

The next step exploits this inequality for $\mathfrak{cs}$.

*Part 3*: Replacing $\mathfrak{cs}$ in (6.14) with the right hand side of (6.17) leads to the inequality

$$\mathfrak{a}(\mathfrak{d}|_s) \le -(r - \pi)\mathrm{w}_A + c_{c}(\mathrm{M} + \mathrm{O} + r + r^{2/3} \mathrm{M}^{4/3}) . \tag{6.18}$$

Replace the function $\mathrm{w}_A$ in (6.18) by $\mathrm{E}$ using Lemma 6.2. Having done so, rearrange terms to obtain the following inequality for $\mathrm{E}$:

$$(r - \pi)\mathrm{E} \le -\mathfrak{a}(\mathfrak{d}|_s) - (r - \pi)\textstyle\sum_{z\in ¥} C_{\mathbb{S},z}\Big( \int_{\gamma^{(z)}} \hat{a}_{A_+} \Big) + c_{c}(\mathrm{O} + r + r^{2/3} \mathrm{M}^{4/3}) . \tag{6.19}$$

As the function $s \to \mathfrak{a}(\mathfrak{d}|_s)$ is non-increasing, the right hand side of (6.19) is no less than

$$-\mathfrak{a}(\mathfrak{c}_+) - (r-\pi)\sum_{z\in\Upsilon} C_{\mathbb{S},z}\Big(\int_{\gamma^{(z)}} \hat{a}_{A_+}\Big) + c_c(O + r + r^{2/3}M^{4/3}) .$$

(6.20)

Use the $A_+$ versions of (6.11), (6.15) and (6.16) to bound the combined two left most terms in (6.20) by $c_0(r\mathrm{M}(\mathfrak{c}_+) + r^{2/3}\mathrm{M}(\mathfrak{c}_+)^{4/3})$. Using this bound leads to the inequality

$$rE \le c_c(r\mathrm{M}(\mathfrak{c}_+) + r^{2/3}\mathrm{M}(\mathfrak{c}_+)^{4/3} + O + r + r^{2/3}M^{4/3})$$

(6.21)

when $r > 2\pi$. The assumed $\mathrm{M}(\mathfrak{c}_+) \le c$ bound and (6.21) imply that $rE \le c_c(O + r + r^{2/3}M^{4/3})$.

*Part 4*: Let F for the moment denote any given non-negative function on [-1, 1] and let $\underline{F}$ denote its integral over this interval. The measure of the set of points where F is less than $8\underline{F}$ must be greater than $\frac{15}{8}$. This being the case, suppose that $F'$ is a second non-negative function. Then there are points in $[s-1, s+1]$ where both F and $F'$ are less than $8\underline{F}$ and $8\underline{F}'$, respectively.

With the preceding in mind, introduce $\underline{O}(s)$ to denote the integral of $O(\cdot)$ over the interval $[s-1, s+1]$. Fix $s' \in [s-1, s+1]$ where $O(s') \le 8\underline{O}(s)$ and $M(s') \le 8\underline{M}(s)$. The $s'$ version of the inequality $rE \le c_c(O + r + r^{2/3}M^{4/3})$ implies that

$$rE(s') \le c_c(\underline{O}(s) + r + r^{2/3}\,\underline{M}(s)^{4/3}) .$$

(6.22)

As explained next, the inequality in (6.12) follows from (6.20) with three additional replacements. The first replacement invokes Lemma 4.2 to substitute $2A_{\mathfrak{d}}$ for $\underline{O}(s)$. The second replacement invokes (6.7) to replace $\underline{M}(s)$ with $\sup_{x\ge s}\underline{E}(x)$.

To explain the final repacement, fix for the moment $s'' \in [s-3, s-1]$ and invoke (6.6) with $t = s' - s''$. With the first and second replacements made, (6.22) and (6.6) imply

$$rE(s'') \le c_c(r + A_{\mathfrak{d}} + r^{2/3}\sup_{x\ge s}\underline{E}(x)^{4/3}).$$

(6.23)

View both sides of (6.23) as functions on $[s-3, s-1]$, with the right hand side being the constant function. Integrate both sides over this interval. The integral of the left hand side is $r\underline{E}(s-2)$ and that of the right is twice what is written in (6.23). The resulting inequality with the assumed $A_{\mathfrak{d}} \le c$ bound leads directly to (6.12) when evaluated at $s+2$ rather than $s$.

**c) Local convergence to pseudoholomorphic subvarieties**

The upcoming Proposition 6.3 describes how certain instanton solutions to (4.1) can be used to determine pseudoholomorphic subvarieties in bounded subsets of $\mathbb{R} \times Y$. Proposition 6.3 and the subsequent two lemmas about pseudoholomorphic subvarieties are used to prove Lemma 6.2 and they are invoked again in Section 7 to prove Theorem 1.5. Proposition 6.3 is the analog of Proposition 4.1 in [T6] and subsequent two lemmas are the respective analogs of Lemma 4.6 and Corollary 4.7 in [T6].

Proposition 6.3 and the two lemmas use $Y_*$ to denote either $M_\delta \cup \mathcal{H}_0$ or $Y$. Their assertions with regards to $Y_* = M_\delta \cup \mathcal{H}_0$ are used in the upcoming proof of Lemma 6.2 and those that concern $Y_* = Y$ are used in Section 7.

Proposition 6.3 introduces the notion of a pseudoholomorphic subvariety in an open set of $\mathbb{R} \times Y$. To define this term, let $U \subset \mathbb{R} \times Y$ denote the open set. A subset $C \subset U$ is said to be a pseudoholomorphic subvariety in U when the conditions below are met.

- C *is the intersection between* U *and a closed subset,* C´, *of a neighborhood of* U.
- *The complement in* C´ *of a finite set of points is a smoothly embedded submanifold of this neighborhood with* J*-invariant tangent space.*
- C´ *has no totally disconnected components.*
- *The integral of w over* C´ *is finite.*
- *There exists an* $s \in \mathbb{R}$ *independent upper bound for the integral of* $ds \wedge \hat{a}$ *over the intersection between* C´ *and* $[s-1, s+1] \times Y$.

(6.24)

The subvariety C is said to be irreducible when the complement in C of any finite set of points is connected.

**Proposition 6.3**: *Fix* $c \geq 1$ *and in the case* $Y_* = Y$, *also* $\mathcal{K} > 1$. *There exists* $\kappa_c > 1$, *and given* $m > 100$, *there exists* $\kappa_{c,m} > \pi$*; these having the following significance: Fix* $r \geq \kappa_{c,m}$ *and fix* $\mu \in \Omega$ *with* $\mathcal{P}$*-norm less than 1. Suppose that* $\mathfrak{d} = (A, \psi = (\alpha, \beta))$ *is an instanton solution to the* $(r, \mu)$ *version of (4.1) with* $A_\mathfrak{d} < cr$. *If* $Y_* = Y$ *assume in addition that* $\sup_{s\in\mathbb{R}} \underline{M}(s) < \mathcal{K}$. *Let* $I \subset \mathbb{R}$ *denote an interval of length* $2m$. *There exists a finite set,* $\vartheta$, *of at most* $\kappa_c$ *elements with each element having the form* $(C, m)$ *with* C *an irreducible, pseudoholomorphic subvariety in* $I \times Y_*$ *and* m *a positive integer no greater than* $\kappa_c$. *Distinct pairs from* $\vartheta$ *have distinct subvariety components. Furthermore,*

- $\sup_{z\in\cup_{(C,m)\in\vartheta}(C\cap(I\times Y_*))} \mathrm{dist}(z, \alpha^{-1}(0)) + \sup_{z\in(\alpha^{-1}(0)\cap(I\times Y_*))} \mathrm{dist}(\cup_{(C,m)\in\vartheta} C, z) < m^{-1}$
- *Let* $\upsilon$ *denote the restriction to* $I \times Y_*$ *of a smooth 2-form on* $\mathbb{R} \times Y$ *with* $\|\upsilon\|_\infty = 1$ *and* $\|\nabla\upsilon\|_\infty \leq m^{-1}$. *Then* $\left| \frac{i}{2\pi}\int_{I\times Y_*} \upsilon \wedge F_{\hat{A}} - \sum_{(C,m)\in\vartheta} m\int_{C\cap I\times Y_*} \upsilon \right| \leq m^{-1}$.

***Proof of Proposition 6.3***: But for one extra remark in the $Y_* = M_\delta \cup \mathcal{H}_0$ case, the arguments for the first bullet and Items a) and b) of the second bullet of Proposition 4.1 in [T6] can be used with only notational changes to prove Proposition 6.3. The extra remark concerns the assumption made by Proposition 4.1 for a bound on $\underline{M}(s)$. The arguments for Proposition 4.1 in [T6] require only a bound on Lemma 4.10's function $M(x,\rho)$ for $\rho = c_0^{-1}$ and for all $x \in \mathbb{R} \times Y_*$. Such a bound comes from Proposition 5.1's bound on $\underline{M}_\Diamond$ when $Y = M_\delta \cup \mathcal{H}_0$.

The next two lemmas are used to say more about the subvarieties that can arise in the context of Proposition 6.3. The lemma that follows directly is an analog of Lemma 4.6 in [T6]. The notation is such that if $s \in \mathbb{R}$ and C is a pseudoholomorphic subvariety that is defined in a neighborhood of $\{s\} \times Y_*$, then $C|_s$ denotes $C \cap (\{s\} \times Y_*)$.

**Lemma 6.4**: *Given $m > 1$ and $\varepsilon > 0$, there exists $\kappa_{m\varepsilon} > 1$ with the following significance: Suppose that C is a closed, irreducible, pseudoholomorphic subvariety in a neighborhood of $[-4, 4] \times Y_*$ with $\int_{C\cap([-4,4]\times Y_*)} w < \kappa_{m\varepsilon}^{-1}$ and $\int_{C\cap([-4,4]\times Y_*)} ds \wedge \hat{a} \le m$. Then each point of $C|_s$ for $|s| \le 1$ has distance along $Y_*$ no greater than $\varepsilon$ from a closed integral curve, $\gamma$, of length less than $m + \varepsilon$. Moreover, there is a positive integer q which is bounded by an m and $\varepsilon$ independent multiple of m, and is such that if $\upsilon$ is a smooth 2-form on $[-1,1] \times Y_*$ with $\|\upsilon\|_\infty = 1$ and $\|\nabla\upsilon\|_\infty \le \varepsilon^{-1}$, then $\left| \int_{C\cap([-1,1]\times Y_*)} \upsilon - q \int_{[-1,1]\times\gamma} \upsilon \right| \le \varepsilon$.*

***Proof of Lemma 6.4***: The proof is the same but for notation of Lemma 4.6 in [T6].

The next lemma is an analog of Corollary 4.7 in [T6]. Note in this regard that Corollary 4.7 in [T6] makes an assumption that is not guaranteed here, this being that all integral curves of $v$ with an a priori length bound are nondegenerate. The upcoming Lemma 6.5 is a version of Corollary 4.7 in [T6] that suffices for the present purposes.

**Lemma 6.5**: *Given $m > 1$ and $\varepsilon > 0$, there exists $\kappa_{m\varepsilon} > 1$ with the following significance: Let $\mathbb{I} \subset \mathbb{R}$ denote an interval of length at least 4, and suppose that C is an irreducible, pseudoholomorphic subvariety in a neighorhood of $\mathbb{I} \times Y_*$ with $\int_{C\cap(I'\times Y_*)} w < \kappa_{m\varepsilon}^{-1}$ and $\int_{C\cap(I'\times Y_*)} ds \wedge \hat{a} < m$ for all intervals $I' \subset \mathbb{I}$ of length 1. Assume in addition that C has intersection number zero with all submanifolds in $\mathbb{R} \times Y$ of the form $\{s\} \times S$ with S being a cross-sectional sphere in $\mathcal{H}_0$. Let $I \subset \mathbb{I}$ denote the subset with distance at least 3 from any boundary point of $\mathbb{I}$. There exists a finite set $\Theta$ consisting of pairs $(\gamma, q)$ with $\gamma$ a* *closed, integral curve of* $\nu$ *and* q *a positive integer. The set* $\Theta$ *is such that no two pair share the same closed integral curve. Moreover,*

- *The intersection of* $\gamma$ *with* $M_\delta$ *is a collection of arcs that begin on the boundary of radius* $\delta$ *coordinate balls about the index* 1 *critical points of f in* M *and end on the boundary of radius* $\delta$ *coordinate balls about the index* 2 *critical points of f in* M.
- $\sum_{(\gamma,q)\in\Theta} q\, \ell_\gamma < m + \varepsilon$.
- *Each point of* $C|_s$ for $s \in I$ *has distance along* Y *less than* $\varepsilon$ *from* $\cup_{(\gamma,q)\in\Theta}\gamma$. *Conversely, each point in* $\cup_{(\gamma,q)\in\Theta}\gamma$ *has distance no greater than* $\varepsilon$ *from* $C|_s$.
- *If* $\upsilon$ *is a smooth 2-form on* $I \times Y$ *with* $\|\upsilon\|_\infty = 1$ *and* $\|\nabla\upsilon\|_\infty \le \varepsilon^{-1}$. *Then*

$$| \int_{C\cap(I\times Y)} \upsilon \ - \ \sum_{(\gamma,q)\in\Theta} q \int_{I\times\gamma} \upsilon \ | < \varepsilon.$$

***Proof of Lemma 6.5***: But for one additional remark, the argument in [T6] that explains how Corollary 4.7 in [T6] follows from Lemma 4.6 in [T6] explains why Lemma 6.5 follows from Lemma 6.4. The additional remark concerns both the first bullet of the lemma and the assumption for Corollary 4.7 in [T6] of nondegenerate Reeb orbits. The assumption that C has intersection number zero with submanifolds of the form $\{s\} \times S$ with $S \subset \mathcal{H}_0$ being a cross sectional sphere replaces the nondegeneracy assumption in Lemma 4.6 of [T6] and it leads directly to the first bullet of Lemma 6.5. To see how this comes about, note that if $S \subset \mathcal{H}_0$ is a constant u sphere, then $\{s\}\times S$ is pseudoholomorphic, so if $C|_s$ is close to a closed integral curve of $\nu$ that crosses $\mathcal{H}_0$, then C will have positive intersection number with $\{s\}\times S$. This is ruled out by assumption. Meanwhile, Section II.2 finds that the only closed integral curves of $\nu$ that don't intersect $\mathcal{H}_0$ are hyperbolic and so non-degenerate. Moreover, those that intersect $M_\delta$ are described by the first bullet of Lemma 6.5.

## d) Proof of Lemma 6.2

The proof has four parts. These parts use $c_c$ to denote a number greater than 1 that depends only on $c$. Its value can be assumed to increase between successive appearances.

*Part 1*: The curvature of the version (1.15) that defines $\hat{A}$ can be written as

$$F_{\hat{A}} = (1 - \wp)(ds \wedge \tfrac{\partial}{\partial s} A + {*B_A}) + \wp' \nabla_A\alpha \wedge \nabla_A\bar{\alpha} \ . \tag{6.25}$$

with it understood that the d$s$ component of $\nabla_A\alpha$ is $\frac{\partial}{\partial s}\alpha$. The notation here uses $\wp$ and $\wp'$ to denote the respective functions on $I_* \times Y_*$ given by $\wp|_{t=|\alpha|^2}$ and $(\frac{d}{dt}\wp)|_{t=|\alpha|^2}$. Use (6.25) to see that $w \wedge F_{\hat{A}}$ can be written as $-i\mathrm{F}ds \wedge \hat{a} \wedge w$ with $\mathrm{F}$ being

$$i(1 - \wp)(\tfrac{\partial}{\partial s}A)_v - i\wp'((\tfrac{\partial}{\partial s}\alpha)(\nabla_A\bar{\alpha})_v - (\nabla_A\alpha)_v(\tfrac{\partial}{\partial s}\bar{\alpha})) , \tag{6.26}$$

where $(\frac{\partial}{\partial s}A)_v$ and $(\nabla_A\alpha)_v$ denote the pairing of these 1-forms with the vector field $v$.

*Part 2*: Let I denote an interval of length 2 and introduce V to denote the subset of $I \times (M \times \mathcal{H}_0)$ where $|\alpha|^2 < \frac{5}{8}$. The support of $F_{\hat{A}}$ in $I \times (M \times \mathcal{H}_0)$ is in V. Use (6.25) and (6.26) to see that

$$c_0 \int_V |\tfrac{\partial}{\partial s}\alpha|^2 \geq |\tfrac{i}{2\pi}\int_{I\times(M_\delta \cup \mathcal{H}_0)} w \wedge F_{\hat{A}}| - c_0 \operatorname{vol}(V) \tag{6.27}$$

where vol(V) denotes the volume of the set V. Proposition 5.1 bounds the integral of $r(1 - |\alpha|^2)$ over V by $c_c$ and this implies that $\operatorname{vol}(V) \leq r^{-1} c_c$. Therefore, (6.27) implies that

$$c_0 r \int_{I\times Y} |\tfrac{\partial}{\partial s}\psi|^2 \geq r\,|\tfrac{i}{2\pi}\int_{I\times(M_\delta \cup \mathcal{H}_0)} w \wedge F_{\hat{A}}| - c_c . \tag{6.28}$$

This last bound leads directly to the following conclusion: If $\varepsilon \in (r^{-1}c_c, 1)$, there are most $\varepsilon^{-1}c_c$ disjoint intervals in I of the form $[s - 1, s + 1]$ with $|\frac{i}{2\pi}\int_{[s-1,s+1]\times(M_\delta \cup \mathcal{H}_0)} w \wedge F_{\hat{A}}| > \varepsilon$.

*Part 3*: Apply the $Y_* = M_\delta \cup \mathcal{H}_0$ version of Proposition 6.3 to intervals of length 200 in $\mathbb{R}$. Use the first bullet of the latter, Lemma 6.5 and the final conclusion in Part 2 to see that there exists a set of at most $c_c$ points in $\mathbb{R}$ with the following property: If $s \in \mathbb{R}$ has distance 1 or more from all points in this set, then $F_{\hat{A}} = 0$ and $\alpha/|\alpha|$ is $\hat{A}$-covariantly constant at points with distance $c_0^{-1}$ or less from any z ∈ ¥ version of $\{s\} \times \gamma^{(z)}$. Let $\mathcal{Q}$ denote this finite set in $\mathbb{R}$.

Suppose that $s \in \mathbb{R}$ has distance less than 2 from some point in $\mathcal{Q}$. The fact that $\mathcal{Q}$ has at most $c_c$ elements implies that there are points in $[s - c_c, s + c_c]$ with distance at least 2 from each point in $\mathcal{Q}$. Let $s'$ denote such a point. Use the $s$ and $s'$ versions of (6.11) with the derivative bound in (6.10) to conclude that Lemma 6.2 is true for $s$ if and only if it is true for $s'$.

Introduce $\mathcal{Q}_* \subset \mathbb{R}$ to denote the set of points with distance less than 2 from some point in $\mathcal{Q}$. Let $(s, s') \subset \mathbb{R}$ denote a connected component of $\mathcal{Q}_*$. Then $|s' - s| < c_c$ because $\mathcal{Q}$ has at most $c_c$ elements. This understood, use the $s$ and $s'$ versions of (7.11) with (7.10) again to see that

$$|\Sigma_{z\in ¥} c_{\mathbb{S},z}\Big(\int_{\{s\}\times\gamma^{(z)}} \hat{a}_A\Big) - \Sigma_{z\in ¥} c_{\mathbb{S},z}\Big(\int_{\{s'\}\times\gamma^{(z)}} \hat{a}_{A_+}\Big)| < \kappa_c \tag{6.29}$$

Given the conclusions of the preceding two paragraphs,the fact that $\mathcal{Q}$ has at most $c_c$ elements implies that Lemma 6.2 holds if it (6.29) is also true when $s$ and $s'$ are any two elements in the same component of $\mathbb{R}-\mathcal{Q}_*$.

*Part 4*: To see about (6.29) when $s$ and $s'$ are in the same component of $\mathbb{R}-\mathcal{Q}_*$, fix for the moment a point $z \in ¥$. Write $\hat{A}$ as $\hat{A} = A_E + \hat{a}_{\hat{A}}$ and use the $\mathbb{R} \times Y$ version of (1.15) with Lemma 4.8 to see that

$$|\int_{\{s\}\times\gamma^{(z)}} \hat{a}_A - \int_{\{s\}\times\gamma^{(z)}} \hat{a}_{\hat{A}}| \le c_0 . \tag{6.30}$$

when $s$ has distance 2 or more from every point in $\mathcal{Q}$. Note that this inequality also holds with A replaced by $A_+$ and with $\hat{A}$ replaced by $\hat{A}_+$, this being the $s \to \infty$ limit of $\hat{A}$.

With (6.30) in mind, suppose that $s' > s$ are in the same component of $\mathbb{R}-\mathcal{Q}_*$. Use Stoke's theorem to see that

$$\int_{\{s'\}\times\gamma^{(z)}} \hat{a}_{\hat{A}} - \int_{\{s\}\times\gamma^{(z)}} \hat{a}_{\hat{A}} = \int_{[s,s']\times\gamma^{(z)}} F_{\hat{A}} . \tag{6.31}$$

The right hand side of (6.31) is zero, so it follows using the $s$ and $s'$ versions of (6.30) that the integral of $\hat{a}_A$ over $\{s\}\times\gamma^{(z)}$ differs by at most $c_0$ from its integral over $\{s'\} \times \gamma_z$. Thus, (6.29) does indeed hold for any pair $s' > s$ in the same component of $\mathbb{R}-\mathcal{Q}_*$.

## 7. Proposition 1.1-1.4 and Theorem 1.5

This last section supplies the proofs for Section 1's propositions and theorem.

### a) Proofs of Propositions 1.1-1.3

Leave out for the moment the second and third bullets of Proposition 1.1 and the assertions of Propositions 1.2 and 1.3 that refer to $\hat{\mathcal{Z}}^{\ge}_{SW,r}$. The remaining assertions of these propositions, those that refer only to $\hat{\mathcal{Z}}_{\mathrm{SW,r}}$, are all special cases of theorems from [KM]. To elaborate, the essential concern is a compactness theorem for the space of instanton solutions (1.20). See in particular the discussion at the beginning of Chapter 29.2 in [KM]. The desired compactness theorem is Proposition 29.2.1 in [KM]. This is because the $r > \pi$ versions of (1.14) and (1.20) are defined by what is said in [KM] to be a *monotone* perturbation.

The second bullet of Proposition 1.1 follows directly from Proposition 2.4. The third bullet of Proposition 1.1 and the assertions about $\hat{\mathcal{Z}}^{\geq}_{\mathrm{SW,r}}$ in Propositions 1.2 and 1.3 follow from a proof that the value of the function $\mathrm{x}$ in (1.16) on the $s \to \infty$ limit of any relevant instanton is no less than its value on the $s \to -\infty$ limit if the instanton contributes to the differential on the chain complex, or to one of the other homomorphisms. This property of $\mathrm{x}$ follows from the upcoming Proposition 7.1 together with Lemmas 2.5, 4.1 and 4.2. Note in this regard that Proposition 7.1 proves this assertion about $\mathrm{x}$ for instanton solutions to (4.1), this being the version of (1.20) that uses $\mathfrak{g} = \mathfrak{e}_\mu$ with $\mu \in \Omega$ having $\mathcal{P}$-norm less than 1. Even so, the fact that $\lim_{s\to\infty} \mathrm{x}(\mathfrak{d}|_s) \geq \lim_{s\to-\infty} \mathrm{x}(\mathfrak{d}|_s)$ for the instanton solutions to (4.1) implies this inequality is also true for any instanton solutions to a $\mathfrak{g} = \mathfrak{e}_\mu + \mathfrak{p}$ version of (1.20) that contributes to the differential or the other relevant homomorphisms if $\mathfrak{p}$ comes from a certain residual set in $\mathcal{P}_\mu$ and has small $\mathcal{P}$-norm. More is said about why this is after the statement of Proposition 7.1.

To set the stage for Proposition 7.1, let $\gamma$ denote a closed, integral curve of $v$. Define the function $\mathrm{x}_\gamma$ on Conn(E) by the rule that assigns to a connection A on E the integral over the curve $\gamma$ of the $i\mathbb{R}$ valued 1-form $\frac{i}{2\pi}(\hat{A} - A_E)$. The $\gamma = \gamma^{(z_0)}$ version of $\mathrm{x}_\gamma$ is the function $\mathrm{x}$ in (1.16).

The proposition assumes that $\gamma \subset M_\delta \cup \mathcal{H}_0$ and that $\gamma$ has a tubular neighborhood of the sort described directly. Let $\ell$ denote the length of $\gamma$ and let $t \in \mathbb{R}/(\ell\mathbb{Z})$ denote an affine parameter for $\gamma$. Use z to denote the complex coordinate for $\mathbb{C}$. The operative assumption is that $\gamma$ has a tubular neighborhood with coordinates (t, z) that are defined for $|z|$ less than a positive constant and are such that

- *The curve $\gamma$ is the $z = 0$ locus.*
- *The vector field $v$, the 2-form $w$ and the 1-form $\hat{a}$ appear as*

$$v = \frac{\partial}{\partial t} + \cdots \ , \quad w = \frac{i}{2} dz \wedge d\bar{z} + \cdots \quad \textit{and} \quad \hat{a} = dt \ ,$$

  *where the unwritten terms are bounded by a constant multiple of $|z|$.*
- *The vector field $v$ annihilates $|z|^2$*

(7.1)

It follows from the constructions in Section II.1c and Section II.1d that each $z \in ¥$ version of $\gamma^{(z)}$ has a tubular neighborhood with coordinates of the sort described by (7.1); and in particular, the curve $\gamma^{(z_0)}$ has such a tubular neighborhood.

**Proposition 7.1**: *There exists* $\kappa \geq \pi$ *and given* $c \geq 1$, *there exists* $\kappa_c > 1$ *with the following significance: Fix* $r \geq \kappa$ *and* $\mu \in \Omega$ *with* $\mathcal{P}$*-norm less than 1 and suppose that* $\mathfrak{d} = (A, \psi)$ *is an instanton solution to the* $(r, \mu)$ *version of (4.1) with* $A_{\mathfrak{d}} < c\, r \ln r$. *Let* $\gamma$ *denote a closed, integral curve of* $v$ *that lives entirely in* $M_\delta \cup \mathcal{H}_0$ *and has a tubular neighborhood with coordinates of the sort described by (8.1). Then* $\lim_{s\to\infty} x_\gamma(A|_s) \geq \lim_{s\to-\infty} x_\gamma(A|_s)$.

Given what is said by the fourth bullet of Proposition 2.7 and Lemma 4.1, the assumption $A_{\mathfrak{d}} < c\, r \ln r$ is satisfied if the difference between the value of $\mathfrak{f}_\mathfrak{s}$ on the $s \to \infty$ limit of $\mathfrak{d}$ and the value of $\mathfrak{f}_\mathfrak{s}$ on the $s \to -\infty$ limit of $\mathfrak{d}$ is no greater than $c \ln r$. The bound on $A_{\mathfrak{d}} < c\, r \ln r$ in Proposition 7.1 is used to invoke Lemma 5.2 so as to bound $\mathfrak{d}$'s version of the function $\underline{M}$ by $c_c r^{6/7}$. This bound on $\underline{M}$ is then used to invoke Lemma 4.7. Lemma 4.7 in turn is used to write $\frac{\partial}{\partial s} A$ and $B_A$ as in (5.7). A crucial point in this regard is that the function $\sigma$ that appears in (5.7) is constrained to obey $\sigma > -c_c r^{-1/q}$ and $1 - \sigma > -c_c r^{-1/q}$ with $q \in (1, c_0)$ a fixed constant and with $c_c$ a constant greater than 1. These bounds on $\sigma$ follow from the fact that the left hand side of (5.9) is non-negative.

The proof of Proposition 7.1 is given in Section 7b. What follows directly explains how Proposition 7.1 is used to prove the assertions in Propositions 1.1-1.3 that refer to $\hat{\mathcal{Z}}^{\geq}_{SW,r}$. To this end, suppose that the conclusions of Proposition 7.1 hold for instanton solutions to a given $\mathfrak{g} = \mathfrak{e}_\mu + \mathfrak{p}$ version of (1.20) if $(r,\mu)$ obey its assumptions and if the perturbation $\mathfrak{p}$ is in $\mathcal{P}_\mu$. If $\mathfrak{d}$ is an instanton solution to this version of (1.20) and if it contributes to either the differential or one of the other relevant homomorphisms of the Seiberg-Witten chain complex, then the $s \to \infty$ limit of $\mathfrak{f}_\mathfrak{s}(\mathfrak{d}|_s)$ is either 1 or 2 more than the $s \to -\infty$ limit of $\mathfrak{f}_\mathfrak{s}(\mathfrak{d}|_s)$. Use this observation with Proposition 2.7 and Lemma 4.1 to see that $\mathfrak{d}$ obeys Proposition 7.1's bound $A_{\mathfrak{d}} \leq c\, r \ln r$.

Given what was said in the preceding paragraph, the assertions about $\hat{\mathcal{Z}}^{\geq}_{SW,r}$ in Propositions 1.1-1.3 hold if the conclusions of Proposition 7.1 hold for instanton solutions to any $\mathfrak{g} = \mathfrak{e}_\mu + \mathfrak{p}$ version of (1.20) if $(r,\mu)$ obey its assumptions and if $\mathfrak{p} \in \mathcal{P}_\mu$ has small $\mathcal{P}$-norm. To see why this is so, assume it to be false so as to derive nonsense. Under this contrary assumption, there is a sequence $\{\mathfrak{p}_n\}_{n=1,2,\ldots}$ with the following two properties: The $\mathcal{P}$-norm of each $n \in \{1, 2, \ldots\}$ version of $\mathfrak{p}_n$ is less than $n^{-1}$ and the conclusions of Proposition 7.1 fail for some instanton solution to the $\mathfrak{g} = \mathfrak{e}_\mu + \mathfrak{p}_n$ version of (1.20) with $A_{\mathfrak{d}} < c\, r \ln r$. Let $\{\mathfrak{d}_{\mathfrak{p}_n}\}_{n=1,2,\ldots}$ denote a corresponding sequence of recalcitrant instantons. This sequence can be chosen so that all its constituent members have the same $s \to \infty$ limit, and all have the same $s \to -\infty$ limit. The latter are denoted respectively by $\mathfrak{c}_+$ and $\mathfrak{c}_-$. Since the function $\mathrm{X}$ takes integer values on the solutions to (1.13), the operative assumption in what follows is that $\mathrm{X}(\mathfrak{c}_+) \le \mathrm{X}(\mathfrak{c}_-) - 1$.

The function $s \to \mathfrak{a}(\mathfrak{d}_{\mathfrak{p}_n}|_s) + \mathfrak{p}_n(\mathfrak{d}_{\mathfrak{p}_n}|_s)$ is a non-increasing function on $\mathbb{R}$ and as the sequence $\{\mathfrak{p}_n\}_{n=1,2,\ldots}$ is bounded and converges to zero, the fact that the set of $\mathfrak{d} = \mathfrak{d}_{\mathfrak{p}_n}$ versions of $A_{\mathfrak{d}}$ is bounded implies that the sequence $\{\mathfrak{d}_{\mathfrak{p}_n}\}_{n=1,2,\ldots}$ has a subsequence that converges in the sense described in Chapter 16 of [KM] to what is said in Definition 16.1.2 of [KM] to be a *broken trajectory*. In the situation here, such a trajectory consists of a non-empty, finite, ordered set $\{\mathfrak{d}_k\}_{k=1,2,\ldots,N}$ of instanton solutions to (4.1) with the following property: The $s \to \infty$ limit of $\mathfrak{d}_k$ is the $s \to -\infty$ limit of $\mathfrak{d}_{k+1}$ for $k < N$. In addition, $\mathfrak{c}_-$ is the $s \to -\infty$ limit of $\mathfrak{d}_1$ and $\mathfrak{c}_+$ is the $s \to \infty$ limit of $\mathfrak{d}_N$. This being the case, $\mathrm{X}(\mathfrak{c}_+) - \mathrm{X}(\mathfrak{c}_-)$ can be written as $\sum_{k=1,2,\ldots N} (\lim_{s\to\infty} \mathrm{X}(\mathfrak{d}_k|_s) - \lim_{s\to-\infty} \mathrm{X}(\mathfrak{d}_k|_s))$. This sum is non-negative if each $\mathfrak{d}_k$ obeys $A_{\mathfrak{d}_k} \le c(r\ln r + 1)$ so as to invoke Proposition 7.1.

To see that this last bound is obeyed, use Lemma 4.2 to see that each $\mathfrak{d} = \mathfrak{d}_k$ version of $A_{\mathfrak{d}}$ is positive, so the desired bound holds if it holds for $\sum_{1\le k\le N} A_{\mathfrak{d}_k}$. Meanwhile, the fact that the function on $\mathbb{R}$ given by the rule $s \to \mathfrak{a}(\mathfrak{d}_{\mathfrak{p}_n}|_s) + \mathfrak{p}_n(\mathfrak{d}_{\mathfrak{p}_n}|_s)$ is non-increasing and the fact that $\{\mathfrak{p}_n\}_{n=1,2,\ldots}$ converges to zero implies that the $n \to \infty$ limit of the set of $\mathfrak{d} = \mathfrak{d}_{\mathfrak{p}_n}$ versions of $A_{\mathfrak{d}}$ exists. Moreover, the manner of convergence of $\{\mathfrak{d}_{\mathfrak{p}_n}\}_{n=1,2,\ldots}$ to $\{\mathfrak{d}_k\}_{k=1,2,\ldots}$ as described in Chapter 16 of [KM] guarantees that the limit of the corresponding set of $\mathfrak{d} = \mathfrak{d}_{\mathfrak{p}_n}$ versions of $A_{\mathfrak{d}}$ is no less than $\sum_{1\le k\le N} A_{\mathfrak{d}_k}$. In fact, the limit equals the sum if all solutions to the $(r,\mu)$ version of (1.13) are non-degenerate.

The conclusion that $\mathrm{X}(\mathfrak{c}_+) - \mathrm{X}(\mathfrak{c}_-) \ge 0$ violates the assumptions and so constitutes the desired nonsense.

**b) Proof of Proposition 7.1**

The proof has six parts. By way of notation, $c_c$ is used in what follows to denotes a constant which is greater than 1 that is detemined solely by $c$, $\gamma$ and the geometry of Y. In particular, $c_c$ does not depend on $\mathfrak{d}$, nor does it depend on the chosen values of r or $\mu$. The value of $c_c$ can be assumed to increase between successive appearances.

*Part 1*: Let $\mathfrak{c}_+ = (A_+, \psi_+)$ denote the $s \to \infty$ limit of $\mathfrak{d}$ and let $\mathfrak{c}_- = (A_-, \psi_-)$ denote the corresponding $s \to -\infty$ limit of $\mathfrak{d}$. Both $\hat{A}_+$ and $\hat{A}_-$ are flat with trivial holonomy on a fixed, but small radius tubular neighborhood of $\gamma$ if r is greater than a purely $\gamma$-dependent constant. This fact implies that $\mathrm{X}_\gamma(A_+) - \mathrm{X}_\gamma(A_-) \in \mathbb{Z}$.

With the preceding in mind fix for the moment a smooth, closed 2-form on Y with compact support in this tubular neighborhood whose De Rham cohomology class is the image of the Poincaré dual of the class in $H_1(Y;\mathbb{Z})$ that is defined by viewing the oriented loop $\gamma$ as a 1-cycle. Use $\upsilon_\gamma$ to denote the chosen 2-form.

Re-introduce $\hat{A}$ to denote the connection that is defined by A using the formula in (1.15) with it understood that $\nabla_A\alpha$ has $ds$ component equal to $\frac{\partial}{\partial s}\alpha$. The curvature of this connection is depicted in (6.25). Stokes' theorem writes $x_\gamma(A_+) - x_\gamma(A_-)$ as the integral over $\mathbb{R} \times Y$ of the curvature 2-form $\frac{i}{2\pi} F_{\hat{A}} \wedge \upsilon_\gamma$.

By way of a parenthetical remark, if $A_0 \leq c\,r$, then Proposition 5.1 and Proposition 6.3 can be invoked if r is greater than a purely $c$ dependent constant. Assume this to be the case. It follows from Proposition 6.3 and what is said in Part 4 of the proof of Lemma 6.2 that the integral of $\frac{i}{2\pi} F_{\hat{A}} \wedge \upsilon_\gamma$ is a weighted, algebraic count with positive weights of the intersections between the submanifold $\mathbb{R} \times \gamma$ and a pseudoholomorphic subvariety that is defined in some neighborhood of $\mathbb{R} \times \gamma$. Thus $x_\gamma(A_+) - x_\gamma(A_-) \geq 0$.

The equality between $x_\gamma(A_+) - x_\gamma(A_-)$ and the integral of $\frac{i}{2\pi} F_{\hat{A}} \wedge \upsilon_\gamma$ does not depend on the chosen version of $\upsilon_\gamma$ as long as its support is in a radius $c_0^{-1}$ tubular neighborhood of $\gamma$. This being the case, the remainder of this Part 1 defines a useful choice. To this end, let T denote a radius $c_0^{-1}$ tubular neighborhood of $\gamma$ with coordinates of the sort that are described in (7.1). Assume that T appears in these coordinates as $\mathbb{R}/(\ell\,\mathbb{Z}) \times D_0$ where $D_0$ is a radius $c_0^{-1}$ disk centered at the origin in $\mathbb{C}$.

The desired version of $\upsilon_\gamma$ is constructed with the help of a nonnegative, non-increasing function on $[0, \infty)$ with support in $[0,2]$. The latter is denoted in what follows by $q$ and it has the following properties

- $q(x) = 1$ *where* $x < 1$ *and* $q(x) = 0$ *where* $x > 2$.
- $q(x) = \dfrac{e^{-1/(2-x)}}{1 - e^{-1/(1-x)} + e^{-1/(2-x)}}$ *where* $x \in [1, 2]$.

(7.2)

Use q to denote the integral of the function $x \to 2\pi x\, q(x)$. With $q$ in hand, fix D so as to be greater than 100 times the inverse of the radius of $D_0$. The value of D can be taken smaller than a constant that is determined ultimately by $c$ and $\gamma$. With D fixed, let $q_D$ denote the function on $D_0$ that is defined by the rule $z \to q^{-1} D^2 q(Dz)$. Use the coordinates for T in (7.1) to view $q_D$ as a function on Y with compact support in T.

The desired version of $\upsilon_\gamma$ is defined to be zero on the complement of T and defined to equal $q_D w$ on T. So defined, the condition in the third bullet of (7.1) guarantees that $\upsilon_\gamma$ is closed. This understood, it follows from the definition that its De Rham cohomology class is the image in De Rham cohomology of the Poincaré dual of $\gamma$'s class in $H_1(Y;\mathbb{Z})$.

*Part 2*: The upcoming Lemma 7.2 refers to certain notions that are defined directly. The first of these is $\rho_r$, this used to denote $(\ln r)^4$. To define the rest, fix $s_0 \in \mathbb{R}$, $t_0 \in \mathbb{R}/(\ell\mathbb{Z})$ and a point $z_0 \in D_0$ with $|z_0|$ at most half the diameter of $D_0$. The lemma uses $Q^{(s_0,t_0,z_0)}$ to denote the set in $\mathbb{R}\times\mathbb{R}/(\ell\mathbb{Z})\times D_0$ whose $(s,t,z)$ coordinates obey the two conditions $|s - s_0| \le 2\rho_r$ and $|t-t_0|+|z-z_0| < 4\rho_r$. The integral of $iF_{\hat{A}} \wedge (ds \wedge \hat{a} + w)$ over $Q^{(s_0,t_0,z_0)}$ is denoted by $\Delta_{(s_0,t_0,z_0)}$. The lemma uses $Q_{(s_0,z_0)}$ to denote $\cup_{t\in\mathbb{R}/(\ell\mathbb{Z})} Q^{(s_0,t,z_0)}$; this is the set whose $(s,t,z)$ coordinates obey the two constraints $|s - s_0| < 2\rho_r$ and $|z-z_0| < 4\rho_r$ with no constraint on t.

**Lemma 7.2**: *Given $c > 1$ there exists $\kappa_c > 1$ with the following significance: Fix $r \ge \kappa_c$ and $\mu \in \Omega$ with $\mathcal{P}$-norm less than 1 and suppose that $\mathfrak{d} = (A, \psi)$ is an instanton solution to the $(r,\mu)$ version of (4.1) with $A_{\mathfrak{d}} < c\, r \ln r$. Fix $s_0 \in \mathbb{R}$ and $z_0 \in D_0$ with $|z_0|$ less than half the radius of $D_0$. If $\sup_{t\in\mathbb{R}/(\ell\mathbb{Z})} \Delta_{(s_0,t,z_0)} > \kappa_c \rho_r^2$, then*

$$\int_{Q_{(s_0,z_0)}} iF_{\hat{A}} \wedge w > \kappa_c^{-1} \rho_r \sup_{t\in\mathbb{R}/(\ell\mathbb{Z})} \Delta_{(s_0,t,z_0)} .$$

***Proof of Lemma 7.2***: The proof has five steps. These steps use $\Delta_*$ for $\sup_{t\in\mathbb{R}/(\ell\mathbb{Z})} \Delta_{(s_0,t,z_0)}$.

<u>Step 1</u>: This step states five observations that play central roles in the subsequent steps. The first observation constitutes the next equation. This writes $\frac{\partial}{\partial s}A$ and $B_A$ as in (5.7). These equations with (6.25), (4.1) and Lemma 4.6 lead to

- $*(iF_{\hat{A}} \wedge w) = (1-\wp)(1-\sigma)r(1-|\alpha|^2) + 2\wp'|(\nabla^{(1,0)}\alpha)_0|^2 + \mathfrak{e}_A$ ,
- $*(iF_{\hat{A}} \wedge ds \wedge \hat{a}) = (1-\wp)\sigma\, r(1-|\alpha|^2) + 2\wp'|(\nabla^{(1,0)}\alpha)_1|^2 + \mathfrak{e}_B$ ,

(7.3)

with the notation as follows. What is denoted by $\sigma$ is the function that appears in (5.7). To define $(\nabla^{(1,0)}\alpha)_0$ and $(\nabla^{(1,0)}\alpha)_1$, first introduce $\nabla^{(1,0)}\alpha$ to denoted $T^{1,0}(\mathbb{R}\times Y)$ part of the covariant derivative of $\alpha$. View the latter as a homomorphism from the (1,0) summand in $T(\mathbb{R}\times Y)\otimes\mathbb{C}$ to E. What is denoted by $(\nabla^{(1,0)}\alpha)_0$ is the restriction of this homomorphism to the span of $\frac{\partial}{\partial s} - iv$; and what is denoted by $(\nabla^{(1,0)}\alpha)_1$ is the restriction of $\nabla^{(1,0)}\alpha$ to the $+i$ eigenspace of J in the $K^{-1}\otimes\mathbb{C}$ summand in $T(\mathbb{R}\times Y)\otimes\mathbb{C}$. Meanwhile, $\mathfrak{e}_A$ and $\mathfrak{e}_B$ are such that their absolute values are bounded by $c_0((1-\wp)+\wp')$. The $(1-\wp)$ contribution to latter bounds follow from the bounds on $\mathfrak{z}_A$ and $\mathfrak{z}_B$ in (6.7), and the $\wp'$ contribution follows by using (4.1) to write the $T^{0,1}$ part of $\nabla_A\alpha$ as a linear combination of covariant derivatives of $\beta$ and then invoking Lemma 4.6.

Adding the two equalities in (7.3) leads to the second observation:

$$*(iF_{\hat{A}} \wedge (ds \wedge \hat{a} + w)) \geq \tfrac{1}{4}(1-\wp)r + 2\wp'|\nabla^{(1,0)}\alpha|^2 + \mathfrak{e}\ , \tag{7.4}$$

where $\mathfrak{e} = 0$ where $\wp = 1$ and $|\mathfrak{e}| \leq c_0 r^{-3/2}$ in any event. By way of an explanation, this inequality follows from the fact that $|\alpha| < \frac{9}{16}$ on the support of $(1-\wp)$ and from the fact that $\wp' \leq c_0(1-\wp)^{3/4}$. To elaborate, use the bound $\wp' \leq c_0(1-\wp)^{3/4}$ to see that the contribution to the right hand side of (7.4) from the $\mathfrak{e}_A$ and $\mathfrak{e}_B$ terms in (7.3) is larger than $r^{1/2}(1-\wp)$ only where $(1-\wp)$ is less than $c_0 r^{-3/2}$. At points where such is the case, the bound $\wp' \leq c_0(1-\wp)^{3/4}$ with the fact that $|\alpha| < \frac{9}{16}$ on the support of $(1-\wp)$ implies that $\frac{1}{16}(1-\wp)r(1-|\alpha|^2) + c_0\wp' r^{1/2}$ is no greater than $c_0 r^{-3/2}$.

The third observation is stated by the next inequality:

$$\int_{Q_{(s_0,z_0)}} iF_{\hat{A}} \wedge w \geq -c_c(r^{-1/(q+1)}\Delta_* + r^{-3/2})\ . \tag{7.5}$$

This is a consequence of (7.4), the first bullet of (7.3) and the fact that (5.7)'s function $\sigma$ is less than $1 + c_c r^{-1/q}$ The point being that the expression on the right side of the first bullet in (7.3) is negative only where $\sigma$ is greater than 1 or the $\mathfrak{e}_A$ term is negative and dominates the others. The derivation of (7.5) also uses the fact that the length of $\gamma$ when measured in units where $\rho_r = 1$ is $\ell(\ln r)^4$; and it then uses the bound $(\ln r)^4 r^{-1/q} < c_0 r^{-1/(1+q)}$.

The fourth observation is a direct corollary to (7.5):

*Fix* $t \in \mathbb{R}/(\ell\mathbb{Z})$. *If* $m > 0$ *and if the integral of* $iF_{\hat{A}} \wedge w$ *over* $Q^{(s_0,t,z_0)}$ *is greater than* $m + c_c r^{-1/(q+1)}(\Delta_* + 1)$*, then* $\int_{Q_{(s_0,z_0)}} iF_{\hat{A}} \wedge w > m$.

(7.6)

The fifth observation is a tautology that comes by writing $iF_{\hat{A}} \wedge (ds \wedge \hat{a} + w)$ as the sum of $iF_{\hat{A}} \wedge w$ and $iF_{\hat{A}} \wedge ds \wedge \hat{a}$. To set the notation, fix $t_0 \in \mathbb{R}/(\ell\mathbb{Z})$ such that $\Delta_{(s_0,t_0,z_0)} = \Delta_*$. Use $Q^*$ to denote $Q^{(s_0,t_0,z_0)}$.

*If the integral of* $iF_{\hat{A}} \wedge w$ *over* $Q^*$ *is less than* $\frac{1}{2}\Delta_*$,
*then the integral of* $i{*}B_{\hat{A}} \wedge ds \wedge \hat{a}$ *over* $Q^*$ *is greater than* $\frac{1}{2}\Delta_*$.

(7.7)

This is so because $iF_{\hat{A}} \wedge ds \wedge \hat{a} = i{*}B_A \wedge ds \wedge \hat{a}$.

Step 2: This step outlines the argument that is used to prove Lemma 7.2. Assume that $\Delta_* > \rho_r^2$. It follows from (7.6) that the assertion of Lemma 7.2 is true if the integral of $iF_{\hat{A}} \wedge w$ over $Q^*$ is greater than $\frac{1}{100}\Delta_*$ so assume that this is not the case. Use (7.7) to see that the integral of $i{*}B_{\hat{A}} \wedge ds \wedge \hat{a}$ over $Q^*$ is greater than $\frac{1}{2}\Delta_*$.

The constant $(s,t)$ slices of $Q^*$ are disks that lie either in a cross-sectional sphere of $\mathcal{H}_0$ or a level set of $f$ in $M_\delta$. The former are compact surfaces without boundary, and so are most of the latter. The integral of $i{*}B_{\hat{A}}$ over a surface of this sort without boundary is $2\pi$ times the pairing of the first Chern class of E with the homology class defined by the surface and this no greater than $2\pi G$. Thus, the integral of $i{*}B_{\hat{A}}$ over such a surface is a priori bounded by $2\pi G$. If the integral of $i{*}B_{\hat{A}}$ over a disk in one of these surfaces is greater than this bound, then there must be other parts of the surface where the corresponding integral is negative.

The second bullet in (7.3) implies that the pull-back of $i{*}B_{\hat{A}}$ to such a surface is no smaller than $-c_c r^{1-1/q}$ times the area form, this being the pull-back of $w$. As a consequence, the area where this pull-back is negative can be no smaller than $c_c r^{-1+1/q}$. Meanwhile, the two bullets of (7.3) imply the following: The pull-back of $i{*}B_{\hat{A}}$ to the surface is negative at a point only if ${*}(ds \wedge \frac{\partial}{\partial s}A \wedge w)$ is of order r. This implies that $|\frac{\partial}{\partial s}A|^2$ is order $r^2$ and so the integral of the latter function on the surface is no less than $r^{1+1/q}$.

As explained below, the extra factor $r^{1/q}$ leads to a violation of Lemma 4.2. This violation is avoided only if $\Delta_*$ is a priori bounded by $\kappa_c \rho_r^2$.

Step 3: Given $(s, t) \in \mathbb{R} \times \mathbb{R}/(\ell\mathbb{Z})$ with $|s - s_0| < 2\rho_r$ and $|t - t_0| < 4\rho_r$, introduce by way of notation $D_{(s,t)}$ to denote the slice at $(s, t)$ of $Q^*$; and use $E(s,t)$ to denote the integral of $\frac{i}{2\pi}{*}B_A$ over the disk $\{(s,t)\} \times D_{(s,t)}$. Given $n \in \{0,1,2,\ldots\}$, let $U_n \subset \mathbb{R} \times \mathbb{R}/(\ell\mathbb{Z})$ denote the set of points with $|s - s_0| < 2\rho_r$, $|t - t_0| < 4\rho_r$ and $E(s,t) \in [nG, 2nG)$. Use $v_n$ to denote the measure of $U_n$.

Suppose for the remainder of this step that $\gamma(t_0)$ is either in $\mathcal{H}_0$ or in the part of $M_\delta$ where $f$ is either less than $1 - 4\delta^2$, or between $1 + 2\delta^2$ and $2 - 2\delta^2$, or greater than $2 + 4\delta^2$. This assumption has the following implication: No point in $Q^*$ is on a level set of $f$ that enters a radius $\delta$ coordinate ball centered on either an index 1 or index 2 critical point of $f$ in M. Keep this fact in mind.

Fix $n \in \{2, \ldots\}$ such that $U_n \neq \emptyset$ and fix $(s, t) \in U_n$. The disk $\{(s, t)\} \times D_{(s,t)}$ lies in a compact surface in $\mathbb{R} \times (M_\delta \cup \mathcal{H}_0)$ whose tangent space is annihilated by $\hat{a}$. This surface is either in a component of a level set of $f$ in $M_\delta$ or a cross-sectional sphere of $\mathcal{H}_0$. Use $S_{(s,t)}$ to denote this surface. The integral of $\frac{i}{2\pi}{*}B_{\hat{A}}$ over $S_{(s,t)}$ is equal to G if $S_{(s,t)}$ is an $f$-level set with $f \in (1 + \delta^2, 2 - \delta^2)$ part of $M_\delta$; it is equal to 0 otherwise.

Since $\mathrm{E}(s,t) \geq 2\mathrm{G}$, the integral of $\frac{i}{2\pi} *B_{\hat{A}}$ over $S_{(s,t)} - D_{(s,t)}$ must be less than -G. To see what this entails, write the pull-back of $\frac{i}{2\pi} *B_{\hat{A}}$ to $S_{(s,t)}$ as $\mathrm{B}\,w$ and use (5.7) to see that B is the function that appears on the right hand side of (7.3). The latter function is no less than $-c_c((1-\wp)\, r^{1-1/q} + r^{-3/2})$. The factor $-c_c r^{-3/2}$ contributes no more than $-c_c r^{-3/2}$ to the integral of $\mathrm{B}\,w$ and this is no more than $-10^{-6}$ if $r > c_c$. Assume this to be the case.

What was just said implies that the measure of the set in $S_{(s,t)} - D_{(s,t)}$ where B is such that $(1-\wp)\sigma < 0$ is no less than $c_c^{-1}(n-1)\mathrm{G}\,r^{-1+1/q}$. Noting that $|\alpha| < \frac{3}{4}$ on this set, it follows from the first bullet in (7.3) that $*(i\,ds \wedge \frac{\partial}{\partial s}A \wedge w) > c_c^{-1} r$ on this same set. This implies in particular that $|\frac{\partial}{\partial s}A|^2 > c_0 r^2$ on a set of measure greater than $c_c^{-1}(n-1)\mathrm{G}\,r^{-1+1/q}$ in $S_{(s,t)}$. It follows as a consequence that the integral of $|\frac{\partial}{\partial s}A|^2$ over $(\cup_{(s,t)\in U_n} S_{(s,t)}) \subset \mathbb{R} \times Y$ is no less than $c_c^{-1} r^{1+1/q}(n-1)\mathrm{G}\,V_n$.

<u>Step</u> <u>4</u>: Suppose that $f(\gamma(t_0))$ is between $1-4\delta^2$ and $1+2\delta^2$ or else between $2-2\delta^2$ and $2+4\delta^2$. Fix $(s,t) \in \mathbb{R} \times \mathbb{R}/(\ell\mathbb{Z})$ with $|s-s_0| < 2\rho_r$ and $|t-t_0| \leq 4\rho_r$. If $t > t_0$, use $T_{(s,t)}$ to denote the set of points in $\mathbb{R} \times T$ that have coordinates $(s, \tau, z)$ with $\tau \in [t, t_0 + 8\delta]$ and with z such that $|z-z_0| + |t-t_0| \leq 4\rho_r$. If $t \leq t_0$, define $T_{(s,t)} \subset \mathbb{R} \times T$ to be the set that have coordinates $(s, \tau, z)$ with $\tau \in [t_0 - 8\delta, t]$ and with z as in the $t > t_0$ case. In either case, $T_{(s,t)}$ is a manifold with corners. In the $t > t_0$ case, there are three codimension 1 faces of $T_{(s,t)}$. These are the disks $D_{(s,t)}$ and $D_{(s,t_0+8\delta)}$, and third is the cylinder consisting of the points $(s,\tau, z)$ with $\tau \in [t, 8\delta + t_0]$ and $|z-z_0| + |t-t_0| = 4\rho_r$. There is a similar story when $t < t_0$; The cylinder face of $T_{(s,t)}$ in this case is the set of points $(s, \tau, z)$ with $\tau \in [-8\delta + t_0, t]$ and with z as in the $t > t_0$ case. In either case, let $C_{(s,t)}$ denote the cylindrical face of $T_{(s,t)}$.

Use Stoke's theorem to see that the absolute value of the difference between the integral of $i*B_{\hat{A}}$ over the two disk faces of $T_{(s,t)}$ is equal to the absolute value of the integral of $i*B_{\hat{A}}$ over $C_{(s,t)}$. This integral involves only the $\mathfrak{r}$ - $\mathfrak{X}$ part of (5.7)'s depiction of $B_A$. Meanwhile, the contribution from the term proportional to $\wp'$ in (6.25) is bounded by $c_0\, \wp'(|(\nabla^{(1,0)}\alpha)_0|\,|(\nabla^{(1,0)}\alpha)_1| + 1)$. With these last facts in mind, introduce by way of notation

$$N = \int_{\{(s,t):|s-s_0|<2\rho_r \text{ and } |t-t_0|<4\rho_r\}} \Big( \int_{C_{(s,t)}} ((1-\wp)(|\mathfrak{X}|+|\mathfrak{r}|) + \wp'(|(\nabla^{(1,0)}\alpha)_0|\,|(\nabla^{(1,0)}\alpha)_1|+1)) \Big)\, ds \wedge dt . \tag{7.8}$$

Let $Q^{\Diamond}$ denote either $Q^{(s_0, t_0+8\delta, z_0)}$ or $Q^{(s_0, t_0-8\delta, z_0)}$. What was said above about Stoke's theorem implies that

$$Q = | \int_{Q^*} (i * B_{\hat A} \wedge ds \wedge \hat a) \ - \int_{Q^\diamond} (i * B_{\hat A} \wedge ds \wedge \hat a) | \ < c_0 N \ .$$

(7.9)

Use $T^*$ to denote the union of the $t \in [t_0 - 8\delta, t_0 + 8\delta]$ versions of $Q^{(s_0, t, z_0)}$. Fix for the moment $e > 1$. Change variables in the integration that defines N and use (5.8) with the fact that $\sigma > -c_c r^{-1/q}$ and $(1 - \sigma) > -c_c r^{-1/q}$ to see that N is no greater than

$$e^{-1} \int_{T^*} ds \wedge \hat a \wedge i * B_{\hat A} \ + \ c_c e \int_{T^*} iF_{\hat A} \wedge w \ + c_c r^{-1/(q+1)} \Delta_* \ .$$

(7.10)

The left most integral in (7.10) is no greater than $c_0 \rho_r^{-1} \Delta_*$ and so the left most term in (7.10) is no greater than $c_0 e^{-1} \rho_r^{-1} \Delta_*$. Therefore, if $e = 1000 c_0 \rho_r^{-1}$ and if $Q > \frac{1}{100} \Delta_*$, then

$$\int_{T^*} iF_{\hat A} \wedge w > c_c^{-1} \rho_r \Delta_* \ .$$

(7.11)

when $r > c_c$. If (7.11) holds with $r > c_c$, then what is said by Lemma 7.2 is true, this being a consequence of (7.5).

Assume that (7.11) does not hold. Then $Q < \frac{1}{100} \Delta_*$ when $r \ge c_c$ and so a repetition of Step 2 with $(s_0, t_0, z_0)$ replaced by either $(s_0, t_0 + 8\delta, z_0)$ or with $(s_0, t_0 - 8\delta, z_0)$ supplies a lower bound for the integral of $|\frac{\partial}{\partial s} A|^2$ over $(\cup_{(s,t) \in U_n} S_{(s,t)}) \subset \mathbb{R} \times Y$ for each $n \in \{2, \ldots\}$.

<u>Step</u> <u>5</u>: Sum the bounds from Step 3 or Step 4 for the integral of $|\frac{\partial}{\partial s} A|^2$ over $(\cup_{(s,t) \in U_n} S_{(s,t)}) \subset \mathbb{R} \times Y$ for $n = 2, 3, \ldots$ to see that the integral over $\mathbb{R} \times Y$ of $|\frac{\partial}{\partial s} A|^2$ is no less than $(c_c^{-1} \Delta_* - 2G V) r^{1+1/q}$ where V is the upper bound for the various $(s,t)$ versions of the sum $V_0 + V_1$. Since V is at most than $16\pi \rho_r^2$, this implies that

$$\int_{\mathbb{R} \times Y} | \frac{\partial}{\partial s} A |^2 \ge (c_c^{-1} \Delta_* - c_0 G \rho_r^2) \, r^{1+1/(2q)} \ .$$

(7.12)

What with Lemma 7.2's assumption about $A_0$, this last inequality runs afoul of what is asserted by Lemma 4.1 if $\Delta_* > c_c \rho_r^2$ and $r > c_c$.

*Part 3*: The next lemma is an analog of Lemma 7.2 for pairs $(s_0, z_0) \in \mathbb{R} \times D_0$ whose corresponding $\Delta_*$ is small. Lemma 7.3 uses the following notation: Given $x \in \mathbb{R}$ and $\rho \in (r^{-1/2}, c_0^{-1})$, the lemma uses $\hat M_{(x,\rho)}$ to denote the integral of $iF_{\hat A} \wedge (ds \wedge \hat a + w)$ over the ball of radius $\rho$ in $\mathbb{R} \times Y$ centered at x. The lemma also introduces $T_{(s_0, z_0)}$ to denote the radius $\rho_r$ tubular neighborhood of the $s = s_0$ and $z = z_0$ slice of $\mathbb{R} \times T$, this being the set of points of the form $(s, t, z)$ with $(s - s_0)^2 + |z - z_0|^2 < \rho_r^2$ and with $t \in \mathbb{R}/(\ell\mathbb{Z})$.

**Lemma 7.3**: *Given* $c > 1$ *there exists* $\kappa_c > 1$ *with the following significance: Fix* $r \geq \kappa_c$ *and* $\mu \in \Omega$ *with* $\mathcal{P}$*-norm less than 1 and let* $\mathfrak{d} = (A, \psi)$ *denote an instanton solution to the* $(r, \mu)$ *version of (4.1) with* $A_{\mathfrak{d}} < c\, r \ln r$. *Fix* $(s_0, t_0) \in \mathbb{R} \times \mathbb{R}/(\ell\mathbb{Z})$ *and a point* $z_0 \in D_0$ *with* $|z_0|$ *less than one fourth the radius of* $D_0$. *If* $\hat{M}_{((s_0, t_0, z_0), \rho_r)} \neq 0$, *then there exists* $s_1 \in \mathbb{R}$ *and* $z_1 \in D_0$ *with* $|s_1 - s_0| < 8\rho_r$ *and* $|z_1 - z_0| < 8\rho_r$ *with* $\int_{T_{(s_1, z_1)}} iF_{\hat{A}} \wedge w > \kappa_c^{-1} \rho_r^4$.

The proof of Lemma 7.3 invokes an $\hat{A}$ analog of Lemma 4.10, this being

**Lemma 7.4**: *There exists* $\kappa \geq \pi$, *and given* $q \geq 1$, *there exists* $\kappa_q \geq 1$ *with the following significance: Fix* $r \geq \kappa$ *and* $\mu \in \Omega$ *with* $\mathcal{P}$*-norm less than 1. Let* $\mathfrak{d} = (A, \psi)$ *denote an instanton solution to the* $(r, \mu)$ *version of (4.1) with* $A_{\mathfrak{d}} < r^2$ *and* $\sup_{s \in \mathbb{R}} \underline{M}(s) \leq r^{1-1/q}$. *Suppose that* $x \in \mathbb{R} \times Y$ *is a point where* $|\alpha| \leq \frac{3}{4}$.

- *If* $\rho_1 > \rho_0$ *are in* $(\kappa_q r^{-1/2}, \kappa_q^{-1})$, *then* $\hat{M}_{(x, \rho_1)} \geq \kappa_q^{-1} \rho_1^2/\rho_0^2\, \hat{M}_{(x, \rho_0)}$.
- *If* $\rho \in (\kappa_q r^{-1/2}, \kappa_q^{-1})$, *then* $\hat{M}_{(x, \rho)} \geq \kappa^{-1} \rho^2$.

It follows from Lemmas 2.5, 4.1 and 5.2 that the assumptions of Lemma 7.3 are met using $q > 6$.

***Proof of Lemma 7.4***: Given Lemmas 4.4 and the first bullet of 4.8 and the formula in (6.25) for $F_{\hat{A}}$, the proof of Proposition 3.1 in the article *SW=>Gr* from [T8] can be used with only cosmetic changes to prove the assertions. A second proof deduces Lemma 7.4 from Lemmas 4.8, 4.10 and 4.12 by proving the following assertion:

*There is a purely* q *dependent constant,* $c_q$, *which is greater than 1 and is such that if* $x \in \mathbb{R} \times Y$ *is a point where* $|\alpha| < \frac{1}{2}$ *and* $\rho \in (r^{-1/2}, c_q^{-1})$, *then*

$$c_q^{-1} M_{(x,\rho)} < \hat{M}_{(x,\rho)} < c_q M_{(x,\rho)}. \tag{7.13}$$

What follows is a sketch of the argument for (7.13). To start, fix $m > 100$ and use Lemma 4.12 with Proposition 4.11 to see that the contribution to $M_{(x,\rho)}$ from the set of points where $|\alpha| \geq 1 - m^{-1}$ is greater than $c_q^{-1} \hat{M}_{(x,\rho)}$ if r is greater than a purely $m$ and q dependent constant. This proves the upper bound in (7.13). To prove the lower bound, use Lemma 4.12 and Proposition 4.11 to see that the contribution to $M_{(x,\rho)}$ from this same set is no less than $c_{m_q} \hat{M}_{(x,\rho)}$ if $r > c_{m_q}$ with $c_{m_q} > 1$ being a contant that depends only on $m$ and q. Meanwhile, the assertions in the second and third bullets of Lemma 4.8 can be used to prove the following: If $m \ge c_q$, then the contribution to $M_{(x,\rho)}$ from the complement of the set where $|\alpha| \le 1 - m^{-1}$ is no greater than $c_{m_q}$ times the contribution to $M_{(x,\rho)}$ from the set where $|\alpha| \ge 1 - m^{-1}$.

The assertion of the second bullet of Lemma 7.3 follows from Lemma 4.12 and Proposition 4.11 as they imply that $|\alpha| < \frac{1}{2}$ on a ball in $\mathbb{R} \times Y$ of radius at least $c_0^{-1} r^{-1/2}$ with distance at most $c_0 r^{-1/2}$ from x.

*Part 4*: This part of the Proposition 7.1's proof supplies a proof of Lemma 7.3. By way of notation, the proof uses $\hat{M}_*$ to denote $\sup_{t\in\mathbb{R}/(\ell\mathbb{Z})} \hat{M}_{((s_0,t,z_0),\rho_r)}$. The proof uses the same conventions about $c_c$ that is used in Lemma 7.2's proof; and it introduces one additional convention of the same sort: Given $m > 1$, the proof uses $c_{cm}$ to denote a number that is greater than 1 and depends only on $m$, $c$, $\gamma$ and the geometry of Y. In particular, this number does not depend on $\mathfrak{d}$ nor on r. The value of $c_{cm}$ can be assumed to increase between successive appearances.

***Proof of Lemma 7.3***: The proof has five steps.

<u>Step 1</u>: There exists $N < c_0$ and a set of N points $\{(s_k, z_k)\}_{k=1,2\ldots,N}$ with the following properties: First, $|s_0 - s_k| < 8\rho_r$ and $|z_0 - z_k| < 8\rho_r$ for each $k \in \{1, \ldots, N\}$. Second, the union of the sets $\{T_{(s_k,z_k)}\}_{k=0,1,\ldots,N}$ contains $Q_{(s_0,z_0)}$. This understood, invoke Lemma 7.2 to see that Lemma 7.3's assertion is true unless $\sup_{t\in\mathbb{R}/(\ell\mathbb{Z})} \Delta_{(s_0,t,z_0)} < c_c \rho_r^2$. Indeed, if this isn't the case, then the integral of $iF_{\hat{A}} \wedge w$ over least one $k \in \{0, 1, \ldots, N\}$ version of $T_{(s_k,z_k)}$ will be greater than $c_c^{-1} \rho_r^3$.

Let $t_0 \in \mathbb{R}/(\ell\mathbb{Z})$ denote a point with $\hat{M}_{((s_0,t_0,z_0),\rho_r)} > 0$. If such is the case, then there must be a point in the radius $\rho_r$ ball centered at $(s_0, t_0, z_0)$ where $\wp < 1$ and so a point in this ball where $|\alpha| < \frac{3}{4}$. Let $(s_1, t_1, z_1)$ denote such a point. The operative assumption that $\sup_{t\in\mathbb{R}/(\ell\mathbb{Z})} \Delta_{(s_0,t,z_0)} < c_c \rho_r^2$ requires that $\sup_{t\in\mathbb{R}/(\ell\mathbb{Z})} \hat{M}_{((s_1,t,z_1),\rho_r)} < c_c \rho_r^2$ also. This being the case, it is enough to prove the following assertion:

$$\textit{Given } m > 10\textit{, there exists } c_{cm} > 1 \textit{ such that if } r \ge c_{cm} \textit{ and if } \sup_{t\in\mathbb{R}/(\ell\mathbb{Z})} \hat{M}_{((s_1,t,z_1),\rho_r)} \in (0, m\rho_r^2)\textit{, then } \int_{T_{(s_1,z_1)}} iF_{\hat{A}} \wedge w > c_{cm}^{-1} \rho_r^4. \tag{7.14}$$

The remaining steps prove (7.14).

<u>Step</u> <u>2</u>: Let $T_{1/4} \subset T_{(s_1, z_1)}$ denote the set of points whose $(s, z)$ coordinates are such that $|s - s_1|^2 + |z - z_1|^2 < \frac{1}{16}\rho_r^2$. The assertion below summarizes the content of this step.

$$\textit{If } r > c_{cm} \textit{ then either (8.14) holds or } F_{\hat{A}} = 0 \textit{ on } T_{1/4} \cap \mathcal{H}_0 . \tag{7.15}$$

The proof of (7.15) is given after a digression that follows directly. The proof invokes two key facts that are supplied by this digression.

To start the digression, let x´ denote the point with $(s, t, z)$ coordinates $(s', t', z')$ with $s'$ and $z'$ constrained so that $|s' - s_1|^2 + |z' - z_1|^2 \le \frac{1}{4}\rho_r^2$. The operative assumption in (7.14) requires that $\hat{M}_{(x', \rho_r/2)} \le m\rho_r^2$. Assume in addition that $|\alpha| < \frac{3}{4}$ at x´. The first bullet of Lemma 8.4 requires the bound $\hat{M}_{(x',\rho)} \le c_c m \rho^2$ for all $\rho \in (c_c r^{-1/2}, c_c^{-1})$ if $r \ge c_c$. Fix $R > m^4$, $\varepsilon \in (0, m^{-4})$ and $k \in \{10, 11, \ldots\}$ to invoke Lemma 4.12 with the given q and value of $m$. As will be apparent in the proof of (7.15), choices for R, ε and k that depend only on $c$ and $m$ are sufficient. In any event, with R, ε and k chosen, Lemma 4.12 with the given bound on $\hat{M}_{(x',\rho)}$ will be invoked with it understood that r is greater than a constant that depends only on $m$, $c$ and the chosen values for R, ε and k. Of particular interest here is the fact that the corresponding solution $(A_0, \alpha_0)$ of (4.6) is described by Items a) and b) of the third bullet of Proposition 4.1. The assertions of these two items imply the following:

FACT 1: *There are zeros of* $\alpha_0$ *with distance less than* $c_{cm}$ *from the origin in* $\mathbb{C}^2$ *and so there are zeros of* $\alpha$ *with distance less than* $c_{cm}\, r^{-1/2}$ *or less from* x´.

FACT 2: *Each zero of* $\alpha_r$ *with distance less than* R *from the origin has distance less than* $c_{cm}\varepsilon$ *from a zero of* $\alpha_0$*; and each zero of* $\alpha(0)$ *with distance less than* R *from the origin has distance less than* $c_{cm}\varepsilon$ *from a zero of* $\alpha_r$.

With regards to FACT 1, the assertion about the distance from origin of zeros of $\alpha_0$ follows from three facts: The equations in (4.6) are elliptic modulo the action of $C^\infty(\mathbb{C}^2; S^1)$; the polynomial that defines the zero locus of $\alpha_0$ has a priori bounded degree; and $|\alpha_0(0)| < \frac{3}{4} + \varepsilon$ because $|\alpha_r(0)| < \frac{3}{4}$. FACT 2 follows from the a priori degree bound for the polynomial that defines the zeros of $\alpha_0$. In particular, this fact has the following consequences: All but at most a finite set of affine lines in $\mathbb{C}^2$ intersects $\alpha_0^{-1}(0)$ in a finite set of points. Those that do not have this property are irreducible components of $\alpha_0^{-1}(0)$. Moreover, if a line intersects $\alpha_0^{-1}(0)$ in a finite set of points, then the local degree of each intersection point is positive and their sum is bounded by a purely $m$-dependent constant. Given that these zeros have positive local degree, each such intersection point must have distance less than $c_{cm}\varepsilon$ from an $\alpha_r = 0$ point on the affine line if all points on the line at distance R from the origin have distance $c_{cm}\varepsilon R$ or more from all zeros of $\alpha_0$.

With the digression now over, what follows is the proof of (7.15). To start the argument, suppose that $x \in T_{1/4} \cap \mathcal{H}_0$ and that $F_{\hat{A}} \neq 0$ at x. As $\wp < 1$ at x, so $|\alpha| < \frac{3}{4}$ at x. It follows from FACT 1 that the integral of $iF_{\hat{A}} \wedge (ds \wedge \hat{a} + w)$ over the radius $c_{cm} r^{-1/2}$ ball centered at x is greater than $c_{cm}^{-1} r^{-1}$ and so it follows from Lemma 7.4 that $\hat{M}_{(x,\rho_r/4)}$ is greater than $c_{cm}^{-1} \rho_r^2$. If the integral of $iF_{\hat{A}} \wedge w$ over the radius $\frac{1}{4}\rho_r$ ball centered at x is greater than $\frac{1}{2}\hat{M}_{(x,\rho_r/4)}$, then the conclusions of Lemma 7.3 follow because the integral of $iF_{\hat{A}} \wedge w$ over $T_{(s_0,z_0)}$ is no less than $-c_{cm} r^{1/q} m^2 \rho_r^2$.

Granted what was just said assume that the integral of $iF_{\hat{A}} \wedge w$ over the radius $\frac{1}{4}\rho_r$ ball centered at x is less than $\frac{1}{2}\hat{M}_{(x,\rho_r/4)}$. It then follows that the integral of $iF_{\hat{A}} \wedge (ds \wedge \hat{a})$ over this same ball is greater than $c_{cm}^{-1} \rho_r^2$.

Let $(s_x, t_x, z_x)$ denote the $(s,t,z)$ coordinates of x. Introduce $Q^x$ to denote the subset of $\mathbb{R} \times \mathbb{R}/(\ell\mathbb{Z}) \times D_0$ whose $(s,t,z)$ coordinates obey $|s - s_x| < 2\rho_r$ and $|t - t_x| + |z - z_x| < 4\rho_r$. Given $(s, t)$ with $|s - s_x| < 2\rho_r$ and $|t - t_x| < 4\rho_r$, use $D_{(s,t)}$ to denote the constant $(s,t)$ slice of $Q^x$ and use $E(s,t)$ to denote the integral of $\frac{i}{2\pi} {*B_{\hat{A}}}$ over $D_{(s,t)}$. What follows is a direct consequence of Lemma 4.12 and FACT 1 and FACT 2:

$$\textit{If } E(s,t) \textit{ is greater than } c_{cm}\varepsilon\textit{, then } E(s,t) \textit{ is greater than } 1 - c_{cm}\varepsilon. \tag{7.16}$$

Given $n \in \{1, \ldots\}$, let $U_n \subset \mathbb{R} \times \mathbb{R}/(\ell\mathbb{Z})$ denote the set where the conditions $|s - s_x| < 2\rho_r$, $|t - t_x| < 4\rho_r$ and $E(s,t) \in [n - \frac{1}{2}, n + \frac{1}{2})$ hold. Use $U_0$ to denote the set of points with $(s, t)$ such that $E(s, t) \in [0, \frac{1}{2})$. Use $v_n$ to denote the measure of $U_n$.

Given $(s, t)$ as just defined, let $S_{(s,t)}$ denote the slice of $\mathcal{H}_0 \cap M_\delta$ containing $D_{(s,t)}$. This is a J-holomorphic 2-sphere that has pairing 0 with the first Chern class of E. In particular, the integral of $\frac{i}{2\pi} {*B_{\hat{A}}}$ over $S_{(s,t)}$ is zero. This being the case, an almost verbatim copy of the arguments from Step 3 of the proof of Lemma 7.2 prove that the integral of $|\frac{\partial}{\partial s} A|^2$ over $(\cup_{(s,t)\in U_n} S_{(s,t)})$ is no smaller than $c_c^{-1} r^{1+1/q} v_n$. What with (7.16), an almost verbatim repeat of the arguments in Step 5 of the proof of Lemma 7.2 prove that

$$\int_{\mathbb{R}\times Y} |\tfrac{\partial}{\partial s} A|^2 \geq (c_c^{-1} - c_{cm}\varepsilon)\, \rho_r^2\, r^{1+1/(2q)}\,. \tag{7.17}$$

This runs afoul of Lemma 4.1 if $\varepsilon < c_{cm}^{-1}$ and the latter happens if $\varepsilon < c_{cm}$ and $r > c_{cm}$.

Step 3: This step states two observations that are used in the subsequent steps. To set the stage, fix $\tau \in [\frac{1}{4}, \frac{1}{8}]$ and introduce $T_\tau$ to denote the subset of points in $T_{(s_1, z_1)}$ whose $(s, z)$ coordinates are such that $|s - s_1|^2 + |z - z_1|^2 \le \tau^2 \rho_r^2$. The first observation here is that

$$c_c^{-1} \rho_r^2 < \int_{T_\tau} iF_{\hat{A}} \wedge (ds \wedge a + w) < c_c m \rho_r \ . \tag{7.18}$$

By way of an explanation, the lower bound follows from the version of the top bullet of Lemma 7.4 that takes $x = (s_1, t_1, z_1)$ and $\rho_0 = c_c r^{-1/2}$ and $\rho_1 = \frac{1}{2} \rho_r$. Meanwhile, the upper bound follows from the bound $\hat{M}_* < m\rho_r^2$ and the fact that $T_{1/4}$ can be covered by $c_0 \rho_r^{-1}$ balls of radius $\frac{1}{16} \rho_r$ with centers in $T_{1/4}$.

The second observation concerns the integral of $iF_{\hat{A}} \wedge w$. To say more, let U denote a given open subset of $T_{1/4}$. Then

$$\int_U iF_{\hat{A}} \wedge w \ge -c_c m \rho_r r^{-1/q} \ . \tag{7.19}$$

Given the upper bound in (7.18), this follows from the fact that the function $1 - \sigma$ that appears in (5.7) is no less than $-c_c r^{-1/q}$.

Step 4: Assume for this and the remaining steps that $F_{\hat{A}} = 0$ on $T_{1/4} \cap \mathcal{H}_0$. This being the case, the $\mathbb{R}/(\ell\mathbb{Z})$ parameter t on $T_{1/4}$ can be lifted to an $\mathbb{R}$ valued parameter on a neighorhood of the support of $|F_{\hat{A}}|$ and nothing is lost by assuming that the now $\mathbb{R}$-valued parameter t is constrained to an interval $I \subset \mathbb{R}$ of the form $[-c_0, c_0]$ at points in $T_{1/4}$ with distance 1 or less from the support of $|F_{\hat{A}}|$. Meanwhile, it follows from (7.1) that the 2-forms $ds \wedge \hat{a}$ and $w$ on $\mathbb{R} \times I \times D$ can be written as

$$ds \wedge \hat{a} = d(-t\,ds) \quad \textit{and} \quad w = \tfrac{i}{4} d(z\,d\bar{z} - \bar{z}\,dz + \cdots) \tag{7.20}$$

where the unwritten terms in the formula for $w$ have no $ds$ component and are bounded in absolute value by $|z|^2$.

Fix $\tau \in [\frac{1}{8}, \frac{1}{4}]$ and use Stoke's theorem with (7.18) so see that

$$\int_{T_\tau} iF_{\hat A} \wedge ds \wedge \hat a = \int_{(\partial I\times Y)\cap T_\tau} i * B_{\hat A} \wedge t\, ds \ . \tag{7.21}$$

Look at (6.25) and (5.7) to see that the absolute value of the right hand side of (7.21) is no greater than

$$c_0 \int_{(\partial I\times Y)\cap T_\tau} ((1-\wp)(|\mathfrak{X}|+|\mathfrak{r}|) + \wp'(|(\nabla^{(1,0)}\alpha)_0||(\nabla^{(1,0)}\alpha)_1|+1)) \ . \tag{7.22}$$

Fix $e > 64$ to be determined shortly and use (7.21), (7.22) with (5.7), (5.9) and the bound on $\hat M_*$ to see that

$$\int_{[\frac{1}{8},\frac{1}{4}]} (\int_{T_\tau} iF_{\hat A} \wedge ds \wedge \hat a)\, d\tau \le c_0 e \int_{T_{1/4}-T_{1/8}} iF_{\hat A} \wedge w + e^{-1} \int_{T_{1/4}-T_{1/8}} iF_{\hat A} \wedge ds \wedge \hat a + c_0 r^{-1/2} \ . \tag{7.23}$$

This last inequality is the input for the final step in the proof of (7.14).

<u>Step</u> <u>5</u>: There are two cases to consider with regards to (7.23). The first is that when the integral on the left hand side of (7.23) is less than

$$e^{-1/4} \int_{[\frac{1}{8},\frac{1}{4}]} (\int_{T_\tau} iF_{\hat A} \wedge (ds \wedge \hat a + w)\, d\tau \ . \tag{7.24}$$

If this is the case, then there exists $\tau \in [\frac{1}{4}, \frac{1}{8}]$ such that

$$\int_{T_\tau} iF_{\hat A} \wedge w > \tfrac{1}{2} \int_{T_\tau} iF_{\hat A} \wedge (ds \wedge a + w) \ . \tag{7.25}$$

Use the lower bound in (7.18) to see that the integral on the left hand side of (7.25) is no less than $c_0^{-1}\rho_r^2$. Thus $\int_{T_\tau} iF_{\hat A} \wedge w \ge c_0^{-1} \rho_r^2$. This with the bound in (7.19) implies what is asserted by (7.14).

The other case to consider is that where the left hand side of (7.23) is greater than what is written in (7.25). It follows from the lower bound in (7.18) that what is written in (7.24) is greater than $c_0 e^{-1} \rho_r^2$. Meanwhile, the term on the right hand side of (7.23) that is proportional to the integral of $iF_{\hat A} \wedge ds \wedge \hat a$ is bounded by $c_0 e^{-1} \int_{T_{1/4}} iF_{\hat A} \wedge (ds \wedge \hat a + w)$. The upper bound in (7.18) implies that this last expression is no greater than $c_0 e^{-1} m \rho_r$. Granted all of this, then (7.23) implies that

$$\int_{T_{1/4}-T_{1/8}} iF_{\hat{A}} \wedge w > c_0^{-1}\, e^{-5/4} \rho_r^{\,2} - c_0 e^{-2} m \rho_r \ .$$

(7.26)

If $e = \rho_r^{-8/5}$, then the right hand side of (7.26) is greater than $c_0 \rho_r^{\,4}$. This with (7.19) implies what is asserted in (7.14).

*Part 5*: Fix $(s_0, t_0, z_0)$ where the function $q_D$ from Part 1 is non-zero and reintroduce Lemma 7.2's set $Q_{(s_0,z_0)}$ so as to consider the integral

$$\int_{Q_{(s_0,z_0)}} iF_{\hat{A}} \wedge q_D w \ .$$

(7.27)

The assertion that follows summarizes what is said here about (7.27).

*Assume that* r *is greater than a constant that depends only on* D*, c,* $\gamma$ *and the geometry of* Y. *Then the integral in (7.27) is non-negative if* $q_D(z_0) > c_0 e^{-(\ln r)^2/c_0}$ .

(7.28)

To explain (7.28), note first that integral in (7.27) is zero if $\sup_{t\in\mathbb{R}/(\ell\mathbb{Z})} \Delta_{(s_0,t,z_0)}$ is zero, so assume that this is not the case. Use $\Delta_*$ to denote this supremum. Write $q_D$ on $Q_{(s_0,z_0)}$ as $q_D(z_0) + \mathfrak{q}$ with $\mathfrak{q}(z_0) = 0$. The integral in (7.27) has the corresponding decomposition as

$$q_D(z_0) \int_{Q_{(s_0,z_0)}} iF_{\hat{A}} \wedge w + \int_{Q_{(s_0,z_0)}} iF_{\hat{A}} \wedge \mathfrak{q} w \ .$$

(7.29)

To see about the right most term in (7.29), let $Q_- \subset Q_{(s_0,z_0)}$ denotes the set of points where $iF_{\hat{A}} \wedge w$ is a negative multiple of the volume form. The inequality in (7.19) has its $Q_{(s_0,z_0)}$ analog; this being the lower bound

$$\int_{Q_-} iF_{\hat{A}} \wedge w \geq -c_c \rho_r^{-1} \Delta_* \ .$$

(7.30)

The proof of (7.30) is identical to that of (7.19) with it understood that the $Q_{(s_0,z_0)}$ version of the upper bound in (7.18) replaces the integration domain with $Q_{(s_0,z_0)}$ and replaces the term $c_c m \rho_r$ on the far right hand side of (7.18) with $c_c \rho_r^{-1} \Delta_*$. Granted (7.30), write $Q_{(s_0,z_0)}$ as $(Q_{(s_0,z_0)} - Q_-) \cup Q_-$ and then use Taylor's theorem with remainder to see that

$$\int_{Q_{(s_0,z_0)}} iF_{\hat A} \wedge q w \ \geq -\,(\sup_{\{z:\ |z-z_0|<4\rho_r\}} |\tfrac{\partial}{\partial z} q_D|)\,\rho_r \int_{Q_{(s_0,z_0)}} iF_{\hat A} \wedge w \ - c_c\,(\sup_{\{z:\ |z-z_0|<4\rho_r\}} |\tfrac{\partial}{\partial z} q_D|)\, r^{-1/q}\,\Delta_* . \tag{7.31}$$

Introduce $\varsigma(z_0)$ to denote

$$q_D(z_0) - (\sup_{\{z:\ |z-z_0|<4\rho_r\}} |\tfrac{\partial}{\partial z} q_D|\,)\,(\rho_r + c_c r^{-1/q}). \tag{7.32}$$

If $\varsigma(z_0)$ is positive, then Lemmas 7.2 and 7.3 with (7.29) and (7.31) find

$$\int_{Q_{(s_0,z_0)}} iF_{\hat A} \wedge q_D w \ \geq \varsigma(z_0)\, c_c^{-1} \rho_r^{\,2}\, \Delta_* \,, \tag{7.33}$$

which is positive. A look at (7.2) finds $\varsigma(z_0)$ to be negative only in the case that

$$2 - D|z_0| \leq c_0\, \rho_r^{\,1/2} \,, \tag{7.34}$$

in which case $q_D$ is less than $c_0\, e^{-(\ln r)^2/c_0}$ because $\rho_r = (\ln r)^4$.

*Part 6*: Define the function $\mathfrak{f}$ on $\mathbb{R}$ by the rule

$$s \to \mathfrak{f}(s) = \int_{[s-2\rho_r,\, s+2\rho_r]\times Y} \tfrac{i}{2\pi} F_{\hat A} \wedge q_D w \,. \tag{7.35}$$

It follows from what is said in (7.28) that $\mathfrak{f}(s) \geq -c_0\, e^{-(\ln r)^2/c_0}$ . Note here and future reference, that the function on [0,1] given by the rule $x \to e^{-(\ln x)^2/c_0}$ is bounded from above by $c_k x^{-k}$ for any $k > 0$ with $c_k$ being an purely k-dependent constant.

With the lower bound on $\mathfrak{f}$ in mind, define the subset $\mathcal{W} \subset \mathbb{R}$ by the following rule: A point $s \in \mathcal{W}$ if and only if $\mathfrak{f}(s)$ is negative. The set $\mathcal{W}$ is open. More to the point, the fact that $\mathrm{X}_\gamma(A_+) - \mathrm{X}_\gamma(A_-)$ is an integer implies that $\mathrm{X}_\gamma(A_+) - \mathrm{X}_\gamma(A_-)$ is negative only if the measure of $\mathcal{W}$ is greater than $c_0^{-1}\, e^{(\ln r)^2/c_0}$ . The paragraphs that follow explain why $\mathcal{W}$ is nowhere near this large.

Let $I \subset \mathbb{R}$ denote a closed interval of length 1 where the total change in the function $s \to \mathfrak{a}(\mathfrak{d}|_s)$ across the interval is less than $r^{-1}$. Invoke Lemma 4.2 to see that

$$\int_{I\times Y} (|\tfrac{\partial}{\partial s} A|^2 \ + |\mathfrak{B}_{(A,\psi)}|^2 \ + 2r(|\tfrac{\partial}{\partial s}\psi|^2 \ + |D_A\psi|^2)) \leq 2r^{-1} \,. \tag{7.36}$$

This fact is used momentarily.

If $\mathcal{W}$ has total length greater than $r^4$, then it can be covered by no less than $c_0^{-1} r^4$ closed intervals of length 1 with center at a point in $\mathcal{W}$ and such that no more than $c_0^{-1}$ of these intervals contain any given point. Given that the total change along $\mathbb{R}$ of the function $s \to \mathfrak{a}(\mathfrak{d}|_s)$ is bounded by $c\,r\ln r$, there are at least $c_0^{-1}r^4$ intervals in the cover where (7.36) holds. Let I denote one of the latter and let $s_0 \in \mathcal{W}$ denote I's center point.

By assumption, $\mathfrak{f}(s)$ is negative, and so there exists $z_0 \in D_0$ with $\sup_{t\in\mathbb{R}/(\ell\mathbb{Z})} \Delta_{(s_0,t,z_0)}$ greater than zero. It then follows from Lemmas 7.2 and 7.3 that the integral of $iF_{\hat{A}} \wedge w$ over $Q_{(s_0,z_0)}$ is greater than $c_c^{-1}\rho_r^2 \sup_{t\in\mathbb{R}/(\ell\mathbb{Z})} \Delta_{(s_0,t,z_0)}$. Given the formula in (6.25) and the bounds in Lemma 4.8 for $|\nabla_A\beta|$, this can occur only if the integral of $|\frac{\partial}{\partial s}A| + |\frac{\partial}{\partial s}\psi|^2$ over $Q_{(s_0,z_0)}$ is likewise greater than $c_c^{-1}\rho_r^2 \sup_{t\in\mathbb{R}/(\ell\mathbb{Z})} \Delta_{(s_0,t,z_0)}$. Since the latter is in any event greater than $c_c^{-1}\rho_r^4$, it follows that the integral of $|\frac{\partial}{\partial s}A|^2 + |\frac{\partial}{\partial s}\psi|^2$ over $Q_{(s_0,z_0)}$ is greater than $c_c^{-1}\rho_r^{12}$. This runs afoul of what is asserted by (7.36).

**c) The proof of Proposition 1.4**

The subsequent argument for Proposition 1.4 differs little from those used in [T2] to prove corresponding assertions that concern the analogs of (1.14) and (1.20) in the case when $\hat{a}$ is replaced by a contact 1-form and $w$ is replaced by the latter's exterior derivative.

Use $D_1$ to denote the data $(r_1, \mu_1, \mathfrak{p}_1)$ and use $D_2$ to denote the corresponding set $(r_2, \mu_2, \mathfrak{p}_2)$. Fix a smooth map, $s \to r_{(s)}$, from $\mathbb{R}$ to $[r_1, r_2]$ which is equal to $r_1$ for $s << -1$ and equal to $r_2$ for $s >> 1$. Fix a smooth map, $s \to \mu_{(s)}$, that is equal to $\mu_1$ for $s << -1$ and equal to $\mu_2$ for $s >> 2$. Finally, fix a smooth map from $\mathbb{R}$ to $\mathcal{P}$ of the form $s \to \mathfrak{p}_{(s)}$ such that $\mathfrak{p}_{(s)} = \mathfrak{p}_1$ for $s << -1$ and $\mathfrak{p}_{(s)} = \mathfrak{p}_2$ for $s >> 1$. The data $(r_{(s)}, \mu_{(s)}\ \mathfrak{p}_{(s)})$ can be used to define a version of (1.20), this being a system of equations that requires a map from $\mathbb{R}$ to $\mathrm{Conn}(E) \times C^\infty(Y;\mathbb{S})$ to obey

- $\frac{\partial}{\partial s}A + B_A - r_{(\cdot)}(\psi^\dagger\tau\psi - i\hat{a}) + \frac{1}{2} B_{A_K} - \mathfrak{T}^{(\cdot)}{}_{(A,\psi)} = 0$ .
- $\frac{\partial}{\partial s}\psi + D_A\psi - \mathfrak{S}^{(\cdot)}{}_{(A,\psi)} = 0,$

(7.37)

where $r_{(\cdot)}$ is the function $s \to r_{(s)}$ and $\mathfrak{T}^{(\cdot)}$ and $\mathfrak{S}^{(\cdot)}$ at any given $s \in \mathbb{R}$ are the gradients of $\mathfrak{e}_{\mu_{(s)}} + \mathfrak{p}_{(s)}$ along the respective $\mathrm{Conn}(E)$ and $C^\infty(Y;\mathbb{S})$ factors of $\mathrm{Conn}(E) \times C^\infty(Y;\mathbb{S})$. Of interest are the solutions to (8.37) with $s \to -\infty$ and $s \to \infty$ limits, these being solutions to the respective $(r_1,\mu_1)$ and $(r_2,\mu_2)$ versions of (1.13). Such a solution is called an instanton.

As explained in Chapter 25 of [KM], if the map $\mathfrak{p}_{(\cdot)}$ is chosen from a suitable residual set, then there will be but a finite number of instantons of the form $s \to \mathfrak{d}|_s$ with $\lim_{s\to\infty} \mathfrak{f}_{\mathfrak{s}}(\mathfrak{d}|_s) - \lim_{s\to-\infty} \mathfrak{f}_{\mathfrak{s}}(\mathfrak{d}|_s) = 0$. Chapter 25 of [KM] explains how to use the latter set to define a homomorphism between the $(r_1,\mu_1)$ version of $\hat{\mathcal{Z}}_{SW,r}$ to the $(r_2,\mu_2)$ version $\hat{\mathcal{Z}}_{SW,r}$ that preserves the relative $\mathbb{Z}/p_M\mathbb{Z}$ gradings and intertwines the endomorphisms that define the respective $D_1$ and $D_2$ differentials. This chapter also explains why the homomorphism intertwines the endomorphisms that define the actions of $\mathbb{Z}[\mathbb{U}]\otimes(\wedge^*(H_1(Y;\mathbb{Z})/\text{torsion}))$ on the respective $D_1$ and $D_2$ homologies.

The relevant homomorphism is defined by its action on the set of generators of the $(r_1,\mu_1)$ version of $\hat{\mathcal{Z}}_{SW,r}$. As such, it has the form

$$[\mathfrak{c}] \to \Sigma_{[\mathfrak{c}']}\ W_{[\mathfrak{c}'],[\mathfrak{c}]}\ [\mathfrak{c}'] \tag{7.38}$$

where the sum is indexed by the elements in the $(r_2,\mu_2)$ version of $\hat{\mathcal{Z}}_{SW,r}$ with any given coefficient $W_{[\mathfrak{c}'],[\mathfrak{c}]}$ being an integer. Only a finite set of these are non-zero. In particular, $W_{[\mathfrak{c}'],[\mathfrak{c}]}$ is non-zero only if $\mathfrak{f}_{\mathfrak{s}}(\mathfrak{c}) = \mathfrak{f}_{\mathfrak{s}}(\mathfrak{c}')$. In the latter case, Chapter 25 in [KM] defines $W_{[\mathfrak{c}'],[\mathfrak{c}]}$ to be a certain $\pm 1$ weighted count of the instanton solutions to (8.37) with $s \to -\infty$ limit equal to $\mathfrak{c}$ and with $s \to \infty$ limit equal to $\mathfrak{c}'$.

What is said in Chapter 25 in [KM] proves that the homomorphism in (7.38) induces an isomorphism between the corresponding $D_1$ and $D_2$ versions of $H^{\infty}_{SW,r}$. Chapter 25.5 in [KM] asserts that the respective isomorphisms that are defined by any two such $\mathbb{R}$-parameterized families are identical. The fact that any two such isomorphism are identical lead directly to the conclusion that there is a canonical isomorphism between the $D_1$ and $D_2$ versions of $H^{\infty}_{SW,r}$.

The proof of the assertions in Proposition 1.4 that concern $H^{-}_{SW,r}$, $H^{+}_{SW,r}$, and the long exact sequence that relates the latter with $H^{\infty}_{SW,r}$ has two parts.

*Part 1*: Let $r_* = \frac{1}{2}(r_1+r_2)$. Fix $\mu_* \in \Omega$ with $\mathcal{P}$-norm less than 1 and such that all instanton solutions to the $(r_*,\mu_*)$ version of (1.13) are non-degenerate. Fix in addition an element $\mathfrak{p}_* \in \mathcal{P}_{\mu_*}$ with small norm that obeys the $(r_*,\mu_*)$ version of (1.22). Use the data set $(r_*,\mu_*,\mathfrak{p}_*)$ to define the corresponding versions of $\mathbb{Z}(\hat{\mathcal{Z}}_{SW,r})$, $\mathbb{Z}(\hat{\mathcal{Z}}^{<}_{SW,r})$ and $\mathbb{Z}(\hat{\mathcal{Z}}_{SW,r})/\mathbb{Z}(\hat{\mathcal{Z}}^{<}_{SW,r})$, the operator $\partial^*_{SW}$ and then $H^{\infty}_{SW,r}$, $H^{-}_{SW,r}$ and $H^{+}_{SW,r}$. Let $\mathbb{L}_{*1}$ denote a homomorphism of the sort described above from the $(r_1,\mu_1)$ version of $\mathbb{Z}(\hat{\mathcal{Z}}_{SW,r})$ to the $(r_*,\mu_*)$ version, and let $\mathbb{L}_{2*}$ denote a homorphism of this sort from the $(r_*,\mu_*)$ version to the $(r_2,\mu_2)$ version.

Assume that $\mathbb{L}_{*1}$ maps the $(r_1,\mu_1)$ version of $\mathbb{Z}(\hat{\mathcal{Z}}^{<}_{SW,r})$ to the $(r_*,\mu_*)$ version and that $\mathbb{L}_{2*}$ maps the $(r_*,\mu_*)$ version to the $(r_2,\mu_2)$ version, then the composition $\mathbb{L}_{2*}\mathbb{L}_{*1}$ will map the $(r_1,\mu_1)$ version of $\mathbb{Z}(\hat{\mathcal{Z}}^{<}_{SW,r})$ to the $(r_2,\mu_2)$ version. The induced isomorphism from the $(r_1,\mu_1)$ version of $H^{\infty}{}_{SW,r}$ to the $(r_2,\mu_2)$ version of $H^{\infty}{}_{SW,r}$ will map the $(r_1,\mu_1)$ version of $H^{-}{}_{SW,r}$ to its $(r_2,\mu_2)$ counterpart, and likewise define a homomorphism between the respective version of $H^{+}{}_{SW,r}$. These homomorphisms will necessarily intertwine the homomorphisms in the two exact sequences.

Chapter 26 in [KM] explains why the composition $\mathbb{L}_{2*}\mathbb{L}_{*1}$ induces the canonical isomorphism from the $(r_1,\mu_1)$ version of $H^{\infty}{}_{SW,r}$ to the $(r_2,\mu_2)$ version. This understood, it follows as a consequence of what was said in the preceding paragraph that the canonical isomorphism from the $(r_1,\mu_1)$ version of $H^{\infty}{}_{SW,r}$ to the $(r_2,\mu_2)$ version induces respective homomorphisms from the $(r_1,\mu_1)$ versions of $H^{+}{}_{SW,r}$ and $H^{-}{}_{SW,r}$ to their $(r_2,\mu_2)$ counterparts that intertwine the associated long exact sequence homomorphisms.

*Part 2*: Given the conclusions of Part 1, the assertions in Proposition 1.4 that concern $H^{-}{}_{SW,r}$ and $H^{+}{}_{SW,r}$ and the long exact sequence relating $H^{\infty}{}_{SW,r}$, $H^{-}{}_{SW,r}$ and $H^{+}{}_{SW,r}$ follow as corollaries of the following lemma.

**Lemma 7.5**: *The versions of* $\kappa$ *in Propositions 1.1-1.3 can be chosen so that what is said below is also true. Fix* $R > \kappa$ *and there exists* $\kappa_R > \kappa$ *with the following significance: Fix* $r_1 \in (\kappa, R)$ *and an element* $\mu_1 \in \Omega$ *with* $\mathcal{P}$*-norm less than 1 such that all solutions to the* $(r_1,\mu_1)$ *version of* (1.13) *are non-degenerate. Fix* $r_2 > r_1$ *with* $|r_1 - r_2| < \kappa_R^{-1}$ *and fix* $\mu_2 \in \Omega$ *with* $\mathcal{P}$*-norm less than 1 such that all solutions to the* $(r_2,\mu_2)$ *version of* (1.13) *are also non-degenerate. Require in addition that* $\mu_2 - \mu_1$ *have* $\mathcal{P}$ *norm less than* $\kappa_R^{-1}$. *Fix respective elements* $\mathfrak{p}_1$ *and* $\mathfrak{p}_2$ *in* $\mathcal{P}$ *with small norm that obey the* $\mu = \mu_1$ *and* $\mu = \mu_2$ *versions of (1.22). There are homomorphisms of the sort described at the outset from the* $(r_1,\mu_1,\mathfrak{p}_1)$ *version of* $\mathbb{Z}(\hat{\mathcal{Z}}_{SW,r})$ *to the* $(r_2,\mu_2,\mathfrak{p}_2)$ *version that maps the* $(r_1,\mu_1,\mathfrak{p}_1)$ *version of* $\mathbb{Z}(\hat{\mathcal{Z}}^{<}_{SW,r})$ *to the* $(r_2,\mu_2,\mathfrak{p}_2)$ *version of* $\mathbb{Z}(\hat{\mathcal{Z}}^{<}_{SW,r})$.

***Proof of Lemma 7.5***: Suppose that no such $\kappa_R$ exists so as to derive nonsense. Given this contrarian assumption, there exist sequence $(r_{n1},\mu_{n1},\mathfrak{p}_{n1})$ and $(r_{n2},\mu_{n2},\mathfrak{p}_{n2})$ that obey the assumptions of the lemma but not the conclusions with $|r_{n1} - r_{n2}| < n^{-1}$, with $|\mu_{1n} - \mu_{n2}| < n^{-1}$ and with the $\mathcal{P}$-norms of $\mathfrak{p}_{n1}$ and $\mathfrak{p}_{n2}$ being less than $n^{-1}$. For each n, fix a corresponding $\mathbb{R}$ parametrized data set $(r_{n(\cdot)},\mu_{n(\cdot)},\mathfrak{p}_{n(\cdot)})$ that gives a version of (7.37) with instanton solutions that can be used to define the homomorphism between the respective $(r_{n1},\mu_{n1},\mathfrak{p}_{n1})$ and $(r_{n2},\mu_{n2},\mathfrak{p}_{n2})$ versions of $\mathbb{Z}(\hat{\mathcal{Z}}_{SW,r})$ in the manner described in the paragraph that surrounds (7.38). Such a path can and should be chosen so that $|\mu_{n(\cdot)} - \mu_{n1}| < n^{-1}$ and such that $\mathfrak{p}_{n(\cdot)}$ has $\mathcal{P}$-norm less than $n^{-1}$. The resulting index n homomorphism will map the $(r_{n1}, \mu_{n1}, \mathfrak{p}_{n1})$ version of $\mathbb{Z}(\hat{\mathcal{Z}}^{<}_{SW,r})$ to the $(r_{n2}, \mu_{n2}, \mathfrak{p}_{n2})$ version if the following is true: Let $\mathfrak{d}$ denote an instanton solution to the index n version of (7.37) with equal $s \to \infty$ and $s \to -\infty$ limits of $\mathfrak{f}_{\mathfrak{s}}(\mathfrak{d}|_s)$. Then the $s \to \infty$ limit of $\mathrm{x}(\mathfrak{d}|_s)$ is no less than the $s \to -\infty$ limit of $\mathrm{x}(\mathfrak{d}|_s)$.

Granted this point, there is for each n, an instanton solution to be denoted by $\mathfrak{d}_n$ with equal $s \to \infty$ and $s \to -\infty$ limits of $\mathfrak{f}_{\mathfrak{s}}(\mathfrak{d}_n)$ and with $\lim_{s\to\infty} \mathrm{x}(\mathfrak{d}|_s) \le \lim_{s\to-\infty} \mathrm{x}(\mathfrak{d}|_s) - 1$. An almost verbatim repetition of the final three paragraphs of Section 7a generates nonsense via Proposition 7.1 from the existence of the sequence $\{\mathfrak{d}_n\}_{n=1,2,\ldots}$.

**d) Proof of Theorem 1.5: Part 1**

This subsection gives a proof of the first and second bullets of Theorem 1.5 and explains how the third bullet follows from an auxilliary proposition that is proved in the next subsection.

To start, fix r large and define Theorem 1.5's map $\hat{\Phi}^r$ to be the L´ version of the map supplied by Proposition 3.1. What is said by Proposition 3.1 implies that the resulting version of $\mathbb{L}^r$ defines a $\mathbb{Z}$-module monomorphism from $\mathbb{Z}(\hat{\mathcal{Z}}_{ech,M}{}^{L'})$ to $\mathbb{Z}(\hat{\mathcal{Z}}_{SW,r})$ obeying the first and second bullets of Theorem 1.5.

The upcoming Proposition 7.6 is used to prove Theorem 1.5's third bullet. To give some background for Item c) of the third bullet in this proposition, suppose for the moment that $\mathfrak{c}$ and $\mathfrak{c}'$ are non-degenerate solutions to some $(r,\mu)$ version of (1.13), and suppose that $\mathfrak{d}$ is a non-degenerate instanton solution to the corresponding version of (4.1) with $s \to -\infty$ limit $\mathfrak{c}'$ and $s \to \infty$ limit $\mathfrak{c}$. The corresponding operator in (1.21) is Fredholm, and this being the case, suppose that its index is 1 and that it has trivial cokernel. These properties are all that is needed to compute the ±1 weight that $\mathfrak{d}$ would contribute to the coefficient $\mathrm{W}_{[\mathfrak{c}'],[\mathfrak{c}]}$ in (1.24) were the pair $(r, \mathfrak{g} = \mathfrak{e}_\mu)$ suitable for defining the $\partial^*_{SW}$ homology. This point is important for the following reason: If $(r, \mathfrak{g} = \mathfrak{e}_\mu + \mathfrak{p})$ with $\mathfrak{p}$ from $\mathcal{P}_\mu$ is suitable for defining the $\partial^*_{SW}$ homology, then the $(r,\mu)$ instanton $\mathfrak{d}$ contributes to the $(r, \mathfrak{g} = \mathfrak{e}_\mu + \mathfrak{p})$ version of $\mathrm{W}_{[\mathfrak{c}'],[\mathfrak{c}]}$, and its contribution to the $(r, \mathfrak{g} = \mathfrak{e}_\mu + \mathfrak{p})$ version of $\mathrm{W}_{[\mathfrak{c}'],[\mathfrak{c}]}$ is this same ±1.

Proposition 7.6 refers to notation that is used in Section 1 to describe the endomorphisms that generate the respective actions of $\mathbb{Z}[\mathbb{U}] \otimes \wedge^* H_1(Y;\mathbb{Z})$/torsion action on the $\partial_{ech}$ homology and on the $\partial^*_{SW}$ homology. By way of a reminder, from Part 4 of Section 1b, the endomorphism of $\mathbb{Z}(\hat{\mathcal{Z}}_{ech,M})$ that defines the action of $\mathbb{U}$ on the $\partial_{ech}$ homology is defined as in (1.9) by a set of integers, $\{N^{\mathbb{U}}{}_{\hat{\Theta}',\hat{\Theta}}\}_{\hat{\Theta}',\hat{\Theta}\in\hat{\mathcal{Z}}_{ech,M}}$, these being either 0 or 1. Part 7 of Section 1c describes an analogous set of integers, $\{W^{\mathbb{U}}{}_{[\mathfrak{c}'],[\mathfrak{c}]}\}_{[\mathfrak{c}'],[\mathfrak{c}]\in\hat{\mathcal{Z}}_{SW,r}}$,

that appear in the version of (1.24) for the endomorphism that gives the action of $\mathbb{U}$ on the $\partial^*_{SW}$ homology. Part 4 of Section 1b defines the set $\{\{[\gamma^{(z)}]\}_{z\in¥}, \{\hat{\imath}_p\}_{p\in\Lambda}\}$, this being a set of 1-cycles that give a basis for $H_1(Y;\mathbb{Z})$/torsion. Let $\hat{\imath}$ denote a cycle from this set. This same Part 4 of Section 1b defines the endomorphism that gives the action of $\hat{\imath}$ on the $\partial_{ech}$ homology. This endomorphism is a version of (1.9) whose coefficients are denoted by $\{N^{\hat{\imath}}_{\hat{\Theta}',\hat{\Theta}}\}_{\hat{\Theta}',\hat{\Theta}\in\hat{\mathcal{Z}}_{ech,M}}$. Meanwhile, Part 7 of Section 1c defines the endomorphism that gives the $\hat{\imath}$'s action on the $\partial^*_{SW}$ homology by a version of (1.24) with the coefficients denoted by $\{W^{\hat{\imath}}_{[\mathfrak{c}'],[\mathfrak{c}]}\}_{[\mathfrak{c}'],[\mathfrak{c}]\in\hat{\mathcal{Z}}_{SW,r}}$.

**Proposition 7.6**: *Given* $\Theta\in\mathcal{Z}_{ech,M}$, *there exist* $\kappa_\Theta\geq 1$ *with the following significance: Fix* $r\geq\kappa_\Theta$ *and an element* $\mu\in\Omega$ *with* $\mathcal{P}$ *norm less than 1. Suppose that* $\hat{\Theta}\in\hat{\mathcal{Z}}_{ech,M}$ *is a lift of* $\Theta$.

- *The element* $\hat{\Theta}$ *is in the image of* $\hat{\Phi}^r$. *Use* $\mathfrak{c}$ *to denote its image in* $\hat{\mathcal{Z}}_{SW,r}$.
- *If* $\mathfrak{c}'\in\hat{\mathcal{Z}}_{SW,r}$ *is such that* $\mathcal{M}_1(\mathfrak{c}',\mathfrak{c})\neq\emptyset$, *then* $\mathfrak{c}'$ *is in the image of* $\hat{\Phi}^r$.
- *If* $\hat{\Theta}'\in\hat{\mathcal{Z}}_{ech,M}$ *is such that* $\mathcal{M}_1(\hat{\Theta}',\hat{\Theta})\neq\emptyset$, *then* $\hat{\Theta}'$ *is in the domain of* $\hat{\Phi}^r$. *Granted that such is the case, use* $\mathfrak{c}'$ *to denote* $\hat{\Phi}^r(\hat{\Theta}')$.
    a) *The space* $\mathcal{M}_1(\mathfrak{c}',\mathfrak{c})$ *has only non-degenerate instantons.*
    b) *There exists a smooth,* $\mathbb{R}$*-equivariant, 1-1, onto map* $\Psi^r$: $\mathcal{M}_1(\hat{\Theta}',\hat{\Theta})\to\mathcal{M}_1(\mathfrak{c}',\mathfrak{c})$.
    c) *The* $\pm 1$ *weight that any given element in* $\mathcal{M}_1(\hat{\Theta}',\hat{\Theta})$ *contributes to the coefficient* $N_{\hat{\Theta}',\hat{\Theta}}$ *in* (1.9)*'s formula for* $\partial_{ech}\hat{\Theta}$ *is the same as the weight that its* $\Psi^r$*-image would contribute to the coefficient* $W_{[\mathfrak{c}'],[\mathfrak{c}]}$ *in (1.24)'s formula for* $\partial^*_{SW}\mathfrak{c}$.
    d) *Let* $\hat{\imath}$ *denote a cycle from the set* $\{\{[\gamma^{(z)}]\}_{z\in¥}, \{\hat{\imath}_p\}_{p\in\Lambda}\}$. *The contribution of any given element in* $\mathcal{M}_1(\hat{\Theta}',\hat{\Theta})$ *to the coefficient* $N^{\hat{\imath}}_{\hat{\Theta}',\hat{\Theta}}$ *is the same as that of its* $\Psi^r$ *image to the coefficient* $W^{\hat{\imath}}_{[\mathfrak{c}'],[\mathfrak{c}]}$.
- *If* $\mathfrak{c}'\in\hat{\mathcal{Z}}_{SW,r}$ *is such that* $\mathcal{M}_{2,p}(\mathfrak{c}',\mathfrak{c})\neq\emptyset$, *then* $\mathfrak{c}'=\hat{\Phi}^r(\hat{\Theta}')$ *with* $\hat{\Theta}'\in\hat{\mathcal{Z}}_{ech,M}$ *being the unique element that contributes to the* $\mathbb{U}$*-map coefficient* $N^{\mathbb{U}}_{\hat{\Theta}',\hat{\Theta}}$. *The corresponding space* $\mathcal{M}_{2,p}(\mathfrak{c}',\mathfrak{c})$ *contains a single instanton which is non-degenerate and the contribution of the latter to* $W^{\mathbb{U}}_{[\mathfrak{c}'],[\mathfrak{c}]}$ *is* 1, *this being* $N^{\mathbb{U}}_{\hat{\Theta}',\hat{\Theta}}$.

This proposition is proved in the next subsection. Accept it as true until then. What follow directly use Proposition 7.6 to prove the third bullet of Theorem 1.5.

Fix L´. The corresponding set $\mathcal{Z}_{ech,M}{}^{L'}$ contains but a finite number of elements. This understood, introduce $\kappa_{L'}$ to denote the largest of the $\Theta\in\mathcal{Z}_{ech,M}{}^{L'}$ versions of the constant $\kappa_\Theta$ supplied by Proposition 7.6. Fix $r>\kappa_{L'}$ and fix an element $\mu\in\Omega$ with $\mathcal{P}$-norm less than 1 and such that all solutions to the $(r,\mu)$ version of (1.13) are non-degenerate. Use these solutions to define the set $\mathcal{Z}_{SW,r}$ and to define $\mathbb{Z}$-module $\mathbb{Z}(\hat{\mathcal{Z}}_{SW,r})$. Fix an element $\mathfrak{p} \in \mathcal{P}_\mu$ with small $\mathcal{P}$-norm that can be used to define $\partial^*_{SW,r}$.

The assertion of the third bullet in Theorem 1.5 to the effect that $\mathbb{L}^r \partial_{ech,M} = \partial^*_{SW}\mathbb{L}^r$ follows directly from what is said by the second and third bullets of Proposition 7.6. Indeed, as elements in $\mathcal{P}_\mu$ vanish to second order on the images in Conn(E) $\times$ $C^\infty(Y;\mathbb{S})$ of non-degenerate instantons with $\iota_{(\cdot)} = 1$, the second bullet of Proposition 7.6 with Items a) and b) of the third bullet of Proposition 7.6 imply that any instanton that is used to compute the action of $\partial^*_{SW}$ on the image of $\Phi^r$ is in the image of some version of $\Psi^r$. Granted that such is the case, Item c) of the third bullet of Proposition 7.6 guarantees that the contribution of such an instanton to $\partial^*_{SW}$ is the same as the contribution of its $\Psi^r$-inverse image to $\partial_{ech}$.

The fact that $\mathbb{L}^r$ intertwines the endomorphisms that define the respective actions of $\mathbb{Z}[\mathbb{U}] \otimes H_1(Y;\mathbb{Z})$/torsion action on the $\partial_{ech}$ homology and on the $\partial^*_{SW}$ homology follows directly from what was the preceding paragraph with Item d) of the third bullet of Proposition 7.6 and the fourth bullet of Proposition 7.6.

## e) Proof of Proposition 7.6

The proof of the proposition has seven parts.

*Part 1*: The first bullet of the proposition follows from Proposition 3.1. It is also the case that if r is large and $\hat{\Theta}'$ is such that $\mathcal{M}_1(\hat{\Theta}', \hat{\Theta}) \neq \emptyset$, then $\hat{\Theta}'$ is also in the image of $\hat{\Phi}^r$ when r is large. The reason being that there are but a finite number of such $\hat{\Theta}'$ in $\hat{\mathcal{Z}}_{ech,M}$, this a fact that is explained in Section II.Ab.

The second bullet of Proposition 7.6 follows from Propositions 3.1 and 6.1 as does the assertion in the fourth bullet to the effect that if r is large, and if $\mathfrak{c}'$ is such that $\mathcal{M}_{2,\mathfrak{p}}(\mathfrak{c}', \mathfrak{c}) \neq \emptyset$, then $\mathfrak{c}'$ is in the image of $\Phi^r$. To explain how this comes about, introduce the function $\mathrm{M}$ on Conn(E) $\times$ $C^\infty(Y;\mathbb{S})$ from (1.30). Proposition 3.1 bounds $\mathrm{M}(\mathfrak{c}_+)$ by a multiple of $\sum_{\gamma\in\Theta} \ell_\gamma$. Suppose that $\mathfrak{c}'$ is a solution to the $(r,\mu)$ version of (1.3) and is such that $\mathcal{M}_1(\mathfrak{c}', \mathfrak{c}) \neq \emptyset$ or $\mathcal{M}_{2,\mathfrak{p}}(\mathfrak{c}', \mathfrak{c}) \neq \emptyset$. Let $\mathfrak{d}$ denote an instanton in either one of these spaces. Use Lemma 4.1 to see that $A_\mathfrak{d} \leq c_0 r(\mathrm{M}(\mathfrak{c}) + 1)$ and thus by $c_0 r(\sum_{\gamma\in\Theta} \ell_\gamma + 1)$. Granted this bound, use Proposition 6.1 to conclude that $\mathfrak{d}$'s version of the function $\underline{\mathrm{M}}$ is bounded by $c_\Theta$ with $c_\Theta > 1$ determined solely by the value of $\Theta$. It follows as a consequence that the $s \to -\infty$ limit of $\mathfrak{d}$'s version of $\underline{\mathrm{M}}$ is bounded by $c_\Theta$ and so $\mathrm{M}(\mathfrak{c}') < c_\Theta$. This being the case, Proposition 3.1 asserts that $\mathfrak{c}'$ is in the image of $\hat{\Phi}^r$ if r is larger than a purely $\Theta$ dependent constant.

*Part 2*: Keep in mind for what follows that the almost complex structure J is such that any $(\hat{\Theta}', \hat{\Theta})$ version of $\mathcal{M}_1(\hat{\Theta}', \hat{\Theta})$ has a finite set of $\mathbb{R}$ orbits and the Fredholm operator associated to each such orbit has trivial cokernel. The next remark is also important for what follows: Given $k \in \mathbb{Z}$, use $\hat{\Theta}_k$ for the moment to denote the translate of $\hat{\Theta}$ by the action of k on $\hat{\mathcal{Z}}_{ech,M}$ that comes about by viewing $\hat{\mathcal{Z}}_{ech,M}$ as a principal $\mathbb{Z}$ bundle over the set $\mathcal{Z}_{ech,M}$. Fix $\hat{\Theta}' \in \hat{\mathcal{Z}}_{ech,M}$. The respective sets of component subvarieties of $\mathcal{M}_1(\hat{\Theta}', \hat{\Theta})$ and of $\mathcal{M}_1(\hat{\Theta}'_k, \hat{\Theta}_k)$ are identical.

Granted these last facts, the construction that is described in Sections 4-7 of [T9] can be employed with only minor alterations if r is large to construct an $\mathbb{R}$-equivariant, injective map from any given $\hat{\Theta}' \in \hat{\mathcal{Z}}_{ech,M}$ version of $\mathcal{M}_1(\hat{\Theta}', \hat{\Theta})$ to the corresponding space $\mathfrak{M}_1(\mathfrak{c}', \mathfrak{c})$. Denote this map by $\Psi^r$. The arguments in Section 2b of [T5] can be used to construct $\Psi^r$ to have the following property: Let C denote any given component subvariety in $\mathcal{M}_1(\hat{\Theta}', \hat{\Theta})$. Write the instanton $\Psi^r(C)$ as $(A, \psi = (\alpha, \beta))$. If $\hat{\imath} \in \{\{[\gamma^{(z)}]\}_{z \in Ⲩ}, \{\hat{\imath}_p\}_{p \in \Lambda}\}$, then the intersections between C and $\mathbb{R} \times \hat{\imath}$ enjoy a 1-1 correspondence between those of $\alpha^{-1}(0)$ and $\mathbb{R} \times \hat{\imath}$ and this correspondence is such that partnered intersections have the same local intersection number. Note in this regard that J is such that C's intersections with $\mathbb{R} \times \hat{\imath}$ are finite in number and transversal. This is also the case for the intersections of $\alpha^{-1}(0)$ and $\mathbb{R} \times \hat{\imath}$. Note in addtion that the distance between any given point in $C \cap (\mathbb{R} \times \hat{\imath})$ and its corresponding partner in $\alpha^{-1}(0) \cap (\mathbb{R} \times \hat{\imath})$ is bounded by a $\Theta$ dependent multiple of $r^{-1/2}$.

What follows is a parenthetical remark with regards to the use here of the constructions in Sections 4-7 of [T9] and in Section 2b of [T5]. The constructions here use the simplest versions of those in the latter references by virtue of three facts, the first being that all integral curves of $v$ from all elements $\mathcal{Z}_{ech,M}$ are hyperbolic. The second is implied by the first: All subvarietes from any $\hat{\Theta}', \hat{\Theta} \in \hat{\mathcal{Z}}_{ech,M}$ version of $\mathcal{M}_1(\hat{\Theta}', \hat{\Theta})$ have the following property: Let C denote an element in $\mathcal{M}_1(\hat{\Theta}', \hat{\Theta})$. If $|s| >> 1$, then distinct components of the any constant $s$ slice of C are in small radius tubular neighborhoods of distinct integral curves of $v$, and each such component is isotopic in this neighborhood to its core integral curve. The final fact constitutes what is asserted by Lemma 3.2.

*Part 3*: The arguments in Section 3 of [T4] can be used with only very minor changes when r is large to prove the following: Fix $\hat{\Theta}' \in \hat{\mathcal{Z}}_{ech,M}$ with $\mathcal{M}_1(\hat{\Theta}', \hat{\Theta}) \neq \emptyset$. The map $\Psi^r$ restricts to each component of $\mathcal{M}_1(\hat{\Theta}', \hat{\Theta})$ as an $\mathbb{R}$ equivariant diffeomorphism onto a smooth component of $\mathfrak{M}_1(\mathfrak{c}', \mathfrak{c})$ with only non-degenerate instantons. Moreover, the contribution of any given component of $\mathcal{M}_1(\hat{\Theta}', \hat{\Theta})$ to the coefficient $N_{\hat{\Theta}',\hat{\Theta}}$ is the same as that of its $\Psi^r$ image to $W_{[\mathfrak{c}'],[\mathfrak{c}]}$. Note in this regard that the assumptions in Equation (1.14) of [T4] are not needed, this being a consequence of the three facts that are stated in the final paragraph of Part 2.

Given what was said in Part 1 about intersections with $\hat{\imath} \in \{\{[\gamma^{(z)}]\}_{z\in ¥}, \{\hat{\imath}_p\}_{p\in\Lambda}\}$ versions of $\mathbb{R} \times \hat{\imath}$, the conclusions of the preceding paragraph lead directly to the following: If $\hat{\imath} \in \{\{[\gamma^{(z)}]\}_{z\in ¥}, \{\hat{\imath}_p\}_{p\in\Lambda}\}$, then the contribution of any given component of $\mathcal{M}_1(\hat{\Theta}', \hat{\Theta})$ to $N^{\hat{\imath}}_{\hat{\Theta}',\hat{\Theta}}$ is identical to that if its $\Psi^r$ image to $W^{\hat{\imath}}_{[\mathfrak{c}'],[\mathfrak{c}]}$.

*Part 4*: Suppose that $\hat{\Theta}' \in \hat{\mathcal{Z}}_{ech,M}$ is such that $N^{\mathbb{U}}_{\hat{\Theta}',\hat{\Theta}} \neq \emptyset$. If r is large, then constructions in Sections 4-7 of [T9] with those in Sections 2d and 4 of [T5] construct a component of $\mathcal{M}_{2,p}(\mathfrak{c}', \mathfrak{c})$ that contains a single instanton which is non-degenerate and suitable for use in the definition of the coefficient $W^{\mathbb{U}}_{[\mathfrak{c}'],[\mathfrak{c}]}$ and contributes +1 to this integer. Note again that only the simplest versions of what is done in [T9] and [T6] are needed because only the J-holomorphic subvariety $(\times_{\gamma\in\Theta}(\mathbb{R}\times\gamma)) \cup (\{0\}\times S)$ is used, and this is the union of disjoint product cylinders and a compact submanifold. Note also that there is no need to introduce the notion of a (δ,L) approximation to use the constructions in [T6], this being yet another consequence of the three facts stated at the end of Part 2.

*Part 5*: It remains yet to prove that $\Psi^r$ maps any given version of $\mathcal{M}_1(\hat{\Theta}', \hat{\Theta})$ onto the corresponding version of $\mathcal{M}_1(\mathfrak{c}, \mathfrak{c}')$ and to prove that $\mathcal{M}_{2,p}(\mathfrak{c}', \mathfrak{c})$ has just the one component that is described in Part 4. The proofs that these assertions are true uses almost verbatim versions of arguments in Sections 4-7 of [T6] and in Section 4e of [T5]. Only the simplest cases of the arguments from [T6] and [T5] are needed, this also a consequence of the three facts stated in the final paragraph of Part 2. What follows directly and in Parts 6 and 7 say more about the analogs here of the relevant parts of [T6] and [T5].

Let $\mathfrak{c}'$ denote a solution to the $(r, \mu)$ version of (1.13) with either $\mathcal{M}_1(\mathfrak{c}', \mathfrak{c})$ or $\mathcal{M}_{2,p}(\mathfrak{c}', \mathfrak{c})$ nonempty. Let $\mathfrak{d}$ denote an instanton in one or the other of these spaces. The applications of the arguments from [T6] require as input the bound $\underline{M} < c_\Theta$ from Part 1 on $\mathfrak{d}$'s version of the function $\underline{M}$. Keep in mind that such a bound exists.

Suppose that there exists for each $n \in \{1, 2, \ldots\}$ a pair $r_n > n$ and an element $\mu_n$ in $\Omega$ with $\mathcal{P}$-norm less than 1 such that the $(r_n, \mu_n)$ version of the map $\Psi^r$ is not onto. If this is the case, there exists $\hat{\Theta}' \in \hat{\mathcal{Z}}_{ech,M}^{L'}$ and for each n, either an instanton solution to the $(r_n,\mu_n)$ version of (4.1) in the corresponding version $\mathcal{M}_1(\mathfrak{c}', \mathfrak{c})$ that is not in the image of the relevant version of $\Psi^r$, or an instanton in $\mathcal{M}_{2,p}(\mathfrak{c}', \mathfrak{c})$ that is not the one from Part 4. Use $\mathfrak{d}_n$ to denote this instanton. The latter is written when needed as $(A_n, \psi_n)$.

The rest of Part 5 and Parts 6 and 7 assume that the sequence $\{\mathfrak{d}_n\}_{n=1,2,\ldots}$ contains an infinite subset from the corresponding version of $\mathcal{M}_1(\mathfrak{c}', \mathfrak{c})$. Granted that this is the case, Lemmas 6.1 in [T6] has the following analog:

**Lemma 7.7**: *There is an element* $C \subset \mathcal{M}_1(\Theta_-, \Theta_+)$*, a subsequence of* $\{\mathfrak{d}_n\}_{n=1,2,\ldots}$ *(hence renumbered consecutively from 1) and a corresponding sequence of constant translations along the* $\mathbb{R}$ *factor of* $\mathbb{R} \times Y$*, all with the following property: For each* n*, write the translated version of* $\psi_n$ *as a pair* $(\alpha_n, \beta_n)$*. The sequence whose n'th element is*

$$sup_{z\in C}\, \mathrm{dist}(z, \alpha_n^{-1}(0)) + \sup_{z\in\alpha_n^{-1}(0)} \mathrm{dist}(C, z)$$

*converges with limit zero. In addition, if* $I \subset \mathbb{R}$ *is an interval of length 1 and* $\upsilon$ *is a 2-form on* $\mathbb{R} \times Y$ *with* $\|\upsilon\|_\infty = 1$ *and support on* $I \times Y$*, then the sequence whose n'th element is* $\frac{i}{2\pi}\int_{\mathbb{R}\times Y} \upsilon \wedge F_{\hat{A}_n} - \int_C \upsilon$ *also converges with limit zero.*

The proof of this lemma is given in Part 6; assume it to be true in the mean time. Lemma 7.7 leads to the analog of Lemma 6.2 in [T6], this has the identical assumptions and adds the following conclusion:

$$\lim_{n\to\infty} r_n^{1/2}\big(sup_{z\in C}\, \mathrm{dist}(z, \alpha_n^{-1}(0)) + \sup_{z\in\alpha_n^{-1}(0)} \mathrm{dist}(C, z)\big) = 0 \ . \tag{7.39}$$

The proof of (7.39) is almost identical to that of Lemma 6.2 in [T6] and so the reader is referred to Section 7 of [T6] for the proof of the latter's Lemma 6.2. By way of a guide to the proof of Lemma 6.2 of [T6], much of what is done in Section 7 of [T6] is of no concern to (7.39) because of the three facts listed at the end of Part 2. In particular, the integer m that enters in Lemmas 7.2–7.5 and Lemma 7.7 of [T6] can be set equal to 1. Moreover, most of the delicate estimates in Section 7d of [T6] are not needed because distinct $s >> 1$ slices, or distinct $s << -1$ slices of any given subvariety from $\mathcal{M}_1(\hat{\Theta}', \hat{\Theta})$ are in tubular neighborhoods of distinct integral curves of $v$ and are isotopic in these neighborhoods to the core integral curve.

Given Lemma 7.7 and (7.39), the argument to prove that $\Psi^r$ is onto is an almost verbatim copy of the arguments in Sections 6b-6e of [T6]. Only the simplest cases of these arguments are needed by virtue of the facts listed in the final paragraph of Part 1. In any event, the modifications of the arguments in Sections 6b-6e of [T6] are minimal and so nothing more will be said.

*Part 6*: The proof of Lemma 7.7 invokes an analog of Proposition 5.5 in [T6], this constituting the lemma that follows.

**Lemma 7.8**: *Given* $c \geq 1$, *there exists* $\kappa \geq 1$, *and given* $m > \kappa$, *there exists* $\kappa_m \geq 1$ *which, with* $\kappa$, *has the following significance: Fix* $r \geq \kappa_m$ *and* $\mu \in \Omega$ *with* $\mathcal{P}$*-norm less than 1. Suppose that* $\mathfrak{d} = (A, \psi = (\alpha, \beta))$ *is an instanton solution to (4.1) with* $A_{\mathfrak{d}} \leq c r$ *and with* $\lim_{s\to\infty} M(\mathfrak{d}|_s) < c$.

- *Each point in* $\mathbb{R} \times Y$ *where* $|\alpha| \leq 1 - m^{-1}$ *has distance at most* $\kappa r^{-1/2}$ *from where* $\alpha = 0$.
- *Moreover, there exists*
  a) *A positive integer* $N \leq \kappa$ *and a cover of* $\mathbb{R}$ *as* $\cup_{1\leq k\leq N} I_k$ *by connected open sets of length at least* $\frac{1}{2} m$. *These are such that* $I_k \cap I_{k'} = \emptyset$ if $|k - k'| > 1$. *In addition, if* $|k - k'| = 1$, *then* $I_k \cap I_{k'}$ *has length between* $\frac{1}{128} m$ *and* $\frac{1}{64} m$.
  b) *For each* $k \in \{1, 2, \ldots, N\}$, *a set* $\vartheta_k$ *whose typical element is a pair* $(C, m)$ *where* m *is a positive integer and where* $C \subset \mathbb{R} \times Y$ *is an irreducible, pseudoholomorphic subvariety. No two pair from* $\vartheta_k$ *contain the same subvariety component and* $\sum_{(C,m)\in\vartheta_k} m\int_C w < \kappa$.

*In addition, these sets* $\{\vartheta_k\}_{k=1,\ldots,N}$ *are such that*

1) $\sup_{z\in\cup_{(C,m)\in\vartheta_k} C \text{ and } s(z)\in I_k} \operatorname{dist}(z, \alpha^{-1}(0)) + \sup_{z\in\alpha^{-1}(0) \text{ and } s(z)\in I_k} \operatorname{dist}(\cup_{(C,m)\in\vartheta_k} C, z) < m^{-1}$.

2) *Let* $k \in \{1, \ldots, N\}$, *let* $I' \subset I_k$ *denote an interval of length 1, and let* $\upsilon$ *denote the restriction to* $I' \times Y$ *of a 2-form on* $\mathbb{R} \times Y$ *with* $\|\upsilon\|_\infty = 1$ *and* $\|\nabla\upsilon\|_\infty \leq m^{-1}$. *Then* $\left|\frac{i}{2\pi}\int_{I'\times Y} \upsilon \wedge F_{\hat{A}} - \sum_{(C,m)\in\vartheta} m\int_C \upsilon\right| \leq m^{-1}$.

This arguments for this lemma are given momentarily. Assume it to be true for the subsequent proof of Lemma 7.7.

***Proof of Lemma 7.7***: The proof has three steps.

Step 1: Pass to a subsequence of $\{(r_n,\mu_n), \mathfrak{d}_n\}_{n=1,2,\ldots}$ and renumber from 1 with the subsequence chosen so that Lemma 8.8 can be invoked with $m = n$ for each index n. Lemma 7.8 provides a corresponding sequence $\{\vartheta_{k,n}\}_{k=1,\ldots,N_n}$ with each $N_n$ a priori bounded by Lemma 7.8's constant $\kappa$. Since the sequence $\{N_n\}_{n=1,2,\ldots}$ is bounded, the sequence $\{(r_n,\mu_n), \mathfrak{d}_n\}_{n=1,2,\ldots}$ can be assumed to have the property that $N_n = N$ for all n.

A subsequence of $\{(r_n,\mu_n), \mathfrak{d}_n\}_{n=1,2,\ldots}$ can be chosen and renumbered from 1 so that the sequence $\{\{\vartheta_{k,n}\}_{k=1,2\ldots N}\}_{n=1,2,\ldots}$ converges to what is said to be a broken pseudoholomorphic subvariety. This is a collection $\{\vartheta_k\}_{k=1,2,\ldots,N'}$ of sets with the properties that are listed next. First, each $\vartheta_k$ is a finite set of pairs with each pair having the form $(C, m)$ with $C \in \mathbb{R} \times Y$ being an irreducible pseudoholomorphic subvariety and with m being a positive integer. Moreover, if $N > 1$, then each $\vartheta_k$ contains at least one pair whose subvariety is not $\mathbb{R}$-invariant.

The second property concerns the large $s$ limits of the elements in each $\vartheta_k$. In particular, the $s \to \infty$ limit of the constant $s$ slices of $\vartheta_k$ determine a finite set denoted by $\Theta_{k+}$ whose elements are pairs of the form $(\gamma, q)$ with $\gamma$ being a closed integral curve of $v$ and q being a positive integer. The manner by which $\vartheta_k$ determines $\Theta_{k,+}$ is as follows: For $s \in \mathbb{R}$, let $C|_s \subset Y$ denote the constant $s$ slice of C. View $\cup_{(C,m)\in\vartheta_k} m\, C|_s$ as a current. This current converges as $s \to \infty$ and the limit is the current $\cup_{(\gamma,q)\in\Theta_{k+}} q\gamma$. By the same token, the $s \to -\infty$ limit of $\cup_{(C,m)\in\vartheta_k} m\, C|_s$ determines a second set, $\Theta_{k,-}$, this having the same form as $\Theta_{k+}$. The collection $\{(\Theta_{k,-}, \Theta_{k,+})\}_{k=1,2,\ldots N'}$ are constrained by the requirement

$$\Theta_{1-} = \Theta' \;,\quad \Theta_{k,+} = \Theta_{k+1,-} \textit{ for } k = 1, \ldots, N'-1, \textit{ and } \Theta_{N,+} = \Theta \;. \tag{7.40}$$

Here and in what follows, $\Theta'$ denotes the image of $\hat{\Theta}'$ via the projection to $\mathcal{Z}_{ech,M}$. The convergence of $\{\{\vartheta_{k,n}\}_{k=1,2\ldots N}\}_{n=1,2\ldots}$ to $\{\vartheta_{k'}\}_{k'=1,2\ldots N'}$ is analogous to that described in the paragraph surrounding Equation (5.38) in [T6] with the only salient modification being the replacement of da in this equation in [T6] with $w$.

Step 2: The constraints on the first Chern class of E, what is said in (7.40) and what is said in Section II.3 about pseudoholomorphic subvarieties in $\mathbb{R} \times Y$ place extra constraints on the pairs from $\{\vartheta_k\}_{k=1,\ldots,N'}$. The first constraint involves the integers components of these pairs: If $(C,m) \in \cup_{k=1,\ldots,N'} \vartheta_k$, then $m = 1$ unless either C is compact or all components of its constant $s$ slices converge as $|s| \to \infty$ to closed integral curves of $v$ in $\cup_{p\in\Lambda} \mathcal{H}_p$. The remaining constraints involves the sets $\{\{\Theta_{k,-}, \Theta_{k+}\}\}_{k=1,\ldots,N'}$:

- *If $\gamma$ comes from a pair in $\cup_{k=1,\ldots,N'}(\Theta_{k,-} \cup \Theta_{k+})$, then $\gamma$ is disjoint from $\mathcal{H}_0$ and as a consequence, $\gamma$ lies entirely in the union of the $f \in (1,2)$ part of* $M_\delta$ *with $\cup_{p\in\Lambda}\mathcal{H}_p$.*
- *Each* $(\gamma, q) \in \cup_{k=1,\ldots,N'}(\Theta_{k,-} \cup \Theta_{k+})$ *has* $q = 1$ *unless* $\gamma \in \cup_{p\in\Lambda}\mathcal{H}_p$.
- *Fix* $k \in \{1, \ldots, N'\}$ *and let* $\Theta_{k,*}$ *denote either* $\Theta_{k,-}$ *or* $\Theta_{k,+}$. *Then* $(\cup_{(\gamma,1)\in\Theta_{k,*}} \gamma) \cap M_\delta$ *consists of* G *arcs that pair the index 1 and index 2 critical points of $f$ in* M *in the sense that distinct arcs start on the respective boundaries of the radius* $\delta$ *coordinate balls about distinct index 1 critical points of $f$ and end on the respective boundaries of the radius* $\delta$ *coordinate balls about distinct index* 2 *critical points of $f$.*

(7.41)

To prove that these constraints must be satisfied, start with $\vartheta_{N'}$ to see that the constraints on the integer components of its pairs are forced by the condition $\Theta_{N',+} = \Theta$. The constraints on $\Theta_{N,-}$ are then forced by the first Chern class considerations and the constraints on the integer components of the pairs in $\vartheta_N$. The constraints on $\Theta_{N',-}$ must be obeyed by $\Theta_{N'-1,+}$, and these constraints require the asserted constraints on the integer components of the pairs in $\vartheta_{N'-1}$. The latter with the first Chern class considerations force the constraints on $\Theta_{N'-1,-}$ and thus on $\Theta_{N'-2,+}$. Continuing in this vein proves that the constraints in the preceding paragraph must hold for all $k \in \{1, \ldots, N'\}$.

Step 3: Fix $k \in \{1, \ldots, N'\}$ and use $Z_k$ to denote the 2-cycle in Y given by the push-forward via the projection of $\sum_{(C,m)\in\vartheta_k} m[C]$ with [C] here denoting the non-compact cycle in Y that is carried by the fundamental class of C. The boundary of $Z_k$ is the 1-cycle $\sum_{(\gamma,q)\in\Theta_{k,+}} q[\gamma] - \sum_{(\gamma,q)\in\Theta_{k,-}} q[\gamma]$.

Definition 2.14 in [Hu2] uses $[Z_k]$ to define the embedded contact homology index, this being an integer that is denoted here by $I(\Theta_{k-}, \Theta_{k+}, Z_k)$. Let Z denote $\sum_{1\le k\le N'} Z_k$. Given Remark 2.16 in [Hu2], what is said by (7.40) implies that

$$I(\Theta', \Theta, Z) = \sum_{1\le k\le N'} I(\Theta_{k,-}, \Theta_{k,+}, Z_k) \; . \tag{7.42}$$

The argument in Part 2 of the proof of Lemma 6.1 in [T6] for the non-torsion case can be copied here with only minor changes to see that $I(\Theta', \Theta, Z) = 1$. This argument uses what is said in Lemma 7.8 about the large n versions of $\alpha_n^{-1}(0)$ and the fact that the instanton $\mathfrak{d}_n$ is in the $(r_n, \mu_n)$ version of $\mathcal{M}_1(\mathfrak{c}', \mathfrak{c})$.

It follows from Hutching's Definition 2.14 in [Hu2], from the description in Propositions II.3.1-II.3.4 of the pseudoholomorphic subvarieties in $\mathbb{R} \times Y$, and from (8.41) that $I(\Theta', \Theta, Z) = 1$ if and only if $N' = 1$, in which case $\vartheta_1$ defines an element in $\mathcal{M}_1(\hat{\Theta}', \hat{\Theta})$. The argument for this uses the inherently positive intersection numbers between pseudoholomorphic subvarieties in much the same way as used in the proof of Lemma III.8.3 to more than off set negative contributions to the sum in (7.42) that come from any given pair $(C, m) \in \cup_{k=1,\ldots,N'} \vartheta_k$. This is illustrated in the proof of Lemma III.8.3 by the formula (III.8.6) .

Given that $N' = 1$ and that $\vartheta_1$ defines an element in $\mathcal{M}_1(\hat{\Theta}', \hat{\Theta})$, then what is said in Lemma 7.7 follows directly from the conclusions of Lemma 7.8.

*Part 7*: This part contains the

***Proof of Lemma 7.8***: The arguments are much like the simplest versions of those used to prove Proposition 5.5 in [T6]. The six steps that follow describe what is needed from [T6] and what parts of these arguments need more than purely cosmetic changes.

Step 1: Given the bound in Proposition 6.1 on $\underline{M}$, Proposition 4.5 in [T6] has a simpler analog here also. This analog is a slightly weaker version of Lemma 7.8 that differs from Lemma 7.8 only to the extent that it does not make the claim that the pseudoholomorphic subvarieties in any given $k \in \{1, \ldots, N\}$ version of $\vartheta_k$ are defined on the whole of $\mathbb{R} \times Y$. The weaker version claims instead that the $\vartheta_k$ subvarieties are defined on a neighborhood of $I_k \times Y$.

The argument that derives Lemma 7.8 from its weak version amounts to little more than a standard application of a local form of the Gromov compactness theorem for pseudoholomorphic subvarieties. This argument differs little from the compactness theorems in [BEHWZ]. Given the a priori bound on $\underline{M}$ from Proposition 6.1, the derivation of Lemma 7.8 from its weak analog differs only in notation from what is said in [T6] to deduce Proposition 5.5 in [T6] from Proposition 5.1 in [T6].

Step 2: What follows here and in the subsequent steps prove the weak version of Lemma 7.8 using a modified version of the argument for Proposition 4.5 in [T6]. The modified version of this proposition is stated by the next lemma.

**Lemma 7.9**: *Given $c \geq 1$, there exists $\kappa \geq 1$, and given $m > \kappa$, there exists $\kappa_m \geq 1$ which, with $\kappa$, has the following significance: Fix $r \geq \kappa_m$ and $\mu \in \Omega$ with $\mathcal{P}$-norm less than 1. Suppose that $\mathfrak{d} = (A, \psi = (\alpha, \beta))$ is an instanton solution to (4.1) with $A_{\mathfrak{d}} \leq c r$ and with $\lim_{s\to\infty} M(\mathfrak{d}|_s) < c$. Let $I \subset \mathbb{R}$ denote an interval of length at least $m$.*

- *Each point in $I \times Y$ where $|\alpha| \leq 1 - m^{-1}$ has distance $\kappa r^{-1/2}$ or less from a zero of $\alpha$.*
- *There exists a finite set, $\vartheta$, whose components are pairs of the form $(C, m)$ where $C$ is a closed, irreducible pseudoholomorphic subvariety in a neighborhood of the closure of $I \times Y$ and where $m$ is a positive integer. Moreover, no two pairs in $\vartheta$ share the same subvariety component. This set is such that*

  a) $\sup_{z\in\cup_{(C,m)\in\vartheta} C \text{ and } s(z)\in I} \operatorname{dist}(z, \alpha^{-1}(0)) + \sup_{z\in\cup_{(C,m)\in\vartheta} C \text{ and } s(z)\in I} \operatorname{dist}(\cup_{(C,m)\in\vartheta} C, z) < m^{-1}$.

  b) *Let $\upsilon$ denote a smooth 2-form on $I \times Y$ with compact support, with $\|\upsilon\|_\infty = 1$ and with $\|\nabla\upsilon\|_\infty \leq m$. Then* $\left| \frac{i}{2\pi}\int_{I\times Y} \upsilon \wedge F_{\hat{A}} - \sum_{(C,m)\in\vartheta} m\int_C \upsilon \right| \leq m$.

  c) $\sum_{(C,m)\in\vartheta} m\int_C w \leq \kappa$.

The proof of this lemma is given momentarily. The next lemma plays a central role in the proof, and in subsequent arguments in this section.

**Lemma 7.10**: *Given $c \geq 1$, there exists $\kappa_c \geq 1$ with the following significance: Fix $r \geq \kappa_c$ and $\mu \in \Omega$ with $\mathcal{P}$-norm less than 1. Suppose that $\mathfrak{d} = (A, \psi = (\alpha, \beta))$ is an instanton solution to (4.1) with $A_{\mathfrak{d}} \leq c r$ and with $\lim_{s\to\infty} M(\mathfrak{d}|_s) < c$. Let $I \subset \mathbb{R}$ denote an open interval. Then* $-\kappa_c r^{-1/3} < \int_{I\times Y} iF_{\hat{A}} \wedge w < \kappa_c$.

***Proof of Lemma 7.9***: Given the bound in Proposition 6.1 on $\underline{\mathrm{M}}$, what is asserted by Lemma 7.9 follows directly from the $\mathrm{Y}_* = \mathrm{Y}$ version of Proposition 6.3 with the help of Lemma 7.10. The latter is needed to deduce Item c) of second bullet of Lemma 7.9 from the second bullet of Proposition 6.3.

***Proof of Lemma 7.10***: Invoke Proposition 6.1 to bound $\underline{\mathrm{M}}$ by $c_c$ with $c_c$ denoting here and in what follows a purely $c$ dependent constant which is greater than 1. Its value can be assumed to increase between successive appearances.

Consider first the upper bound. To this end, write I as $(s_1, s_2)$. Fix $\mathrm{s} \in (0, 1)$ and suppose that $m > 0$ is such that the integral of $iF_{\hat{A}} \wedge w$ over $[s_1 - \mathrm{s}, s_1 + \mathrm{s}] \times \mathrm{Y}$ is bounded from above by $m$. As explained directly, this implies that the integral of $iF_{\hat{A}} \wedge w$ is bounded by $m + c_c r^{-1}$. To see why this is, use the bound on $\underline{\mathrm{M}}$ to first invoke Lemma 4.9 and then use Lemma 4.9 to see that the function $\sigma$ in (5.7) is such that $1 - \sigma > -c_c r^{-1}$. Given this lower bound, use the top bullet in (7.3) with the bound $\underline{\mathrm{M}} < c_c$ to see that the respective integrals of $iF_{\hat{A}} \wedge w$ over $[s_1 - \mathrm{s}, s_1] \times \mathrm{Y}$ and over $[s_2, s_2 + \mathrm{s}] \times \mathrm{Y}$ are bounded from below by $-r^{-1} c_c$.

Since $\underline{\mathrm{M}}$ is bounded by $c_c$, there exists $\mathrm{s} \in (0, 1)$ such that $\mathrm{M}(\mathfrak{d}|_{s_1 - \mathrm{s}})$ and $\mathrm{M}(\mathfrak{d}|_{s_2 + \mathrm{s}})$ are both bounded by $c_c$. Given what is said in the preceding paragraph, it is sufficient to bound the integral of $iF_{\hat{A}} \wedge w$ over $[s_1 - \mathrm{s}, s_2 + \mathrm{s}] \times \mathrm{Y}$. This is done by comparing this integral to the integral of $iF_A \wedge w$ over $[s_1 - \mathrm{s}, s_2 + \mathrm{s}] \times \mathrm{Y}$. The comparison is made momentarily. What follows directly studies the integral of $iF_A \wedge w$ over $[s_1 - \mathrm{s}, s_2 + \mathrm{s}] \times \mathrm{Y}$.

Write $A = A_E + \hat{a}_A$ and integrate by parts to see that

$$\int_{[s_1 - \mathrm{s},\, s_2 - \mathrm{s}] \times \mathrm{Y}} iF_A \wedge w = i \int_{\{s_2 + \mathrm{s}\} \times \mathrm{Y}} \hat{a}_A \wedge w - i \int_{\{s_1 - \mathrm{s}\} \times \mathrm{Y}} \hat{a}_A \wedge w \,. \tag{7.43}$$

Meanwhile, use (1.26)-(1.28) and the fact that $\frac{i}{2\pi}(F_{A_E} + \frac{1}{2} F_{A_K})$ can be written as $w + d\mathfrak{b}$ to see that

$$r^{-1}(\mathfrak{a}(\mathfrak{d}|_{s_1 - \mathrm{s}}) - \mathfrak{a}(\mathfrak{d}|_{s_2 + \mathrm{s}})) = (1 - 2r^{-1})\Big(i \int_{\{s_2 + \mathrm{s}\} \times \mathrm{Y}} \hat{a}_A \wedge w - i \int_{\{s_1 - \mathrm{s}\} \times \mathrm{Y}} \hat{a}_A \wedge w\Big) + \mathfrak{e}\,, \tag{7.44}$$

where $\mathfrak{e}$ has absolute value bounded by $c_c r^{-1/3}$. To see why this is, note first that $\mathfrak{e}$ is bounded by a sum of two terms, the first being

$$c_0 r^{-1}\Big|\Big(\int_{\{s_2+s\}\times Y} \hat{a}_A \wedge *B_A\Big)\Big| + c_0(r^{-1/3}\mathrm{M}(\mathfrak{d}|_{s_2+s})^{4/3} + 1)$$

(7.45)

and the second having the same form but for the replacement of $s_2 + s$ with $s_1 - s$. By way of an explanation, the integral term comes from the part of $\mathfrak{cs}$ that involves $\hat{a}_A \wedge d\hat{a}_A$ and the term $c_0\, r^{-1/2}\, \mathrm{M}$ bounds the sum of the absolute values of three integrals, these being the respective integrals of $r^{-1}(i{*}B_A \wedge \mathfrak{b})$, of $r^{-1}(i{*}B_A \wedge \mu)$ and of $2\psi^\dagger D_A\psi$. The absolute values of these three integrals are bounded by $c_0 r^{-1/2}\mathrm{M}$; this thanks to Lemmas 4.8 and 4.9.

The integral of $\hat{a}_A \wedge *B_A$ over $[s_1 - s, s_2 + s] \times Y$ is the same as that with $\hat{a}_A$ replaced by $\hat{a}_A - u^{-1}du$ with u being any smooth map from Y to $S^1$. This being the case, no generality is lost by assuming that the $L^2$ orthogonal projection of $\hat{a}_A$ can be written as $\hat{a}_A^{\perp} + \mathfrak{p}$ where $\mathfrak{p}$ is a harmonic 1-form on Y with pointwise norm bounded by $c_0$ and where $\hat{a}_A^{\perp}$ is coclosed and orthogonal to the space of harmonic 1-forms on Y. Since $d\hat{a}_A^{\perp} = *B_A$, the argument used to prove Lemma 2.4 in [T1] can be used to prove that $|\hat{a}_A|^{\perp}$ is no greater than $c_0(r^{2/3}|\mathrm{M}|^{1/3} + 1)$. This understood, it follows from Lemma 4.9 that the integral of $r^{-1}(\hat{a}_A \wedge *B_A)$ that appears in (7.45) is bounded by $c_0(r^{-1/3}|\mathrm{M}|^{4/3} + 1)$.

To finish the story on (7.43), use Lemma 4.2 with the bound $\mathrm{A}_{\mathfrak{d}} < c\,r$ to bound the left hand side of (7.44) by $c$. Granted this bound and that on what is denoted by $\mathfrak{e}$ in (7.44), it then follows that the right hand side of (7.43) is no greater that $(1 - c_c r^{-1/3})\,\mathfrak{c}$.

Write $\hat{A}$ as $A_E + \hat{a}_{\hat{A}}$ and use integration by parts to write the $\hat{A}$ analog of the formula in (7.43). The latter has the integral of $iF_{\hat{A}} \wedge w$ over $[s_1 - s, s_2 + s] \times Y$ on the left hand side and has the same right hand side as the original version but for the replacement of $\hat{a}_A$ by $\hat{a}_{\hat{A}}$. This understood, use (1.15) with the bounds in Lemma 4.8 to see that the absolute value of the difference between right hand sides the respective $\hat{a}_A$ and $\hat{a}_{\hat{A}}$ versions of (7.43) is no greater than $c_c r^{-1/2}$.

To prove Lemma 7.10's lower bound assertion, use the bound $1 - \sigma > -c_c r^{-1}$ and the bound on $\underline{\mathrm{M}}$ by $c_c$ to see that the assertion holds for intervals of length less than 2. This being the case, take I to have length greater than 2. Write $I = (s_1, s_2)$ and suppose that $s \in (0, 1)$ and $m > 0$ are such that the integral of $iF_{\hat{A}} \wedge w$ over $[s_1 + s, s_2 - s] \times Y$ is greater than $-m$. The aforementioned bounds on $1 - \sigma$ and $\underline{\mathrm{M}}$ imply that the integral of $iF_{\hat{A}} \wedge w$ over both $[s_1, s_1 + s] \times Y$ and $[s_2 - s, s_2] \times Y$ is no less than $-c_c r^{-1}$. As a consequence, the integral of $iF_{\hat{A}} \wedge w$ over $I \times Y$ is no less than $-(m + c_c r^{-1})$. With the preceding understood, use the fact that $\underline{\mathrm{M}}$ is bounded by $c_c$ to choose s so that both $\mathrm{M}(\mathfrak{d}|_{s_1+s})$ and $\mathrm{M}(\mathfrak{d}|_{s_2-s})$ are bounded by $c_c$. The plan is to compare the integral of $iF_{\hat{A}} \wedge w$ over $[s_1 + s, s_2 - s] \times Y$ with that of $iF_A \wedge w$ over $[s_1 + s, s_2 - s] \times Y$. Use (7.43)-(7.45) with s replaced by -s to see that the latter integral is no less than $-c_c r^{-1/3}$. Meanwhile, the right hand side of this -s version of (7.43) differs from the right hand side of its Â counter part by at most $c_c r^{-1/2}$. The argument for this is identical but for the change s → -s as that given in the preceding paragraph.

Step 3: The $Y_* = Y$ versions of Lemmas 6.4 and 6.5 play the role here of that played by Lemma 4.6 in [T6] and Corollary 4.7 in [T6]. The next lemma is a replacement for Lemma 4.8 in [T6].

**Lemma 7.11**: *Given $m > 1$, there exists $\kappa_m > 1$, and given $\varepsilon > 0$, there exists $R_\varepsilon > 16$; and these have the following significance: Let $\mathbb{I} \subset \mathbb{R}$ denote an interval of length at least $2R_\varepsilon$, and suppose that C is a closed, irreducible, pseudoholomorphic subvariety in a neighorhood of $\mathbb{I} \times Y$ with $\int_{C\cap(I'\times Y)} w < \kappa_m^{-1}$ and $\int_{C\cap(I'\times Y)} ds \wedge \hat{a} \le m$ for all intervals $I' \subset \mathbb{I}$ of length 1. Assume in addition that C has intersection number zero with all submanifolds in $\mathbb{R} \times Y$ of the form $\{s\} \times S$ with S being a cross-sectional sphere in $\mathcal{H}_0$. Let $I \subset \mathbb{I}$ denote the subset with distance at least $R_\varepsilon$ from any boundary point of $\mathbb{I}$. There exists a finite set $\Theta$ consisting of pairs $(\gamma, q)$ with $\gamma$ a closed, integral curve of $v$ and q a positive integer. The set $\Theta$ is such that two pair share the same closed integral curve. Moreover,*

- $\sum_{(\gamma,q)\in\Theta} q\, \ell_\gamma \le m$.
- *Each point of* $C|_s$ for $s \in I$ *has distance along* Y *less than* $\varepsilon$ *from* $\cup_{(\gamma,q)\in\Theta} \gamma$. *Conversely, each point in* $\cup_{(\gamma,q)\in\Theta} \gamma$ *has distance no greater than* $\varepsilon$ *from* $C|_s$.
- *If* $\upsilon$ *is a smooth 2-form on* $I \times Y$ *with* $\|\upsilon\|_\infty = 1$ *and* $\|\nabla\upsilon\|_\infty \le \varepsilon^{-1}$. *Then*

$$\left| \int_{C\cap(I\times Y)} \upsilon - \sum_{(\gamma,q)\in\Theta} q \int_{I\times\gamma} \upsilon \right| < \varepsilon .$$

***Proof of Lemma 7.11***: The proof of Lemma 4.8 in [T6] can be copied with only the replacement of M with Y and with the references to Corollary 4.7 in [T6] replaced by references to Lemma 6.5.

The lower bound in Lemma 7.10 for integrals of $iF_{\hat{A}} \wedge w$ serves as a replacement for Lemma 4.9 in [T6].

Step 4: The remaining arguments for the weak version of Lemma 7.8 are similar in most respects, and simpler, than those given in Parts 4 and 5 of Section 5d in [T6] to prove Proposition 4.5 in [T6].

To complete the proof of the weak version of Lemma 7.8, fix $\varepsilon' > 0$ and define the subset $\mathcal{I} \subset \mathbb{Z}$ by the rule that places a given integer k in $\mathcal{I}$ if and only if the integral of $iF_{\hat{A}} \wedge w$ over $[k, k+1] \times Y$ is greater than $\varepsilon'$. It follows from the asserted upper bound from Lemma 7.10 in the case $I = \mathbb{R}$ and from the lower bound as applied to the components of $\mathbb{R} - (\cup_{k\in\mathcal{I}}[k,k+1])$ that $\mathcal{I}$ is a finite set with the number of components bounded by a constant that depends solely on $c$ and $\varepsilon'$. Use $n_{\varepsilon'}$ to denote this number.

Introduce the number $R_{\varepsilon'}$ from Lemma 7.11. There is a set, $\mathcal{V}$, of at most $n_{\varepsilon'}$ intervals in $\mathbb{R}$ and a pair of numbers, $c_{m*\varepsilon'}$ and $c_{m\varepsilon'}$, with the properties listed below.

- $c_{m*\varepsilon'}$ *and* $c_{m\varepsilon'}$ *are determined soley by* $n_{\varepsilon'}$ *and* $m$. *In any event,* $c_{m*\varepsilon'} > c_{m\varepsilon'} > 100\, n_{\varepsilon'}$.
- $\cup_{I\in\mathcal{V}} I$ *contains* $\cup_{k\in\mathcal{I}}[k,k+1]$.
- *Suppose that* $I \in \mathcal{V}$.
  a) $I$ *has length greater than* $c_{m\varepsilon'}(m + R_{\varepsilon'})$ *but less than* $c_{m*\varepsilon'}(m + R_{\varepsilon'})$.
  b) *If* $I \subset \mathcal{V}$ *then* $I$ *contains at least one* $k \in \mathcal{I}$ *version of* $[k, k+1]$.
  c) *If* $I \subset \mathcal{V}$*, then* $I \cap (\cup_{k\in\mathcal{I}}[k,k+1])$ *has distance at least* $10 R_{\varepsilon'}$ *from* $I$*'s boundary*.
- *If* $I$ *and* $I'$ *are distinct intervals from* $\mathcal{V}$ *with non-empty intersection, then* $I \cap I'$ *has length greater than* $\frac{1}{128} m$.
- *Each component of* $\mathbb{R} - (\cup_{I\in\mathcal{V}} I)$ *has length greater than* $4m$.

(7.46)

Given $I \subset \mathcal{V}$, use $I_* \subset I$ to denote the set of points with distance $\frac{1}{64} m$ or greater from any boundary point of I. The assertion of the weak version of Lemma 7.8 follow directly by using Step 2's analog of Proposition 4.1 in [T6] for its interval $\mathbb{I}$ taken in turn to be the intervals from $\mathcal{V}$ and using Lemma 7.11 for each component of $\mathbb{R} - (\cup_{I\in\mathcal{V}} I_*)$ with the constant in both replaced by $\varepsilon'$ and with the latter being a suitable function of $m$.

**f) Proof of Theorem 1.5: Part 2**

The five parts of this subsection complete the proof of Theorem 1.5. Part 1 proves the assertion at the very end of the theorem about the versions of $\mathbb{L}^r$ that are defined by distinct data sets. Part 2 of the subsection talks about some points in the proof given in Part 1 that are used implicitly in Parts 3-5. Parts 3 and 4 of the subsection prove the fifth bullet of Theorem 1.5; and Part 5 uses what is done in Parts 3 and 4 to prove the fourth bullet of Theorem 1.5.

*Part 1*: A proof is given momentarily for the final assertion of Theorem 1.5. What follows directly spells out what need proving. Fix $L' > 1$ and suppose that $(r, \mu, \mathfrak{p})$ and $(r', \mu', \mathfrak{p}')$ are data sets that satisfy the conditions demanded by Theorem 1.5. This is to say that the solutions to the respective $(r, \mu)$ and $(r', \mu')$ versions of (1.13) are nondegenerate and holonomy nondegenerate, and that the respective instanton solutions to the $(r, \mathfrak{g} = \mathfrak{e}_\mu + \mathfrak{p})$ and $(r', \mathfrak{g} = \mathfrak{e}_{\mu'} + \mathfrak{p}')$ versions of (1.20) are also nondegenerate. Granted these assumptions, the pair $(r, \mu)$ can be used to define $\mathbb{Z}(\hat{\mathcal{Z}}_{SW,r})$ and $\mathbb{Z}(\hat{\mathcal{Z}}^{<}_{SW,r})$, and $(r, \mu, \mathfrak{p})$ can be used to define the endomorphism $\partial^*_{SW}$. By the same token, $(r', \mu', \mathfrak{p}')$ can be used to define $\mathbb{Z}(\hat{\mathcal{Z}}_{SW,r'})$ and $\mathbb{Z}(\hat{\mathcal{Z}}^{<}_{SW,r'})$, and with $\mathfrak{p}'$ they define the corresponding version of $\partial^*_{SW}$. As noted in the Proposition 1.4, there is a canonical homomorphism between $H^{\infty}{}_{SW,r}$, $H^{-}{}_{SW,r}$ and $H^{+}{}_{SW,r}$ and the corresponding primed triad that intertwines the respective long exact sequences. Now suppose in addition that r and r′ are such that Proposition 3.1 can be used to define the $\hat{\Phi}^r$ and $\hat{\Phi}^{r'}$ on $\mathcal{Z}_{ech,M}{}^{L'}$. Theorem 1.5 asserts that this canonical homomorphism between homology groups can be lifted to a chain complex homomorphism between the respective $\mathbb{Z}(\hat{\mathcal{Z}}_{SW,r})$ and $\mathbb{Z}(\hat{\mathcal{Z}}_{SW,r'})$ which has the properties demanded by Proposition 1.4 and also intertwines the two versions of $\mathbb{L}^r$.

To prove this assertion of Theorem 1.5, return for the moment to Section 7c. The existence of a lift of Proposition 1.4's canonical homology homomorphism to a homomorphism from $\mathbb{Z}(\hat{\mathcal{Z}}_{SW,r})$ to $\mathbb{Z}(\hat{\mathcal{Z}}_{SW,r'})$ that satisfied Proposition 1.4's requirements is proved in Section 7c. In particular, it follows from what is said in Section 7c that such a lift can be found that factors as $\hat{\imath}_N \circ \hat{\imath}_{N-1} \circ \cdots \circ \hat{\imath}_1$ with $\hat{\imath}_1$ mapping the $(r, \mu, \mathfrak{p})$ version of $\mathbb{Z}(\hat{\mathcal{Z}}_{SW,r})$ to an $(r_1, \mu_1, \mathfrak{p}_1)$ version, with $\hat{\imath}_2$ mapping the latter version to an $(r_2, \mu_2, \mathfrak{p}_2)$ version, and so on, and with $\hat{\imath}_N$ mapping an $(r_{N-1}, \mu_{N-1}, \mathfrak{p}_{N-1})$ version to the $(r', \mu', \mathfrak{p}')$ version. These various data sets are such that each $k \in \{1, \ldots, N-1\}$ version of $|r_{k+1} - r_k|$ is very small as are the $\mathcal{P}$-norms of $\mu_{k+1} - \mu_k$ and $\mathfrak{p}_{k+1} - \mathfrak{p}_k$. Likewise, both $|r_1 - r|$ and $|r' - r_{N-1}|$ are small as are the $\mathcal{P}$-norms of $\mu_1 - \mu$, $\mathfrak{p}_1 - \mathfrak{p}$ and also $\mu' - \mu_{N-1}$ and $\mathfrak{p}' - \mathfrak{p}_{N-1}$. Note that 'small' here means as small as desired (but not zero) at the expense of increasing N.

Of particular import is that the sequence $\{r_k\}_{k=1,2...}$ can be chosen so that each element obeys the requirements set forth by Proposition 3.1 to define the corresponding version of $\hat{\Phi}^{(\cdot)}$ on $\hat{\mathcal{Z}}_{ech,M}{}^{L'}$. According to Proposition 3.1, each such version of $\hat{\Phi}^{(\cdot)}$ maps to nondegenerate and holonomy nondegenerate instantons. It is also the case that each element can be assumed to obey the conditions set forth in Proposition 7.6 for all pairs $(\Theta, \Theta') \in \mathcal{Z}_{ech,M}{}^{L'}$. Proposition 7.6 supplies for each such pair a corresponding map $\Psi^{(\cdot)}$ and of particular import is that the image of the latter consists of non-degenerate solutions to (4.1). Granted these last remarks, the final assertion of Theorem 1.5 follows from what is said in Parts 3-5 of Section 3h in [T2]. Part 2 of this subsection says more about these parts of [T2].

*Part 2*: The appeal to Parts 3-5 of Section 3h in [T2] uses only the fact that sufficiently large versions of $\hat{\Phi}^{(\cdot)}$ map to $\mathcal{G}_{M_\Lambda}$ orbits of nondegenerate solutions to the relevant version of (1.13). What follows elaborates on what nondegeneracy implies. Let $\mathfrak{c}$ denote a nondegenerate solution to a given $(r, \mu)$ version of (1.13). The nondegeneracy assumption is used in two related ways. The first uses $\mathfrak{c}$'s nondegeneracy with the implicit function theorem to build a smooth map from a neighborhood of r in $(\pi, \infty) \times \Omega$ into $\mathrm{Conn}(E) \times C^\infty(Y; \mathbb{S})$ with two salient properties: The map sends $(r, \mu)$ to $\mathfrak{c}$ and it maps any $(r', \mu')$ in its domain to a solution to the $(r', \mu')$ version of (1.13). In addition if $\mathfrak{c}'$ is a solution to the $(r', \mu')$ version of (1.13) and if the $C^\infty(Y; S^1)$ orbit of $\mathfrak{c}'$ is sufficiently close to $\mathfrak{c}$, then the image of $(r', \mu')$ via the map lies on this orbit. This map is denoted by $\hat{\mathfrak{c}}_{\mathfrak{c}}$ in what follows.

The second use of the nondegeneracy assumption concerns instanton solutions to (7.37). To say more, suppose that $(r, \mu)$ and $(r', \mu')$ are very close in $(-\infty, \pi) \in \Omega$, and that $\mathfrak{p} \in \mathcal{P}_\mu$ and $\mathfrak{p}' \in \mathcal{P}_{\mu'}$ are likewise very near each other in $\mathcal{P}$. Consider (7.37) when the data set $(r_{(\cdot)}, \mu_{(\cdot)}, \mathfrak{p}_{(\cdot)})$ has $s \to -\infty$ given by $(r, \mu, \mathfrak{p})$, $s \to \infty$ limit given by $(r', \mu', \mathfrak{p}')$; and when it is such that $(r_{(\cdot)}, \mu_{(\cdot)}, \mathfrak{p}_{(\cdot)})$ is nearly constant as $s$ varies in $\mathbb{R}$. Because $\mathfrak{c}$ is nondegenerate, standard pertubative techniques will prove the following: There is a unique instanton solution to (7.37) with Fredholm index equal to 0 whose $s \to -\infty$ limit is $\mathfrak{c}$. This instanton is very close to $\mathfrak{c}$ at each $s \in \mathbb{R}$ and its $s \to \infty$ limit is a solution to the $(r', \mu')$ version of (1.13) that is very close to $\mathfrak{c}$. In particular, this limit is the translate of $\hat{\mathfrak{c}}_{\mathfrak{c}}(r', \mu')$ by a map from Y to $S^1$ that is very close to the constant map to $1 \in S^1$. Let $\mathfrak{c}'$ denote this limit.

Looking ahead to Parts 3-5, the fact that the map from Y to $S^1$ is almost the constant map has the following implications: The values on $\mathfrak{c}'$ of the functions $\mathfrak{cs}$ and $\mathrm{W}$ in (1.26) and (1.27) are identical to their values on $\hat{\mathfrak{c}}_{\mathfrak{c}}(r', \mu')$. Likewise, the value on $\mathfrak{c}'$ of the $(r', \mu')$ version of (1.28)'a function $\mathfrak{a}$ is the same as its value on $\hat{\mathfrak{c}}_{\mathfrak{c}}(r', \mu')$. The $r'$ version of the function $\mathrm{M}$ automatically has the same values on $\mathfrak{c}'$ and $\hat{\mathfrak{c}}_{\mathfrak{c}}(r', \mu')$.

*Part 3*: The proof of the fifth bullet of Theorem 1.5 is given here and in Part 4. To start the proof, note that the residual set that is described by the second bullet in (1.18) can be chosen so as to have the properties listed below:

*There is a countable, non-accumulating bad set in* $(\pi, \infty)$ *such that if* r *avoids it, then*

- *The corresponding* $(r, \mathfrak{g} = \mathfrak{e}_\mu)$ *version of* $\mathcal{Z}_{\mathrm{SW},r}$ *is a finite set of* $C^\infty(Y; S^1)$ *orbits in* $\mathrm{Conn}(E) \times C^\infty(Y; \mathbb{S})$.
- *Each solution to (1.13) is non-degenerate and holonomy non-degenerate pairs.*
- *If* $\mathfrak{c}$ *and* $\mathfrak{c}'$ *are solutions to (1.13) in distinct* $C^\infty(Y; S^1)$ *orbits, then* $\mathfrak{a}^f(\mathfrak{c}) \neq \mathfrak{a}^f(\mathfrak{c}')$.

(7.47)

The arguments in Section 7b and 7c of [T1] can be used almost verbatim to prove this.

Fix an element $\mu \in \Omega$ with $\mathcal{P}$-norm less than 1 that is described (7.47). The latter describes a certain countable, non-accumulating subset of $(\pi, \infty)$. Denote this set by $\mathcal{U}$. If $r > \pi$ and is not $\mathcal{U}$, then the solutions to the $(r, \mu)$ version of (1.13) are suitable for defining the $\mathbb{Z}$-module $\mathbb{Z}(\hat{\mathcal{Z}}_{SW,r})$. Fix $r > \pi$ in the complement of $\mathcal{U}$ and choose a suitably generic element $\mathfrak{p} \in \mathcal{P}_\mu$ to define the differential $\partial^*_{SW}$ on $\mathbb{Z}(\hat{\mathcal{Z}}_{SW,r})$. Let $\iota$ denote a given class in either $\mathbb{H}^\infty$, $\mathbb{H}^-$ or $\mathbb{H}^+$. The class $\iota$ is then represented by a $\partial^*_{SW}$ cycle in $\mathbb{Z}(\hat{\mathcal{Z}}_{SW,r})$. This is to say that any given representative of $\iota$ can be written as

$$\mathfrak{z} = \sum_{[\mathfrak{c}]\in\hat{\mathcal{Z}}_{SW,r}} z_{[\mathfrak{c}]}\, [\mathfrak{c}] \,, \tag{7.48}$$

where each $z_{[\mathfrak{c}]} \in \mathbb{Z}$ and where only finitely many of these integers are non-zero. Associate to such a representative cocycle the number

$$\mathfrak{a}^f[\mathfrak{z},r] = \inf_{[\mathfrak{c}]\in\hat{\mathcal{Z}}_{SW,r} \text{ and } z_{\mathfrak{c}}\neq 0} \{\mathfrak{a}^f[\mathfrak{c}]\} \,; \tag{7.49}$$

and associate to the class $\iota$ the number

$$\mathfrak{a}^f_\iota[r] = \sup\{\mathfrak{a}^f[\mathfrak{z},r] \text{: } \mathfrak{z} \in \mathbb{Z}(\hat{\mathcal{Z}}_{SW,r}) \textit{ represents } \iota\} \,. \tag{7.50}$$

There are but a finite number of $C^\infty(Y; S^1)$ equivalence classes of solutions to (1.13) and so there is at least one cycle in $\mathbb{Z}(\hat{\mathcal{Z}}_{SW,r})$ that represents $\iota$ with $\mathfrak{a}^f[\cdot,r]$ equal to $\mathfrak{a}^f_\iota[r]$. Proposition 2.7 finds $\mathfrak{a}^f_\iota[r] < c_0 r$ and Proposition 2.7 with Lemma 2.5 find $\mathfrak{a}^f_\iota[r] > -c_0 r \ln r$.

Propositions 4.7 and 4.8 in [T2] have the following analog:

**Lemma 7.12**: *Choose $\mu \in \Omega$ with $\mathcal{P}$-norm less than 1 and described (7.47). Denote by $\mathcal{U} \subset (\pi,\infty)$ the bad set. Given $r \in (\pi, \infty) - \mathcal{U}$, use the solutions to the $(r,\mu)$ version of (1.13) to define $\mathbb{Z}(\hat{\mathcal{Z}}_{SW,r})$. There is a smooth map, $r \to \mathfrak{p}_r$, from $(\pi, \infty) - \mathcal{U}$ to $\mathcal{P}$ such that*

- *For each $r \in (\pi, \infty) - \mathcal{U}$, the element $\mathfrak{p}_r$ vanishes to second order on all solutions to the $(r,\mu)$ version of (1.13).*
- *The pair $(r, \mathfrak{g} = \mathfrak{e}_\mu + \mathfrak{p}_r)$ is suitable for defining $\partial^*_{SW}$ on $\mathbb{Z}(\hat{\mathcal{Z}}_{SW,r})$ if $r \in (\pi, \infty) - \mathcal{U}$ is chosen from a discrete set, $\mathcal{V}$, that accumulates only on the points in $\mathcal{U}$.*
- *Proposition 1.4's canonical isomorphism between the various $r \in (\pi, \infty) - (\mathcal{U} \cup \mathcal{V})$ versions of the $\partial^*_{SW}$ homology groups $H^\infty_{SW,r}$, $H^-_{SW,r}$ and $H^+_{SW,r}$ is such that following is true: If $\iota$ is any given non-zero homology class in $H^\infty_{SW}$, $H^-_{SW}$ or $H^+_{SW}$, then the assignment $r \to \mathfrak{a}^f_\iota[r]$ as defined above for $r \in (\pi, \infty) - (\mathcal{U} \cup \mathcal{V})$ is the restriction of a continuous and piecewise differentiable function on $(\pi, \infty)$.*

***Proof of Lemma 7.12***: But for notation and interchanging min with max, the proof mimics the arguments for Proposition 2.5 in [T7] and for Proposition 4.2 in [T1]. Note in this regard that the arguments in these papers use a homomorphism between the $(r, \mu, \mathfrak{p}_r)$ and $(r', \mu, \mathfrak{p}_{r'})$ versions of the $\partial^*_{SW}$ homology that is not obviously the canonical isomorphism. Even so, what is said in Part 2 of this subsection with arguments much like those in Section 7c can be used with the arguments [T7] and [T1] to obtain a proof of Lemma 7.12's assertion about the canonical isomorphism. See also Proposition 10.7 in [Hu3] and its proof.

Let $I$ denote a given class in either $H^{\infty}_{SW}$, $H^{-}_{SW}$ or $H^{+}_{SW}$. The function $\mathfrak{a}^f{}_I$ is important only to the extent that it can be used to analyze a second function of r, this denoted by $M_I(\cdot)$. To define the latter, fix $\mathfrak{p}_{(\cdot)}$ as in Lemma 7.12 so as to define $\partial^*_{SW}$ on $\mathbb{Z}(\hat{\mathcal{Z}}_{SW,r})$ when $r \in (\pi,\infty)-(\mathcal{U} \cup \mathcal{V})$. Given such a value for r, let $\mathfrak{z}$ denote a $\partial^*_{SW}$ cycle that represents $I$. Write $\mathfrak{z}$ as in (8.47) and define $M[\mathfrak{z},r] = \sup_{[\mathfrak{c}]\in\hat{\mathcal{Z}}_{SW,r} \text{ and } Z_{\mathfrak{c}}\neq 0} M(\mathfrak{c})$. Define

$$M_I[r] = \inf\{M[\mathfrak{z},r]: \mathfrak{z} \in \mathbb{Z}(\hat{\mathcal{Z}}_{SW,r}) \textit{ represents } I \textit{ and } \mathfrak{a}^f[\mathfrak{z},r] = \mathfrak{a}^f{}_I\} \ . \tag{7.51}$$

It follows from what is said in Part 2 of this section that $M_I$ is a priori a smooth function on $(\pi,\infty)-(\mathcal{U}\cup\mathcal{V})$.

*Part 3*. Theorem 1.5's fifth bullet follows from Proposition 3.1 if the following assertion is true:

$$\textit{Fix a class } I \textit{ in } H^{\infty}_{SW}, H^{-}_{SW} \textit{ or } H^{+}_{SW}. \textit{ Then the corresponding function } M_I(\cdot) \textit{ is bounded.} \tag{7.52}$$

Of course, (7.52) makes sense only when $\mu$ is described by (7.47); but Theorem 1.5 follows in any event using the fact that the set described in (7.47) is dense in $\Omega$ and what is said in Part 2 of this subsection.

To see about (7.52), fix an interval component of $(\pi,\infty)-(\mathcal{U} \cup \mathcal{V})$ and differentiate the expression in (1.28) on this interval to see that

$$\frac{d}{dr}\left(-\frac{2}{r}\,\mathfrak{a}^f{}_I\right) = \frac{1}{r^2}\,\mathfrak{cs}^f(\mathfrak{c}) \ , \tag{7.53}$$

with $\mathfrak{c}$ being a particular solution to (1.13) whose equivalence class has non-zero coefficient in some representative cycle for $I$ with $\mathfrak{a}^f[\cdot,r] = \mathfrak{a}^f{}_I$. Use Lemma 2.5 and Proposition 2.7 to see that the right hand side of (7.52) is no greater than $c_0 r^{-4/3}(\ln r)^{4/3}$. This being the case, integrate (7.52) on the components of $(\pi,\infty)-(\mathcal{U} \cup \mathcal{V})$ and use the fact that $\mathfrak{a}^f{}_I[\cdot]$ is continuous to see that $-\mathfrak{a}^f{}_I[r] \le c_I r + r^{1/3}(\ln r)^{4/3}$ with $c_I$ being a constant that depends on I but is independent of r. This last bound plus the bound implied by Lemma 2.5 and Proposition 2.7 for $|\mathfrak{cs}^f|$ requires $W^f[\mathfrak{c}] \le c_I + r^{-2/3}(\ln r)^{4/3}$ with $\mathfrak{c}$ being a particular solution to (1.13) whose equivalence class has non-zero coefficient in some representative cycle for I with $\mathfrak{a}^f[\cdot, r] = \mathfrak{a}^f{}_I$. Thus, $W^f[\mathfrak{c}]$ is bounded by an r-independent constant determined by the class I. This understood, (7.52) follows directly from the second bullet of Proposition 2.7.

*Part 4*: This part proves the fourth bullet of Theorem 1.5 in four steps.

<u>Step 1</u>: Fix $\mu$ from the set described by the second bullet of (1.18) and use $\mathcal{U}$ to denote the denote the corresponding countable, non-accumulating set $(\pi, \infty)$. Suppose that $r \in (\pi, \infty) - \mathcal{U}$ and is sufficiently large. In particular, require that r with $\mu$ and a suitably generic, small normed element $\mathfrak{p}$ from $\mathcal{P}$ can be used to define $\mathbb{L}^r$ on $\mathcal{Z}_{ech,M}{}^L$ so as to satisfy the first three bullets of Theorem 1.5 and the fifth bullet. Require in addition that the final assertion of Theorem 1.5 hold for $(r, \mu, \mathfrak{p})$ and data sets $(r', \mu, \mathfrak{p}')$ with $r' \ge r$.

<u>Step 2</u>: Introduce $\mathbb{Q}_{ech}{}^L$ to denote either $\mathbb{Z}(\hat{\mathcal{Z}}_{ech,M}{}^L)$, $\mathbb{Z}(\hat{\mathcal{Z}}_{ech,M}{}^{L,<})$ or $\mathbb{Z}(\hat{\mathcal{Z}}_{ech,M}{}^L)/\mathbb{Z}(\hat{\mathcal{Z}}_{ech,M}{}^{L,<})$ and let $\mathbb{Q}_{SW}$ denote the corresponding $\mathbb{Z}(\hat{\mathcal{Z}}_{SW,r})$, $\mathbb{Z}(\hat{\mathcal{Z}}^<_{SW,r})$ or $\mathbb{Z}(\hat{\mathcal{Z}}_{SW,r})/\mathbb{Z}(\hat{\mathcal{Z}}^<_{SW,r})$. Let $\varsigma$ denote an element in $\mathbb{Q}_{ech}{}^L$ such that $\mathbb{L}^r\varsigma = \partial^*{}_{SW}\mathfrak{z}$ with $\mathfrak{z}$ being an element in $\mathbb{Q}_{SW}$. If the fourth bullet of Theorem 1.5 holds for $\varsigma$, then it holds for $\varsigma + \partial_{ech}\varsigma'$ for any $\varsigma' \in \mathbb{Q}_{ech}{}^L$, this being a consequence of the second and third bullets of Theorem 1.5 and the fact that $\mathbb{L}^r$ is a monomorphism.

Since $(\partial^*{}_{SW})^2 = 0$, the second and third bullets of Theorem 1.5 require $\partial_{ech,M}\varsigma = 0$. Granted this, then $\varsigma$ defines a class in the homology of the chain complex $(\mathbb{Q}_{ech}{}^L, \partial_{ech})$. Use $I_\varsigma$ to denote this class. What is said in the preceding paragraph implies that the question of whether the fourth bullet of Theorem 1.5 holds for a given element $\varsigma$ depends only on the class $I_\varsigma$.

<u>Step 3</u>: The map $\mathbb{L}^r$ induces a homomorphism from the $(\mathbb{Q}_{ech}{}^L, \partial_{ech})$ homology to the $(\mathbb{Q}_{SW}, \partial^*{}_{SW})$ homology. The conclusions of Step 2 mean that the fourth bullet of Theorem 1.5 is asking about the kernel of this map. Let $\mathbb{K}_L$ denote this kernel. As explained next, the $\mathbb{Z}$-module $\mathbb{K}_L$ does not depend on the value of r if sufficiently large.

To start the explanation, suppose that $I_\varsigma \in \mathbb{K}_L$. This is to say that $\mathbb{L}^r\varsigma = \partial^*{}_{SW}\mathfrak{z}$ with $\mathfrak{z} \in \mathbb{Q}_{SW}$. As $\partial_{ech}\varsigma = 0$, the chain $\mathbb{L}^{r'}\varsigma$ is annihilated by the $(r', \mu, \mathfrak{p}')$ version of $\partial^*{}_{SW}$ when $r' \ge r$ is disjoint from $\mathcal{U}$. It follows as a consequence that $\mathbb{L}^{r'}\varsigma$ defines a class in either $H^\infty{}_{SW,r'}$, $H^-{}_{SW,r'}$ or $H^+{}_{SW,r'}$ as the case may be. It therefore defines a class in $H^\infty{}_{SW}$, $H^-{}_{SW}$ or $H^+_{SW}$ Given what is said in the last assertion of Theorem 1.5, this class in $H^\infty_{SW,r'}$, $H^-_{SW,r'}$ or $H^+_{SW,r'}$ corresponds to the class in $H^\infty_{SW,r}$, $H^-_{SW,r}$ or $H^+_{SW,r}$ that is defined by $\mathbb{L}^r\varsigma$.

The homology of $(\mathbb{Q}_{ech}{}^L, \partial_{ech})$ is finitely generated, and so $\mathbb{K}_L$ is finitely generated. This understood, the assertion made by the fourth bullet of Theorem 1.5 holds for all $\varsigma$ with $\mathbb{L}^r\varsigma \subset \text{Image}(\partial^*_{SW})$ if it holds for a judiciously chosen finite set of such elements. This understood, the fourth bullet of Theorem 1.5 holds if the following assertion is true:

$$\textit{Fix } \varsigma \in \mathbb{Q}_{ech}{}^L \textit{ with } \mathbb{L}^r(\varsigma) = \partial^*_{SW}\mathfrak{z} \textit{ for some } \mathfrak{z} \in \mathbb{Q}_{SW}.$$
$$\textit{There exists } L' > L \textit{ such that } \varsigma = \partial_{ech}\varsigma' \textit{ for some } \varsigma' \in \mathbb{Q}_{ech}{}^{L'}. \tag{7.54}$$

The proof of this assertion is given in Step 4.

<u>Step 4</u>: Arguments much like those that prove Lemma 7.12 prove the following:

**Lemma 7.13**: *Choose* $\mu \in \Omega$ *with* $\mathcal{P}$*-norm less than 1 and described (7.47). Denote by* $\mathcal{U} \subset (\pi, \infty)$ *the associated, non-accumulation set. Given* $r \in (\pi, \infty) - \mathcal{U}$*, use the solutions to the* $(r, \mu)$ *version of (1.13) to define the* $\mathbb{Z}(\hat{\mathcal{Z}}_{SW,r})$*. There is a smooth map,* $r \to \mathfrak{p}_r$*, from* $(\pi, \infty) - \mathcal{U}$ *to* $\mathcal{P}$ *such that*

- *For each* $r \in (\pi, \infty) - \mathcal{U}$*, the element* $\mathfrak{p}_r$ *vanishes to second order on all solutions to the* $(r, \mu)$ *version of (1.13).*
- *Data* $(r, \mathfrak{g} = \mathfrak{e}_\mu + \mathfrak{p}_r)$ *is suitable for defining* $\partial^*_{SW}$ *on* $\mathbb{Z}(\hat{\mathcal{Z}}_{SW,r})$ *if* $r \in (\pi, \infty) - \mathcal{U}$ *is chosen from the complement of a discrete set,* $\mathcal{V}$*, with accumulation only at points in* $\mathcal{U}$*.*
- *The assignment* $r \to \mathfrak{a}^f{}_\varsigma[r]$ *as defined above for* $r \in (\pi, \infty) - (\mathcal{U} \cup \mathcal{V})$ *is the restriction of a continuous and piecewise differentiable function on* $(\pi, \infty)$*.*

Fix $\mathfrak{p}_{(\cdot)}$ as in Lemma 7.13 so as to define $\partial^*_{SW}$ when $r \in (\pi, \infty) - (\mathcal{U} \cup \mathcal{V})$. Given such a value for r, let $\mathfrak{z}$ denote an element which is in either $\mathbb{Z}(\hat{\mathcal{Z}}_{SW,r})$ or $\mathbb{Z}(\hat{\mathcal{Z}}^<_{SW,r})$ with $\mathfrak{z}$ in the latter when $\mathbb{Q}_{ech}{}^L$ is $\mathbb{Z}(\hat{\mathcal{Z}}_{ech,M}{}^{L,<})$. Assume that $\mathfrak{z}$ obeys $\partial^*_{SW}\mathfrak{z} = \mathbb{L}^r\varsigma$ in $\mathbb{Q}_{SW}$. Write $\mathfrak{z}$ as in (7.48) and set $\mathrm{M}[\mathfrak{z}, r] = \sup_{[\mathfrak{c}] \in \hat{\mathcal{Z}}_{SW,r} \text{ and } Z_\mathfrak{c} \neq 0} \mathrm{M}(\mathfrak{c})$. Now define

$$\mathrm{M}_\varsigma[r] = \inf\{\mathrm{M}[\mathfrak{z}, r] : \mathfrak{z} \textit{ is such that } \partial^*_{SW}\mathfrak{z} = \mathbb{L}^r\varsigma \textit{ in } \mathbb{Q}_{SW} \textit{ and } \mathfrak{a}^f[\mathfrak{z}, r] = \mathfrak{a}^f{}_\varsigma\}\ . \tag{7.55}$$

An almost verbatim copy of the argument in Part 4 of this subsection proves that $\mathrm{M}_\varsigma[r]$ is bounded. This being the case, (7.54) follows from what is said in Proposition 3.1 about the map $\hat{\Phi}^r$.

## Appendix: The proof of Proposition 2.6

This appendix supplies a proof for Proposition 2.6. Much of what is done here mirrors similar constructions in Section 3 of [T9] and Section 2 of [T4]. Even so, the reworking can be justified for two reasons. First, the spectral flow function is not invariant under the action on Conn(E) $\times$ $C^\infty(Y;\mathbb{S})$ of the whole of the group $C^\infty(Y;\mathbb{S})$; and so care must be taken so as to not introduce a spurious gauge transformation in any given step of the proof. Care must also be taken so as not to introduce spurious factors of $\ln r$ in any given step. Such factors are easy to come by because there are solutions to (1.13) with (1.30)'s function $M$ being greater than $c_0^{-1}\ln r$. The need to avoid spurious gauge transformations and spurious factors of $\ln r$ accounts for the much of the length of the proof.

This appendix has three sections, these labeled as A, B and C.

## A. The eigenvalue equation $\mathfrak{L}_{c,r}\mathfrak{b} = \lambda\mathfrak{b}$

This section of the appendix supplies some necessary back ground for the proof of Proposition 2.6. Much of what is done here borrows heavily from Sections 3a-c of [T9] and Section 2a in [T4].

### a) Pairs in Conn(E) $\times$ $C^\infty(Y;\mathbb{S})$ and solutions to the vortex equation

This subsection uses solutions to (2.8) to construct pairs of connection on E and section of $\mathbb{S}$ over the complement in Y of tubular neighborhoods of a chosen subset of curves from $\cup_{p\in\Lambda}\{\hat\gamma_p^+ \cup \hat\gamma_p^-\}$. These constructions mimic those in Section 3 of [T9]. There are six parts to what follows.

*Part 1*: The input from the vortex equations includes first a pair $(A_0,\alpha_0)$ that obeys (2.8) and is such that $\frac{1}{2\pi}(1-|\alpha_0|^2)$ is integrable. As noted in (3.1), the integral of this function is a non-negative integer and of interest here is the case when the integer in question is 1. This is to say that $(A_0,\alpha_0)$ define a point in the vortex moduli space $\mathfrak{C}_1$. Use the pair $(A_0,\psi_0)$ to construct the square integrable solution $\varsigma$ to (3.27). Meanwhile, let y denote the square integrable, real valued function on $\mathbb{C}$ that solves the equation

$$-\partial\bar\partial y + \tfrac{1}{2}|\alpha_0|^2 y = -2^{-1/2}(1 - |\alpha_0|^2)\,. \tag{A.1}$$

The pair $(y, \varsigma)$ can be written explicitly in terms of $\alpha_0$ and its covariant derivative in the given case when $(A_0, \alpha_0)$ determine an element in $\mathfrak{C}_1$. For example, if $\alpha_0^{-1}(0) = 0$, then

$$\varsigma = -\bar z\bar\alpha_0^{-1}(1-|\alpha_0|^2) \quad \textit{and} \quad y = 2^{1/2} z\alpha_0^{-1}\partial_{A_0}\alpha_0\,. \tag{A.2}$$

In general, if $m \geq 1$ and if $(A_0, \alpha_0)$ defines a point in $\mathfrak{C}_m$, then there is a unique square integrable solution to (A.1) and a unique, square integrable solution to (3.27). This pair $(y, \varsigma)$ obeys $|y| + |\varsigma| \leq c_m\, e^{-\text{dist}(\cdot, \alpha_0^{-1}(0))/c_0}$ . Note for future reference that $\varsigma$ obeys Equation (2.35) of [T4].

Suppose that $(A_0, \alpha_0)$ defines an element in $\mathfrak{C}_1$ and is such that $\alpha_0^{-1}(0)$ is the origin in $\mathbb{C}$. Any two such solutions differ by the action of $C^\infty(\mathbb{C}; S^1)$ as they correspond to a single point in $\mathfrak{C}_1$; this is the point with $\sigma_1$ in (3.2) equal to zero. This point in $\mathfrak{C}_1$ is called the *symmetric vortex* as it is the unique fixed point in $\mathfrak{C}_1$ of the action by $S^1$ that is induced by the latters action on $\mathbb{C}$ as the group of rotations about the origin. There is a unique solution to (2.8) that maps to the symmetric vortex in $\mathfrak{C}_1$ and is such that

$$A_0 = \theta_0 - a_0\, \tfrac{1}{2}(z^{-1}dz - \bar{z}^{-1}d\bar{z}) \quad \textit{and} \quad \alpha_0 = |\alpha_0|\, \tfrac{z}{|z|}\ , \tag{A.3}$$

where the notation has $\theta_0$ denoting the product flat connection on the product line bundle over $\mathbb{C}$ and $a_0$ denoting a real valued function on $\mathbb{C}$. Lemma 3.3 can be used to prove that $a_0$ and $|\alpha_0|$ obey

$$|1 - a_0| \leq c_0(1 - |\alpha_0|^2) \quad \textit{and} \quad 1 - |\alpha_0| \leq c_0\, e^{-|z|/c_0}\ , \tag{A.4}$$

Note that if $m > 1$, then the point in $\mathfrak{C}_m$ with (3.2)'s coordinate functions all zero also corresponds to solutions $(A_0, \alpha_0)$ with $\alpha_0^{-1}(0) = 0$. There is in this case a unique solution with $\alpha_0^{-1}(0) = 0$ that has $\alpha_0 = |\alpha_0|\,(\tfrac{z}{|z|})^m$.

*Part 2*: This part of the subsection defines various terms and notions that are employed in the subsequent parts.

The constructions that are described below require the a priori choice of constants $c_v \geq 10^6$, $z \geq c_v^{\,6}$ and $\rho_* \geq c_v^{\,2} z^{-1/2}$. The lower bound for $c_v$ is increased to some $c_0 \geq 10^6$ in the applications to follow. The constants $c_v$ and $\rho_*$ are also constrained so that $c_v^{\,2}\rho_* \leq c_0^{-1}$ and in particular, $c_v^{\,2}\rho_*$ is smaller than $\frac{1}{100}$ times the maximum allowed radius of any transverse disk.

Proposition 3.1's map $\Phi^r$ uses the pairs constructed below with $z = r$, with $c_v$ on the order of 1, with $\rho_*$ constrained to be greater than $r^{-1/2+\delta}$ for some fixed $\delta > 0$. The proof of Proposition 2.6 uses versions of the constants with $z \in (c_0, r)$, with $c_v = c_0$, and with $\rho_*$ no larger than $c_0 z^{-1/2}$.

Introduce $Y_{*\Lambda}$ to denote the subset of Y with distance at least $c_v^{\,2}\rho_*$ from $\cup_{p\in\Lambda}(\hat{\gamma}_p^+ \cup \hat{\gamma}_p^-)$. Given that $c_v^{\,2}\rho_* \leq c_0^{-1}$, this subset $Y_{*\Lambda}$ is a smooth manifold with boundary whose boundary components are tori and whose complement is a disjoint union of solid tori tubular neighborhoods of the curves from the set $\{\hat{\gamma}_p^+ \cup \hat{\gamma}_p^-\}_{p\in\Lambda}$.

*Part 3*: A set of union of components of $Y - Y_{*\Lambda}$ must be specified in advance before starting the construction. This can be the empty set. The chosen subset of $Y - Y_{*\Lambda}$ is denoted in what follows by $T_{*\Lambda}$. The constructions that follow define a connection on E's restriction to $Y_{*\Lambda} \cup T_{*\Lambda}$ and a section of $\mathbb{S}$ over $Y_{*\Lambda} \cup T_{*\Lambda}$. Proposition 3.1's map $\Phi^r$ uses only the case when $T_{*\Lambda} = Y - Y_{*\Lambda}$. The proof of Proposition 2.6 uses all possible versions of $T_{*\Lambda}$.

The construction of a connection and section of $\mathbb{S}$ over $Y_{*\Lambda} \cup T_{*\Lambda}$ requires the choice of a finite set $\Theta$ whose elements are described below.

- *Each element in* $\Theta$ *is either a curve from* $\{\hat{\gamma}_p^+ \cup \hat{\gamma}_p^-\}_{p\in\Lambda}$ *that lies in* $T_{*\Lambda}$ *or a properly embedded, 1-dimensional, oriented submanifold in* $Y_{*\Lambda}$. *The curves from* $\Theta$ *are distinct.*
- *The* $T_{*\Lambda}$ *boundary components of* $Y - Y_{*\Lambda}$ *are disjoint from the* $Y_{*\Lambda}$ *curves from* $\Theta$. *Any other boundary component of* $Y_{*\Lambda}$ *contains no more then two endpoints of arcs from* $\Theta$, *and if two, then one has* $u < 0$ *while the other has* $u > 0$.
- *Suppose the* $\gamma$ *is a curve from* $\Theta$ *in* $Y_{*\Lambda}$.
  a) *The unit length, oriented tangent vector to* $\gamma$ *has distance at most* $c_v z^{-1/2}$ *from* $v$.
  b) *The curve* $\gamma$ *intersects any given* $\mathfrak{p} \in \Lambda$ *version of* $\mathcal{H}_p$ *where* $1 - 3\cos^2\theta > 0$.
  c) *If* $\gamma$ *is disjoint from a given boundary torus of* $Y_{*\Lambda}$, *then it has distance greater than* $3c_v\rho_*$ *from this torus.*
  d) *If* $\gamma$ *intersects a boundary torus of* $Y_{*\Lambda}$, *then it does so only at its endpoints. These intersections are transversal. Moreover, one end point of* $\gamma$ *lies where* $u < 0$ *on some boundary component of* $Y_{*\Lambda}$ *and the other where* $u > 0$ *on some boundary component of* $Y_{*\Lambda}$.
- *The intersection of* $\cup_{\gamma\in\Theta}\gamma$ *with* $M_\delta$ *sits in the* $f^{-1}(1,2)$ *part of* $M_\delta$. *This intersection consists of* G *properly embedded segments that pair the index 1 and index 2 critical points of f in the sense that distinct segments start on the boundary of the radius* $\delta$ *coordinate balls about distinct index 1 critical points of f and end on the boundary of the radius* $\delta$ *coordinate balls about distinct index 2 critical points.*

(A.5)

The proof of Proposition 2.6 uses only versions of $\Theta$ that lack curves from the set $\cup_{p\in\Lambda}\{\hat{\gamma}_p^+ \cup \hat{\gamma}_p^-\}$. Versions with curves from this set are needed to define Proposition 3.1's map $\Phi_r$.

Let $\gamma$ denote a 1-manifold in $Y_{*\Lambda}$ from an element in $\Theta$. Introduce $U_\gamma$ to denote the union of the radius $4\rho_*$ transverse disks centered at the points in $\gamma$. Use $U_\gamma' \subset U_\gamma$ to denote the union of the radius $\rho_*$ transverse disks centered at the points in $\gamma$. If $\gamma$ is a $\cup_{p\in\Lambda}\{\hat{\gamma}_p^+ \cup \hat{\gamma}_p^-\}$ curve from Θ, use $U_\gamma$ to denote the union of the radius $4\rho_{*\Lambda}$ transverse disks centered on γ and use $U_\gamma' \subset U_\gamma$ to denote the union of the concentric radius $\rho_{*\Lambda}$ transverse disks. Keep in mind the following consequence of the formula in Section 1a for *v*: If $c_v \ge c_0$, and if $\gamma \subset Y_{*\Lambda}$ is from Θ, then $U_\gamma$ is an open solid torus with γ the core circle. Moreover, if γ and γ′ are in $Y_{*\Lambda}$ and come from distinct elements in Θ, then $U_\gamma \cap U_{\gamma'} = \emptyset$. It is assumed in what follows that $c_v$ is such as to guarantee this.

*Part 4*: This part of the subsection describes certain set of preferred coordinates for the various versions of $U_\gamma$. Each element in this set is determined in a canonical fashion by an isometric isomorphism from $K^{-1}|_\gamma$ to $\gamma \times \mathbb{C}$.

To define these coordinates, introduce $T_\gamma$ to denote the union of the transverse disks of radius $c_0^{-1}$ centered at the points of γ. Choose this radius so that the union of the transverse disks with centers on any length less than $c_0^{-1}$ segment of γ is a solid, embedded cylinder with the segment as the core arc. The desired coordinates for $U_\gamma$ are obtained by restricting the domain of a set of functions on $T_\gamma$ that define local coordinates on each such solid cylinder.

Let $\ell_\gamma$ denote the length of γ. The first of these functions is a parameter, denoted by t, with values in $\mathbb{R}/(\ell_\gamma \mathbb{Z})$ when γ is a closed loop and with values in an interval of length $\ell_\gamma$ otherwise. The coordinate t is constant along each transverse disk in $T_\gamma$ with center on γ. The other coordinate is denoted by z; it is a $\mathbb{C}$-valued function that identifies each transverse disk with the radius $c_0^{-1}$ disk in $\mathbb{C}$ centered at the origin. The coordinate identification is such that the origin in $\mathbb{C}$ corresponds to γ.

THE DEFINITION: *Fix a $\mathbb{C}$-linear isomorphism between $K^{-1}|_\gamma$ and $\gamma \times \mathbb{C}$. This defines an orthonormal, oriented frame for the kernel of â along γ. Use this isomorphism with the metric's exponential map to identify a tubular neighborhood of γ with $\gamma \times \mathbb{C}$. Use* t *to denote an affine coordinate long γ with the property that the corresponding tangent vector field has unit length and positive pairing with â. The coordinate* t *and the standard complex coordinate* z *for $\mathbb{C}$ are the desired coordinate functions.*

These coordinates are such that the z = 0 locus is γ and $\partial_z$ along γ is in the kernel of *â* and has norm $2^{-1/2}$. The first order Taylor's expansion writes *v* and *w* as

- $v = \frac{\partial}{\partial t} + 2i(\nu z + \mu \bar{z} - x_\gamma)\frac{\partial}{\partial z} - 2i(\nu \bar{z} + \bar{\mu} z - \bar{x}_\gamma)\frac{\partial}{\partial \bar{z}} + \cdots$,
- $w = \frac{i}{2} dz\wedge d\bar{z} - (\nu z + \mu \bar{z} - x_\gamma)\, d\bar{z} \wedge dt - (\nu \bar{z} + \bar{\mu} z - \bar{x}_\gamma) dz\wedge dt + \cdots$,

(A.6)

where $\nu$ is a real valued function of t while $\mu$ and $x_\gamma$ are $\mathbb{C}$-valued functions of t with $x_\gamma$ such that $|x_\gamma| \le c_\nu z^{-1/2}$. The unwritten terms are bounded by $c_0(|\nu| + |\mu|)|z|(z^{-1/2} + |z|)$ with $c_0$ here dependent on $\nu$ and $\mu$. Note that $\nu$ must be real so as to have $dw = 0$. What is said in Item a) of the third bullet of (A.5) leads to asserted the bound on $|x_\gamma|$.

Changing the isomorphism between $K^{-1}|_\gamma$ and $\gamma \times \mathbb{C}$ writes z as $z = u(t) z'$ with u being a smooth map from the domain of t to $S^1$. The resulting version of (A.6) replaces $\nu$ with $\nu' = \nu + \frac{i}{2} u^{-1} \frac{d}{dt} u$, it replaces $\mu$ with $\mu' = u^{-2}\mu$ and it has $x_\gamma$ repaced by $x_\gamma' = u^{-1} x_\gamma$.

As is explained momentarily, this last observation has the following important consequence: Coordinates of the sort just described can be found with the property that the functions $\nu$ and $\mu$ in (A.6) obey $|\nu| + |\mu| \le c_0$. To see why this is, fix a point p in the interior of $\gamma$ and a unitary frame for $K^{-1}|_p$. Parallel transport this frame along a small length interval in $\gamma$ containing p. Use the latter frame with the exponential map to define the coordinates (t, z) for a solid cylinder with this interval as the core arc. Use $T_p$ to denote this cylinder. The Lie derivative of $w$ by $\frac{\partial}{\partial \bar z}$ is bounded by $c_0$ at p because the covariant derivative of $\frac{\partial}{\partial \bar z}$ is zero at p. Use the (t, z) version of (A.6) to see that this Lie derivative at p is $-\nu\, dz \wedge dt - \mu\, d\bar z \wedge dt$. This implies in particular that $|\mu| \le c_0$ at p. The fact that $|\mu|$ is independent of the chosen orthonormal frame for $K^{-1}|_\gamma$ implies that $|\mu| \le c_0$ along the whole of $\gamma$ no matter what frame is used to the coordinates. Meanwhile, the freedom to change $\nu$ to $\nu - \frac{i}{2} u^{-1} \frac{d}{dt} u$ can be exploited to obtain a version of the coordinates with the function $\nu$ such that $|\nu|$ is also bounded along $\gamma$. Indeed, if $\gamma$ is not closed, then this equation can be integrated so as to obtain a version with $\nu = 0$. If $\gamma$ is closed, then a version can be found with $\nu$ constant and less than $c_0\, \ell_\gamma^{-1}$.

Unless told otherwise, assume in here and in Appendix B and C that any chosen coordinate system of the sort described above has $|\nu| + |\mu| < c_0$. A second convention with regards to these coordinates concerns the case when $\gamma$ is an integral curve from the set $\cup_{p\in\Lambda}\{\hat\gamma_p^+ \cup \hat\gamma_p^-\}$. As explained in Part 5 of Section 3c there is a version of these coordinates with both $\nu$ and $\mu$ contant, with $\mu$ real and such that $\mu > |\nu|$. These constant values for $\nu$ and $\mu$ are denoted at times by $\nu_0$ and $\mu_0$. This $(\nu_0, \mu_0)$ version of the coordinates should be assumed in what follows when $\gamma \in \cup_{p\in\Lambda}\{\hat\gamma_p^+ \cup \hat\gamma_p^-\}$.

A coordinate system of the sort described above should be chosen for each set from the collection $\{U_\gamma\}_{\gamma\in\Theta}$. These chosen coordinate systems are used in what follows.

*Part 5*: Fix a set $\Theta$ as described in Part 3. The corresponding pair of connection and spinor is defined with the help of the open cover of $Y_{*\Lambda} \cup T_{*\Lambda}$ that consists of the collection $\{U_\gamma \cap Y_{*\Lambda}\}_{\gamma\in\Theta}$ and the set $U_0 = (Y_{*\Lambda} \cup T_{*\Lambda}) - (\cup_{\gamma\in\Theta}(U_\gamma' \cap Y_{*\Lambda})$. Use $\mathfrak{U}$ to denote the collection of sets consisting of $U_0$ and $\{U_\gamma\}_{\gamma\in\Theta}$. For each $U \in \mathfrak{U}$, an isometric isomorphism must first be chosen to identify $E|_U$ with $U \times \mathbb{C}$. The bundle E can be reconstructed from these isomorphisms using the corresponding transition functions.

The pair $(A,\psi)$ on any given set U from the cover is written using the isomorphism between $E|_U$ and $U \times \mathbb{C}$ as $A = \theta_0 + a_U$ where $\theta_0$ denotes the product connection on $U \times \mathbb{C}$ and where $a_U$ is an $i\mathbb{R}$ valued 1-form on U. Meanwhile, $\psi$ is written as $(\alpha_U, \beta_U)$ where $\alpha_U$ and $\beta_U$ denote a respective $\mathbb{C}$-valued function and section of $K^{-1}$ on U. With regards to the section $\beta_U$, the $U \neq U_0$ versions of $K^{-1}|_U$ come with an isomorphism $K^{-1}|_U = U \times \mathbb{C}$ that is defined using the chosen coordinates for U. This isomorphism is defined so that its inverse maps the constant section 1 of $U \times \mathbb{C}$ to the section that can written as $(1+\mathfrak{r})(\partial_z + \mathfrak{q})$ with $\mathfrak{r}$ being real and $\mathfrak{q}$ being orthogonal to $\partial_z$. Moreover, $|\mathfrak{r}| + |\mathfrak{q}| \leq c_0|z|$. Granted this isomorphism, the $U \neq U_0$ versions of $\beta_U$ are also viewed as functions on U.

What follows specifies the various $U \in \mathfrak{U}$ versions of $(a_U, (\alpha_U, \beta_U))$ on the complement in U of $\cup_{U' \in \mathfrak{U}-U} U \cap U'$.

WHEN $U = U_0$: *Fix an isomorphism* $E|_{U_0} = U_0 \times \mathbb{C}$. *Set* $(a_{U_0} = 0, (\alpha_{U_0} = 1, \beta_{U_0} = 0))$.

(A.7)

The definition of $(a_U, (\alpha_U, \beta_U))$ for $U \in \mathfrak{U}-U_0$ requires first the introduction of the rescaling map from $\mathbb{C}$ to itself that multiplies the coordinates by $z^{1/2}$. The latter map is denoted here by $r_z$. The definition refers to the functions y and $\varsigma$ on $\mathbb{C}$ that are depicted in (A.2) and the function $a_0$ on $\mathbb{C}$ given in (A.3).

WHEN $U = U_\gamma$: *Fix an isomorphism between* $E|_U$ *and* $U \times \mathbb{C}$. *Use Part 4's coordinates to identify* U *with the product of either* $S^1$ *or the appropriate interval with the radius* $\rho_*$ *concentric disk in* $D_0$. *Use* $(A_0, \alpha_0)$ *to denote the symmetric solution to (2.8) from* $\mathfrak{C}_1$ *with* $\alpha_0 = |\alpha_0| \frac{z}{|z|}$. *Set* $(a_U = i 2^{1/2} \nu r_z^* y dt - \frac{1}{2} r_z^* a_0 (z^{-1} dz - \bar{z}^{-1} d\bar{z}),\ \alpha_U = r_z^* \alpha_0,\ \beta_U = i\mu z^{-1/2} r_z^* \varsigma))$.

(A.8)

By way of a reminder, the pair $(y, \varsigma)$ are defined in (A.2).

*Part 6*: This part describes each $U \subset \mathfrak{U}$ version of $(a_U, (\alpha_U, \beta_U))$ on the intersection between U and $\cup_{U' \in \mathfrak{U}} (U \cap U')$. The transition function between a given set $U \in \mathfrak{U}-U_0$ and $U_0$ are as follows:

- *Suppose that* $\gamma \subset Y_{*\Lambda}$. *The transition function for* $U_0 \cap U_\gamma$ *identifies the constant section* 1 *of the bundle* $U_0 \times \mathbb{C}$ *with the section* $\frac{z}{|z|}$ *of the bundle* $U_\gamma \times \mathbb{C}$.
- *Suppose that* $\gamma \in \cup_{p \in \Lambda} \{\hat{\gamma}_p^+ \cup \hat{\gamma}_p^-\}$. *The transition function for* $U_0 \cap U_\gamma$ *on the part where* $|z| \leq 2c_\nu \rho_*$ *maps the section* 1 *of* $U_0 \times \mathbb{C}$ *to the section* $\frac{z}{|z|}$ of $U_\gamma \times \mathbb{C}$.

(A.9)

What is said in (A.9) is all that is needed for Proposition 3.1's map $\Phi^r$ because the latter has no cases where three sets from the cover intersect.

Suppose that $U \in \mathfrak{U} - U_0$. Use $\hat{U}$ to denote the part of $U \cap U_0$ that is considered in (A.9). The definition of $(a_U, (\alpha_U, \beta_U))$ on $\hat{U}$ follows. The definition introduces $\chi_{\hat{U}}$ to denote $\chi(\rho_*^{-1}|z| - 1)$.

- $a_U = \nu\, \chi_{\hat{U}}\, i 2^{1/2} r_z^* y\, dt - \frac{1}{2}(1 - \chi_{\hat{U}} + \chi_{\hat{U}} r_z^* a_0)(z^{-1}dz - \bar{z}^{-1} d\bar{z})$,
- $\alpha_U = (1 - \chi_{\hat{U}}(1 - r_z^*|\alpha_0|)) \frac{z}{|z|}$ ,
- $\beta_U = i \mu z^{-1/2} \chi_{\hat{U}} r_z^* \varsigma$ .

(A.10)

Equations (A.7)-(A.10) together define a smooth pair of connection on E's restriction to the $U_0 \cup (\cup_{\gamma \in \Theta} U_\gamma)$ and section $\mathbb{S}$ over this same set. In particular, they define a pair of connection on E's restriction to $Y_{*\Lambda} \cup T_{*\Lambda}$ and section of $\mathbb{S}$ over $Y_{*\Lambda} \cup T_{*\Lambda}$.

## b) Constraints

The operator in (1.17) will be analyzed in the case where the relevant version of $(A, \psi)$ is assumed to have five properties that are given momentarily. In particular, these properties are satisfied both by solutions to a given $(r, \mu)$ version of *(1.13)*. The upcoming Lemma A.1 asserts that these properties are also satisfied by the pairs that are constructed in Section Aa.

The upcoming properties refer to constants $c_0 \geq 100$ and $z \geq c_0^{10}$, These lower bounds are increased in subsequent subsections. The properties refer to a given pair $(A, \psi) \in \text{Conn}(E) \oplus C^\infty(Y; \mathbb{S})$. By way of a look ahead, a pair $(A, \psi)$ with these properties looks much like a solution to the $r = z$ and $\mu$ version of (1.13) with $\mu \in \Omega$ having $\mathcal{P}$-norm bounded by 1. The properties listed below are such that $(A, \psi)$ resembles such a solution in as much as

$$z^{-1/2}|B_A - z(\psi^\dagger \tau \psi - \hat{a})| + |D_A \psi| \leq c_0^{-6} z^{1/2} + c_0 .$$

(A.11)

The list of properties follows directly.

PROPERTY 1: *The section* $\psi = (\alpha, \beta)$ *is such that*
- $|\alpha| \leq 1 + c_0 z^{-1}$ *and* $|\beta| \leq c_0 z^{-1/2}$.
- $|\nabla_A \alpha|^2 \leq c_0(z|1 - |\alpha|^2| + 1)$.
- $|\nabla_A \beta| \leq c_0$.

The second property introduces the following notation: The section $D_A\psi$ of $\mathbb{S}$ is written as $([D_A\psi]_0, [D_A\psi]_1)$ with respect to the splitting $\mathbb{S} = E \oplus EK^{-1}$. As always, the Hodge dual of the curvature 2-form of A is denoted by $B_A$.

PROPERTY 2: *The 1-form* $B_A$ *and the section* $D_A\psi$ *of* $\mathbb{S}$ *are such that*

- $|\langle \hat{a}, B_A\rangle + i\,z(1 - |\alpha|^2)| \le c_0^{-20} z + c_0$.
- $|\hat{a} \wedge B_A| \le c_0^{6} z^{1/2} |1 - |\alpha|^2|^{1/2} + c_0$.
- $|[D_A\psi]_0| \le c_0^{6}$.
- $|[D_A\psi]_1| \le c_0^{-10} z^{1/2} + c_0$.

Note for future reference that the third bullets of PROPERTIES 1 and 2 have the following consequence: Let $(\nabla_A\alpha)_v$ denote the section of E that is obtained by pairing with $\nabla_A\alpha$ with $v$. The norm of this section $(\nabla_A\alpha)_v$ is bounded by $c_0\, c_0$ because $[D_A\psi]_0$ is the sum of $i(\nabla_A\alpha)_v$ with a linear combination of covariant derivatives of $\beta$.

The third property introduces $Y_{\Diamond z}$ to denote the subset of Y with distance at least $c_0^{4} z^{-1/2}$ from $\cup_{p\in\Lambda}(\hat{\gamma}_p^+ \cup \hat{\gamma}_p^-)$. This subset $Y_{\Diamond z}$ is a smooth manifold with boundary whose boundary components are tori. The third property also refers to the $(A, \psi)$ version of the connection $\hat{A}$ that is defined in (1.15).

PROPERTY 3: *The zero locus of* $\alpha$ *in* $Y_{\Diamond z}$ *is transversal and it consists of the disjoint union of at most* G *components with each a properly embedded arc or circle. The zero locus of* $\alpha$ *in* $Y_{\Diamond z}$ *has the following additional properties*:

- *Any given boundary component of* $Y_{\Diamond z}$ *contains either zero or two endpoints of the arc components of* $\alpha$*'s zero locus in* $Y_{\Diamond z}$. *If two, then the distance between them is at least* $100c_0^{2} z^{-1/2}$. *Moreover,* $u < 0$ *on one and* $u > 0$ *on the other.*
- *Suppose that* $\gamma$ *is a component of the zero locus of* $\alpha$ *in* $Y_{\Diamond z}$.
  a) *The unit length, oriented tangent vector to* $\gamma$ *has distance at most* $c_0 z^{-1/2}$ *from* $v$.
  b) *The curve* $\gamma$ *intersects any given* $\mathfrak{p} \in \Lambda$ *version of* $\mathcal{H}_{\mathfrak{p}}$ *where* $1 - 3\cos^2\theta > 0$.
  c) *If* $\gamma$ *is disjoint from a given boundary torus of* $Y_{\Diamond z}$, *then it has distance greater than* $3c_0^{3} z^{-1/2}$ *from this torus.*
  d) *If* $\gamma$ *intersects a boundary torus of* $Y_{\Diamond z}$, *then it does so only at its endpoints and these intersections are transversal.*
- *The intersection of* $\alpha$*'s zero locus with* $M_\delta$ *lies in the* $f^{-1}(1,2)$ *part of* $M_\delta$. *This intersection consists of* G *properly embedded segments that pair the index 1 and index 2 critical points of f in the sense that distinct segments start on the boundary of the* *radius δ coordinate balls about distinct index 1 critical points of f and end on the boundary of the radius δ coordinate balls about distinct index 2 critical points.*

- *The 2-form* $\frac{i}{2\pi} F_{\hat{A}}$ *has compact support and integral* 1 *on any transverse disk in* Y *with radius* $c_0 z^{-1/2}$ *and center at a zero of* α *in* $Y_{\lozenge z}$.

The fourth property constrains α away from its zero locus.

PROPERTY 4: *The absolute value of* $1 - |\alpha|^2$ *is less than* $c_0^{-10}$ *at all points with distance greater than* $c_0 z^{-1/2}$ *from the zero locus of* α *in* Y.

The final property is not strictly speaking required; it is imposed solely to avoid some extra effort. To set the notation, let p ∈ Y denote a given point. Fix a $\mathbb{C}$-linear isometry between $\mathbb{C}$ and the kernel of $\hat{a}$ at p. With $z$ given, use $\varphi_z$: $\mathbb{C} \to Y$ to denote the composition of first multipication on $\mathbb{C}$ = kernel($\hat{a}$) by $z^{-1/2}$ followed by the metric's exponential map. With $(A, \psi = (\alpha, \beta)) \in \text{Conn}(E) \times C^\infty(Y; \mathbb{S})$ given, use $(A_z, \alpha_z)$ to denote the pull-back of (A, α) using $\varphi_z$.

PROPERTY 5: *Fix* p ∈ Y. *The pair* $(A_z, \alpha_z)$ *has distance at most* $c_0^{-10}$ *in the* $C^4$*-topology from a solution to the vortex equations when restricted to the disk of radius* $c_0$ *with center at the origin in* $\mathbb{C}$.

To see an example of a pair with these properties, fix $r \geq c_0$ and $\mu \in \Omega$ with $\mathcal{P}$-norm bounded by 1. Every solution to the corresponding $(r, \mu)$ version of (1.13) has PROPERTIES 1-5 if $r \geq c_0$ with $c_0$ chosen less than $r^{1/6}$ and so as to avoid at most G intervals of a priori bounded length. By way of an explanation, the top three bullets in PROPERTY 1 are asserted by Lemma 2.1 and the fourth bullet has zero on its right hand side. Meanwhile (1.13) guarantees PROPERTY 2, Lemma 2.3 guarantees PROPERTY 4 and Lemma 2.9 guarantees PROPERTY 5. The first bullet of Proposition 2.4 gives Item a) of the second bullet of PROPERTY 3, its second bullet guarantees Item b) of PROPOSITION 3, its third bullet guarantees the third bullet of PROPERTY 3 and its fifth bullet guarantees the fourth bullet of PROPERTY 3. If the first bullet of PROPERTY 3 or Items c) or d) of the second bullet of PROPERTY 3 are not obeyed for a given choice of $c_0$, then at most G replacements of the form $c_0 \to c_0 + c_0$ will satisfy all of them. That this is so follows directly from the first three bullets of Proposition 2.4 and the formula for ν in (1.3).

The lemma that follows asserts that certain versions of the pairs (A, ψ) that are described in Section Aa also have the five properties listed above.

**Lemma A.1**: *There exists* $\kappa \geq 100$ *with the following significance: Fix parameters* $c_v \geq \kappa$ *and* $z \geq \kappa c_v^{10}$, *and then set* $\rho_* = c_v^2 z^{-1/2}$. *Fix a set* $T_{*\Lambda}$ *and then a set* $\Theta$ *as described by (A.5) which obeys the first and second bullets of the* $(z, c_0 = c_v)$ *version of* PROPERTY 3. *Suppose that* $(A, \psi) \in \mathrm{Conn}(E) \times C^\infty(Y; \mathbb{S})$ *is given by the* $(z, c_v, \rho_*)$ *version of (A.7)-(A.10) on* $Y_{*\Lambda} \cup T_{*\Lambda}$ *and that the* $(z, c_0 = c_v)$ *version of* PROPERTIES *1,2,4 and 5 hold on* $Y-(Y_{*\Lambda} \cup T_{*\Lambda})$. *Then* $(A, \psi)$ *obeys the* $(z, c_0 = c_v)$ *version of* PROPERTIES *1-5 on the whole of* Y.

***Proof of Lemma A.1***: This is a straightforward calculation given what is said in Section 3d and (3.3) about solutions to (2.8). Much the same calculation is done in Sections 2e and 2f of the article *Gr => SW* in [T8]. The specifics of the calculation are omitted.

## c) Bounds on eigenvectors

Suppose in what follows that $\mathfrak{c} = (A, \psi)$ satisfies PROPERTIES 1-5 in the previous subsection. The two lemmas in this subsection give some preliminary information about eigenvectors of the associated version of the operator $\mathfrak{L}_{\mathfrak{c},z}$, this being the $z = r$ version of the operator that is depicted in (1.17). The notation is such that $\nabla\mathfrak{b}$ is used to denote the covariant derivative of a given $\mathfrak{b} \in C^\infty(Y; iT^*Y \oplus \mathbb{S} \oplus i\mathbb{R})$ that is defined by the Levi-Civita covariant derivative on the $iT^*Y$ summand, the covariant derivative on sections of $\mathbb{S}$ that is define by the Levi-Civita connection with the connection A, and the exterior derivative on the sections of $i\mathbb{R}$. These lemmas also use $\|\cdot\|_2$ to denote the $L^2$ norm on Y.

**Lemma A.2**: *There exists* $\kappa \geq 100$ *and, given* $c_0 \geq \kappa$, *there exists* $\kappa_{c0} \geq \kappa$ *with the following significance: Fix* $z \geq \kappa_{c0} c_0^{10}$ *and suppose that* $\mathfrak{c} = (A, \psi) \in \mathrm{Conn}(E) \times C^\infty(Y; \mathbb{S})$ *obeys the* $(c_0, z)$ *version of* PROPERTIES 1 *and* 2 *in Section Ab.*

- *Let* $\mathfrak{b} = (b, \eta, \phi)$ *denote the an eigenvector of* $\mathfrak{L}_{\mathfrak{c},z}$. *Use* $\lambda$ *to denote the corresponding eigenvalue. Then* $\|\nabla\mathfrak{b}\|_2 \leq \kappa(\lambda + c_0 z^{1/2})\|\mathfrak{b}\|_2$.
- *Suppose in addition that* $\mathfrak{c}$ *obeys* PROPERTY 4 *in Section Ab and that* $|\lambda| \leq c_0^{-\kappa} z^{1/2}$. *Fix* $m > 2c_0$. *The* $L^2$ *norm of* $\mathfrak{b}$ *over the subset in* Y *with distance greater than* $mz^{-1/2}$ *from* $\alpha^{-1}(0)$ *is no greater than* $\kappa m^{-1}$.

To set the notation for the next lemma suppose that $\mathfrak{b} = (b, \eta, \phi)$ is a section of $iT^*Y \oplus \mathbb{S} \oplus i\mathbb{R}$. The lemma writes b as $b = b_0\hat{a} + b^\perp$ where $b^\perp$ annihilates $v$; and it writes $\eta$ with respect to the splitting $\mathbb{S} = E \oplus EK^{-1}$ as $\eta = (\eta_0, \eta_1)$. Lemma A.3 also uses $(\nabla b^\perp)_v$ and $(\nabla_A \eta_0)_v$ to denote the directional covariant derivatives along the vector field $v$.

**Lemma A.3**: *There exists* $\kappa \geq 100$*, and given* $c_0$*, there exists* $\kappa_{c0} \geq \kappa$ *with the following significance: Fix* $z \geq \kappa_{c0}\, c_0^{10}$ *and let* $\mathfrak{c} = (A, \psi) \in \mathrm{Conn}(E) \times C^\infty(Y; \mathbb{S})$ *denote a pair that obeys the* $(c_0, z)$ *version of* PROPERTIES 1-4 *in Section Ab. Suppose that* $\mathfrak{b} = (b, \eta, \phi)$ *is an eigenvector of* $\mathfrak{L}_{\mathfrak{c},z}$ *with* $L^2$ *norm equal to 1, and use* $\lambda$ *to denote* $\mathfrak{b}$*'s eigenvalue. Assume that* $|\lambda| \leq c_0^{-\kappa} z^{1/2}$.

- *The* $L^2$ *norms of* $b_0$*,* $\eta_1$ *and* $\phi$ *are bounded by* $c_0^{\kappa} z^{-1/2}$ *and the* $L^2$ *norms of their covariant derivatives are bounded by* $\kappa_{c0}$.
- *The* $L^2$ *norms of* $(\nabla b^\perp)_v$ *and* $(\nabla \eta_0)_v$ *are bounded by* $(\kappa_{c0} + \kappa|\lambda|)$.

The proofs of Lemmas A.2 and A.3 are given momentarily. A Bochner-Weitzenbock formula for $\mathfrak{L}_{\mathfrak{c},r}^{\;2}$ plays the central role in the arguments for these lemmas. To state this formula, fix $z \geq 1$ and let $\mathfrak{c} = (A, \psi)$ denote a pair in $\mathrm{Conn}(E) \times C^\infty(Y; \mathbb{S})$ and let $\mathfrak{L}_{\mathfrak{c},z}$ denote the corresponding version of (1.17). The respective $iT^*Y$, $\mathbb{S}$ and $i\mathbb{R}$ components of $\mathfrak{L}_{\mathfrak{c},z}^{\;2}\mathfrak{b}$ are:

- $\nabla^\dagger\nabla b + 2z|\psi|^2 b - 2^{1/2} z^{1/2}(\nabla\psi^\dagger\eta - \eta^\dagger\nabla\psi) + 2^{-1/2} z^{1/2}((D\psi)^\dagger\tau\eta + \eta^\dagger\tau D\psi) + \mathrm{Ric}(b)$ ,
- $D_A^{\;2}\eta + z[(\psi^\dagger\eta - \eta^\dagger\psi)\psi - cl(\psi^\dagger\tau\eta + \eta^\dagger\tau\psi)\psi] - 2^{3/2} z^{1/2}\langle b, \nabla\psi\rangle - 2^{1/2} z^{1/2}(cl(b) + \phi) D\psi$ ,
- $d^\dagger d\phi + 2z|\psi|^2\phi + 2^{-1/2} z^{1/2}((D\psi)^\dagger\eta - \eta^\dagger D\psi)$ .

(A.12)

To explain the notation, $\mathrm{Ric}(\cdot)$ denotes the endomorphism of $T^*Y$ defined by the Ricci tensor. Meanwhile, $\langle\ ,\ \rangle$ is defined as follows: Let $V \to Y$ denote any given vector bundle. Given V, then $\langle\ ,\ \rangle$ is the homomorphism from the bundle $T^*Y \otimes (V \otimes T^*Y)$ to V that is defined by the Riemannian metric.

***Proof of Lemma A.2***: To prove the first bullet, take the inner product between $\mathfrak{b}$ and $\mathfrak{L}_{\mathfrak{c},z}^{\;2}\,\mathfrak{b}$ and integrate the result over Y. Use (A.11) and (A.12) with an integration by parts to obtain the asserted bound. Note in this regard that the bounds on $|\psi|$ by $c_0$ and that on $|\nabla\psi|$ by $c_0 z^{1/2}$ are needed.

To prepare the stage for the proof of the second bullet, let $\mathfrak{c} \in \mathrm{Conn}(E) \times C^\infty(Y; \mathbb{S})$ denote any given element. What is written in (A.12) can be depicted schematically as

$$\mathfrak{L}_{\mathfrak{c},z}^{\;2}\,\mathfrak{b} = \nabla^\dagger\nabla\mathfrak{b} + 2z\mathfrak{b} + \mathfrak{e}(\mathfrak{b}) ,$$

(A.13)

where the endomorphism $\mathfrak{e}$ obeys $|\mathfrak{e}| \leq c_0(|B_A| + z^{1/2}|\nabla\psi| + c_0)$. If $\mathfrak{b}$ is an eigenvector of $\mathfrak{L}_{\mathfrak{c},z}$ with eigenvalue $\lambda$, then (A.13) leads to the inequality

$$\mathrm{d}^\dagger \mathrm{d}|\mathfrak{b}| + 2z|\mathfrak{b}| - |\mathfrak{e}||\mathfrak{b}| \leq \lambda^2|\mathfrak{b}| \,.$$

(A.14)

Now suppose that $\mathfrak{b}$ is as described in Lemma A.2's second bullet. The assumption that $\lambda \leq c_0^{-1}z^{1/2}$ and what is said in PROPERTIES 1, 2 and 4 in Section Ab have the following consequence: The inequality in (A.14) implies the more straightforward inequality

$$\mathrm{d}^\dagger \mathrm{d}|\mathfrak{b}| + z|\mathfrak{b}| \leq 0 \ .$$

(A.15)

at points with distance $c_0 z^{-1/2}$ or more from $\alpha^{-1}(0)$. To exploit this inequality, let $\Delta_\alpha(\cdot)$ denote for a moment the function $\mathrm{dist}(\cdot, \alpha^{-1}(0))$. Given $m \geq 2c_0$, convolve the function $\chi(2 - 2m^{-1}z^{1/2}\Delta_\alpha(\cdot))$ with a suitably chosen smoothing kernel to construct a non-negative function on Y with the following properties: Let $g_m$ denote this function. Then $g_m = 1$ where the distance to $\alpha^{-1}(0)$ is greater than $mz^{-1/2}$ and $g_m = 0$ where the distance to Y is less than $\frac{1}{2}mz^{-1/2}$. Furthermore, $|dg_m| \leq c_0 m^{-1}z^{1/2}$. Multiply both sides of (A.22) by $g_m^2|\mathfrak{b}|$ and integrate by parts. The resulting inequality implies the bound $z\| g_m\mathfrak{b}\|_2 \leq c_0 m^{-1}z\| \mathfrak{b}\|_2$. Divide both sides of this by $z$ to obtain what is asserted by Lemma A.2.

***Proof of Lemma A.3***: The bounds in the second bullet follow from those in the first from the form of $\mathfrak{L}_{c,z}$. Indeed, the relevant version of the equation $\mathfrak{L}_{c,z}\,\mathfrak{b} = \lambda\mathfrak{b}$ equates $(\nabla b^\perp)_v$ and $(\nabla\eta_0)_v$ with linear combinations of the following: First, covariant derivatives of $b_0$, $\eta_1$ and $\phi$. Second linear combinations of $z^{1/2}b_0$, $z^{1/2}\eta_1$, $z^{1/2}\phi$ times factors of $\alpha$ or its complex conjugate. Third, linear combinations of factors of $z^{1/2}b^\perp$ and $z^{1/2}\eta_0$ times factors of $\beta$ or its complex conjugate. Finally, components of $\lambda\,\mathfrak{b}$. This property of $\mathfrak{L}_{c,z}$ is directly evident from its depiction in the upcoming (A.26) and (A.27).

The proof of the first bullet has six steps.

<u>Step 1</u>. The asserted bounds are proved with the help of (A.12). The bounds for $\phi$ will use the formula in (A.12) for the $i\mathbb{R}$ component of $\mathfrak{L}_{c,z}^2\,\mathfrak{b}$. Those for $b_0$ are obtained with the help of the formula in (A.12) for the $iT^*Y$ component of $\mathfrak{L}_{c,z}^2\,\mathfrak{b}$ by projecting the latter onto the span in of $\hat{a}$. Those for $\eta_1$ are obtained using the formula in (A.12) for the $\mathbb{S}$ component of $\mathfrak{L}_{c,z}^2\,\mathfrak{b}$ by projecting the latter onto the $EK^{-1}$ summand of $\mathbb{S}$. In this regard, the projection of the $iT^*Y$ component of $\mathfrak{L}_{c,z}^2\,\mathfrak{b}$ to the span of $\hat{a}$ can be written as

$$d^{\dagger}d\, b_0 + 2z\, |\psi|^2 b_0 - 2^{1/2} z^{1/2}((\nabla\psi)_v{}^{\dagger}\eta - \eta^{\dagger}(\nabla\psi)_v) + 2^{-1/2} z^{1/2}((D\psi)^{\dagger} cl(\hat{a})\eta + \eta^{\dagger} cl(\hat{a}) D\psi) + \mathfrak{R}_0(\nabla b) + \langle \hat{a}, Ric(b)\rangle \,, \tag{A.16}$$

where $(\nabla\psi)_v$ denotes the directional covariant derivative along $v$ and where $\mathfrak{R}_0$ denotes a linear form on $T^*Y \otimes T^*Y$ that is defined by the covariant derivatives of $\hat{a}$. In particular, the latter is bounded in absolute value by $c_0$. Meanwhile, the projection of the $\mathbb{S}$ component of $\mathfrak{L}_{c,z}{}^2\, \mathfrak{b}$ to the $EK^{-1}$ summand of $\mathbb{S}$ can be written as

$$\nabla_A{}^{\dagger}\nabla_A\eta_1 + i\langle \hat{a}, B_A\rangle \eta_1 + 2z(|\alpha|^2 + |\beta|^2)\eta_1 + cl(B_A{}^{\perp})\eta_0 - 2^{3/2} z^{1/2}\langle b, \nabla\beta\rangle - 2^{1/2} z^{1/2}[(cl(b) + \phi) D\psi]_1 + \mathfrak{R}_1(\nabla\eta) + \mathfrak{r}_1(\eta) \,, \tag{A.17}$$

where the notation use $B_A{}^{\perp}$ to denote $B_A - \hat{a}\langle \hat{a}, B_A\rangle$, it uses $[\cdot]_1$ to denote the $EK^{-1}$ component of the given section of $\mathbb{S}$; and it uses $\mathfrak{R}_1$ and $\mathfrak{r}_1$ to denote endomorphisms that depend only on the Riemannian metric.

Step 2: The notation that follows uses $\xi$ to denote $(b_0, \eta_1, \phi)$ and it uses $\nabla\xi$ to denote the 3-tuple whose first and third entries are $db_0$ and $d\phi$, and whose second entry is $\nabla_A\eta_1$. Fix a constant $m_\Lambda > 8c_0{}^4$ to be determined shortly. Suppose in this step that the $L^2$ norm of $\xi$ over the part of $Y_{\Diamond z}$ with distance greater than $m_\Lambda z^{-1/2}$ from the boundary of $Y_{\Diamond z}$ is less than $m_\Lambda{}^{-1/4} \|\xi\|_2$.

Introduce $\theta_\Lambda$ to denote the characteristic function for the set of points in $Y_{\Diamond z}$ with distance at least $m_\Lambda z^{-1/2}$ from the boundary of $Y_{\Diamond z}$. Meanwhile, use the function $\chi$ to construct a smooth, non-negative function which is 1 where the distance to $Y - Y_{\Diamond z}$ is less than $2m_\Lambda z^{-1/2}$ and zero where the distance to this set is greater than $4m_\Lambda z^{-1/2}$. Use $\chi_\Lambda$ to denote this function. The function $\chi_\Lambda$ can and should be constructed so that its differential obeys $|d\chi_\Lambda| \le 16\, m_\Lambda{}^{-1} z^{1/2}$. Note that $|d\chi_\Lambda|$ has support where $\theta_\Lambda$ is equal to 1.

Take the $L^2$ inner product of the components of $\chi_\Lambda{}^2\xi$ with the relevant parts of the eigenvalue equation $\mathfrak{L}_{c,z}{}^2\, \mathfrak{b} = \lambda\mathfrak{b}$. Use the third bullet of (A.12), (A.16) and (A.17) with an integration by parts derive the following inequality:

$$\|\nabla(\chi_\Lambda\xi)\|_2{}^2 \le c_0\lambda^2 \|\chi_\Lambda\xi\|_2{}^2 + c_0 m_\Lambda{}^{-2} z \|\theta_\Lambda\xi\|_2{}^2 + c_0((c_0{}^{-10} z + c_0)\|\chi_\Lambda\xi\|_2{}^2 + c_0{}^6 z^{1/2} \|\chi_\Lambda\xi\|_2 \|\mathfrak{b}\|_2 + \|\mathfrak{b}\|_2{}^2). \tag{A.18}$$

The rest of this step explains how the various terms on in this inequality come about.

The term $\|\nabla(\chi_\Lambda\xi)\|_2{}^2$ on the left hand side and the term $c_0 m_\Lambda{}^{-2} z \,\|\theta_\Lambda\xi\|_2{}^2$ on the right hand side arise from the integration by parts that rewrites the $L^2$ inner product

between $\chi_\Lambda{}^2\xi$ and $\nabla^\dagger\nabla\xi$ as the square of the $L^2$ norm of $\nabla(\chi_\Lambda\xi)$ and a term with derivatives of $\chi_\Lambda$. The former accounts for the term on the left hand side of (A.18) and the latter accounts for the appearance of $c_0\, m_\Lambda{}^{-2} z\, \|\theta_\Lambda\xi\|_2{}^2$ on the right hand side of (A.18). These two terms with the term $c_0\|\mathfrak{b}\|_2{}^2$ on the right hand side of (A.18) also account for the $L^2$ inner product between $\chi_\Lambda{}^2\xi$ and the $\mathfrak{R}_0(\nabla\mathfrak{b})$ and $\mathfrak{R}_1(\nabla\mathfrak{b})$ terms in (A.12) and (A.13). The term $\lambda^2\|\chi_\Lambda\xi\|_2{}^2$ comes from the $L^2$ inner product between $\chi_\Lambda{}^2\xi$ and $\lambda\xi$.

The terms $(c_0{}^{-10}z + c_0)\|\chi_\Lambda\xi\|_2{}^2$ and $c_0{}^6 z^{1/2}\|\chi_\Lambda\xi\|_2\|\mathfrak{b}\|_2$ and $\|\mathfrak{b}\|_2{}^2$ on the right hand side of (A.18) account for the $L^2$ inner product between components of $\chi_\Lambda{}^2\xi$ and the various terms in the third bullet of (A.12), (A.16) and (A.17) that lack covariant derivatives of components of $\xi$ To elaborate, there are, first of all, the terms that have $2z|\psi|^2$ multiplying $\phi$ in (A.16) and $b_0$ in (A.12). These are discarded when writing (A.18) as they contribute nonpositive terms to the right hand side of (A.18). There is also a nonpositive contribution to the right hand side of (A.18) from the $2z(|\alpha|^2+|\beta|^2)\eta_1$ term in (A.17) and from $i\langle\hat{a}, B_A\rangle\eta_1$. PROPERTY 1 is used to rewrite this last term as $z(1-|\alpha|^2)\eta_1$ plus a remainder term that is bounded by $(c_0{}^{-10}z + c_0)|\eta_1|$. The remainder term is accounted for by a part of the $(c_0{}^{-10}z+c_0)\|\chi_\Lambda\xi\|_2{}^2$ term on the right hand side of (A.18).

The other terms without covariant derivatives of $\xi$ in the third bullet of (A.12), (A.16) and (A.17) are bounded by either

$$c_0|[D_A\psi]_1|\,|\xi| \quad or \quad c_0(|B_A{}^\perp| + |(\nabla\psi)_v| + |\nabla\beta| + |[D_A\psi]_0| + 1)|\mathfrak{b}|. \tag{A.19}$$

With (A.19) understood, what follows is a consequence of PROPERTIES 1 and 2: The terms without covariant derivatives of $\xi$ that are bounded by $c_0|[D_A\psi]_1|\,|\xi|$ are accounted for by the term $(c_0{}^{-10}z+c_0)\|\chi_\Lambda\xi\|_2{}^2$ on the right most right hand side of (A.18). Meanwhile, the terms without covariant derivatives of $\xi$ that are bounded by the right most expression in (A.19) are accounted for by the term $c_0{}^6 z^{1/2}\|\chi_\Lambda\xi\|_2\|\mathfrak{b}\|_2 + \|\mathfrak{b}\|_2{}^2$ in (A.18). Note in this regard that PROPERTIES 1 and 2 imply the bound $|(\nabla\psi)_v| \le c_0\, c_0{}^6$. This is stated explicitly with regards to $(\nabla\beta)_v$; and $|(\nabla_A\alpha)_v| \le c_0\, c_0{}^6$ because $[D_A\psi]_0$ is a sum of $i(\nabla_A\alpha)_v$ and linear combinations of covariant derivatives of $\beta$.

<u>Step 3</u>: Fix $m_\Lambda = 100\, c_0{}^4$ and use the assumption $\|\theta_\Lambda\xi\|_2 \le m_\Lambda{}^{-1/4}\|\xi\|_2$ to see that the right hand expression in (A.18) is at most $c_0(\lambda + c_0{}^{-10}z)\|\chi_\Lambda\xi\|_2{}^2 + c_0{}^k$ with $k \le c_0$. Meanwhile, the left hand side of (A.18) is no less than $c_0\, m_\Lambda{}^{-2} z\|\chi_\Lambda\xi\|_2$. Indeed, this follows from a standard Dirichelet eigenvalue inequality given that $|\chi_\Lambda\xi|$ has compact support in the radius $(m_\Lambda + c_0{}^4) z^{-1/2}$ tubular neighborhood of $\cup_{p\in\Lambda}\{\hat{\gamma}_p^+ \cup \hat{\gamma}_p^-\}$. These upper and lower bounds find $\|\xi\|_2{}^2 \le c_0\, c_0{}^k z^{-1}$ if $c_0 \ge c_0$ and $\lambda \le c_0{}^{-5} z^{1/2}$. This gives the first assertion of the first bullet of Lemma A.3.

Step 4: Fix $m_\Lambda = 100c_0^4$ so as to invoke the conclusions of Step 2. With $m_\Lambda$ fixed, use $\chi$ to construct a smooth, non-negative function on Y which is equal to 1 at distances greater than $m_\Lambda z^{-1/2}$ from Y−$Y_{\Diamond z}$ and equal to zero on Y−$Y_{\Diamond z}$. This function is denoted by $\chi_{\Diamond z}$. The function $\chi_{\Diamond z}$ can be constructed so that $|d\chi_{\Diamond z}| \le 32\, m_\Lambda^{-1} z^{-1/2}$. Fix a second constant $m \in (c_0, c_0^2)$. This step makes the following two assumptions:

- *The* $L^2$ *norm of* $\xi$ *over the part of* $Y_{\Diamond z}$ *with distance greater than* $m_\Lambda z^{-1/2}$ *from the boundary of* $Y_{\Diamond z}$ *is not less than* $m_\Lambda^{-1/4} \|\xi\|_2$.
- *The* $L^2$ *norm of* $\chi_{\Diamond z}\xi$ *over the part of* $Y_{\Diamond z}$ *with distance* $mz^{-1/2}$ *or more from the zero locus of* $\alpha$ *is greater than* $m^{-1/4}\|\chi_{\Diamond z}\xi\|_2$.

(A.20)

Use $\chi$ once more, now to define a smooth, non-negative function which is 1 where the distance to $\alpha^{-1}(0)$ is greater than $mz^{-1/2}$ and zero where the distance is less than $\frac{1}{2}mz^{-1/2}$. Let $\chi_m$ denote this function. Given that $m \le c_0^2$, what is said by the first bullet of PROPERTY 3 and what is said by Item c) of the second bullet of PROPERTY 3 imply that the function $\chi_m$ can be constructed so that its differential obeys $|d\chi_m| \le 16\, m^{-1} z^{1/2}$. This bound is assumed in what follows. Introduce $\theta_m$ to denote the characteristic function for the support of $|d\chi_m|$ and $\theta_{\Diamond z}$ to denote the characteristic function for the support of $|d\chi_{\Diamond z}|$.

Take the $L^2$ inner product of $(\chi_m\chi_{\Diamond z})^2\xi$ with the eigenvalue equation $\mathfrak{L}_{c,z}^2\,\mathfrak{b} = \lambda\,\mathfrak{b}$ and use either the third bullet of (A.13) or (A.16) or (A.17) with an integration by parts and (A.11) to derive from these integrals the inequality that follows:

$$\|\nabla(\chi_m\chi_{\Diamond z}\xi)\|_2^2 + \tfrac{1}{2} z\|\chi_m\chi_{\Diamond z}\xi\|_2^2 \le c_0\lambda^2\|\chi_m\chi_{\Diamond z}\xi\|_2^2 + c_0 z(m^{-2}\|\theta_m\xi\|_2^2 + m_\Lambda^{-2}\|\theta_{\Diamond z}\xi\|_2^2) + c_0 c_0^6 z^{1/2}\|\chi_m\chi_{\Diamond z}\xi\|_2 .$$

(A.21)

This proof of this inequality invokes PROPERTY 4, the bounds for the norms of the components of $B_A$ and $D_A\psi$ that are asserted in PROPERTY 2 and the bound for $|\nabla_A\beta|$ that is asserted by PROPERTY 1. As noted previously, these imply that $|(\nabla\psi)_\nu| \le c_0 c_0^6$.

To make something of (A.21), use the first bullet in (A.20) to conclude that $\|\theta_{\Diamond z}\xi\|_2 \le m_\Lambda^{1/4}\|\chi_{\Diamond z}\xi\|_2$ and then use the second to see that $\|\theta_{\Diamond z}\xi\|_2 \le (m\,m_\Lambda)^{1/4}\|\chi_m\chi_{\Diamond z}\xi\|_2$. Use this bound in (A.21) to conclude that

$$\tfrac{1}{4} z\|\chi_m\chi_\Lambda^c\xi\|_2^2 \le c_0\lambda^2\|\chi_m\chi_\Lambda^c\xi\|_2^2 + c_0(m_\Lambda^{1/4}m^{-3/4} + m^{1/4}m_\Lambda^{-3/4})z\|\xi\|_2^2 + c_0 c_0^k ,$$

(A.22)

with $k < c_0$. By assumption, $m \le m_\Lambda$ and so $m^{1/4} m_\Lambda^{-3/4} \le \frac{1}{100}$. Meanwhile, $m_\Lambda^{1/4} m^{-3/4} \le \frac{1}{100}$ if $m \ge c_0 c_0^{4/3}$. Assume this to be the case. If it is also the case that $\lambda \le c_0^{-1} z^{1/2}$, then (A.18) finds $\|\chi_m \chi_\Lambda^c \xi\|_2 \le c_0 c_0^{k'} z^{-1/2}$ with $k' \le c_0$ This last bound with (A.20) gives the bound $\|\xi\|_2^2$ by $c_0 c_0^{k''} z^{-1/2}$ with $k'' \le c_0$ if $m = c_0 c_0^{4/3}$ and if $\lambda \le c_0^{-1} z^{1/2}$.

Step 5: This step assumes that the top bullet in (A.20) is satisfied but that the lower bullet is violated. Use $\chi$ yet again, this time to construct a smooth, non-negative function which is 0 where the distance to $\alpha^{-1}(0)$ is greater than $2mz^{-1/2}$ and 1 where the distance is less than $mz^{-1/2}$. Denote this function by $\chi^c_m$. Given that $m \le c_0^2$, the function $\chi^c_m$ can and should be constructed so that $|d\chi^c_m| \le 32 m^{-1} z^{1/2}$. Use $\theta^c_m$ to denote the characteristic function for the support of $d\chi^c_m$.

Take the $L^2$ inner product of $(\chi^c_m \chi_{\diamond z})^2 \xi$ with the two sides of the equation in either the third bullet of (A.12) or (A.16) or (A.17) and use an integration by parts with (A.11) to see from these integrals that

$$\|\nabla(\chi^c_m \chi_{\diamond z} \xi_m)\|_2^2 \le c_0 (\lambda^2 + c_0^{-10} z) \|\chi^c_m \chi_{\diamond z} \xi\|_2^2 + m^{-2} z \|\theta^c_m \chi_{\diamond z} \xi\|_2 + m_\Lambda^{-2} z \|\theta_{\diamond z} \xi\|_2^2 + c_0^6 z^{1/2} \|\chi^c_m \chi_{\diamond z} \xi\|_2) . \tag{A.23}$$

The problematic terms on the right hand side of (A.23) those with $\|\theta_{\diamond z}\xi\|_2$ and $\|\theta^c_m \chi_{\diamond z}\xi\|_2$. The former is dealt with as follows: The top bullet in (A.20) asserts that $\|\theta_{\diamond z}\xi\|_2 \le m_\Lambda^{1/4} \|\chi_{\diamond z}\xi\|_2$. Hold on to this for the moment. The triangle inequality finds $\|\chi_{\diamond z}\xi\|_2 \le \|\chi^c_m \chi_{\diamond z}\xi\|_2 + \|(1-\chi^c_m)\chi_{\diamond z}\xi\|_2$, and thus $\|\chi_{\diamond z}\xi\|_2 \le \|\chi^c_m \chi_{\diamond z}\xi\|_2 + m^{-1/4} \|\chi_{\diamond z}\xi\|_2$ because the lower bullet in (A.20) is violated. Thus $\|\chi_{\diamond z}\xi\|_2 \le 2\|\chi^c_m \chi_{\diamond z}\xi\|_2$ and so the bound $\|\theta_{\diamond z}\xi\|_2 \le m_\Lambda^{1/4} \|\chi_{\diamond z}\xi\|_2$ implies that $\|\theta_{\diamond z}\xi\|_2 \le 2 m_\Lambda^{1/4} \|\chi^c_m \chi_{\diamond z}\xi\|_2$. Meanwhile, the problematic term $\|\theta^c_m \chi_{\diamond z}\xi\|_2$ is bounded by $c_0 m^{-1/4} \|\chi^c_m \chi_{\diamond z}\xi\|_2$ because the lower bullet in (A.20) is violated.

Insert the bounds in the preceding paragraph for $\|\theta_{\diamond z}\xi\|_2$ and $\|\theta^c_m \chi_{\diamond z}\xi\|_2$ in (A.23) and use the fact that $m_\Lambda = 100 c_0^4$ and $m = c_0 c_0^{4/3}$ to see that

$$\|\nabla(\chi^c_m \chi_{\diamond z}\xi_m)\|_2^2 \le c_0 (1 + (z^{-1}\lambda^2 + c_0^{-4} + m^{-5/2}) z) \|\chi^c_m \chi_{\diamond z}\xi_m\| + c_0 c_0^k , \tag{A.24}$$

with $k < c_0$. Consider now the left hand side of (A.24). To this end, let $D \subset Y_{\diamond\zeta}$ denote a transverse disk centered at a point in $\alpha^{-1}(0)$ with radius $2mz^{-1/2}$. The points in the support of $\chi^c_m \chi_{\diamond z}$ have distance at most $2mz^{-1/2}$ from $\alpha^{-1}(0)$, and so the Dirichelet inequality implies that the $L^2$ norm of $|\nabla(\chi^c_m \chi_{\diamond z}\xi_m)|$ over D is no less than $c_0^{-1} m^{-1} z^{1/2}$ times that of $\chi^c{}_m\chi_{\Diamond z}\xi_m$ over D. This being the case, (A.24) has the following consequence: Assume that $\lambda \le c_0^{-1}m^{-1}z^{1/2}$. Then $c_0^{-1}m^{-2}z\|\chi^c{}_m\chi_{\Diamond z}\xi_m\|_2^{\,2} \le c_0\,c_0^{\,k}$. This last bound leads directly to the desired upper bound on $\|\xi\|_2$.

Step 6: Steps 2-5 established Lemma A.3's claim about the $L^2$ norms of $\phi$, $b_0$ and $\eta_1$. Granted this claim, take the $L^2$ inner product of both sides of the eigenvalue equation $\mathfrak{L}_{c,z}^{\,2}\,\mathfrak{b} = \lambda\mathfrak{b}$ with $\xi$ and use the third bullet of (A.12), (A.16) and (A.17) with the bounds in PROPERTIES 1 and 2 to derive the bound $\|\nabla\xi\|_2^{\,2} \le \lambda^2\|\xi\|_2^{\,2} + c_0\,c_0^{\,k}(z\|\xi\|_2^{\,2} + 1)$ with $k \le c_0$. This last bound implies what Lemma A.3 asserts about the $L^2$ norm of $\nabla\xi$.

## d) The vortex operator

This subsection constitutes a digression to supply various observations that are used subsequently to say more about eigenvectors of $\mathfrak{L}_{c,z}$ near the zero locus of $\alpha$. The discussion here is given in four parts.

*Part 1*: Assume in what follows that $(A, \psi = (\alpha, \beta))$ obeys the constraints in Section Ab. Fix a point $p \in Y$ on the zero locus of $\alpha$ or on one of the curves from the set $\cup_{p\in\Lambda}\{\hat{\gamma}_p^{+} \cup \hat{\gamma}_p^{-}\}$ in a component of $Y - Y_{\Diamond z}$ that contains zeros of $\alpha$. In the former case, set $c_1 = 20c_0$ and in the latter case, set $c_1 = c_0^{\,4}$. Fix an isometric, $\mathbb{C}$-linear identification between $\text{kernel}(\hat{a})|_p$ and $\mathbb{C}$. With this identification understood, let $\varphi_z$ denote the map from $\mathbb{C}$ to Y that is obtained by composing first multiplication by $z^{-1/2}$ and then the metric's exponential map. Use $(A_z, \alpha_z)$ to denote the $\varphi_z$-pull-back of $(A, \alpha)$.
Use $\vartheta_z$ in what follows to denote the $(A_z, \alpha_z)$ version of (3.4)'s operator $\vartheta$. Of particular interest is this operator on concentric disks about the origin in $\mathbb{C}$ with radius $c_1$ or less.

The analysis of $\vartheta_z$ uses the following consequence of PROPERTIES 1 and 2 in Section Ab: The pair $(A_z, \alpha_z)$ comes close to solving (2.8)'s vortex equations on the radius $c_1$ disk centered at the origin in $\mathbb{C}$ in the sense that

$$|{*}F_{A_z} + i(1 - |\alpha_z|^2)| + |\bar{\partial}_{A_z}\alpha_z| \le c_0(c^{-1}|1 - |\alpha_z|^2| + c_0^{\,-10}). \tag{A.25}$$

The ramifications with regards to $\vartheta_z$ stem from the fact that the right most term in (3.6) is bounded by $c_0\,c_0^{\,-1}$ on the radius $20\,c_1$ disk about the origin if $z \ge c_0^{\,10}$. This suggests that $\vartheta_z\vartheta_z^{\dagger}$ is uniformly positive in a suitable sense. The constructions that follow are used to make a precise statement to this affect. These constructions assume that $z \ge \kappa_c c_0^{\,10}$ with $\kappa_c$ larger than the versions of $\kappa_{c0}$ that appear in Lemmas A.2 and A.3.

*Part 2*: Suppose that $k \in \{0, 1, \ldots, 7\}$ and that $\alpha$ lacks zeros in the open, concentric annulus in the transverse disk centered at p with respective inner and outer radii equal to $(c_1 + k c_0)z^{-1/2}$ and $(c_1 + (k+3) c_0)z^{-1/2}$. Zeros of $\alpha$ on the transverse disk through p with radius $c_1 z^{-1/2}$ correspond via $\varphi_z$ to the zeros of $\alpha_z$ on the corresponding $c_1$ disk about the origin in $\mathbb{C}$. In any event, use $\mathcal{A} \subset \mathbb{C}$ to denote the concentric annulus with inner radius $(c_1 + (k+1) c_0)$ and outer radius $(c_1 + (k+2) c_0)$.

It is a consequence of PROPERTY 4 that $|\alpha_z| \geq 1 - c_0^{-10}$ in $\mathcal{A}$. This being the case, there is an isomorphism over $\mathcal{A}$ between $\varphi_z{}^*E$ and $\mathcal{A} \times \mathbb{C}$ that maps $\alpha_z$ to $|\alpha_z|$ with the latter viewed at any given point as a complex number with zero imaginary part. Let $\theta_*$ denote the product connection on the trivial line bundle $\mathcal{A} \times \mathbb{C}$. This isomorphism pulls back $A_z$ as $\theta_* + a_z$ where $a_z$ is an $i\mathbb{R}$ valued 1-from on $\mathcal{A}$. The second bullet of PROPERTY 1 with PROPERTY 4 imply that $|a_z| \leq c_0 (c_0^{-9} + z^{-1/2})$.

Fix a non-negative, radial function on $\mathbb{C}$ which is equal to 1 where the distance to the origin is less than $(c_1 + (k + \frac{4}{3}) c_0)$ and equal to zero where the distance to the origin is greater than $(c_1 + (k + \frac{5}{3}) c_0)$. Choose a function whose derivative is bounded in absolute value by $10 c_0^{-1}$. Denote the chosen function by $\chi_*$.

Define a complex hermitian line bundle $E_z \to \mathbb{C}$ by identifying it with E on the radius $(c_1 + (k + \frac{5}{4}) c_0)$ disk about the origin in $\mathbb{C}$ and with the product bundle on the complement of the radius $(c_1 + (k+1) c_0)$ disk about the origin in $\mathbb{C}$. Use the isomorphism between $\varphi_z{}^*E|_{\mathcal{A}}$ and $\mathcal{A} \times \mathbb{C}$ to define the necessary clutching function. A unitary connection, $A_{z*}$, is defined on $E_z$ by setting $A_{z*} = A_z$ on the disk about the origin in $\mathbb{C}$ with radius $(c_1 + (k + \frac{5}{4}) c_0)$ and by setting $A_{z*} = \theta_* + \chi_* a_z$ on the complement of the disk about the origin with radius $(c_1 + (k+1) c_0)$. Use $\alpha_{z*}$ to denote the section of $E_z$ given by $\alpha$ over the radius $(c_1 + (k + \frac{5}{4}) c_0)$ disk centered at the origin and given by $(1 - \chi_*) + \chi_* |\alpha_z|$ over the complement of the radius $(c_1 + (k+1) c_0)$ disk. The connection $A_{z*}$ is flat and $\alpha_{z*}$ has norm 1 and is also $A_{z*}$-covariantly constant on the complement of the radius $(c_1 + \frac{2}{3} c_0)$ disk about the origin in $\mathbb{C}$. The pair $(A_{z*}, \alpha_{z*})$ also comes close to solving the vortex equations on the whole of $\mathbb{C}$ in the sense that (A.25) still holds.

*Part 3*: Use $n_*$ to denote the integral of $\frac{i}{2\pi} F_{A_{z*}}$ over $\mathbb{C}$. This is equal to 1 if the point p is in $Y_{\Diamond z}$, but can be greater than 1 if $p \in \cup_{p \in \Lambda} \{\hat{\gamma}_p^+ \cup \hat{\gamma}_p^-\}$. It follows from (A.25) with PROPERTIES 1 and 2 that $n_*$ is a positive integer.

The local Euler number of the zeros of $\alpha_{z*}$ sum to $n_*$ because $\alpha_{z*}$ has norm 1 and is $A_{z*}$ covariantly constant on the complement of a compact set in $\mathbb{C}$. The following lemma says more about these zeros.

**Lemma A.4**: *There exists* $\kappa > 100$ *such that if* $c_0 \geq \kappa$ *and* $z \geq \kappa c_0^{10}$, *then* $n_* \leq \kappa c_0^4$. *Moreover there exists an open set in* $\mathbb{C}$ *with the following three properties*:

- *The set can be covered by* $n_*$ *disks of radius* 4.
- $|\alpha_z| \geq \kappa^{-1}$ *on its complement.*
- *The sum of the local Euler numbers of the zeros of* $\alpha_z$ *in each component of this set is non-zero and positive.*

***Proof of Lemma A.4***: The bound on $n_*$ follows from PROPERTIES 3 and 4; and the other bullets follow from PROPERTY 5 using Lemma 2.9.

Fix a set in $\mathbb{C}$ that obeys the three bullets of Lemma A.4. Denote the set of components of this set by $Z_p$. Given $U \in Z_p$, use $\mathfrak{m}_U \in \{1, 2, \ldots, n_*\}$ to denote the sum of the local Euler numbers of $\alpha_z$ on U.

*Part 4*: Use $\vartheta_{z*}$ to denote the $(A_{z*}, \alpha_{z*})$ version of (3.4)'s operator $\vartheta$. The upcoming Lemma A.5 list some salient features of $\vartheta_{z*}$. This lemma uses $L^2(\mathbb{C}; \mathbb{C} \oplus E_z)$ to denote the completion of the vector space of smooth and compactly supported sections of the bundle $\mathbb{C} \times (\mathbb{C} \oplus E_z)$ using the norm whose square sends a given compactly supported section $\mathfrak{z} = (x, \iota)$ to the integral of $|\mathfrak{z}|^2$. This norm is denoted by $\|\cdot\|_2$. Use $\nabla\mathfrak{z}$ to denote $(dx, \nabla_{A_z}\iota)$. Lemma A.5 uses $L^2_1(\mathbb{C}; \mathbb{C} \oplus E_z)$ to denote the completion of this same space using the inner product whose square sends the given element $\mathfrak{z}$ to the integral of over $\mathbb{C}$ of $|\nabla\mathfrak{z}|^2 + |\mathfrak{z}|^2$. This defining norm for $L^2_1(\mathbb{C}; \mathbb{C} \oplus E_z)$ is denoted by $\|\cdot\|_{2,1}$.

**Lemma A.5**: *There exists* $\kappa > 100$ *such that what follows is true if* $c_0 \geq \kappa$ *and* $z \geq \kappa c_0^{10}$.

- *The operator* $\vartheta_{z*}$ *extends as a bounded, Fredholm operator from* $L^2_1(\mathbb{C}; \mathbb{C} \oplus E_z)$ *to* $L^2(\mathbb{C}; \mathbb{C} \oplus E_z)$ *with index equal to* $n_*$ *and with trivial cokernel.*
- *If* $\mathfrak{z} \in L^2_1(\mathbb{C}; \mathbb{C} \oplus E_z)$, *then* $\|\vartheta_{z*}^{\dagger}\mathfrak{z}\|_2 \geq \kappa^{-1}\|\mathfrak{z}\|_{2,1}$; *and if* $\mathfrak{z}$ *is* $L^2$*-orthogonal to the kernel of* $\vartheta_{z*}$, *then* $\|\vartheta_{z*}\mathfrak{z}\| \geq \kappa^{-1}\|\mathfrak{z}\|_{2,1}$.
- *Square integrable elements in the kernel of* $\vartheta_{z*}$ *are smooth and in* $L^2_1(\mathbb{C}; \mathbb{C} \oplus E_z)$.
- *If* $\mathfrak{z}$ *is an* $L^2(\mathbb{C}; \mathbb{C} \oplus E_z)$ *and in the kernel of* $\vartheta_{z*}$, *then* $|\mathfrak{z}| \leq \kappa \sum_{U \in Z_p} \mathfrak{m}_U e^{-\mathrm{dist}(\cdot, U)/2}$. *Moreover, if* $U \in Z_p$ *is a component with distance greater than* $\kappa$ *from the others and* *such that* $\mathfrak{m}_U = 1$*, then the* $L^2(\mathbb{C};\mathbb{C}\oplus E_z)$ *kernel of* $\vartheta_{z*}$ *has a non-zero element with the properties listed below. The list uses* $\mathfrak{z}_U$ *to denote this element.*

a) $|\mathfrak{z}_U| \le \kappa\, e^{-\mathrm{dist}(\cdot,U)/2}\, \|\mathfrak{z}_U\|_2$

b) *Any* $L^2(\mathbb{C};\mathbb{C}\oplus E_z)$ *element in the kernel of* $\vartheta_{z*}$ *can be written as* $x\mathfrak{z}_U + \mathfrak{z}'$ *with* $x \in \mathbb{C}$ *and with* $\mathfrak{z}'$ *such that* $|\mathfrak{z}'| \le \sum_{U'\in Z_p - U} \mathfrak{m}_{U'} e^{-\mathrm{dist}(\cdot,U')/2}$.

***Proof of Lemma A.5***: The fact that $\vartheta_{z*}$ defines a bounded map from $L^2_1(\mathbb{C};\mathbb{C}\oplus E_z)$ to $L^2(\mathbb{C};\mathbb{C}\oplus E_z)$ follows from the appearance of the $L^2$ norm of the covariant derivative in the definition of $L^2_1(\mathbb{C};\mathbb{C}\oplus E_z)$. To prove it Fredholm, it is necessary to prove that the kernel in $L^2_1(\mathbb{C};\mathbb{C}\oplus E_z)$ is finite dimensional, that the range is closed and that the cokernel is finite dimensional. The finite dimensional kernel and the closed range follow as a consequence of the Rellich lemma with the verification of the following:

*There exists* $\varepsilon > 0$ *and* $R \ge 1$ *such that* $\|\vartheta_{z*}\mathfrak{z}\|_2 \ge \varepsilon\|\mathfrak{z}\|_2$ *if the support of* $\mathfrak{z}$ *has compact support in the complement of the radius* R *disk in* $\mathbb{C}$ *about the origin.*

This follows by virtue of the fact that $\vartheta_{z*}^\dagger\vartheta_{z*}(x,\iota) = ((-\partial\bar\partial + \frac{1}{2})x, (-\partial_{A_{z*}}\bar\partial_{A_{z*}} + \frac{1}{2})\iota)$ where $A_{z*}$ is flat, $\alpha_{z*}$ has norm 1 and is also $A_{z*}$ covariantly constant.

The fact that the range is closed implies that the cokernel is isomorphic to the kernel of the adjoint. Standard elliptic regularity identifies the latter with the kernel of $\vartheta_{z*}^\dagger$. The fact that the latter is trivial can be seen using the $(A_{z*},\alpha_{z*})$ version of (3.6). This is done by invoking the bounds in (A.25) after commuting covariant derivatives to $\bar\partial_{A_0}\partial_{A_0}$ with $\partial_{A_{z*}}\bar\partial_{A_{z*}} + \frac{1}{2} + \mathfrak{e}$ with $|\mathfrak{e}| \le c_0\, c_0^{-1}(|1 - |\alpha_{z*}|^2 + 1)$

The fact that the dimension of the kernel is $n_*$ can be seen by comparing $\vartheta_{z*}$ with the version of $\vartheta$ that is defined by a pair $(A_0,\alpha_0)$ that obeys (2.8)'s vortex equations and is such $1 - |\alpha_0|^2$ is integrable and with integral equal to $2\pi n_*$. Such a comparison can be made by using what is said in Section 2a of [T9] to construct a [0,1]-parametrized path of pairs in $\mathrm{Conn}(E_z) \times C^\infty(\mathbb{C};E_z)$ that starts at $(A_{z*},\alpha_{z*})$, ends at such a solution to (2.8) and is such that each member of the family defines a Fredholm version of $\vartheta$. The construction of such a path amounts to little more than an exercise with cut-off functions and so no more will be said.

Granted the first bullet, the assertions of the second bullet are straightforward consequences of two facts, the first being that $\vartheta_{z*}$ is Fredholm with trivial cokernel and the second being (3.6). As for the third bullet, standard elliptic regularity arguments prove that the elements in the kernel of $\vartheta_{z*}$ are smooth. Meanwhile, the fact that the $L^2$ kernel of $\vartheta_{z*}$ coincides with its $L^2{}_1$ kernel follows from what was said above about $\vartheta_{z*}{}^{\dagger}\vartheta_{z*}$ where $A_{z*}$ is flat, $\alpha_{z*}$ has norm one and $\alpha_{z*}$ is $A_{z*}$-covariantly constant.

The assertions of the fourth bullet can be proved using the same sorts of arguments in Part 5 from Section 2a in [T9]. The modifications to these arguments are straightforward given that the properties listed in Section Ab imply that $(A_{z*}, \alpha_{z*})$ looks very much like a solution to the vortex equations with $1 - |\alpha_0|^2$ integrable and with integral equal to $2\pi n_*$. Note in particular what is said by (A.25). Note that Part 5 of Section 2a of [T9] states a stronger version of what is asserted by the fourth bullet for the version of $\vartheta$ that is defined using just such a solution to the vortex equations. As nothing fundamentally new is needed for the arguments in the case of $\vartheta_{z*}$, the details of the proof of the fourth bullet are omitted.

## e) The definition of $\mathrm{Ker}_\vartheta$ and $\Pi_\vartheta$

Assume here that $\mathfrak{c} = (A, \psi) \in \mathrm{Conn}(E) \oplus C^\infty(Y; \mathbb{S})$ obeys PROPERTIES 1-5 in Section Ab as defined with parameters $c_0$, $z$ and $c$ with $c_0$ and $z$ chosen so as to satisfy the requirements of Lemmas A.2-A.5. Parts 1 and 2 of this subsection uses versions of $\vartheta_{z*}$ to construct a complex line bundle over each component of $\alpha^{-1}(0)$ and a complex vector bundle over the sets that comprise an open cover of certain components of $\cup_{p\in\Lambda}\{\hat{\gamma}_p^+ \cup \hat{\gamma}_p^-\}$. This bundle is denoted in each case by $\mathrm{Ker}_\vartheta$. Part 3 defines a $\mathbb{C}$-linear homomorphism from the space of sections of $\bar{K} \oplus E \to Y$ to the space of sections of each version of $\mathrm{Ker}_\vartheta$, this denoted by $\Pi_\vartheta(\cdot)$. Part 4 defines a norm on the direct sum of these spaces of sections. This map is used to say more about eigenvectors of the operator $\mathfrak{L}_{c,z}$.

The construction of the associated complex line bundle and the associated homomorphism from $C^\infty(Y; \bar{K} \oplus E)$ for a component of $\alpha^{-1}(0)$ in $Y_{\Diamond z}$ mimics constructions in Section 3 of [T9]. The construction for a curve in $\cup_{p\in\Lambda}\{\hat{\gamma}_p^+ \cup \hat{\gamma}_p^-\}$ mimics constructions in Section 5 of [T9]. The definition of the norm also mimics what is done in Section 5 of [T9].

*Part 1*: Use $\gamma$ to denote a component of the zero locus of $\alpha$ in $Y_{z\Diamond}$. Given $p \in \gamma$, define the pair $(A_{z*}, \alpha_{z*})$ on $\mathbb{C}$. The $L^2$ kernel of the corresponding operator $\vartheta_{z*}$ is 1-dimensional since the integral of $\frac{i}{2\pi} F_{\hat{A}}$ over the radius $c_1 z^{-1/2}$ transverse disk with center p is equal to 1. The association to each point in $\gamma$ of the $L^2$ kernel of the corresponding version of $\vartheta_{z*}$ defines a complex line bundle over $\gamma$, this being $\mathrm{Ker}_\vartheta$.

*Part 2*: Use $\gamma$ now to denote an element in $\cup_{p\in\gamma}\{\hat{\gamma}_p^+ \cup \hat{\gamma}_p^-\}$. Consider first the case when there are no zeros of $\alpha$ on the nearby boundary component of $Y_{\Diamond z}$. The associated version of $\mathrm{Ker}_\vartheta$ is the zero dimensional bundle if the corresponding component of $Y-Y_{\Diamond z}$ has no zeros of $\alpha$. Suppose next that this component has zeros of $\alpha$. It follows from Item c) of the third bullet of PROPERTY 3 that any given $p \in \gamma$ version of the pair $(A_{z*},\alpha_{z*})$ can be defined using $k = 0$. This understood, let $n_*$ denote the integral over $\mathbb{C}$ of $\frac{i}{2\pi}$ times the curvature 2-form of $A_{z*}$. This positive integer does not depend on the chosen point in $\gamma$. Lemma A.5 asserts that any given $p \in \gamma$ version of $\vartheta_{z*}$ has $L^2$ kernel dimension equal to $n_*$. As p varies in $\gamma$, these $L^2$ kernels define a rank $n_*$ complex vector bundle over $\gamma$. This is the bundle $\mathrm{Ker}_\vartheta$.

Suppose next that $\gamma \in \cup_{p\in\Lambda}\{\hat{\gamma}_p^+ \cup \hat{\gamma}_p^-\}$ and that the nearby boundary component of $Y_{\Diamond z}$ has zeros of $\alpha$. Fix $p \in \gamma$ and let $D_0 \subset Y$ denote for the moment the transverse disk centered at p with radius $(c_0^4 + 10c_0)\, z^{-1/2}$. Granted that $c_0 \geq c_0$, use the second bullet of PROPERTY 3 with the formula for $v$ in (1.3) to find $k \in \{0, 1, \ldots, 7\}$ such that the following is true: The concentric, closed annulus in $D_0$ with inner radius $(c_0^4 + k\, c_0)\, z^{-1/2}$ and with outer radius $(c_0^4 + (k+3)\, c_0)\, z^{-1/2}$ has no zeros of $\alpha$. To say more about why this is so, suppose that $\upsilon$ is a connected, closed segment of an integral curve of $v$ with each endpoint having distance either $c_0^4 z^{-1/2}$ or $(c_0^4 + 10\, c_0)\, z^{-1/2}$ from $\gamma$. The formula in (1.3) implies that the $\phi$ angle changes monotonically on $\upsilon$ with total change being much less than $2\pi$ if $c_v > c_0$.

If $p \in \gamma$ and if $k \in \{0, 1, \ldots, 7\}$ and there are no zeros of $\alpha$ in the transverse disk centered at p with distance from p between $(c_0^4 + k\, c_0)\, z^{-1/2}$ and $(c_0^4 + (k+3)\, c_0)\, z^{-1/2}$, then such is the case for transverse disk centered at all points in some open neighborhood of p in $\gamma$. This being the case, $\gamma$ can be written as the union of 8 open sets, $\{\gamma_k\}_{k=0,1,\ldots,7}$ where the $\gamma_k$ corresponds to the subset of points in $\gamma$ where k has the property just described. The formula for $v$ in (1.3) implies that $\gamma_k$ will have at most 2 components.

Fix $k \in \{0, \ldots, 7\}$ and $p \in \gamma_k$. Use the chosen value for k to construct the pair $(A_{z*},\alpha_{z*})$ and the operator $\vartheta_{z*}$. The association to a point $p \in \gamma_k$ of the corresponding $L^2$ kernel of $\vartheta_{z*}$ defines a finite rank, complex vector bundle over $\gamma_k$. This bundle is $\mathrm{Ker}_\vartheta$.

*Part 3*: Let $\gamma_*$ denote either a component in $Y_{\Diamond z}$ of the zero locus of $\alpha$ or else a given $\gamma \in \cup_{p\in\Lambda}\{\hat{\gamma}_p^+ \cup \hat{\gamma}_p^-\}$ and $k \in \{0,\ldots,7\}$ version of $\gamma_k$. The associated C-linear map $\Pi_\vartheta$: $C^\infty(Y; \bar{K} \oplus E) \to C^\infty(\gamma_*; \mathrm{Ker}_\vartheta)$ is defined as follows: Fix $p \in \gamma_*$ and re-introduce the map $\varphi_z$ and the function $\chi_*$ from Part 2 of Section Ad that is used to define the corresponding pair $(A_{z*},\alpha_{z*})$. If $\mathfrak{f}_0$ is a section of $\bar{K} \oplus E$, then $\chi_*\varphi_z{}^*(\mathfrak{f}_0)$ defines an element in $C^\infty(\mathbb{C}; \mathbb{C} \oplus E_{z*})$ with compact support. The $L^2$ orthogonal projection of $\chi_*\varphi_z{}^*(\mathfrak{f}_0)$ to the $L^2$ kernel of $\vartheta_{z*}$ is the value of the section $\Pi_\vartheta(\mathfrak{f}_0)$ at p.

Introduce $\Theta_*$ to denote the set whose elements are the components of $\alpha$'s zero locus in $Y_{\Diamond z}$ and the various $\gamma \in \cup_{p\in\Lambda}\{\hat\gamma_p^+ \cup \hat\gamma_p^-\}$ and $k \in \{0, 1, \ldots, 7\}$ version of $\gamma_k$ with it understood in the latter case that $\gamma = \gamma_{k=0}$ and $\gamma_k = \emptyset$ for $k > 0$ if the nearby boundary component of $Y_{\Diamond z}$ lacks zeros of $\alpha$. The map $\Pi_\vartheta$ is viewed in what follows as a $\mathbb{C}$-linear map from $C^\infty(Y;\mathbb{V}_0)$ to $\oplus_{\gamma_*\in\Theta_*} C^\infty(\gamma_*;\mathrm{Ker}_\vartheta)$.

*Part 4*: This last part of the subsection defines a version of the $L^2$ norm on $\oplus_{\gamma_*\in\Theta_*} C^\infty(\gamma_*;\mathrm{Ker}_\vartheta)$. To this end, let $\mathfrak{q}$ denote an element in this vector space. The corresponding norm is denoted by $\|\mathfrak{q}\|_2$. The definition that follows writes a given $\gamma_* \in \Theta_*$ component of $\mathfrak{q}$ as $\mathfrak{q}_{\gamma_*}$ and it writes the integral over $\gamma_*$ of $|\mathfrak{q}_{\gamma_*}|^2$ as $\|\mathfrak{q}_{\gamma_*}\|_2$. Granted this notation, $\|\mathfrak{q}\|_2^2 = \sum_{\gamma_*\in\Theta_*} \|\mathfrak{q}_{\gamma_*}\|_2^2$.

## f) Rewriting $\mathfrak{L}_{(\cdot)}$

What follows in this subsection is used subsequently to bring what is said by Sections Ad and Ae into the $\mathfrak{L}_{c,z}$ story. It is necessary to start by introducing some new notation. The annihilator of $v$ in $T^*Y$ is defined to be the 2-dimensional subbundle of $T^*Y$ that is orthogonal to the 1-form $\hat a$. This subbundle is dual to kernel($\hat a$). The almost complex structure J splits its complexification as $K \oplus \bar K$ with it understood that K annihilates the -i eigenbundle of J's action on the complexification of the kernel of $\hat a$.

Introduce $I_{\mathbb{C}}$ to denote the product bundle $Y\times\mathbb{C}$. Write the complexification of the direct sum of the line $\mathbb{R}\hat a \subset T^*Y$ with $I_{\mathbb{R}}$ as $I_{\mathbb{C}} \oplus \bar I_{\mathbb{C}}$ with it understood that the projection to the $I_{\mathbb{C}}$ of point $(b_0\hat a, \phi) \in \mathbb{R}\hat a \oplus I_{\mathbb{R}}$ is $-b_0 + i\phi$.

Introduce $\mathbb{V}_0$ to denote $\bar K \oplus E$ and $\mathbb{V}_1$ to denote $I_{\mathbb{C}} \oplus EK^{-1}$. Define an $\mathbb{R}$ linear isomorphism from $iT^*Y \oplus \mathbb{S} \oplus iI_{\mathbb{R}}$ to $\mathbb{V}_0 \oplus \mathbb{V}_1$ as follows: Let $(b,\eta,\phi)$ denote a given point in $iT^*Y \oplus \mathbb{S} \oplus iI_{\mathbb{R}}$. Write b as $b_0\hat a + b^\perp$ and use $q$ to denote the orthogonal projection of $b^\perp$ to the subbundle $\bar K$ of the complexification of the annihilator of $v$. Use $p$ to denote $(-b_0 + i\phi)$ Write $\eta$ as $(\eta_0, \eta_1)$ using the identification of $\mathbb{S}$ with $E \oplus EK^{-1}$. The desired isomorphism sends $(b,\eta,\phi)$ to $((q, \eta_0), (p, \eta_1)) \in \mathbb{V}_0 \oplus \mathbb{V}_1$. This isomorphism is used in what follows to write any given section of $iT^*Y \oplus \mathbb{R} \oplus iI_{\mathbb{R}}$ as a section of $\mathbb{V}_0 \oplus \mathbb{V}_1$ and vice versa.

To continue setting notation, suppose that $q$ is a section of $\bar K$. Use $\partial^K q$ to denote the projection of the covariant derivative of $q$ to the $\bar K \otimes K$ summand in $\bar K \otimes (T^*Y)_{\mathbb{C}}$. Meanwhile, use $(\nabla q)_v$ to denote the corresponding projection to the $\bar K \otimes (\mathbb{C}\hat a)$ summand. When $p$ denotes a section of $I_{\mathbb{C}}$, use $\bar\partial^K p$ and $(\nabla p)_v$ to denote the respective projections of $dp$ to the $\bar K$ and $\mathbb{C}\hat a$ summands of $(T^*Y)_{\mathbb{C}}$. When $\eta = (\eta_0,\eta_1)$ denotes a section of $\mathbb{S} =$ $E\oplus EK^{-1}$, write the directional covariant derivatives of $\eta_0$ and $\eta_1$ along $v$ as $(\nabla_A\eta_0)_v$ and $(\nabla_A\eta_1)_v$, write the $\bar{K}$ part of the covariant derivative of $\eta_0$ as $\bar{\partial}_A^K\eta_0$ and write the K part of the covariant derivative of $\eta_1$ as $\partial_A^K\eta_1$.

With this notation understood, the operator $\mathfrak{L}_{c,v}$ can be veiwed in the manner of (3.13) and (3.14) in [T?SE2] as an operator mapping $C^\infty(Y;\mathbb{V}_0\oplus\mathbb{V}_1)$ to itself. Viewed in this light, the operator is denoted by $\mathfrak{L}_{\mathbb{V}}$. Let $\mathfrak{f}$ denote a given section of $C^\infty(Y;\mathbb{V}_0\oplus\mathbb{V}_1)$. To write $\mathfrak{L}_{\mathbb{V}}\mathfrak{f}$, first write the $\mathbb{V}_0 = \bar{K}\oplus E$ component of $\mathfrak{f}$ as $(q,\eta_0)$ and the $\mathbb{V}_1 = I_{\mathbb{C}}\oplus EK^{-1}$ component as $(p,\eta_1)$. The $\bar{K}$ and E summands of the $\mathbb{V}_0$ component of $\mathfrak{L}_{\mathbb{V}}\mathfrak{f}$ are

- $i(\nabla q)_v - 2i(-\bar{\partial}^K p + \frac{1}{\sqrt{2}} z^{1/2}\bar{\alpha}\eta_1) - \sqrt{2}\, iz^{1/2}\,\bar{\eta}_0\,\beta + t_{0q}q$ ,
- $i(\nabla_A\eta_0)_v - 2i(-\partial_A^K\eta_1 + \frac{1}{\sqrt{2}} z^{1/2}\alpha p) - \sqrt{2}\, iz^{1/2}\,\bar{q}\,\beta + t_{0\eta}\eta_1$ ;

(A.26)

and the respective $I_{\mathbb{C}}$ and $K^{-1}E$ summands of the $\mathbb{V}_1$ component of $\mathfrak{L}_{\mathbb{V}}\mathfrak{f}$ are

- $-i(\nabla p)_v + 2i(\partial^K q + \frac{1}{\sqrt{2}} z^{1/2}\bar{\alpha}\eta_0) + \sqrt{2}\, iz^{1/2}\bar{\eta}_1\,\beta + t_{1q}q + t_{1p}p$ ,
- $-i(\nabla_A\eta_1)_v + 2i(\bar{\partial}_A^K\eta_0 + \frac{1}{\sqrt{2}} z^{1/2}\alpha q) + \sqrt{2}\, iz^{1/2}\,\bar{p}\,\beta + t_{1\eta}\eta_0$ .

(A.27)

Here, each $t_{**}$ denotes an $\mathbb{R}$-linear homomorphism between summands of $\mathbb{V}_0\oplus\mathbb{V}_1$ that depends only on the metric and has norm bounded by $c_0$.

The description of the $\mathbb{V}_1$ component of $\mathfrak{L}_{\mathbb{V}}\mathfrak{f}$ given in (A.27) proves sufficient for what is to come. The description in (A.26) of the $\mathbb{V}_0$ component requires some additional rewriting. To begin this task, use $\gamma$ now to denote a small length open segment of $\alpha$'s zero locus in $Y_{\Diamond z}$ or a curve from $\cup_{p\in\Lambda}\{\hat{\gamma}_p^+\cup\hat{\gamma}_p^-\}$ whose corresponding component of $Y-Y_{\Diamond z}$ has a zero of $\alpha$. Fix a point $p\in\gamma$ and use the coordinates from Part 4 of Section Aa to parametrize a neighborhood of $\gamma$ in Y. This neighborhood is denoted for now by T. Nothing is lost by taking the t = 0 point to be the value of $\gamma$'s affine parameter at p. The segment $\gamma$ is assumed in what follows to be parametrized by $t\in(-\rho,\rho)$ with $\rho$ a constant that will be specified in the applications to come. This constant $\rho$ in any event obeys $\rho\in(c_0z^{-1/2}, c_0^{-1})$.

Each constant t slice of $(-\rho,\rho)$ is the intersection between the transverse disk through the corresponding point in $\gamma$ and the tubular neighborhood of $\gamma$. An identification of $E|_T$ with E's restriction to the t = 0 slice of T writes the directional covariant derivative along $\partial_t$ of A as $(\nabla_A)_{\partial_t} = \partial_t + a_{A0}$ where $a_{A0}$ is an $i\mathbb{R}$ valued function on T. With an identification of this sort chosen, then the terms $(\nabla q)_v$ and $(\nabla_A\eta_0)_v$ in (A.26) can be written using (A.6) as

- $\frac{\partial}{\partial t} q + 2i(\nu(z-x_\gamma)+\mu(\bar z - \bar x_\gamma))\frac{\partial}{\partial z} q - 2i(\nu(\bar z - \bar x_\gamma) + \bar\mu(z-x_\gamma))\frac{\partial}{\partial \bar z} q + \mathfrak{r}_q \cdot dq$
- $\frac{\partial}{\partial t}\eta_0 + a_{A0}\eta_0 + 2i(\nu(z-x_\gamma)+\mu(\bar z - \bar x_\gamma))\partial_A\eta_0 - 2i(\nu(\bar z - \bar x_\gamma)+\bar\mu(z-x_\gamma))\bar\partial_A\eta_0 + \mathfrak{r}_\eta\cdot\nabla_A q$

(A.28)

where the notation is such that $\partial_A$ is the covariant version of $\frac{\partial}{\partial z}$ and $\bar\partial_A$ is the covariant version of $\frac{\partial}{\partial \bar z}$. What are denoted by $\mathfrak{r}_q$ and $\mathfrak{r}_\eta$ obey $|\mathfrak{r}_q| + |\mathfrak{r}_\eta| \le c_0|z|(z^{-1/2} + |z|)$. Meanwhile, the terms just to the right of $(\nabla q)_\nu$ and $(\nabla_A\eta_0)_\nu$ in (A.26) can be written as

- $-2i(-\bar\partial^K p + \frac{1}{\sqrt2} z^{1/2}\bar\alpha\eta_1) = -2i(-\bar\partial p + \frac{1}{\sqrt2} z^{1/2}\bar\alpha\eta_1) + \mathfrak{e}_q\cdot\nabla p$,
- $-2i(-\partial_A^K\eta_1 + \frac{1}{\sqrt2} z^{1/2}\alpha p) = -2i(-\partial_A\eta_1 + \frac{1}{\sqrt2} z^{1/2}\alpha p) + \mathfrak{e}_\eta\cdot\nabla_A\eta_1$,

(A.29)

where $|\mathfrak{e}_q| + |\mathfrak{e}_\eta| \le c_0|z|^2$. The remaining terms in (A.26) can be written as

- $-\sqrt2 iz^{1/2}\bar\eta_0\beta + 2\nu q + 2\mu\bar q + \mathfrak{w}\cdot q$
- $-\sqrt2 iz^{1/2}\bar q\beta + \mathfrak{t}_{0\eta}\eta_1$

(A.30)

where $|\mathfrak{w}| \le c_0|z|$ and $|\mathfrak{t}_{0\eta}| \le c_0$.

### g) The operator $\mathfrak{L}_{\mathbb{V}}$ and $\Pi_\vartheta$

The next lemma hints at the role played by $\Pi_\vartheta$ as it talks about the $\Pi_\vartheta$ image of the $\mathbb{V}_0$ part of an eigenvector of $\mathfrak{L}_{\mathbb{V}}$. This lemma and subsequent discussions abuse notation to some extent by using $\Pi_\vartheta$ to denote two maps to $\oplus_{\gamma_*\in\Theta_*}\mathrm{Ker}_\vartheta$. The first is Section Ae's map from $C^\infty(Y; \mathbb{V}_0)$ and the second is the map from $C^\infty(Y; \mathbb{V}_0\oplus\mathbb{V}_1)$ that is obtained Section Ae's map by first projecting to the $\mathbb{V}_0$ summand.

**Lemma A.6**: *There exists $\kappa \ge 1$ and given $c_0 \ge \kappa$, there exists $\kappa_{c0} \ge \kappa$ with the following significance: Fix $c_0 \ge \kappa$ and $z \ge \kappa_{c0} c_0^{10}$ Suppose that $\mathfrak{c} = (A,\psi)$ obeys the corresponding version of* PROPERTIES 1-5 *in Section Ab. The assumptions of Lemmas A.2-A.5 are satisfied, and this understood, let $\mathfrak{f}$ denote an eigenvector of the operator $\mathfrak{L}_{\mathbb{V}}$ with eigenvalue bounded in absolute value by $c_0^{-\kappa} z^{1/2}$. Then $\|\Pi_\vartheta\mathfrak{f}\|_2 \ge (1-\kappa c_0^{-1})\|\mathfrak{f}\|_2$.*

***Proof of Lemma A.6***: Choose $c_0$ and $z$ so that Lemmas A.2-A.5 can be invoked. Write $\mathfrak{f}$ as $(\mathfrak{f}_0,\mathfrak{f}_1)$. Lemma A.3 finds $\|\mathfrak{f}_1\|_2 \le c_0 c_0^k z^{-1/2}\|\mathfrak{f}\|_2$ with $k \le c_0$ so $\mathfrak{f}_0$ accounts for most of the $L^2$ norm of $\mathfrak{f}$. Meanwhile, the second bullet of Lemma A.2 asserts that the $L^2$ norm of $\mathfrak{f}$ on the set of points with distance $2c_0 z^{-1/2}$ or more from $\alpha^{-1}(0)$ is bounded by $c_0 c_0^{-1}\|\mathfrak{f}\|_2$. As a consequence, the bulk of the $L^2$ norm of $\mathfrak{f}_0$ is accounted for by its $L^2$ norm on the radius $2c_0^{-1}z^{-1/2}$ tubular neighborhood of α's zero locus. The contribution to the $L^2$ norm from this part of Y is analyzed in the four steps that follow.

Step 1: Reintroduce the set $\Theta_*$ from Part 3 of Section Ae and let $\gamma_*$ denote a given element in $\Theta_*$. This is to say that $\gamma_*$ is either a component of α's zero locus in $Y_{\Diamond z}$ or some $\gamma \in \cup_{p\in\Lambda}\{\hat{\gamma}_p^+ \cup \hat{\gamma}_p^-\}$ and $k \in \{0, \ldots, 7\}$ version of $\gamma_k$. Each $p \in \gamma_*$ has an associated version of map $\varphi_z$ and function $\chi_*$ on $\mathbb{C}$ as described in Part 2 of Section Ad. In particular, the assignment to each point in $\gamma_*$ of the $L^2$ norm over $\mathbb{C}$ of the corresponding version of $\chi_*\varphi_z{}^*(\mathfrak{f}_0)$ defines a function on $\gamma_*$. The second bullet Lemma A.2 implies that

$$\Sigma_{\gamma_*\in\Theta_*}\int_{\gamma_*}\|\chi_*\varphi_z{}^*(\mathfrak{f}_0)\|_2^2 \geq (1-c_0c_0^{-1})\|\mathfrak{f}\|_2 . \tag{A.31}$$

This inequality is exploited in Step 4.

Step 2: Fix $\gamma_* \in \Theta_*$. If γ is a component of α's zero locus in $Y_{\Diamond z}$, set $c_1 = c_0$ and if not, set $c_1 = c_0^4$, this being the definition of $c_1$ that is used in Part 2 of Section Ad to construct the versions of $\vartheta_{z*}$ that are associated to the points in $\gamma_*$. Use $T_{\gamma_*}$ in what follows to denote the union of the radius $(c_1+10c_0)\,z^{-1/2}$ transverse disks with centers on $\gamma_*$.

Write $\mathfrak{f}_0 = (q,\eta_0)$ and assign to each $p \in \gamma_*$ the element $((\varphi_z^{-1})^*(\chi_*))\,\mathfrak{f}_0$ this being a section of $\mathbb{V}_0$ over the transverse disk through p whose components are written as $(q_*,\eta_{0*})$ These sections define smooth a section of $\mathbb{V}_0$ over $T_{\gamma_*}$ and they are viewed in this way. Use (A.27) with Lemmas A.3 and PROPERTY 1 to see that

$$\|(\partial^K q_* + \tfrac{1}{\sqrt{2}} z^{1/2}\bar{\alpha}\eta_{0*})\|_2^2 + \|(\bar{\partial}_A^K\eta_{0*} + \tfrac{1}{\sqrt{2}} z^{1/2}\alpha q_*)\|_2^2 \leq c_0(\lambda^2 + c_0^{-2}z)\,\|\mathfrak{f}\|_2^2 \tag{A.32}$$

with it understood that the $L^2$ norms on the right hand side denote integrals over $T_{\gamma_*}$

Step 3: Use $\mathfrak{z}_*$ in what follows to denote any given $p \in \gamma$ version of $\chi_*\varphi_z{}^*(\mathfrak{f}_0)$. The operator $\vartheta_{z*}$ enters the story by virtue of the fact that

$$\vartheta_{z*}\mathfrak{z}_* = z^{-1/2}\varphi_z{}^*(\partial^K q_* + \tfrac{1}{\sqrt{2}} z^{1/2}\bar{\alpha}\eta_{0*},\ \bar{\partial}_A^K\eta_{0*} + \tfrac{1}{\sqrt{2}} z^{1/2}\alpha q_*) + \mathfrak{e} \tag{A.33}$$

where the term denoted by $\mathfrak{e}$ has compact support in the radius $c_1+(k+2)c_0$ disk about the origin in $\mathbb{C}$ and $|\mathfrak{e}|$ is no greater than than the $\varphi_z$-pull-back of $c_0(c_0z^{-1}|\nabla\mathfrak{f}_0| + c_0^{-1}|\mathfrak{f}_0|)$. The argument to prove (A.33) is identical but for notation to that used in [T9] to derive the latter's (3.14) and (3.15).

Step 4: By definition, $\Pi_\vartheta\mathfrak{f}$ at points on $\gamma_*$ is the $L^2$-orthogonal projection of $\mathfrak{z}_*$ to the $L^2$-kernel of $\vartheta_{z*}$. This understood, write $\mathfrak{z}_* = \Pi_\vartheta\mathfrak{f} + \mathfrak{z}_*^\perp$. As the point in $\gamma_*$ varies, so $\mathfrak{z}_*^\perp$ varies and this understood, view $\|\mathfrak{z}_*^\perp\|_2$ as a function on $\gamma_*$. Lemma A.5 asserts that $\|\vartheta_{z*}\mathfrak{z}_*^\perp\|_2 \geq c_0^{-1}\|\mathfrak{z}_*^\perp\|_2$. This bound with Lemma A.2's first bullet and (A.32) and (A.33) imply that

$$(1 - c_0\mathcal{c}_0^{-1})\int_{\gamma_*}\|\chi_*\varphi_z{}^*(\mathfrak{f}_0)\|_2^2 - \int_{\gamma_*}|\Pi_\vartheta(\mathfrak{f}_0)|_2^2 \leq c_0(\lambda^2 z^{-1} + \mathcal{c}_0^{-2})\|\mathfrak{f}\|_2^2 \tag{A.34}$$

when $\mathcal{c}_0 > c_0$ and $z \geq c_{\mathcal{c}0}$ with the latter constant depending only on $\mathcal{c}_0$.

The inequalities in (A.31) and (A.34) imply that $\|\Pi_\vartheta\mathfrak{f}\|_2^2 \geq (1 - c_0\mathcal{c}_0^{-1})\|\mathfrak{f}\|_2^2$ if it is the case that $\mathcal{c}_0 > c_0$ and $|\lambda| \leq c_0^{-1}\mathcal{c}_0^{-1}z^{1/2}$.

**h) The equation $\Pi_\vartheta\mathfrak{L}_{\mathbb{V}}\mathfrak{f} = \lambda\Pi_\vartheta\mathfrak{f}$ on $Y - Y_{\Diamond z}$**

What is asserted by Lemma A.6 implies that $\Pi_\vartheta\mathfrak{f}$ determines $\mathfrak{f}$ for the most part if $\mathfrak{f}$ is an eigenvector of $\mathfrak{L}_{\mathbb{V}}$ whose corresponding eigenvalue is greater than $-c_0^{-1}\mathcal{c}_0^{-1}z^{1/2}$ but less than $c_0^{-1}\mathcal{c}_0^{-1}z^{1/2}$. This fact lies behind the focus in this subsection and the next on the $\Pi_\vartheta$ projection of the eigenvalue equation $\mathfrak{L}_{\mathbb{V}}\mathfrak{f} = \lambda\mathfrak{f}$. By way of a look at what is to come, (A.26) with (A.28)-(A.30) are used here to rewrite a given $\gamma_* \in \Theta_*$ component of the projected equation $\Pi_\vartheta(\mathfrak{L}_{\mathbb{V}}\mathfrak{f}) = \lambda\Pi_\vartheta\mathfrak{f}$ as

$$\tfrac{i}{2}\,\partial_t(\Pi_\vartheta\mathfrak{f}) + \mathcal{R}\cdot\Pi_\vartheta\mathfrak{f} + \mathfrak{e}(\mathfrak{f}) = \lambda\Pi_\vartheta\mathfrak{f} \tag{A.35}$$

where $\mathcal{R}$ is an $\mathbb{R}$-linear section of the bundle of endomorphisms of $\mathrm{Ker}_\vartheta|_{\gamma_*}$ and where $\mathfrak{e}$ is an $\mathbb{R}$-linear functional of $\mathfrak{f}$ that has small norm when $\mathfrak{f}$ has $L^2$ norm equal to 1. An equation of this sort appears because the $\Pi_\vartheta$-image of the right hand of (A.29) at any given $t \in \gamma_*$ can be written schematically as $\vartheta_{z_t}{}^\dagger\mathfrak{z} + \mathfrak{r}$ where $\mathfrak{z}$ depends on $(p, \eta_1)$ and $\mathfrak{r}$ is small in a suitable sense. Meanwhile, $\vartheta_{z_t}{}^\dagger\mathfrak{z}$ projects to 0 in $\mathrm{Ker}_\vartheta|_t$ and so the lack of an a priori small bound for the norm of $\vartheta_{z_t}{}^\dagger\mathfrak{z}$ is of no concern.

This subsection use (A.26) with (A.28)-(A.30) to say more about (A.35) when the given element $\gamma_* \in \Theta_*$ is some $\gamma \in \cup_{p\in\Lambda}\{\hat{\gamma}_p^+ \cup \hat{\gamma}_p^-\}$ and $k \in \{0,1,\ldots,7\}$ version of $\gamma_k$. The salient points are summarized by Lemma A.7. The first two parts of this subsection set up the background for Lemma A.7; the third part contains the lemma and its proof.

*Part 1*: To set the notation for what is to come, introduce $\kappa_\Diamond$ and $\kappa_{c\Diamond}$ to denote the larger of the respective versions of $\kappa$ and $\kappa_{c0}$ that are supplied by Lemmas A.2-A.6. Fix $c_0 \geq \kappa_\Diamond$ and $z \geq \kappa_{c\Diamond}\, c_0^{10}$ for use in Section Ab. Let $(A, \psi = (\alpha, \beta))$ denote a pair that satisfies PROPERTIES 1-5 in Section Ab using the given values of $c_0$ and $z$. Define the set $\Theta_*$ as in Section Ae, and focus attention on a given element in $\Theta_*$ that has the form $\gamma_k$ with $\gamma \in \cup_{p\in\Lambda}\{\hat{\gamma}_p^+ \cup \hat{\gamma}_p^-\}$ and integer $k \in \{0, \ldots, 7\}$. If this element is the whole of $\gamma$, write $\gamma$ as the union of two open sets of length $\frac{3}{4}\ell_\gamma$ with distance $\frac{1}{2}\ell_\gamma$ between their respective midpoints. These sets are denoted in what follows by $\gamma_+$ and $\gamma_-$. Introduce $\gamma_*$ to denote $\gamma_+$ or $\gamma_k$ if $\gamma_k = \gamma$ and to denote $\gamma_k$ if $\gamma_k \neq \gamma$.

Let $T \subset Y$ denote the radius $c_0^{-1}$ tubular neighborhood of $\gamma$ with radius chosen so as to use $\gamma$'s version of the coordinates from Part 4 of Section Aa for T with $\nu$ and $\mu$ constant and real, and with $\mu$ greater than $|\nu|$. To spare notation, suppose that the $t = 0$ point is in $\gamma_*$. Use $T_* \subset T$ to denote the set of points with coordinate $t \in \gamma_*$ and $|z| \leq 2c_v^4$. Fix once and for all an isomorphism between E's restriction to the transverse disk in T through the $t = 0$ point and the product bundle over this same disk. Parallel transport along the constant $z \in \mathbb{C}$ slices of $T_*$ from the $t = 0$ transverse disk defines an isomorphism between $E|_{T_*}$ and $T_* \times \mathbb{C}$. This isomorphism writes A on $T_*$ as

$$A = \theta_0 + \tfrac{1}{2}(\mathrm{A}\, d\bar{z} - \bar{\mathrm{A}}\, dz) \tag{A.36}$$

with $\mathrm{A}$ being a $\mathbb{C}$-valued function on $T_*$. This isomorphism makes $\alpha$ a $\mathbb{C}$-valued function.

*Part 2*: Given $t \in \gamma_*$, use $\vartheta_{zt}$ to denote the $z = r$ version of Section Ad's operator $\vartheta_{z*}$ that is defined by the restriction of $(A, \alpha)$ to the transverse disk through t. The isomorphism between $E|_{T_*}$ and $T_* \times \mathbb{C}$ writes the family $\{\vartheta_{zt}\}_{t\in\gamma_*}$ as a smooth, 1-parameter family of operators on $\mathbb{C}$ and it identifes the bundle $\mathrm{Ker}_\vartheta|_{\gamma_*}$ with a subbundle of the product bundle $\gamma_* \times C^\infty(\mathbb{C}; \mathbb{C}\oplus\mathbb{C})$. In particular, this isomorphism writes any given section of $\mathrm{Ker}_\vartheta$ over $\gamma_*$ as a map from $\gamma_*$ to $C^\infty(\mathbb{C}; \mathbb{C}\oplus\mathbb{C})$. Viewed in this way, the $L^2$ orthogonal projection on $L^2(\mathbb{C}; \mathbb{C}\oplus\mathbb{C})$ induces a covariant derivative on sections of $\mathrm{Ker}_\vartheta|_{\gamma_*}$ as follows: Let $t \to \zeta = (x_t, \iota_t)$ denote a smooth map from $\gamma_*$ to $C^\infty(\mathbb{C}; \mathbb{C}\oplus\mathbb{C})$ such that $\zeta$ at each $t \in \gamma_*$ is a square integrable element in the kernel of $\vartheta_{zt}$. The covariant derivative of this section is denoted by $\mathrm{D}\zeta$ and it is defined at any given $t \in \vartheta$ by the rule

$$\mathrm{D}\zeta = \tfrac{d}{dt}\zeta + \vartheta_{zt}^{\dagger}\varpi\,, \tag{A.37}$$

where $\varpi$ is the square integrable solution to the equation

$$\vartheta_{z_t}\vartheta_{z_t}^{\dagger}\varpi + (\tfrac{\partial}{\partial t}\vartheta_{z_t})\zeta = 0 .$$

(A.38)

Note in this regard that $\frac{\partial}{\partial t}\vartheta_{z_t}$ is an endomorphism of the product bundle $\mathbb{C} \times (\mathbb{C}\oplus\mathbb{C})$ whose coefficients are defined by the t-derivative of the function A and the corresponding covariant derivative of $\alpha$ on the union of the radius $(c_v^4 + 10c_v)\,r^{-1/2}$ transverse disks in $T_*$ with centers at the points in $\gamma_*$. This covariant derivative on $\mathrm{Ker}_\vartheta|_{\gamma_*}$ is metric compatible with it understood that the metric on this bundle is that induced by the $L^2$ inner product on the space of square integrable maps from $\mathbb{C}$ to $\mathbb{C}\oplus\mathbb{C}$.

Let m denote the fiber dimension of $\mathrm{Ker}_\vartheta|_{\gamma_*}$, this being the dimension of the $L^2$ kernel of any $t \in \gamma_*$ version of $\vartheta_{z_t}$. Fix an $L^2$ orthonormal basis for $\mathrm{Ker}_\vartheta|_{\gamma_*}$ at $t = 0$. Parallel transport this basis along $\gamma_*$ using the connection defined by (A.37) and (A.38) to define an isomorphism from $\mathrm{Ker}_\vartheta|_{\gamma_*}$ to $\gamma_* \times \mathbb{C}^m$. This isomorphism is used in the upcoming Lemma A.7 to view a section of $\mathrm{Ker}_\vartheta|_{\gamma_*}$ as a map from $\gamma_*$ to $\mathbb{C}^m$.

*Part 3*: The stage is now set for Lemma A.7

**Lemma A.7**: *There exists $\kappa \ge \kappa_\Diamond$ and given $c_0 \ge \kappa$, there exists $\kappa_{c0} \ge \kappa_{c\Diamond}$ with the following significance: Fix $c_0 \ge \kappa$ and $z \ge \kappa_{c0}\, c_0^{10}$ Suppose that $\mathfrak{c} = (A,\psi)$ obeys the corresponding version of* PROPERTIES 1-5 *in Section Ab. Let $\mathfrak{f}$ denote an eigenvector of the operator $\mathfrak{L}_{\mathbb{V}}$ with eigenvalue bounded in absolute value by $c_0^{-\kappa} z^{1/2}$. Define $\gamma_*$ as in Part 1 and view both $\Pi_\vartheta\mathfrak{f}$ and $\Pi_\vartheta(\mathfrak{L}_{\mathbb{V}}\mathfrak{f})$ along $\gamma_*$ as maps from $\gamma_*$ to $\mathbb{C}^m$ as instructed in Part 2. Viewed in this way, the equation $\Pi_\vartheta(\mathfrak{L}_{\mathbb{V}}\mathfrak{f}) = \lambda\Pi_\vartheta\mathfrak{f}$ has the form*

$$\tfrac{i}{2}\tfrac{d}{dt}(\Pi_\vartheta\mathfrak{f}) + \mathfrak{r}(\mathfrak{f}) = \lambda\,\Pi_\vartheta\mathfrak{f}$$

*with the endomorphism $\mathfrak{r}$ being an $\mathbb{R}$-linear functional of $\mathfrak{f}$ that obeys* $\int_{\gamma_*} |\,\mathfrak{r}(\mathfrak{f})\,| \le c_0^{\kappa}\,\|\mathfrak{f}\|_2$.

***Proof of Lemma A.7***: Use $\langle\cdot,\cdot\rangle$ to denote the $\mathrm{Ker}_\vartheta|_{\gamma_*}$ inner product at given $t \in \gamma_*$. With $\mathrm{Ker}_\vartheta|_{\gamma_*}$ viewed as $\gamma_* \times \mathbb{C}^m$, this is just the Hermitian inner product on $\mathbb{C}^m$; and with $\mathrm{Ker}_\vartheta|_{\gamma_*}$ viewed as a subbundle of $\gamma_* \times C^\infty(\mathbb{C};\mathbb{C}\oplus\mathbb{C})$, this same Hermitian inner product is the $L^2$ inner product on the subspace of square integrable maps from $\mathbb{C}$ to $\mathbb{C}\oplus\mathbb{C}$. Let $\zeta$ denote a covariantly constant section of $\mathrm{Ker}_\vartheta|_{\gamma_*}$ with unit $L^2$-norm. Use the two views of $\langle\cdot,\cdot\rangle$ with (A.26), (A.28)-(A.30) to write $\langle\zeta,\Pi_\vartheta\mathfrak{L}_{\mathbb{V}}\mathfrak{f}\rangle$ as $\frac{i}{2}\frac{d}{dt}\langle\zeta,\Pi_\vartheta\mathfrak{f}\rangle + \mathfrak{r}_{\mathbb{V}}(\mathfrak{f}) + \mathfrak{r}_\omega(\mathfrak{f})$ where the function $t \to \mathfrak{r}_{\mathbb{V}}(\mathfrak{f})|_t$ comes from the inner product of $\zeta$ with the $\Pi_\vartheta$ images of all but the term $\frac{\partial}{\partial t}\mathfrak{f}_0 = (\frac{\partial}{\partial t}q, \frac{\partial}{\partial t}\eta_1)$ in (A.28), and with the $\Pi_\vartheta$ images of all of the terms in (A.29) and (A.30). Meanwhile, the function $t \to \mathfrak{r}_\omega(\mathfrak{f})|_t$ is the right most term in the identity

$$\frac{i}{2}\int_{\mathbb{C}}\zeta^\dagger\chi_*\varphi_z^*(\tfrac{\partial}{\partial t}\mathfrak{f}_0) = \frac{i}{2}\frac{\partial}{\partial t}\int_{\mathbb{C}}\zeta^\dagger\chi_*\varphi_z^*(\mathfrak{f}_0) - \frac{i}{2}\int_{\mathbb{C}}(\tfrac{\partial}{\partial t}\zeta)^\dagger\chi_*\varphi_z^*(\mathfrak{f}_0)\,. \tag{A.39}$$

The term $\mathfrak{r}_{\mathbb{V}}(\mathfrak{f})$ is such that

$$\int_{\gamma_*}(\int_{\mathbb{C}}|\mathfrak{r}_{\mathbb{V}}(\mathfrak{f})|^2) \le c_0 c_0^k \|\mathfrak{f}\|_2^2 \tag{A.40}$$

with $k \le c_0$. This can be seen by using the first bullet of Lemma A.2 to bound the contributions from (A.28), by using (A.33) with Lemma A.3 to bound those from (A.29), and by using PROPERTY 1 in Section Aa to bound the contribution to the terms with $\beta$ in (A.30). Note with regards to (A.28) that the function $a_{A0}$ is zero because there is no dt component on the right hand side of (A.36).

To obtain the desired bound on the term $\mathfrak{r}_\omega(\mathfrak{f})$, use (A.37) to rewrite this term as

$$\frac{i}{2}\int_{\mathbb{C}}(\vartheta_{zt}^\dagger\varpi)\chi_*\varphi_z^*(\mathfrak{f}_0)\,. \tag{A.41}$$

Use the Minkowski inequality to bound (A.41) by the product of the $L^2$ norm of $\vartheta_{z_t}{}^\dagger\omega$ on $\mathbb{C}$ and that of $\chi_*\varphi_r^*(\mathfrak{f}_0)$. The latter norm is bounded by a uniform multiple of the $L^2$ norm of $\mathfrak{f}$ over the transverse disk centered at the given point on $\gamma_*$. Use (A.38) with Lemma A.6 to bound the $L^2$ norm on $\mathbb{C}$ of $\vartheta_{z_t}{}^\dagger\omega$ by a multiple of the $L^2$ norm on $\mathbb{C}$ of $|\frac{\partial}{\partial t}\vartheta_{z_t}||\zeta|$. This in turn is bounded by

$$c_0 \sup_{\{(t,z)\in\gamma_*\times\mathbb{C}:\ |z|<c_0^4+10c_0\}} (r^{-1/2}|\tfrac{\partial}{\partial t}A| + |\tfrac{\partial}{\partial t}\alpha|) \tag{A.42}$$

because $\zeta$ has unit $L^2$ norm.

To say something about the size of (A.42), note first that $\frac{\partial}{\partial t}\alpha$ is the covariant derivative of $\alpha$ along the coordinate vector field $\frac{\partial}{\partial t}$ because Part 1's isomorphism between E over $T_*$ writes A as depicted in (A.36). Meanwhile this covariant derivative differs from $(\nabla_A\alpha)_\nu$ by at most $c_0 c_0^4 z^{-1/2}|\nabla_A\alpha|$ because the vector fields $\frac{\partial}{\partial t}$ and $\nu$ differ on $T_*$ by at most $c_0 c^4 z^{-1/2}$. This being the case, use PROPERTIES 1 and 2 in Section Aa to see that the derivative of $\alpha$ in (A.42) is no greater than $c_0 c_0^6$. Meanwhile, $\frac{1}{2}\frac{\partial}{\partial t}A$ is the $dt\, d\bar{z}$ component of the curvature 2-form of A because of the lack of a term in (A.36) that is proportional to dt. Use this fact with the aforementioned bound on $|\frac{\partial}{\partial t} - \nu|$ and PROPERTY 2 in Section Aa to see that term $z^{-1/2}|\frac{\partial}{\partial t} A|$ in (A.42) is likewise no greater than $c_0\, c_0^{\,6}$.

**i) Lemma A.1 and the equation $\Pi_\vartheta \mathfrak{L}_{\mathbb{V}} \mathfrak{f} = \lambda \Pi_\vartheta \mathfrak{f}$**

The lemma that follows talks about (A.35) in the case when $(A, \psi)$ is described by Lemma A.1. This lemma refers to the functions $\nu$ and $\mu$ that appear in a given version of (A.6) for the case when the relevent curve from $\Theta$ lies in $Y_{*\Lambda}$. Choose a coordinate system of the sort described in Part 4 of Section Aa for each such curve from $\Theta$ with a bound by $c_0$ on the corresponding versions of $|\nu|$ and $|\mu|$. Such a choice is assumed implicitly in the lemma.

**Lemma A.8**: *There exists* $\kappa \geq 100$ *and, given* $c_v \geq \kappa$, *there exists* $\kappa_{c_v} \geq \kappa$ *both greater than their incarnations in Lemmas A.1-A.6 and with the following additional property: Fix* $z \geq \kappa_{c_v} c_v^{\,10}$ *and a pair* $(A, \psi) \in \mathrm{Conn}(E) \times C^\infty(Y; \mathbb{S})$ *that is described by Lemma A.1 using the chose value of* $c_v$. *Let* $\gamma$ *denote a curve from* $\Theta$ *in* $Y_{*\Lambda}$. *Suppose that* $\lambda$ *is an eigenvalue of the corresponding version of* $\mathfrak{L}_{\mathbb{V}}$ *with* $|\lambda| \leq c_v^{-\kappa} z^{1/2}$ *and let* $\mathfrak{f}$ *denote the corresponding eigenvector. Use* $\zeta$ *to denote the section of the line bundle* $\mathrm{Ker}_\vartheta|_\gamma \to \gamma$ *given by* $(\Pi_\vartheta \mathfrak{f})|_\gamma$. *There is an isomorphism of* $\mathrm{Ker}_\vartheta|_\gamma$ *with* $\gamma \times \mathbb{C}$ *that writes* $\Pi_\vartheta \mathfrak{f}$ *as a map* $\zeta: \gamma_* \to \mathbb{C}$ *and the* $\Pi_\vartheta$*-image along* $\gamma$ *of the eigenvalue equation* $\mathfrak{L}_{\mathbb{V}} \mathfrak{f} = \lambda \mathfrak{f}$ *as*

$$\frac{i}{2}\frac{d}{dt}\zeta + \nu\zeta + \mu\bar{\zeta} = \lambda\zeta + \mathfrak{e}(\mathfrak{f}) \,,$$

*where* $\mathfrak{e}(\mathfrak{f})$ *is an* $\mathbb{R}$*-linear functional of* $\mathfrak{f}$ *that has* $L^2$ *norm bounded by* $\kappa c_v^{-1} \|\mathfrak{f}\|_2$.

***Proof of Lemma A.8***: Use (A.28)-(A.30) with the conclusions of Lemmas A.2-A.6 as input for the arguments that are used in Steps 9 and 10 from Section 2a in [T4]. These arguments with one addition give a proof. Steps 9 and 10 in Section 2a of [T4] prove the latter's Lemma 2.1 which is the analog of Lemma A.8 for the case with $\hat{a}$ is replaced by a contact 1-form. The one addition concerns the terms in (A.28) that involve $x_\gamma$. To say more about these terms, remark first that they appear only when $\gamma \in Y_{\Diamond z}$. The relevant vortex solution defines the centered solution in $\mathfrak{C}_1$ and if $(A_0, \alpha_0)$ denotes such a solution, then the $L^2$-kernel of the corresponding version of $\vartheta$ is 1-dimensional and spanned by $\frac{1}{\sqrt{\pi}}(\frac{1}{\sqrt{2}}(1 - |\alpha_0|^2, \partial_{A_0}\alpha_0)$. This fact with an integration by parts shows that the $x_\gamma$ terms contribute only to the $\mathfrak{e}(\mathfrak{f})$ term in the statement of the lemma.

The next lemma states a stronger version of what is asserted by Lemma A.7 for cases when $(A, \psi)$ on a given component of $Y-(Y_{*\Lambda} \cup T_{*\Lambda})$ is also given by the constructions in Section Aa using $\rho_* = c_v^4$. To set the stage for the lemma, introduce $\gamma$ to denote the curve from $\cup_{p\in\Lambda}\{\hat{\gamma}_p^+ \cup \hat{\gamma}_p^-\}$ in the given component of $Y-(Y_{*\Lambda} \cup T_{*\Lambda})$. Assume that the curves in $Y_{*\Lambda}$ from $\Theta$ have distance at least $(c_v^4 + 3c_v^3)\, z^{-1/2}$ from $\gamma$. Use T here to denote the set of points with distance less than $(c_v^4 + c_v^3)\, z^{-1/2}$ from $\gamma$. Coordinates for T are given by $\gamma$'s version of the coordinates from Part 4 of Section Aa with $\nu$ and $\mu$ constant and real with $\mu > |\nu|$.

The definition of $(A,\psi)$ on T refers to a function, $\chi_{\diamond\diamond}$, of the radial coordinate $|z|$ on T. This function is non-negative, it is equal to 1 where $|z|$ is less than $(c_v^4 - \frac{7}{4} c_v^2)\, z^{-1/2}$, it is equal to zero where $|z|$ is greater than $(c_v^4 - \frac{5}{4} c_v^2)\, z^{-1/2}$, and the norm of its derivative has absolute value bounded by $32\, c_v^{-3}\, z^{1/2}$. Note in particular that the function $\chi_{\diamond\diamond}$ is zero on $T-(T \cap Y_{\diamond\diamond})$ and it is equal to zero on $T \cap Y_{*\Lambda}$. Such a function can be readily constructed using the function $\chi$.

Let m denote a given positive integer. There is a unique solution to (2.8) with (3.1) equal to m and having the following properties: Write this solution as $(A_{m0}, \alpha_{m0})$. Then $\alpha_{m0} = |\alpha_{m0}|(\frac{z}{|z|})^m$. Meanwhile, $A_{m0}$ can be written in terms of the product connection $\theta_0$ as $A_{m0} = \theta_0 - a_{m0}\, \frac{m}{2}\, (z^{-1}dz - \bar{z}^{-1} d\bar{z})$. Note that both $|a_{m0}|$ and $|\alpha_{m0}|$ are functions only of the radial distance to the origin in $\mathbb{C}$. The m = 1 version of $a_{m0}$ is denoted by $a_0$ in (A.3). Any given $m \geq 1$ version of $a_{m0}$ obeys the analog of the the m = 1 bound in (A.4), this being $|1 - a_{m0}| \leq c_0 (1 - |\alpha_{m0}|)$. The pair $(A_{m0}, \alpha_{m0})$ defines the point in the space $\mathfrak{C}_m$ from Part 1 of Section 3a that maps via the coordinates in (3.2) to the origin in $\mathbb{C}^m$. Let $y_m$ and $\varsigma_m$ denote the $(A_{m0}, \alpha_{m0})$ versions of the functions y and $\varsigma$ that are described in Section Aa.

Fix an isomorphism between $E|_T$ and $T \times \mathbb{C}$ and use this isomorphism to view A as a connection on $T \times \mathbb{C}$ and the component $\alpha$ of $\psi$ as a complex valued function on T. Use this isomorphism with the coordinates from Part 4 of Section Aa to view $\beta$ as a complex valued function also. With this view understood, the connection A is written as $\theta + a_U$ where $a_U$ is an $i\mathbb{R}$ valued 1-form on T. The 1-form $a_U$, $\alpha$ and $\beta$ are defined as follows:

- $a_U = \nu\, \chi_{\diamond\diamond}\, i\, 2^{1/2}\, r_z^{*} y_m\, dt - \frac{m}{2}\, (1 - \chi_{\diamond\diamond} + \chi_{\diamond\diamond}\, r_z^{*} a_{m0})(z^{-1}dz - \bar{z}^{-1} d\bar{z})$,
- $\alpha_U = (1 - \chi_{\diamond\diamond}(1 - r_z^{*}|\alpha_{m0}|)\, (\frac{z}{|z|})^m$,
- $\beta_U = i\, \mu\, z^{-1/2} \chi_{\diamond\diamond}\, r_z^{*} \varsigma_m$ .

(A.44)

Let $\mathfrak{c}_m$: $S^1 \to \mathfrak{C}_m$ denote the constant map to the point given by $(A_{m0}, \alpha_{m0})$. The upcoming lemma also refers to the $(A_{m0}, \alpha_{m0})$ version of the linear operator that is depicted in (3.10).

**Lemma A.9**: *Fix* $m \ge 1$*, there exists* $\kappa \ge 100$ *and, given* $c_0 \ge \kappa$*, there exists* $\kappa_{c0} \ge \kappa$*, both greater than their incarnations in Lemmas A.1-A.6 and with the following significance: Fix* $z \ge \kappa c_0^{10}$*. Set* $c_v = c_0$ *and then set* $\rho_* = c_0^2 z^{-1/2}$*. Fix a set* $T_{*\Lambda}$ *and then a set* $\Theta$ *as described by (A.5) which obeys the first and second bullets of the* $(z, c_0)$ *version of* *PROPERTY 3.*

- *Suppose that* $(A, \psi) \in \text{Conn}(E) \times C^\infty(Y; \mathbb{S})$ *is given by the* $(z, c_0, \rho_* = c_0^2 z^{-1/2})$ *version of (A.7)-(A.10) on* $Y_{*\Lambda} \cup T_{*\Lambda}$*. Fix a component of* $Y - (Y_{*\Lambda} \cup T_{*\Lambda})$ *and suppose that* $(A, \psi)$ *is given by (A.44) on this component. Assume that PROPERTIES 1,2,4 and 5 hold on the rest of* $Y - (Y_{*\Lambda} \cup T_{*\Lambda})$*. Then* $(A, \psi)$ *obeys PROPERTIES 1-5 on the whole of* $Y$*.*
- *Suppose that* $\lambda$ *is an eigenvalue of the corresponding version of* $\mathfrak{L}_{\mathbb{V}}$ *with* $|\lambda| \le c_0^{-\kappa} z^{1/2}$ *and let* $\mathfrak{f}$ *denote the corresponding eigenvector. Let* $\gamma$ *denote the curve from* $\cup_{p \in \Lambda}\{\hat{\gamma}_p^+ \cup \hat{\gamma}_p^-\}$ *in the given component of* $Y - (Y_{*\Lambda} \cup T_{*\Lambda})$*. Use* $\zeta$ *to denote the section of the line bundle* $\text{Ker}_\vartheta|_\gamma \to \gamma$ *given by* $(\Pi_\vartheta \mathfrak{f})|_\gamma$*. There is an isomorphism of* $\text{Ker}_\vartheta|_\gamma$ *with* $\gamma \times \mathbb{C}$ *that writes* $\Pi_\vartheta \mathfrak{f}$ *as a map* $\zeta: \gamma_* \to \mathbb{C}$ *and the* $\Pi_\vartheta$*-image along* $\gamma$ *of the eigenvalue equation* $\mathfrak{L}_{\mathbb{V}} \mathfrak{f} = \lambda \mathfrak{f}$ *as* $\zeta \to \frac{i}{2} \nabla_t \zeta + (\nabla_{\zeta_{\mathbb{R}}} \nabla^{1.0} \hbar)|_{\mathfrak{c}_m} + \mathfrak{e}(\zeta)$ *where* $\mathfrak{e}(\mathfrak{f})$ *is an* $\mathbb{R}$*-linear functional of* $\mathfrak{f}$ *that has* $L^2$ *norm bounded by* $\kappa c_0^{-1} \|\mathfrak{f}\|_2$*.*

***Proof of Lemma A.9***: The proof of the first bullet is a version of what is done in Sections 2e and 2f of the article *Gr => SW* in [T8]. The proof of the second bullet is a version of what is done in Steps 9 and 10 in Section 2a of [T4].

## B. Vortex equation solutions and $(A, \psi)$.

This section of the appendix supplies additional material for the proof of Proposition 2.6. To give a look ahead, suppose that $(A, \psi = (\alpha, \beta))$ is a solution to a given $(r, \mu)$ version of (1.13) with $\mu$ being a 1-form from $\Omega$ whose $\mathcal{P}$-norm is less than 1. The replacement of $(A, \psi)$ with a pair made from vortex solutions facilitates the upcoming analysis of the r-dependence of the spectrum of $\mathfrak{L}_{\mathbb{V}}$. By way of a reminder, the value of $\mathfrak{f}_s$ at $(A, \psi)$ requires comparing the spectrum of the $(z = r, (A, \psi))$ version of $\mathfrak{L}_{\mathbb{V}}$ with that of a version defined using $z = 1$. If the $z = r$ version of the operator $\mathfrak{L}_{\mathbb{V}}$ is defined not by $(A, \psi)$ but by a pair made from vortex solutions, then Lemmas A.8 and A.9 can used to analyze the spectrum of $\mathfrak{L}_{\mathbb{V}}$. The input from these lemmas are used in Appendix C to study the spectral flow for the versions of $\mathfrak{L}_{\mathbb{V}}$ along 1-parameter family that is defined by the value of $z$ and a corresponding $z$-dependent pair in $\mathrm{Conn(E)} \times C^\infty(Y;\mathbb{S})$ that is built from vortex solutions.

The constructions that follow in this appendix use $(A,\psi)$ to construct a new pair in $\mathrm{Conn(E)} \times C^\infty(Y;\mathbb{S})$ that is defined on all of Y using solutions to the vortex equations in (2.8). This new pair is denoted by $(A_\Diamond,\psi_\Diamond)$. The norm of the difference between the values of $\mathfrak{f}_s$ as defined using the $(z = r, (A,\psi))$ version of $\mathfrak{L}_{\mathbb{V}}$ and using the $(z = r, (A_\Diamond, \psi_\Diamond))$ version of $\mathfrak{L}_{\mathbb{V}}$ is shown to be bounded by an $(A,\psi)$ and r independent constant. It proves convenient to construct the desired pair $(A_\Diamond, \psi_\Diamond)$ in two stages. The first stage constructs a pair that is denoted by $(A_*,\psi_*)$. This pair is defined on most, but not all of Y using solutions to the vortex equations in (2.8). In particular, the definition does not use vortex solutions near certain curves from the set $\cup_{p\in\Lambda}\{\hat\gamma_p^+ \cup \hat\gamma_p^-\}$. The second stage modifes $(A_*,\psi_*)$ near these curves to obtain the desired pair $(A_\Diamond, \psi_\Diamond)$.

**a) The construction of $(A_*,\psi_*)$.**

This subsection constructs the desired pair $(A_*,\psi_*)$ from data supplied by the given solution to (1.13). The first four parts of this subsection construct $(A_*,\psi_*)$. The fifth part of the subsection explains why $(A_*,\psi_*)$ does not depend on the coordinates from Part 4 of Section Aa that are chosen in Part 2. The sixth and final part of the subsection constructs a path in $\mathrm{Conn(E)}\times C^\infty(Y;\mathbb{S})$ between $(A_*,\psi_*)$ and the given solution to (1.13).

*Part 1.* The constructions in Section Aa are used to define $(A_*, \psi_*)$ over most of Y. These constructions require as input the specification of parameters $c_v$, $z$, and $\rho_*$. The parameter $c_v$ is chosen in a two step process as follows. A preliminary step chooses a parameter $c_{v1}$ so as to be larger than the various incarnations of the constant $\kappa$ that are given by Proposition 2.4 and Lemmas A.1-A.9. With $c_{v1}$ chosen, let $\kappa_\Diamond$ denote the largest of the various $c_0 \in [c_{v1},2c_{v1}]$ versions of the constant $\kappa_{c0}$ that are given by these same lemmas. Assume that $r > 2^{20}\kappa_\Diamond c_{v1}^{10}$ and suppose that $(A,\psi = (\alpha, \beta))$ is a solution to the $(r,\mu)$ version of (1.13) with $\mu$ a given element in $\Omega$ with $\mathcal{P}$-norm smaller than 1. The second step to choosing $c_v$ depends specifically on the form of the zero locus of $\alpha$ and thus on the chosen solution $(A,\psi)$ of (1.13). The constant $c_v$ should be chosen from the interval $[c_{v1},2c_{v1}]$ so as to satisfy certain conditions that are stated momentarily. These conditions refer to the subset $Y_{\Diamond\Diamond}$ of Y that consists of those points with distance at least $(c_v^4 - 3c_v^3)\, r^{-1/2}$ from each curve in the set $\cup_{p\in\Lambda}\{\hat\gamma_p^+ \cup \hat\gamma_p^-\}$.

- *If a component of $\alpha$'s zero locus in $Y_{\Diamond\Diamond}$ is disjoint from a given boundary torus of $Y_{\Diamond\Diamond}$, then all of its points have distance greater than $6c_v^3 r^{-1/2}$ from this torus.*

- *If a component of* $\alpha$*'s zero locus in* $Y_{\Diamond\Diamond}$ *intersects a boundary torus of* $Y_{\Diamond\Diamond}$ *then this intersection point is an endpoint of the component and it is a transversal intersection. One end point of such a component lies where* $u < 0$ *on some boundary component of* $Y_{\Diamond\Diamond}$ *and the other where* $u > 0$ *on some boundary component of* $Y_{\Diamond\Diamond}$.
- *If a given boundary component of* $Y_{\Diamond\Diamond}$ *intersects the zero locus of* $\alpha$*, then it does so at two points. The distance between these points is at least* $100c_v^2 r^{-1/2}$*; and one lies where* $u < 0$ *and the other where* $u > 0$.

(B.1)

Use Proposition 2.4 with the formula for $\nu$ in (1.3) to see that (B.1) will hold if $c_v$ is chosen from the complement of at most G intervals of length $c_0$ in $[c_{v1}, 2c_{v1}]$. These intervals are determined by $\alpha$ and thus by the chosen $(A, \psi)$.

Take $z = r$ and $\rho_* = c_v^2 r^{-1/2}$ to complete the specification of Section Aa's required parameters.

*Part 2*: Use $Y_{*\Lambda} \subset Y_{\Diamond\Diamond}$ to denote the set of points with distance at least $c_v^4 r^{-1/2}$ from each curve in the set $\cup_{p\in\Lambda}\{\hat{\gamma}_p^+ \cup \hat{\gamma}_p^-\}$. The constructions in Section Aa require as additional input the choice of a union of components of $Y - Y_{*\Lambda}$, this denoted by $T_{*\Lambda}$. Define $T_{*\Lambda}$ as follows: A component of $Y - Y_{*\Lambda}$ is in $T_{*\Lambda}$ if and only if the component lacks zeros of $\alpha$.

Having specified $T_{*\Lambda}$, the next order of business is to specify a set $\Theta$ that consists of embedded 1-manifolds in $Y_{*\Lambda} \cup T_{*\Lambda}$. These are the components of $\alpha^{-1}(0) \cap Y_{*\Lambda}$. In particular, $\Theta$ has no curves from $\cup_{p\in\Lambda}\{\hat{\gamma}_p^+ \cup \hat{\gamma}_p^-\}$. The constraint in (B.1) guarantees that the requirements of bullets two and three of (A.5) and bullets one and two of PROPERTY 3 of Section Ab are met by the curves in $\Theta$. That such is the case can be seen using Proposition 2.4 with the formula for $\nu$ in (1.3).

Having specified $T_{*\Lambda}$ version of $\Theta$, Part 3 of Section Aa introduces a set denoted by $U_0$ and sets $\{U_\gamma\}_{\gamma\in\Theta}$. The collection of $\{U_0\} \cup \{U_\gamma\}_{\gamma\in\Theta}$ is denoted by $\mathfrak{U}$. Keep in mind that the union of the sets from $\mathfrak{U}$ contains $Y_{*\Lambda} \cup T_{*\Lambda}$. The constructions in Section Aa require choosing coordinates of the sort described in Part 4 of Section Aa for each $U_\gamma$. Make such a choice once and for all.

Section Aa also requires isomorphisms between the various $U \in \mathfrak{U}$ versions of $E|_U$ and $U \times \mathbb{C}$. Consider first the case of $U_0$. The chosen lower bound for r implies that $|\alpha|$ is nearly 1 on $U_0$ and in particular $|\alpha| > \frac{3}{4}$. This being the case, there is an isomorphism between $E|_U$ and $U \times \mathbb{C}$ that sends $\alpha$ to the map from $U_0$ to $\mathbb{C}$ given by $|\alpha|$. This is the isomorphism to use for $U_0$. Consider next the case for $U_\gamma$ with $\gamma$ a given curve from $\Theta$. The chosen coordinates for $U_\gamma$ supply an isomorphism from $E|_{U_\gamma}$ to $U_\gamma \times \mathbb{C}$ that makes $\alpha$ appear as the map form $U_\gamma$ to $\mathbb{C}$ given by $|\alpha| \frac{z}{|z|}$. Use this isomorphism for $U_\gamma$.

*Part 3*: Section Aa uses the data supplied by Parts 1 and 2 to construct a pair of connection on E and section of E over $U_0 \cup (\cup_{\gamma\in\Theta} U_\gamma)$. The desired $(A_*, \psi_*)$ is defined so as to equal this pair from Section Aa over $Y_{*\Lambda} \cup T_{*\Lambda}$. This understood, this part of the subsection and Part 4 define $(A_*, \psi_*)$ over the components of $Y-(Y_{*\Lambda} \cup T_{*\Lambda})$.

Reintroduce the set $Y_{\diamond\diamond}$ from Part 1, this being the subset of Y whose points have distance $(c_v{}^4 - 3c_v{}^3)\, r^{-1/2}$ from the curves in $\cup_{p\in\Lambda}\{\hat\gamma_p^+ \cup \hat\gamma_p^-\}$. Fix a component of $Y-(Y_{*\Lambda} \cup T_{*\Lambda})$ and use T to denote the radius $(c_v{}^4 + c_v{}^3)\, r^{-1/2}$ tubular neighorhood of the corresponding curve from the set $\cup_{p\in\Lambda}\{\hat\gamma_p^+ \cup \hat\gamma_p^-\}$. This set T is open and the given component is an open subset of T with compact closure.

The definition to come of $(A_*, \psi_*)$ on T uses the coordinates from Part 4 of Section Aa that are defined by T's curve from $\cup_{p\in\Lambda}\{\hat\gamma_p^+ \cup \hat\gamma_p^-\}$. The definition also refers to the function, $\chi_{\diamond\diamond}$, that was introduced in the discussions just prior to Lemma A.9. By way of a reminder, this is a non-negative function of the radial coordinate |z| on T that is equal to 1 where |z| is less than $(c_v{}^4 - \frac{7}{4} c_v{}^2)\, r^{-1/2}$, and equal to zero where |z| is greater than $(c_v{}^4 - \frac{5}{4} c_v{}^2)\, r^{-1/2}$. The norm of its derivative has absolute value bounded by $32 c_v{}^{-3} r^{1/2}$. This function $\chi_{\diamond\diamond}$ is equal to 1 on $T-(T\cap Y_{\diamond\diamond})$ and it is equal to 0 on $T \cap Y_{*\Lambda}$.

The definition of $(A_*, \psi_*)$ over T when $\alpha$ lacks zeros on $T \cap Y_{\diamond\diamond}$ occupies the remainder of Part 3. To start the definition in this case, define $(A_*, \psi_*)$ over $T-(T\cap Y_{\diamond\diamond})$ to be $(A, \psi)$'s restriction to this same subset of T. To define $(A_*, \psi_*)$ on the rest of T, use the first bullet of (B.1) and Proposition 2.4 so conclude that all points in $T\cap Y_{\diamond\diamond}$ have distance at least $2c_v{}^3 r^{-1/2}$ from a zero of $\alpha$ if $c_v \ge c_0$. This last observation has two immediate and not unrelated consequences, the first being that $T \cap Y_{*\Lambda}$ is contained in Section Aa's open set $U_0$. The second consequence comes via Lemma 2.3 which guarantees that $|\alpha| \ge 1 - e^{-c_v^2}$ on $T\cap Y_{\diamond\diamond}$ if $c_v \ge c_0$. Granted that these facts, Part 2's isomorphism from $E|_{U_0}$ to $U_0 \times \mathbb{C}$ sending $\alpha$ to $|\alpha|$ extends over $T \cap Y_{\diamond\diamond}$ using this same rule to identify $E|_{T\cap Y_{**}}$ with $(T\cap Y_{\diamond\diamond})\times\mathbb{C}$. This isomorphism depicts A on $T\cap Y_{\diamond\diamond}$ as a connection on $(T\cap Y_{\diamond\diamond}) \times \mathbb{C}$; and viewed as such, A can be written as $A = \theta_0 + a_{A,U_0}$ where $\theta_0$ denotes the product connection and where $a_{A,U_0}$ is an $i\mathbb{R}$ valued 1-form on $T\cap Y_{\diamond\diamond}$. Use the isomorphism to write $(\alpha, \beta)$ as $(|\alpha|, \beta_{U_0})$ with $\beta_{U_0}$ being a section over $T\cap Y_{\diamond\diamond}$ of the bundle $K^{-1}$.

Write $\psi_*$ as $(\alpha_*, \beta_*)$ with respect to the $E\oplus EK^{-1}$ splitting of $\mathbb{S}$. Granted this notation, use the isomorphism from the preceding paragraph to define $(A_*, \psi_* = (\alpha_*, \beta_*))$ over $T\cap Y_{\diamond\diamond}$ by declaring

$$A_* = \theta_0 + \chi_{\diamond\diamond}\, a_{A,U_0}\,, \quad \alpha_* = (1-\chi_{\diamond\diamond}) + \chi_{\diamond\diamond}|\alpha|\,, \ \textit{and} \ \ \beta_* = \chi_{\diamond\diamond}\, \beta_{U_0}\,. \tag{B.2}$$

The definition given in (B.2) smoothly extends $(A_*, \psi_*)$ from $U_0$ to $U_0 \cup T$ because the pair $(A_*, \psi_*)$ on $U_0$ is defined in Section Aa using Part 2's isomorphism between $E|_{U_0 \cap T}$ and $(U_0 \cap T) \times \mathbb{C}$ as $(A_* = \theta_0, \psi_* = (1, 0))$.

*Part 4*: This part assumes that $\alpha$ has zeros in $T \cap Y_{\Diamond\Diamond}$. The definition in this case also sets $(A_*, \psi_*)$ equal to $(A, \psi)$ on $T-(T \cap Y_{\Diamond\Diamond})$. Four steps are used to define $(A_*, \psi_*)$ on $T \cap Y_{\Diamond\Diamond}$.

Step 1: To set the stage for the definition on $T \cap Y_{\Diamond\Diamond}$ use (B.1), Proposition 2.4 and the depiction of $v$ in (1.3) to see that $\alpha$'s zero locus in the $Y_{\Diamond\Diamond}$ closure of $T \cap Y_{\Diamond\Diamond}$ consists of two embedded, closed arcs, each with one endpoint on the boundary torus of the closure of T and the other on $Y_{\Diamond\Diamond}$'s boundary torus in T. Moreover,

- *The oriented unit tangent vector to each arc differs from* $v$ *by at most* $c_0 r^{-1/2}$.
- *Each arc has transversal intersections with the level sets of* $|z|$.
- *One arc sits where* $u < 0$ *and the other where* $u > 0$ *and both where* $1 - 3\cos^2\theta > 0$.
- *The distance between any given point in one arc from any given point in the other is at least* $100 c_v^2 r^{-1/2}$.

(B.3)

Each arc from (B.3) extends a curve from the set $\Theta$ into $Y_{*\Lambda} \cup (T \cap Y_{\Diamond\Diamond})$ so as to move a boundary point on T's boundary component of $Y_{*\Lambda}$ to T's boundary component of $Y_{\Diamond\Diamond}$. Let $\gamma$ denote such an extended curve. The open set $U_\gamma$ from Section Aa likewise extends into T with the same definition as the union of the radius $4c_v^2 r^{-1/2}$ transverse disks centered on the extension of $\gamma$. This is an open solid torus with core circle $\gamma$. The open set $U_0$ also extends into $T \cap Y_{\Diamond\Diamond}$ as the complement of the union of the radius $c_v^2 r^{-1/2}$ disks centered on the relevant two arcs from (B.3).

Step 2: Granted what is said in Step 1, then Section Aa's definitions can be used to extend $(A_*, \psi_*)$ into $T \cap Y_{\Diamond\Diamond}$. The extended pair over $T \cap Y_{\Diamond\Diamond}$ is denoted by $(A_{*T}, \psi_{*T})$. By way of a reminder, the extension over the complement of the radius $3c_v^2 r^{-1/2}$ tubular neighborhoods of (B.3)'s arcs is written using the isomorphism of E with the product $\mathbb{C}$ bundle that sends $\alpha$ to $|\alpha|$. Meanwhile, $(A_{*T}, \psi_{*T})$ is written over the radius $4c_v^2 r^{-1/2}$ tubular neighborhood of either of (B.3)'s arcs using the coordinates from Part 4 of Section Aa and the isomorphism of E with the product $\mathbb{C}$ bundle that sends $\alpha$ to $|\alpha| \frac{z}{|z|}$. The respective formula on these sets are given below. These formulae write $\psi_{*T}$ in two component form with respect to the splitting of $\mathbb{S}$ as $E \oplus EK^{-1}$. The formulae use $\theta_0$ to denote the product connection on the product $\mathbb{C}$ bundle.

- $A_{*T} = \theta_0$ *and* $\psi_{*T} = (1, 0)$ ,
- $A_{*T} = \theta_0 + i2^{1/2}\nu r_r^* y\,dt - \frac{1}{2} r_r^* a_0(z^{-1}dz - \bar{z}^{-1}d\bar{z})$ *and* $\psi_{*T} = (r_r^*\alpha_0, i\mu r^{-1/2} r_r^* \varsigma)$ .

(B.4)

By way of comparison, the isomorphisms used in (B.4) writes (A, ψ) over the complement of the radius $3c_v^2 r^{-1/2}$ tubular neighborhoods of (B.3)'s arcs and over the radius $4c_v^2 r^{-1/2}$ tubular neighborhood of either arc as

- $A = \theta_0 + a_{A,U_0}$ *and* $\psi = (|\alpha|, \beta_{U_0})$ ,
- $A = \theta_0 + a_{A,U_\gamma}$ *and* $\psi = (|\alpha|\frac{z}{|z|}, \beta_{U_\gamma})$ ,

(B.5)

where $a_{A,U_0}$ and $a_{A,U_\gamma}$ are $i\mathbb{R}$ valued 1-forms. Keep in mind that

$$a_{A,U_0} = a_{A,U_\gamma} + \tfrac{1}{2}(z^{-1}dz - \bar{z}^{-1}d\bar{z}) \quad \textit{and} \quad \beta_{U_0} = \tfrac{\bar{z}}{|z|}\, a_{A,U_\gamma}$$

(B.6)

on the intersection of the respective domains of definition.

The pair $(A_{*T}, \psi_{*T})$ is not the desired extension of $(A_*, \psi_*)$ because it is observedly not the same as (A, ψ) near the boundary torus in T of $Y_{\Diamond\Diamond}$.

<u>Step 3</u>: Let $\gamma \in \Theta$ denote a component with it understood that γ extends into $Y_{\Diamond\Diamond}$. Let $U_\gamma' \subset U_\gamma$ denote the radius $c_v^2 r^{-1/2}$ tubular neighborhood of γ. This step defines a smooth map $u_\gamma: U_\gamma \to S^1$ so as to define.

$$a_{A,\gamma} = a_{A,U_\gamma} - u_\gamma^{-1}du_\gamma \ , \ \alpha_\gamma = |\alpha|\, u_\gamma \tfrac{z}{|z|} \quad \textit{and} \quad \beta_\gamma = u_\gamma \beta_{U_\gamma} \ .$$

(B.7)

The map $u_\gamma$ is constructed so as to obey $u_\gamma = 1$ on $U_\gamma – U_\gamma'$. This being the case, then the pair $(\theta_0 + a_{A,\gamma}, (\alpha_\gamma, \beta_\gamma))$ is gauge equivalent to (A, ψ) on $U_\gamma$ and it extends as (A, ψ) to the whole of $U_0$.

The definition of $u_\gamma$ requires the introduction of a function $\chi_\gamma$ which is given on the $Y_{*\Lambda}$ part of $U_\gamma$ by the rule $z \to \chi(2c_v^{-2} r^{1/2}|z| - 1)$. The function $\chi_\gamma$ on the $Y_{\Diamond\Diamond} – Y_{*\Lambda}$ part of $U_\gamma$ is the product of the function $z \to \chi(2c_v^{-2} r^{1/2}|z| - 1)$ with a second non-negative function. The latter is also constructed using χ and it has the following properties: It is a function of the distance to the nearby component of $\cup_{p\in\Lambda}\{\hat{\gamma}_p^+ \cup \hat{\gamma}_p^-\}$. This second function equals 1 where the distance to these curves is greater than $c_v^4 - \frac{9}{4}c_v^2$ and it equals 0 where the distance is less than $c_v^4 - \frac{11}{4}c_v^2$. The derivative of this second function should have absolute value no greater than $100c_v^{-2} r^{1/2}$. Note in particular that this definition of $\chi_\gamma$ makes it zero on $U_\gamma$'s intersection with a neighborhood of the boundary of $Y_{\Diamond\Diamond}$.

With $\chi_\gamma$ in hand write, $a_{A,U_\gamma}$ as $a_{A,U_\gamma} = a_{A0}\,dt + \frac{1}{2}(A\,d\bar{z} - \bar{A}\,dz)$ with A being a $\mathbb{C}$-valued function on $\mathbb{C}$ and $a_{A0}$ being an $i\mathbb{R}$ valued function on $\mathbb{C}$. The map $u_\gamma$ is defined by the rule

$$u_\gamma = e^{\hat{o}_\gamma} \quad \textit{where} \quad \hat{o}_\gamma(t,z) = \chi_\gamma \tfrac{1}{2} \int_0^1 (\bar{z}A - z\bar{A})|_{(t,sz)}\,ds \ . \tag{B.8}$$

By way of explanation, the map $u_\gamma$ is designed in part so that the 1-form $a_\gamma$ annihilates the radial vector field $z\frac{\partial}{\partial z} + \bar{z}\frac{\partial}{\partial \bar{z}}$ where $\chi_\gamma = 1$. Note in addition that $u_\gamma$ extends to the whole of Y as a smooth map to $S^1$ that is equal to 1 on the complement of a compact set in $U_\gamma$.

<u>Step 4</u>: The desired pair $(A_*, \psi_*)$ is written below using the isomorphisms of E with the product bundle that are used in (B.4) and (B.5). The formula over the complement of the radius $3c_v^2 r^{-1/2}$ tubular neighborhoods of (B.3)'s arcs is

$$A_* = \theta_0 + \chi_{\diamond\diamond}\, a_{A,U_0} \quad \textit{and} \quad \psi_* = ((1 - \chi_{\diamond\diamond}) + \chi_{\diamond\diamond}|\alpha|, \chi_{\diamond\diamond}\, \beta_{U_0}) \ . \tag{B.9}$$

Meanwhile, $A_*$ and $\psi_*$ are written over the radius $4c_v^2 r^{-1/2}$ tubular neighborhood of the extension to $Y_{\diamond\diamond} \cap T$ of an arc $\gamma \in \Theta$ as $A_* = \theta_0 + a_*$ and $\psi_* = (\alpha_*, \beta_*)$ where

- $a_* = \chi_{\diamond\diamond}\, a_{A,\gamma} + (1 - \chi_{\diamond\diamond})\,[\chi_{\hat{0}}\, i\, 2^{1/2} \nu\, r_r^* y\, dt - \frac{1}{2}(1 - \chi_{\hat{0}} + \chi_{\hat{0}} r_r^* a_0)(z^{-1}dz - \bar{z}^{-1}d\bar{z})]$
- $\alpha_* = ((1 - \chi_{\diamond\diamond})(1 - \chi_{\hat{0}}(1 - r_r^*|\alpha_0|)\frac{z}{|z|} + \chi_{\diamond\diamond}\alpha_\gamma)$
- $\beta_* = (1 - \chi_{\diamond\diamond})(i\mu r^{-1/2}\chi_{\hat{0}}\, r_r^* \varsigma) + \chi_{\diamond\diamond}\,\beta_\gamma$.

$$\tag{B.10}$$

The formulae in (B.6)-(B.8) guarantee that (B.9) and (B.10) define a smooth connection on E and section of $\mathbb{S}$ over $Y \cap T$ because $\frac{z}{|z|}$ is the transition function between the relevant product $\mathbb{C}$ bundles.

*Part 5*: This part of the subsection explains why $(A_*, \psi_*)$ does not depend on the choices made in Part 2 of coordinate charts from Part 4 of Section Aa. What follows is the short explanation: A change in the coordinate chart for any given $\gamma \in \Theta$ also changes the product structure for the bundle E over the corresponding set $U_\gamma$. The change in the product structure must be taken into account when comparing versions of $(A_*, \psi_*)$ that are defined by two different choices from Part 4 of Section Aa. The changed product structure compensates for the apparent coordinate dependence in the formula for $(A_*, \psi_*)$. The next two paragraphs say somewhat more about how this comes about.

Recall that a change in the coordinate chart writes the coordinate z on $U_\gamma$ as $u(t)\,z'$ with u being a smooth map from $\gamma$ to $S^1$. To see the effect, consider first the formula in the second bullet of (B.4) that depicts $(A_*, \psi_*)$ on $U_\gamma' \cap Y_{*\Lambda}$. Write the pull-back of the expressions on the right hand side of the equations in the lower bullet of (B.4) via the map $(t, z') \to (t, z = u(t) z')$ in terms of $\nu' = \nu + \frac{i}{2} u^{-1} \frac{d}{dt} u$ and $\mu' = u^{-2}\mu$. Use (A.2) to write $y = -2^{-1/2}(1 - a_0)$ and use the formula for $\varsigma$ in (A.2) to see that the $(t, z') \to (t, u(t) z')$ pull-back of the expression for $A_{*T}$ in the lower bullet of (B.4) is obtained from the $(z', \nu', \mu')$ version of the expression by subtracting $(u^{-1} \frac{d}{dt} u)\, dt$. Meanwhile, the pull-back of the formula for $\psi_{*T}$ in the lower bullet of (B.4) is obtained from the $(z', \nu', \mu')$ version by multiplying the latter by u. These changes are precisely off-set by the change in the product structure.

The invariance of (B.10) with respect to coordinate change can be seen by writing $a_{A,U_\gamma}$ as $\frac{1}{2}(\alpha^{-1}\nabla_A\alpha - \bar{\alpha}^{-1}\nabla_A\bar{\alpha}) - \frac{1}{2}(z^{-1}dz - \bar{z}^{-1}d\bar{z})$ so as to compare $a_{A,U_\gamma}$ with its pull-back via the map $(t, z') \to (t, u(t) z')$.

*Part 6*: Part 3 defined various $\gamma \in \Theta$ versions of a map $u_\gamma$ to $S^1$ from the corresponding $Y_{\Diamond\Diamond}$ extension of $U_\gamma$. As noted at the end of Part 3, such a map extends to the whole of Y as a smooth map that sends the complement of $U_\gamma$ to 1. Let u denote the product of these extended maps, this a smooth map from Y to $S^1$. This part of the subsection describes a path in $\mathrm{Conn}(E) \times C^\infty(Y; \mathbb{S})$ between $(A - u^{-1}du, u\psi)$ and $(A_*, \psi_*)$. This path is parametrized by $\tau \in [0, 1]$ with the $\tau = 0$ member being $(A - u^{-1}du, u\psi)$ and the $\tau = 1$ member being $(A_*, \psi_*)$. The $\tau \in [0, 1]$ member of this path is denoted in what follows by $(A_{*\tau}, \psi_{*\tau})$ and $\psi_{*\tau}$ is written as $(\alpha_{*\tau}, \beta_{*\tau})$ with respect to the splitting of $\mathbb{S}$ as $E \oplus E^{-1}$. The pair $(A_{*\tau}, \psi_{*\tau})$ on $Y - Y_{\Diamond\Diamond}$ is defined to be $(A, \psi)$. The pair $(A_{*\tau}, \psi_{*\tau})$ on the $Y_{*\Lambda}$ part of $U_0 - (\cup_{\gamma \in \Theta} U_\gamma)$ is defined using the $\alpha \to |\alpha|$ isomorphism from $E|_{U_0}$ to $U_0 \times \mathbb{C}$ by the rules

$$A_{*\tau} = \theta_0 + (1 - \tau)\, a_{A,U_0}, \quad \alpha_{*\tau} = \tau + (1 - \tau)|\alpha| \quad \textit{and} \quad \beta_{*\tau} = (1 - \tau)\, \beta_{U_0}\ . \tag{B.11}$$

Meanwhile, the definition on any given $Y_{\Diamond\Diamond} \cap T$ part of $U_0 - (\cup_{\gamma \in \Theta} U_\gamma)$ is obtained from the formula in (B.9) by replacing $\chi_{\Diamond\Diamond}$ with $(1 - \tau) + \tau\chi_{\Diamond\Diamond}$. The pair $(A_{*\tau}, \psi_{*\tau})$ on the $Y_{*\Lambda}$ part of any given $\gamma \in \Theta$ version of $U_\gamma$ is defined using the $\alpha \to |\alpha| \frac{z}{|z|}$ isomorphism from $E|_{U_\gamma}$ to $U_\gamma \times \mathbb{C}$ by the rules

- $A_{*\tau} = \theta_0 + \tau(i 2^{1/2} \nu\, r_r^* y\, dt - \frac{1}{2} r_r^* a_0 (z^{-1}dz - \bar{z}^{-1} d\bar{z})) + (1 - \tau)\, a_{A\gamma}$
- $\alpha_{*\tau} = \tau\, r_r^* \alpha_0 + (1 - \tau)\, \alpha_\gamma$, *and* $\beta_{*\tau} = (1 - \tau)\, \beta_\gamma + \tau\, i \mu r^{-1/2} r_r^* \varsigma)$

(B.12)

The definition over any given $Y_{\Diamond\Diamond} \cap T$ part of $U_\gamma$ is obtained from the formula in (B.10) by repacing $\chi_{**}$ with $(1 - \tau) + \tau\chi_{\Diamond\Diamond}$ and replacing $\chi_{\hat{U}}$ by $\tau\chi_{\hat{U}}$.

By way of a parenthetical remark, this path in $\mathrm{Conn}(E) \times C^{\infty}(Y;\mathbb{S})$ does not depend on the chosen coordinates from Part 4 of Section Aa.

**b) $(A_*, \psi_*)$ and PROPERTIES 1-5**

The upcoming Lemma B.1 asserts that $(A_*, \psi_*)$ and each $\tau \in [0,1]$ member of the path $\tau \to (A_{*\tau}, \psi_{*\tau})$ have all five of the properties that are listed in Section Ab. This lemma is proved using the apriori bounds on the various components of $(A, \psi)$ and $(A_*, \psi_*)$ that are supplied by Lemma B.2.

**Lemma B.1**: *There exists* $\kappa > 1$ *and given* $c_v \geq \kappa$, *there exists* $\kappa_{c_v} \geq \kappa$ *with the following significance: Suppose that* $r \geq \kappa_{c_v} c_v^{10}$ *and suppose that* $(A, \psi = (\alpha, \beta))$ *is a solution to the* $(r, \mu)$ *version of (1.13) with μ a given element in Ω with $\mathcal{P}$-norm smaller than 1. Then the corresponding* $(A_*, \psi_*)$ *satisfies the* $c_0 = c_v$ *and* $z = r$ *version of* PROPERTIES 1-5 *in Section Ab as do all* $\tau \in [0,1]$ *members of the path* $\tau \to (A_{*\tau}, \psi_{*\tau})$.

The proof of this lemma is given momentarily.

Lemma B.2 talks about various components of $(A, \psi)$ on the $Y_{\Diamond\Diamond}$ extensions of $U_0$ and the various $\gamma \in \Theta$ versions of $U_\gamma$. To set the stage for this lemma, use the $\alpha \to |\alpha|$ isomorphism between $E|_{U_0}$ and $U_0 \times \mathbb{C}$ to write $(A, (\alpha, \beta))$ on $U_0$ as $(\theta + a_{A,U_0}, (|\alpha|, \beta_{U_0}))$. The lemma also uses $(a_{A,U_0})_v$ to denote the pairing of the 1-form $a_{A,U_0}$ with $v$.

With $\gamma \in \Theta$ fixed, Lemma B.2 uses the coordinates from Part 4 of Section Aa for $U_\gamma$. Lemma B.2 uses the coordinates from Part 4 of Section Aa, the map $u_\gamma$ from (B.8) and the $\alpha \to |\alpha|\frac{z}{|z|}$ isomorphism between $E|_{U_\gamma}$ and $U_\gamma \times \mathbb{C}$ to write $(A - u_\gamma^{-1} du_\gamma, (u_\gamma \alpha, u_\gamma \beta))$ as $(\theta + a_{A,\gamma}, (\alpha_\gamma, \beta_\gamma))$; and it writes $a_{A,\gamma}$ as $a_{A0,\gamma} dt + \frac{1}{2}(A_\gamma d\bar{z} - \bar{A}_\gamma dz)$. Lemma B.2 also borrows the functions $a_0$ and $\alpha_0$ from (A.3).

**Lemma B.2**: *Fix* $m \geq 1$. *There exists an* m-*dependent* $\kappa > 1$ *and given* $c_v \geq \kappa$, *there exists* $\kappa_{c_v} \geq \kappa$ *with the following significance: Take* $r \geq \kappa_{c_v} c_v^{10}$ *and suppose that* $(A, \psi = (\alpha, \beta))$ *is a solution to the* $(r, \mu)$ *version of (1.13) with μ a given element in Ω with $\mathcal{P}$-norm smaller than 1. Define* $(A_*, \psi_*)$ *as instructed in Section Ba. Then*

- $r^{-1/2}|a_{A,U_0}| + |(a_{A,U_0})_v| + |1 - |\alpha|| + r^{1/2}|\beta_{U_0}| \leq c_v^{-m}$ *on the* $Y_{\Diamond\Diamond}$ *extension of* $U_0$.
- $r^{-1/2}|A_\gamma - r_r^* a_0| + |\alpha_\gamma - r_r^* \alpha_0| \leq c_v^{-m}$ *and* $|a_{A0,\gamma}| \leq \kappa c_v^2$ *on the part of the* $Y_{\Diamond\Diamond}$ *extension of any given* $U_\gamma$ *where the distance to* $\cup_{p \in \Lambda}\{\hat{\gamma}_p^+ \cup \hat{\gamma}_p^-\}$ *is greater than* $(c_v^4 - 2c_v^2) r^{-1/2}$.

The proof of Lemma B.1 assumes that Lemmas B.2 is true.

***Proof of Lemma B.1***: The two steps that follow verify the five properties. These steps use $\kappa_c$ to denote a constant whose value is greater than 1 and depends only on an upper bound for Lemma B.2's constant m and $c_v$, but not on the particulars of $(A, \psi)$ nor on r. This constant can be assumed to increase between subsequent appearances.

Step 1: Given the definition of $\Theta$, what is said in Proposition 2.4 implies PROPERTY 3. The others PROPERTIES hold at where the distance to $\cup_{p\in\Lambda}\{\hat{\gamma}_p^+ \cup \hat{\gamma}_p^-\}$ less than $c_v^4 - 2c_v^2$ if they hold for $(A,\psi)$ which is the case when $c_v \geq c_0$ and $r \geq \kappa_c$.

The remainder of this step verifies PROPERTIES 1, 2, 4 and 5 on $U_0-(\cup_{\gamma\in\Theta}U_\gamma)$. Consider first the $Y_{*\Lambda}$ part of this set. The inequality asserted by the first bullet of PROPERTY 1 and by the first two bullets of PROPERTY 2 follow directly from Lemmas 2.1 and 2.3. Lemma 2.3 also leads directly to PROPERTY 4 and Lemma 2.9 to PROPERTY 5. To verify the remaining parts of PROPERTIES 1 and 2, take $m \geq 100$ and use Lemma B.2's bound $|a_{A,U_0}| \leq c_v^{-m} r^{1/2}$ and the bound $|(a_{A,U_0})_\nu| \leq c^{-m}$ with Lemmas 2.1 and 2.3 to see that $r^{-1/2}|\nabla_{A_{*\tau}}\alpha_{*\tau}|$ and both $|(\nabla_{A_{*\tau}}\alpha_{*\tau})_\nu|$ and $|\nabla_{A_{*\tau}}\beta_{*\tau}|$ are bounded by $c_0$ on the $Y_{*\Lambda}$ part of $U_0-(\cup_{\gamma\in\Theta}U_\gamma)$. The latter set of bounds lead directly to the bounds on the $Y_{*\Lambda}$ part of $U_0-(\cup_{\gamma\in\Theta}U_\gamma)$ that are stated by the second and third bullets of PROPERTY 1 and by the third and fourth bullets of PROPERTY 2 on $Y_{*\Lambda}$.

Given (B.9) and its $(A_{*\tau},\psi_{*\tau})$ analog, the arguments from the preceding paragraph with but one additional comment establish PROPERTIES 1, 2, 4 and 5 on the $Y_{\Diamond\Diamond}-Y_{*\Lambda}$ part of $U_0-(\cup_{\gamma\in\Theta}U_\gamma)$. The additional comment concerns the function $\chi_{\Diamond\Diamond}$ that appears in (B.9), this being the fact that $|(\nabla\chi_{\Diamond\Diamond})_\nu| \leq c_0 c_v^6$ on this part of $Y_{\Diamond\Diamond}$. Indeed, such a bound follows because $\chi_{\Diamond\Diamond}$ is independent of t and because (1.3) finds $|\nu - \frac{\partial}{\partial t}| \leq c_0 c_v^4 r^{-1/2}$ at all points with distance $c_v^4 r^{-1/2}$ or less from $\cup_{p\in\Lambda}(\hat{\gamma}_p^+ \cup \hat{\gamma}_p^-)$.

Step 2: This step verifies PROPERTIES 1, 2, 4 and 5 on any $\gamma \in \Theta$ version of $U_\gamma$. To start, take $m \geq 100$ in the second bullet of Lemma B.2. Its assertions about $A_\gamma$ and $\alpha_\gamma$ imply directly the first bullet of PROPERTY 1 and PROPERTIES 4 and 5. Introduce $\nabla_{A_{*\tau}}^{\perp}$ to denote the $A_{*\tau}$ covariant derivative along the constant t slices of $U_\gamma$. This same part of Lemma B.2 implies that

$$|\nabla_{A_{*\tau}}^{\perp}(\alpha_\gamma - r_r^*\alpha_0)| \leq c_0 c_v^{-50} r^{1/2}(1-|\alpha_\gamma|^2)^{1/2} + c_0 . \tag{B.13}$$

Lemma 2.1 and Lemma B.2's assertions about $A_\gamma$ give the bound $|\nabla_{A_{*\tau}}^{\perp}\beta_\gamma| \leq c_0$ and they imply that (A.2)'s function $\varsigma$ is such that $r^{-1/2}|\nabla_{A_{*\tau}}^{\perp}(r_r^*\varsigma)| \leq c_0$. These bounds with (B.13) imply part of what is required by the second and third bullets of PROPERTY 1 and part of what is required by PROPERTY 2. The remaining parts of PROPERTY 1 and PROPERTY 2 follow directly from Lemma B.2's bound on $|a_{0A,\gamma}|$ given that the absolute value of the directional derivative of $\chi_{\Diamond\Diamond}$ along $\nu$ is bounded by $c_0 c_v^6$ and that of $\chi_{\hat{U}}$ is bounded by $c_0 c_v^2$; the latter being a consequence of what is asserted in the first bullet of (B.3).

***Proof of Lemma B.2***: The proof has four steps. The proof also uses $\kappa_c$ to denote a constant greater than 1 that depends only on a given upper bound for m and $c_v$.

Step 1: This step proves the first bullet of Lemma B.2. The bounds for $1 - |\alpha|$ and for $r^{1/2}\beta_{U_0}$ come directly from Lemmas 2.1 and 2.3 as they bound both by $e^{-c_v}$ if $c_v \geq c_0$ and $r \geq \kappa_c$ To obtain the other bounds, write $\nabla_A\alpha$ on $U_0$ as $d|\alpha| + a_{A,U_0}|\alpha|$. Given that $a_{A,U_0}$ is $i\mathbb{R}$ valued, Lemmas 2.1 and 2.3 imply that $|a_{A,U_0}| \leq e^{-c_v} r^{1/2}$ if $c_v \geq c_0$ and $r \geq \kappa_c$. Meanwhile, these same lemmas together with the vanishing of the $EK^{-1}$ component of $D_A\psi$ imply the bound $|(a_{A,U_0})_\nu| \leq e^{-c_v}$ These bounds lead directly to what is asserted by the first bullet.

Step 2: To start the proof of the second bullet of Lemma B.2, use Proposition 2.4 and Lemma 2.9 to see that

$$\big|\,|\alpha|/|r_r^*\alpha_0| - 1\big| \leq e^{-c_v} \tag{B.14}$$

on $U_\gamma$ when $c_v \geq c_0$ and $r \geq \kappa_c$. Keeping this in mind, use the isomorphism between E over $U_\gamma$ and the product bundle that sends $\alpha$ to $|\alpha|\frac{z}{|z|}$ to write $\alpha$ over $U_\gamma$ as $|\alpha|\frac{z}{|z|}$. Having done so, (B.14) implies that $|\alpha - r_r^*\alpha_0|$ is bounded by $e^{-c_v}$ on $U_\gamma$. Use the same isomorphism to write A over $U_\gamma$ as $A = \theta_0 + a_{A,U_\gamma}$; and use the coordinates from Part 4 of Section Aa on $U_\gamma$ to again write $a_{A,U_\gamma}$ as $a_{A,U_\gamma} = a_{A0}\,dt + \frac{1}{2}(A\,d\bar{z} - \bar{A}\,dz)$. The bound in (B.14) together with Lemma 2.9 have the following additional consequence: For any given t, the $(r_r)^{-1}$ pull-back of $(A, \alpha)|_t$ to the radius $4c_v^2$ ball about the origin in $\mathbb{C}$ differs from $(a_0|z|^2, \alpha_0)$ in the $C^6$ topology by less than $c_v^4 e^{-c_v}$ if $c_v \geq c_0$ and $r \geq \kappa_c$.

The function $a_0$ is a function of $|z|^2$ because $(A_0,\alpha_0)$ gives the symmetric vortex in $\mathfrak{C}_1$ and $\alpha_0 = |\alpha_0|\frac{z}{|z|}$. It follows as a consequence that $\bar{z}a_0 - z\bar{a}_0 = 0$. This understood, the fact that $|z|^2 a_0 - (r_r^{-1})^* A$ has small $C^6$ norm implies that the $r_r^{-1}$ pull-back of the function $\hat{o}_\gamma$ in (B.8) from any given constant t slice of $U_\gamma$ has $C^6$ norm bounded by $c_0 c_v^4 e^{-c_v}$ on the ball of radius $4c_v^2$ centered at the origin in $\mathbb{C}$.

The preceding observation about $\hat{o}_\gamma$ has the following consequence: Write the 1-from $a_{A,\gamma}$ now as $a_{A0,\gamma}dt + \frac{1}{2}(A_\gamma d\bar{z} - \bar{A}_\gamma dz)$. For any given t, the $(r_r)^{-1}$ pull-back of the pair $(A_\gamma, \alpha_\gamma)|_t$ to the radius $4c_v^2$ ball about the origin in $\mathbb{C}$ differs from $(a_0|z|^2, \alpha_0)$ in the $C^6$ topology by less than $c_v^6 e^{-c_v}$ if $c_v \geq c_0$ and $r \geq \kappa_c$. These bounds lead directly to Lemma B.2's assertions about $A_\gamma$ and $\alpha_\gamma$.

Step 3: The bound on $a_{A0,\gamma}$ requires first a bound on $a_{A0,\gamma}$ on the $|z| \geq r^{-1/2}$ part of $U_\gamma$. By way of a parenthetical remark, $(\theta_0 + a_{A,\gamma}, (\alpha_\gamma, \beta_\gamma))$ are used in (B.10) and (B.12) in lieu of $(\theta_0 + a_{A,U_\gamma}, (|\alpha|, \beta_{U_\gamma}))$ in part because no bound of the form $|a_{A0}| \leq c_0 c_v^2$ has been found for the whole of $U_\gamma$. As explained below, a bound of this sort does exist on the complement of any given radius tubular neighborhood of $\gamma$ and the latter bound is needed to derive the desired bound for $|a_{A0,\gamma}|$.

To bound $|a_{A0}|$, use the depiction of A on $U_\gamma$ as $A = \theta_0 + a_{A0}dt + \frac{1}{2}(A d\bar{z} - \bar{A} dz)$ and that of $\alpha$ as $|\alpha| \frac{z}{|z|}$ to write the A-covariant derivative of $\alpha$ along $\frac{\partial}{\partial t}$ as

$$(\nabla_A\alpha)_{\partial/\partial t} = \left(\frac{\partial}{\partial t}|\alpha| + a_{A0}|\alpha|\right)\frac{z}{|z|}. \tag{B.15}$$

Since $a_{A0}$ is $i\mathbb{R}$ valued, the norm of this directional covariant derivative is greater that $|a_{A0}|\,|\alpha|$. Meanwhile, $|v - \frac{\partial}{\partial t}| \leq c_0 c_v r^{-1/2}$ on $U_\gamma$; and as $|(\nabla_A\alpha)_v| \leq c_0$ it follows from the bound in Lemma 2.1 that $|a_{A0}|\,\alpha \leq c_0 c_v$. Thus, $|a_{A0}| \leq c_0 c_v |z|^{-1}$ at any $z \neq 0$ point on $U_\gamma$.

Use $d^\perp a_{A0}$ to denote the differential of $a_{A0}$ along the constant t slices of $U_\gamma$. The identity in (B.15) is used next to obtain a bound by $c_0 c_v (r^{1/2} + |z|^{-1})|z|^{-1}$ on $|d^\perp a_{A0}|$. To get this bound, first write $(\nabla_A\alpha)_{\partial/\partial t}$ as $(\nabla_A\alpha)_v + \mathfrak{R}\cdot\nabla_A\alpha$ where $\mathfrak{R}$ is an endomorphism with norm bounded by $c_0 r^{-1/2}$ and with derivative bounded by $c_0$. Use the $EK^{-1}$ component of the equation $D_A\psi = 0$ to write $(\nabla_A\alpha)_v$ as a linear combination of covariant derivates of $\beta$. Meanwhile, (B.15) finds $a_{A0} = \mathrm{im}(\alpha^{-1}(\nabla_A\alpha)_{\partial/\partial t})$ and so

$$|d^\perp a_{A0}| \leq c_0|\alpha|^{-1}\left(|\alpha|^{-1}|\nabla_A\alpha|\,|(\nabla_A\alpha)_{\partial/\partial t}| + |\nabla\mathfrak{R}|\,|\nabla_A\alpha| + |\mathfrak{R}||\nabla_A{}^2\alpha| + |\nabla_A{}^2\beta|\right). \tag{B.16}$$

The desired bound for $|d^\perp a_{A0}|$ follows from (B.16) and Lemma 2.1.

Step 4: The bounds for $|a_{A0}|$ and $|d^\perp a_{A0}|$ in Step 3 are used first to bound $|a_{A0,\gamma}|$ on the $|z| \geq r^{-1/2}$ part of $U_\gamma$ under the henceforth implicit assumption that the distance to $\cup_{p\in\Lambda}\{\hat{\gamma}_p^+ \cup \hat{\gamma}_p^-\}$ is greater than $c_v^4 - 2c_v^2$. To this end, note first that $|a_{A0,\gamma}| \leq |a_{A0}| + |\partial_t\hat{o}_\gamma|$ and so what is needed is a suitable bound on $|\partial_t\hat{o}_\gamma|$. To obtain one, use (B.8) to see that $|\partial_t\hat{o}_\gamma| \leq c_0|z|\,|\partial_t A|$. Meanwhile, $|\partial_t A| \leq c_0(|d^\perp a_{A0}| + |F_A(\frac{\partial}{\partial t}, \cdot)|$ where $F_A$ denotes the curvature

2-form of A. Use the top bullet in (1.13) with the fact that $\frac{\partial}{\partial t}$ is very close to $\nu$ to see that $|F_A(\frac{\partial}{\partial t}, \cdot)|$ is bounded by $c_0 r^{1/2} c_v$ on $U_\gamma$. What with Step 3's bound for $|d^\perp a_{A0}|$, the latter bound implies that $|\partial_t \hat{o}_\gamma| \le c_0 c_v^2$ on the $|z| \ge r^{-1/2}$ part of $U_\gamma$. This with Step 3's bound for $|a_{A0}|$ leads directly to the desired $|a_{A0,\gamma}| \le c_0 c_v^2$ bound on the $|z| \ge r^{-1/2}$ part of $U_\gamma$.

To obtain the desired bound for $|a_{A0,\gamma}|$ on the $|z| \le r^{-1/2}$ part of $U_\gamma$, fix z with $|z| = r^{-1/2}$ and for any given $\rho \in [0, 4c_v^2 r^{-1/2}]$, write

$$a_{A0,\gamma}|_{\rho z} = a_{A0,\gamma}|_z - \int_{[\rho, r^{-1/2}]} \partial_s (a_{A0,\gamma}|_{sz})\, ds \ . \tag{B.17}$$

Meanwhile, the function $\hat{o}_\gamma$ was chosen specifically so as to guarantee that the 1-form $\frac{1}{2}(A_\gamma d\bar{z} - \bar{A}_\gamma dz)$ annihilate the radial vector field on $\mathbb{C}$, and this implies the identity

$$\partial_s(a_{A0,\gamma}|_{sz}) = -F_A(\tfrac{\partial}{\partial t}, \tfrac{\partial}{\partial |z|})|_{sz} \ . \tag{B.18}$$

As noted in the preceding paragraph, the norm of the right hand side of (B.18) is bounded by $c_0 c_v r^{1/2}$. Use this bound for $|\partial_s(a_{A0,\gamma}|_{sz/|z|})|$ in (B.17) to obtain bound for $|a_{A0,\gamma}|$ by $c_0 c_v^2$ on the $|z| \le r^{-1/2}$ part of $U_\gamma$ when $c_v \ge c_0$ and $r \ge \kappa_c$.

### c) The difference between $\mathfrak{f}_s$ at $(A, \psi)$ and at $(A_*, \psi_*)$

Both the $(A,\psi)$ and the $(A_*,\psi_*)$ version of $\mathfrak{L}_{(\cdot),r}$ might have non-trivial kernel. What follows first defines what is meant by the norm of the spectral flow difference if this is the case. The subsequent Proposition B.3 asserts that this difference is bounded by a purely $c_v$ dependent constant.

Let $\mathfrak{c}_0$ and $\mathfrak{c}_1$ denote a given pair in $\mathrm{Conn}(E) \times C^\infty(Y; \mathbb{S})$. Fix $z_0, z_1 \ge 1$ and introduce $\mathfrak{L}_0$ and $\mathfrak{L}_1$ to denote the respective $(\mathfrak{c}_0, z_0)$ and $(\mathfrak{c}_1, z_1)$ versions of $\mathfrak{L}_{(\cdot)}$. The norm of the spectral flow between $\mathfrak{L}_0$ and $\mathfrak{L}_1$ is denoted here by $|\mathfrak{f}_{s1} - \mathfrak{f}_{s0}|$ and it is defined as follows: Fix $\varepsilon > 0$ and introduce $\mathcal{C}_{0\varepsilon} \subset (\mathrm{Conn}(E) \times C^\infty(Y; \mathbb{S})) \times (0, \infty)$ to denote the set of pairs $(\mathfrak{c}', z')$ such that $\mathfrak{c}'$ has $C^2$ distance less than $\varepsilon$ from $\mathfrak{c}_0$ and $|z' - z_0| < \varepsilon$. Require in addition that the $(\mathfrak{c}', z')$ version of $\mathfrak{L}_{(\cdot)}$ have trivial kernel. Define $\mathcal{C}_{1\varepsilon}$ likewise using $(\mathfrak{c}_1, z_1)$. Granted this notation, define

$$|\mathfrak{f}_{s1} - \mathfrak{f}_{s0}| = \lim_{\varepsilon \to 0} \sup \{ |\mathfrak{f}_s(\mathfrak{c}_1', z_1') - \mathfrak{f}_s(\mathfrak{c}_0', z_0')| : (\mathfrak{c}_0', z_0') \in \mathcal{C}_{0\varepsilon} \textit{ and } (\mathfrak{c}_1', z_1') \in \mathcal{C}_{1\varepsilon} \}. \tag{B.19}$$

Perturbation theory can by used to prove that the lim-sup on the right in (B.19) is finite; and that it is equal to the norm of the honest spectral flow difference when both the $(\mathfrak{c}_1, z_1)$ and $(\mathfrak{c}_2, z_2)$ versions of $\mathfrak{L}_{(\cdot)}$ have trivial kernel. The limit in (B.9) is said in what follows to be the *norm of the difference between the values* $\mathfrak{f}_{\mathfrak{s}}$.

**Proposition B.3**: *There exists* $\kappa \geq 100$*, and given* $c_v \geq \kappa^2$*, there exists* $\kappa_c \geq \kappa$ *with the following significance: Suppose that* $c_v \geq \kappa$*, that* $r \geq \kappa_c c_v^{10}$ *and that* $\mu \in \Omega$ *has* $\mathcal{P}$*-norm bounded by 1. Let* $(A, \psi)$ *denote a solution to the* $(r, \mu)$ *version of (1.13). Use* $(A, \psi)$ *as directed in Section Ba to construct the pair* $(A_*, \psi_*)$*. The norm of the difference between the values of* $\mathfrak{f}_{\mathfrak{s}}$ *at* $(A, \psi)$ *and at* $(A_*, \psi_*)$ *is bounded by* $\kappa$.

***Proof of Proposition B.3***: The $\tau = 0$ point on the path $\tau \to (A_{*\tau}, \psi_{*\tau})$ is $(A - u^{-1}du, u\psi)$ where $u: Y \to S^1$ is a homotopically trivial map. This being the case, it is sufficient to exhibit an r and $(A, \psi)$ independent bound for the absolute value of the difference between $\mathfrak{f}_{\mathfrak{s}}$ at $(A - u^{-1}du, u\psi)$ and at $(A_*, \psi_*)$. Such a bound is derived in the subsequent five parts of the proof.

*Part 1*: Suppose that $\mathbb{L}$ is a given Hilbert space and that $\mathcal{L}$ is an unbounded, self-adjoint operator on $\mathbb{L}$. Assume that $\mathcal{L}$ has pure point spectrum with no accumulation points and such that each eigenvalue has finite multiplicity. Let $e$ denote a bounded, self adjoint operator on $\mathbb{L}$ and suppose that $\{e_\tau\}_{\tau \in [0,1]}$ is a real analytic family of bounded, self adjoint operators on $\mathbb{L}$ with $e_0 = 0$ and $e_1 = e$. Section 2 of [T3] explains how to label the eigenvalues of each $\tau \in [0,1]$ version of $\mathcal{L} + e$ by the integers so that the following is true: Given an integer n, let $\{\lambda_{n\tau}\}_{\tau \in [0,1]}$ denote the corresponding 1-parameter family of eigenvalues. Then the map $\tau \to \lambda_{n\tau}$: $[0, 1] \to \mathbb{R}$ is continuous and piecewise real analytic. Moreover, the corresponding 1-parameter family of eigenvectors varies in a real analytic fashion where $\lambda_{n(\cdot)}$ does. Let $\{\mathfrak{f}_{(\tau)}\}_{\tau \in [0,1]}$ denote a corresponding 1-parameter family of unit length eigenvectors. The map $\tau \to \mathfrak{f}_{(\tau)}$ can be assumed real analytic on the open subsets in $[0,1]$ where $\lambda_{n(\cdot)}$ is real analytic. As noted in [T3], the derivative of $\lambda_{n(\cdot)}$ where it is real analytic is given by

$$\frac{d}{d\tau}\lambda_{n\tau} = \langle \mathfrak{f}_{(\tau)}, (\tfrac{d}{d\tau} e_\tau)\, \mathfrak{f}_{(\tau)} \rangle_{\mathbb{L}} \,, \tag{B.20}$$

where $\langle \, , \rangle_{\mathbb{L}}$ denotes here the inner product on $\mathbb{L}$.

*Part 2*: Let $\mathbb{L}$ denote the Hilbert space $L^2(Y; \mathbb{V}_0 \oplus \mathbb{V}_1)$, let $\mathcal{L}$ denote the $z = r$ and $(A_{*0}, \psi_{*0})$ version of the operator $\mathfrak{L}_{\mathbb{V}}$, and $\mathcal{L} + e$ denote the corresponding $(A_*, \psi_*)$ version of this operator. The next lemma implies in part that what is said in Part 1 can be invoked for this version of $\mathbb{L}$, $\mathcal{L}$ and $e$. This lemma uses $\kappa_\lozenge$ to denote a number that is greater than the versions of the constant $\kappa$ that appear in Lemmas A.1-A.8 and B.1-B.2.

**Lemma B.4**: *Fix* $m \geq \kappa_\lozenge$. *There exists an* m*-dependent* $\kappa > 1$ *and given* $c_v \geq \kappa$, *there exists* $\kappa_{c_v} \geq \kappa$ *with the following significance: Suppose that* $r \geq \kappa_{c_v} c_v^{10}$ *and suppose that* $(A, \psi = (\alpha, \beta))$ *is a solution to the* $(r, \mu)$ *version of (1.13) with* $\mu$ *a given element in* $\Omega$ *with* $\mathcal{P}$*-norm smaller than 1. Construct as in Lemmas B.6 and B.7 the family of operators* $\{\mathfrak{L}_{\mathbb{V},\tau}\}_{\tau\in[0,1]}$. *Fix* $n \in \mathbb{Z}$ *and let* $\{\lambda_{n\tau}\}_{\tau\in[0,1]}$ *denote the corresponding family of eigenvalues. Then* $|\lambda_{n\tau}| = 0$ *for some* $\tau \in [0,1]$ *only if* $|\lambda_{n\tau'}| \leq \kappa c_v^{-m}$ *for all* $\tau' \in [0,1]$.

Lemma B.4 is proved in the upcoming Section Bd of this appendix.

Let $\kappa_\lozenge$ now denote a constant that is greater than the various versions of $\kappa$ that appear in Lemmas A.2-A.8 and Lemmas B.1-B.2 and B.4. Fix $c_v$ and r so that the assumptions of these lemmas are met. Let $\Delta$ denote the dimension of the span of the eigenvectors of $\mathfrak{L}_{\mathbb{V},1}$ with eigenvalue beween $-\kappa_\lozenge c_v^{-1}$ and $\kappa_\lozenge c_v^{-1}$. It is a consequence of Lemma B.4 that the norm of the spectral flow difference between $(A_{*0}, \psi_{*0})$ and $(A_*, \psi_*)$ is no greater than $\Delta$. This being the case, Proposition B.3 follows if $\Delta$ has an r and $(A, \psi)$-independent bound given a suitable r and $(A, \psi)$ independent choice of m and then $c_v$. A choice for $c_v$ that yields such a bound $\Delta$ is derived in Parts 4 and 5 of the proof. The subsequent Part 3 of the proof supplies two observations in the form of lemmas that are used in the derivation.

*Part 3*: To set the stage for the first lemma, use $\kappa_{\lozenge\lozenge}$ to denote the version of $\kappa$ given by Proposition 2.4. Fix $r \geq \kappa_{\lozenge\lozenge}$ and $\mu \in \Omega$ with $\mathcal{P}$-norm bounded by 1. Let $(A, \psi)$ denote a solution to the corresponding $(r, \mu)$ version of (1.13). Let $\gamma \subset Y$ denote a closed, connected segment in $\alpha$'s zero locus whose points have distance at least $100\,\kappa_{\lozenge\lozenge}^2 r^{-1/2}$ from all curves in $\cup_{p\in\Lambda}\{\hat{\gamma}_p^+ \cup \hat{\gamma}_p^-\}$. Use the coordinates from Part 4 of Section Aa to define the functions $\nu$ and $\mu$ on $\gamma$ and having done so, use $L_\gamma$ to denote the operator on $C^\infty(\gamma; \mathbb{C})$ that is defined by the rule $\zeta \to \frac{i}{2}\frac{d}{dt}\zeta + \nu\zeta + \mu\bar{\zeta}$. This operator defines a function on $C^\infty(\gamma; \mathbb{C})$ by the rule $\zeta \to \|L_\gamma\zeta\|_2$ where $\|\cdot\|_2$ denotes here the $L^2$ norm on $C^\infty(\gamma; \mathbb{C})$. This function can be restricted to any given linear subspace in $C^\infty(\gamma; \mathbb{C})$. Given $T > 0$, there is always an integer that is greater than or equal to the dimension of any linear subspace in $C^\infty(\gamma; \mathbb{C})$ on which the function $\zeta \to \|L_\gamma\zeta\|_2$ obeys $\|L_\gamma\zeta\|_2 \leq T\|\zeta\|_2$.

The upcoming Lemma B.5 concerns $L_\gamma$ and a least upper bound of the sort just described. By way of a parenthetical remark, the versions of $L_\gamma$ that appear in Lemma A.8 are of particular interest with regards to the proof of Proposition B.3.

**Lemma B.5**: *There exists* $\kappa > \kappa_{\diamond\diamond}$ *with the following significance: Fix* $r \geq \kappa$ *and* $\mu \in \Omega$ *with* $\mathcal{P}$*-norm bounded by* 1. *Suppose that* $(A, \psi = (\alpha, \beta))$ *is a solution to the corresponding* $(r, \mu)$ *version of (1.13). Let* $\gamma$ *denote a closed, connected segment of the zero locus of* $\alpha$ *whose points have distance at least* $100\kappa_{\diamond\diamond} r^{-1/2}$ *from all of the curves in the set* $\cup_{p\in\Lambda}\{\hat{\gamma}_p^+ \cup \hat{\gamma}_p^-\}$. *Use* $\Delta_\gamma$ *to denote the least upper bound for the dimensions of the linear subspaces in* $C^\infty(\gamma;\mathbb{C})$ *on which the function* $\zeta \to \|L_\gamma\zeta\|_2$ *obeys* $\|L_\gamma\zeta\|_2 \leq \kappa^{-1}\|\zeta\|_2$. *This least upper bound obeys* $\Delta_\gamma \leq \kappa$.

The argument for this lemma would be straight-forward were there an r-independent upper bound on $\gamma$'s length, but such bound does not exist. In any event, the proof is given momentarily. The next lemma states an analog of Lemma B.5 with the straight-forward argument for its proof. This lemma enters the proof of Proposition B.3 in conjunction with Lemma A.7. Lemma B.6 also plays a role in Lemma B.5's proof.

**Lemma B.6**: *There exists* $\kappa \geq 1$ *with the following significance: Fix* $R > 0$ *and* $T > 0$. *The least upper bound for the dimensions over* $\mathbb{C}$ *of the linear subspaces in* $C^\infty([0,R];\mathbb{C})$ *on which the function* $\zeta \to \|\frac{d}{dt}\zeta\|_2$ *obeys* $\|\frac{d}{dt}\zeta\|_2 \leq T\|\zeta\|_2$ *is bounded by* $(2 + \pi^{-1}RT)$.

***Proof of Lemma B.6***: The subset of elements in $C^\infty([0,R];\mathbb{C})$ that vanish at both endpoints has complex codimension 2. This understood, the least upper bound in question is no greater than 2 plus the number of linearly independent eigenvectors for the operator $-\frac{d^2}{dt^2}$ on $C^\infty([0,R];\mathbb{C})$ that vanish at both endpoints and have eigenvalue less than $T$. This number is $\pi^{-1}RT$.

***Proof of Lemma B.5***: The four steps that follow constitute the proof.

Step 1: The lemma is proved by cutting $\gamma$ into at a concatenation of $c_0$ closed, connected segments, and then bounding a version of $\Delta_{(\cdot)}$ on each segment. To explain why such a cutting strategy works, suppose for the moment that $\gamma_0 \subset \gamma$ is a closed, connected segment. Fix $T > 0$ and introduce $\Delta_{\gamma_0,T}$ to denote the least upper bound for the dimensions of the linear subspaces in $C^\infty(\gamma_0;\mathbb{C})$ on which the function $\|L_{\gamma_0}(\cdot)\|_2$ is bounded by $T^{-1}\|\cdot\|_2$ with $\|\cdot\|_2$ denoting here the $L^2$ norm on $C^\infty(\gamma_0;\mathbb{C})$. Suppose that $\gamma_1, \gamma_2$ are two such segments that share at least one endpoint. Then $\Delta_{\gamma_1\cup\gamma_2,T} \leq 4 + \Delta_{\gamma_1,T} + \Delta_{\gamma_2,T}$. This is because the subspace in $C^\infty(\gamma;\mathbb{C})$ that vanishes at the common end points of $\gamma_1$ and $\gamma_2$ has codimension 2 if they share one endpoint and codimension 4 if they share two endpoints.

With the preceding understood, suppose that $\gamma$ is written as the concatentation of N segments $\{\gamma_k\}_{1\le k\le N}$. Iterate the bound given in the previous paragraph to see that $\Delta_{\gamma,T}$ is no greater than $4N + \sum_{1\le k\le N} \Delta_{\gamma_k,T}$.

Step 2: Fix $\varepsilon > 0$ and let $\gamma^\varepsilon \subset \gamma$ denote the part of $\gamma$ with distance at least $\varepsilon$ from the curves in the set $\cup_{p\in\Lambda}\{\hat\gamma_p^+ \cup \hat\gamma_p^-\}$. As explained directly, $\Delta_{\gamma^\varepsilon,T} \le c_0(1 + T^{-2})|\ln\varepsilon|$.

To see why this bound holds, keep in mind that $L_\gamma$ is defined by the pair $(\nu,\mu)$ and the latter are defined by the chosen coordinates from Part 4 of Section Aa. Granted that such is the case, any version of $L_\gamma$ can be obtained from a given version by conjugating the given version with a map from $\gamma$ to $S^1$. This implies, in particular, that $\Delta_{\gamma^\varepsilon,T}$ does not depend on the choice of coordinates. This being the case, choose the coordinates from Part 4 of Section Aa so that the resulting pair $\nu$ and $\mu$ are such that $|\nu| + |\mu| \le c_0$. Use B to denote to denote an upper bound for $|\nu|$ and $|\mu|$ on $\gamma^\varepsilon$

Let L denote the length of $\gamma^\varepsilon$. Let $\mathcal{V} \subset C^\infty(\gamma^\varepsilon;\mathbb{C})$ denote a linear subspace of the sort under consideration. If $\gamma^\varepsilon$ has no endpoints, then the dimension of $\mathcal{V}$ is no greater than the dimension of the span of the eigenvectors of $-\frac{d^2}{dt^2}$ with eigenvalue no greater than $(B+T^{-1})^2$. As noted in the proof of Lemma B.6, if $\gamma^\varepsilon$ has endpoints, then $\mathcal{V}$'s dimension is at most 4 more than the span of the Dirichelet eigenvectors of $-\frac{d^2}{dt^2}$ with eigenvalue no greater than $(B+T^{-1})^2$. In both cases, there are at most $c_0(1+B+T^{-1})^2L$ linearly independent eigenvectors with this eigenvalue bound. Meanwhile, Proposition 2.4 with Proposition II.2.7 and Lemma II.2.2 imply among other things that the length of $\gamma^\varepsilon$ is no greater than $c_0|\ln\varepsilon|$, and that both $|\nu|$ and $|\mu|$ are bounded by $c_0$.

Step 3: Let $\hat\gamma \in \cup_{p\in\Lambda}\{\hat\gamma_p^+ \cup \hat\gamma_p^-\}$ denote a given curve. As explained in Part 4 of Section 1a, there is a version of the coordinates from Part 4 of Section Aa for $\hat\gamma$ with both $\nu$ and $\mu$ constant, with $\mu$ real and such that $\mu > |\nu|$. This version is assumed in what follows. The corresponding constant values for $\nu$ and $\mu$ are denoted by $\nu_0$ and $\mu_0$.

Fix $\varepsilon > 0$ such that the radius $\varepsilon$ tubular neighborhood of $\hat\gamma$ is well inside the coordinate chart just described. Let T denote such a tubular neighborhood, and suppose that $\upsilon \subset T$ is a closed, connected segment in T of an integral curve of $v$. Taylor's theorem with remainder can be used with the formulae in (A.3) to see that $\upsilon$ has a tubular neighborhood with coordinates from Part 4 of Section Aa with $|\nu - \nu_0| + |\mu - \mu_0| < c_0\varepsilon$.

Reintroduce from Proposition 2.4 the subset $Y_r \subset Y$. By way of a reminder, the points in $T \cap Y_r$ have distance no less than $c_0 r^{-1/2}$ from $\hat\gamma$. Let $\gamma^T$ denote a properly embedded, connected component of $\alpha$'s zero locus in the closure of $T \cap Y_r$. Thus, $\gamma^T$ has two boundary points, either both on the boundary of the closure of T, or one on the latter and one on $Y_r$'s boundary torus in T.

<u>Step 4</u>: Define the operator $L_0$ on $C^\infty(\mathbb{R};\mathbb{C})$ by the rule $\zeta \to \frac{i}{2}\frac{d}{dt}\zeta + \nu_0\zeta + \mu_0\bar{\zeta}$. Fix $\ell > 0$ and restrict $L_0^2$ to the subspace of elements in $C^\infty([0, \ell\,];\mathbb{C})$ that vanish at the boundary points. The corresponding Dirichelet eigenvalues of $L_0^2$ on this domain are of the form $\kappa^2 + \nu_0^2 + \mu_0^2 \pm 2\nu_0(\kappa^2+\mu_0)^2$ with $\kappa = \frac{\pi(2k+1)}{4\ell}$ for $k \in \mathbb{Z}$. Note in particular that the smallest eigenvalue is greater than $(\mu_0-\nu_0)^2$ when $\mu_0 > \nu_0$, thus greater than $c_0^{-1}$.

Let $L_{\gamma^T}$ denote the restriction of $L_\gamma$ to $C^\infty(\gamma^T;\mathbb{C})$. What was said in the preceding paragraph and what was said in the final paragraph of Step 3 have the following implication: Let $\zeta \in C^\infty(\gamma^T;\mathbb{C})$ denote an element that vanishes at both boundary points of $\gamma^T$. Then $\|L_{\gamma^T}\zeta\|_2 \geq ((\mu_0-\nu_0)-c_0\varepsilon)\|\zeta\|_2$. Thus, if $\varepsilon < c_0^{-1}$, then

$$\|L_{\gamma^T}\zeta\|_2 \geq \tfrac{1}{2}(\mu_0-\nu_0)\|\zeta\|_2\,; \tag{B.21}$$

and this implies that $\Delta_{\gamma^T,(\nu_0-\mu_0)/2} \leq 4$. Note that (B.21) and the preceding bound on $\Delta_{\gamma^T,(\nu_0-\mu_0)/2}$ hold even if the pair $(\nu,\mu)$ that define L along $\gamma^T$ are not $\varepsilon$ close to $(\nu_0, \mu_0)$. This is because any two versions of $(\nu,\mu)$ that arise from different choices for coordinates from Part 4 of Section Aa define corresponding versions of L that are obtained from each other by conjugating with a map from $\gamma$ to $S^1$.

Lemma B.5 follows from the $\Delta_{\gamma^T,(\nu_0-\mu_0)/2} \leq 4$ bound and those given in Step 2.

*Part 4*: Fix m and then $c_v$ and r so as to invoke the conclusions of the Lemmas A.2-A.8 and Lemmas B.4-B.6. Sum Lemma B.5's integers $\{\Delta_\gamma\}_{\gamma\in\Theta}$ and use $\Delta_\Theta$ to denote the result. Let $\kappa_\bullet$ denote the larger of the version of the constant $\kappa$ that appears in Lemmas A.4, A.7 and B.5. Let N denote the number of components of $Y-Y_{*\Lambda}$ with zeros of $\alpha$. Each such component has the same length, this denoted by $\ell_*$. Lemma A.7 associates an integer m to each such component. As noted in Part 3 of Section Ad, no version of m is greater than $\kappa_\bullet c_v^4$.

Use N to denote the number of linearly independent eigenvectors of $\mathfrak{L}_{\mathbb{V},1}$ with eigenvalue between $-c_v^{-1}$ and $c_v^{-1}$. As is explained in the subsequent paragraphs, N is no greater than $10^4(\Delta_\Theta + N\kappa_* c_v^4(1 + \ell_*\kappa_\bullet c_v^{\kappa_\bullet}))$ if $c_v \geq c_0$ and r is larger than a constant that depends only on $c_v$.

To see why this bound holds, suppose in what follows that N is larger than what is claimed so as to generate nonsense. If N is larger than the asserted bound, then there exists a section, $\mathfrak{f}$, of $\mathbb{V}_0 \oplus \mathbb{V}_1$ with four properties that are described next. Lemmas B.5 and B.6 guarantee that the third and fourth properties can be satisfied. First, $\mathfrak{f}$ is a linear combination of an orthonormal set of eigenvectors of $\mathfrak{L}_{\mathbb{V},1}$ with eigenvalue between $-c_v^{-1}$ and $c_v^{-1}$. Second, $\mathfrak{f}$ has unit $L^2$ norm. The third property concerns the curves in $\Theta$. Let $\gamma$ denote such a curve. Let $\zeta_\gamma$ denote $\gamma$'s component of $\Pi_\vartheta\mathfrak{f}$. The function $\zeta \to \|L_\gamma\zeta\|_2$ from Lemma B.5 is such that $\|L_\gamma\zeta_\gamma\|_2 \geq \kappa_\bullet^{-1}\|\zeta_\gamma\|_2$. The final property concerns the components of $Y-Y_{*\Lambda}$ that contain zeros of $\alpha$. Let $\gamma \in \cup_{p\in\Lambda}\{\hat{\gamma}_p^+ \cup \hat{\gamma}_p^-\}$ denote a curve from such a component. Part 1 of Section Ah describes a cover of $\gamma$ by open sets $\{\gamma_k\}_{0\leq k\leq 7}$ or by open sets $\gamma_+$ and $\gamma_-$. Let $\gamma_*$ denote a component of this cover. Denote by m the fiber dimension of $\mathrm{Ker}_\vartheta|_{\gamma_*}$. Use the isomorphism between $\mathrm{Ker}_\vartheta|_{\gamma_*}$ and $\gamma_* \times \mathbb{C}^m$ in Lemma A.7 to view the component of $\Pi_\vartheta\mathfrak{f}$ in $\mathrm{Ker}_\vartheta|_{\gamma_*}$ as a map from $\gamma_*$ to $\mathbb{C}^m$. Let $\zeta_{\gamma_*}$ denote this map. Then $\|\frac{i}{2}\frac{d}{dt}\zeta_{\gamma_*}\|_2$ is greater than $2\kappa_\bullet c_v^{\kappa_\bullet}\|\zeta_{\gamma_*}\|$.

What with the third and fourth properties, Lemmas A.7 and A.8 imply that $\|\Pi_\vartheta\mathfrak{L}_{\mathbb{V},1}\mathfrak{f}\|_2 \geq \kappa_\bullet^{-1}\|\Pi_\vartheta\mathfrak{f}\|_2$. This being the case, Lemma A.6 finds $\|\Pi_\vartheta\mathfrak{L}_{\mathbb{V},1}\mathfrak{f}\|_2 \geq \frac{1}{2}\kappa_\bullet^{-1}\|\mathfrak{f}\|_2$ if $c_v \geq c_0$ and if r is creater than a purely $c_v$-dependent constant. Meanwhile, it follows from the definitions that $\|\mathfrak{q}\|_2 \geq c_0^{-1}\|\Pi_\vartheta\mathfrak{q}\|_2$ for any given section $\mathfrak{q}$ of $\mathbb{V}_0\oplus\mathbb{V}_1$. Take $\mathfrak{q}$ to be $\mathfrak{L}_{\mathbb{V},1}\mathfrak{f}$ to see that $\|\mathfrak{L}_{\mathbb{V},1}\mathfrak{f}\|_2 \geq (c_0\kappa_\bullet)^{-1}\|\mathfrak{f}\|_2$. Even so, the first property listed in the preceding paragraph requires the bound $\|\mathfrak{L}_{\mathbb{V},1}\mathfrak{f}\|_2 \leq c_v^{-1}\|\mathfrak{f}\|_2$ and so $c_v^{-1} \geq (c_0\kappa_\bullet)^{-1}$ unless $\mathfrak{f}$ is identically zero; and $\mathfrak{f} \neq 0$ because of the second of the listed properties.

This $(A,\psi)$ and r-independent lower bound for $c_v^{-1}$ is the required nonsense because $c_v$ has no a priori upper bound.

## d) Proof of Lemma B.4

The proof of Lemma B.4 has three parts.

*Part 1*: This part of the proof states an auxilliary lemma that augments what is said by Lemma B.2. By way of a reminder, Lemma B.2 concerns the $Y_{\Diamond\Diamond}$ extension of a given $\gamma \in \Theta$ version of $U_\gamma$. With coordinates from Part 4 of Section Aa chosen, Lemma B.2 uses the $\alpha \to |\alpha|\frac{z}{|z|}$ isomorphism from $E|_{U_\gamma}$ to $U_\gamma \times \mathbb{C}$ to write $(A - u_\gamma^{-1}du_\gamma, (u_\gamma\alpha, u_\gamma\beta))$ as $(\theta + a_{A,\gamma}, (\alpha_\gamma, \beta_\gamma))$ with the map $u_\gamma$: $U_\gamma \to S^1$ defined in (B.8). It goes on to write the $i\mathbb{R}$ valued 1-form $a_{A,\gamma}$ as $a_{A0,\gamma}dt + \frac{1}{2}(A_\gamma d\bar{z} - \bar{A}_\gamma dz)$. Whereas Lemma B.2 talks about the functions A and $\alpha_\gamma$, the upcoming Lemma B.7 talks about the functions $a_{A0,\gamma}$ and $\beta_\gamma$. This lemma brings in the functions $\varsigma$ and y on $\mathbb{C}$ from (A.2); and it uses $u$ to denote the function of the radial coordinate $|z|$ on $\mathbb{C}$ that given by integrating the function $1 - |\alpha_0|^2$ along any ray from the origin starting at distance $|z|$ from the origin. Thus

$$u(|z|) = \int_{|z|}^{\infty} (1-|\alpha_0|^2)|_s ds .$$

(B.22)

The lemma also invokes the coordinates from Part 4 of Section Aa to bring in the corresponding version of the function $t \to x_\gamma(t)$.

**Lemma B.7**: *Fix* m ≥ 1. *There exists an* m*-dependent* $\kappa > 1$ *and given* $c_v \geq \kappa$*, there exists* $\kappa_{c_v} \geq \kappa$ *with the following significance: Take* $r \geq \kappa_{c_v} c_v^{10}$ *and suppose that* $(A, \psi = (\alpha, \beta))$ *is a solution to the* $(r, \mu)$ *version of (1.13) with* $\mu$ *a given element in* $\Omega$ *with* $\mathcal{P}$*-norm smaller than 1. Use* $(A, \psi)$ *as directed in Section Ba to construct the pair* $(A_*, \psi_*)$. *Let* $\gamma$ *denote a component of* $\alpha$*'s zero locus in* $Y_{*\Lambda}$.

- $|r^{1/2}\beta_\gamma - i\mu r_r^* \varsigma| \leq c_v^{-m}$ *at all points in the* $Y_{*\Lambda}$ *part of* $U_\gamma$.
- $|a_{A0,\gamma} - i\nu 2^{1/2} r_r^* y + i r (x_\gamma \frac{\bar z}{|z|} + \bar x_\gamma \frac{z}{|z|}) r_r^* u| \leq c_v^{-m}$ *at all points in the* $Y_{*\Lambda}$ *part of* $U_\gamma$.

***Proof of Lemma B.7***: The proof has six steps. The first four prove the top bullet and the last two prove the lower bullet. As in the proofs of Lemmas B.1 and B.2, what is denoted by $\kappa_c$ is a constant with value greater than 1 that depends only on m and $c_v$; in particular, it has no r and (A, ψ) dependence.

<u>Step 1</u>: The $EK^{-1}$ component of the equation $D_A{}^2\psi = 0$ can be written schematically as $-((\nabla_A)_v)^2\beta - \bar\partial_A \partial_A \beta + \frac{1}{2} r|\alpha|^2 \beta = -\mu \partial_A \alpha + \mathfrak{r}$ where $|\mathfrak{r}| \leq c_0$. This equation is used in what follows on the extension of $U_\gamma$ into $Y_{\diamond\diamond}$ with the coordinates from Part 4 of Section Aa. The section β of $EK^{-1}$ is viewed as a $\mathbb{C}$-valued function on $U_\gamma$ using these coordinates and the chosen isomorphism on $U_\gamma$ between E and $U_\gamma \times \mathbb{C}$. Meanwhile, connection $A - u_\gamma^{-1}du_\gamma$ is written as $A_\gamma = a_{a0,\gamma}dt + \frac{1}{2}(A_\gamma d\bar z - \bar A_\gamma dz)$. Use the derivative bounds for β given by Lemma 2.1 to replace the derivative $(\nabla_A)_v$ with $\partial_t$ at so as to obtain from (A.36) the equation

$$-\frac{1}{4}\partial_t^2\beta_\gamma - \bar\partial_A \partial_A \beta_\gamma + \frac{1}{2} r|\alpha|^2 \beta_\gamma = -\mu \partial_A \alpha_\gamma + \mathfrak{r}$$

(B.23)

with $\mathfrak{r}$ different and now such that $|\mathfrak{r}| \leq c_0 c_v^2$ when $c_v \geq c_0$ and $r \geq \kappa_c$. The notation here uses $\bar\partial_A = \frac{\partial}{\partial \bar z} + \frac{1}{2} A_\gamma$ and it uses $\partial_A$ for the complex conjugate operator.

Introduce by way of notation $\beta_*$ to denote $\mu r^{-1/2} r_r^* \varsigma$. It follows from (3.27) that the section $\beta_*$ obeys an $(A_*, \alpha_*)$ analog of (B.23):

$$-\frac{1}{4}\partial_t^2\beta_* - \bar\partial_{A_*} \partial_{A_*} \beta_* + \frac{1}{2} r|\alpha_*|^2 \beta_* = -\mu \partial_{A_*} \alpha_* + \mathfrak{r} ,$$

(B.24)

with $\mathfrak{r}$ here different from its (B.23) incarnation but such that $|\mathfrak{r}| \le c_0 c_v^2$.

Step 2: Use $\Delta_A$ to denote $A_\gamma - A_*$ and use $\Delta_\alpha$ to denote $\alpha_\gamma - \alpha_*$. Their absolute values on $U_\gamma$ are such that $r^{-1/2}|\Delta_A| + |\Delta_\alpha| \le r^{1/2} e^{-c_v}$ when $c_v \ge \kappa$ and $r \ge \kappa_c$. Write the connection $A_*$ as $A - \Delta_A$, and write $\alpha_*$ as $\alpha - \Delta_\alpha$ to rewrite (B.22) as

$$-\tfrac{1}{4}\partial_t^2\beta_* - \bar\partial_A\partial_A\beta_* + \tfrac{1}{2}r|\alpha|^2\beta_* = -\mu\,\partial_A\alpha + \mu\,\partial_A\Delta_\alpha - \mu\,\Delta_A^{(1,0)}\Delta_\alpha - \mathfrak{R}\cdot\beta_* + \mathfrak{r} \tag{B.25}$$

where the notation uses $\Delta_A^{(1,0)}$ to denote the (1,0) part of $\Delta_A$. Meanwhile, what is written as $\mathfrak{R}\cdot\beta_*$ is linear in $\beta_*$ and can be written as

$$(\bar\partial_A\Delta_A^{(1,0)})\beta_* + \Delta_A^{(0,1)}\partial_A\beta_* + \Delta_A^{(0,1)}\bar\partial_A\beta_* - \mathfrak{w}_0(\Delta_A,\Delta_A)\beta_* + r\,\mathfrak{w}_1(\Delta_\alpha)\beta_* \tag{B.26}$$

where $|\mathfrak{w}_{0,1}| \le c_0$ and with $\mathfrak{r}$ different but still obeying $|\mathfrak{r}| \le c_0 c_v^2$.

Let $\Delta_\beta = \beta - \beta_*$. The two equations (B.23) and (B.25) imply that $\Delta_\beta$ obeys

$$-\tfrac{1}{4}\partial_t^2\Delta_\beta - \bar\partial_A\partial_A\Delta_\beta + \tfrac{1}{2}r|\alpha|^2\Delta_\beta = -\mu\,\partial_A\Delta_\alpha + \mu\,\Delta_A^{(1,0)}\Delta_\alpha + \mathfrak{R}\cdot\beta_* + \mathfrak{r} \tag{B.27}$$

where $\mathfrak{r}$ is again a term with absolute value bounded by $c_0 c_v^2$. Use $o_\beta$ to denote the function $\tfrac{1}{2}|\Delta_\beta|^2$. Take the Hermitian inner product of both sides of (B.27) with $\Delta_\beta$ and commute covariant derivatives of A to obtain an equation for $o_\beta$, this being the next equation. This upcoming equation uses $\nabla_A^{\perp}$ to denote the covariant derivative along the constant t slices of $U_\gamma$ and it uses $re[\cdot]$ to denote the real part of indicated expression. What follows is the promised equation for $o_\beta$:

$$\begin{aligned}&-\tfrac{1}{4}\partial_t^2 o_\beta - \bar\partial\partial\, o_\beta + \tfrac{1}{2}r(1+|\alpha|^2)o_\beta = \\ &\qquad -\tfrac{1}{4}|\partial_t\Delta_\beta|^2 - \tfrac{1}{4}|\nabla_A^{\perp}\Delta_\beta|^2 + re[-\mu\,\bar\Delta_\beta\,\partial_A\Delta_\alpha + \mu\,\bar\Delta_\beta\,\Delta_A^{(1,0)}\Delta_\alpha + \bar\Delta_\beta\,\mathfrak{R}\cdot\beta_*] + \mathfrak{r}\,,\end{aligned} \tag{B.28}$$

What is denoted by $\mathfrak{r}$ here signifies a term that with absolute value bounded by $c_0 c_v^2$.

Step 3: Let $p = (t, z)$ denote a given point in $\mathbb{R}\times\mathbb{C}$ and introduce $G_p(\cdot)$ to denote the Green's function for the operator $-(\partial_t^2 + 4\bar\partial\partial) + 2r$ on $\mathbb{R}\times\mathbb{C}$ with pole at p. This Green's function is positive and is such that

$$G_p \le c_0\,\tfrac{1}{|p-(\cdot)|}e^{-\sqrt{r}|p-(\cdot)|} \quad \textit{and} \quad |dG_x| \le c_0\left(\tfrac{1}{|p-(\cdot)|^2} + \sqrt{r}\right)e^{-\sqrt{r}|p-(\cdot)|}\,. \tag{B.29}$$

where d here denotes the full exterior derivative on $\mathbb{R} \times \mathbb{C}$. Introduce $\chi_\diamond$ to denote the function on the $Y_{**}$ extended $U_\gamma$ that is defined as follows: Take $\chi_\diamond = \chi_{\hat{U}}$ on $Y_{*\Lambda}$ and take it equal to $\chi_{\hat{U}}(1 - \chi_{\diamond\diamond})$ on $U_\gamma$'s intersection with a given component of $Y_{\diamond\diamond} - Y_{*\Lambda}$. Thus, $\chi_\diamond$ has compact support on the interior of $U_\gamma$ and it is equal to one at the points in $Y_{*\Lambda}$ with distance $c_v^2 r^{-1/2}$ or less from $\gamma$. Take p to be a point in $U_\gamma$ where $\chi_{\hat{U}}$ is equal to 1. Multiply both sides of (B.29) by $(\chi_\diamond)^2 G_x$ and then integrate both sides to obtain an equality between two integrals. The left hand side integral is denoted by $I_L$ and the right hand side by $I_R$. As explained in the next step of the proof, the equality $I_L = I_R$ leads directly to the bound $|\Delta_\beta| \le c_0 r^{-1/2} c_v^{10} e^{-c_v/2}$ at p. This proves the assertion of the top bullet of Lemma B.7 at points in $U_\gamma \cap Y_{*\Lambda}$ with distance $c_v^2 r^{-1/2}$ or less from $\gamma$ if $c_v \ge c_0$ and $r \ge \kappa_c$. Meanwhile, Lemmas 2.1 and 2.2 with (3.3) imply what is asserted by the top bullet of Lemma B.7 on the points in $U_\gamma \cap Y_{*\Lambda}$ with distances greater than $c_v^2 r^{-1/2}$ from $\gamma$ if $c_v \ge c_0$ and $r \ge \kappa_c$.

Step 4: The asserted bound on $I_L$ is derived by integrating by parts twice in the relevant integral.. The result can be written as $I_L = \frac{1}{2} |\Delta_\beta|^2(x) + \mathfrak{e}$ where $\mathfrak{e}$ is a function with $|\mathfrak{e}| \le c_0 r^{-1} e^{-c_v^2/c_0}$. By way of explanation, the function $\mathfrak{e}$ comes from an integral whose integrand has a term that is bounded in absolute value by $c_0 |\Delta_\beta|^2 (|d\chi_\diamond|^2 + |d^\dagger d\chi_\diamond|) G_x$ and one that is bounded in absolute value by $c_0 |\Delta_\beta|^2 |d\chi_x| |dG_x|$. Since these terms are supported at distances no less than $(1 + c_0^{-1}) c_v^2 r^{-1/2}$ from x, and since $|\Delta_\beta| \le |\beta| + |\beta_*| \le c_0 r^{-1/2}$ in any event, the claim about $I_L$ is a consequence of the exponential factors in (B.29).

Meanwhile, the integral $I_R$ is bounded by $c_0 r^{-1} c_v^{10} e^{-c_v}$. By way of explanation, some judiciously chosen applications of integration by parts will remove derivatives along the $\mathbb{C}$ factor of $\mathbb{R} \times \mathbb{C}$ from $\Delta_\alpha$ and $\Delta_A$ and replace them with terms that have derivatives of either $G_x$, $\chi_\diamond$ or covariant derivatives of $\beta_*$. A covariant derivative of $\beta_*$ is bounded by $c_0(|\nabla_A^\perp \Delta_\beta| + |\nabla_A \beta|)$. Lemma 2.1 has $|\nabla_A \beta| \le c_0$ and this with (B.29) together with the bounds from Step 2 for $|\Delta_A|$ and $|\Delta_\alpha|$ can be used in a straightforward fashion to bound the integrals that result by $c_0 r^{-1} c_v^{20} e^{-c_v}$.

Step 5: This step begins the proof of the lower bullet of Lemma B.7. As is explained in this step and Step 6, the second bullet's assertion is a consequence of the identity in (B.18) and what is asserted by a version of Lemma B.7's top bullet that uses a suitable $m' > m$. To exploit (B.18), first write $\frac{\partial}{\partial |z|}$ as $\frac{z}{|z|} \frac{\partial}{\partial z} + \frac{\bar{z}}{|z|} \frac{\partial}{\partial \bar{z}}$. Next write $\frac{\partial}{\partial z}$ as $\frac{1}{2}(\hat{e}_1 - i\hat{e}_2) + \mathfrak{r}_t \frac{\partial}{\partial t} + \mathfrak{r}_2$ where $\{\hat{e}_1, \hat{e}_2\}$ are an oriented, orthonormal frame for the kernel of $\hat{a}$, where $|\mathfrak{r}_t| \le c_0 |z|$, and where $|\mathfrak{r}_2| \le c_0 |z|^2$. Doing so writes $\frac{\partial}{\partial |z|}$ in terms of $\{\hat{e}_1, \hat{e}_2\}$. Use this depiction in (A.6). Meanwhile, use (A.6) to write $\frac{\partial}{\partial t}$ in terms of $v$. The result of all of this rewriting replaces the the curvature component on the right hand side of (B.18) by

$$- |z|^{-1}(z_1 F_A(\nu,\hat{e}_1) + z_2 F_A(\nu,\hat{e}_2) + (2\nu|z|^2 + \mu\,\bar{z}^2 + \bar{\mu}\,z^2 - x_\gamma \bar{z} - \bar{x}_\gamma z) F_A(\hat{e}_1,\hat{e}_2)) + \mathfrak{r}\ ,$$

(B.30)

where the notation is such that $z_1$ and $z_2$ denote the repective real and imaginary parts of z, and where $\mathfrak{r}$ denotes a term with absolute value bounded by $|z|^2|F_A|$.

Since the 2-form $F_A$ is the Hodge dual of $B_A$, the equation in (1.13) can be invoked to replace (B.30) with

$$-r\,|z|^{-1}(\bar{\alpha}\beta\,\bar{z} - \alpha\bar{\beta}\,z + i(2\nu|z|^2 + \mu\,\bar{z}^2 + \bar{\mu}\,z^2 - x_\gamma \bar{z} - \bar{x}_\gamma z)(1 - |\alpha|^2) + \mathfrak{r}'\ ,$$

(B.31)

where $\mathfrak{r}'$ has norm obeying $|\mathfrak{r}'| \le c_0\, c_v\, e^{-\sqrt{r}|z|/c_0}$ , this due to Lemmas 2.1 and 2.2. A further rewriting uses the top bullet in Lemma B.4 to replace $\beta$ in (B.31) by $i\,\mu\, r^{-1/2}\, r_r^{*}\varsigma$ plus a term with small norm. Make this substitution and then invoke the formula for $\varsigma$ in (A.2) to write (B.31) as

$$-i\,r\,(2\nu|z| - x_\gamma \tfrac{\bar{z}}{|z|} - \bar{x}_\gamma \tfrac{z}{|z|})(1 - |\alpha|^2) + \mathfrak{r}''$$

(B.32)

where $|\mathfrak{r}''| \le c_0\,(c_v^{20-m'} r^{1/2} + c_v)\, e^{-\sqrt{r}|z|/c_0}$ .

<u>Step</u> <u>6</u>: Granted (B.32), use the formula for y in (A.2) and the formula for for $A_0$ in (A.3) to write $y = -2^{1/2}(a_0 - 1)$. Keeping in mind that $a_0$ is real and a function of $|z|^2$, use the formula in (A.3) and the top bullet in (2.8) to see that $\partial_{|z|} a_0 = |z|(1 - |\alpha_0|^2)$. This understood, it follows from (B.18) and (B.32) that

$$\partial_s(a_{A0,\gamma} - i\nu\, 2^{1/2}\, r_r^{*} y)|_{sz} = i\,r\,(x_\gamma \tfrac{\bar{z}}{|z|} + \bar{x}_\gamma \tfrac{z}{|z|})(1 - |r_r^{*}\alpha_0|^2) + \mathfrak{r}''\ .$$

(B.33)

To exploit (B.33), first look again at what is said in the first paragraph from Step 4 of the proof of Lemma B.2 to see that $|a_{A0,\gamma}| \le c_0\, c_v^2\, e^{-\sqrt{r}|z|/c_0}$ where $|z| \ge r^{-1/2}$ on $U_\gamma \cap Y_{*\Lambda}$. Given (3.3), such a bound also holds for $r_r^{*}y$. These bounds imply what is asserted by the second bullet of Lemma B.7 at the points in $U_\gamma \cap Y_{*\Lambda}$ with distance greater than $c_v^2 r^{-1/2}$ from $\gamma$. Given this last observation, integrate both sides of (B.33) from a given value of s to $c_v^2 r^{-1/2}$ with a choice for $m' > 2m + 100$ to obtain the second bullet's assertion on the part $U_\gamma \cap Y_{*\Lambda}$ at points with distance less then $c_v^2 r^{-1/2}$ from $\gamma$ when $c_v \ge c_0$ and $r \ge \kappa_c$.

*Part 2*: The next lemma writes the various versions of $\mathfrak{L}_{\mathbb{V},\tau}$ as $\mathfrak{L}_{\mathbb{V},1} + e_\tau$ and gives a bound for the norm of the $\tau$-derivative of $e_\tau$.

**Lemma B.8**: *Fix* $m \geq 1$. *There exists an* m*-dependent* $\kappa > 1$ *and given* $c_v \geq \kappa$*, there exists* $\kappa_{c_v} \geq \kappa$ *with the following significance: Take* $r \geq \kappa_{c_v} c_v^{10}$ *and suppose that* $(A, \psi = (\alpha, \beta))$ *is a solution to the* $(r, \mu)$ *version of (1.13) with* $\mu$ *a given element in* $\Omega$ *with* $\mathcal{P}$*-norm smaller than 1. Define* $\{(A_{*\tau}, \psi_{*\tau})\}_{\tau \in [0,1]}$ *as instructed in Section Ba. Given* $\tau \in [0,1]$*, let* $\mathfrak{L}_{\mathbb{V},\tau}$ *denote the* $(A_{*\tau}, \psi_{*\tau})$ *and* $z = r$ *version of the operator* $\mathfrak{L}_{\mathbb{V}}$*. Write* $\mathfrak{L}_{\mathbb{V},\tau}$ *as* $\mathfrak{L}_{\mathbb{V},1} + e_\tau$*. Then the resulting map* $\tau \to e_\tau$ *from* $[0, 1]$ *to* $C^\infty(Y; \mathbb{V}_0 \oplus \mathbb{V}_1)$ *is real analytic with derivative bound* $|\frac{d}{d\tau} e_\tau| \leq c_v^{-m} r^{1/2}$ *at all points in* Y.

***Proof of Lemma B.8***: Given Lemmas B.2 and B.7, the assertion is a direct consequence of the formula for $(A_{*\tau}, \psi_{*\tau})$ in Section Ba and the formula for $\mathfrak{L}_{\mathbb{V}}$ in (A.26) and (A.27).

The lemma that follows uses the $|\frac{d}{d\tau} e_\tau| \leq c_v^{-m} r^{1/2}$ bound from Lemma B.8 to give a version of Lemma B.4 with an r-dependent eigenvalue bound.

**Lemma B.9**: *Fix* $m \geq 1$. *There exists an* m*-dependent* $\kappa > 1$ *and given* $c_v \geq \kappa$*, there exists* $\kappa_{c_v} \geq \kappa$ *with the following significance: Suppose that* $r \geq \kappa_{c_v} c_v^{10}$ *and that* $(A, \psi = (\alpha, \beta))$ *is a solution to the* $(r, \mu)$ *version of (1.13) with* $\mu$ *a given element in* $\Omega$ *with* $\mathcal{P}$*-norm smaller than 1. Construct as in Lemma B.8 the family of operators* $\{\mathfrak{L}_{\mathbb{V},\tau}\}_{\tau \in [0,1]}$ *and introduce* $\{\lambda_{n\tau}\}_{\tau \in [0,1]}$ *to denote the resulting family of eigenvalues. If* $|\lambda_{n\tau}| = 0$ *for some* $\tau \in [0,1]$ *then* $|\lambda_{n\tau'}| \leq c_v^{-m} r^{1/2}$ *for all* $\tau' \in [0,1]$.

***Proof of Lemma B.9***: Return momentarily to the context in Part 1. Introduce T to denote $\sup_{\tau \in [0,1]} \|\frac{d}{d\tau} e_\tau\|$. It follows from (B.20) that any $n \in \mathbb{Z}$ version of the map $\tau \to \lambda_{n\tau}$ is such that $|\lambda_{n\tau'} - \lambda_{n\tau}| \leq T$ for any pair $\tau, \tau' \in [0,1]$. This implies in particular that $|\lambda_{n\tau}| > 0$ for all $\tau$ if $|\lambda_{n\tau'}| > T$ for any $\tau' \in [0,1]$. Given Lemma B.8, this last observation leads directly to the assertion in Lemma B.9 when applied to the family $\{\mathfrak{L}_{\mathbb{V},\tau}\}_{\tau \in [0,1]}$ with $T = c_v^{-m} r^{1/2}$.

*Part 3*: The three steps that follow complete the proof of Lemma B.4.

<u>Step 1</u>: If $m > c_0$, then Lemma B.9 can be invoked. With m so chosen, suppose that $\tau \in [0,1]$ and that $\lambda_{n\tau} = 0$. Let $\{\mathfrak{f}_{(\tau')}\}_{\tau' \in [0,1]}$ denote the corresponding family of eigenvectors. Fix $\tau' \in [0, 1]$. If $c_v \geq \kappa_\Diamond$ and if r is greater than a purely $c_v$-dependent constant, then Lemma B.9's bound on $|\lambda_{n\tau'}|$ implies that the assumptions of Lemmas A.2, A.3 and A.6 are met with $z = r$, with $(A_{*\tau'}, \psi_{*\tau'})$ used in lieu of $(A, \psi)$, with $\lambda = \lambda_{n\tau'}$, and with $\mathfrak{f} = \mathfrak{f}_{(\tau')}$. In particular, what is asserted by the first bullet of Lemma A.3 holds using $c_v$ in lieu of $c_0$. This is to say that when $\mathfrak{f}$ is written in terms of its $\mathbb{V}_0$ and $\mathbb{V}_1$ components as $(\mathfrak{f}_0, \mathfrak{f}_1)$, then the component $\mathfrak{f}_1$ has $L^2$ norm bounded by $c_0 c_v^k r^{-1/2}$ with $k \le c_0$.

Step 2: Assume that $c_v$ and r are chosen so that the assumptions of Lemmas B.1-B.2 and B.7-B.9 are met. Suppose that $\tau' \in [0,1]$ is a point where the map $\lambda_{n(\cdot)}$ is real analytic. Let $a'_v$ to denote endomorphism of $\mathbb{S}$ given by the derivative $(\frac{d}{d\tau}\nabla_{A_{*(\cdot)}})_v$ at the point $\tau'$. It follows from Lemmas B.2 and B.7 that $|a'_v| \le c_0 c_v^{-m}$. Write $\psi_{*(\cdot)}$ as $(\alpha_{*(\cdot)}, \beta_{*(\cdot)})$ and let $\beta'$ denote the derivative $\frac{d}{d\tau}\beta_{*(\cdot)}$ at $\tau'$. Lemma B.7 implies that $|\beta'| \le c_0 c_v^{-m} r^{-1/2}$. Meanwhile, Lemma B.8 has $|\frac{d}{d\tau} e_{(\cdot)}| \le c_0 c_v^{-m} r^{-1/2}$

Step 3: Look at (A.26) and (A.27) to see that the absolute value of the inner product between $\mathfrak{f}$ and $(\frac{d}{d\tau} e_{(\cdot)})|_{\tau'}\mathfrak{f}$ at any given point in Y is no greater than

$$c_0 c_v^{-m} r^{1/2} |\mathfrak{f}_0||\mathfrak{f}_1| + c_0(|a'_v| + r^{1/2}|\beta'|)|\mathfrak{f}|^2 . \tag{B.34}$$

The integral of the expression in (B.34) over Y is no less than $|\frac{d}{d\tau}\lambda_{n(\cdot)}|$ at $\tau'$, this being a consequence of (B.20). Meanwhile, what is said by Steps 1 and 2 imply that the integral of the expression in (B.34) is no greater than $c_0 c_v^{-m+k}$.

Use $m'$ to denote $m - k$. The argument used in the proof of Lemma B.9 proves that the bound by $c_0 c_v^{-m'}$ on $|\frac{d}{d\tau}\lambda_{n(\cdot)}|$ implies that $|\lambda_{n\tau}| = 0$ for some $\tau$ only if $|\lambda_{n(\tau')}| \le c_0 c_v^{-m'}$ for all $\tau' \in [0,1]$. Since $m'$ can be any positive number greater than $\kappa_\Diamond$, this last bound implies what is asserted by Lemma B.4.

**e) The pair $(A_\Diamond, \psi_\Diamond)$**

This subsection modifes $(A_*, \psi_*)$ on the components of $Y-(Y_{*\Lambda} \cup T_{*\Lambda})$ so that the resulting pair is given on these components by solutions to the vortex equations. The five parts of this subsection describe this modification.

*Part 1*: This first part describes the modification in the simplest case. To this end, fix attention on a component of $Y-(Y_{*\Lambda} \cup T_{*\Lambda})$ whose boundary is disjoint from the zero locus of $\alpha$. This is the simplifying assumption. Let $\gamma$ denote the curve in this component from $\cup_{p\in\Lambda}\{\hat{\gamma}_p^+ \cup \hat{\gamma}_p^-\}$ and let $T \subset Y$ denote the subset of points with distance $(c_v^4 + c_v^3)\, z^{-1/2}$ or less from $\gamma$. Fix coordinates for T from Part 4 of Section Aa with $\nu$ constant and $\mu$ both constant, real valued and positive.

The next lemma supplies a particular sort of isomorphism from $E|_T$ to the product bundle $T\times\mathbb{C}$.

**Lemma B.10**: *There exists* $\kappa > 1$ *and given* $c \geq \kappa$*, there exists* $\kappa_c \geq \kappa$ *with the following significance: Take* $r \geq \kappa_c c^{10}$ *and let* $(A, \psi = (\alpha, \beta))$ *denote a solution to the* $(r, \mu)$ *version of (1.13) with* $\mu$ *a given element in* $\Omega$ *with* $\mathcal{P}$*-norm smaller than 1. Suppose that* $\gamma \in \cup_{p\in\Lambda}\{\hat{\gamma}_p^+ \cup \hat{\gamma}_p^-\}$*, that* $\alpha$ *has zeros at distances less than* $\kappa r^{-1/2}$ *from* $\gamma$ *but no zeros at distance between* $\kappa r^{-1/2}$ *and* $(c^4+3c^3)r^{-1/2}$ *from* $\gamma$*. Fix coordinates* from *Part 4 of Section Aa for the radius* $(c^4+3c^3)r^{-1/2}$ *tubular neighborhood of* $\gamma$*. There is an isomorphism on the concentric, radius* $(c^4+c^3)r^{-1/2}$ *tubular neighborhood of* $\gamma$ *between* E *and the product bundle with the properties listed below.*

- *The isomorphism writes* $A = \theta + a_{A0}dt + \frac{1}{2}(\mathcal{A}\,d\bar{z} - \bar{\mathcal{A}}\,dz)$ *with* $|a_{A0}| \leq c_0$ *and* $|\mathcal{A}| \leq c_0 r^{1/2}$.
- *Use* m *to denote the sum of the local Euler numbers of the zeros of* $\alpha$ *on any radius* $c_v{}^4 r^{-1/2}$ *transverse disk centered on* $\gamma$*. The isomorphism writes* $\alpha$ *as* $|\alpha|(\frac{z}{|z|})^m$ *at points with distances between* $2\kappa r^{-1/2}$ *and* $(c^4+c^3)$ *from* $\gamma$*.*

This lemma is proved momentarily.

Take $c = c_v$ in Lemma B.10 and use the lemma's isomorphism between $E|_T$ and $T\times\mathbb{C}$ to write $A_\Diamond$ as $\theta + a_\Diamond$ with $a_\Diamond$ being in $i\mathbb{R}$ valued 1-form on T. Write $\psi_\Diamond$ as $(\alpha_\Diamond, \beta_\Diamond)$ and use the isomorphism to view $\alpha_\Diamond$ as a $\mathbb{C}$ valued function on T. Use the isomorphism and the chosen coordinate system to view $\beta_\Diamond$ as a $\mathbb{C}$-valued function also. Let m denote the rank of the complex bundle $\mathrm{Ker}(\vartheta)|_\gamma$. The data $a_\Diamond$, $\alpha_\Diamond$, and $\beta_\Diamond$ are given by what is written on the right hand side of the respective top, middle and bottom bullets in (A.44).

***Proof of Lemma B.10***: According to Proposition 2.4, the sum of the local Euler numbers of the zeros of $\alpha$ with distance at most $c_0 r^{-1/2}$ from $\gamma$ is bounded by $c_0$. Meanwhile, Lemma 2.3 finds $|\alpha| \geq 1 - \frac{1}{100}$ at distances greater than $c_0 r^{-1/2}$ from $\gamma$ but less than $(c^4+2c^3)r^{-1/2}$ from $\gamma$. Granted the $\alpha$'s local Euler number is m, and granted this lower bound on $|\alpha|$, there exists $c_1 \leq c_0$ and an isomorphism from E on the $|z| \geq c_1 r^{-1/2}$ part of the radius $(c^4+2c^3)r^{-1/2}$ tubular neighborhood of $\gamma$ that writes $\alpha$ as $|\alpha|(\frac{z}{|z|})^m$. Use this isomorphism to write A as in the first bullet of the lemma. The isomorphism writes the dt part of $\nabla_A\alpha$ as $(\partial_t|\alpha| + a_{A0}|\alpha|)(\frac{z}{|z|})^m$. Given that $D_A\psi = 0$, it follows that the dt part of $\nabla_A\alpha$ is bounded by $|\nabla_A\beta| + c_0 r^{-1/2}|\nabla_A\alpha|$. This understood, then Lemma 2.1's bound imply that $|a_{A0}| \leq c_0$ on the $|z| \geq c_1 r^{-1/2}$ part of the radius $(c^4+2c^3)r^{-1/2}$ tubular neighborhood of $\gamma$. The bound for $|\mathcal{A}|$ on this same part of the radius $(c^4+2c^3)r^{-1/2}$ tubular neighborhood of $\gamma$ follows from Lemma 2.1's bound for $|\nabla_A\alpha|$.

The isomorphism just described will be modified on the $|z| \leq \frac{7}{4}c_1 r^{-1/2}$ tubular neighborhood of $\gamma$ to obtain an isomorphism between E and the product bundle on the whole $(c^4+2c^3)\,r^{-1/2}$ tubular neighborhood of $\gamma$ that obeys the first bullet of the lemma. To this end, remark first that there is an isomorphism between $E|_\gamma$ and $\gamma\times\mathbb{C}$ that writes the pull-back of A along $\gamma$ as $\hat{a}_{A2}dt$ with $\hat{a}_{A2}$ constant with absolute value less than $\frac{2\pi}{\ell_\gamma}$. Fix such an isomorphism, and then use parallel transport by A along the rays through the origin in the constant t slices of the tubular neighborhood to extend this isomorphism to the $|z| \geq 2c_1 r^{-1/2}$ part of the tubular neighborhood. An isomorphism of this sort writes A as $\theta + a_{A2}dt + \frac{1}{2}\hat{A}(z\,d\bar{z} - \bar{z}\,dz)$ where $\hat{A}$ is a real valued function on the radius $2c_1 r^{-1/2}$ tubular neighborhood of $\gamma$. This function obeys the equation $\rho^{-1}\frac{\partial}{\partial\rho}(\rho^2\hat{A}) = F_A(\frac{\partial}{\partial z}, \frac{\partial}{\partial\bar{z}})$ with $\rho$ denoting $|z|$. Integrate this identity starting at $|z| = 0$ using (1.13) to see that $\rho|\hat{A}| \leq c_0 r^{1/2}$. Meanwhile, $\frac{\partial}{\partial\rho}a_{A2} = F_A(\frac{\partial}{\partial\rho}, \frac{\partial}{\partial t})$ because $A - \theta$ has no $d\rho$ component. Integrate the latter identity using (1.13) with the fact $|\partial_t - v| \leq c_0|z|$ to see that $|a_{A2}| \leq c_0$ where $|z| \leq 2c_1 r^{-1/2}$.

The preceding paragraphs describe two isomorphisms between E and the product bundle that are both defined on the $c_1 r^{-1/2} \leq |z| \leq 2c_1 r^{-1/2}$ part of the radius $(c^4+2c^3)\,r^{-1/2}$ tubular neighborhood of $\gamma$. The corresponding transition function is a map from this part of the tubular neighborhood to $S^1$. Use $\hat{u}$ to denote this map. The bounds on $a_{A1}$ and $a_{A2}$ imply that $|\frac{\partial}{\partial t}\hat{u}| \leq c_0$ and those on A and $\rho|\hat{A}|$ imply that $|\frac{\partial}{\partial z}\hat{u}| \leq c_0 r^{1/2}$. This understood, the map $\hat{u}$ with a cut-off function constructed from $\chi$ can be used in a straighforward manner to modify the outer isomorphism where $\frac{5}{4}c_1 r^{-1/2} < |z| \leq \frac{7}{4}c_1 r^{-1/2}$ so that the result agrees with the inner isomorphism where $c_1 r^{-1/2} < |z| < \frac{5}{4}c_1 r^{-1/2}$ and is such that the norms of the new versions of $a_{A1}$ and A are still bounded by $c_0$ and $c_0 r^{1/2}$ respectively.

*Part 2*: Fix attention on a component of $Y-(Y_{*\Lambda}\cup T_{*\Lambda})$ whose boundary intersects the zero locus of $\alpha$ and let $\gamma$ again denote the corresponding curve from the set $\cup_{p\in\Lambda}\{\hat{\gamma}_p^+ \cup \hat{\gamma}_p^-\}$. Introduce the coordinates from Part 4 of Section Aa for $\gamma$ with $\nu$ constant and with $\mu$ constant, real and positive. It follows as a consequence of what is said in (B.1) that the zero locus of $\alpha$ extends as two properly embedded arcs in the part of the radius $c_v^4 r^{-1/2}$ tubular neighborhood of $\gamma$ where $(c_v^4 - 3c_v^3)\,r^{-1/2} \leq |z| \leq c_v^4 r^{-1/2}$. The following lemma describes an extension of these two arcs as the end segments of a single, properly embedded arc in the radius $c_v^4 r^{-1/2}$ tubular neighborhood of $\gamma$.

**Lemma B.11**: *There exists $\kappa > 1$ and given $c_v \geq \kappa$, there exists $\kappa_{c_v} \geq \kappa$ with the following significance: Suppose that $r \geq \kappa_{c_v} c_v^{10}$ and suppose that $(A, \psi = (\alpha, \beta))$ is a solution to the $(r, \mu)$ version of (1.13) with $\mu$ a given element in $\Omega$ with $\mathcal{P}$-norm smaller than 1. Let $\gamma$ denote a curve from $\cup_{p\in\Lambda}\{\hat{\gamma}_p^+ \cup \hat{\gamma}_p^-\}$ with points at distance $c_v^4\, r^{-1/2}$ from the zero locus of $\alpha$. Let T denote the set of points with distance $(c_v^4 + c_v^3)\,r^{-1/2}$ or less from $\gamma$. There exists a smooth, properly embedded arc in T with the properties listed below.*

- *The arc is the zero locus of* $\alpha$ *at points with distance* $c_v r^{-1/2}$ *or more from* $\gamma$.
- *The arc lies in the* $1 - 3\cos^2\theta > 0$ *part of* T.
- *Each point in the arc has distance greater than* $\kappa^{-1} c_v r^{-1/2}$ *from* $\gamma$.
- *A unit length tangent vector to the arc has distance at most* $\kappa c_v r^{-1/2}$ *from* $v$.
- *Fix a closed, transverse disk in* T *with radius* $c_v^{\,4} r^{-1/2}$, *center point on* $\gamma$, *and no zeros of* $\alpha$ *on its boundary. Let* $m_*$ *denote the intersection number between this disk and the arc and let* $m_\alpha$ *denote the sum of the local Euler numbers of the zeros of* $\alpha$ *on the disk. Then* $m_\alpha - m_*$ *is non-negative and independent of the chosen disk.*

***Proof of Lemma B.11***: There are various ways to construct an arc with the desired properties. The construction that follows is perhaps more complicated than is needed for now, but resulting arc has certain extra properties that are exploited later. The construction has four steps.

Step 1: Fix $m \in (1, c_0^{-1} c_v)$. If $\alpha^{-1}(0) \cap T$ has a component whose points all have distance greater than $m^{-4} c_v r^{-1/2}$ from $\gamma$, then the desired arc is this component of $\alpha^{-1}(0) \cap T$. Proposition 2.4 guarantees that the conditions of the lemma are met if $\kappa$ is greater than a purely $m$-dependent constant. The remaining steps assume that the components of $\alpha^{-1}(0)$ that intersects the boundary of T have points with distance $m^{-4} c_v r^{-1/2}$ from $\gamma$. This understood, keep in mind the following consequence of Proposition 2.4: The zero locus of $\alpha$ in the part of T with distance at least $m^{-4} c_v r^{-1/2}$ from $\gamma$ consists of two arcs, one where u is everywhere positive and the other where u is every where negative. These are denoted respectively by $\upsilon_+$ and $\upsilon_-$ in what follows. Note also that the unit tangent vector to either arc has distance at most $c_0 r^{-1/2}$ from $v$.

The subsequent three steps take $\cos\theta = \frac{1}{\sqrt{3}}$ on $\gamma$. The construction for the case when $\cos\theta = -\frac{1}{\sqrt{3}}$ on $\gamma$ is identical but for certain sign changes and will not be given.

Step 2: Let $\theta_* \in (0, \pi)$ denote the angle that obeys $\cos\theta_* = \frac{1}{\sqrt{3}}$. It follows from the formula for $v$ in (1.3) that any integral curve of $v$ in T can be parametrized as a map from an interval $I \subset \mathbb{R}$ to $(\mathbb{R}/(2\pi\mathbb{Z})) \times \mathbb{R}^2$ of the form $t \to (\phi = -t, u = b\, x_0(t), \theta = \theta_* + y_0(t))$ where $b = \frac{\sqrt{3}}{2\sqrt{2}} e^R (x_0 + 4e^{-2R})^{1/2}$ and where x and y are smooth functions that obey

$$\frac{d}{dt} x_0 = \lambda y_0 + \mathfrak{e}_x(x_0, y_0) \quad \textit{and} \quad \frac{d}{dt} y_0 = \lambda x_0 + \mathfrak{e}_y(x_0, y_0) \tag{B.35}$$

with $\lambda = 4\sqrt{6}\, e^{-R} (x_0 + 4e^{-2R})^{1/2}$ and with the pair $\mathfrak{e}_x$ and $\mathfrak{e}_y$ being smooth functions of the coordinates (x, y) on $\mathbb{R}^2$ that obey $|\mathfrak{e}_x| + |\mathfrak{e}_y| \le c_0(x^2 + y^2)$.

This parametrization of integral curves of $\nu$ in T suggests the introduction of coordinates (x, y) for T by writing $u = bx$ and $y = \theta_* + y$. These coordinates are such that if $m \geq c_0$ and if $p \in [0, 4)$, then the points in T with $(x^2 + y^2)^{1/2} = m^{-p} c_v r^{-1/2}$ have distance less the $c_0 m^{-p} c_v r^{-1/2}$ from $\gamma$ and distance greater than $c_0^{-1} m^{-p} c_v r^{-1/2}$ from $\gamma$.

It follows from what is said in Proposition 2.4 that $\upsilon_+$ can be parametrized as a map from an interval $I_+ \subset \mathbb{R}$ to $(\mathbb{R}/(2\pi\mathbb{Z})) \times \mathbb{R}^2$ of the form

$$t \to (\varphi = -t,\ u = b\, x_+(t)\,,\ \theta = \theta_* + y_+(t)) \tag{B.36}$$

where $x_+$ and $y_+$ obey a modified version of (B.35) that adds respective terms $\mathfrak{r}_{x+}$ and $\mathfrak{r}_{y+}$ of t to the left and right hand equations. These are smooth functions of t with absolute value bounded by $c_0 r^{-1/2}$. The arc $\upsilon_-$ can be parametrized in a similar fashion as a map with domain an interval $I_- \subset \mathbb{R}$ by a pair of functions $(x_-, y_-)$ that obey a modified version of (B.35) that adds respective terms $r_{x-}$ and $r_{y-}$ to the left hand equations, both being functions of t with absolute value bounded by $c_0 r^{-1/2}$.

Use S to denote this torus $(x^2 + y^2)^{1/2} = m^{-2} c_v r^{-1/2}$. This torus is intersected transversely in one point by $\upsilon_+$ and likewise by $\upsilon_-$. No generality is lost by parametrizing $I_+$ and $I_-$ so these points $S \cap \upsilon_+$ and $S \cap \upsilon_-$ occur at respective parameter values $t = 2\pi + t_\lozenge$ and $t = -2\pi - t_\lozenge$ with $t_\lozenge \in [0, 2\pi)$. More is said about $I_+$ and $I_-$ momentarily.

<u>Step</u> <u>3</u>: Introduce coordinates p and q on $\mathbb{R}^2$ by the rules $p = y + x$ and $q = y - x$. Writing (B.35) or its $\upsilon_+$ or $\upsilon_-$ analogs in terms p and q gives an equation for a pair of maps $t \to p_*(t)$ and $t \to q_*(t)$. Here $*$ denotes either 0, + or - as the case may be. The equation in question has the form

$$\tfrac{d}{dt} p_* = \lambda p_* + \mathfrak{e}_p(p_*, q_*) + \mathfrak{r}_{p*} \quad \textit{and} \quad \tfrac{d}{dt} q_* = -\lambda q_* + \mathfrak{e}_q(p_*, q_*) + \mathfrak{r}_{q*}\ , \tag{B.37}$$

where $\mathfrak{e}_p$ and $\mathfrak{e}_q$ are smooth functions of p and q that obey $|\mathfrak{e}_p| + |\mathfrak{e}_q| \leq c_0 (p^2 + q^2)$; and where $\mathfrak{r}_{p0}$ and $\mathfrak{r}_{q0}$ are zero (this the case of (B.35)), while $\mathfrak{r}_{p+}$, $\mathfrak{r}_{q+}$, $\mathfrak{r}_{p-}$ and $\mathfrak{r}_{q-}$ are functions of t with absolute value bounded by $c_0 r^{-1/2}$.

It follows from (B.37) that $p_+$ and $q_+$ where $(p_+^2 + q_+^2)^{1/2} \leq c_v r^{-1/2}$ have the form

$$p_+(t) = p_{\lozenge+} e^{\lambda(t - 2\pi - t_\lozenge)} + \mathfrak{w}_{p+} \quad \textit{and} \quad q_+(t) = q_{\lozenge+} e^{-\lambda(t - 2\pi - t_\lozenge)} + \mathfrak{w}_{q+}, \tag{B.38}$$

where $p_{\lozenge+} = p(2\pi + t_\lozenge)$ and $q_{\lozenge+} = q(2\pi + t_\lozenge)$, and where $\mathfrak{w}_{p+}$ and $\mathfrak{w}_{q+}$ are functions of t with absolute value bounded by $c_0 r^{-1/2}$ where $|t| \leq c_0$. The t-derivatives of $\mathfrak{w}_{p+}$ and $\mathfrak{w}_{q-}$ for such t are also bounded by $c_0 r^{-1/2}$.

If $m > c_0$, then the fact that $\upsilon_+$ intersects the locus where $(p^2 + q^2)^{1/2} = c_0 m^{-4} c_v r^{-1/2}$ leads to the following observations:

- *The interval* $[-2\pi - t_\Diamond, 2\pi + t_\Diamond] \subset I_+$.
- $|p_{\Diamond+} - m^{-2} c_v r^{-1/2}| + m^2 |q_{\Diamond+}| \le c_0 m^{-4} c_v r^{-1/2}$.

(B.39)

To see why this is, note first that the torus S is the locus $(p^2 + q^2)^{1/2} = 2^{1/2} m^{-2} c_v r^{-1/2}$ and so neither $p_{\Diamond+}$ nor $q_{\Diamond+}$ is greater than $2^{1/2} m^{-2} c_v r^{-1/2}$. However $p_{\Diamond+} \ge 2^{-1/2} m^{-2} c_v r^{-1/2}$ because $p \ge |q|$ where $u > 0$. Granted that $p_{\Diamond+} \ge 2^{-1/2} m^{-2} c_v r^{-1/2}$ and granted that $\upsilon_+$ intersects the locus where $(p^2 + q^2)^{1/2} = c_0 m^{-4} c_v r^{-1/2}$, then the left hand identity in (B.37) requires that the parameter t at this intersection is less than $\lambda^{-1} \ln(m^{-2}) + c_0$. This gives the top bullet in (B.39) if $m > c_0$. The fact that t on $I_+$ has values less than $\lambda^{-1} \ln(m^{-2}) + c_0$ with the right hand identity in (B.38) finds $q_+$ at such values of t greater than $q_{\Diamond+} c_0 m^2$ and so $m^2 |q_+|$ must be less than $c_0 m^{-4} c_v r^{-1/2}$. This implies the lower bullet in (B.39).

What was just said about $p_+$ and $q_+$ has its $p_-$ and $q_-$ analog. By way of a summary, these functions can be written as

$$p_-(t) = p_{\Diamond-} e^{\lambda(t + 2\pi + t_\Diamond)} + \varpi_{p-} \quad \textit{and} \quad q_-(t) = q_{\Diamond-} e^{-\lambda(t + 2\pi + t_\Diamond)} + \varpi_{q-}, \tag{B.40}$$

where $p_{\Diamond-} = p_-(-2\pi - t_\Diamond)$ and $q_{\Diamond-} = q(-2\pi - t_\Diamond)$, and where $\varpi_{p-}$ and $\varpi_{q-}$ are functions of t with absolute value bounded by $c_0 r^{-1/2}$ where $|t| \le 100\pi$. Their respective t-derivatives are also bounded by $c_0 r^{-1/2}$ for such t. The $(p_-, q_-)$ analog of (B.39) reads

- *The interval* $[-2\pi - t_\Diamond, 2\pi + t_\Diamond] \subset I_-$.
- $m^2 |p_{\Diamond-}| + |q_{\Diamond-} - m^{-2} c_v r^{-1/2}| \le c_0 m^{-4} c_v r^{-1/2}$.

(B.41)

The proof of (B.41) differs only cosmetically from that of (B.39).

<u>Step</u> <u>4</u>: The arc $\upsilon$ coincides with the $t > 2\pi + t_\Diamond$ part of $\upsilon_+$ and the $t < -2\pi - t_\Diamond$ part of $\upsilon_-$. The remaining part of $\upsilon$ is parametrized by $[-2\pi - t_\Diamond, 2\pi + t_\Diamond]$. The definition that follows for this part of $\upsilon$ refers to a non-negative function on $\mathbb{R}$ that is denoted by $\sigma$ and defined by the rule $t \to \sigma(t) = \chi(1 - \frac{4}{\pi} t)$. This function is equal to 1 where $t < \frac{\pi}{4}$ and it is equal to 0 where $t > \frac{\pi}{2}$.

The $t \in [-2\pi - t_\Diamond, 2\pi + t_\Diamond]$ point of $\upsilon$ is written as $(\varphi = -t, u = b\, x_\upsilon(t)\ \theta = \theta_* + y_\upsilon(t))$ with $x_\upsilon$ and $y_\upsilon$ being functions on the interval $[-2\pi - t_\Diamond, 2\pi + t_\Diamond]$. The functions $x_\upsilon$ and $y_\upsilon$ are written here as $x_\upsilon = \frac{1}{2}(p_\upsilon - q_\upsilon)$ and $y_\upsilon = \frac{1}{2}(p_\upsilon + q_\upsilon)$ with $p_\upsilon$ and $q_\upsilon$ given by the following rule:

- $p_\upsilon(t) = \sigma(-t)(p_{\lozenge+}\, e^{\lambda(t - 2\pi - t_\lozenge)} + \varpi_{p+}(t)) + \sigma(t)\, p_-(t)$ ,
- $q_\upsilon(t) = \sigma(t)(q_{\lozenge-}\, e^{-\lambda(2\pi + t_\lozenge + t)} + \varpi_{q-}(t)) + \sigma(-t)\, q_+(t)$ .

(B.42)

It is a straightforward matter to check that the arc $\upsilon$ has all of the required properties.

*Part 3*: Use $\kappa_\lozenge$ to denote a constant that is greater than the versions of $\kappa$ that appear in Lemmas A.2-A.9 and in Lemmas B.10 and B.11. Assume in what follows that $c_v \ge \kappa_\lozenge$ and that r is greater than the $c_0 = c_v$ lower bounds in given in Lemmas A.2-A.9 and the lower bounds given in Lemmas B.10 and B.11. Fix a component of $Y-(Y_{*\Lambda} \cup T_{*\Lambda})$ whose boundary has a zero of $\alpha$. Let $\gamma$ denote the nearby curve from $\cup_{p\in\Lambda}\{\hat{\gamma}_p^+ \cup \hat{\gamma}_p^-\}$ and let T denote the set of points in Y with distance $(c_v^4 + c_v^3)\, r^{-1/2}$ or less $\gamma$. This part of the subsection defines $(A_\lozenge, \psi_\lozenge)$ on T. The definition has four steps. These steps use $\upsilon$ denote the arc that is supplied by Lemma B.11. These steps also use $\gamma$'s version of the coordinates (t, z) from Part 4 of Section Aa for T.

Step 1: Define $U_{T0} \subset T$ as follows: The $|z| \ge (c_v^4 - 2c_v^2)\, r^{-1/2}$ part of $U_{T0}$ consists of the points in T with distance greater than $c_v^2 r^{-1/2}$ from $\upsilon$. The $|z| < (c_v^4 - 2c_v^2)\, r^{-1/2}$ part of $U_{T0}$ consists of the points with distance greater than $c_v^{1/4} r^{-1/2}$ from $\upsilon$ and with distance greater than $c_v^{1/4} r^{-1/2}$ from $\gamma$. Note that the $|z| \ge c_v^4\, r^{-1/2}$ part of $U_{T0}$ coincides with the $Y_{*\Lambda} \cap T$ part of Section Ba's set $U_0$. This understood, fix an isomorphism over $U_{T0}$ between E and $U_{T0} \times \mathbb{C}$ that sends $\alpha$ to $|\alpha|$ on the part of $U_{T0}$ where $|\alpha| \ge \frac{1}{2}$ . Such an isomorphism extends Section Ba's isomophism from the $|z| \ge (c_v^4 - 2c_v^2)\, r^{-1/2}$ part of $U_{T0}$ to the whole of $U_{T0}$. This isomorphism identifies $A_\lozenge$ on $U_{T0}$ with the product connection and it identifies $\alpha_\lozenge$ with the constant $1 \in \mathbb{C}$. The component $\beta_\lozenge$ is everywhere zero on $U_{T0}$.

What follows is, for now, just a parenthetical remark: Suppose that $\lambda$ and $\lambda'$ are two isomorphisms from $E|_{U_{T0}}$ to $U_{T0} \times \mathbb{C}$ that agree where $|\alpha| > \frac{1}{2}$ . Then $\lambda' = e^{ix}\lambda$ with x being a real valued function which is 0 where $|\alpha| > \frac{1}{2}$ . That this is so is a consequence of what is said in Proposition 2.4 about the zero locus of $\alpha$ in T.

Step 2: Let $U_{\upsilon-} \subset T$ denote the subset of points with $|z| < (c_v^4 - \frac{7}{4} c_v^3)\, r^{-1/2}$ and distance less than $4c_v^{1/4} r^{-1/2}$ from $\upsilon$. To this end, keep in mind that the $|z| \ge (c_v^4 - 2c_v^3)\, r^{-1/2}$ part of $\upsilon$ coincides with this same part of $\alpha^{-1}(0) \cap T$. This understood, fix coordinates for $U_{\upsilon-}$ from Part 4 of Section Aa and coincide on the $|z| \ge (c_v^4 - 2c_v^3)\, r^{-1/2}$ part of $U_{\upsilon-}$ with those used in Section Ba. Denote these coordinates by $(t_\upsilon, z_\upsilon)$ so as to distinguish them from the coordinates t and z that are used for T. The restriction of E to the $|z| \ge (c_v^4 - 2c_v^3)\, r^{-1/2}$ part of $U_{\upsilon-}$ has its $\alpha \to |\alpha|\frac{z_\upsilon}{|z_\upsilon|}$ isomorphism with the product bundle. Extend this to the product bundle so as to give an isomorphism over the whole of $U_{\upsilon-}$ between E and the product bundle. This extension should be such that the corresponding transition function for Step 1's isomorphism from $E|_{U_{T0}}$ to $U_{T0}\times \mathbb{C}$ sends the constant section 1 of $U_{T0}\times \mathbb{C}$ to the section $\frac{z_\upsilon}{|z_\upsilon|}$ of $U_{\upsilon-}\times \mathbb{C}$ on the $|z_\upsilon| \ge c_v^{1/4} r^{-1/2}$ part of $U_{\upsilon-}$. The formula in the next equation defines $A_\Diamond$ and $\alpha_\Diamond$ on $U_{\upsilon-}$ by viewing them via this isomorphism as a connection on the product bundle and map to $\mathbb{C}$. This isomorphism with the coordinates $(t_\upsilon, z_\upsilon)$ are used to view $\beta_\Diamond$ as a map to $\mathbb{C}$ also.

The upcoming equation uses $\chi_{\hat{U}\upsilon-}$ to denote the function of $|z_\upsilon|$ given by the rule $\chi(c_v^{-1/4} r^{1/2}|z_\upsilon| - 1)$. The equation also uses $(\nu_\upsilon, \mu_\upsilon)$ to denote the $(t_\upsilon, z_\upsilon)$ version of the functions $\nu$ and $\mu$ from $\upsilon$'s version of (A.6). This equation once again brings in the functions $\alpha_0$ and $a_0$ from (A.3) and the corresponding versions of y and $\varsigma$ from (A.2).

- $A_\Diamond = \theta + \nu_\upsilon\, \chi_{\hat{U}\upsilon-}\, i\, 2^{1/2} r_r^* y\, dt_\upsilon - \frac{1}{2}(1 - \chi_{\hat{U}\upsilon-} + \chi_{\hat{U}\upsilon-} r_r^* a_0)(z_\upsilon^{-1} dz_\upsilon - \bar{z}_\upsilon^{-1} d\bar{z}_\upsilon)$,
- $\alpha_\Diamond = (1 - \chi_{\hat{U}\upsilon-}(1 - r_r^*|\alpha_0|))\frac{z_\upsilon}{|z_\upsilon|}$ ,
- $\beta_\Diamond = i\, \mu_\upsilon r^{-1/2} \chi_{\hat{U}\upsilon-} r_r^* \varsigma$ .

(B.43)

<u>Step 3</u>: This step defines $A_\Diamond$, $\alpha_\Diamond$ and $\beta_\Diamond$ on the part of T with $|z| \ge (c_v^4 - 2c_v^3) r^{-1/2}$ where the distance to $\upsilon$ is less than $4c_v^2 r^{-1/2}$. To this end, use the $\alpha \to \alpha\frac{z_\upsilon}{|z_\upsilon|}$ isomorphism between E and the product bundle over this part of T to view $A_\Diamond$ as a connection on the product bundle and $\alpha_\Diamond$ as a map to $\mathbb{C}$. Use this same isomorphism with the coordinates $(t_\upsilon, z_\upsilon)$ to view $\beta_\Diamond$ as a map to $\mathbb{C}$ also.

Reintroduce $\chi_{\Diamond\Diamond}$ to denote the function that appears in (B.2). This function equals 1 where $|z| < (c_v^4 - \frac{7}{4}c_v^3) r^{-1/2}$ and it equals 0 where $|z| > (c_v^4 - \frac{5}{4}c_v^3) r^{-1/2}$. The definition uses $\chi_{\hat{U}\upsilon+}$ to denote the function of $|z_\upsilon|$ given by $(1 - \chi_{\Diamond\Diamond})\, \chi(c_v^{-2} r^{1/2}|z_\upsilon| - 1) + \chi_{\Diamond\Diamond}\chi(c_v^{-1/4} r^{1/2}|z_\upsilon| - 1)$. This function is equal to (B.35)'s function $\chi_{\hat{U}\upsilon-}$ where $|z| \le (c_v^4 - \frac{7}{4}c_v^3) r^{-1/2}$ and it is equal to the (B.10)'s function $\chi_{\hat{U}}$ where $|z| \ge (c_v^4 - \frac{5}{4}c_v^3) r^{-1/2}$.

Replace $\chi_{\hat{U}\upsilon-}$ in (B.43) with $\chi_{\hat{U}\upsilon+}$ to obtain the formulas for $A_\Diamond$, $\alpha_\Diamond$ and $\beta_\Diamond$ on the part of T with $|z| \ge (c_v^4 - 2c_v^3) r^{-1/2}$ and with distance less than $4c_v^2 r^{-1/2}$ to $\upsilon$.

<u>Step 4</u>: This last step defines $A_\Diamond$, $\alpha_\Diamond$ and $\beta_\Diamond$ on the $|z| < \frac{3}{4} c_v^{1/2} r^{-1/2}$ part of T. The definition requires Lemma B.11's integer m. The definition also requires the choice of an isomorphism between E on this part of T and the product bundle. A choice for such an isomorphism should be made subject to the following constraint: The resulting transition

function on the $\frac{1}{2}\,c_v^{1/2}r^{-1/2} < |z| < \frac{3}{4}\,c_v^{1/2}r^{-1/2}$ part of T between the product product bundle over $|z| < \frac{3}{4}\,c_v^{1/2}r^{-1/2}$ part of T and the product bundle $U_{T0}\times\mathbb{C}$ sends the latter's constant section 1 to the former's section $z \to (\frac{z}{|z|})^m$. An isomorphism of this sort exist because γ represents the class 0 in $H_1(Y;\mathbb{Z})$. Moreover, the space of isomorphisms that obey this constraint is contractible. The chosen isomorphism is used to view $A_\Diamond$ as a connection on the product bundle over this part of T and $\alpha_\Diamond$ here as a $\mathbb{C}$ valued function. This isomorphism with the coordinates (t, z) is used to view $\beta_\Diamond$ as a $\mathbb{C}$ valued function as well.

In the case m = 0, the transition function is the constant function. In this case, $A_\Diamond$ on the $|z| < \frac{3}{4}\,c_v^{1/2}r^{-1/2}$ part of T is set equal to the product connection, the function $\alpha_\Diamond$ is the constant function 1, and $\beta_\Diamond$ is zero.

Assume next that m > 0. Introduce $\chi_{**}$ the function $z \to \chi(4c_v^{-1/2}r^{-1/2}|z| - 1)$. This function equals 1 where $|z| \le \frac{1}{4}\,c_v^{1/2}r^{-1/2}$ and it equals 0 where $|z| \ge \frac{1}{2}\,c_v^{1/2}r^{-1/2}$. The definition of $A_\Diamond$, $\alpha_\Diamond$ and $\beta_\Diamond$ in the next equation uses $\chi_{**}$; it uses ν and μ to denote the pair of functions that are given by γ's version of (A.6); and uses the functions $\alpha_{m0}$, $a_{m0}$, $y_m$ and $\varsigma_m$ that appear in (A.44). What follows defines $A_\Diamond$, $\alpha_\Diamond$ and $\beta_\Diamond$ on $|z| < \frac{3}{4}\,c_v^{1/2}r^{-1/2}$ part of T.

- $A_\Diamond = \theta + \nu\,\chi_{**}\,i\,2^{1/2}\,r_r^*y_m dt - \frac{m}{2}(1 - \chi_{**} + \chi_{**}\,r_r^*a_{m0})(z^{-1}dz - \bar{z}^{-1}d\bar{z})$,
- $\alpha_\Diamond = (1 - \chi_{**}(1 - r_r^*|\alpha_{m0}|)(\frac{z}{|z|})^m$,
- $\beta_\Diamond = i\,\mu\,r^{-1/2}\chi_{**}\,r_r^*\varsigma_m$.

(B.44)

*Part 4*: The next lemma asserts two important features of the large $c_v$ and r versions of $(A_\Diamond,\psi_\Diamond)$.

**Lemma B.12**: *There exists* $\kappa \ge 100$, *and given* $c_v \ge \kappa$, *there exists* $\kappa_{c_v} > \kappa$ *with the following significance: Suppose that* $r \ge \kappa_{c_v} c_v^{10}$ *and suppose that* $(A, \psi = (\alpha, \beta))$ *is a solution to the* $(r, \mu)$ *version of (1.13) with* $\mu$ *a given element in* Ω *with* $\mathcal{P}$*-norm smaller than 1.*

- *The corresponding* $(A_\Diamond,\psi_\Diamond)$ *does not depend on the coordinates from Part 4 of Section Aa that are chosen from the various* $\gamma \in \Theta$ *versions of* $U_\gamma$.
- *The corresponding* $(A_\Diamond,\psi_\Diamond)$ *satisfies the* $c_0 = c_v$ and $z = r$ *versions of* PROPERTIES 1-5 *in Section Ab.*

***Proof of Lemma B.12***: The fact that $(A_\Diamond,\psi_\Diamond)$ does not depend on the chosen coordinates from Part 4 of Section Aa follows directly from the fact that $(A_*,\psi_*)$ does not depend on these choices. The assertion in the second bullet follows from Lemma A.1 if PROPERTIES 1, 2, 4 and 5 are obeyed on each component of $Y-(Y_{*\Lambda} \cup T_{*\Lambda})$. The fact that these properties are obeyed follows from (3.3) and the fact that $y_m$ and $\varsigma_m$ and their derivatives obey similar bounds.

**f) A path from $(A_*, \psi_*)$ to $(A_\Diamond, \psi_\Diamond)$**

This subsection derives an $(A, \psi)$ and r independent bound for the norm of the difference between the values of $\mathfrak{f}_\mathfrak{s}$ at $(A_*, \psi_*)$ and $(A_\Diamond, \psi_\Diamond)$. The proposition that follows makes the formal assertion.

**Proposition B.13**: *There exists* $\kappa \geq 100$*, and given* $c_v \geq \kappa$*, there exists* $\kappa_{c_v} > \kappa$ *with the following significance: Suppose that* $r \geq \kappa_{c_v} c_v^{10}$ *and suppose that* $(A, \psi = (\alpha, \beta))$ *is a solution to the* $(r, \mu)$ *version of (1.13) with* $\mu$ *a given element in* $\Omega$ *with* $\mathcal{P}$*-norm smaller than 1. Then the norm of the difference between the respective values of* $\mathfrak{f}_\mathfrak{s}$ *at* $(A_*, \psi_*)$ *and* $\mathfrak{f}_\mathfrak{s}$ *at* $(A_\Diamond, \psi_\Diamond)$ *is bounded by* $\kappa$.

The proof of this proposition is given in Part 8 of the subsection. The intervening parts define a certain path in $\mathrm{Conn}(E) \times C^\infty(Y; \mathbb{S})$ from $(A_*, \psi_*)$ to $(A_\Diamond, \psi_\Diamond)$ that is used in the proof.

*Part 1*: The path is parametrized by [0, 1] and a given $\tau \in [0,1]$ member is denoted by $(A_{\Diamond\tau}, \psi_{\Diamond\tau})$ with $\tau = 0$ member $(A_*, \psi_*)$ and $\tau = 1$ member $(A_\Diamond, \psi_\Diamond)$. As defined, the pairs $(A_*, \psi_*)$ and $(A_\Diamond, \psi_\Diamond)$ agree on $Y_{*\Lambda} \cup T_{*\Lambda}$ and this will be the case for all pairs along the path between them. The definition of the path $\{(A_{\Diamond\tau}, \psi_{\Diamond\tau})\}_{\tau \in [0,1]}$ on a given component of $Y-(Y_{*\Lambda} \cup T_{*\Lambda})$ is supplied momentarily. The definition uses $\gamma$ to denote the corresponding curve from the set $\cup_{p \in \Lambda} \{\hat{\gamma}_p^+ \cup \hat{\gamma}_p^-\}$, and it uses T to denote the set of points with distance $(c_v^4 + c_v^3) r^{-1/2}$ or less from $\gamma$. The definition also uses $\gamma$'s version of the coordinates from Part 4 of Section Aa for T that has $\nu$ constant and $\mu$ constant, real and greater than $|\nu|$.

By way of an overview of what is to come, the path $\tau \to (A_{\Diamond\tau}, \psi_{\Diamond\tau})$ first moves $(A_*, \psi_*)|_T$ to a pair with two salient features: It is very close to $(A_*, \psi_*)$ in a large k version of the $C^k$ topology; and it is constructed from a 1-parameter family of vortex solutions to (2.8) with the parameter being the points in $\gamma$. The parametrization is such that the pull-back via the scaling map $z \to r^{1/2} z$ of a given $t \in \gamma$ solution to the vortex equation defines the restriction of the new pair to the constant t slice of T. A homotopy of this $\gamma$-parameterized family through $\gamma$-parametrized families of solutions to (2.8) is used to define the second part of the path $\tau \to (A_{\Diamond\tau}, \psi_{\Diamond\tau})$. The end member of this second part of the path is very close to $(A_\Diamond, \psi_\Diamond)$ in a large k version of the $C^k$ topology. The third part of the path moves this end pair to $(A_\Diamond, \psi_\Diamond)$.

*Part 2*: Suppose for the moment $\alpha$ has no zeros on the boundary of the closure of a given transverse disk in T with center on $\gamma$. If this is the case, then the sum of the local Euler numbers of the zeros of $\alpha$ can be defined, and this sum is a positive integer. If $\alpha$ has no zeros on T's boundary torus, then this sum is the same for all transverse disks of this sort. If $\alpha$ has zeros on this torus, then there are two values that occur unless both zeros of $\alpha$ on the boundary of T have the same value of the parameter t. These two values differ by 1. In any event, use $m_\alpha$ to denote larger of the possible values for the sum of the local Euler numbers.

*Part 3*: The construction requires a suitable isomorphism between $E_T$ and $T \times \mathbb{C}$. To obtain one, fix an isomorphism between $E|_\gamma$ and $\gamma \times \mathbb{C}$ that writes the pull-back of A on $\gamma$ as $\theta + a_{A0}\,dt$ with $a_{A0}$ being constant and having absolute value less than $\frac{2\pi}{\ell_\gamma}$. Use parallel transport by A along the rays from the origin in each constant t disk to define an isomorphism between $E|_T$ and $T \times \mathbb{C}$. View A and $\alpha$ using this isomorphism as a pair of connection on the product bundle and map to $\mathbb{C}$, and use the coordinates (t, z) with this isomorphism to view $\beta$ likewise as a map to $\mathbb{C}$. Write A as $\theta + a_{A0}dt + a_A^\perp$ with $a_{A0}$ being an $i\mathbb{R}$ valued function and where $a^\perp$ has the form $\frac{1}{2}\hat{A}(z\,d\bar{z} - \bar{z}\,dz)$ with $\hat{A}$ being a real valued function. Use this isomorphism to view $\alpha$ as a map from T to $\mathbb{C}$. As explained in the next paragraph, the functions $\alpha$, $a_{A0}$ and the 1-form $a_A^\perp$ are such that

$$|\tfrac{\partial}{\partial t}\alpha| + |a_{A0}| + r^{-1/2}\,(|a_A^\perp| + |\tfrac{\partial}{\partial t}a_A^\perp|) \le c_0 c_v^{\,4}\,. \tag{B.45}$$

To justify these bounds, introduce polar coordinates on $\mathbb{C}$ by writing $z = \rho\, e^{i\sigma}$. The pull-back of $A - \theta$ to a given constant t slice of T is $a_A^\perp$. When written using $d\rho$ and $d\sigma$ this i-valued 1-form appears as $a_A^\perp = -i\hat{A}\rho^2 d\sigma$ with $\hat{A}$ being an $\mathbb{R}$-valued function. The fact that this pull-back of $A - \theta$ lacks a $d\rho$ component implies that $\frac{\partial}{\partial\rho} a_{A0} = F_A(\frac{\partial}{\partial\rho}, \frac{\partial}{\partial t})$ where $F_A$ is the curvature 2-form of A. Keeping in mind that $\frac{\partial}{\partial t}$ and $v$ differ on T by no more than $c_0|z|$, integrate this identity using the bounds in Lemmas 2.1, B.2 and B.7 to obtain the asserted bound for $|a_{A0}|$. Use this bound on $|a_{A0}|$, the aforementioned bound on $|\frac{\partial}{\partial t} - v|$, the fact that $|(\nabla_A\alpha)_v| \le c_0$ and $|\nabla_A\alpha| \le c_0 r^{1/2}$ to see that $|\frac{\partial}{\partial t}\alpha|$ is bounded by $c_0 c_v^{\,4}$. The asserted bound on $a_A^\perp$ follows by integrating the curvature identity $\rho^{-1}\frac{\partial}{\partial\rho}(\rho^2\hat{A}) = F_A(\frac{\partial}{\partial z}, \frac{\partial}{\partial\bar{z}})$ using again Lemmas 2.1, B.2 and B.7. To obtain the bound for the t-derivative of $a_A^\perp$, first differentiate the curvature identity $\frac{\partial}{\partial\rho} a_{A0} = F_A(\frac{\partial}{\partial\rho}, \frac{\partial}{\partial t})$ to obtain an equation for $\frac{\partial}{\partial\rho}(\frac{\partial}{\partial\sigma} a_{A0})$. Meanwhile, a third curvature identity has $\rho^2 \frac{\partial}{\partial t}\hat{A} + i\frac{\partial}{\partial\sigma} a_{A0} = iF_A(\frac{\partial}{\partial t}, \frac{\partial}{\partial\sigma})$. Differentiate this last equation to obtain an equation for $\frac{\partial}{\partial\rho}(\rho^2 \frac{\partial}{\partial t}\hat{A})$ that involves derivatives of the components of $F_A(\frac{\partial}{\partial t}, \cdot)$. Given what was said previously about $\frac{\partial}{\partial t}$ and $\nu$, and given the bounds in Lemmas 2.1, B.2 and B.7, integration of this last equation finds $|\frac{\partial}{\partial t} a_A{}^{\perp}| \le c_0 c_v{}^4 r^{1/2}$.

*Part 4*: This part defines $(A_{\Diamond\tau}, \psi_{\Diamond\tau})$ for $\tau \in [0, \frac{1}{3}]$. The definition requires a preliminary lemma, Lemma B.14. To set the notation, reintroduce $\varphi_r: \mathbb{C} \to \mathbb{C}$, the rescaling map given by the rule $z \to r^{-1/2} z$; and introduce $D_* \subset \mathbb{C}$ to denote the radius $c_v{}^4$ disk with center at the origin.

**Lemma B.14**: *There exists* $\kappa \ge 100$*, and given* $m \ge 1$ *and* $c_v \ge \kappa$*, there exists* $\kappa_{c,m} > \kappa$ *with the following significance: Suppose that* $r \ge \kappa_{c,m} c_v{}^{10}$ *and suppose that* $(A, \psi = (\alpha, \beta))$ *is a solution to the* $(r, \mu)$ *version of (1.13) with* $\mu$ *a given element in* $\Omega$ *with* $\mathcal{P}$*-norm smaller than 1. Use* $(A, \psi)$ *to define* T *as above. There exists a smooth map* $t \to \mathfrak{c}(t)$ *from* $\gamma$ *to the space of solutions to (2.8) on* $\mathbb{C}$ *with the properties listed below.*

- *The integral in (3.1) is finite, independent of* t *and either* $m_\alpha$ *or* $m_\alpha + 1$.
- *Any given version of* $\mathfrak{c}(t)$ *has the form* $(\theta + \mathcal{A}_{*t}, \alpha_{*t})$ *and for each* t*, the pair* $(\mathcal{A}_{*t}, \alpha_{*t})$ *on* $D_*$ *differs from* $(\varphi_r{}^* a_A{}^{\perp}, \varphi_r{}^*\alpha)|_t$ *in the* $C^0$ *topology by at most* $c_v{}^{-m}$.
- *The assignment* $t \to (\mathcal{A}_{*t}, \alpha_{*t})$ *is such that* $|\frac{\partial}{\partial t}\alpha_{*t}| + |\frac{\partial}{\partial t}\mathcal{A}_{*t}| \le \kappa c_v{}^4$ *on* $D_*$.

This lemma is proved in Part 5; assume it for now. The $\tau \in [0, \frac{1}{3}]$ version of $A_{\Diamond\tau}$ on the $|z| \le (c_v{}^4 - 2c_v{}^3) r^{-1/2}$ part of T is $A_{\Diamond\tau} = \theta + (1 - 3\tau) a_{A0} dt + a_A{}^{\perp} + 3\tau(r_r{}^*\mathcal{A}_* - a_A{}^{\perp})$. Meanwhile, the respective E and $EK^{-1}$ components of $\psi_{\Diamond\tau}$ on this same portion of T are defined by the rule $\alpha_{\Diamond\tau} = \alpha + 3\tau(r_r{}^*\alpha_* - \alpha)$ and $\beta_{\Diamond\tau} = (1 - 3\tau)\beta$. The definition on the rest of T is given by using the connection $\theta + (1 - 3\tau) a_{A0} dt + a_A{}^{\perp} + 3\tau(r_r{}^*\mathcal{A}_* - a_A{}^{\perp})$ in lieu of A and the sections $\alpha + 3\tau(r_r{}^*\alpha_* - \alpha)$ and $(1 - 3\tau)\beta$ in lieu of $(\alpha, \beta)$ to define the various functions and 1-forms that appear in (B.8)-(B.10). Keep in mind when doing so that the various isomorphisms between E and the product bundle that are invoked when writing (B.8)-(B.10) are not the isomorphisms that are used here.

*Part 5*: This part contains the proof of Part 4's lemma.

***Proof of Lemma B.14***: The proof has five steps.

$\underline{\text{Step 1}}$: Let $D_T$ denote the centered, radius $c_v^4 + c_v^3$ disk in $\mathbb{C}$. Fix $t \in \gamma$ and use Lemma 2.9 with the pair $(A,\alpha)|_t$ to obtain a solution to (2.8)'s vortex equations on $\mathbb{C}$ that can be written as $(\theta + \mathcal{A}_{1t}, \alpha_{1t})$ and is such that $\alpha_{1t} - \varphi_r{}^*\alpha|_t$ and $\mathcal{A}_{1t} - \varphi_r{}^* a_A{}^{\perp}|_t$ on $D_T$ have $C^{3m+1}$ norm bounded by $c_v^{-3m-1}$. The $\gamma$-parameterized family $t \to (\mathcal{A}_{1t}, \alpha_{1t})$ need not be continuous, let alone differentiable; nor must it obey the first bullet's requirement at any given $t \in \gamma$.

$\underline{\text{Step 2}}$: To obtain a $\gamma$-parametrized family that obeys the first bullet's requirements, consider first the case where $\alpha$ lacks zeros on the boundary of T. In this case, $\alpha_{1t}$ has $m_\alpha$ zeros counting multiplicities that lie where $|z| \le c_0$ and no zeros where $|z|$ is greater than $c_0$ and less than $c_v^4 + c_v^3 + 1$. Let $z \to \wp_t(z)$ denote the monic, degree $m_\alpha$ polynomial on $\mathbb{C}$ whose zeros with their corresponding multiplicity are those of $\alpha_{1t}$. Write this polynomial as $z^{m_\alpha} + \sigma_{1t} z^{m_\alpha - 1} + \cdots + \sigma_{m_\alpha t}$ and use the values of $\{\sigma_{qt}\}_{1 \le q \le m_\alpha}$ as the coordinates for an element in the vortex moduli space $\mathfrak{C}_{m_\alpha}$. Let $\mathfrak{c}_t$ denote this element. It follows from what is said in Section 2a of [T9] that there exists a purely $m$-dependent constant $c_m > 1$, and given $c_v > c_m$, there exists a purely $m$ and $c_v$ dependent constant $c_{m,c}$ with the following significance: If $c_v > c_m$ and $r > c_{m,c}$, then there exists a vortex solution on $\mathbb{C}$ that maps to $\mathfrak{c}_t$ which when written as $(\theta + \mathcal{A}_{2t}, \alpha_{2t})$ is such that the pair $(\mathcal{A}_{2t} - \varphi_r{}^* a_A{}^{\perp}, \alpha_{2t} - \varphi_r{}^*\alpha)$ on $D_*$ has $C^{3m}$ norm bounded by $2c_v^{-3m}$.

$\underline{\text{Step 3}}$: Consider next the case when $\alpha$ has two zeros on the boundary of T. Let $t_+$ and $t_-$ denote the t-values of the points where these zeros occur. One of these zeros will lie where $u > 0$ and the other where $u < 0$. Use $(t_{\alpha+}, z_+)$ to denote the coordinates of the former and $(t_{\alpha-}, z_-)$ to denote those of the latter. Let $I_\alpha \subset \gamma$ denote the oriented segment that starts at $t_{\alpha+}$ and ends at $t_{\alpha-}$ with $I_\alpha$ being the single point $t_+$ when $t_{\alpha+} = t_{\alpha-}$. The significance of $I_\alpha$ is as follows: Fix a transverse disk of radius $(c_v^4 + c_v^3) r^{-1/2}$ with center in the interior of $I_\alpha$. Then the sum of the local Euler number of the zeros of $\alpha$ on such a disk is equal to $m_\alpha - 1$. Meanwhile, this sum for a transverse disk with center on $\gamma - I_\alpha$ is equal to $m_\alpha$. Keeping in mind that the coordinate t is $\mathbb{R}/\ell_\gamma$-valued, let $t_+ \in [0, \ell_\gamma)$ denote the lift to $\mathbb{R}$ of $t_{\alpha+}$ and introduce $t_-$ to denote the lift to $\mathbb{R}$ of $t_{\alpha-}$ with $\frac{1}{2}\ell_\gamma \le t_- - t_+ < \frac{3}{2}\ell_\gamma$. Introduce $\mathfrak{p}: [t_+, t_-] \to \gamma$ to denote the projection map. The inverse image of any given point in $\gamma - I_\alpha$ is empty if $t_- - t_+ < \ell_\gamma$ and it contains a single point if $t_- - t_+ \ge \ell_\gamma$. The inverse image of any given point in $I_\alpha$ has a single point if $t_- - t_+ \le \ell_\gamma$ and two points otherwise. Fix a smooth map $z_1: [t_+, t_-] \to \mathbb{C}$ with $|\frac{\partial}{\partial t} z_1|$ constant, with $|z_1| \le (c_v^4 + \frac{9}{8} c_v^3)$ for all t, and such that $z_1(t_+) = z_+$ and $z_1(t_-) = z_-$. Require in addition that the image of $(t_+, t_-)$ lie where $|z| > c_v^4 + c_v^3$.

For each $t \in \gamma$, define monic polynomials $z \to \wp_{1t}(z)$ and $z \to \wp_{2t}(z)$ as follows:

- *Suppose that* $t \notin I_\alpha$: *The zeros of* $\wp_{1t}$ *with their corresponding multiplicity are those of* $\alpha_{1t}$ *with distance* 1 *or less from some* $|z| \le c_v^4 + c_v^3$ *zero of* $\varphi_r^*\alpha$. *If* $\mathfrak{p}^{-1}(t) = \emptyset$*, then* $\wp_{2t} = 1$; *and* $\wp_{2t} = z - z_I(\mathfrak{p}^{-1}(t))$ *if* $\mathfrak{p}^{-1}(t) \ne \emptyset$ .
- *Suppose that* $t \in I_\alpha$: *The zeros of* $\wp_{1t}$ *with their corresponding multiplicity are those of* $\alpha_{1t}$ *with distance* 1 *or less from some* $|z| \le c_v^4 + c_v^3 - 2$ *zero of* $\varphi_r^*\alpha$. *Meanwhile* $\wp_{2t}$ *is* $\prod_{t' \in \mathfrak{p}^{-1}(t)} (z - z_I(t'))$

(B.46)

For each $t \in \gamma$, use $\wp_t$ to denote the product $\wp_{1t}\wp_{2t}$. This is a monic polynomial with t-independent degree, either $m_\alpha$ or $m_\alpha + 1$. Let $m_*$ denote this degree. Use the coefficients of $\wp_t$ to specify a point in the vortex moduli space $\mathfrak{C}_{m_*}$. The observations in Section 2a of [T9] can be used to derive a purely $m$-dependent constant $c_m > 1$, and given $c_v > c_m$, a purely $m$ and $c_v$ dependent constant $c_{m,c}$ with the following significance: If $c_v > c_m$ and if $r > c_{m,c}$, then there is a solution on $\mathbb{C}$ to (2.8) that maps to $\wp_t$'s point in $\mathfrak{C}_{m_*}$ and can be written as $(\theta + \mathcal{A}_{2t}, \alpha_{2t})$ with $(\mathcal{A}_{2t}, \alpha_{2t})$ such that $(\mathcal{A}_{2t} - \varphi_r^* a_A{}^\perp, \alpha_{2t} - \varphi_r^*\alpha)$ on $D_*$ has $C^{3m}$ norm bounded by $2c_v^{-3m}$ .

<u>Step</u> <u>4</u>: The map $t \to (\mathcal{A}_{2t}, \alpha_{2t})$ satisfies the requirements of the first and second bullets of Lemma B.14, but it need not be smooth and if smooth, it need not satisfy the requirements of the third bullet. To remedy this defect, first introduce $c_\chi$ to denote the integral of the function $t \to \chi(|t| - 1)$ over $\mathbb{R}$. Fix for the moment $L \ge 1$ and define the map from $\gamma$ into $C^\infty(D_T; iT^*\mathbb{C} \oplus \mathbb{C})$ given by the rule $t \to (\mathcal{A}^L{}_t, \alpha^L{}_t)$ where

$$(\mathcal{A}^L{}_t, \alpha^L{}_t) = c_\chi^{-1} \int L\chi(L\,|\,t-s\,|-1)(\mathcal{A}_{1s}, \alpha_{1s})ds\,.$$

(B.47)

This is a smooth map. What follows are two consequences of (B.47). If $L \le c_v^{3m/2+4}$ then the $C^0$ norms of $\alpha^L{}_t - \varphi_r^*\alpha|_t$ and $\mathcal{A}^L{}_t - \varphi_r^* a_A{}^\perp|_t$ on $D_*$ are bounded by $c_v^{-3m/2}$. Moreover, the map $t \to (\mathcal{A}^L{}_t, \alpha^L{}_t)$ is such that $|\frac{\partial}{\partial t}\alpha^L{}_t| + |\frac{\partial}{\partial t}\mathcal{A}^L{}_t| \le c_0 c_v^4$ on $D_*$. However, any given t version $(\theta + \mathcal{A}^L{}_T, \alpha^L{}_t)$ need not obey the vortex equations. Even so, the pair comes very close to doing so.

<u>Step</u> <u>5</u>: To obtain a pair that obeys the vortex equation, introduce the $(\theta + \mathcal{A}^L{}_t, \alpha^L{}_t)$ version of (3.4)'s operator $\vartheta$. Denote this operator by $\vartheta_{tL}$. As explained momentarily, there is a smooth map $t \to \mathfrak{h}_t =$ from $\gamma$ to $C^\infty(\mathbb{C}; \mathbb{C} \oplus \mathbb{C}) \cap L^2(\mathbb{C}; \mathbb{C} \oplus \mathbb{C})$ with the following properties: Write $\vartheta_{tL}{}^\dagger \mathfrak{h}_t$ as $(2^{-1/2} e_{0t}, e_{1t})$. Then $(e_{0t}, e_{1t})$ has $C^0$ norm bounded by $c_0 c_v^{-3m/2}$ on $D_*$, its t-derivative on $D_*$ has pointwise norm bounded by $c_0 c_v^{4-3m/2}$ and the pair of connection and map to $\mathbb{C}$ given by $(\theta + \mathcal{A}^L{}_T + e_{0t} d\bar{z} - \bar{e}_{0t} dz,\ \alpha^L{}_t + e_{1t})$ obeys the vortex equations on $\mathbb{C}$ and defines a point in $\mathfrak{C}_{m_*}$. To explain, note that the vortex equations are obeyed if $\mathfrak{h}_t$ obeys an equation having the schematic form

$$\vartheta_{tL}\vartheta_{tL}{}^{\dagger}\mathfrak{h}_t + (\vartheta_{tL}{}^{\dagger}\mathfrak{h}_t)\#(\vartheta_{tL}{}^{\dagger}\mathfrak{h}_t) = \mathfrak{q}_t \tag{B.48}$$

where $\mathfrak{q}_t$ has $C^0$ and $L^2$ norm bounded by $c_0 c_v^{-3m/2}$. Given (3.6) and this small norm for $\mathfrak{q}_t$, the contraction mapping theorem on a suitable Hilbert space can be used to find a smooth solution $t \to \mathfrak{h}_t$ from $\gamma$ to $C^\infty(\mathbb{C};\mathbb{C}\oplus\mathbb{C}) \cap L^2(\mathbb{C};\mathbb{C}\oplus\mathbb{C})$ with $C^1$ norm bounded by $c_0 c_v^{-3m/2}$. The contraction mapping theorem construction will guarantee a pointwise norm bound by $c_0 c_v^4$ for the t-derivative of $\vartheta_{Lt}{}^{\dagger}\mathfrak{h}_t$ on D.

Granted the preceding, set $\mathcal{A}_{*t} = \mathcal{A}^L{}_T + e_{0t} d\bar{z} - \bar{e}_{0t} dz$ and set $\alpha_{*t} = \alpha^L{}_t + e_{1t}$. The resulting map $t \to (\mathcal{A}_{*t}, \alpha_{*t})$ obeys all of the lemma's requirements.

*Part 6*: This part defines the pair $(A_{\Diamond\tau}, \psi_{\Diamond\tau})$ for $\tau \in [\frac{1}{3}, \frac{2}{3}]$. To this end, let D denote a given constant t slice of T, this being a transverse disk with center on $\gamma$ and radius $(c_v^4 + c_v^3) r^{-1/2}$. If the sum of the local Euler numbers of $\alpha$ on D is defined, then it is also defined for $\alpha_\Diamond$ and these sums are the same. Note that all of the local Euler numbers of $\alpha_\Diamond$ are positive. Let $t \in \gamma$ denote the center point of such a disk. Define a monic polynomial $z \to \wp_{\Diamond t}(z)$ on $\mathbb{C}$ using the rules that follow. If $\alpha$ has no zeros on the boundary of T, then $\wp_{\Diamond t}(z) = z^{m_\alpha}$. If $\alpha$ has zeros on the boundary of T, define $\wp_{\Diamond t}$ by using $\alpha_{\Diamond t}$ in lieu of $\alpha_{1t}$ in (B.46). Meanwhile, let $\wp_{*t}$ denote the monic polynomial on $\mathbb{C}$ whose roots with their corresponding multiplicity are the zeros of the function $\alpha_{*t}$ from Lemma B.3. Note that all such zeros have positive local Euler number. The polynomials $\wp_{\Diamond t}$ and $\wp_{*t}$ have the same degree. Denote this degree by $m_*$.

Given $\tau \in [\frac{1}{3}, \frac{2}{3}]$, set $\wp_{\tau t}$ to be the monic polynomial $(2-3\tau)\, \wp_{*t} + (3\tau - 1)\, \wp_{\Diamond\tau}$. The resulting 1-parameter family of polynomials interpolates between $\wp_{*t}$ and $\wp_{\Diamond t}$. For any given pair $(\tau, t)$, the coefficients of $\wp_{\tau t}$ defines a point in $\mathfrak{C}_{m_*}$ that varies smoothly with variations in $\tau$ and t with the variation in $\tau$ being real analytic. With $\tau$ fixed for the moment, let $t \to \mathfrak{c}_\tau(t)$ denote the corresponding map from $\gamma$ to $\mathfrak{C}_m$. Lemma B.14 describes a lift of the map $t \to \mathfrak{c}_{\tau=1/3}(t)$ to a smooth map $t \to (\theta + \mathcal{A}_{*t}, \alpha_{*t})$ from $\gamma$ into the space of solutions to (2.8) on $\mathbb{C}$. The next lemma describes a corresponding smooth lift of the two variable map $(\tau, t) \to \mathfrak{c}_\tau(t)$ from $[\frac{1}{3}, \frac{2}{3}] \times \gamma$ to $\mathfrak{C}_{m_*}$.

**Lemma B.15**: *There exists* $\kappa \geq 100$*, and given* $c_v \geq \kappa$*, there exists* $\kappa_{c_v} > \kappa$ *with the following significance: Suppose that* $r \geq \kappa_{c_v} c_v^{10}$ *and suppose that* $(A, \psi = (\alpha, \beta))$ *is a solution to the* $(r, \mu)$ *version of (1.13) with* $\mu$ *a given element in* $\Omega$ *with* $\mathcal{P}$*-norm smaller than 1. Fix a component of the corresponding version of* $Y-(Y_{*\Lambda} \cup T_{*\Lambda})$ *and introduce the latter's version of the integer* m *and the map* $(\tau, t) \to \mathfrak{c}_\tau(t)$ *from* $[\frac{1}{3}, \frac{2}{3}] \times \gamma$ *to* $\mathfrak{C}_m$. *There is a smooth map* $(\tau,t) \to (\mathcal{A}_{\tau t}, \alpha_{\tau t})$ *from* $[\frac{1}{3}, \frac{2}{3}] \times \gamma$ *to* $C^\infty(\mathbb{C}; iT^*\mathbb{C} \oplus \mathbb{C})$ *which is real analytic with respect to variations in* $\tau$ *and such that at each* $(\tau, t) \in [\frac{1}{3}, \frac{2}{3}] \times \gamma$*, the pair of connection on the product bundle* $\mathbb{C} \times \mathbb{C}$ *and map to* $\mathbb{C}$ *given by* $(\theta + \mathcal{A}_{\tau t}, \alpha_{\tau t})$ *satisfies (2.8) and projects to* $\mathfrak{c}_\tau(t)$. *In addition,*

- $|\frac{\partial}{\partial t} \alpha_{\tau t}| + |\frac{\partial}{\partial t} \mathcal{A}_{\tau t}| \leq \kappa c_v^4$ *on* $D_*$ .
- $|\frac{\partial}{\partial \tau} \alpha_{\tau t}| + |\frac{\partial}{\partial \tau} \mathcal{A}_{\tau t}| \leq \kappa c_v^4$ *on* $D_*$ .

This lemma is proved momentarily. By way of notation, any given $\tau \in [\frac{1}{3}, \frac{2}{3}]$ version of the map $t \to (\mathcal{A}_{\tau t}, \alpha_{\tau t})$ from $\gamma$ into $C^\infty(\mathbb{C}; iT^*\mathbb{C} \oplus \mathbb{C})$ is denoted in what follows by $(\mathcal{A}_\tau, \alpha_\tau)$.

The $\tau \in [\frac{1}{3}, \frac{2}{3}]$ version of the connection $A_{\Diamond\tau}$ on the $|z| \leq (c_v^4 - 2c_v^3) r^{-1/2}$ part of T is $A_{\Diamond\tau} = \theta + r_r^* \mathcal{A}_\tau$ and the respective E and $EK^{-1}$ components of $\psi_{\Diamond\tau}$ on this same portion of T are defined by the rule $\alpha_{\Diamond\tau} = r_r^* \alpha_\tau$ and $\beta_{\Diamond\tau} = 0$. The definition on the rest of T is given by using the connection $\theta + r_r^* \mathcal{A}_\tau$ in lieu of A and the sections $r_r^* \alpha_\tau$ and 0 in lieu of $(\alpha, \beta)$ to define the various functions and 1-forms that appear in (B.8)-(B.10). Keep in mind when doing so that the various isomorphisms between E and the product bundle that are invoked when writing (B.8)-(B.10) are not the isomorphisms that are used here.

***Proof of Lemma B.15***: The existence of a lift of the map $(\tau,t) \to \mathfrak{c}_\tau(t)$ follows from what is said in Section 2c of [T9]. The existence of a lift with t and $\tau$ derivatives bounded by $c_0 c_v^4$ follows from what is said in this same Section 2c of [T9] using (2.5), (2.11), (2.12) and (2.19) in [T9].

*Part 7*: This part defines $(A_{\Diamond\tau}, \psi_{\Diamond\tau})$ for $\tau \in [\frac{2}{3}, 1]$. This definition is given below by (B.49). To set the notation, view $A_\Diamond$ and the pair $(\alpha_\Diamond, \beta_\Diamond)$ as a respective connection on $T \times \mathbb{C}$ and pair of maps from T to $\mathbb{C}$ using the same isomorphism of $E|_T$ with $T \times \mathbb{C}$ that is used to define the $\tau = \frac{2}{3}$ version of $(A_{\Diamond\tau}, \psi_{\Diamond\tau})$. The definition writes this depiction of $A_\Diamond$ as $\theta + A_\Diamond'$ and it writes this depiction of $(\alpha_\Diamond, \beta_\Diamond)$ as $(\alpha_\Diamond', \beta_\Diamond')$. The connection $A_{\Diamond\tau=\frac{2}{3}}$ is written below as $\theta + A_{\Diamond\frac{2}{3}}$ . Equation (B.49) refers to a map $\hat{u}: \gamma \times \mathbb{C} \to S^1$ that is described below by Lemma B.16. Fix $\tau \in [\frac{2}{3}, 1]$ and what follows defines $(A_{\Diamond\tau}, \psi_{\Diamond\tau})$ on T.

- $A_{\Diamond\tau} = \theta + (3\tau - 1)(A_\Diamond' - \hat u^{-1}d\hat u) + (3 - 3\tau)A_{\Diamond\frac{2}{3}}$ ,
- $\alpha_{\Diamond\tau} = (3\tau - 1)\hat u\alpha_\Diamond' + (3 - 3\tau)\alpha_{\Diamond\frac{2}{3}}$ *and* $\beta_{\Diamond\tau} = (3\tau - 1)\hat u\beta_\Diamond'$.

(B.49)

The map û is constructed in the proof of Lemma B.16.

**Lemma B.16**: *There exists* $\kappa \geq 100$*, and given* $m \geq 1$ *and* $c_v \geq \kappa$*, there exists* $\kappa_{c,m} > \kappa$ *with the following significance: Suppose that* $r \geq \kappa_{c,m} c_v^{10}$ *and let* $(A, \psi = (\alpha, \beta))$ *denote a solution to the* $(r, \mu)$ *version of (1.13) with* $\mu$ *an element in* $\Omega$ *with* $\mathcal{P}$*-norm less than 1. Use* $(A, \psi)$ *to define* T *and the corresponding versions of* $(A_\Diamond', \alpha_\Diamond')$ and $(A_{\Diamond\frac{2}{3}}, \alpha_{\Diamond\frac{2}{3}})$.
*There exists a smooth map* $\hat u$: $T \to S^1$ *such that*

$$r^{-1/2}|A_\Diamond' - \hat u^{-1}d\hat u - A_{\Diamond\frac{2}{3}}| + |\hat u\alpha_\Diamond' - \alpha_{\Diamond\frac{2}{3}}| \leq c_v^{-m} \text{ and } |A_\Diamond'(\tfrac{\partial}{\partial t}) - \hat u^{-1}\tfrac{\partial}{\partial t}\hat u| \leq \kappa c_v^{4}.$$

***Proof of Lemma B.16***: The two steps that follow constructs û on the $|z| \leq \frac{1}{2} c_v^{1/2} r^{-1/2}$ portion of T. But for cosmetic changes, the same construction supplies û on the rest of T.

Step 1: To define û where $|z| \leq \frac{1}{2} c_v^{1/2} r^{-1/2}$, recall from Part 3 that the pull-back of $(A_\Diamond, \psi_\Diamond)$ to the $|z| \leq \frac{1}{4} c_v^{1/2} r^{-1/2}$ portion of a transverse disk centered on any given $t \in \gamma$ is the solution to the vortex equations in (2.8) given by $(\theta - \frac{m}{2} r_r^* a_{m0}(z^{-1}dz - \bar z^{-1}d\bar z), r_r^*\alpha_{m0})$. Meanwhile, the pull-back of $(A_{\Diamond\tau=\frac{2}{3}}, \alpha_{\Diamond\tau=\frac{2}{3}})$ to the same part of the transverse disk centered at $t \in \gamma$ is a solution to (2.8) that was written as $(\theta + r_r^*\mathcal{A}_{\frac{2}{3}t}, r_r^*\alpha_{\frac{2}{3}t})$. The two $\mathbb{C}$ valued functions $\alpha_{m0}$ and $\alpha_{\frac{2}{3}t}$ have the same zero locus on the $|z| \leq \frac{3}{4} c_v^{1/2}$, this being the origin. Moreover, they have the same local degree at 0. What follows is a consequence: There exists a smooth map, denoted here by u, from the $|z| \leq \frac{9}{16} c_v^{1/2}$ part of $\gamma \times \mathbb{C}$ to $S^1$ such that $u\alpha_{\frac{2}{3}} = |\alpha_{\frac{2}{3}}|z^m$.

Fix a positive integer $m$. Granted the preceding, use what is said in Part 4 of Section 2a in [T9] about solutions to (2.8) to find a purely $m$ dependent lower bound for $c_v$ such that the subsequent assertion is true when $c_v$ exceeds this bound. Introduce $d^\perp u$ to denote the exterior derivative of u along the constant t slices of $\gamma \times \mathbb{C}$. For $t \in \gamma$, the pairs $(\mathcal{A}_{\frac{2}{3}t} - u^{-1}d^\perp u, u\alpha_{\frac{2}{3}t})$ and $(-\frac{m}{2} a_{m0}(z^{-1}dz - \bar z^{-1}d\bar z), \alpha_{m0})$ differ by at most $c_v^{-3m}$ in the $C^{2m}$ topology on the $|z| \leq \frac{5}{8} c_v^{1/2}$ disk in $\mathbb{C}$.

This last conclusion has the following consequence: If $c_v$ is greater than a purely $m$ dependent lower bound, then the map $u_1 = u^{-1}(\varphi_r^*\alpha_\Diamond')^{-1}\alpha_{m0}$, from the $|z| \leq \frac{9}{16} c_v^{1/2}$ part of $\gamma \times \mathbb{C}$ to $S^1$ is such that for any $t \in \gamma$, the pair $(\mathcal{A}_{\frac{2}{3}t}, \alpha_{\frac{2}{3}t})$ and the pull-back to $\{t\} \times \mathbb{C}$ of the pair $(\varphi_r^*A_\Diamond' - u_1^{-1}du_1, u_1\varphi_r^*\alpha_\Diamond')$ differ by less than $c_v^{-5m/2}$ in the $C^{2m}$ topology on the disk $|z| \le \frac{5}{8} c_v^{1/2}$.

Step 2: The map $u_1$ can be replaced by a map $u_2$: $\gamma\times\mathbb{C} \to S^1$ such that if $c_v$ is greater than a purely $m$-dependent constant, then

- *For* $t \in \gamma$, *the the pull-back to* $\{t\} \times \mathbb{C}$ *of* $(\varphi_r^*A_\Diamond' - u_2^{-1}du_2, u_2\varphi_r^*\alpha_\Diamond')$ *and* $(\mathcal{A}_{\frac{2}{3}t}, \alpha_{\frac{2}{3}t})$ *differ pointwise on the* $|z| \le \frac{5}{8} c_v^{1/2}$ *part of* $\mathbb{C}$ *by at most* $c_v^{-m}$.
- $|A_\Diamond'(\frac{\partial}{\partial t}) - u_2^{-1}\frac{\partial}{\partial t}u_2| \le c_0 c_v^4$ *on the* $|z| \le \frac{5}{8} c_v^{1/2}$ *part of* $\mathbb{C}$.

(B.50)

The map $\hat{u}$ on the $|z| \le \frac{5}{8} c_v^{1/2} r^{-1/2}$ part of T is defined to be $r_r^*u_2$.

To construct $u_2$ write $u_1$ as $u^{-1}(\varphi_r^*\alpha_\Diamond')^{-1}\alpha_{m0}$ and write u as $e^{i(2\pi nt/\ell_\gamma + x)}$ with n being an integer and x being a real valued function on $\gamma\times\mathbb{C}$. The map $u_2$ has the form $e^{i(2\pi n/\ell_\gamma + x_2)}$ with $x_2$ being a real valued function on $\gamma \times \mathbb{C}$. The function $x_2$ is the smoothing of the function x that is given by the rule $x_2|_t = c_\chi^{-1}\int L\chi(L\,|\,t\text{-}s\,|\text{-}1)\, x\,|_s ds$ with $c_\chi$ being the constant that appears in (B.47) and with $L = c_v^{5m/4}$. The resulting map $u_2$ obeys the inequality in the first bullet of (B.50) if $c_v$ is greater than a purely $m$-dependent constant. This is direct consequence of the $c_v^{-5m/2}$ bound obtained in Step 1. Meanwhile, $|\frac{\partial}{\partial t} x_2| \le c_0 c_v^4$, this being a consequence of this same $c_v^{-5m/2}$ bound and the bound in the top bullet of Lemma B.15. Granted all of this, then the Lemma B.16's right most inequality is obeyed if the integer n is such that $|n| \le c_0 c_v^4$.

To obtain such a bound for n, fix a constant z circle in T with $|z| = \frac{1}{4} \kappa_\Diamond^{-1} c_v$ and with distance at least $\frac{1}{100} \kappa_\Diamond^{-1} c_v^{1/4}$ or more from $\alpha$'s zero locus. Proposition 2.4 guarantees the existence of such circles $c_v \ge c_0$ and if r is greater than a purely $c_v$ dependent constant. The integral over the chosen circle of $-i\,u^{-1}\frac{\partial}{\partial t}u$ is equal to $2\pi\, n\, \ell_\gamma^{-1}$, and so upper and lower bounds on this integral gives a bound for |n|. A suitable bound is obtained by writing $\alpha_{\frac{2}{3}}$ as $|\alpha_{\frac{2}{3}}|\, u^{-1} z^m$ to derive the identity

$$i(\alpha_{\frac{2}{3}})^{-1}\frac{\partial}{\partial t}\alpha_{\frac{2}{3}} = 2\pi n\, \ell_\gamma^{-1} + \frac{\partial}{\partial t}x + i\,\frac{\partial}{\partial t}\ln(|\alpha_{\frac{2}{3}}|)\,.$$

(B.51)

Integrate both sides of this identity on the given circle. The integral of the right hand side is $2\pi n$; and the top bullet in Lemma B.15 bounds the absolute value of the integral of the left hand side by $c_0 c_v^4$.

*Part 8*: The promised proof of Proposition B.13 is given below. By way of a look ahead, the proof uses the results from Appendex A in much the same way as does the proof of Proposition B.3. Most of what is said in Appendex A requires PROPERTIES 1-5 in Section Ab; the fact that each $\tau \in [0,1]$ version of $(A_{\Diamond\tau}, \psi_{\Diamond\tau})$ has these properties is asserted by the next lemma.

**Lemma B.17**: *There exists* $\kappa \ge 100$*, and given* $c_v \ge \kappa$*, there exists* $\kappa_{c_v} > \kappa$ *with the following significance: Suppose that* $r \ge \kappa_{c_v} c_v^{10}$ *and suppose that* $(A, \psi = (\alpha, \beta))$ *is a solution to the* $(r, \mu)$ *version of (1.13) with* $\mu$ *a given element in* $\Omega$ *with* $\mathcal{P}$*-norm smaller than 1. Each element in the corresponding path* $\{(A_\tau, \psi_\tau)\}_{\tau\in[0,1]}$ *obeys* PROPERTIES 1-5 *in Section Ab.*

***Proof of Lemma B.17***: The assertion follows from Lemma A.1 if it is the case that PROPERTIES 1-5 hold on $Y-(Y_{*\Lambda} \cup T_{*\Lambda})$. To verify that this is indeed the case, focus attention now on a given component of this set. Let $\gamma$ denote the corresponding curve from $\cup_{p\in\Lambda}\{\hat{\gamma}_p^+ \cup \hat{\gamma}_p^-\}$ and let T denote the radius $(c_v^4 + c_v^3) r^{-1/2}$ tubular neighborhood of $\gamma$. The fact that PROPERTIES 4 and 5 hold on T when $c_v \ge c_0$ and r is larger than a purely $c_v$-dependent constant follows from (A.4) and Lemmas B.14-B.16. The fact that PROPERTIES 1 and 2 hold on T follows from (A.4), and Lemmas B.14 –B.15 given that the vectors fields $\frac{\partial}{\partial t}$ and $v$ differ on T by at most $c_0 c_v^4$. The details of the argument are straightforward and left to the reader but for the remark that the verification of the second and third bullets of PROPERTY 2 require the third bullet of Lemma B.14, the first bullet of Lemma B.15 and the bound for $|A_\Diamond{}'(\frac{\partial}{\partial t}) - \hat{u}^{-1}\frac{\partial}{\partial t}\hat{u}|$ in Lemma B.16.

***Proof of Proposition B.13***: The assertion of the proposition follows if there is a purely $c_v$ dependent $\kappa_c \ge 1$ with the following property: Assume that $c_v \ge c_0$ and $r \ge \kappa_c$. Fix any interval $[\tau, \tau'] \subset [0,1]$ of length at most $\kappa_c^{-1}$. Then the norm of the difference between the values of $\mathfrak{f}_\mathfrak{s}$ at $(A_{\Diamond\tau}, \psi_{\Diamond\tau})$ and at $(A_{\Diamond\tau'}, \psi_{\Diamond\tau'})$ is bounded by $c_0$. The six steps of the proof exhibit purely $c$ dependent $\kappa_c$ with this property.

Step 1: What is said in Part 1 of the proof of Proposition B.3 applies to the family $\{\mathfrak{L}_{\mathbb{V}\tau}\}_{\tau\in[0,1]}$ where any given $\tau \in [0,1]$ member is the $(A_{\Diamond\tau}, \psi_{\Diamond\tau})$ version of the operator $\mathfrak{L}_{\mathbb{V}}$ that is depicted in (A.26) and (A.27). This being the case, there is the corresponding set of eigenvalue families $\{\lambda_{n\tau}\}_{n\in\mathbb{Z},\tau\in[0,1]}$. Keep in mind that all $\tau \in [0,1]$ versions of $\mathfrak{L}_{\mathbb{V}\tau}$ are identical on $Y_{*\Lambda} \cup T_{*\Lambda}$. This has the following consequence: Fix $n \in \mathbb{Z}$ and an interval in $[0,1]$ where the map $\tau \to \lambda_{n\tau}$ is differentiable. Let $\tau \to \mathfrak{f}_{(\tau)}$ denote the corresponding family of unit $L^2$-norm eigenvectors. Then the relevant version of (B.20) has the form

$$\tfrac{d}{d\tau}\lambda_{n\tau} = \int_{Y-(Y_{*\Lambda}\cup T_{*\Lambda})} \mathfrak{f}_{(\tau)}^{\dagger}(\tfrac{d}{d\tau}\mathfrak{L}_{\mathbb{V}\tau})\mathfrak{f}_{(\tau)} \,. \tag{B.52}$$

Fix $m > c_0$ and take $c_v$ and r so as to invoke Lemmas B.14 and B.16. If I is either of the intervals $[0, \frac{1}{3}]$ or $[\frac{2}{3}, 1]$, then these lemmas imply that the $\tau \in$ I versions of $\frac{d}{d\tau}\mathfrak{L}_{\mathbb{V}(\cdot)}$ is an endomorphism of $\mathbb{V}$ with pointwise norm bounded by $c_0 c_v^{-m} r^{1/2}$. This being the case, integrate (B.52) to draw the following conclusion: Fix $n \in \mathbb{Z}$. If $\lambda_{n(\cdot)}$ has a zero in I, then $|\lambda_{n\tau}| \le c_0 c_v^{-m} r^{1/2}$ for all $\tau \in$ I.

Suppose in addition that $c_v$ and r are such that Lemma B.15 can also be invoked. Lemma B.15 implies that the $\tau \in [\frac{1}{3}, \frac{2}{3}]$ version of $\frac{d}{d\tau}\mathfrak{L}_{\mathbb{V}(\cdot)}$ has pointwise norm bounded by $c_0 r^{1/2}$. This understood, fix an interval $I \subset [\frac{1}{3}, \frac{2}{3}]$ with length at most $c_v^{-m}$. Integrate (B.52) on the interval I to deduce the analog of what is said at the end of the preceding paragraph: Fix $n \in \mathbb{Z}$. If $\lambda_{n(\cdot)}$ has a zero in I, then $|\lambda_{n\tau}| \le c_0 c_v^{-m} r^{1/2}$ for all $\tau \in$ I.

<u>Step 2</u>: With $m \ge 1$ fixed, take $c_v$ and r large enough to invoke the preceding lemmas in this Appendix B and the $c_0 = c_v$ versions of the lemmas in Appendix A. Let I now denote any given interval in $[0,1]$ of length at most $c_v^{-m}$. If $m > c_0$, then Lemma A.6 can be invoked to draw the following conclusion: Let $n \in \mathbb{Z}$ be such that $\lambda_{n(\cdot)}$ has a zero in I. Fix $\tau \in$ I and use $\mathfrak{f}_{(\tau)}$ to denote an eigenvector of $\mathfrak{L}_{\mathbb{V}\tau}$ with eigenvalue $\lambda_{n\tau}$. Then $\|\Pi_{\vartheta}\mathfrak{f}_{(\tau)}\| \ge (1 - c_0 c_v^{-1})\|\mathfrak{f}_{(\tau)}\|_2$.

Supposing that $m \ge c_0$, that $c_v$ is greater than a purely $m$-dependent constant, and that r is greater than a purely $m$ and $c_v$ dependent constant, then Lemmas A.7 and A.8 can be invoked to conclude the following: Let $I \subset [0, 1]$ denote an interval of length $c_0^{-1} c_v^{-m}$. If $n \in \mathbb{Z}$ and $\lambda_{n(\cdot)}$ has a zero in I, then $|\lambda_{n\tau}| \le c_v^{-m}$ for all $\tau \in$ I.

<u>Step 3</u>: Let $I \subset [0, 1]$ denote an interval of length at most $c_0^{-1} c_v^{-m}$. Write I as $[\tau, \tau']$. Granted that the conclusion of the preceding step hold for I, then the argument used in Part 4 of the proof of Proposition B.3 can be repeated with only notational changes to see that the norm of the difference between the respective values of $\mathfrak{f}_{\mathfrak{s}}$ at $(A_{\Diamond\tau}, \psi_{\Diamond\tau})$ and at $(A_{\Diamond\tau'}, \psi_{\Diamond\tau'})$ is at most $c_0$.

## C. Paths in Conn(E) $\times$ $C^\infty$(Y; $\mathbb{S}$) from vortex solutions

This last section of the appendix first constructs a deformation of ($A_\Diamond$, $\psi_\Diamond$) through a family of pairs in Conn(E) $\times$ $C^\infty$(Y; $\mathbb{S}$), all made from vortex solutions as in Section Aa and (A.44) using $z$ = r. The end result is then deformed through a family that is defined using vortex solutions as done in Section Aa and (A.44) using ever *increasing* values of $z$. The end result of this deformation is a pair whose resulting version of $\mathfrak{L}_{\mathbb{V}}$ as defined using $z$ >> r can be compared with that of a $z$ = $\mathcal{O}(1)$ version using a strategy from [T4]. These comparisons are used in Section Ce of this appendix to complete the proof of Propostion 2.6

### a) Deforming the zero locus of $\alpha_\Diamond$

The zero locus of $\alpha_\Diamond$ is a disjoint union of two sorts of embedded circles. The first are curves from the set $\cup_{\mathfrak{p}\in\Lambda}\{\hat{\gamma}^+_{\mathfrak{p}} \cup \hat{\gamma}^-_{\mathfrak{p}}\}$. The remainder consist of a finite set of at most G embedded circles that look very much like the subset of curves from a generator of the embedded contact homology chain complex that intersect the $f^{-1}(1,2)$ part of $M_\delta$. With this in mind, this subsection constructs a path in Conn(E)$\times C^\infty$(Y; $\mathbb{S}$) from ($A_\Diamond$, $\psi_\Diamond$) that ends at a pair of connection on E and section of $\mathbb{S}$ with the following property: The section of $\mathbb{S}$ when written with respect to the decomposition $\mathbb{S}$ = E$\oplus$E$K^{-1}$ has E component whose zero locus consists entirely of closed integral curves of $v$.

The construction of the desired path occupies the first three parts of the subsection. The fourth part of the subsection states and then proves a proposition that supplies an r-independent bound for the absolute value of the difference between $\mathfrak{f}_{\mathfrak{s}}$ at ($A_\Diamond$, $\psi_\Diamond$) and at the end member of the path.

*Part 1*: Given that r $\geq$ $c_0$, it follows from Proposition 2.4 and Proposition II.2.7 that there exists a set of closed integral curves of $v$ whose intersection with $M_\delta$ is every where very close to $\alpha^{-1}(0) \cap M_\delta$. This set of curves is denoted here by $\Theta^\alpha$; it is parametrized as in Proposition II.2.7 as $\Theta^\alpha$ = ($\hat{\upsilon}^\alpha$ , $(\mathfrak{k}^\alpha{}_{\mathfrak{p}})_{\mathfrak{p}\in\Lambda}$). The component $\hat{\upsilon}^\alpha$ from $\Theta^\alpha$ describes how the curves from $\Theta^\alpha$ intersect $M^\delta$, and each $\mathfrak{p} \in \Lambda$ version of $\mathfrak{k}^\alpha{}_{\mathfrak{p}}$ is an integer that describes how the curves from $\Theta^\alpha$ intersect $\mathcal{H}_{\mathfrak{p}}$. The paragraphs that follow say more about the significance of the parameterization that is used by [KLTII].

What is denoted by $\hat{\upsilon}^\alpha$ signifies a certain set of G segments of integral curves of $v$ in the $f^{-1}(1,2)$ part of $M_\delta$; these being integral curves that extend into M as integral curves of the pseudogradient vector field for $f$ that was used in Section II.1 to define the geometry of Y. The segments that comprise $\hat{\upsilon}^\alpha$ define a pairing between the index 1 critical points of the incarnation of $f$ as a function on M and the latter's index 2 critical points in the following sense: Each arc from this set starts on the boundary of the radius $\delta$ coordinate ball in $M_\delta$ corresponding to an index 1 critical point of $f$, and each ends on the boundary of the radius $\delta$ coordinate ball in $M_\delta$ of an index 2 critical point of $f$. Moreover, distinct arcs start on distinct radius $\delta$ coordinate balls and end on distinct radius $\delta$ coordinate balls. The section $\alpha$ determines $\hat{\upsilon}^\alpha$ in the following way: The pairing of index 1 critical points of $f|_M$ with index 2 critical points that is determined via $\alpha$ as described in the third bullet of Proposition 2.4 is the same pairing given by $\hat{\upsilon}^\alpha$. Moreover, the respective components of $\alpha^{-1}(0) \cap M_\delta$ and $\hat{\upsilon}^\alpha$ that pair the same index 1 and index 2 critical points of $f|_M$ are in each other's radius $c_0^{-1}\delta$ tubular neighborhoods.

As noted above, the component $(\mathfrak{k}^\alpha_{\mathfrak{p}})_{\mathfrak{p}\in\Lambda}$ of $\Theta^\alpha$ consist of a set of integers that are labeled by the pairs in $\Lambda$. The remainder of Part 1 explains how $\alpha$ determines this set. To this end, let $\upsilon$ denote a component of the zero locus of $\alpha_\Diamond$ that intersects $M_\delta$ and let $\hat{\upsilon}^{\alpha,\upsilon} \subset \hat{\upsilon}^\alpha$ denote the subset which corresponds to $\upsilon \cap M_\delta$ in the sense that corresponding arcs label the same index 1 and index 2 critical points of $f|_M$. Introduce $\Lambda_\upsilon$ to denote the subset of $\mathfrak{p} \in \Lambda$ with $\upsilon \cap \mathcal{H}_{\mathfrak{p}} \neq \emptyset$ and suppose for the moment that $\mathfrak{k} = (\mathfrak{k}_{\mathfrak{p}})_{\mathfrak{p}\in\Lambda_\upsilon}$ is a given set of integers parametrized by $\Lambda_\upsilon$. Proposition II.2.7 uses the sets $\hat{\upsilon}^{\alpha,\upsilon}$ and $\mathfrak{k}$ to define a closed integral curve of $v$. Let $\upsilon^{\mathfrak{k}}$ to denote this integral curve of $v$. The next paragraph summarizes some facts about $\upsilon^{\mathfrak{k}}$ that follow from Proposition II.2.7.

The label $\mathfrak{k}$ makes a significant difference with regards to the behavior of $\upsilon^{\mathfrak{k}}$ on the various $\mathfrak{p} \in \Lambda_\upsilon$ versions of $\mathcal{H}_{\mathfrak{p}}$. To say more, fix an element $\mathfrak{p} \in \Lambda_\upsilon$. Then $\upsilon^{\mathfrak{k}} \cap \mathcal{H}_{\mathfrak{p}}$ is an arc that crosses $\mathcal{H}_{\mathfrak{p}}$ where $1 - 3\cos^2\theta > 0$ starting from the $u < 0$ boundary of $\mathcal{H}_{\mathfrak{p}}$ and ending on the $u > 0$ boundary of $\mathcal{H}_{\mathfrak{p}}$. These endpoints have distance at most $c_0^{-1}\delta$ from the corresponding endpoints of $\upsilon \cap \mathcal{H}_{\mathfrak{p}}$. This understood, define a continuous and piecewise smooth loop in $\mathcal{H}_{\mathfrak{p}}$ as follows: Start on the $u < 0$ boundary point of $\upsilon^{\mathfrak{k}} \cap \mathcal{H}_{\mathfrak{p}}$ and travel along $\upsilon^{\mathfrak{k}} \cap \mathcal{H}_{\mathfrak{p}}$ to its $u > 0$ boundary. Take the short geodesic arc from this boundary point of $\upsilon^{\mathfrak{k}} \cap \mathcal{H}_{\mathfrak{p}}$ to the nearby boundary point of $\upsilon \cap \mathcal{H}_{\mathfrak{p}}$. Having done so, travel in the reverse direction along $\upsilon \cap \mathcal{H}_{\mathfrak{p}}$ to its boundary point on the $u < 0$ boundary of $\mathcal{H}_{\mathfrak{p}}$. Then take the short geodesic arc to the starting point on $\upsilon^{\mathfrak{k}} \cap \mathcal{H}_{\mathfrak{p}}$. The result is an oriented, piecewise smooth loop in the $1 - 3\cos^2\theta$ part of $\mathcal{H}_{\mathfrak{p}}$ and thus a class in the first homology of the $1 - 3\cos^2\theta$ part of $\mathcal{H}_{\mathfrak{p}}$. Meanwhile, the first homology of this part of $\mathcal{H}_{\mathfrak{p}}$ is isomorphic to $\mathbb{Z}$ with generator being the $u = 0$, $\cos\theta = 0$ circle. The loop just constructed from $\upsilon \cap \mathcal{H}_{\mathfrak{p}}$ and $\upsilon^{\mathfrak{k}} \cap \mathcal{H}_{\mathfrak{p}}$ defines an element in this homology class, thus an integer multiple of the generator. This integer can be written as $\mathfrak{m}_{\upsilon,\mathfrak{p}} + \mathfrak{k}_{\mathfrak{p}}$ with $\mathfrak{m}_{\upsilon,\mathfrak{p}}$ depending on $\upsilon \cap \mathcal{H}_{\mathfrak{p}}$ but not on $\mathfrak{k}$.

Granted the preceding, any given $\mathfrak{p} \in \Lambda_\upsilon$ version of the integer $\mathfrak{k}^\alpha_{\mathfrak{p}}$ coming from $\Theta^\alpha$ is $-\mathfrak{m}_{\upsilon,\mathfrak{p}}$. This is to say that $\mathfrak{k}_{\mathfrak{p}} = \mathfrak{k}^\alpha_{\mathfrak{p}}$ version of the loop in $\mathcal{H}_{\mathfrak{p}}$ described in the preceding paragraph is null-homotopic.

The subsequent parts of this subsection use $\upsilon^\alpha \subset \Theta^\alpha$ to denote the loop that is defined by the subsets $\hat{\upsilon}^{\alpha,\upsilon} \subset \hat{\upsilon}^\alpha$ and components $(\mathfrak{k}^\alpha_{\mathfrak{p}})_{\mathfrak{p}\in\Lambda_\upsilon} \subset (\mathfrak{k}_{\mathfrak{p}})_{\mathfrak{p}\in\Lambda}$.

*Part 2*: The introduction promises a path in $\mathrm{Conn}(E)\times C^\infty(Y;\mathbb{S})$ from $(A_\Diamond, \psi_\Diamond)$ that ends at a pair whose section of $\mathbb{S} = E\oplus EK^{-1}$ has E component with zero locus consisting entirely of closed integral curves of $v$; these being the curves from $\Theta^\alpha$ and the curves from $\cup_{p\in\Lambda}\{\hat{\gamma}_p^+ \cup \hat{\gamma}_p^-\}$ that lie in $\alpha_\Diamond^{-1}(0)$. The path in $\mathrm{Conn}(E)\times C^\infty(Y;\mathbb{S})$ is parametrized by $[0,1]$ and a given $\tau \in [0,1]$ member of this path denoted in what follows by $(A_{\Diamond 1\tau}, \psi_{\Diamond 1\tau})$. The definition of this element in $\mathrm{Conn}(E) \times C^\infty(Y;\mathbb{S})$ is given momentarily. The lemma that follows directly supplies input for the definition.

**Lemma C.1**: *Fix* $m \geq 1$. *There an m-dependent constant* $\kappa \geq 100$, *and given* $c_v \geq \kappa$, *there exists* $\kappa_{c_v} > \kappa$ *with the following significance: Suppose that* $r \geq \kappa_{c_v} c_v^{10}$ *and suppose that* $(A, \psi = (\alpha, \beta))$ *is a solution to the* $(r, \mu)$ *version of (1.13) with* $\mu$ *a given element in* $\Omega$ *with $\mathcal{P}$-norm smaller than 1. The parameters* $\kappa$, $c_v$ *and* $r$ *are suitable for use Lemma B.11 and in particular for constructing* $(A_\Diamond, \psi_\Diamond = (\alpha_\Diamond, \beta_\Diamond))$ *and the corresponding set* $\Theta^\alpha$. *Let* $\upsilon$ *denote a component of* $\alpha_\Diamond^{-1}(0)$ *that intersects* $M_\delta$ *and let* $\upsilon^\alpha$ *denote the corresponding element in* $\Theta^\alpha$. *There exists an isotopy from* $[0, 1] \times \upsilon$ *into* $Y$ *starting from* $\upsilon$, *ending at* $\upsilon^\alpha$ *and with the properties listed below. The list uses* $\upsilon^\alpha{}_\tau$ *to denote the* $\tau \in [0,1]$ *curve of the isotopy.*

- *Each point in* $\upsilon^\alpha{}_\tau$ *has distance at most* $m^{-1}$ *from the corresponding point in* $\upsilon$.
- *Each point in* $\upsilon^\alpha{}_\tau$ *has distance at least* $\kappa^{-1} c_v r^{-1/2}$ *from each curve in* $\cup_{p\in\Lambda}\{\hat{\gamma}_p^+ \cup \hat{\gamma}_p^-\}$.
- *The unit tangent vector to* $\upsilon^\alpha{}_\tau$ *has distance at most* $c_v r^{-1/2}$ *from* $v$; *and it has distance at most* $\kappa r^{-1/2}$ *from* $v$ *at the points where the distance is at least* $c_v r^{-1/2}$ *from each curve in* $\cup_{p\in\Lambda}\{\hat{\gamma}_p^+ \cup \hat{\gamma}_p^-\}$.
- *The push-forward via this isotopy of* $\frac{\partial}{\partial\tau}$ *is bounded by* $\kappa_{c_v}$.

This lemma is proved in Section Cb.

*Part 3*: Granted Lemma C.1, fix $\tau \in [1,2]$ so as to define $(A_{\Diamond 1\tau}, \psi_{\Diamond 1\tau})$. The definition of this pair is identical to that of $(A_\Diamond, \psi_\Diamond)$ given in Section Be but for the one change and one added remark. What follows directly is the one change to Section Be's definition. Let $\upsilon$ denote a given component of the zero locus of $\alpha_\Diamond$ that intersects $M_\delta$. By way of a reminder, $\upsilon$'s intersection with $Y_{*\Lambda} \cup T_{*\Lambda}$ is a union of components of $\alpha$'s zero locus in $Y_{*\Lambda} \cup T_{*\Lambda}$; and $\upsilon$'s intersection with any given component of $Y-(Y_{*\Lambda} \cup T_{*\Lambda})$ is described by Lemma C.1. This understood, replace $\upsilon$ in the formula that appear in Section Be with the corresponding curve $\upsilon^\alpha{}_\tau$ that is supplied by Lemma B.18.

The added remark addresses the need to specify an isomorphism between E and the product bundle over a certain neighborhood of each curve $\upsilon^\alpha{}_\tau$ and over the complement of the union of a certain smaller neighborhood about $\cup_{\upsilon^\alpha \in\Theta^\alpha} \upsilon^\alpha{}_\tau$ and a neighborhood of the components of the zero locus of $\alpha_\Diamond$ from $\cup_{p\in\Lambda}\{\hat{\gamma}_p^+ \cup \hat{\gamma}_p^-\}$. The required isomorphisms are already specified for the $\tau = 0$ case, these being the ones needed to define $(A_*, \psi_*)$ and $(A_\Diamond, \psi_\Diamond)$. The three steps that follow describe the $\tau > 0$ versions of these isomorphisms.

Step 1: Let $\upsilon^\alpha$ denote a given component of $\Theta^\alpha$. The definition of $(A_\Diamond, \psi_\Diamond)$ refered to larger and smaller neighborhoods of the corresponding curve $\upsilon$. Each $\tau \in [0,1]$ version of $\upsilon^\alpha{}_\tau$ has its analogous neighborhoods, these denoted by $U_{\upsilon,\tau}$ and $U_{\upsilon,\tau}{}'$. The set $U_{\upsilon,\tau}$ is a neighborhood of $\upsilon^\alpha{}_\tau$ that is defined as follows: Its intersection with $Y_{*\Lambda} \cup T_{*\Lambda}$ is the radius $4c_v{}^2 r^{-1/2}$ tubular neighborhood of $\upsilon^\alpha{}_\tau$. To describe the remainder of $U_{\upsilon,\tau}$, fix a component of $Y-(Y_{*\Lambda} \cup T_{*\Lambda})$ and let $\gamma$ denote for the moment the corresponding curve from $\cup_{p\in\Lambda}\{\hat{\gamma}_p^+ \cup \hat{\gamma}_p^-\}$. Reintroduce $\gamma$'s version of the coordinates $(t, z)$ that are used on this component to define $(A_\Diamond, \psi_\Diamond)$. Let T denote the $|z| \le (c_v{}^4 + c_v{}^3) r^{-1/2}$ part of the coordinate chart. The set $U_{\upsilon,\tau}$ intersects the $|z| \ge (c_v{}^4 - 2c_v{}^3) r^{-1/2}$ part of T as the radius $4c_v{}^2 r^{-1/2}$ tubular neighborhood of this part of $\upsilon^\alpha{}_\tau$. The intersection of $U_{\upsilon,\tau}$ with the rest of T is the radius $4c_v{}^{1/2} r^{-1/2}$ tubular neighborhood of this part of $\upsilon^\alpha{}_\tau$. The set $U_{\upsilon,\tau}{}'$ is defined analogously, but with the factor of 4 missing.

Introduce for each $\tau \in [0, 1]$ the neighborhood of $\upsilon^\alpha{}_\tau$ that is defined as follows: This neighborhood intersects the complement in Y of the radius $(c_v{}^4 - 3c_v{}^3) r^{-1/2}$ tubular neighborhoods of the curves from $\cup_{p\in\Lambda}\{\hat{\gamma}_p^+ \cup \hat{\gamma}_p^-\}$ as the tubular neighborhood of $\upsilon^\alpha{}_\tau$ with radius $8c_v{}^2 r^{-1/2}$ ; and it intersects the remaining part of Y as the concentric tubular neighborhood of $\upsilon^\alpha{}_\tau$ with radius $8\, c_v{}^{1/2} r^{-1/2}$. This neighborhood is denoted by $U_{\upsilon,\tau*}$. The set $U_{\upsilon,\tau}$ is a proper subset of $U_{\upsilon,\tau*}$.

Step 2: Fix an isomorphism between $K^{-1}|_\upsilon$ and $\upsilon \times \mathbb{C}$ that gives a version of the coordinates from Part 4 of Section Aa for $\upsilon$ with $|\nu| + |\mu| \le c_0$. The push forward of $\frac{d}{d\tau}$ by the map that defines Lemma C.1's isotopy gives a vector field along the image of the isotopy. Parallel transport along the integral curves of this vector field defines an isomorphism over any given $\tau \in [0,1]$ version of $\upsilon^\alpha{}_\tau$ between $K^{-1}$ and the product bundle. Use this isomorphism to define a $\upsilon^\alpha{}_\tau$ version of the coordinates from Part 4 of Section Aa. The associated pair $(\nu, \mu)$ is such that $|\nu| + |\mu| \le c_0$, this being a consequence of the fourth bullet in Lemma C.1.

Fix $\tau \in [0, 1]$. The pushforward of $\frac{d}{d\tau}|_\tau$ appears with respect to the $\upsilon^\alpha{}_\tau$ version of the $(t, z)$ coordinate chart as a vector that is defined at $z = 0$. View this vector as a vector field on $U_{\upsilon,\tau*}$ whose coefficients have no z-dependence. Use the function $\chi$ to extend the the latter vector field from $U_{\upsilon,\tau}$ to the rest of Y so as to be equal to 0 on the complement of $U_{\upsilon,\tau*}$ and so that its commutator with $\frac{\partial}{\partial z}$ on $U_{\upsilon,\tau*}$ is bounded by $\kappa_{c_v} r^{1/2}$ with $\kappa_{c_v}$ denoting the constant from Lemma C.1. The existence of an extension with this property follows from what is said by the fourth bullet of Lemma C.1. This extension is denoted in what follows by $v_{\upsilon,\tau}$. Use $v^\alpha$ to denote the vector field on $[0,1] \times Y$ that is defined by the rule $v_\alpha|_\tau = \frac{d}{d\tau} + \sum_{\upsilon\in\Theta_\alpha} v_{\upsilon,\tau}$.

Define $\pi_\alpha$: $[0,1] \times Y \to Y$ to be the map that sends any given point $(\tau, p)$ to the point in $\{0\} \times Y$ that lies on the integral curve through p of the vector field $v^\alpha$. The map $\pi_\alpha$ is a surjection that restricts to any $\tau \in [0,1]$ version of $\upsilon^\alpha{}_\tau$ as a diffeomorphism onto $\upsilon$.

<u>Step</u> <u>3</u>: Let $\pi$: $[0,1] \times Y \to Y$ denote the standard projection. The respective pull-backs of E by $\pi$ and $\pi_\alpha$ are isomorphic. Use $\varphi_\alpha$: $\pi^*E \to \pi_\alpha{}^*E$ to denote the isomorphism that is defined by parallel transport along the fibers of $\pi_\alpha$ by the connection $\pi_\alpha{}^*A_{\Diamond 1}$. The pull-back $\pi_\alpha{}^*\alpha_{\Diamond 1}$ defines a section of $\pi_\alpha{}^*E$ and so $\varphi_\alpha{}^{-1}(\pi_\alpha{}^*\alpha_{\Diamond 1})$ defines a section of $\pi^*E$. This section is denoted by $\hat\alpha_\Diamond$ and its restriction, $\hat\alpha_\Diamond|_\tau$, to any given constant $\tau$ slice is a section of E. The zero locus of the latter is $\cup_{\upsilon^\alpha\in\Theta^\alpha} \upsilon^\alpha{}_\tau$; it vanishes transversely with degree one on each component curve.

<u>Step</u> <u>4</u>: Fix $\tau \in [0,1]$ and introduce $U_{0,\tau}$ to denote the $\{\upsilon^\alpha{}_\tau\}_{\upsilon^\alpha\in\Theta^\alpha}$ version of the set $U_0$. This is the complement of $\cup_{\upsilon^\alpha\in\Theta^\alpha} U_{\upsilon,\tau}{}'$ and the union of the radius $c_0{}^{-1}c_v{}^{1/2}r^{-1/2}$ tubular neighborhoods of the curves from $\cup_{p\in\Lambda}\{\hat\gamma_p^+ \cup \hat\gamma_p^-\}$ in the zero locus of $\alpha_\Diamond$. The constructions that define $(A_{\Diamond 1\tau},\psi_{\Diamond 1\tau})$ require an isomorphism over $U_{0,\tau}$ between E and the product bundle. Use the isomorphism that sends $\hat\alpha_\Diamond|_\tau$ to $|\hat\alpha_\Diamond|_\tau|$.

The constructions that define $(A_{\Diamond 1\tau},\psi_{\Diamond 1\tau})$ also require an isomorphism between the bundle E and the product bundle over each set from the collection $\{U_{\upsilon,\tau}\}_{\upsilon^\alpha\in\Theta^\alpha}$. This isomorphism is defined using the $\upsilon^\alpha{}_\tau$ version of the coordinates (t, z). The desired isomorphism sends the section $\hat\alpha_\Diamond|_\tau$ to $|\hat\alpha_\Diamond|_\tau| \frac{z}{|z|}$.

*Part 4*: The next lemma compares $\mathfrak{f}_s$ at $(A_\Diamond,\psi_\Diamond)$ with $\mathfrak{f}_s$ at $(A_{\Diamond 1\tau=1},\psi_{\Diamond 1\tau=1})$.

**Proposition C.2**: *There exists* $\kappa \ge 100$, *and given* $c_v \ge \kappa$, *there exists* $\kappa_{c_v} > \kappa$ *with the following significance: Suppose that* $r \ge \kappa_{c_v} c_v{}^{10}$ *and suppose that* $(A,\psi = (\alpha, \beta))$ *is a solution to the* $(r, \mu)$ *version of (1.13) with* $\mu$ *a given element in* $\Omega$ *with* $\mathcal{P}$*-norm smaller than 1. The parameters* $\kappa$, $c_v$ *and* r *are suitable for use Lemma B.18 and in particular for constructing the path* $\{(A_{\Diamond 1\tau},\psi_{\Diamond 1\tau})\}_{\tau\in[0,1]}$. *The norm of the difference between the values of* $\mathfrak{f}_s$ *at* $(A_\Diamond,\psi_\Diamond)$ *and* $\mathfrak{f}_s$ *at* $(A_{\Diamond 1\tau=1},\psi_{\Diamond 1\tau=1})$ *is no greater than* $\kappa_{c_v}$.

***Proof of Proposition C.2***: The proof is much like that of Lemma B.17. In any event, there are five steps.

Step 1: Use the same arguments that prove Lemma B.12 to prove that Lemma B.12's assertion also holds for each $\tau \in [0,1]$ version of $(A_{\Diamond\tau}, \psi_{\Diamond\tau})$.

Step 2: For any given $\tau \in [0,1]$, use $\mathfrak{L}_{\mathbb{V}\tau}$ for the moment to denote the $(A_{\Diamond 1\tau}, \psi_{\Diamond\tau 1})$ version of the operator $\mathfrak{L}_{\mathbb{V}}$ that is depicted in (A.26) and (A.27). This family is not real analytic, but there are as small as desired perturbations that make it so, and this being the case, what is said in Part 1 of the proof of Proposition B.3 can be assumed to apply. Let $\{\lambda_{n\tau}\}_{n\in\mathbb{Z},\tau\in[0,1]}$ denote the corresponding set of eigenvalue families. The analog of (B.52) in this case reads

$$\frac{d}{d\tau}\lambda_{n\tau} = \Sigma_{\upsilon^\alpha\in\Theta^\alpha} \int_{U_{\upsilon,\tau}} \mathfrak{f}_{(\tau)}^{\dagger}(\tfrac{d}{d\tau}\mathfrak{L}_{\mathbb{V}\tau})\mathfrak{f}_{(\tau)} , \tag{C.1}$$

this because the $\tau$-derivative of $(A_{\alpha 1\tau}, \psi_{\alpha 1\tau})$ has support only in $\cup_{\upsilon^\alpha\in\Theta^\alpha} U_{\upsilon,\tau}$.

Step 3: It follows from what is said in the fourth bullet of Lemma B.18, in Lemma C.1 and in Part 3 that the endomomorphism $\frac{d}{d\tau}\mathfrak{L}_{\mathbb{V}\tau}$ of $\mathbb{V}_0\oplus\mathbb{V}_1$ has pointwise norm bounded by $\kappa_{c1} r^{1/2}$ where $\kappa_{c1}$ is a constant that is purely $c_v$ dependent. With this in mind, fix an integer $m \geq 1$ and let $I \in [0,1]$ denote an interval of length at most $m^{-1}$. The formula in (C.1) implies that $\lambda_{n(\cdot)}$ has a zero on I only if $|\lambda_{n\tau}| \leq m^{-1}\kappa_{c1} r^{1/2}$ for each $\tau \in I$.

This understood, it follows from the lemmas in Section Aa that if $c_v \geq c_0$ and r is greater that a purely $c_v$ dependent constant, then there is a second purely $c_v$ dependent constant $\kappa_{c2} > \kappa_{c1}$ with the following significance: Take $m > \kappa_{c2}$ and suppose that $n \in \mathbb{Z}$ is such that $\lambda_{n(\cdot)}$ has a zero in I. Fix $\tau \in I$ and use $\mathfrak{f}_{(\tau)}$ to denote an eigenvector of $\mathfrak{L}_{\mathbb{V}\tau}$ with eigenvalue $\lambda_{n\tau}$. Then $\|\Pi_\vartheta \mathfrak{f}_{(\tau)}\| \geq (1 - c_0 c_v^{-1})\|\mathfrak{f}_{(\tau)}\|_2$.

Granted the preceding, use Lemmas A.2, A.3, A.7 and A.8 with (A.28)-(A.30) to deduce the following: Suppose that $c_v \geq c_0$, r is greater that a purely $c_v$ dependent constant, and that $m$ is greater than yet another purely $c_v$ dependent constant. Suppose that $n \in \mathbb{Z}$ and $\lambda_{n(\cdot)}$ has a zero on I. Then $|\lambda_{n\tau}| \leq c_v^{-4}$ for all $\tau \in I$.

Step 4: Take $c_v$, r and $m$ so as to use what is said in Steps 1-3. Fix $\tau' > \tau \in [0,1]$ with $\tau' - \tau < m^{-1}$. Since $m$ need only be greater than a purely $c_v$ dependent constant, assume that it is no greater than this constant plus 1. The argument used in Part 4 of the proof of Proposition B.3 can be repeated with only notational changes to see that the norm of the difference between the values of $\mathfrak{f}_{\mathfrak{s}}$ at $(A_{\Diamond 1\tau}, \psi_{\Diamond 1\tau})$ and at $(A_{\Diamond 1\tau'}, \psi_{\Diamond 1\tau'})$ is at most $\kappa_c$ with $\kappa_c$ being a purely $c_v$ dependent constant. This conclusion implies what is asserted by Proposition B.19 as $[0,1]$ can be covered by $2m$ intervals of length less than $m^{-1}$.

**b) The proof of Lemma C.1**

The proof has fourteen steps.

Step 1: Fix $\mathfrak{p} \in \Lambda$ such that $\upsilon$ crosses $\mathcal{H}_{\mathfrak{p}}$. The curve $\upsilon$ crosses the $u = R + \ln\delta$ sphere in $\mathcal{H}_{\mathfrak{p}}$ and shortly intersect the $f = 1 + \delta^2$ surface in $\mathcal{H}_{\mathfrak{p}}$ as it continues out of $\mathcal{H}_{\mathfrak{p}}$ to cross $M_\delta$. Let $z_{\mathfrak{p}+}$ denote this intersection point. By way of a reminder, the function $f$ where $u \geq R + \ln\delta$ in $\mathcal{H}_{\mathfrak{p}}$ is given by $f = 1 + e^{2(R-u)}(1 - 3\cos^2\theta)$. The point $z_{\mathfrak{p}+}$ is the starting point of a component of the $f \in (1 + \delta^2, 2 - \delta^2)$ part of $\upsilon \cap M_\delta$. The ending point of this component lies on the $f = 2 - \delta^2$ surface in one of the handles $\{\mathcal{H}_{\mathfrak{p}'}\}_{\mathfrak{p}' \in \Lambda}$. Let $\mathfrak{p}' \in \Lambda$ denote the relevant pair and let $z_{\mathfrak{p}'-}$ denote this ending point on the $f = 2 - \delta^2$ surface in $\mathcal{H}_{\mathfrak{p}'}$.

By way of a reminder from Part 2 in Section II.1c, the index 1 critical point from $\mathfrak{p}$ has an ascending disk in $M_\delta$ that intersects the Heegaard surface $\Sigma$ as a smoothly embedded circle, this denoted by $C_{\mathfrak{p}+}$; and the index 2 critical point from $\mathfrak{p}'$ has a descending disk in $M_\delta$ that intersects the Heegaard surface $\Sigma$ as a smoothly embedded circle, this denoted by $C_{\mathfrak{p}'-}$. The segment of $\upsilon$ that starts at $z_{\mathfrak{p}+}$ and ends at $z_{\mathfrak{p}'-}$ intersects $\Sigma$ at a point with distance $c_0^{-1}$ or less from a point in $C_{\mathfrak{p}+} \cap C_{\mathfrak{p}-}$. Use $z_\upsilon$ for this point in $\upsilon \cap \Sigma$ and use $z_*$ for the nearby point in $C_{\mathfrak{p}+} \cap C_{\mathfrak{p}'-}$. The point $z_\upsilon$ is well inside a certain coordinate neighborhood of $z_*$. This neighborhood has coordinates $(\varphi, \hbar)$ which are defined where $|\varphi|^2 + |\hbar|^2$ is bounded by a constant that depends only on the geometry of M. The pair $(\varphi, \hbar)$ are such that $z = \varphi + i\hbar$ is a holomorphic coordinate for the neighborhood.

Lie transport by $v$ of the functions $(\varphi, \hbar)$ along $v$'s integral curves defines coordinates $(t, \varphi, \hbar)$ for a closed cylinder in $M_\delta$ with t being the value of $f$ along the integral curves of $v$. The coordinate t is restricted to the interval $[1 + \delta^2, 2 - \delta^2]$. The corresponding $t = 1 + \delta^2$ boundary disk of the cylinder is a disk in the $u > R + \ln\delta$ part of $\mathcal{H}_{\mathfrak{p}}$. The function $\varphi$ on this boundary disk is such that $d\varphi = d\phi$. The function $\hbar$ on this disk is the function $e^{-2(R-u)}\cos\theta\sin^2\theta$. The $t = 2 - \delta^2$ boundary disk of this coordinate cylinder is a disk in the $u \leq -R - \ln\delta$ part of $\mathcal{H}_{\mathfrak{p}'}$. The function $\varphi$ on this boundary disk is is either $e^{-2(R+u)}\cos\theta\sin^2\theta$ or it is $-e^{-2(R+u)}\cos\theta\sin^2\theta$. In the former case, $d\hbar = d\phi$ on this boundary disk; and $d\hbar = -d\phi$ in the latter case.

The segment of $\upsilon$ between $z_{\mathfrak{p}+}$ and $z_{\mathfrak{p}-}$ is in this coordinate cylinder and as such, it appears as the graph of the form $t \to (t, z = z_\upsilon(t))$. The function $z_\upsilon(\cdot)$ solves the $\tau = 0$ version of the $\tau \in [0, 1]$ family of differential equations depicted in the upcoming (C.2). The depiction of this family introduces a certain $\mathbb{C}$ valued function, $x_\upsilon$, on $[1 + \delta^2, 2 - \delta^2]$ with norm bounded by $c_0 r^{-1/2}$. A given $\tau \in [0, 1]$ member of the family requires a $\mathbb{C}$ valued function of t to obey

$$\tfrac{i}{2}\tfrac{d}{dt} z + (1 - \tau)\, x_\upsilon = 0 \tag{C.2}$$

for $t \in [1 + \delta^2, 2 - \delta^2]$. Given $z_0 \in \mathbb{C}$ with norm bounded by $c_0^{-1}$, integration finds a unique solution to (C.2) with $z(1 + \delta^2) = z_0$. There is also a unique solution with $z(2 - \delta^2) = z_0$. In either case, the solution obeys $|z(\cdot) - z_0| \le (1 - \tau)\, c_0 r^{-1/2}$.

Step 2: Fix $\varepsilon_1 \in (0, c_0^{-1})$. This section uses $c_\varepsilon$ to denote a purely $\varepsilon_1$ dependent constant that is greater than 1 and whose value can be assumed to increase on successive appearances.

Fix $r \ge c_\varepsilon$. Suppose that $\mathfrak{p} \in \Lambda$ is such that $\upsilon$ crosses $\mathcal{H}_p$. Assume in addition that each point of $\upsilon$ has distance $\varepsilon_1$ or greater from both $\hat{\gamma}_p^+$ and $\hat{\gamma}_p^-$. This being the case, $\upsilon$ coincides with a segment in $\mathcal{H}_p$ of $\alpha$'s zero locus and so its tangent vector here has distance at most $c_0 r^{-1/2}$ from $v$. Reintroduce $z_{p-}$ to denote the point on $\upsilon \cap \mathcal{H}_p$ where $\upsilon$ intersects the $e^{-2(R+u)}(1 - 3\cos^2\theta) = \delta^2$ surface in $\mathcal{H}_p$ and introduce again $z_{p+}$ to denote the point where $\upsilon$ intersects the $e^{-2(R-u)}(1 - 3\cos^2\theta) = \delta^2$ surface in $\mathcal{H}_p$.

Let $\gamma$ denote the segment in $\mathcal{H}_p$ of the integral curve of $v$ that starts at $z_{p-}$ and lies in the $e^{-2(R-|u|)}(1 - 3\cos^2\theta) \le \delta^2$ part of $\mathcal{H}_p$. Given that $\upsilon$'s tangent vector differs from $v$ by at most $c_0 r^{-1/2}$, the $e^{-2(R-|u|)}(1 - 3\cos^2\theta) \le \delta^2$ part of $\upsilon$ in $\mathcal{H}_p$ lies entirely in the radius $c_\varepsilon r^{-1/2}$ tubular neighborhood of $\gamma$. The function $1 - 3\cos^2\theta$ is positive on $\gamma$ if $r \ge c_\varepsilon^{-1}$ and the segment $\gamma$ ends on the $e^{-2(R-u)}(1 - 3\cos^2\theta) = \delta^2$ surface in $\mathcal{H}_p$. In fact, the radius $c_\varepsilon^{-1}$ tubular neighborhood of $\gamma$ lies entirely in the $1 - 3\cos^2\theta > 0$ part of $\mathcal{H}_p$ and its boundary consists of one disk on the $e^{-2(R+u)}(1 - 3\cos^2\theta)$ surface and the other on the $e^{-2(R-u)}(1 - 3\cos^2\theta)$ surface. This neighborhood has coordinates (t, z) as described in Part 4 of Section Aa with $|\nu| + |\mu| \le c_0$ and with the $t = 0$ point being the $e^{-2(R+u)}(1 - 3\cos^2\theta)$ point on $\gamma$.

The segment of $\upsilon$ in $\mathcal{H}_p$ between $z_{p-}$ and $z_{p+}$ appears in these coordinates as the graph $t \to (t, z = z_\upsilon(t))$ where $z_\upsilon(0) = 0$. The function $z_\upsilon(\cdot)$ is a solution to the $\tau = 0$ member of a $\tau \in [0,1]$ family of differential equations for a $\mathbb{C}$ valued function of t; this being an equation of the form

$$\tfrac{i}{2}\tfrac{d}{dt} z + \nu z + \mu \bar{z} + (1 - \tau)\, x_\upsilon + \mathfrak{e}(z) = 0\ , \tag{C.3}$$

where $\mathfrak{e}$ is a smooth function on the radius $c_0^{-1}$ ball in $\mathbb{C}$ centered at the origin with the property that $|\mathfrak{e}| \le c_0|z|^2$ and $|d\mathfrak{e}| \le c_0 |z|$. Meanwhile, $x_\upsilon$ is a smooth function of t that obeys $|x_\upsilon| \le c_0 r^{-1/2}$.

Fix $\tau \in [0, 1]$ and $z_0 \in \mathbb{C}$ in the domain of $\mathfrak{e}$ with $|z_0| \le c_0^{-1}$. Then there is a unique solution to $\tau$'s version of (C.3) that is defined on a neighborhood of 0 with $t = 0$ value $z_0$. Let $z(\cdot)$ denote this solution. If $|z_0| \le c_\varepsilon^{-1}$, then this $z(\cdot)$ will be defined for all values of t and it will obey $|z(\cdot)| \le c_\varepsilon |z_0|$. The solution depends smoothly on the data $(\tau, z_0)$. The solutions to the $\tau = 1$ version of (C.3) are the segments of the integral curves of $v$ in the $e^{-2(R-|u|)}(1 - 3\cos^2\theta)$ part of $\mathcal{H}_p$ that start at distances less than $c_\varepsilon^{-1}$ from $z_{p-}$.

Step 3: The observations in Steps 1 and 2 suggest the lemma that follows.

**Lemma C.3**: *Given* $\varepsilon > 0$ *there exists* $\kappa_\varepsilon > 1$ *with the following significance: Fix* $r \ge \kappa_\varepsilon$ *and suppose that* $(A, \psi = (\alpha, \beta))$ *is a solution to the* $(r, \mu)$ *version of (1.13) with* $\mu$ *a given element in* $\Omega$ *with* $\mathcal{P}$*-norm smaller than 1. Let* $\upsilon$ *denote a component of* $\alpha^{-1}(0)$ *whose points have distance* $\varepsilon$ *or more from each curve in* $\cup_{p\in\Lambda}\{\hat{\gamma}_p^+ \cup \hat{\gamma}_p^-\}$*. Then* $\upsilon$ *is in the radius* $\kappa_\varepsilon r^{-1/2}$ *tubular neighborhood of closed, integral curve of* $v$*.*

***Proof of Lemma C.3***: The proof also uses $c_\varepsilon$ to denote a purely $\varepsilon$-dependent constant that is greater than 1. The value of $c_\varepsilon$ can be assumed to increase between successive appearances. Fix a point $p \in \upsilon \cap \Sigma$ and use what is said in Steps 1 and 2 to construct a segment of an integral curve of $v$ that starts at p, ends at a point $p' \in \Sigma$ with distance at most $c_\varepsilon r^{-1/2}$ from p and is such that $\upsilon$ lies in its radius $c_\varepsilon r^{-1/2}$ tubular neighborhood. Let $\gamma_p$ denote this segment. With this in mind, the arguments used in Step 4 of the proof of Proposition II.2.7 can be used with only cosmetic modifications to prove that $\gamma_p$ has distance at most $c_\varepsilon r^{-1/2}$ from a closed integral curve of $v$.

What follows directly is a proof of Lemma C.1 in the case when $\upsilon$ obeys the assumptions of Lemma C.3 for a given $\varepsilon$. To start, let $\gamma_0$ denote now the closed integral curve of $v$ that is supplied by Lemma C.3. If $r \ge c_\varepsilon$, then the radius $4\kappa_\varepsilon$ tubular neighborhood of $\gamma_0$ will intersect each $\mathfrak{p} \in \Lambda$ version of $\mathcal{H}_p$ only where $1 - 3\cos^2\theta > 0$. Keeping this in mind, use Part 4 of Section Aa to define coordinates (t, z) for this tubular neighborhood of $\gamma_0$ with $|\nu|$ and $|\mu|$ bounded by $c_0$. The curve $\upsilon$ appears in these coordinates as the graph of a map $t \to z(t)$ with $|z| \le \kappa_\varepsilon r^{-1/2}$ and with $|\frac{d}{dt} z| \le c_0 r^{-1/2}$. Define the family $\{\upsilon^\alpha{}_\tau\}_{\tau\in[0,1]}$ by writing any given member as the graph of the map $t \to (1 - \tau) z(\cdot)$.

Step 4: Let $\gamma$ denote $\hat{\gamma}_p^+$. Use $\theta_* \in (0, \frac{\pi}{2})$ in what follows to denote the angle with $\cos\theta_* = \frac{1}{\sqrt{3}}$. Use b to denote $\frac{3}{2\sqrt{2}} e^R (x_0 + 4e^{-2R})^{1/2}$. Fix $\varepsilon \in (0, c_0^{-1})$ with the upper bound chosen so that the $\mathbb{R}/(2\pi\mathbb{Z})$ valued function $\phi$ and the pair $(x = b^{-1}u, y = \theta - \theta_*)$ define coordinates on the radius $\varepsilon$ tubular neighborhood of $\gamma$. Let $p = y + x$ and $q = y - x$ and fix $\varepsilon_0 < \varepsilon$ so that the locus where $(p^2 + q^2)^{1/2} = \varepsilon_0$ lies in the radius $\varepsilon$ tubular neighborhood of $\gamma$. The notation that follows uses $c_\varepsilon$ to denote a constant that is greater than 1 and depends solely on $\varepsilon_0$. Its value can be assumed to increase between successive appearances.

Fix $\varepsilon_1$ and r as in Step 2 with $\varepsilon_1$ chosen so that $(p^2 + q^2) \le \frac{1}{8}\varepsilon_0$ on the radius $\varepsilon_1$ tubular neighborhood of $\gamma$. Suppose that $\mathfrak{p} \in \Lambda$ is such that $\upsilon \cap \mathcal{H}_p \neq \emptyset$ but assume in this case that $\upsilon$ has boundary points on the radius $\varepsilon_1$ tubular neighborhood of either $\hat{\gamma}_p^+$. Much the same argument holds if $\upsilon$ has boundary points on $\hat{\gamma}_p^-$ and so the latter case will not be discussed.

The part of $\upsilon \cap \mathcal{H}_p$ where $e^{-2(R-|u|)}(1 - 3\cos^2\theta) \le \delta^2$ but not in the $(p^2 + q^2)^{1/2} < \frac{1}{2}\varepsilon_0$ part of radius $\varepsilon$ tubular neighborhood of $\gamma$ consists of two segments, these denoted $\upsilon_-$ and $\upsilon_+$ in what follows. The function u is negative on $\upsilon_-$ and positive on $\upsilon_+$. Both segments lie in the zero locus of $\alpha$ and have transversal intersection with the $(p^2+q^2)^{1/2} = \frac{1}{2}\varepsilon_0$ locus. The starting point of $\upsilon_-$ is $z_{p-}$. Use $z_-$ to denote the point of $\upsilon_-$ on the $(p^2+q^2)^{1/2} = \varepsilon_0$ surface in $\gamma$'s radius $\varepsilon$ tubular neighborhood.

The u < 0 part of the segment of the integral curve of $v$ in $\mathcal{H}_p$ that contains $z_-$ and lies where $e^{-2(R+u)}(1 - 3\cos^2\theta) \le \delta^2$ will start on the $e^{-2(R+u)}(1 - 3\cos^2\theta) = \delta^2$ surface at distance $c_\varepsilon r^{-1/2}$ or less from $z_{p-}$ and it will intersect the $(p^2+q^2)^{1/2} = \frac{1}{2}\varepsilon_0$ surface in the radius $\varepsilon$ tubular neighborhood of $\gamma$. Introduce $\gamma_-$ to denote the segment of the u < 0 part of this integral curve that runs between its $e^{-2(R+u)}(1 - 3\cos^2\theta) = \delta^2$ point and its intersection with the $(p^2+q^2)^{1/2} = \frac{1}{2}\varepsilon_0$ surface in the radius $\varepsilon$ tubular neighborhood of $\gamma$. The part of $\upsilon_-$ that lies outside the locus where $(p^2+q^2) < \frac{1}{2}(1 + c_\varepsilon r^{-1/2})\varepsilon_0$ is in the radius $c_\varepsilon r^{-1/2}$ tubular neighborhood of $\gamma_-$.

Fix coordinates (t, z) for the radius $c_0^{-1}$ tubular neighborhood of $\gamma_-$ from Part 4 of Section Aa with the z = 0 locus being $\gamma_-$ and the t = 0 point being $z_-$. Require in addition that $|\nu|$ and $|\mu|$ are bounded by $c_0$. The intersection of this tubular neighbhorhood with the $e^{-2(R+u)}(1 - 3\cos^2\theta) = \delta^2$ surface is a disk neighborhood of boundary point of $\gamma_-$ in this surface. If the radius of this tubular neighborhood is less than $c_\varepsilon$, then its intersection with the surfaces in the radius $\varepsilon$ tubular neighborhood of $\gamma$ where $(p^2 + q^2)^{1/2}$ is constant and between $\frac{3}{4}\varepsilon_0$ and $2\varepsilon_0$ are disks that lie in the u < 0 part of these surfaces.

The segment $\upsilon_-$ appears in the coordinates (t, z) as a graph $t \to (t, z(t))$ with $z(t)$ obeying the $\tau = 0$ version of a $\tau \in [0,1]$ family of equations that has the same form as that depicted in (C.3). This solution has $z(0) = 0$ and $|z(\cdot)| \le c_\varepsilon r^{-1/2}$. Note that the solutions to the $\gamma_-$ and $\tau = 1$ version of (C.3) are integral curves of $v$.

Solutions to (C.3) for all values of $\tau$ can readily be found. In particular, there exists a purely $\varepsilon_0$ dependent constant, $c_{p\varepsilon}$, that is greater than 1 and has the following significance: Fix $\tau \in [0,1]$ and a point $z_0$ in the $(p^2+q^2) = \varepsilon_0$ surface with distance less that $c_{p\varepsilon}^{-1}$ from $z_-$. Use $\Delta$ to denote this distance. There is a unique solution to the $\gamma_-$ version of (C.3) for the chosen value of $\tau$ that contains $z_0$ and with norm bounded for all t by $c_\varepsilon(\Delta + r^{-1/2})$. Moreover, varying the data $(\tau, z_0)$ changes the corresponding solution in a smooth fashion; and the three parameter family of solutions so defined is such that the derivative of $z(\cdot)$ with respect to changes of $(\tau, z_0)$ is bounded at each t by $c_\varepsilon$.

Step 5: This step uses the same notation as Step 4. Use $z_+$ to denote the point of $\upsilon_+$ on the locus $(p^2+q^2)^{1/2} = \varepsilon_0$. The end point of $\upsilon_+$ is on the $e^{-2(R-u)}(1-3\cos^2\theta) = \delta^2$ surface; this being the point $z_{p+}$. If $r \geq c_\varepsilon$ then the segment of the integral curve of $v$ in $\mathcal{H}_p$ that contains $z_+$ will intersect the surface $(p^2+q^2)^{1/2} = \frac{1}{2}\varepsilon_0$ where $u > 0$ at a point with distance at most $c_\varepsilon r^{-1/2}$ from the point where $\upsilon_+$ intersects this surface. It will also intersect the surface where $e^{-2(R-u)}(1-3\cos^2\theta) = \delta^2$. Introduce $\gamma_+$ to denote the segment of this integral curve of $v$ that starts on the $(p^2+q^2)^{1/2} = \frac{1}{2}\varepsilon_0$ surface in the radius $\varepsilon$ tubular neighborhood of $\gamma$ and ends on the $e^{-2(R-u)}(1-3\cos^2\theta) = \delta^2$ surface. The radius $c_\varepsilon^{-1}$ tubular neighborhood of $\gamma_+$ will intersect this surface in a disk, and it will intersect each surface in the radius $\varepsilon$ tubular neighborhood of $\gamma$ where $(p^2 + q^2)^{1/2}$ is constant and between $\frac{3}{4}\varepsilon_0$ and $2\varepsilon_0$ as a disk in the $u > 0$ part of the surface in question.

Fix coordinates (t, z) for the radius $c_\varepsilon^{-1}$ tubular neighborhood of $\gamma$ from Part 4 of Section Aa with $|\nu|$ and $|\mu|$ bounded by $c_0$ and with the $t = 0$ point being the point $z_+$. The segment $\upsilon_+$ appears in the coordinates (t, z) as a graph $t \to (t, z(t))$ with $z(t)$ obeying the $\tau = 0$ version of a $\tau \in [0,1]$ family of equations that has the same form as that depicted in (C.3). This solution has $z(0) = 0$ and $|z(\cdot)| \leq c_\varepsilon r^{-1/2}$. The solutions to the $\gamma_+$ and $\tau = 1$ version of (C.3) are integral curves of $v$.

The constant $c_{p\varepsilon}$ from Step 4 can be chosen so that there is a $\gamma_+$ analog of what is said in the final paragraph of Step 4. This is to say that the following is true: Choose any $\tau \in [0,1]$ and a point $z_0$ in the $(p^2+q^2) = \varepsilon_0$ surface with distance less that $c_{p\varepsilon}^{-1}$ from $z_+$. Use $\Delta$ to denote this distance. There is a unique solution to the $\gamma_+$ version of (C.3) for the chosen value of $\tau$ that contains $z_0$ and has norm bounded for all t by $c_\varepsilon(\Delta + r^{-1/2})$. Varying the data $(\tau, z_0)$ changes the corresponding solution in a smooth fashion and the three parameter family of solutions so defined is such that the derivative of $z(\cdot)$ with respect to changes of $(\tau, z_0)$ is bounded at each t by $c_\varepsilon$.

Step 6: Fix $\varepsilon \in (0, c_0^{-1})$ and suppose that $\upsilon$ intersects the radius $\varepsilon$ tubular neighborhood of either $\hat{\gamma}_p^+$ or $\hat{\gamma}_p^-$. What follows considers the case when the curve in question is $\hat{\gamma}_p^+$. As the same analysis holds for the other case modulo some sign changes, the latter case is not discussed. Use $\gamma$ now to denote $\hat{\gamma}_p^+$ and let $\theta_* \in (0, \frac{\pi}{2})$ denote the angle with $\cos\theta_* = \frac{1}{\sqrt{3}}$, this the value of $\theta$ on $\gamma$. Coordinates for a neighborhood of $\gamma$ are given by the $\mathbb{R}/(2\pi\mathbb{Z})$ function $\phi$ and $\mathbb{R}$-valued functions (x, y) that are defined by the rules ($u = bx$, $\theta = \theta_* + y$) where $b = \frac{3}{2\sqrt{2}} e^R (x_0 + 4e^{-2R})^{1/2}$ is the constant that appears in (B.36) and in Steps 4 and 5. Introduce as in these same steps functions p and q given by $p = y + x$ and $q = y - x$.

If υ's intersection with the radius ε tubular neighborhood of γ lies entirely in α's zero locus, then the part of the curve υ where $(p^2 + q^2)^{1/2} \le c_0^{-1}\varepsilon$ can be parametrized by an interval $I \subset \mathbb{R}$ as a map $t \to (\phi = -t, p = p_\upsilon(t), q = q_\upsilon(t))$. No generality is lost in this case by taking I to contain the origin $0 \in \mathbb{R}$ and to take $t = 0$ to be the $u = 0$ point on υ. Thus $p_\upsilon(0) = q_\upsilon(0)$. If $(p^2+q^2)^{1/2} \le m^{-2} c_v r^{-1/2}$ on $\upsilon \cap \mathcal{H}_p$, then it follows from what is said in Step 3 of the proof of Lemma B.11 that the part of the curve υ where $(p^2 + q^2)^{1/2} \le c_0^{-1}\varepsilon$ can also be parametrized by and interval $I \subset \mathbb{R}$ containing 0; this map has again the form $t \to (\phi = -t, p = p_\upsilon(t), q = q_\upsilon(t))$. In this case the 0 point in I is taken as in Step 2 of the proof of Lemma B.11. It follows from (B.42) that $|p_\upsilon(0) - q_\upsilon(0)| \le c_0 m^{-4} c_v r^{-1/2}$.

In all cases, the functions $p_\upsilon$ and $q_\upsilon$ obey an equation of the form that is depicted in (B.37). This equation is reproduced below:

$$\frac{d}{dt} p_\upsilon = \lambda p_\upsilon + \mathfrak{e}_p(p_\upsilon, q_\upsilon) + \mathfrak{r}_{p\upsilon} \quad \textit{and} \quad \frac{d}{dt} q_\upsilon = -\lambda q_\upsilon + \mathfrak{e}_q(p_\upsilon, q_\upsilon) + \mathfrak{r}_{q\upsilon} \,. \tag{C.4}$$

By way of a reminder, $\lambda = 4\sqrt{6}\, e^{-R} (x_0 + 4 e^{-2R})^{1/2}$ and the functions $\mathfrak{e}_p$ and $\mathfrak{e}_q$ are smooth and have absolute value bounded by $c_0 (p^2 + q^2)$. Meanwhile, $\mathfrak{r}_{p\upsilon}$ and $\mathfrak{r}_{q\upsilon}$ are smooth functions of t. Their absolute values are bounded at times $t \in I$ where $\upsilon(t) \subset \alpha^{-1}(0)$ by $c_0 r^{-1/2}$, in particular this occurs where $(p^2+q^2)^{1/2} \ge c_v r^{-1/2}$. In general, their absolute values are bounded by $c_0 c_v r^{-1/2}$. A smaller upper bound is given the upcoming (C.5).

What follows says more about $\mathfrak{r}_{p\upsilon}$ and $\mathfrak{r}_{q\upsilon}$ at times $t \in I$ where their absolute value is greater than $c_0 r^{-1/2}$. To this end, reintroduce the constant $m$ that is used to define υ; and reintroduce $t_\Diamond \in [0, 2\pi)$ from Step 3 of Lemma B.11. As done in this same Step 3, take the parametrization for I so that it is only necessary to consider times $t \in [-2\pi - t_\Diamond, 2\pi + t_\Diamond]$. The pair $p_\upsilon$ and $q_\upsilon$ on this interval are given in (B.42). Differentiate (B.42) and compare with (C.4) to see that $\mathfrak{r}_{p\upsilon}$ and $\mathfrak{r}_{q\upsilon}$ obey

- $\mathfrak{r}_{p\upsilon} \le c_0\, m^{-6} c_v r^{-1/2}$ *for* $t \in [0, 2\pi + t_\Diamond]$ *and* $\mathfrak{r}_{p\upsilon} \le c_0\, m^{-2} c_v r^{-1/2}$ *for* $t \in [-2\pi - t_\Diamond, 0]$.
- $\mathfrak{r}_{q\upsilon} \ge -c_0\, m^{-6} c_v r^{-1/2}$ *for* $t \in [-2\pi - t_\Diamond, 0]$ *and* $\mathfrak{r}_{q\upsilon} \ge -c_0\, m^{-2} c_v r^{-1/2}$ *for* $t \in [0, 2\pi + t_\Diamond]$ .
- $\mathfrak{r}_{p\upsilon} \ge -c_0\, m^{-6} c_v r^{-1/2}$ *and* $r_{q\upsilon} \le c_0\, m^{-6} c_v r^{-1/2}$ *for all* $t \in [-2\pi - t_\Diamond, 2\pi + t_\Diamond]$.

(C.5)

The constant $m$ is left unspecified for now but ultimately chosen to be less than $c_0$. The choice of $m$ determines in part a lower bound for Lemma C.1's constant κ.

To say something about the respective lengths of the $t > 0$ and $t < 0$ parts of I, introduce Δ to denote value of $(p_\upsilon^2 + q_\upsilon^2)^{1/2}$ at the $0 \in I$. Fix $\varepsilon_0 \in (0, \varepsilon)$ so that the coordinates functions p and q are defined where $(p^2+q^2)^{1/2} < 2\varepsilon_0$. Let $t_+$ and $t_-$ denote the respective values of t in I where $(p_\upsilon{}^2 + q_\upsilon{}^2)^{1/2} = \varepsilon_0$. As explained in Step 6,

$$|t_+ - \lambda^{-1}\ln(\varepsilon_0\Delta^{-1})| \le c_0 \quad \textit{and} \quad |t_- + \lambda^{-1}\ln(\varepsilon_0\Delta^{-1})| \le c_0 \tag{C.6}$$

if $\varepsilon_0 \le c_0{}^{-1}$. These bounds imply in part that the lengths of the $t > 0$ and $t < 0$ parts of I differ by at most $c_0$.

Step 7: Introduce $Y_0$ and $X_0$ to denote the value of $y = \frac{1}{2}(p+q)$ and $x = \frac{1}{2}(p-q)$ at the $t = 0$ point on $\upsilon$. Note that $X_0 = 0$ if $(p^2+q^2) \ge m^{-2}c_v r^{-1/2}$ on $\upsilon$; and (B.39) and (B.41) imply that $|X_0| \le c_0 m^{-6} c_v r^{-1/2}$ otherwise. Meanwhile, $Y_0 = 2^{-1/2}\Delta$ if $(p^2+q^2) \ge m^{-2}c_v r^{-1/2}$ on $\upsilon$ and $|Y_0 - 2^{-1/2}m^{-2}c_v r^{-1/2}| \le c_0 m^{-4} c_v r^{-1/2}$ otherwise.

Fix $Y \in (\frac{1}{4}Y_0, 4Y_0)$. Given $\tau \in [0, 1]$, there is a unique map $t \to (p_{Y,\tau}(t), q_{Y,\tau}(t))$ from a maximal interval $I_{Y,\tau} \subset \mathbb{R}$ to $\mathbb{R}^2$ that obeys the equation

$$\frac{d}{dt}p_{Y,\tau} = \lambda p_{Y,\tau} + \mathfrak{e}_p(p_{Y,\tau},q_{Y,\tau}) + (1-\tau)\mathfrak{r}_{p\upsilon} \quad \textit{and} \quad \frac{d}{dt}q_{Y,\tau} = -\lambda q_{Y,\tau} + \mathfrak{e}_q(p_{Y,\tau},q_{Y,\tau}) + (1-\tau)\mathfrak{r}_{q\upsilon} \tag{C.7}$$

with $p_{Y,\tau}(0) = Y+X_0$ and $q_{Y,\tau}(0) = Y-X_0$ and with $(p_{Y,\tau}{}^2 + q_{Y,\tau}{}^2)^{1/2} \le \varepsilon_0$ for $t \in I_{Y,\tau}$ with equality only at each boundary point of $I_{Y,\tau}$. A proof of existence and uniqueness can be had using standard techniques from the theory of ordinary differential equations. Use $t_{Y,\tau+}$ and $t_{Y,\tau-}$ to denote the respective negative and positive endpoints of $I_{Y,\tau}$.

The first implication of (C.7) concerns the size of $q_{Y,\tau}$ relative to $p_{Y,\tau}$ where $t \ge 0$: As explained momentarily,

$$|q_{Y,\tau}| \le c_0(\varepsilon_0 p_{Y,\tau} + m^{-2} c_v r^{-1/2}) \tag{C.8}$$

for $t \in [0, t_{Y,\tau+}]$ when $\varepsilon_0 \le c_0{}^{-1}$ and $m > c_0$. To prove this, fix $\varsigma > 0$ and set $w = |q| - \varsigma|p|$. It follows from (C.5) and (C.7) that

$$\frac{d}{dt}w \le -\lambda w - 2\lambda\varsigma|p_{Y,\tau}| + c_0\varepsilon_0(|w| + |p_{Y,\tau}|) + c_0 m^{-2} c_v r^{-1/2}\ . \tag{C.9}$$

Take $\varsigma = c_0\lambda^{-1}\varepsilon_0$ to see that $\frac{d}{dt}w \le -c_0{}^{-1}w + c_0 m^{-2}c_v r^{-1/2}$ if $\varepsilon_0 \le c_0{}^{-1}$. Multiply both sides of this last inequality by $e^{t/c_0}$ and integrate to obtain (C.8).

With regards to $p_{Y,\tau}$, note first that (C.7) with (C.5) imply that $p_{Y,\tau}$ is an increasing function of t when t is positive. To say more about the size of $p_{Y,\tau}$ it proves useful to introduce the norm on $C^\infty([0, t]; \mathbb{R})$ for $t \le t_{Y,\tau+}$ given by $h \to \|h\|_t = \sup_{s\in[0,t]} e^{-\lambda s}|h(s)|$. Use (C.8) with (C.5) and the right hand equation in (C.7) to see that

$$|p_{Y,\tau} - e^{\lambda t}(Y+X_0)| \le c_0(e^{\lambda t}\|p_{Y,\tau}\|_t^2 + m^{-6}c_v r^{-1/2}) .$$

(C.10)

Given that $Y+X_0 \ge c_\varepsilon^{-1} m^{-2} c_v r^{-1/2}$, this last equation implies that

$$(1 - c_0\varepsilon_0) e^{\lambda t}(Y+X_0 - c_0 m^{-6} c_v r^{-1/2}) \le p_{Y,\tau}(t) \le (1+c_0\varepsilon) e^{\lambda t}(Y+X_0 + c_0 m^{-6} c_v r^{-1/2})$$

(C.11)

for $t \in [0, t_{Y,\tau+}]$ when $\varepsilon_0 \le c_0^{-1}$ and $m \ge c_0$. Note that (C.11) with (C.8) implies that the function $1 - 3\cos^2\theta$ is positive along the trajetory $t \to (p_{Y,\tau}(t), q_{Y,\tau}(t))$ where $(p_{Y,\tau}^2 + q_{Y,\tau}^2)^{1/2}$ is greater than $c_0 m^{-2} c_v r^{-1/2}$.

The same sort of arguments can be used for $t \in [-2\pi - t_{D,\tau-}, 0]$ to see that

- $|p_{Y,\tau}| \le c_0(\varepsilon_0 |q_{Y,\tau}| + m^{-2} c_v r^{-1/2})$ ,
- $(1 - c_0\varepsilon_0) e^{-\lambda t}(Y - X_0 - c_0 m^{-6} c_v r^{-1/2}) \le q_{Y,\tau}(t) \le (1+c_0\varepsilon) e^{-\lambda t}(Y - X_0 + c_0 m^{-6} c_v r^{-1/2})$ ,

(C.12)

for $t \le [t_{Y,\tau-}, 0]$ if $\varepsilon_0 \le c_0^{-1}$ and $m \ge c_0$.

These bounds for $p_{Y,\tau}$ and $q_{Y,\tau}$ have the following implication with regards to the times $t_{Y,\tau+}$ and $t_{Y,\tau-}$. To say more, let $t_{Y,\tau*}$ denote either $t_{Y,\tau+}$ or $-t_{Y,\tau-}$. Then

$$|t_{Y,\tau*} - \lambda^{-1}\ln(\varepsilon_0 Y^{-1})| \le c_0\varepsilon_0 + c_\varepsilon \Delta m^{-6} c_v r^{-1/2}.$$

(C.13)

Given that $\Delta m^{-6} c_v r^{-1/2} \le c_0 m^{-4}$, the right hand side of (C.12) is at most $c_0\varepsilon_0 + c_\varepsilon m^{-4}$.

<u>Step</u> <u>8</u>: Suppose next that $(Y, \tau)$ and $(Y', \tau')$ are as described in Step 7. Introduce $P$ to denote $p_{Y,\tau} - p_{Y',\tau'}$ and $Q$ to denote $q_{Y,\tau} - q_{Y',\tau}$. Subtract their respective versions of (C.7) to derive equations for $P$ and $Q$ that can be written as

$$\tfrac{d}{dt}P = \lambda P + \mathfrak{z}_{pp}P + \mathfrak{z}_{pq}Q + (\tau - \tau')\mathfrak{r}_{p\upsilon} \quad \textit{and} \quad \tfrac{d}{dt}Q = -\lambda Q + \mathfrak{z}_{qp}P + \mathfrak{z}_{qq}Q + (\tau - \tau')\mathfrak{r}_{q\upsilon} ,$$

(C.14)

where each $\mathfrak{z}_{\bullet\bullet}$ is a function of t with norm bounded by $c_0(|p_{Y,\tau}| + |p_{Y',\tau'}| + |q_{Y,\tau}| + |q_{Y',\tau'}|)$. These equations can be analyzed using the same tools used in Part 7 to draw the conclusions expressed by following inequalities. The analysis for $t \ge 0$ finds

- $|Q| \le c_0(\varepsilon_0 |P| + |Y - Y'| + m^{-2} c_v r^{-1/2}|\tau - \tau'|)$ .
- $|P - e^{\lambda t}(Y - Y')| \le c_0 e^{\lambda t}((\varepsilon_0 + \Delta|\ln(\varepsilon_0^{-1}\Delta)|)|Y - Y'| + m^{-6} c_v r^{-1/2}|\tau - \tau'|)$.

(C.15)

Meanwhile, the analysis for $t \le 0$ leads to

- $|P| \le c_0(\varepsilon_0 Q|_t + |Y - Y'| + m^{-2} c_v r^{-1/2}|\tau - \tau'|)$
- $|Q - e^{-\lambda t}(Y - Y')| \le c_0 e^{-\lambda t}((\varepsilon_0 + \Delta|\ln(\varepsilon_0^{-1}\Delta)|)|Y - Y'| + m^{-6} c_v r^{-1/2}|\tau - \tau'|)$.

(C.16)

The bounds in (C.10)-(C.13) and (C.15), (C.16) play central roles in what follows.

Step 9: The bounds in (C.10)-(C.13) and (C.15) and (C.16) can be used to say something about $t_{Y,\tau+} - t_{Y',\tau'+}$ and $t_{Y,\tau-} - t_{Y',\tau'-}$. To this end, suppose for the sake of argument that $t_{Y,\tau+} \ge t_{Y',\tau'+}$. Write

$$p_{Y',\tau'}^2 + q_{Y',\tau'}^2 = p_{Y,\tau}^2 + q_{Y,\tau}^2 - 2(P p_{Y,\tau} + Q q_{Y,\tau}) + P^2 + Q^2 ,$$

(C.17)

and set $t = t_{Y,\tau}$. Use the fact that $(p_{Y,\tau}^2 + q_{Y,\tau}^2)^{1/2} = \varepsilon_0$ at $t = t_{Y,\tau+}$ with (C.8), (C.10) and (C.15) to see that $p_{Y',\tau'}^2 + q_{Y',\tau'}^2$ at $t = t_{Y,\tau+}$ is

$$p_{Y',\tau'}^2 + q_{Y',\tau'}^2 = \varepsilon_0^2(1 - 2Y^{-1}(Y - Y') + \mathfrak{e})$$

(C.18)

where $|\mathfrak{e}| \le c_0(\varepsilon_0 Y^{-1}|Y - Y'| + Y^{-2}|Y - Y'|^2 + m^{-4}|\tau - \tau'|)$. This last inequality with (C.10) and (C.15) allows $t_{Y',\tau'+}$ to be written as

$$t_{Y',\tau'+} = t_{Y,\tau+} + \lambda^{-1} Y^{-1}(Y - Y') + \mathfrak{e} ,$$

(C.19)

where $\mathfrak{e}$ has the same absolute value bound as its namesake in (C.18). The same sort of arguments write $t_{Y',\tau'-}$ as $t_{Y',\tau'-} = t_{Y,\tau-} - \lambda^{-1} Y^{-1}(Y - Y') + \mathfrak{e}$ with $\mathfrak{e}$ being different then its namesakes in (C.18) and (C.19) but obeying the same absolute value bound.

Step 10: The functions p and q are convenient to use on the radius $\varepsilon$ tubular neighborhood of $\gamma$, but less so elsewhere on $\mathcal{H}_p$ and in particular, less so near the boundary of $\mathcal{H}_p$. The function $\hbar = f(u)\cos\theta\sin^2\theta$ is far more convenient, this in part because the final arguments for Lemma C.1's proof are much the same as those used in the proof of Proposition II.2.7. At distance $\varepsilon$ or less from $\gamma$, the function $\hbar$ can be readily written in terms of p and q; and doing so leads to the following formula:

$$\hbar = \tfrac{2}{3\sqrt{3}}(x_0 + 4e^{-2R}) + \tfrac{2}{\sqrt{3}}(x_0 + 4e^{-2R})pq + \mathfrak{h} ,$$

(C.20)

where $|\mathfrak{h}| \le c_0(p^2 + q^2)^{3/2}$ and $|\frac{\partial}{\partial p}\hbar| + |\frac{\partial}{\partial q}\hbar| \le c_0(p^2 + q^2)$.

Fix $\tau \in [0,1]$ and fix Y as in Step 7 so as to define the interval $I_{Y,\tau}$ and the corresponding pair of functions $p_{Y,\tau}$ and $q_{Y,\tau}$ on $I_{Y,\tau}$. Of interest here is the function on $I_{Y,\tau}$ given by the rule $t \to \hbar(p_{Y,\tau}(t), q_{Y,\tau}(t))$. This function is denoted in what follows by $\hbar_{Y,\tau}$.

particular interest are the values $\hbar_{Y,\tau}$ at the $t = t_{Y,\tau+}$ and at $t = t_{Y,\tau-}$. In particular, (C.8) and (C.11) with (C.20) imply that its values at these times differ from $\frac{2}{3\sqrt{3}}(x_0+4e^{-2R})$ by at most $c_0\,\varepsilon_0^3$.

Consider now the functions $\hbar_{Y,\tau}$ and $\hbar_{Y',\tau'}$ with $\tau, \tau' \in [0, 1]$ and with $Y$ and $Y'$ as in Step 8. Of interest is $\hbar_{Y,\tau}(t) - \hbar_{Y',\tau'}(t')$ with $t$ and $t'$ being $t_{Y,\tau+}$ and $t_{Y',\tau'+}$ or else $t$ and $t'$ being $t_{Y,\tau-}$ and $t_{Y',\tau'-}$. Use $t_*$ and $t_*'$ to denote either of these pair of values for t. The absolute value $\hbar_{Y,\tau}(t_*) - \hbar_{Y',\tau'}(t_*')$ obeys the a priori bound

$$|\hbar_{Y,\tau}(t_*) - \hbar_{Y',\tau'}(t_*')| \le c_0\varepsilon_0^3\,(Y^{-1}|Y - Y'| + m^{-4}|\tau - \tau'|)\ ; \tag{C.21}$$

this follows from (C.20) with (C.8), (C.11)-(C.13), (C.15), (C.16) and Step 9's assertions.

Step 11: The arguments that follow in this step and Steps 12 and 13 focus almost entirely on the case when $\upsilon$ intersects but one $\mathfrak{p} \in \Lambda$ version of $\mathcal{H}_{\mathfrak{p}}$. The arguments in the general case are only outlined as they differ from those used for this simplest case in a straightforward fashion; and, in any event, they are much the same as those used for Proposition II.2.7. This step is a guide of sorts for Step 12.

Assume now that $\upsilon$ crosses only one handle from the set $\cup_{\mathfrak{p}\in\Lambda}\mathcal{H}_{\mathfrak{p}}$. Let $\mathfrak{p} \in \Lambda$ denote the relevant pair. With an r-independent $\varepsilon > \varepsilon_0 > 0$ fixed in advance, it is sufficient given what is said in Step 3 to consider only the case where $\upsilon$ intersects the radius $\frac{1}{8}\varepsilon_0$ tubular neighborhood of either $\hat{\gamma}_{\mathfrak{p}}^{+}$ or $\hat{\gamma}_{\mathfrak{p}}^{-}$. The discussion that follows considers only the case of $\hat{\gamma}_{\mathfrak{p}}^{+}$ as the arguments for the other case are identical but for some sign changes and notation. This understood, the notation from Steps 4-9 will be used when referring to the radius $\varepsilon$ tubular neighborhood of this curve. In particular, the curve $\hat{\gamma}_{\mathfrak{p}}^{+}$ is denoted below as $\gamma$. The constant $\varepsilon_0$ is chosen so that the locus in the radius $\varepsilon$ tubular neighborhood where the coordinates p and q obey $(p^2+q^2)^{1/2} \le \varepsilon_0$ lies well inside this tubular neighborhood. The portion of $\upsilon$ in the radius $\varepsilon$ tubular neighborhood of $\gamma$ where $(p^2+q^2)^{1/2} \le \varepsilon_0$ is parametrized by the interval $I \subset \mathbb{R}$ as described in Step 6.

Fix $\tau \in [0, 1]$. The next step constructs a 2-parameter family of continuous and piecewise smooth arcs in $M_\delta \cup \mathcal{H}_{\mathfrak{p}}$ that all start and end on the $f = \frac{3}{2}$ Heegaard surface $\Sigma$ in $M_\delta$. The starting and ending points are both very near $\upsilon$'s intersection with this surface. Any given member of this family starts near $\upsilon$'s intersection point with $\Sigma$ follows $\upsilon$ to the $u = -R - \ln\delta$ sphere in $\mathcal{H}_{\mathfrak{p}}$. The arc stays close to $\upsilon$ through $\mathcal{H}_{\mathfrak{p}}$ so as to exit $\mathcal{H}_{\mathfrak{p}}$ through its $u = R + \ln\delta$ sphere in $\mathcal{H}_{\mathfrak{p}}$. It then follows $\upsilon$ in $M_\delta$ so as to end on the surface $\Sigma$. Each arc from the family is the end-to-end concatenation of 5 smooth segments. The parameter space for the family of arcs is $[(1 - R)Y_0, (1 + R)Y_0] \times [-R, R]$ with $R \le \frac{1}{4}$ to be determined ultimately by $\varepsilon_0$.

Step 12: Fix $\tau \in [0,1]$ and a pair $D = (Y,\sigma) \in [(1-R)Y_0,(1+R)Y_0]\times[-R,R]$. The corresponding continous and piecewise smooth arc in $M_\delta \cup \mathcal{H}_p$ is denoted by $\upsilon_{D,\tau}$. As noted in the previous step, this arc is the end-to-end concatenation of five segments. It proves useful in this regard to describe the middle segment first, then the second and fourth segments and at the end, the first and fifth segments. By way of notation, $c_\varepsilon$ is used in what follows to denote a purely $\varepsilon_0$ dependent constant that is greater than 1. This constant can be assumed to increase between consecutive appearances.

THE MIDDLE SEGMENT: The middle segment crosses the $(p^2+q^2)^{1/2} \le \varepsilon_0$ portion of the radius $\varepsilon$ tubular neighborhood of $\gamma$. This segment is parametrized as the map from the interval $I_{Y,\tau}$ given by the rule $t \to (\phi = -t+\sigma, p = p_{Y,\tau}(t), q = q_{Y,\tau}(t))$ with $I_{Y,\tau}$ and the functions $p_{Y,\tau}$ and $q_{Y,\tau}$ as defined in Step 7. Use what is said in Steps 8 and 9 to see that the $t_{Y,\tau+}$ and $t_{Y\tau-}$ endpoints of this segment on the $(p^2+q^2)^{1/2} = \varepsilon_0$ surface have distance at most $c_0\,\varepsilon_0\,(R+m^{-4})$ from the $u > 0$ and $u < 0$ points where $\upsilon$ intersects this surface. Use $z_{D,\tau-}$ and $z_{D,\tau+}$ to denote these respective endpoints.

THE SECOND AND FOURTH SEGMENTS: Let $z_-$ and $z_+$ denote the respective $u < 0$ and $u > 0$ points where $\upsilon$ intersects the $(p^2+q^2)^{1/2} = \varepsilon_0$ surface in the radius $\varepsilon$ tubular neighorhood of $\gamma$. Introduce as in Steps 4 and 5 the segments $\upsilon_-$ and $\upsilon_+$ of $\upsilon$. By way of a reminder, $\upsilon_-$ starts on the $e^{-2(R+u)}(1-3\cos^2\theta) = \delta^2$ surface in $\mathcal{H}_p$ and ends at $z_-$; and $\upsilon_+$ starts at $z_+$ and ends on the $e^{-2(R-u)}(1-3\cos^2\theta) = \delta^2$ surface in $\mathcal{H}_p$.

Steps 4 and 5 introduce the segments of integral curves of $v$, these being $\gamma_-$ and $\gamma_+$. The former has $u < 0$, contains $z_-$ and starts on the $e^{-2(R+u)}(1-3\cos^2\theta) = \delta^2$ surface in $\mathcal{H}_p$ and the latter constains $z_+$ and ends on the $e^{-2(R-u)}(1-3\cos^2\theta) = \delta^2$ surface. Steps 4 and 5 describe parametrizations of respective tubular neighborhoods of $\gamma_-$ and $\gamma_+$ using coordinates $(t, z)$ with $|z| \le c_\varepsilon^{-1}$ and with the $z = 0$ locus being $\gamma_-$ or $\gamma_+$ as the case may be. The point $(0, 0)$ is $z_-$ in the former case and $z_+$ in the latter. Reintroduce from the final paragraphs of Steps 4 and 5 the constant $c_{p\varepsilon}$. Take $R$ so that $R < (c_0 c_{p\varepsilon})^{-1}$ and take $m$ so that $m > c_0 c_{p\varepsilon}$. Granted these bounds, then the $t = t_{Y,\tau-}$ endpoint $z_{D,\tau-}$ of the middle segment has distance less than $c_{p\varepsilon}^{-1}$ from $z_-$ and the $t = t_{Y,\tau+}$ endpoint $z_{D,\tau+}$ of the middle segment has distance less than $c_{p\varepsilon}^{-1}$ from $z_+$. Let $\Delta_{D,\tau-}$ and $\Delta_{D,\tau+}$ denote these distances.

Step 4 finds a solution to the $(\gamma_-,\tau)$ version of (C.3) that is defined for all $t \in \gamma_-$, contains $z_{Y,\tau-}$ and has pointwise norm bounded by $c_\varepsilon(\Delta_{D,\tau-}+r^{-1/2})$ This solution defines a smoothly embedded arc in the $u < 0$ part of $\mathcal{H}_p$ that starts on the $e^{-2(R+u)}(1-3\cos^2\theta) = \delta^2$ surface. Use $z_{D,\tau-}$ to denote the segment of this arc that starts on this surface and ends at $z_{D,\tau-}$. This arc is the second segment.

Step 5 finds a solution to the $(\gamma_+,\tau)$ version of (C.3) that is defined for each $t \in \gamma_+$, contains $z_{Y,\tau+}$ and has pointwise norm bounded by $c_\varepsilon(\Delta_{D,\tau+}+r^{-1/2})$. This solution defines a smoothly embedded arc in the $u > 0$ part of $\mathcal{H}_p$ that ends on the $e^{-2(R-u)}(1-3\cos^2\theta) = \delta^2$ surface. Use $z_{D,\tau+}$ in what follows to denote the segment of this arc that starts at $z_{D,\tau+}$ and ends on the $e^{-2(R-u)}(1-3\cos^2\theta) = \delta^2$ surface. This arc is the fourth segment.

THE FIRST AND FIFTH SEGMENTS: The starting point of $z_{D,\tau-}$ on the surface in $\mathcal{H}_p$ where $e^{-2(R+u)}(1-3\cos^2\theta) = \delta^2$ has distance at most $c_\varepsilon(\Delta_{D,\tau-} + r^{-1/2})$ from υ's intersection with this surface, the latter denoted in Step 1 by $z_{p-}$. This being the case, it has distance at most $c_\varepsilon(\Delta_{D,\tau-} + r^{-1/2})$ from the $(t = 2-\delta^2, z = 0)$ point in Step 1's coordinate cylinder. Meanwhile, the ending point of $z_{D,\tau+}$ on the $e^{-2(R-u)}(1-3\cos^2\theta) = \delta^2$ surface in $\mathcal{H}_p$ has distance at most $c_\varepsilon(\Delta_{D,\tau+} + r^{-1/2})$ from υ's intersection point with this same surface, and so this ending point of $z_{D,\tau+}$ has distance at most $c_\varepsilon(\Delta_{D,\tau-} + r^{-1/2})$ from $(t = 1+\delta^2, z = 0)$ point in Step 1's coordinate cylinder.

With the preceding understood, suppose that $R < (c_\varepsilon c_{p\varepsilon})^{-1}$ and $m > c_\varepsilon c_{p\varepsilon}$. Granted these bounds, then the starting point of $z_{D,\tau-}$ will be well inside the $t = 2-\delta^2$ boundary disk of Step 1's coordinate cylinder. Use $z_{D,\tau-}$ to denote the z-coordinate of this starting point of $z_{D,\tau-}$. Meanwhile, the ending point of $z_{D,\tau+}$ will be well inside the $t = 1+\delta^2$ boundary disk of Step 1's coordinate cylinder centered on $\gamma_*$. Use $z_{D,\tau+}$ to denote the z-coordinate of this ending point of $z_{D,\tau+}$.

According to Step 1, there is a solution $t \to z(t)$ to (C.2) with $z(2-\delta^2) = z_{D,\tau-}$. Denote this solution by $z^M{}_{D,\tau-}$. The first segment is the $t \in [\frac{3}{2}, 2-\delta^2]$ part of the arc in $M_\delta$ given by the graph $t \to (t, z^M{}_{D,\tau-}(t))$.

There is also a solution $t \to z(t)$ to (C.2) with $z(1+\delta^2) = z_{D,\tau+}$. Denote the latter solution by $z^M{}_{D,\tau+}$. The fifth segment is the $t \in [1+\delta^2, \frac{3}{2}]$ part of the arc in $M_\delta$ given by the graph $t \to (t, z^M{}_{D,\tau+})$.

<u>Step</u> <u>13</u>: Introduce $\phi_+$ and $\phi_-$ to denote the respective $\phi$ coordinates of the $(t=1+\delta^2, z = 0)$ and $(t=2-\delta^2, z = 0)$ points on the two boundary disks of Step 1's coordinate cylinder. The $\phi$ coordinate of the ending point of $z_{D,\tau+}$ on the surface in $\mathcal{H}_p$ where $e^{-2(R-u)}(1-3\cos^2\theta) = \delta^2$ can be written as $\phi_+ + \varphi_{D,\tau+}$ with $|\varphi_{D,\tau+}| \le c_0^{-1}$; and the $\phi$ coordinate of the starting point of $z_{D,\tau-}$ on the $e^{-2(R+u)}(1-3\cos^2\theta) = \delta^2$ surface can be written as $\phi_- + \varphi_{D,\tau-}$ with $|\varphi_{D,\tau-}| \le c_0^{-1}$. Write the respective values of the function $\hbar$ at these boundary points of $z_{D,\tau+}$ and $z_{D,\tau-}$ as $\hbar_{D,\tau+}$ and $\hbar_{D,\tau-}$.

It follows from what is said in Step 1 that the five segment concatenated arc define in Step 12 is a piecewise embedded loop in $M_\delta \cup \mathcal{H}_p$ if

$$\varphi_{D,\tau+} = (-1)^{\hat{o}}(\hbar_{D,\tau-} + (1-\tau)u_1) \quad \textit{and} \quad \hbar_{D,\tau+} = (-1)^{\hat{o}}(\varphi_{D,\tau-} + (1-\tau)u_2) \tag{C.22}$$

where the notation is such that $\hat{o} \in \{0,1\}$ is determined by point $z_* \in C_{p+} \cap C_{p'}$; and where $u_1$ and $u_2$ are the respective real and imaginary parts of

$$-2i \left( \int_{1+\delta^2}^{3/2} x_\upsilon - \int_{3/2}^{2-\delta^2} x_\upsilon \right) \tag{C.23}$$

with $x_\upsilon$ being the function in (C.2). As is explained next, next, there exists a smooth map, $D(\cdot)$: $[0,1] \to ((1-R)Y_0, (1+R)Y_0) \times (-R,R)$ such that for each $\tau \in [0, 1]$, the $D = D(\tau)$ version of (C.22) holds when the following conditions are met:

- $\varepsilon_0 \le c_0^{-1}$.
- $m \ge c_\varepsilon$ *with* $c_\varepsilon > 1$ *being a purely* $\varepsilon$ *dependent constant.*
- $R \le c_\varepsilon m^{-1}$ *with* $c_\varepsilon > 1$ *being a purely* $\varepsilon$ *dependent constant.*
- $c_v \ge c_{\varepsilon,m,R}$ *with* $c_{\varepsilon,m,R} > 1$ *being a constant that depend only on* $\varepsilon$, $m$, *and* $R$.
- $r \ge \kappa_c$ *with* $\kappa_c > 1$ *being a constant that depends only on* $c_v$.

(C.24)

To construct $D(\cdot)$, take $D(0)$ to be the pair $(Y = Y_0, \sigma = 0)$ which obeys (C.22) because the arc defined in Step 12 from $(Y = Y_0, \sigma = 0)$ is the smooth, embedded circle. The construction of $D(\tau)$ for $\tau > 0$ requires a rewriting of (C.22). To set the stage for this fix $D = (Y,\sigma) \in [(1-R)Y_0, (1+R)Y_0] \times [-R,R]$ and $\tau \in [0,1]$. Use (C.19) to write the difference between $\varphi_{D,\tau+}$ and its $(Y_0, \sigma = 0)$, $\tau = 0$ analog as

$$\sigma - \lambda^{-1} Y_0^{-1} (Y - Y_0) + \mathfrak{e} \,, \tag{C.25}$$

where $|\mathfrak{e}| \le c_0((\varepsilon_0 + R) Y_0^{-1} |Y - Y_0| + m^{-4})$. This last formula and its $t_{Y,\tau-}$ analog allow (C.22) to be rewritten as

- $\sigma - \lambda^{-1} Y_0^{-1} (Y - Y_0) - (-1)^{\hat{o}} (\hbar_{Y,\tau}(t_{Y,\tau-}) - \hbar_{Y_0,0}(t_{Y_0,0-})) + \mathfrak{e}_1 = 0$ ,
- $\sigma + \lambda^{-1} Y_0^{-1} (Y - Y_0) - (-1)^{\hat{o}} (\hbar_{Y,\tau}(t_{Y,\tau+}) - \hbar_{Y_0,0}(t_{Y_0,0+})) + \mathfrak{e}_2 = 0$ ,

(C.26)

where $\mathfrak{e}_1$ and $\mathfrak{e}_2$ are functions of $\tau$, $Y$ and $\sigma$ with whose absolute values are bounded by $c_0((\varepsilon_0 + R) Y_0^{-1} |Y - Y_0| + m^{-4})$. The left hand side of (C.26) defines a smooth map, $\mathcal{F}$, from $[0,1] \times [(1-R)Y_0, (1+R)Y_0] \times [-R,R]$ to $\mathbb{R}^2$ with the property that $\mathcal{F} = 0$ if and only if (C.22) is obeyed. What follows is a crucial observation about this map: The differential of $\mathcal{F}$ along the domain's factor $[(1-R)Y_0, (1+R)Y_0] \times [-R,R]$ is surjective if (C.24) holds, this being a consequence of (C.21) and what is said about $\mathfrak{e}$ in (C.25).

Suppose that $\tau_0 \in [0, 1]$ is such that $D(\cdot)$ has been defined on $[0, \tau_0)$. To extend $D(\cdot)$ to a larger interval, use the fact that $[(1-R)Y_0,(1+R)Y_0]\times[-R,R]$ is compact to see that there is a $\tau = \tau_0$ limit point $D(\tau_0)$ for $\{D(\tau)\}_{\tau\to\tau_0}$. It follows from (C.26) that this limit point is unique, that the extension of $D(\cdot)$ to $[0, \tau_0]$ is continuous, and that $D(\tau_0)$ obeys the $\tau = \tau_0$ version of (C.26). Write $D(\tau_0)$ as $(Y,\sigma)$. Use (C.19), (C.21) and the fact that (C.26) is obeyed to conclude that $|\sigma| + |Y - Y_0| \le c_0\, m^{-4}$. This implies that $D(\tau_0)$ lies in the interior of the parameter space $[(1-R)Y_0,(1+R)Y_0]\times[-R,R]$ when (C.24) holds. Since $D(\tau_0)$ is not a boundary point of $[(1-R)Y_0,(1+R)Y_0]\times[-R,R]$, the fact that $\mathcal{F}$ has surjective differential along the $[(1-R)Y_0,(1+R)Y_0]\times[-R,R]$ factor of its domain implies via the inverse function theorem that $D(\cdot)$ has a smooth extension to an open interval in $[0,1]$ that contains $[0, \tau_0]$.

Step 14: Fix $\tau \in [0,1]$. Let $\upsilon^{\alpha}{}_{\tau*}$ denote the continuous, piecewise smooth loop that is defined in Step 12 by the data set $D(\tau)$ from the previous step. The loop $\upsilon^{\alpha}{}_{\tau*}$ is smooth on the interior of each of its five concatenating segments. The implicit function theorem construction implies that the assignment of $\tau \in [0,1]$ to each of the five concatentating segments defines a smoothly varying arc in $M_\delta \cup \mathcal{H}_p$. Moreover, the assertions made by the four bullets in Lemma C.1 holds for each of these five $[0,1]$-parametrized families of arcs. This follows directly from the implicit function theorem construction given what is said at the very end of Steps 4 and 5; and given (C.15), (C.16), the bound for $\mathfrak{e}$ in (C.19) and the bound in (C.21).

The loop $\upsilon^{\alpha}{}_{\tau*}$ for $\tau \in (0,1)$ is continuous, but its derivative may be discontinuous at four points, these being the loop's intersection points with the boundary spheres of the $|u| \le R+\ln\delta$ part of $\mathcal{H}_p$ and the $(p^2+q^2)^{1/2} = \varepsilon_0$ surface in the radius $\varepsilon$ tubular neighborhood of $\gamma$. Even so, the two concatenating segments near these points are smooth up to their endpoints on the relevant surface, and the corresponding tangent vectors differ by at most $c_0 r^{-1/2}$ at these endpoints. This is because the tangent vector to each segment near these junctions differs from $\nu$ by at most $c_0 r^{-1/2}$. The preceding fact implies that the loop $\upsilon^{\alpha}{}_{\tau*}$ can be smoothed near the junctions of segments so that the result is a smoothly embedded loop that obeys the first three bullets of Lemma C.1. Moreover, it s a straightforward task to define this smoothing without changing the already smooth $\tau = 0$ and $\tau = 1$ versions so that the resulting [0, 1]-parameter family of smooth loops is smoothly parametrized and obeys the assertion of Lemma C.1's fourth bullet at each point. The details of this are straightforward and thus omitted.

Use $\{\upsilon^{\alpha}{}_{\tau}\}_{\tau\in[0,1]}$ to denote the resulting $[0,1]$ parametrized family of smooth loops. This family obeys all of the requirements for Lemma C.1.

### c) Increasing r

This part takes Proposition C.2's pair $(A_{\Diamond 11},\psi_{\Diamond 11})$ as the starting point of a path in $\mathrm{Conn}(E) \times C^\infty(Y;\mathbb{S})$ whose end member is constructed from the same vortex solutions and loops in Y that are used to construct $(A_{\Diamond 11}, \psi_{\Diamond 11})$ but with the given choice of r replaced by a far larger choice. A result from [T4] is brought to bear in the next section; it requires the larger value of r. This larger value of r is denoted by R. There is no upper bound to the value chosen but a lower bound $R > c_0 c_v^4 r$ is imposed.

The path is parametrized by [0,1], and a given $\tau \in [0,1]$ member denoted by $(A_{\bullet\tau}, \psi_{\bullet\tau})$. The definition of $(A_{\bullet\tau}, \psi_{\bullet\tau})$ is identical to that of $(A_{\Diamond 11}, \psi_{\Diamond 11})$ given in Section Ca but for the replacement of r with $r(\tau) = (1-\tau)r + \tau R$.

Keep in mind in what follows that the zero locus of the E summand of any $\tau \in [0,1]$ version of $\psi_{\bullet\tau}$ is identical to that of $\alpha_{\Diamond 11}$ as is the degree of vanishing along any transverse disk centered on the zero locus. By way of a reminder, the zero locus of $\alpha_{\Diamond 11}$ consists solely of closed integral curves of $v$ in $M_\delta \cup (\cup_{p \in \Lambda} \mathcal{H}_p)$, these coming from two sets. The first set consisted of curves that intersect $M_\delta$; this set was denoted by $\Theta^\alpha$. The second set is a subset of $\cup_{p \in \Lambda}\{\hat{\gamma}_p^+ \cup \hat{\gamma}_p^-\}$.

Of interest is the spectral flow between the $\mathfrak{c} = (A_{\bullet 0}, \psi_{*0}) = (A_{\Diamond 11}, \psi_{\Diamond 11})$ version of $\mathfrak{L}_{\mathfrak{c},r}$ and the $\mathfrak{c} = (A_{\bullet 1}, \psi_{\bullet 1})$ version of $\mathfrak{L}_{\mathfrak{c},R}$. Note in particular that the latter operator is defined using R rather than r. The proposition that follows asserts that there is an a priori upper bound for the norm of the spectral flow between these two operators that is independent of the original pair $(A, \psi)$ and r and also R. This proposition uses $\mathfrak{c}(\tau)$ to denote the pair $(A_{\bullet\tau}, \psi_{\bullet\tau})$.

**Proposition C.4**: *There exists* $\kappa \geq 100$, *and given* $c_v \geq \kappa$, *there exists* $\kappa_{c_v} > \kappa$ *with the following significance: Suppose that* $r \geq \kappa_{c_v} c_v^{10}$ *and suppose that* $(A, \psi = (\alpha, \beta))$ *is a solution to the* $(r, \mu)$ *version of (1.13) with* $\mu$ *a given element in* $\Omega$ *with* $\mathcal{P}$*-norm smaller than 1. The values of* $\kappa$, $c_v$ *and* r *are suitable for defining* $(A_{\Diamond 11}, \psi_{\Diamond 11})$ *and any* $R \geq \kappa c_v^6 r$ *version of the family* $\{(A_{\bullet\tau}, \psi_{\bullet\tau})\}_{\tau \in [0,1]}$. *The norm of the spectral flow between the end members of the corresponding family* $\{\mathfrak{L}_{\mathfrak{c}(\tau), r(\tau)}\}_{\tau \in [0,1]}$ *is bounded by* $\kappa$.

***Proof of Proposition C.4***: It is assumed in what follows that $\kappa$, $c_v$ and r are large enough to invoke the various results in the preceding subsections of Appendix C and those in Appendix A and B. The proof that follows has five parts.

*Part 1*: As explained directly, each member of the family $\{(A_{\bullet\tau}, \psi_{\bullet\tau})\}_{\tau \in [0,1]}$ obeys a version of PROPERTIES 1-5 in Section Ab. The proof that such is the case distinguishes between values of $\tau$ near 0 and larger values. To elaborate, note first that the various properties require the specification of constants $c_0$ and $z$. It is a straightforward matter to check that PROPERTIES 1, 2, 4 and 5 are obeyed using $c_v$ in lieu of $c_0$ and $r(\tau)$ in lieu of $z$. It is also a straight forward matter to verify the third and fourth bullets of PROPERTY 3.

Meanwhile, Items a) and b) of the second bullet of PROPERTY 3 follow directly by virtue of the fact that the zero locus of the E summand of $\psi_{\bullet\tau}$ is a union of integral curves of $v$. The story with regards to the first bullet and Items c) and d) of the second bullet of PROPERTY 3 is not so straightforward. The point being that $Y - Y_{\Diamond z}$ is the union of tubular neighborhoods of the curves from the set $\cup_{p\in\Lambda}\{\hat{\gamma}_p^+ \cup \hat{\gamma}_p^-\}$ with radius proportional to $z^{-1/2}$. If a component of the zero locus of $\alpha_{\Diamond 11}$ from $\Theta^\alpha$ interect the $z = r$ version of $Y - Y_{\Diamond z}$, then the first bullet and Items c) and d) of the second bullet of PROPERTY 3 will not hold when $z = r(\tau)$ for values of $\tau$ in certain subsets of $[0,1]$.

To deal with this issue, fix $\tau \in [0,1]$. The $c_0 = c_v$ and $z = r(\tau)$ version of PROPERTY 3 can fail if there exists a curve from $\Theta^\alpha$ and a curve from $\cup_{p\in\Lambda}\{\hat{\gamma}_p^+ \cup \hat{\gamma}_p^-\}$ with the following property: Let $\upsilon^\alpha$ denote the curve from $\Theta^\alpha$ and let $\gamma$ denote the curve from $\cup_{p\in\Lambda}\{\hat{\gamma}_p^+ \cup \hat{\gamma}_p^-\}$. Then the minimum distance between the points in $\upsilon^\alpha$ and $\gamma$ is no less than $(c_v^4 - 3c_v^3)\, r(\tau)^{-1/2}$ and no greater than $(c_v^4 + 3c_v^3)\, r(\tau)^{-1/2}$. This last observation has two immediate consequences, the first being that PROPERTY 3 can fail only in the case when $r(\tau) \le c_0 c_v^6 r$ and thus only if $\tau \le c_0 c_v^6\, r/R$. This is so because the minimum distance between $\upsilon^\alpha$ and $\gamma$ is in any event greater than $c_0^{-1} c_v r^{-1/2}$. To state the second consequence, introduce $c_1 = c_v - 2c_v^{-1}$. If the distance between $\upsilon^\alpha$ and $\gamma$ is no less than $(c_v^4 - 3c_v^3)\, r(\tau)^{-1/2}$, then it is greater than $(c_1^4 + 3c_1^3)\, r(\tau)^{1/2}$ if $c_v \ge c_0$.

Given what was just said, it is a straightforward task to use the pointwise bounds give in Appendix Aa for the absolute values of $\alpha_0$, $a_0$, y and $\varsigma$ to verify the following assertion: If $c_v \ge c_0$, then $(A_{\bullet\tau}, \psi_{\bullet\tau})$ obeys the $c_0 = c_v - c_0 c_v^{-1}$ and $z = r(\tau)$ version of PROPERTIES 1-5 if it does not obey the $c_0 = c_v$ version.

*Part 2*: A suitable bound for the absolute value of the spectral flow is obtained by studying the variation with $\tau$ of the spectrum of the family of operators $\{\mathfrak{L}_{c(\tau),r(\tau)}\}_{\tau\in[0,1]}$. This part of the subsection considers the values of $\tau$ when some closed integral curve of $v$ from $\Theta^\alpha$ and some curve from $\cup_{p\in\Lambda}\{\hat{\gamma}_p^+ \cup \hat{\gamma}_p^-\}$ have minimum distance at most $(c_0^4 + 3c_0^3)\, r(\tau)^{-1/2}$ where it is understood that $c_0 \in (c_v - c_0 c_v^{-1}, c_v]$ and that $(A_{\bullet\tau}, \psi_{\bullet\tau})$ obeys the $c_0$ and $z = r(\tau)$ version of PROPERTIES 1-5 in Section Ab. As noted in Part 1, this condition can hold only if $\tau \le c_0 c_v^6\, r/R$ and so $r(\tau) \le c_0 c_v^6 r$. This understood, suppose in what follows that this minimum distance condition holds for $\tau \le c_0 c_v^6\, r/R$ and that this minimum distance condition does not hold for $\tau \ge c_0 c_v^6\, r/R$.

An almost verbatim repetition of the arguments used to prove Proposition C.1 find an $(A, \psi)$ and r independent bound for the absolute value of the spectral flow for the $\tau \le c_0 c_v^6\, r/R$ part of the family $\{\mathfrak{L}_{c(\tau),r(\tau)}\}_{\tau\in[0,1]}$. The bound for the absolute value of the spectral flow does, however, depend on the choice for $c_v$. The only salient changes to the arguments from the proof of Proposition C.1 involve Steps 1 and 2. Step 1 is replaced by Part 1 above. The change to Step 2 adds extra terms to the right hand side of (C.1) to account for the fact that relevant version of $\frac{d}{d\tau}\mathfrak{L}_{\mathbb{V}\tau}$ is non-zero on the radius $c_0^{-1} c_v r^{-1/2}$ tubular neighborhood of certain curves from the set $\cup_{p\in\Lambda}\{\hat{\gamma}_p^+ \cup \hat{\gamma}_p^-\}$. In any event, the absolute value of the homomorphism $\frac{d}{d\tau}\mathfrak{L}_{\mathbb{V}\tau}$ is bounded by $\kappa_{c1} r^{1/2}$ (R/r). Steps 3 and 4 can be repeated with the only change being that the interval $[0,1]$ is replaced by $[0, c_0 c_v^6 r/R]$ and the latter is divided into some $m \le c_0 c_v^6 \kappa_c$ segments of length at most $\kappa_c^{-1} r/R$.

*Part 3*: Assume that $\tau_0 \in [0,1]$ is such that the following is true: Let $\upsilon^\alpha$ denote a component of the zero locus of $\alpha_{\Diamond 11}$ from $\Theta^\alpha$. Then $\upsilon^\alpha$ has distance greater than $(c_0^4 + 3c_0^3) r(\tau)^{1/2}$ from all curves from $\cup_{p\in\Lambda}\{\hat{\gamma}_p^+ \cup \hat{\gamma}_p^-\}$ with $c_0 \in (c_v - c_0 c_v^{-1}, c_v]$. Part 4 of the proof derives an $(A,\psi)$, r and R independent upper bound for the absolute value of the spectral flow between along the $[\tau_0, 1]$ part of the 1-parameter family of operators $\{\mathfrak{L}_{\mathbb{V}\tau}\}_{\tau\in[0,1]}$ under the assumption that $c_v \ge c_0$ and $r \ge \kappa_c$ with $\kappa_c \ge 1$ denoting a purely $c_v$ dependent constant. Such a bound with the bound in Part 3 implies what is asserted by Proposition C.4.

The arguments in Part 4 invoke the following auxilliary lemma.

**Lemma C.5**: *There exists* $\kappa > 1$ *with the following significance: Let* $\upsilon \in Y$ *denote a closed, integral curve of* $v$ *that lies entirely in* $M_\delta \cup (\cup_{p\in\Lambda}\mathcal{H}_p)$. *Fix coordinates from Part 4 of Section Aa for a tubular neighborhood of* $\upsilon$. *The corresponding version of the operator* $\eta \to \frac{i}{2}\frac{d}{dt}\eta + \nu\eta + \mu\bar{\eta}$ *on* $C^\infty(\gamma, \mathbb{C})$ *has no eigenvalue between* $-\kappa^{-1}$ *and* $\kappa^{-1}$.

This lemma is prove in Part 5.

The arguments in Part 4 require a second auxilliary observation, this concerning the spectrum of operators that are associated to components of the zero locus of $\alpha_{\Diamond 1}$ from the set $\cup_{p\in\Lambda}\{\hat{\gamma}_p^+ \cup \hat{\gamma}_p^-\}$. These operators are versions of those depicted in (3.10) with the pair $(\nu, \mu)$ in (3.8) being that from any $\gamma \in \cup_{p\in\Lambda}\{\hat{\gamma}_p^+ \cup \hat{\gamma}_p^-\}$ version of (A.6) with both functions constant, $\mu$ real and greater than $|\nu|$. With a positive integer, m, chosen, the relevant equivalence class from $\mathfrak{C}_m$ is that defined by the solution to (2.8) and (3.1) with $\alpha_0 = |\alpha_0|(\frac{z}{|z|})^m$. This operator is denoted by $\mathcal{L}_m$. What follows is the second observation.

*Given* $m_* \ge 1$ *there exists* $\kappa > 1$ *such that if* $m \le m_*$, *then the operator* $\mathcal{L}_m$ *has no eigenvalues with absolute value in the interval* $(0, \kappa^{-1}]$.

(C.27)

Such $\kappa$ exists because there are only $m_*$ versions of (3.10) involved and each has discrete spectrum with no accumulation points. By way of a parenthetical remark, it is likely that these versions of (3.10) have trivial kernel and so lack eigenvalues in $[-\kappa^{-1}, \kappa^{-1}]$.

*Part 4*: This part assumes Lemma C.5 to complete the proof of Proposition C.4. This is done in the four steps that follow. These steps use $c_c$ to denote a constant that is greater than 1 and depends only on $c_v$. Its value can increase between successive appearances. These steps also use $\kappa_*$ to denote the smaller of the versions of $\kappa$ that appear in Lemma C.5 and in (C.27).

Step 1: Lemmas A.8 and A.9 can be invoked if $c_v \geq c_0$ and $r \geq c_c$ because the integer m that appears in Lemma A.9 is a priori bounded by $c_0$. This understood, what follows is a direct consequence of what is said by Lemmas A.6, A.8 and A.9: If $c_v \geq c_0$ and $r \geq c_v$, then the number of linearly independent eigenvalues of any $\tau \in [0,1]$ version of $\mathfrak{L}_{c(\tau),r(\tau)}$ with eigenvalue between $-\frac{1}{100}\kappa_*^{-1}$ and $\frac{1}{100}\kappa_*^{-1}$ is bounded by $c_0$.

Step 2: As in the proof of Proposition C.2, no generality is lost by assuming that the parametrization of the family $\{\mathfrak{L}_{c(\tau),r(\tau)}\}_{\tau\in[0,1]}$ is real analytic so as to apply what is said in Part 1 of the proof of Proposition B.3. This understood, let $\{\lambda_{n\tau}\}_{n\in\mathbb{Z},\tau\in[0,1]}$ denote the corresponding family of eigenvalues. Let $c_1$ denote the dimension bound given in Step 1. Given what is said in Step 1, the absolute value of the spectral flow for the $[\tau_0, 1]$ part of the family $\{\mathfrak{L}_{c(\tau),r(\tau)}\}_{\tau\in[0,1]}$ is no greater than $c_1$ unless some $\tau \in [\tau_0, 1]$ version of $\mathfrak{L}_{c(\tau),r(\tau)}$ has an eigenvalue between $\frac{1}{100}\kappa_*^{-1}$ and $\frac{1}{50}\kappa_*^{-1}$. Suppose for the sake of argument that such is the case. Let $\mathfrak{f}$ denote the corresponding eigenfunction and $\lambda$ its eigenvalue.

Let $\upsilon^\alpha$ denote a given component of the zero locus of $\alpha_\Diamond$ from $\Theta^\alpha$ and let $\zeta$ denote the section of the $\gamma = \upsilon^\alpha$ version of line bundle $\mathrm{Ker}_\vartheta|_\gamma \to \gamma$ that is described in Lemma A.8. Lemmas A.8 and Lemma C.5 are not mutually compatible if the $L^2$ norm of $\zeta$ is greater than $c_0 c_v^{-1} \|\mathfrak{f}\|_2$.

Step 3: Let $\gamma \in \cup_{p\in\Lambda}\{\hat{\gamma}_p^+ \cup \hat{\gamma}_p^-\}$ denote a component of the zero locus of $\alpha_{\Diamond 1}$ and let $\zeta$ denote the section of the bundle $\mathrm{Ker}_\vartheta|_\gamma \to \gamma$ that is described in Lemma A.9. Use m to denote the integer for $\gamma$'s version of Lemma A.9. Note in particular that $m \leq c_0$. Introduce $\zeta_-$ to denote the $L^2$-orthogonal projection of $\zeta$ onto the span of the eigenvalues $\mathcal{L}_m$ with eigenvalue 0 or less, and use $\zeta_+$ to denote the $L^2$-orthogonal projection of $\zeta$ onto the span of the eigenvalues of $\mathcal{L}_m$ with eigenvalue greater than $\kappa_*^{-1}$. Note that $\zeta = \zeta_- + \zeta_+$. Lemma A.9 and (C.27) are not mutually compatible if the $L^2$ norm of either $\zeta_-$ or $\zeta_+$ is greater than $c_0 c_v^{-1} \|\mathfrak{f}\|_2$.

Step 4: It follows from what is said in Steps 2 and 3 that $\|\Pi_\vartheta \mathfrak{f}\|_2 \le c_0 c_v \|\mathfrak{f}\|_2$. But if $c_v \ge c_0$, then this last conclusion is incompatible with what is said by Lemma A.6 if $\mathfrak{f}$ is not identically zero.

*Part 5*: This last part of the subsection contains the proof of Lemma C.5.

***Proof of Lemma C.5***: Let $L_\upsilon$ denote the operator in question. Proposition II.2.7 of asserts that $\gamma$ is hyperbolic and such is the case if and only if $L_\upsilon$ has trivial kernel. This understood, the only issue is that of the size of the neighborhood of 0 that lacks eigenvalues. The six steps that follow momentarily prove that such a neighborhood contains an interval of the form $(-c_0^{-1}, c_0)$.

Keep in mind when reading the proof that $L_\upsilon$ is defined by the pair $(\nu, \mu)$ and that the latter are defined by the choice of a unitary frame for $K^{-1}|_\upsilon$. This last fact has the following implication: Any two versions of $(\nu, \mu)$ that arise from $\upsilon$'s version of (A.6) give isospectral versions of $L_\upsilon$. This being the case, no generality is lost by choosing the coordinates so that $|\nu| + |\mu| \le c_0$.

Step 1: Fix $\varepsilon > 0$ so that the radius $\varepsilon$ tubular neighborhood of any given curve in the set $\cup_{p\in\Lambda}\{\hat{\gamma}_p^+ \cup \hat{\gamma}_p^-\}$ has coordinates $(\phi, x, y)$ with $x = b^{-1}u$ and $y = \theta_* + \theta$ with b denoting $\frac{2}{3\sqrt{3}} e^R (x_0 + 4e^{-2R})^{1/2}$ and with $\theta_*$ such that $\cos\theta_* = \pm\frac{1}{\sqrt{3}}$ as the case may be. These are the coordinates used in Step 4 and the subsequent steps of the proof of Lemma C.1. Set $p = y + x$ and $q = y - x$. Fix $\varepsilon_0 \in (0, c_0^{-1}\varepsilon)$ so that the surface $(p^2+q^2)^{1/2} = \varepsilon_0$ lies well inside the radius $\varepsilon$ tubular neighborhood.

Suppose that $\upsilon$ enters the $(p^2+q^2)^{1/2} < \frac{1}{2}\varepsilon_0$ part of the radius $\varepsilon$ tubular neighborhood about a given $\gamma \in \cup_{p\in\Lambda}\{\hat{\gamma}_p^+ \cup \hat{\gamma}_p^-\}$. The discussion that follows considers only the case when $\gamma$'s version of $\cos\theta_*$ is equal to $\frac{1}{\sqrt{3}}$ as the discussion for the other case is identical but for some sign changes. Use $\lambda$ to denote $4\sqrt{6}\, e^{-R}(x_0 + 4e^{-2R})^{1/2}$. The part of $\upsilon$ in the $(p^2+q^2) \le \varepsilon_0$ part of the radius $\varepsilon$ tubular neighborhood of $\gamma$ can be written in terms of the coordinates $(\phi, p, q)$ as the image of a map $t \to (\phi = -t, p_\upsilon(t), q_\upsilon(t))$ with the domain being an interval $I_\upsilon \subset \mathbb{R}$ containing the origin and with the pair $(p_\upsilon, q_\upsilon)$ obeying the version of (C.4) with $\mathfrak{r}_{p\upsilon} = \mathfrak{r}_{q\upsilon} = 0$. They also obey analogs of (C.8), (C.11) and (C.12) with no terms proportional to $r^{-1/2}$ and with $x_0 = 0$.

As in the proof of Lemma C.1, no generality is lost by taking the $t = 0$ point so that $p_\upsilon = q_\upsilon$ at $t = 0$; this being the point where $\upsilon$ crosses the $u = 0$ sphere. With this choice understood, write $I_\upsilon$ as $[t_-, t_+]$. The pair $t_-$ and $t_+$ obey (C.6) with $\Delta$ denoting the value of $(p_\upsilon^2 + q_\upsilon^2)^{1/2}$ at $t = 0$, this being the minimal value of $(p^2+q^2)^{1/2}$ on $\gamma$.

Step 2: Fix T > 1 and let $\zeta$ denote an eigenvector of $L_\upsilon$ whose eigenvalue has absolute value no greater than $T^{-1}\lambda$. The eigenvalue equation for $\zeta$ on the part of $\upsilon$ in the radius $\varepsilon$ tubular neighborhood of $\gamma$ where $(p^2+q^2)^{1/2} \le \varepsilon_0$ can be written as an equation for a pair of $\mathbb{R}$ valued functions $t \to (\zeta_1(t), \zeta_2(t))$ on the interval $I_\upsilon$. This equation has the form

$$\tfrac{d}{dt}\zeta_1 = \lambda\zeta_1 + \mathfrak{e}_{11}\zeta_1 + \mathfrak{e}_{12}\zeta_2 \quad \textit{and} \quad \tfrac{d}{dt}\zeta_2 = -\lambda\zeta_2 + \mathfrak{e}_{21}\zeta_1 + \mathfrak{e}_{22}\zeta_2\,, \tag{C.28}$$

where $\mathfrak{e}_{11}$, $\mathfrak{e}_{12}$, $\mathfrak{e}_{21}$ and $\mathfrak{e}_{22}$ are smooth functions on $I_\upsilon$ that are bounded by $c_0(\varepsilon_0+T^{-1}\lambda)$. The fact that $|\nu| + |\mu| \le c_0$ implies that $c_0^{-1}|\zeta|^2 \le |\zeta_1|^2 + |\zeta_2|^2 \le c_0|\zeta|^2$.

Step 3: The same argument that proves (C.8) can be used with (C.28) to prove that $|\zeta_2| \le c_0(\varepsilon_0+T^{-1}\lambda)|\zeta_1|$ where $t \ge 0$; and it can be used to prove that $|\zeta_1| \le c_0(\varepsilon_0+T^{-1}\lambda)|\zeta_2|$ where $t \le 0$. Granted these bounds, multiply the left hand equation by $\zeta_1$ and the right hand by $\zeta_2$. Integrate the resulting equalities to see that $|\zeta_1|^2 + |\zeta_2|^2$ at $t_+$ and $t_-$ are at most $(\lambda + c_0(\varepsilon_0+T^{-1}\lambda))\,\|\zeta\|_2$. It then follows from (C.28) that

$$(|\zeta_1|^2 + |\zeta_2|^2)(t) \le c_0\|\zeta\|_2\,(e^{-\lambda(t_+-t)/c_0} + e^{-\lambda(|t_-|+t)/c_0}) \tag{C.29}$$

at each $t \in I_\upsilon$.

Fix $L \ge 1$ and suppose that both $t_+$ and $|t_-|$ are greater than $c_0\lambda^{-1}2L$. If such is the case, then (C.29) implies that

$$(|\zeta_1|^2 + |\zeta_2|^2)(t) \le c_0\|\zeta\|_2\, e^{-L} \tag{C.30}$$

at times $t \in I$ with distance $L$ or more from $t_-$ and $t_+$.

Step 4: Use $\chi$ to construct a smooth, non-negative function on $\gamma$ that is equal to 1 except at points in $I_\upsilon$ with distance $L$ or less from either $t_-$ or $t_+$. This function should equal 0 at points on $I_\upsilon$ with distance greater than $L+1$ from both $t_-$ and $t_+$; and the absolute value of its should be bounded by 4. Use $\chi_{\gamma,L}$ to denote this function. What follows is a consequence of (C.30):

$$\|\chi_{\gamma,L}\zeta\|_2 \ge (1 - c_0 e^{-L})\|\zeta\|_2 \quad \textit{and} \quad \|L_\upsilon(\chi_{\gamma,L}\zeta)\|_2 \le T^{-1}\lambda\|\chi_{\gamma,L}\zeta\|_2\,. \tag{C.31}$$

The function $\chi_{\gamma,L}$ can be defined for each $\gamma \in \cup_{p\in\Lambda}\{\hat{\gamma}_p^+ \cup \hat{\gamma}_p^-\}$ of the sort under consideration. Multiply $\zeta$ by all such functions and the result is a section of $K^{-1}|_\upsilon$ with compact support on the part of $\upsilon$ with distance greater than $c_0^{-1}\varepsilon_0\, e^{-c_0L}$ from all curves in the set $\cup_{p\in\Lambda}\{\hat{\gamma}_p^+ \cup \hat{\gamma}_p^-\}$. This section is denoted by $\zeta_{\varepsilon,L}$. What is said by (C.31) implies that $\|L_\upsilon \zeta_{\varepsilon,L}\|_2 \le T^{-1}\lambda\|\zeta_{\varepsilon,L}\|_2$ .

<u>Step</u> <u>5</u>: Use $\upsilon_{\varepsilon,L}$ to denote the part of $\upsilon$ with distance $c_0^{-1}\varepsilon_0\, e^{-c_0 L}$ or more from all curves in the set $\cup_{p\in\Lambda}\{\hat{\gamma}_p^+ \cup \hat{\gamma}_p^-\}$. It follows from what is said in Step 2 that $\upsilon_{\varepsilon,L}$ has length at most $c_0(L+|\ln\varepsilon_0|)$. The fact that $L_\upsilon$ is first order, that its coefficients are bounded by $c_0$ and that $\upsilon_{\varepsilon,L}$ has length at most $c_0(L+|\ln\varepsilon_0|)$ has the following consequence: Let $\eta$ denote a section of $K^{-1}|_\upsilon$ with compact support on $\upsilon_{\varepsilon,L}$. Then $\|L_\upsilon\eta\|_2 \ge c_{\varepsilon,L}^{-1}\|\eta\|_2$ with $c_{\varepsilon,L}$ being a constant that is greater than 1 and depends only on $\varepsilon_0$ and L, but not on $\upsilon$.

This last bound on $\|L_\upsilon\eta\|_2$ runs afoul of the inequality $\|L_\upsilon \zeta_{\varepsilon,L}\|_2 \le T^{-1}\lambda\|\zeta_{\varepsilon,L}\|_2$ unless T is less than $c_0 c_{\varepsilon,L}\lambda$.

<u>Step</u> <u>6</u>: Choose $\varepsilon_1 << c_0^{-1}\varepsilon_0\, e^{-c_0 L}$. Suppose that $\upsilon$ is a closed, integral curve of $v$ whose points have distance $\varepsilon_1$ or more from all curves in the set $\cup_{p\in\Lambda}\{\hat{\gamma}_p^+ \cup \hat{\gamma}_p^-\}$. It follows as a consequence that $\upsilon$'s version of $L_\upsilon$ has no eigenvalues between $-c_{1\varepsilon}^{-1}$ and $c_{1\varepsilon}^{-1}$ with $c_{1\varepsilon} \ge 1$ depending only on $\varepsilon_1$. Such a constant exists because $L_\upsilon$ has trivial kernel, and because there is but a finite set of closed orbits of $v$ in $M_\delta \cup (\cup_{p\in\Lambda}\mathcal{H}_p)$ that have distance $\varepsilon_1$ from $\cup_{p\in\Lambda}\{\hat{\gamma}_p^+ \cup \hat{\gamma}_p^-\}$.

The bound in Step 5 and the bound in the preceding paragraph give a $\upsilon$-independent, strictly positive lower bound to the absolute value of any eigenvalue of $L_\upsilon$.

**d) Decreasing r**

A unique set of closed integral curves of $v$ are defined by three properties, the first three following directly. By way of notation, the set in question is denoted here by $\Theta^0$. The first property requires that all curves from $\Theta^0$ lie in $M_\delta \cup (\cup_{p\in\Lambda}\mathcal{H}_p)$ and that none are from $\cup_{p\in\Lambda}\{\hat{\gamma}_p^+ \cup \hat{\gamma}_p^-\}$. The second property requires that the union of the curves from $\Theta^0$ intersects $M_\delta$ as G segments give the same pairing of the index 1 and index 2 critical points of $f$ as that given by the third bullet of Proposition 2.4 using the zero locus of $\alpha$.

The statement of the third property requires introducing notation from Proposition II.2.7. This proposition characterizes a segment of an integral curve of $v$ in a version of $\mathcal{H}_p$ that starts on the $u < 0$ boundary and ends on the $u > 0$ boundary. Proposition II.2.7 characterizes such a segment by an integer, denoted by $\mathfrak{k}_p$. This $\mathfrak{k}_p$ is such that the total change in the $\phi$ angle along the segment in $\mathcal{H}_p$ can be written as $\sigma + 2\pi\mathfrak{k}_p$ with $\sigma \in [0,2\pi)$.

The first two properties imply that the union of the curves from $\Theta^0$ intersect each $\mathfrak{p} \in \Lambda$ version of $\mathcal{H}_p$ as a single segment of the sort just described. This understood, the third property requires that each of the corresponding $\mathfrak{p} \in \Lambda$ versions of $\mathfrak{k}_p$ be 0.

Let $\Theta^1$ denote the subset of pairs of the form $(\gamma, m)$ where $\gamma \in \cup_{p\in\Lambda}\{\hat{\gamma}_p^+ \cup \hat{\gamma}_p^-\}$ is a component of the zero locus of $\alpha_{\lozenge 11}$ and m is the integer that is used to define $(A_{\lozenge 11}, \psi_{\lozenge 11})$ near $\gamma$ via (A.44).

Part 1 of what follows uses the sets $\Theta^0$ and $\Theta^1$ and a real number $z > c_0$ to specify a pair in $\mathrm{Conn}(E) \times C^\infty(Y; \mathbb{S})$. This pair is denoted in what follows by $\mathfrak{c}(z)$. Each such pair has its corresponding operator $\mathfrak{L}_{\mathfrak{c}(z),z}$. Part 2 of this subsection states and then proves two lemmas that supply an a priori upper bound for the absolute value of the spectral flow between any $z = z_0$ and $z = z_1$ version of $\mathfrak{L}_{\mathfrak{c}(z),z}$. Part 3 states and then proves a proposition that compares the absolute value of the spectral flow between any of the latter versions of $\mathfrak{L}_{\mathfrak{c},z}$ and the version that is defined by taking R very large, $z = R$ and $\mathfrak{c}$ to be $(A_{\bullet 1}, \psi_{\bullet 1})$ as defined using the chosen value for R.

*Part 1*: This part of the subsection defines the pair $\mathfrak{c}(z) \in \mathrm{Conn}(E) \times C^\infty(Y; \mathbb{S})$ for a given $z > c_0$. The definition of $\mathfrak{c}(z)$ on the radius $(c_v{}^4 + 3\, c_v{}^3)\, z^{-1/2}$ tubular neighborhood of any given curve from $\Theta^1$ is given by (A.44) with the integer m coming from the relevant pair in $\Theta^1$.

The four steps that follow define $\mathfrak{c}(z)$ on the complement in Y of the union of the radius $c_v{}^4 z^{-1/2}$ tubular neighborhoods of the curves from $\Theta^1$. By way of a look ahead, Section Aa's construction is used to define $\mathfrak{c}(z)$ on this part of Y.

<u>Step 1</u>: Let $c_v$ denote the constant that is used to define $(A_{\lozenge 11}, \psi_{\lozenge 11})$. Take $z \geq c_0$ and introduce $Y_{*\Lambda}$ to denote the complement in Y of the union of the radius $c_v{}^4 z^{-1/2}$ tubular neighborhoods of the curves from the set $\cup_{p\in\Lambda}\{\hat{\gamma}_p^+ \cup \hat{\gamma}_p^-\}$. Define $T_{*\Lambda}$ to be the subset of $Y - Y_{*\Lambda}$ that consists of the components that do not contain curves from $\Theta^1$. The data consisting of $c_v$, $\rho_* = c_v{}^2 z^{-1/2}$, $T_{*\Lambda}$ and $\Theta = \Theta^0$ supply most but not all of what is needed in Section Aa to define a pair consisting of a Hermitian connection on $E|_{Y_{*\Lambda} \cup T_{*\Lambda}}$ and a section of $\mathbb{S}$ over $Y_{*\Lambda} \cup T_{*\Lambda}$.

The definitions in Section Aa requires the specification of coordinates from Part 4 of Section Aa for each curve in $\Theta^0$. The latter are defined from a chosen isometric isomorphism over each such curve between $K^{-1}$ and the product bundle. Make these choices.

Section Aa also requires isomorphisms between E and the product bundle over certain subsets of $Y_{*\Lambda} \cup T_{*\Lambda}$. These isomorphisms are defined in Step 4. Steps 2 and 3 supply necessary input for the definition in Step 4.

The pair $\mathfrak{c}(z)$ on $Y_{*\Lambda} \cup T_{*\Lambda}$ is the pair that is supplied by Section Aa using the data $c_v$, $\rho_* = c_v{}^2 z^{-1/2}$, $T_{*\Lambda}$, $\Theta = \Theta^0$, the chosen isomorphisms over the curves in $\Theta^0$ between $K^{-1}$ and the product bundle, and the promised isomorphisms between E and the product bundle over the relevant subsets of $Y_{*\Lambda} \cup T_{*\Lambda}$.

Step 2: Appendix Aa introduces an open cover of $Y_{*\Lambda} \cup T_{*\Lambda}$ consisting of a set $U_0$ and a collection of sets $\{U_\gamma\}_{\gamma\in\Theta^0}$. The set $U_0$ is the complement of the union of the radius $c_v^2 z^{-1/2}$ tubular neighborhoods of the curves in $\Theta^0$. Meanwhile, each $\gamma \in \Theta^0$ version of $U_\gamma$ is the radius $4c_v^2 z^{-1/2}$ tubular neighborhood of $\gamma$. The construction of $\mathfrak{c}(z)$ requires an isomorphism between E and the product bundle over $U^0$ and an isomorphism between E and the product bundle over each set from the collection $\{U_\gamma\}_{\gamma\in\Theta^0}$.

Fix $\gamma \in \Theta^0$ to define the isomorphism between E and the product bundle over $U_\gamma$. To do this, remark that the sets $\Theta^0$ and $\Theta^\alpha$ enjoy a 1-1 correspondence with partnered elements being homotopic in $M_\delta \cup (\cup_{p\in\Lambda} \mathcal{H}_p)$. Moreover, the partners intersect $M_\delta$ as arcs that are isotopic via an isotopy that moves points a distance at most $c_0\delta$, this being a consequence of Lemma II.2.5. Let $\upsilon^\alpha$ denote $\gamma$'s partner from $\Theta^\alpha$.

Choose a smoothly embedded, oriented surface in $[0, 1] \times Y$ with the properties listed below.

- *The surface intersects* $[0, c_0^{-1}) \times Y$ as $[0, c_0^{-1}) \times \upsilon^\alpha$
- *The surface intersects* $(c_0^{-1}, 1] \times Y$ as $(c_0^{-1}, 1] \times \gamma$.
- *The surface intersects* $[0,1] \times M_\delta$ *as an embedded rectangle of width less than* $c_0\delta$ *that intersects each constant f surface transversely as a single arc.*
- *The surface intersects the boundary of any radius* $\delta$ *coordinate ball in* $M_\delta$ *transversely as a single arc.*
- *The projection of the surface to* Y *intersects only the* $p \in \Lambda$ *versions of* $\mathcal{H}_p$ *that are crossed by* $\upsilon^\alpha$ *and* $\gamma$*; and its projection in any such* $\mathcal{H}_p$ *is disjoint from* $\hat{\gamma}_p^+$ *and* $\hat{\gamma}_p^-$.

(C.32)

Such a surface can be constructed by mimicking what is done in Step 3 of the proof of Lemma II.5.3 to construct the latter's surface $Z_+$. Use $S_\gamma$ to denote the chosen surface.

Fix $R > c_0 c_v^6 r$ suitable for defining the path $\{(A_{\bullet\tau}, \psi_{\bullet\tau})\}_{\tau\in[0,1]}$ and in any event such that all points in $S_\gamma$ have distance at least $(c_v^4 + 3c_v^3)R^{-1/2}$ from each curve in the set $\cup_{p\in\Lambda}\{\hat{\gamma}_p^+ \cup \hat{\gamma}_p^-\}$. Let $U_S \subset [0,1] \times Y$ denote a tubular neighborhood of $S_\gamma$ that intersects $\{0\} \times Y$ as the radius $4c_v^2 R^{1/2}$ tubular neighborhood of $\upsilon^\alpha$ and that it intersect $\{1\} \times Y$ as $U_\gamma$. Require in addition that points in $U_S$ have distance $(c_v^4 + 3c_v^3)R^{-1/2}$ from $\hat{\gamma}_p^+$ and $\hat{\gamma}_p^-$.

Step 3: Let $\pi$: $[0,1] \times Y \to Y$ denote the projection to the Y factor. As explained in the subsequent paragraphs, the section $\alpha_{\lozenge 11}$ extends over $[0,1] \times Y$ as a section of $\pi^*E$ with zero locus $(\cup_{\gamma\in\Theta^0} S_\gamma) \cup (\cup_{(\gamma,m)\in\Theta^1} [0,1] \times \gamma)$ and with transversal zeros along each $S_\gamma$.

The explanation starts with the 1-cycle $\sum_{\upsilon^\alpha \in \Theta^\alpha} [\upsilon^\alpha] + \sum_{(\gamma,m)\in\Theta^1} m[\gamma]$ where $[\cdot]$ denotes the cycle defined by the fundamental class of the indicated loop. This sum is the weighted sum of the components of the zero locus of $\alpha_{\Diamond 11}$ with the weight of a component being the degree of vanishing of $\alpha_{\Diamond 11}$ on a small radius transverse disk centered on the given component. The class of this cycle in $H_1(Y;\mathbb{Z})$ is Poincaré dual to the first Chern class of E because $\alpha_{\Diamond 11}$ is a section of E.

The first Chern class of E is also Poincaré dual to the class defined by the 1-cycle $\sum_{\gamma\in\Theta^0} [\gamma] + \sum_{(\gamma,m)\in\Theta^1} m[\gamma]$, and as a consequence, the class of the relative 2-cycle $\sum_{\gamma\in\Theta^0} [S_\gamma] + \sum_{(\gamma,m)\in\Theta^1} m[([0,1]\times\gamma]$ on $[0,1]\times Y$ is Poincaré dual to the first Chern class of $\pi^*E$. This being the case, there is a section of $\pi^*E$ whose zero locus defines this same relative 2-cycle. Moreover, there exists such a section with transverse zeros along each $S_\gamma$ and the same local behavior as $\alpha_{\Diamond 11}$ near the origin of any transverse disk in $\{0\}\times Y$ with center on a curve from $\Theta^1$. Use $\hat\alpha$ to denote such a section and use $\hat\alpha|_0$ to denote its restriction to $\{0\}\times Y$. The latter can be written as $u\cdot\alpha_{\Diamond 11}$ with u being a smooth map from the complement in Y of $(\cup_{\upsilon^\alpha\in\Theta^\alpha} \upsilon^\alpha) \cup (\cup_{(\gamma,m)\in\Theta^1} \gamma)$ to $\mathbb{C}-\{0\}$. The section $\alpha_{\Diamond 11}$ has the desired extension if u extends as a map to $\mathbb{C}-\{0\}$ from the complement $[0,1]\times Y$ of $(\cup_{\gamma\in\Theta^0} S_\gamma) \cup (\cup_{(\gamma,m)\in\Theta^1} [0,1]\times\gamma)$.

Let $Y^\alpha$ denote the complement in Y of $(\cup_{\upsilon^\alpha\in\Theta^\alpha} \upsilon^\alpha) \cup (\cup_{(\gamma,m)\in\Theta^1} \gamma)$ and let $X^\alpha$ denote the complement in $[0,1]\times Y$ of $(\cup_{\gamma\in\Theta^0} S_\gamma) \cup (\cup_{(\gamma,m)\in\Theta^1} [0,1]\times\gamma)$. The map u will extend if the restriction homomorphism from $H^1(X^\alpha;\mathbb{Z})$ to $H^1(Y^\alpha;\mathbb{Z})$ is surjective; and this is the case if the inclusion homomorphism from $H_1(Y^\alpha;\mathbb{Z})$/torsion to $H_1(X^\alpha;\mathbb{Z})$/torsion is injective. To prove that this is so, note that its composition with the inclusion homomorphism $H_1(X^\alpha;\mathbb{Z})$ to $H_1([0,1]\times Y;\mathbb{Z})$ is the same as the composition of the homomorphism from $H_1(Y^\alpha;\mathbb{Z})$ to $H_1(Y;\mathbb{Z})$ with the isomorphism given by the push-forward of $\pi$. This understood, the claimed injectivity follows from the fact that the kernel of the inclusion homomorphism from $H_1(Y^\alpha;\mathbb{Z})$ to $H_1(Y;\mathbb{Z})$ is generated by the linking circles of the transverse disks centered on the various curves from $\Theta^1$.

Step 4: Let $\hat\alpha$ now denote an extension of $\alpha_{\Diamond 11}$ to a section of $\pi^*E$ with zero locus $(\cup_{\gamma\in\Theta^0} S_\gamma) \cup (\cup_{(\gamma,m)\in\Theta^1} [0,1]\times\gamma)$ that vanishes transversely along each $S_\gamma$ and is equal to $\pi^*\alpha_{\Diamond 11}$ near each curve from $\Theta^1$. The restriction of this section to $\{1\}\times Y$ is denoted in what follows by $\hat\alpha|_1$. This is a section of E. The required isomorphism over $U_0$ between E and $U_0\times\mathbb{C}$ sends $\hat\alpha|_1$ to its absolute value, $|\hat\alpha|_1|$.

Fix a curve $\gamma\in\Theta^0$. The definition of the required isomorphism between $E|_{U_\gamma}$ and $U_\gamma\times\mathbb{C}$ uses the chosen isomorphism between $K^{-1}|_\gamma$ and $\gamma\times\mathbb{C}$ to define the coordinates (t, z) on $U_\gamma$ from Part 4 of Section Aa. Granted these coordinates, the desired isomorphism over $U_\gamma$ between E and the product bundle takes $\hat{\alpha}|_1$ to $|\hat{\alpha}|_1|\frac{z}{|z|}$.

*Part 2*: This part of the subsection supplies two lemmas that summarize some salient features of the pairs defined in Part 1.

**Lemma C.6**: *There exists* $\kappa > 1$ *and given* $c_v > \kappa$, *there exists* $\kappa_{c_v} > \kappa$ *with the following significance: Suppose that* $r \geq \kappa_{c_v} c_v^{10}$ *and suppose that* $(A, \psi = (\alpha, \beta))$ *is a solution to the* $(r, \mu)$ *version of (1.13) with* $\mu$ *a given element in* $\Omega$ *with* $\mathcal{P}$*-norm smaller than 1. These values of* $c_v$ *and* r *are suitable for defining* $\mathfrak{c}(z)$ *for* $z \geq \kappa_{c_v}$ *given the choice of an isometric isomorphism between* $K^{-1}$ *and the product bundle over each curve from* $\Theta^0$ *and given the choice of a surface of the sort described by (C.32) for each curve from* $\Theta^0$.

- *The resulting* $\mathfrak{c}(z)$ *does not depend on the chosen set of isometric isomorphisms.*
- *The resulting* $\mathfrak{c}(z)$ *depends on the chosen surface and then the extension* $\hat{\alpha}$ *as follows:*
  a) *Respective versions of* $\mathfrak{c}(z)$ *that are defined by different sets of surfaces and extensions differ by the action of a map from* Y *to* $S^1$.
  b) *The homology class of this map defines a class in* $H^1(Y; \mathbb{Z})$ *that is Poincaré dual to a class from the* $\oplus_{p \in \Lambda} H_2(\mathcal{H}_p; \mathbb{Z})$ *summand in (1.4).*

This lemma is proved momentarily.

The next lemma supplies an a priori bound for the absolute value of the spectral flow between versions of $\mathfrak{L}_{\mathfrak{c}(z),z}$ that are defined by distinct choices for $z$.

**Lemma C.7**: *There exists* $\kappa > 1$ *and given* $c_v > \kappa$, *there exists* $\kappa_{c_v} > \kappa$ *with the following significance: Suppose that* $r \geq \kappa_{c_v} c_v^{10}$ *and suppose that* $(A, \psi = (\alpha, \beta))$ *is a solution to the* $(r, \mu)$ *version of (1.13) with* $\mu$ *a given element in* $\Omega$ *with* $\mathcal{P}$*-norm smaller than 1. Use the data from* $\alpha$ *to define* $\mathfrak{c}(z)$ *for* $z \geq \kappa_{c_v}$. *The absolute value of the spectral flow between the* $z = \kappa_{c_v}$ *and any* $z = z_1 \geq \kappa_{c_v}$ *version of* $\mathfrak{L}_{\mathfrak{c}(z),z}$ *is bounded by* $\kappa$.

The proof of this lemma is given directly. The proof assumes that the first bullet of Lemma C.6 is true.

***Proof of Lemma C.7***: Except for one added remark, the proof is identical to that used to prove Proposition C.4. The added remark concerns the use of Lemma A.8 in the proof. In particular, the bounds given for the term $\mathfrak{e}(\mathfrak{f})$ in this lemma depend implicitly on bounds for the functions $\nu$ and $\mu$. Meanwhile, the latter are defined by the coordinates from from Part 4 of Section Aa and thus by the chosen isomorphism over the curve in question between $K^{-1}$ and the product bundle. As Lemma C.6 asserts that $\mathfrak{c}(z)$ does not depend on the chosen isomorphism, choose one with $|\nu| + |\mu| \le c_0$.

***Proof of Lemma C.6***: To prove the first bullet, assume that a choice of surfaces has been made for each curve from $\Theta^0$. Fix $\gamma \in \Theta^0$ and choose an isometric isomorphism between $K^{-1}|_\gamma$ and $\gamma \times \mathbb{C}$ to define $\mathfrak{c}(z)$ on $U_\gamma$. The formulas for $\mathfrak{c}(z)$ are given in (A.8) and (A.9). Granted these formulas, the observations made in the first two paragraphs of Part 5 in Section Ba apply and prove that $\mathfrak{c}(z)$ does not change when the isomorphism changes.

To see about the second bullet, suppose that $\{S_\gamma\}_{\gamma\in\Theta^0}$ and $\{S_\gamma'\}_{\gamma\in\Theta^0}$ are two sets of surfaces of the sort described in (C.32). Let $\varphi_0$ and $\varphi_0'$ denote the corresponding isomorphism between E and the product bundle over $U_0$. Write $\varphi_0'$ as $u_0\varphi_0$ with $u_0$ being a map from $U_0$ to $S^1$. Fix a coordinates from Part 4 of Section Aa for each curve in $\Theta^0$. Given $\gamma \in \Theta^0$, let $\varphi_\gamma$ and $\varphi_\gamma'$ denote the corresponding isomorphisms between E and the product bundle over $U_\gamma$. Write $\varphi_\gamma'$ as $u_\gamma\varphi_\gamma$.

The respective primed and unprimed transition maps from $U_0 \cap U_\gamma$ to $S^1$ that identify the product structure for E over $U_0$ with that over $U_\gamma$ are identical because the same coordinates for $U_\gamma$ are used for the two cases. Use this fact with the formulas in Section Aa to conclude that $u_0 = u_\gamma$ on $U_0 \cup U_\gamma$. This being the case, the collection of maps consisting of $u_0$ and $\{u_\gamma\}_{\gamma\in\Theta^0}$ define a smooth map from Y to $S^1$ that relates the primed and unprimed versions of $\mathfrak{c}(z)$. Let $u\colon Y \to S^1$ denote this map.

Consider now the class defined by u in $H^1(Y;\mathbb{Z})$. This class is determined by the integral of $-\frac{i}{2\pi}u^{-1}du$ over a basis of cycles in Y that generate the free $\mathbb{Z}$-module $H_1(Y;\mathbb{Z})$/torsion. Part 4 of Section 1b describes the set $\{\gamma^{(z)}\}_{z\in ¥}$ of $1 + b_1(M)$ integral curves of $v$ in $M_\delta \cup \mathcal{H}_0$ with the following property: The integral of $-\frac{i}{2\pi}u^{-1}du$ over these cycles detects the image in the summand $H_2(M;\mathbb{Z}) \oplus H_2(\mathcal{H}_0;\mathbb{Z})$ of the Poincaré dual in $H_2(Y;\mathbb{Z})$ of u's cohomology class. To prove that this image is zero, introduce $\hat{\alpha}$ and $\hat{\alpha}'$ to denote the corresponding $\{S_\gamma\}_{\gamma\in\Theta^0}$ and $\{S_\gamma'\}_{\gamma\in\Theta^0}$ extensions of $\alpha_{\lozenge 11}$. Both $\hat{\alpha}$ and $\hat{\alpha}'$ are non-zero on the product of $[0,1]$ with the complement in $M_\delta \cup \mathcal{H}_0$ of the union of the radius $c_0\delta$ tubular neighborhoods of the component segments of $\cup_{\gamma\in\Theta^0}(\gamma \cap M_\delta)$. Use T denote this small radius tubular neighborhood of $\cup_{\gamma\in\Theta^0}(\gamma \cap M_\delta)$. Keep in mind that this set T is disjoint from the set $\cup_{z\in ¥}([0,1] \times \gamma^{(z)})$. This fact can be used to exhibit a homotopy on a neighborhood of $\cup_{z\in ¥}\gamma^{(z)}$ between u and the constant map to $1 \in S^1$: The desired homotopy is parametrized by $[0,1]$ with the $\tau \in [0,1]$ member of the homotopy being the restriction to $\{\tau\} \times ((M_\delta \cup \mathcal{H}_0) - T)$ of $(\hat{\alpha}/|\hat{\alpha}|)(\hat{\alpha}'/|\hat{\alpha}'|)^{-1}$.

*Part 3*: This part of the subsection concerns the spectral flow difference between very large R versions of $\mathfrak{L}_{\mathfrak{c},R}$ as defined using $\mathfrak{c} = (A_{\bullet 1}, \psi_{\bullet 1})$ and the corresponding $z = R$ version of the operator $\mathfrak{L}_{\mathfrak{c}(z),z}$. The proposition that follows says what is needed about this difference.

**Proposition C.8**: *There exists* $\kappa \geq 100$, *and given* $c_v \geq \kappa$, *there exists* $\kappa_{c_v} > \kappa$ *with the following significance: Suppose that* $r \geq \kappa_{c_v} c_v^{10}$ *and suppose that* $(A, \psi = (\alpha, \beta))$ *is a solution to the* $(r, \mu)$ *version of (1.13) with* $\mu$ *a given element in* $\Omega$ *with* $\mathcal{P}$*-norm smaller than 1. The values of* $\kappa$, $c_v$ *and* r *are suitable for defining* $(A_{\Diamond 11}, \psi_{\Diamond 11})$ *and any* $R \geq \kappa c_v^6 r$ *version of* $(A_{\bullet 1}, \psi_{\bullet 1})$. *Fix any sufficiently large* R *and use it define the pair* $(A_{\bullet 1}, \psi_{\bullet 1})$. *Define the* $z = R$ *version of* $\mathfrak{c}(z)$ *using any chosen set of isomorphisms between* $K^{-1}$ *and the product bundle over the curves from* $\Theta^0$*; and using any chosen set of surfaces of the sort described in (C.32) for the curves in* $\Theta^0$. *The norm of the difference between the respective values of the spectral flow function* $\mathfrak{f}_s$ *at* $(A_{\bullet 1}, \psi_{\bullet 1})$ *and at the* $z = R$ *version of* $\mathfrak{c}(z)$ *is bounded by* $\kappa$.

By way of a parenthetical remark, what is said in the second bullet of Lemma C.6 is consistent with what is said in Proposition C.8. This follows from three facts. Here is the first: The function $\mathfrak{f}_s$ is invariant under the action on Conn(E) $\times$ $C^\infty(Y; \mathbb{S})$ of the subgroup of maps from Y to $S^1$ whose corresponding class in $H^1(Y; \mathbb{Z})$ has zero cup product with the first Chern class of the line bundle det($\mathbb{S}$). The second fact concerns the cup product pairing between this first Chern class and a given class $\sigma \in H^1(Y; \mathbb{Z})$: This is the same as the pairing between the first Chern class of det($\mathbb{S}$) and the Poincaré dual of $\sigma$ in $H_2(Y; \mathbb{Z})$. Here is the final fact: The first Chern class of det($\mathbb{S}$) annihilates the $\oplus_{p \in \Lambda} H^2(\mathcal{H}_p; \mathbb{Z})$ summand of $H_2(Y; \mathbb{Z})$.

***Proof of Proposition C.8***: If R is sufficiently large, then the arguments that are used in Section 2b of [T?E=S3] can be imported with only notational changes to prove the proposition.

## e) Proof of Proposition 2.6

This section gives the Proposition 2.6. The argument has five steps.

Step 1: It is convenient to choose a finite set of surrogates for the pair $(A_E, \psi_E)$. This set of surrogates is indexed by the set of all possible pairs of the form $(\Theta^0, \Theta^1)$ that can arise in the previous subsection from large r and $(r, \mu)$ version of (1.13). This indexing set is denoted by $\mathcal{Z}_{HF} \times \mathcal{Z}^1$.

By way of a precise definition, the set $\mathcal{Z}_{HF}$ has distinct elements of the following sort: Let $\Theta^0$ denote a given element. This set $\Theta^0$ consists of at most G distinct, closed integral curves of $v$. All curves in the set $\Theta^0$ lie in $M_\delta \cup (\cup_{p \in \Lambda} \mathcal{H}_p)$ and none are from $\cup_{p \in \Lambda}\{\hat{\gamma}_p^+ \cup \hat{\gamma}_p^-\}$. The union of the curves from $\Theta^0$ intersects $M_\delta$ as G segments that give the same pairing of the index 1 and index 2 critical points of $f$. Finally, the union of the integral curves in $v$ intersect each $\mathfrak{p} \in \Lambda$ version of $\mathcal{H}_p$ as a single segment that runs from the $u < 0$ boundary of $\mathcal{H}_p$ to the $u > 0$ boundary. The intersection is characterized by an integer, $\mathfrak{k}_p$, as in Proposition II.2.7, and the segment in question has $\mathfrak{k}_p = 0$.

As explained in [KLTI] and [KLTII], the set $\mathcal{Z}_{HF}$ determines a set of generators for the Heegaard-Floer homology on M. In any event, $\mathcal{Z}_{HF}$ has finitely many elements.

The set $\mathcal{Z}^1$ consists of elements of the following sort: Let $\Theta^1$ denote a given element. This set $\Theta^1$ consists of pairs of the form $(\gamma, m)$ where $\gamma \in \cup_{p \in \Lambda}\{\hat{\gamma}_p^+ \cup \hat{\gamma}_p^-\}$ and where m is a positive integer. No two pair share the same integral curve component. The integer m is bounded by $c_0 c_v^3$. The set $\Theta^1$ is also finite.

Take each $\Theta^0$ in $\mathcal{Z}_{HF}$ and assign once and for all an isometric isomorphism between $K^{-1}$ and the product bundle over each curve from $\Theta_0$. Let $\kappa_{c*}$ denote the larger of the versions of $\kappa_{c_v}$ that appear in Lemmas C.6 and C.7. Fix $z_0 = \kappa_{c*}^2$ and assign once and for all a product structure for E over the radius $4c_v^2 z_0^{-1/2}$ tubular neighborhood of each curve from $\Theta_0$. Fix once and for all a product structure for E over the complement in Y of the union of the radius $c_v^4 z_0^{-1/2}$ tubular neighborhoods of the curves from $\Theta^0$.

Take each pair $\hat{o} = (\Theta^0, \Theta^1) \in \mathcal{Z}_{HF} \times \mathcal{Z}^1$ and use the data $c_v$, $z = z_0$ with the product structures chosen in the preceding paragraph to construct the corresponding version of the pair $\mathfrak{c}(z = z_0)$ as instructed in Step 1 of Part 1 of Section Cd. Denote this pair by $\mathfrak{c}_{\hat{o}}$. Write this pair as $(A_{\hat{o}}, \psi_{\hat{o}})$ and use $\alpha_{\hat{o}}$ to denote the E summand component of $\psi_{\hat{o}}$.

Since $\mathcal{Z}_{HF} \times \mathcal{Z}^1$ is a finite set, there exists a purely $c_v$ dependent $\kappa_c > 1$ with the following property: Fix $\hat{o} \in \mathcal{Z}_{HF} \times \mathcal{Z}^1$. Then the connection $A_{\hat{o}}$ can be written as $A_E + \hat{a}_{\hat{o}}$ with $\hat{a}_{\hat{o}}$ being an $i\mathbb{R}$ valued 1-form with $|\hat{a}_{\hat{o}}| \le \kappa_c$. In addition, the norm of the difference between the respective values of $\mathfrak{f}_{\mathfrak{s}}$ at $(A_E, \psi_E)$ and $\mathfrak{c}_{\hat{o}}$ is bounded by $\kappa_c$.

<u>Step 2</u>: Fix $c_v \ge c_0$ and $r \ge \kappa_c c_v^{10}$ with $\kappa_c$ being a purely $c_v$ dependent constant. Assume that $c_v$ and $\kappa_c$ are suitably for invoking the results in Appendix A, B and the previous subsections of this Appendix C. Suppose that $(A, \psi = (\alpha, \beta))$ is a solution to the $(r, \mu)$ version of (1.13) with $\mu$ a given element in $\Omega$ with $\mathcal{P}$-norm smaller than 1. Assume in addition that $|x_{\mathbb{S}}(A)| \le c_0$. By assumption, the values of $c_v$ and r are suitable for defining from $(A, \psi)$ the pair $(A_{\Diamond 11}, \psi_{\Diamond 11})$ and any given $R \ge \kappa c_v^6 r$ version of $(A_{\bullet 1}, \psi_{\bullet 1})$. Fix any sufficiently large R and use it define both $(A_{\bullet 1}, \psi_{\bullet 1})$ and the $z \in (\kappa_{c*}, R]$ versions of $\mathfrak{c}(z)$. Use $\mathfrak{c}_{(A,\psi)}(z)$ to denote such a version. Write the pair $\mathfrak{c}_{(A,\psi)}(z)$ as $(A_z, \psi_z)$ and write the E summand of $\psi_z$ as $\alpha_z$. Use $\hat{o} \in \mathcal{Z}_{HF} \times \mathcal{Z}^1$ in what follows to denote the element $(\Theta^0, \Theta^1)$ that is used to define $\mathfrak{c}_{(A,\psi)}(z)$. This element is determined by $(A, \psi)$.

<u>Step 3</u>: Let T denote the union of the radius $c_0 \delta$ tubular neighborhoods of the intersection between $M_\delta$ and the curves from $\Theta^0$. This set T has distance at least $c_0^{-1}$ from the curves in the set $\{\gamma^{(z)}\}_{z \in ¥}$. Moreover, it contains the $M_\delta$ part of the zero locus of $\alpha$ and the $M_\delta$ part of the zero locus of $\alpha_{\hat{o}}$. With this understood, the section $\alpha$ on $(M_\delta \cup \mathcal{H}_0)$–T can be written as $\alpha = |\alpha| u \alpha_{\hat{o}}$ with u being a map from $(M_\delta \cup \mathcal{H}_0)$–T to $S^1$. Write $\nabla_A \alpha$ on $(M_\delta \cup \mathcal{H}_0)$–T as $(d|\alpha| + (u^{-1}du + \hat{a}_A - \hat{a}_{\hat{o}})|\alpha|)\, \alpha_{\hat{o}}$.

Lemma 2.1 asserts that $1 - |\alpha| \le c_0 r^{-1}$ and that $|\nabla_A \alpha| \le c_0$ on $(M_\delta \cup \mathcal{H}_0)$–T. Given that $A = A_E + \hat{a}_A$ and $A_{\hat{o}} = A_E + \hat{a}_{\hat{o}}$, it follows that $|u^{-1}du + \hat{a}_A - \hat{a}_{\hat{o}}| \le \kappa_c$ on $(M_\delta \cup \mathcal{H}_0)$–T with $\kappa_c \ge 1$ being a purely $c_v$ dependent constant. The latter bound has the following consequence: The absolute value of integral of $-\frac{i}{2\pi} u^{-1}du$ over any curve from the set $\{\gamma^{(z)}\}_{z \in ¥}$ is bounded by $\kappa_c$ with $\kappa_c$ again denoting a purely $c_v$ dependent constant.

Fix $z \in (\kappa_{c*}, R)$. The zero locus of $\alpha_z$ in $M_\delta$ also lies in T. This understood, write $\alpha^z$ on $(M_\delta \cup \mathcal{H}_0)$–T as $\alpha_z = u_z |\alpha|^{-1} \alpha$ with $u_z$ being a smooth map to $S^1$. It follows from what is said by Step 4 of Part 3 in Section Ca and by Step 4 of Part 1 in Section Cd that the integral of the 1-form $-\frac{i}{2\pi} u_z^{-1} du_z$ is zero over any curve from the set $\{\gamma^{(z)}\}_{z \in ¥}$. In fact, the extension $\hat{\alpha}$ used in Step 4 of Part 1 in Section 1c can be chosen so that $u_z = 1$.

<u>Step 4</u>: The zero locus of $\alpha_z$ and that of $\alpha_{\hat{o}}$ are identical, it being the union of the curves from $\Theta^0$ and $\Theta^1$. The latter fact implies that $\alpha_z$ can be written on the complement of this zero locus as $\alpha_z = |\alpha_z| |\alpha_{\hat{o}}|^{-1} \hat{u}_z \alpha_{\hat{o}}$ with $\hat{u}_z$ being a smooth map to $S^1$ from the complement in Y of the union of the curves from $\Theta^0$ and $\Theta^1$. It follows from what was said in Step 3 that the integral of $-\frac{i}{2\pi} \hat{u}_z^{-1} d\hat{u}_z$ is zero over any curve from the set $\{\gamma^{(z)}\}_{z \in ¥}$.

Take $z = z_0$ now, this being the value of $z$ that is used to define $(A_{\hat{o}}, \psi_{\hat{o}})$. It follows from the first bullet of Lemma C.6 that the $z = z_0$ version of $\hat{u}_z$ extends to define a smooth map from the whole of Y to $S^1$ and that the $z = z_0$ version of the pair $(A_z, \psi_z)$ can be written as $(A_{z0} = A_{\hat{o}} - \hat{u}_{z0}^{-1} d\hat{u}_{z0}, \psi_{z0} = \hat{u}_{z0} \psi_{\hat{o}})$ on the whole of Y.

<u>Step 5</u>: The integral of $-\frac{i}{2\pi} \hat{u}_{z0}^{-1} d\hat{u}_{z0}$ over the curves from $\{\gamma^{(z)}\}_{z \in ¥}$ is zero, and this implies that the Poincaré dual in $H_2(Y; \mathbb{Z})$ of the class in $H^1(Y; \mathbb{Z})$ defined by $\hat{u}_{z0}$ lies in the $\oplus_{p \in \Lambda} H_2(\mathcal{H}_p; \mathbb{Z})$ summand of $H_2(Y; \mathbb{Z})$. As noted previously, this summand has zero pairing with the first Chern class of $\det(\mathbb{S})$. It follows as a consequence that the spectral flow function $\mathfrak{f}_{\mathfrak{s}}$ has the same value on $\mathfrak{c}_{(A,\psi)}(z_0)$ as it has on $\mathfrak{c}_{\Diamond}$. This being so, the absolute value of $\mathfrak{f}_{\mathfrak{s}}$ on $\mathfrak{c}_{(A,\psi)}(z_0)$ is bounded by $\kappa_c$ with $\kappa_c$ being a purely $c_v$ dependent constant.

Lemma C.7 asserts that the norm of the difference between the values of $\mathfrak{f}_{\mathfrak{s}}$ at $\mathfrak{c}_{(A,\psi)}(z_0)$ and at $\mathfrak{c}_{(A,\psi)}(z = R)$ is bounded by a purely $c_v$ dependent constant, and so the absolute value of $\mathfrak{f}_{\mathfrak{s}}$ at $\mathfrak{c}_{(A,\psi)}(z=R)$ is also bounded by such a constant. Proposition C.8 asserts that the norm of the difference between values of $\mathfrak{f}_{\mathfrak{s}}$ at $\mathfrak{c}_{(A,\psi)}(z=R)$ and $(A_{\bullet 1}, \psi_{\bullet 1})$ is also bounded by a purely $c_v$ dependent constant. Proposition C.3 asserts that such is also the case for the norm of the difference between the values of $\mathfrak{f}_{\mathfrak{s}}$ at $(A_{\bullet 1}, \psi_{\bullet 1})$ and $(A_{\Diamond 11}, \psi_{\Diamond 11})$. Proposition C.2 says the same thing for the norm of the difference between values of $\mathfrak{f}_{\mathfrak{s}}$ at $(A_{\Diamond 11}, \psi_{\Diamond 11})$ and at $(A_{\Diamond}, \psi_{\Diamond})$. Proposition B.13 says this about the norm of the difference between values of $\mathfrak{f}_{\mathfrak{s}}$ at $(A_{\Diamond}, \psi_{\Diamond})$ and $(A_*, \psi_*)$; and Proposition B.3 says this about the norm of the difference between the values of $\mathfrak{f}_{\mathfrak{s}}$ at $(A, \psi)$ and at $(A_*, \psi_*)$.

Adding all of these absolute value bounds to verifies that the absolute value of $\mathfrak{f}_{\mathfrak{s}}$ at $(A, \psi)$ is bounded by a purely $c_v$ dependent constant.

Cagatay Kutluhan; Department of Mathematics, Harvard University, Cambridge MA 02138. *Email address*: kutluhan@math.harvard.edu.

Yi-Jen Lee; Department of Mathematics, Purdue University, West Lafayette, IN 47907 *Email address*: yjlee@math.purdue.edu.

Clifford Henry Taubes; Department of Mathematics, Harvard University, Cambridge, MA 02138. *Email address*: chtaubes@math.harvard.edu.